\documentclass[12pt,a4paper]{amsart}
\usepackage{amssymb,eucal}
\usepackage{mathrsfs}
\usepackage[cmtip,all]{xy}
\usepackage{hyperref}

\calclayout

\newcommand{\+}{\protect\nobreakdash-}
\renewcommand{\:}{\colon}

\newcommand{\rarrow}{\longrightarrow}
\newcommand{\larrow}{\longleftarrow}
\newcommand{\ot}{\otimes}
\newcommand{\ocn}{\odot}
\newcommand{\st}{\star}

\newcommand{\lrarrow}{\mskip.5\thinmuskip\relbar\joinrel\relbar\joinrel
 \rightarrow\mskip.5\thinmuskip\relax}
\newcommand{\llarrow}{\mskip.5\thinmuskip\leftarrow\joinrel\relbar
 \joinrel\relbar\mskip.5\thinmuskip\relax}

\newcommand{\hth}{\mskip 0.5\thinmuskip}

\newcommand{\bu}{{\text{\smaller\smaller$\scriptstyle\bullet$}}}

\DeclareFontFamily{U}{mathx}{\hyphenchar\font45}
\DeclareFontShape{U}{mathx}{m}{n}{
      <5> <6> <7> <8> <9> <10>
      <10.95> <12> <14.4> <17.28> <20.74> <24.88>
      mathx10
      }{}
\DeclareSymbolFont{mathx}{U}{mathx}{m}{n}
\DeclareFontSubstitution{U}{mathx}{m}{n}
\DeclareMathAccent{\widecheck}{0}{mathx}{"71}

\DeclareMathAlphabet{\matholdswab}{U}{yswab}{m}{n}   

\DeclareMathOperator{\Spec}{Spec}
\DeclareMathOperator{\Spf}{Spf}
\DeclareMathOperator{\Hom}{Hom}
\DeclareMathOperator{\Ext}{Ext}
\DeclareMathOperator{\Tot}{Tot}
\DeclareMathOperator{\Mor}{Mor}

\DeclareMathOperator{\fHom}{\mathfrak{Hom}}
\DeclareMathOperator{\Cohom}{\mathfrak{Cohom}}
\DeclareMathOperator{\cHom}{\mathcal H \mskip-.3\thinmuskip
  \text{\rmfamily\mdseries\fontshape{ui}\selectfont om}}

\DeclareMathOperator{\coker}{coker}
\DeclareMathOperator{\rsd}{rsd}
\DeclareMathOperator{\crd}{crd}

\newcommand{\Modl}{{\operatorname{\mathsf{--Mod}}}}
\newcommand{\modl}{{\operatorname{\mathsf{--mod}}}}
\newcommand{\Qcoh}{{\operatorname{\mathsf{--Qcoh}}}}
\newcommand{\Tors}{{\operatorname{\mathsf{--Tors}}}}
\newcommand{\Contra}{{\operatorname{\mathsf{--Contra}}}}
\newcommand{\Sh}{{\operatorname{\mathsf{--Sh}}}}
\newcommand{\Cosh}{{\operatorname{\mathsf{--Cosh}}}}
\newcommand{\Lcth}{{\operatorname{\mathsf{--Lcth}}}}
\newcommand{\Ctrh}{{\operatorname{\mathsf{--Ctrh}}}}
\newcommand{\Proqc}{{\operatorname{\mathsf{--Proqc}}}}
\newcommand{\Indct}{{\operatorname{\mathsf{--Indct}}}}
\newcommand{\Indlc}{{\operatorname{\mathsf{--Indlc}}}}

\newcommand{\dtors}{{\operatorname{\mathsf{-tors}}}}
\newcommand{\dctra}{{\operatorname{\mathsf{-ctra}}}}

\newcommand{\ct}{{\operatorname{\mathrm{-ct}}}}
\newcommand{\qc}{{\operatorname{\mathrm{-qc}}}}

\newcommand{\Ab}{\mathsf{Ab}}
\newcommand{\Hot}{\mathsf{Hot}}
\newcommand{\Ac}{\mathsf{Ac}}
\newcommand{\Com}{\mathsf{Com}}

\newcommand{\AffSch}{\mathsf{AffSch}}
\newcommand{\CSch}{\mathsf{CSch}}
\newcommand{\NoethAffSch}{\mathsf{NoethAffSch}}
\newcommand{\NoethFormSch}{\mathsf{NoethFormSch}}
\newcommand{\NoethSch}{\mathsf{NoethSch}}
\newcommand{\Sch}{\mathsf{Sch}}
\newcommand{\Sets}{\mathsf{Sets}}

\newcommand{\ilim}
 {\mathop{\text{\normalfont``$\varinjlim$''\!\!}}\nolimits}

\newcommand{\sop}{\mathsf{op}}
\newcommand{\tors}{\mathsf{tors}}
\newcommand{\ctra}{\mathsf{ctra}}
\newcommand{\inj}{\mathsf{inj}}
\newcommand{\prj}{\mathsf{prj}}
\newcommand{\fl}{\mathsf{fl}}
\newcommand{\vfl}{\mathsf{vfl}}
\newcommand{\cta}{\mathsf{cta}}
\newcommand{\coa}{\mathsf{coa}}
\renewcommand{\cot}{\mathsf{cot}}
\newcommand{\lct}{\mathsf{lct}}
\newcommand{\lin}{\mathsf{lin}}
\newcommand{\qs}{\mathsf{qs}}
\newcommand{\alf}{\mathsf{alf}}
\newcommand{\al}{\mathsf{al}}
\newcommand{\fq}{\mathsf{fq}}
\newcommand{\cfq}{\mathsf{cfq}}
\newcommand{\bb}{\mathsf{b}}
\newcommand{\abs}{\mathsf{abs}}
\newcommand{\co}{\mathsf{co}}
\newcommand{\ctr}{\mathsf{ctr}}
\newcommand{\bco}{\mathsf{bco}}
\newcommand{\bctr}{\mathsf{bctr}}
\newcommand{\dil}{\mathsf{dil}}
\newcommand{\cdl}{\mathsf{cdl}}

\newcommand{\fX}{\mathfrak X}
\newcommand{\fY}{\mathfrak Y}
\newcommand{\fZ}{\mathfrak Z}
\newcommand{\fU}{\mathfrak U}
\newcommand{\fV}{\mathfrak V}
\newcommand{\fW}{\mathfrak W}
\newcommand{\fO}{\mathfrak O}
\newcommand{\fP}{\mathfrak P}
\newcommand{\fQ}{\mathfrak Q}
\newcommand{\fR}{\mathfrak R}
\newcommand{\fS}{\mathfrak S}
\newcommand{\fT}{\mathfrak T}
\newcommand{\fI}{\mathfrak I}
\newcommand{\fJ}{\mathfrak J}
\newcommand{\fK}{\mathfrak K}
\newcommand{\fL}{\mathfrak L}
\newcommand{\fM}{\mathfrak M}
\newcommand{\fN}{\mathfrak N}
\newcommand{\fA}{\mathfrak A}
\newcommand{\fB}{\mathfrak B}
\newcommand{\fC}{\mathfrak C}
\newcommand{\fE}{\mathfrak E}
\newcommand{\fF}{\mathfrak F}
\newcommand{\fG}{\mathfrak G}
\newcommand{\fH}{\mathfrak H}

\newcommand{\ff}{\mathfrak f}
\newcommand{\fg}{\mathfrak g}
\newcommand{\fh}{\mathfrak h}
\newcommand{\fii}{\mathfrak i}
\newcommand{\fj}{\mathfrak j}
\newcommand{\fk}{\mathfrak k}
\newcommand{\fm}{\mathfrak m}
\newcommand{\fp}{\mathfrak p}
\newcommand{\fqq}{\mathfrak q}

\newcommand{\br}{\mathbf r}
\newcommand{\bs}{\mathbf s}
\newcommand{\bt}{\mathbf t}

\newcommand{\bB}{\mathbf B}
\newcommand{\bC}{\mathbf C}
\newcommand{\bD}{\mathbf D}
\newcommand{\bW}{\mathbf W}
\newcommand{\bT}{\mathbf T}

\newcommand{\cA}{\mathcal A}
\newcommand{\cB}{\mathcal B}
\newcommand{\cC}{\mathcal C}
\newcommand{\cD}{\mathcal D}
\newcommand{\cE}{\mathcal E}
\newcommand{\cF}{\mathcal F}

\newcommand{\cI}{\mathcal I}
\newcommand{\cJ}{\mathcal J}
\newcommand{\cK}{\mathcal K}
\newcommand{\cL}{\mathcal L}
\newcommand{\cM}{\mathcal M}
\newcommand{\cN}{\mathcal N}
\newcommand{\cO}{\mathcal O}
\newcommand{\cP}{\mathcal P}
\newcommand{\cQ}{\mathcal Q}

\newcommand{\rB}{\mathscr B}

\newcommand{\rF}{\mathscr F}
\newcommand{\rG}{\mathscr G}
\newcommand{\rH}{\mathscr H}
\newcommand{\rK}{\mathscr K}
\newcommand{\rL}{\mathscr L}
\newcommand{\rP}{\mathscr P}
\newcommand{\rQ}{\mathscr Q}
\newcommand{\rR}{\mathscr R}

\newcommand{\gJ}{\mathalpha{\text{\small$\matholdswab J$}}}
\newcommand{\gK}{\matholdswab K}
\newcommand{\gL}{\matholdswab L}
\newcommand{\gP}{\matholdswab P}
\newcommand{\gQ}{\matholdswab Q}

\newcommand{\rop}{\mathrm{op}}
\newcommand{\id}{\mathrm{id}}
\newcommand{\Id}{\mathrm{Id}}
\newcommand{\mx}{\mathrm{mx}}
\newcommand{\mn}{\mathrm{mn}}

\newcommand{\sA}{\mathsf A}
\newcommand{\sB}{\mathsf B}
\newcommand{\sC}{\mathsf C}
\newcommand{\sD}{\mathsf D}
\newcommand{\sE}{\mathsf E}
\newcommand{\sF}{\mathsf F}
\newcommand{\sJ}{\mathsf J}
\newcommand{\sK}{\mathsf K}
\newcommand{\sP}{\mathsf P}

\newcommand{\boZ}{\mathbb Z}
\newcommand{\boQ}{\mathbb Q}
\newcommand{\boL}{\mathbb L}
\newcommand{\boR}{\mathbb R}

\newcommand{\Section}[1]{\bigskip\section{#1}\medskip}
\theoremstyle{plain}
\newtheorem{thm}{Theorem}[subsection]
\newtheorem{mainthm}[thm]{Main Theorem}
\newtheorem{prop}[thm]{Proposition}
\newtheorem{lem}[thm]{Lemma}
\newtheorem{cor}[thm]{Corollary}
\theoremstyle{definition}
\newtheorem{rem}[thm]{Remark}
\newtheorem{rems}[thm]{Remarks}

\newtheorem{ex}[thm]{Example}
\newtheorem{exs}[thm]{Examples}

\begin{document}

\title{Contraherent cosheaves of contramodules \\
on Noetherian formal schemes}

\author{Leonid Positselski}

\address{Institute of Mathematics, Czech Academy of Sciences \\
\v Zitn\'a~25, 115~67 Prague~1 \\ Czech Republic}

\email{positselski@math.cas.cz}

\begin{abstract}
 There are two main results in this paper: the MGM equivalence between
the conventional derived categories of abelian/exact categories, and
the derived co-contra correspondence (i.~e., an equivalence between
the coderived and the contraderived category).
 Firstly, for any Noetherian formal scheme that is either of finite
Krull dimension or semi-separated, we construct a triangulated
equivalence between the conventional unbounded derived categories of
the abelian category of quasi-coherent torsion sheaves and the exact
category of contraherent cosheaves of contramodules.
 A natural dedualizing complex of quasi-coherent torsion sheaves is
constructed and used for that purpose.
 Secondly, for any Noetherian formal scheme with a dualizing complex
of quasi-coherent torsion sheaves, we construct a triangulated
equivalence between the coderived category of quasi-coherent torsion
sheaves and the contraderived category of contraherent cosheaves of
contramodules.
 The exact category of contraherent cosheaves of contramodules
on a locally Noetherian formal scheme, as well as the exact categories
of locally contraherent cosheaves of contramodules (with respect to
a given open covering) are defined and discussed in this paper.
 The exposition in the section of preliminaries in adic commutative
algebra is worked out in the greater generality of arbitrary commutative
rings with adic topologies (of finitely generated ideals).
\end{abstract}

\maketitle

\tableofcontents

\section{Introduction}
\medskip

\subsection{{}}
 One of the main original goals of the theory of contraherent cosheaves,
developed in~\cite{Pcosh,Pdomc} and discussed in~\cite{Pphil}, was to
globalize contramodules over formal schemes and ind-schemes.
 In the present paper, we accomplish this goal in the modest generality
of (locally) Noetherian formal schemes with adic topologies on
the rings of functions on affine open formal subschemes.
 These are the formal schemes in the sense of Hartshorne's
textbook~\cite[Section~II.9]{HarAG}, as distinguished from the more
general approaches of Grothendieck's
EGAI~\cite[Sections~0.7 and~I.10]{EGAI} and
\emph{The Stacks project}~\cite[Chapter Tag~0AHW]{SP}.

 Why \emph{contraherent cosheaves of contramodules}?
 Why not \emph{quasi-coherent sheaves of contramodules}?
 This question was raised in~\cite[Section~2]{Pphil} and
answered in~\cite[Section~1.10]{Pcs}, with a definitive counterexample
in~\cite[Remark~10.7 and Example~10.8]{Pcs}.
 Simply put, the localization functors in contramodule categories
do \emph{not} have the requisite exactness properties (which they have
in the usual module categories).
 The exactness properties of the colocalization functors in
the contramodule categories are better.
 Another explanation is that one wants to work with infinite direct
products of contramodules, and the localization functors do not
commute with direct products.
 The main results of the paper~\cite{Pcs} establish the existence,
an explicit construction, and the exactness properties of
the colocalization functors in contramodule categories in a generality
far surpassing the needs of the present paper.

\subsection{{}} \label{introd-sections-2-3}
 This paper starts with a section ``Preliminaries in adic commutative
algebra'', where we discuss the affine geometry of contramodules and
quotseparated contramodules in the generality of commutative rings
with the adic topology of a finitely generated ideal.
 Then follows a section
(Section~\ref{contraherent-cosheaves-of-contramods-secn}) where we
discuss cosheaves of contramodules and define the exact category of
locally contraherent cosheaves of contramodules on a locally Noetherian
formal scheme $\fX$ with respect to its open covering~$\bW$.
 This construction builds upon the foundations of the contraherent
cosheaf theory developed in the book manuscript~\cite{Pcosh}.

 We also consider the dual-analogous abelian category of quasi-coherent
torsion sheaves on $\fX$, which is well-known in
the literature~\cite{AJL,AJPV}.
 A more general construction of the abelian category of quasi-coherent
torsion sheaves on a ``reasonable'' ind-scheme was suggested
in~\cite[Section~7.11.4]{BD2} and developed in detail
in~\cite[Chapter~2]{Psemten}.

 The definitions of the direct and inverse images of quasi-coherent
torsion sheaves and locally contraherent cosheaves of contramodules
under morphisms of locally Noetherian formal schemes are worked out
further along in
Section~\ref{contraherent-cosheaves-of-contramods-secn}.
 Specifically, we define the direct images of locally contraherent
cosheaves of contramodules under suitable affineness conditions, and
the inverse images of locally contraherent cosheaves of contramodules
under suitable tightness and (very) flatness conditions on
the morphism of formal schemes $\ff\:\fY\rarrow\fX$.
 The locality properties of the classes of tight and (very) flat
morphisms of formal schemes are also established, based on an extensive
discussion in the preliminaries.

 Then we construct and discuss the functors of contraherent
$\fHom_\fX(\cM,\cJ)$ and contratensor product $\cM\ocn_\fX\fP$ for
a quasi-coherent torsion sheaf $\cM$, an injective quasi-coherent
torsion sheaf $\cJ$, and a cosheaf of $\fO_\fX$\+modules $\fP$ on
a locally Noetherian formal scheme~$\fX$.
 Here $\fHom_\fX(\cM,\cJ)$ is a locally cotorsion globally contraherent
cosheaf of contramodules on $\fX$, while $\cM\ocn_\fX\fP$ is
a quasi-coherent torsion sheaf on~$\fX$.
 Notice that the construction of the contraherent $\fHom$ functor was
worked out in several alternative settings in~\cite[Section~2.5]{Pcosh}
and~\cite[Section~5.5]{Pdomc} (see the discussion
in~\cite[Sections~5.6 and~8.7]{Pphil} and~\cite[Section~0.10]{Pcosh}).
 Over a formal scheme, there are \emph{no} well-behaved notions of
either flat/very flat quasi-coherent torsion sheaves, or of
contraadjusted/cotorsion quasi-coherent torsion sheaves.
 Hence we are left with only one such setting, viz., that of
injective quasi-coherent torsion sheaves.

\subsection{{}}
 A discussion of injective quasi-coherent torsion sheaves and
projective contraherent cosheaves of contramodules, flat and antilocally
flat contraherent cosheaves of contramodules, flasque quasi-coherent
torsion sheaves and coflasque contraherent cosheaves of contramodules,
dilute quasi-coherent torsion sheaves and codilute contraherent 
cosheaves of contramodules is offered in
Section~\ref{injective-projective-secn}.
 Notice that, \emph{unlike} in~\cite{Pcosh} and~\cite{Pdomc}, we do
\emph{not} consider antilocal contraherent cosheaves in the present
paper, as this concept seems to be based on the notion of a locally
injective contraherent cosheaf, which is \emph{not} applicable to
contraherent cosheaves of contramodules (as there are usually
\emph{no} nonzero injective objects in contramodule categories).
 So we introduce and discuss \emph{codilute} contraherent cosheaves
of contramodules instead, dualizing a concept due to
Murfet~\cite[Section~4.1]{M-n}.

 We extend to locally Noetherian formal schemes Hartshorne's
classification of injective quasi-coherent sheaves on locally
Noetherian schemes~\cite[Section~II.7]{HarRD}.
 A comparable result can be found in Yekutieli's
paper~\cite[Proposition~4.2]{Yek0}.
 Both Hartshorne's classification and our generalization of it are
based on Matlis' classification of injective modules over Noetherian
commutative rings~\cite[Section~3]{Mat}.
 Furthermore, we extend to locally Noetherian formal schemes our
classification of projective locally cotorsion contraherent cosheaves
on locally Noetherian schemes~\cite[Section~6.1]{Pcosh}.
 Both in~\cite{Pcosh} and in the present paper, our approach to such
a classification is based on Enochs' classification of flat cotorsion
modules over Noetherian commutative rings~\cite[Section~2]{En2}.

 Section~\ref{injective-projective-secn} ends with a discussion of
conventional and exotic derived categories of quasi-coherent torsion
sheaves and locally contraherent cosheaves of contramodules, and
the constructions of derived direct image functors acting between
various kinds of derived categories.
 For the inverse images, we do not have enough adjusted objects, so
we restrict ourselves to a discussion of underived (but triangulated)
induced functors of inverse image acting between some derived
categories under suitable assumptions.

\subsection{{}}
 The main results of the paper are stated and proved in
Sections~\ref{dedualizing-secn} and~\ref{dualizing-secn}.
 The now-classical topic of \emph{Matlis--Greenlees--May}
(\emph{MGM}) \emph{duality}, going back to the seminal
papers~\cite{Mat2,GM}, was developed in the subsequent
publications~\cite{AJL0,AJL2,DG,PSY}.
 In particular, the paper~\cite{AJL0} offered an approach to
MGM duality over schemes, and \cite{AJL2}, over formal schemes.
 The approach in the present paper is different.

 The traditional approach to the MGM duality/equivalence involves
an ambient triangulated category $\sD$ with two triangulated
endofunctors $\gamma$ and $\lambda\:\sD\rarrow\sD$.
 The functor $\gamma\:\sD\rarrow\sD$ is some version of a derived adic
torsion functor, while $\lambda$~is a derived adic completion functor.
 The category $\sD$ can be the unbounded derived category
$\sD=\sD(R\Modl)$ of modules over a commutative ring $R$ (in
the context of a commutative ring with a finitely generated ideal
$I\subset R$), or some version of a quasi-coherent derived category
of a scheme or a formal scheme.

 Moreover, the functor $\gamma\:\sD\rarrow\sD$ is an \emph{idempotent
comonad} on $\sD$, i.~e., $\gamma$~is endowed with a natural
transformation $c\:\gamma\rarrow\Id_\sD$ from~$\gamma$ into the identity
functor $\Id_\sD$ such that the two compositions $\gamma\circ c$ and
$c\circ\gamma\:\gamma^2\rarrow\gamma$ are both isomorphisms and equal
to each other.
 Dual-analogously, the functor $\lambda\:\sD\rarrow\sD$ is
an \emph{idempotent monad} on $\sD$, i.~e., $\lambda$~is endowed with
a natural transformation $u\:\Id_\sD\rarrow\lambda$ such that the two
compositions $\lambda\circ u$ and $u\circ\lambda\:\lambda\rarrow
\lambda^2$ are both isomorphisms and equal to each other.
 The key observation is that the adic completion functor~$\lambda$
is right adjoint to the adic torsion functor~$\gamma$; this adjunction
is called the \emph{MGM duality}.

 Then it follows that the functor~$\gamma$ is the coreflector onto
its essential image, i.~e., onto the full triangulated subcategory
$\gamma(\sD)\subset\sD$.
 In other words, the functor $\gamma\:\sD\rarrow\gamma(\sD)$ is
right adjoint to the inclusion $\gamma(\sD)\rarrow\sD$; this is
a consequence of the idempotent comonad structure on~$\gamma$.
 Similarly, the functor~$\lambda$ is the reflector onto the full
triangulated subcategory $\lambda(\sD)\subset\sD$, i.~e., the functor
$\lambda\:\sD\rarrow\lambda(\sD)$ is left adjoint to the inclusion
functor $\lambda(\sD)\rarrow\sD$.
 The full subcategory $\gamma(\sD)\subset\sD$ is called the full
subcategory of \emph{derived torsion complexes}, while the full
subcategory $\lambda(\sD)\subset\sD$ is the full subcategory of
\emph{derived complete complexes}.
 Finally, the restrictions of the functors $\lambda$ and~$\gamma$ are
mutually inverse equivalences $\lambda\:\gamma(\sD)\simeq
\lambda(\sD):\!\gamma$.
 The resulting triangulated equivalence $\gamma(\sD)\simeq\lambda(\sD)$
is called the \emph{MGM equivalence}.

 From our perspective, the value of having such a triangulated
equivalence increases when one has explicit, independent descriptions
of the triangulated categories $\gamma(\sD)$ and~$\lambda(\sD)$ at
one's disposal.
 One can obtain such descriptions with varying degree of explicitness,
depending on one's assumptions.

 In the context of a commutative ring $R$ with an arbitrary finitely
generated ideal $I\subset R$, it is better to take for $\gamma$
and~$\lambda$ the ``sequential'' derived $I$\+adic torsion and
completion endofunctors on $\sD=\sD(R\Modl)$ (in the terminology
of~\cite[Section~2]{Yek3}, see also~\cite[Section~2]{Pmgm}).
 Then the formalism of the preceding paragraphs in this section is
applicable.
 Furthermore, one can describe $\gamma(\sD)\subset\sD$ as the full
subcategory of all complexes with $I$\+torsion cohomology modules in
$\sD(R\Modl)$, and $\lambda(\sD)\subset\sD$ as the full subcategory of
all complexes with $I$\+contramodule cohomology modules in
$\sD(R\Modl)$.
 Here the \emph{$I$\+torsion $R$\+modules} and
the \emph{$I$\+contramodule $R$\+modules} are two full abelian 
subcategories in $R\Modl$ associated with an ideal $I\subset R$.
 These results, going back to~\cite[Section~6]{DG}, can be found in
the form stated above in~\cite[Theorem~3.4]{Pmgm}.
 Under the \emph{weak proregularity} assumption~\cite{PSY} on
the finitely generated ideal $I\subset R$, one can say more.

\subsection{{}}
 The interpretation of the MGM duality over a commutative ring
as an equivalence between the (bounded or unbounded, conventional
or absolute) derived categories of two abelian categories of
torsion modules and contramodules is apparently due to the present
author~\cite[Corollary~3.5 and Theorem~5.10]{Pmgm} (see also
the more recent~\cite[Corollary~8.7]{Ppt}).
 Actually, even within the paper~\cite{Pmgm}, in the context of
a weakly proregular finitely generated ideal in a commutative ring,
two somewhat different approaches to the MGM equivalence were suggested:
a more elementary one in~\cite[Sections~1\+-3]{Pmgm}, and a more
advanced one in~\cite[Sections~4\+-5]{Pmgm}.
 The approach of~\cite[Section~5]{Pmgm} was designed for being
globalized over formal schemes.
 It is the approach of~\cite[Section~5]{Pmgm} that the present
paper aims to develop.

 Thus, our aim in this paper is to work over a formal scheme that is
not embedded, and may not admit any embedding, into a conventional
scheme.
 In this context, our goal is to formulate the MGM duality as
a \emph{derived}, rather than merely triangulated, equivalence;
i.~e., \emph{an equivalence of two derived categories}, rather than
any full subcategories in a derived category.
 This is achieved in
Main Theorems~\ref{fin-krulldim-mgm-main-theorem-I},
\ref{fin-krulldim-mgm-main-theorem-II}
and~\ref{semi-separated-mgm-main-theorem}, which, taken together,
form the first main result of this paper.

 In our extension of the MGM equivalence from commutative rings
to formal schemes, the abelian category of torsion modules is
replaced by the abelian category of quasi-coherent torsion sheaves.
 The abelian category of contramodules, whose role was emphasized
in~\cite{Pmgm}, is replaced by the exact category of contraherent
(or locally contraherent) cosheaves of contramodules.
 Specifically, for any Noetherian formal scheme $\fX$ that is either
of finite Krull dimension or else semi-separated, we construct
a triangulated equivalence
\begin{equation} \label{introd-mgm-formula}
 \boR\fHom_\fX(\cB_\fX^\bu,{-})\:\sD(\fX\Tors)\,\simeq\,
 \sD(\fX\Ctrh):\!\cB_\fX^\bu\ocn_\fX^\boL{-}
\end{equation}
between the conventional unbounded derived category of the abelian
category $\fX\Tors$ of quasi-coherent torsion sheaves on $\fX$ and
the conventional unbounded derived category of the exact category
$\fX\Ctrh$ of contraherent cosheaves of contramodules on~$\fX$.
 Some bounded and absolute derived versions
of~\eqref{introd-mgm-formula} (for bounded or one-sided bounded
derived categories, as well as for absolute derived categories) are
also obtained, particularly in the case of finite Krull dimension.

 The notation $\cB_\fX^\bu$ in~\eqref{introd-mgm-formula} stands for
the natural \emph{dedualizing complex of quasi-coherent torsion
sheaves}, $\cB_\fX^\bu\in\sD^\bb(\fX\Tors)$, which is constructed
in Section~\ref{dedualizing-complexes-on-formal-schemes-subsecn}.
 The notation $\boR\fHom_\fX$ and $\ocn_\fX^\boL$ means (suitably
constructed) derived functors of contraherent $\fHom$ and
contratensor product, cf.\ Section~\ref{introd-sections-2-3} above.

 We also prove that the triangulated
equivalences~\eqref{introd-mgm-formula} commute with the derived
direct images with respect to flat morphisms of Noetherian formal
schemes.

 There are three significant differences between the affine (commutative
algebra) MGM theory of~\cite[Section~5]{Pmgm} or~\cite[Section~8]{Ppt}
and the nonaffine (formal scheme) MGM theory of
Section~\ref{dedualizing-secn} of the present paper:
\begin{itemize}
\item The category of contramodules over a commutative ring with
an ideal is abelian, but the category of contraherent cosheaves of
contramodules on a formal scheme is only exact.
\item The setting in~\cite[Section~5]{Pmgm} and~\cite{Ppt},
following~\cite{PSY}, was that of a weakly proregular finitely generated
ideal, but in the present paper we work under a more restrictive
assumption of Noetherianity (and sometimes even assume finite Krull
dimension, or otherwise semi-separatedness).
\item The result of~\cite[Theorem~5.10]{Pmgm}
or~\cite[Corollary~8.7]{Ppt} is applicable to all the absolute derived
categories alongside with all the conventional derived categories.
 Our Main Theorems~\ref{fin-krulldim-mgm-main-theorem-I}
and~\ref{fin-krulldim-mgm-main-theorem-II}, which are stated in
the context of a Noetherian formal scheme of finite Krull dimension,
taken together, also accept all the conventional derived category
symbols $\st=\bb$, $+$, $-$, or~$\varnothing$, as well as the absolute
derived category symbols $\st=\abs+$, $\abs-$, or~$\abs$.
 However, our other
Main Theorem~\ref{semi-separated-mgm-main-theorem}, which is for
semi-separated Noetherian formal schemes, only handles
the derived category symbols $\st=+$ or~$\varnothing$.
 So we always have an unbounded derived equivalence, but \emph{not}
always a bounded one, depending on the assumptions on the formal
scheme~$\fX$.
\end{itemize}

\subsection{{}}
 The term \emph{co-contra correspondence} (meaning originally
the \emph{comodule-con\-tra\-mod\-ule correspondence}) was coined by
the present author in the monograph~\cite{Psemi} and the subsequent
memoir~\cite{Pkoszul} (see~\cite{Prev} and~\cite{Pksurv} for
the surveys).
 A discussion of the philosophy of co-contra correspondence can be
found in the introductions to the papers~\cite{Pmgm} and~\cite{Pps}.

 In particular, the MGM equivalence, in its formulations suggested
in~\cite{Pmgm,Ppt}, is interpreted as an instance of ``na\"\i ve
co-contra correspondence'' (together with such cosheaf-theoretic
results as~\cite[Theorems~4.8.1 and~6.8.1]{Pcosh}
or~\cite[Theorems~5.26 and~5.27]{Pdomc}).
 The word ``na\"\i ve'' in this context means that the triangulated
equivalence connects the conventional or absolute derived categories.
 The proper (or ``non-na\"\i ve'') co-contra correspondence is supposed
to be an equivalence between the \emph{coderived} and
the \emph{contraderived} category.

 We construct such a triangulated equivalence, for a Noetherian formal
scheme $\fX$ with a dualizing complex $\cD^\bu$ of quasi-coherent
torsion sheaves on $\fX$, in our
Main Theorem~\ref{nonssep-derived-co-contra-main-thm}, which is
the second main result of this paper.
 Specifically, we have
\begin{equation} \label{introd-co-contra-formula}
 \boR\fHom_\fX(\cD^\bu,{-})\:\sD^\co(\fX\Tors)\,\simeq\,
 \sD^\ctr(\fX\Ctrh):\!\cD^\bu\ocn_\fX^\boL{-}.
\end{equation}
 The superindices $\co$ and $\ctr$
in~\eqref{introd-co-contra-formula} stand for the coderived and
the contraderived category (a historical and philosophical
discussion of those can be found in~\cite[Remark~9.2]{PS4}
and~\cite[Section~7]{Pksurv}).

 We also prove that the triangulated
equivalences~\eqref{introd-co-contra-formula} commute with the derived
direct images with respect to smooth morphisms of Noetherian formal
schemes.

 Preceding formulations of the co-contra correspondence in
the context of torsion modules and contramodules can be found
in~\cite[Theorems~D.1.4, D.5.10, and~E.2.9]{Pcosh}
and~\cite[Corollaries~11.4 and~11.5]{Ppt}.
 For previously known manifestations of the co-contra correspondence
phenomenon in the cosheaf-theoretic context,
see~\cite[Theorem~6.8.2]{Pcosh} and~\cite[Corollaries~12.5
and~12.8]{Pdomc}.

 The term \emph{dedualizing complex} was coined by the present
author in~\cite{Pmgm} (while the classical notion of a dualizing
complex, due to Grothendieck, goes back to~\cite{HarRD}).
 Once again, we refer to the introductions to
the papers~\cite{Pmgm} and~\cite{Pps} for a philosophical discussion
of dedualizing versus dualizing complexes and their relevance to
the na\"\i ve versus non-na\"\i ve co-contra correspondence.

\subsection{{}}
 As pointed out in~\cite[Sections~7.11.3--7.11.4]{BD2} and
elaborated in~\cite[Chapters~2\+-3]{Psemten}, the passage from
schemes to ind-schemes (or to formal schemes) splits the theory
of quasi-coherent sheaves into two branches.
 In the terminology of~\cite{Psemten}, there are quasi-coherent
torsion sheaves, and there are \emph{pro-quasi-coherent pro-sheaves}.
 The \emph{flat} pro-quasi-coherent pro-sheaves are well-behaved, in
that they form an exact category in a natural way.

 The objects dual-analogous to flat (or very flat) pro-quasi-coherent
pro-sheaves, called \emph{locally injective} (\emph{locally})
\emph{ind-contraherent ind-cosheaves}, are introduced in
the final Section~\ref{ind-cosheaves-secn} of this paper.
 We construct the functors of direct and inverse image of locally
injective locally ind-contraherent ind-cosheaves with respect to
suitable classes of morphisms of Noetherian formal schemes.

 So the theory of quasi-coherent sheaves and contraherent cosheaves,
as presented in~\cite{Pcosh}, splits roughly into two halves when
one passes from schemes to formal schemes.
 One half generalizes naturally to quasi-coherent torsion sheaves
and contraherent cosheaves of contramodules, while the other half
finds its place in the realm of flat pro-quasi-coherent pro-sheaves
and locally injective ind-contraherent ind-cosheaves.

 For example, injectivity is a property of quasi-coherent torsion
sheaves, while flatness and very flatness are properties of
pro-quasi-coherent pro-sheaves.
 Local contraadjustedness and local cotorsion are properties of
contraherent cosheaves of contramodules, while local injectivity
is a feature of ind-contraherent ind-cosheaves.

 Section~\ref{reminder-secn}, whose aim is to prepare ground for
Section~\ref{ind-cosheaves-secn}, is a reminder of those aspects of
the theory developed in~\cite[Sections~2\+-4]{Pcosh} which do not seem
to be generalizable to quasi-coherent torsion sheaves and contraherent
cosheaves of contramodules, and accordingly, have no counterparts in
Sections~\ref{contraherent-cosheaves-of-contramods-secn}\+-%
\ref{injective-projective-secn} of the present paper.

 In Section~\ref{ind-cosheaves-secn}, we also discuss contraadjusted
and cotorsion flat pro-quasi-coherent pro-sheaves on a semi-separated
Noetherian formal scheme $\fX$, defined straightforwardly as the flat
pro-quasi-coherent pro-sheaves whose inverse images to all closed
subschemes of definition $Z$ in $\fX$ are contraadjusted/cotorsion
flat quasi-coherent sheaves on~$Z$.
 This definition works very well in this context: the cotorsion flat
pro-quasi-coherent pro-sheaves are the injective objects of the exact
category of flat pro-quasi-coherent pro-sheaves on $\fX$, and there
are enough of them; while the pair of classes (very flat
pro-quasi-coherent pro-sheaves, flat contraadjusted pro-quasi-coherent
pro-sheaves) is a heredirary complete cotorsion pair.
 Furthermore, the class of contraadjusted flat pro-quasi-coherent
pro-sheaves on semi-separated Noetherian formal schemes has the usual
antilocality property in the spirit of the discussion
in~\cite[Section~4]{Pal}.
 This provides a partial Noetherian formal scheme analogue of
the theory developed in~\cite[Section~4.1]{Pcosh}.

 Dual-analogously, we define antilocal locally injective
ind-contraherent ind-cosheaves as the locally injective
ind-contraherent ind-cosheaves whose inverse images to all closed
subschemes of definition $Z$ in $\fX$ are antilocal locally injective
contraherent cosheaves on~$Z$.
 This also works quite well: the antilocal locally injective
ind-contraherent ind-cosheaves are the projective objects of the exact
category of locally injective locally ind-contraherent ind-cosheaves
on $\fX$, and there are enough of them.
 This is a partial Noetherian formal scheme analogue of
the theory from~\cite[Section~4.3]{Pcosh}.  {\hbadness=1425\par}

 Furthermore, we also prove a formal scheme version of the theorem of
Esrada and Sl\'avik~\cite[Corollary~1.1 or~6.2]{ES} about a triangulated
equivalence of the derived categories of the exact categories of
flat and very flat quasi-coherent sheaves on a quasi-compact
semi-separated scheme.
 Our formal scheme analogue of this result claims a triangulated
equivalence of the conventional unbounded derived categories of
the exact categories of flat and very flat pro-quasi-coherent
pro-sheaves on a semi-separated Noetherian formal scheme.

\subsection*{Acknowledgement}
 I~am grateful to Jan \v St\!'ov\'\i\v cek, Michal Hrbek, and
Amnon Yekutieli for helpful discussions.
 The author is supported by the GA\v CR project 26-22734S
and the Institute of Mathematics, Czech Academy of Sciences
(research plan RVO:~67985840).

\Section{Preliminaries in Adic Commutative Algebra}

\subsection{Adic topological rings and their maps}
\label{adic-rings-subsecn}
 Let $R$ be a commutative ring.
 By an \emph{adic topology} on $R$ we mean a ring topology for
which there exists a finitely generated ideal $I\subset R$ such that
the ideal $I^n\subset R$ is open for every $n\ge1$ and the descending
sequence $I^n$, \,$n\ge1$, of powers of the ideal $I$ is a base of
neighborhoods of zero in~$R$.
 If this is the case, then the topology on $R$ is said to be
\emph{$I$\+adic}, the ring $R$ is said to be an \emph{adic
topological ring}, and the ideal $I$ is said to be the \emph{ideal of
definition} of the adic topological ring~$R$.

 Let $R$ be an adic topological ring and $I\subset R$ be an ideal of
definition of~$R$.
 Then, for every $n\ge1$, the ideal $I^n\subset R$ is an ideal of
definition, too.
 Another given finitely generated ideal $J\subset R$ is an ideal of
definition if and only if there exist integers $n$, $m\ge1$ such that
$I^n\subset J$ and $J^m\subset I$.
 For any two ideals of definition $I$, $J\subset R$, the ideals
$I+J$ and $IJ\subset R$ are also ideals of definition in~$R$.

\begin{exs}
 The discrete topology on any commutative ring $R$ is adic.
 The zero ideal $(0)\subset R$ is an ideal of definition of $R$
in the discrete topology.
 The indiscrete topology on $R$ is adic as well.
 The unit ideal $(1)\subset R$ is an ideal of definition of $R$
in the indiscrete topology.
\end{exs}

\begin{lem} \label{adically-Noetherian-well-defined-lemma}
 Let $I$ and $J$ be two ideals of definition of an adic topological
ring~$R$.
 Then the quotient ring $R/I$ is Noetherian if and only if
the quotient ring $R/J$ is Noetherian.
\end{lem}

\begin{proof}
 The assertion holds due to the assumption that the ideals $I$
and $J$ are finitely generated.
 If so, then the ring $R/I$ is Noetherian if and only if the ring
$R/I^n$ is (for all, or equivalently, any given $n\ge1$).
\end{proof}

 We will say that an adic topological ring $R$ is \emph{adically
Noetherian} (cf.~\cite[Section~9]{Ppt}) if it satisfies the equivalent
conditions of Lemma~\ref{adically-Noetherian-well-defined-lemma}.

 If the adic topological ring $R$ is adically Noetherian, then
the nilradical of $R/I$ is a nilpotent ideal.
 In other words, the ideal $I$ contains a suitable power of its
radical: there exists $n\ge1$ such that $(\sqrt I)^n\subset I$.
 Consequently, there exists a unique \emph{maximal ideal of
definition} $I_\mx\subset R$, characterized by the properties that
$I_\mx$ is an ideal of definition of $R$ and the quotient ring
$R/I_\mx$ is reduced (i.~e., contains no nonzero nilpotent elements).
 The maximal ideal of definition $I_\mx$ can be constructed as
the radical of any other ideal of definition,
$I_\mx=\sqrt I\subset R$.
 
\begin{lem} \label{continuous-ring-map-lemma}
 Let $R$ and $S$ be two commutative rings with adic topologies,
and let $f\:R\rarrow S$ be a ring homomorphism.
 Then the following four conditions are equivalent:
\begin{enumerate}
\item the ring homomorphism~$f$ is continuous;
\item there exist an ideal of definition $I\subset R$ and
an ideal of definition $J\subset S$ such that $f(I)\subset J$;
\item for every ideal of definition $J\subset S$ there exists
an ideal of definition $I\subset R$ such that $f(I)\subset J$;
\item for every ideal of definition $I\subset R$ there exists
an ideal of definition $J\subset S$ such that $f(I)\subset J$.
\end{enumerate}
 If the topological rings $R$ and $S$ are adically Noetherian, then
conditions~\textup{(1\+-4)} are also equivalent to the following
condition:
\begin{enumerate}
\setcounter{enumi}{4}
\item denoting by $I_\mx\subset R$ and $J_\mx\subset S$ the maximal
ideals of definition of $R$ and $S$, one has $f(I_\mx)\subset J_\mx$.
\end{enumerate}
\end{lem}

\begin{proof}
 (1)~$\Longrightarrow$~(3) The ideal of definition $J\subset S$ is
open, so there exists an open ideal $I'\subset R$ such that
$f(I')\subset J$.
 Now any open ideal $I'\subset R$ contains some ideal of definition,
i.~e., there exists an ideal of definition $I$ such that $I\subset I'$.

 (3)~$\Longrightarrow$~(2) Obvious.
 
 (2)~$\Longrightarrow$~(1) Let $I'\subset R$ and $J'\subset S$ be
ideals of definition for which $f(I')\subset J'$, and let $J\subset S$
be an open ideal.
 Then there exists $n\ge1$ such that $(J')^n\subset J$.
 Hence $f((I')^n)\subset (J')^n\subset J$, so
$(I')^n\subset f^{-1}(J)$.
 Thus $f^{-1}(J)$ is an open ideal in~$R$.
 
 (2)~$\Longrightarrow$~(4) Let $I'\subset R$ and $J'\subset S$ be
ideals of definition for which $f(I')\subset J'$, and let $I\subset R$
be another given ideal of definition.
 Then there exists $n\ge1$ such that $I^n\subset I'$, hence
$f(I^n)\subset J'$.
 Put $J=Sf(I)+J'\subset S$, where $Sf(I)$ is the extension of
the ideal $I\subset R$ in the ring~$S$.
 Then we have $J'\subset J$ and $J^n\subset J'$, hence $J$ is
an ideal of definition of~$S$.
 By construction, $f(I)\subset J$.

 (4)~$\Longrightarrow$~(2) Obvious.
 
 (2)~$\Longrightarrow$~(5) We have $I_\mx=\sqrt I$ and $J_\mx=
\sqrt J$.
 Thus $f(I)\subset J$ implies $f(I_\mx)\subset J_\mx$.

 (5)~$\Longrightarrow$~(2) Take $I=I_\mx$ and $J=J_\mx$.
\end{proof}

\begin{lem} \label{tight-ring-map-lemma}
 Let $R$ and $S$ be two commutative rings with adic topologies,
and let $f\:R\rarrow S$ be a ring homomorphism.
 Then the following seven conditions are equivalent:
\begin{enumerate}
\item there exist an ideal of definition $I\subset R$ and
an ideal of definition $J\subset S$ such that $J\subset S f(I)$;
\item for every ideal of definition $J\subset S$ there exist
an ideal of definition $I\subset R$ and an integer $n\ge1$
such that $J^n\subset S f(I)$;
\item for any two ideals of definition $I\subset R$ and $J\subset S$
there exists an integer $n\ge1$ such that $J^n\subset S f(I)$;
\item for every ideal of definition $I\subset R$ there exists
an ideal of definition $J\subset S$ such that $J\subset S f(I)$;
\item for every open ideal $I\subset R$, the ideal $S f(I)$ is open
in~$S$;
\item the open ideals $I\subset R$ for which the ideal $S f(I)$ is
open in $S$ form a base of neighborhoods of zero in~$R$;
\item there exists an ideal of definition $I\subset R$ such that
$S f(I)$ is an open ideal in~$S$.
\end{enumerate}
\end{lem}

\begin{proof}
 (1)~$\Longrightarrow$~(2) Let $I'\subset R$ and $J'\subset S$ be
ideals of definition for which $J'\subset Sf(I')$, and let $J\subset S$
be another given ideal of definition.
 Then there exists $n\ge1$ such that $J^n\subset J'$, hence
$J^n\subset Sf(I')$.
 It remains to put $I=I'$.

 (2)~$\Longrightarrow$~(3) Let $I'\subset R$ be an ideal of definition
and $n'\ge1$ be an integer such that $J^{n'}\subset Sf(I')$.
 Pick an integer $m\ge1$ such that $(I')^m\subset I$.
 Then $J^{n'm}\subset Sf((I')^m)\subset Sf(I)$.
 It remains to put $n=n'm$.

 (3)~$\Longrightarrow$~(4)~$\Longrightarrow$~(1)
and (3)~$\Longrightarrow$~(2)~$\Longrightarrow$~(1) Obvious.

 (4)~$\Longrightarrow$~(5) Given an open ideal $I\subset R$,
there exists an ideal of definition $I'\subset R$ such that
$I'\subset I$.
 By~(4), there exists an ideal of definition $J'\subset S$ such that
$J'\subset Sf(I')$.
 Then $J'\subset Sf(I)$, and it follows that $Sf(I)$ is an open ideal
in~$S$.

 (5)~$\Longrightarrow$~(6) $\Longrightarrow$~(7)
$\Longleftrightarrow$~(1) Obvious.
 
 (6)~$\Longrightarrow$~(4) Let $I\subset R$ be an ideal of definition.
 By~(6), there exists an open ideal $I'\subset R$ such that
$I'\subset I$ and the ideal $Sf(I')$ is open in~$S$.
 Then there exists an ideal of definition $J\subset S$ such that
$J\subset Sf(I')$.
 It follows that $J\subset Sf(I)$.
\end{proof}

\begin{rem} \label{remark-to-tight-lemma}
 Notice that the equivalent conditions of
Lemma~\ref{tight-ring-map-lemma} do \emph{not} imply that for
every ideal of definition $J\subset S$ there exists an ideal of
definition $I\subset R$ such that $J\subset Sf(I)$.
 For a counterexample, it suffices to take $R=k$ to be a field and
$S=k[\epsilon]/(\epsilon^2)$ to be the $k$\+algebra of dual numbers.
 Endow both $R$ and $S$ with the discrete topologies, and let
$f\:R\rarrow S$ be the identity inclusion.
 Then $I=0$ is the only ideal of definition in $R$, and
$J=(\epsilon)$ is an ideal of definition in $S$, but
$J\not\subset Sf(I)$.
 On the other hand, $J'=(0)$ is also an ideal of definition in $S$,
and $J'=Sf(I)$.
\end{rem}

 We will say that a ring homomorphism $f\:R\rarrow S$ acting between
two rings with adic topologies is a \emph{tight} ring map if
the equivalent conditions of Lemma~\ref{tight-ring-map-lemma} hold.
 Notice that a tight ring map \emph{need not} be continuous in our
definition.
 Using the criterion of Lemma~\ref{tight-ring-map-lemma}(4), one
can easily check that the composition of any two tight ring maps of
adic topological rings is tight.

 Given an $I$\+adic topological ring $R$ and an $R$\+module $M$,
the \emph{$I$\+adic topology on $M$} is defined by the property
that the submodules $I^nM\subset M$, \,$n\ge1$, form a base of
open neighborhoods of zero.
 Clearly, if $I'$ and $I''$ are two ideals of definition of $R$,
then the $I'$\+adic and the $I''$\+adic topologies on $M$ coincide.

\begin{lem} \label{continuous-and-tight-ring-map-lemma}
 Let $R$ and $S$ be two commutative rings with adic topologies,
and let $f\:R\rarrow S$ be a ring homomorphism.
 Then the following five conditions are equivalent:
\begin{enumerate}
\item the ring homomorphism~$f$ is continuous and tight;
\item there exists an ideal of definition $I\subset R$ such that
the topology on $S$ is the $I$\+adic topology;
\item for every ideal of definition $I\subset R$, the topology on $S$
is the $I$\+adic topology;
\item there exists an ideal of definition $I\subset R$ such that
$J=Sf(I)$ is an ideal of definition in~$S$;
\item for every ideal of definition $I\subset R$, the ideal $J=Sf(I)
\subset S$ is an ideal of definition in~$S$.
\end{enumerate}
\end{lem}

\begin{proof}
 (2)~$\Longleftrightarrow$~(3) Follows from the discussion in
the paragraph preceding the lemma.

 (1)~$\Longrightarrow$~(3) Since $f$~is continuous,
Lemma~\ref{continuous-ring-map-lemma}(3) says that for every
ideal of definition $J\subset S$ there exists an ideal of
definition $I\subset R$ such that $Sf(I)\subset J$.
 Since $f$~is tight, Lemma~\ref{tight-ring-map-lemma}(4) says that
for every ideal of definition $I\subset R$ there exists an ideal
of definition $J\subset S$ such that $J\subset Sf(I)$.
 Thus the $J$\+adic and the $I$\+adic topologies on $S$ coincide.

 (3)~$\Longrightarrow$~(1) Under~(3), both the topologies on $R$
and $S$ are the $I$\+adic topologies.
 As every $R$\+module map $M\rarrow N$ is continuous with respect
to the $I$\+adic topologies on $M$ and $N$, it follows that $f$~is
continuous.

 To check that $f$~is tight, let $I\subset R$ and $J\subset S$ be
any two ideals of definition.
 In the $I$\+adic topology, $Sf(I)$ is an open subgroup in~$S$.
 Thus (2)~implies that there exists $n\ge1$ for which
$J^n\subset Sf(I)$.
 So the condition of Lemma~\ref{tight-ring-map-lemma}(3) is satisfied.

 (1)~$\Longrightarrow$~(5) By Lemma~\ref{tight-ring-map-lemma}(4),
there exists an ideal of definition $J'\subset S$ such that
$J'\subset Sf(I)$.
 By Lemma~\ref{continuous-ring-map-lemma}(3), there exists
an ideal of definition $I'\subset R$ such that $Sf(I')\subset J'$.
 Pick $n\ge1$ such that $I^n\subset I'$; then we have
$(Sf(I))^n=Sf(I^n)\subset J'\subset Sf(I)$.
 Since $J'$ is an ideal of defintion in $S$, it follows that
$Sf(I)$ is an ideal of definition in $S$, too.

 (5)~$\Longrightarrow$~(4)~$\Longrightarrow$~(2) Obvious.
\end{proof}

 It is clear from Lemma~\ref{continuous-and-tight-ring-map-lemma} that,
for any homomorphism of commutative rings $f\:R\rarrow S$ and any
adic topology on the ring $R$, there exists a unique adic topology on
the ring $S$ making the ring map~$f$ continuous and tight.
 It also follows that any tight continuous ring map that is
an isomorphism of abstract commutative rings is an isomorphism
of topological rings.

\begin{exs} \label{non-continuous-and-or-non-tight-examples}
 Let $k$~be a field and $R=S=k[x]$ be the ring of polynomials in
one variable~$x$ over~$k$.
 Let $f\:R\rarrow S$ be the identity map.
 Consider the ideals $I=(x)\subset R$, \ $J=(x)\subset S$,
$I_0=(0)\subset R$, and $J_0=(0)\subset S$.
 Then the map~$f$ is continuous and tight with respect to
the $I$\+adic topology on $R$ and the $J$\+adic topology on~$S$.
 The map~$f$ is also continuous and tight with respect to
the $I_0$\+adic topology on $R$ and the $J_0$\+adic topology on~$S$.
 With respect to the $I_0$\+adic topology on $R$ and the $J$\+adic
topology on $S$, the map~$f$ is continuous but not tight.
 With respect to the $I$\+adic topology on $R$ and the $J_0$\+adic
topology on $S$, the map~$f$ is tight but not continuous (in our
terminology).

 To give another example, let $R=k[x]$ be the ring of polynomials
in one variable and $S=k[x,y]$ be the ring of polynomials in
two variables.
 Consider the ideals $I=(x)\subset R$ and $J=(x,y)\subset S$.
 Let $g\:R\rarrow S$ be the identity inclusion map and $h\:S\rarrow R$
be the $k$\+algebra homomorphism taking $x$ to~$x$ and $y$ to~$0$.
 Then the map~$g$ is continuous but not tight, while the map~$h$ is
continuous and tight.
 For a perhaps more complicated example of a tight continuous map,
see Remark~\ref{remark-to-tight-lemma} (in fact, any morphism of
discrete commutative rings is continuous and tight).
\end{exs}

\begin{lem} \label{flat-homomorphism-of-adic-rings}
 Let $R$ and $S$ be two commutative rings with adic topologies, and
let $f\:R\rarrow S$ be a tight continuous ring homomorphism.
 Then the following two conditions are equivalent:
\begin{enumerate}
\item there exists an ideal of definition $I\subset R$
such that the $R/I^n$\+module $S/Sf(I^n)$ is flat for every $n\ge1$;
\item for every open ideal $I\subset R$, the $R/I$\+module $S/Sf(I)$
is flat.
\end{enumerate}
\end{lem}

\begin{proof}
 The point is that if $I\subset I'\subset R$ are two ideals in $R$
such that the $R/I$\+module $S/Sf(I)$ is flat, then the $R/I'$\+module
$S/Sf(I')=R/I'\ot_{R/I}S/Sf(I)$ is flat, too.
\end{proof}

 We will say that a tight continuous homomorphism of adic topological
rings $f\:R\rarrow S$ is \emph{flat} if it satisfies the equivalent
conditions of Lemma~\ref{flat-homomorphism-of-adic-rings}.
 Clearly, if a tight continuous map of adic topological rings~$f$
is flat as a map of abstract rings (i.~e., the $R$\+module $S$ is
flat), then $f$~is also flat as as a tight continuous map of adic
topological rings.
 But the converse is not true in general.

 Using the criterion of
Lemma~\ref{flat-homomorphism-of-adic-rings}(2), one can easily check
that the composition of any two flat tight continuous ring maps of
adic topological rings is flat.

\subsection{Torsion modules} \label{prelim-torsion-modules-subsecn}
 Let $R$ be a commutative ring and $s\in R$ be an element.
 Following the notation in~\cite[Section~1.1]{Pcosh},
\cite[Section~2.1]{Pcta}, etc., we denote by $R[s^{-1}]$ the ring $R$
with the element~$s$ formally inverted.
 In other words, $R[s^{-1}]=S^{-1}R$ is the localization of $R$
at the multiplicative subset $S=\{1,s,s^2,s^3,\dotsc\}\subset R$.

 Let $M$ be an $R$\+module.
 An element $x\in M$ is said to be \emph{$s$\+torsion} if
there exists $n\ge1$ such that $s^nx=0$ in~$M$.
 The $R$\+module $M$ is said to be \emph{$s$\+torsion} if all
the elements of $M$ are $s$\+torsion.
 Equivalently, an $R$\+module $M$ is $s$\+torsion if and only if
$R[s^{-1}]\ot_RM=0$.

 Let $I\subset R$ be an ideal.
 We say that an $R$\+module $\cM$ is \emph{$I$\+torsion} if $\cM$ is
$s$\+torsion for all $s\in R$.
 One can easily see that it suffices to check this condition for
the element~$s$ ranging over any chosen set of generators of
the ideal~$I$ \,\cite[Lemma~1.2]{Ptd}.

 Other definitions of the term ``$I$\+torsion $R$\+module'' exist
in the literature (see, e.~g., \cite[Section~3]{PSY}), but
the difference between them does \emph{not} manifest itself for
\emph{finitely generated} ideals~$I$, which we are really
interested in.
 So, for a finitely generated ideal $I\subset R$, an $R$\+module $\cM$
is $I$\+torsion if and only if, for every element $x\in \cM$,
there exists an integer $n\ge1$ such that $I^nx=0$ in~$\cM$
\,\cite[Lemma~1.3]{Ptd}.

 For any $R$\+module $M$ and any ideal $J\subset R$, we denote by
${}_JM\subset M$ the submodule of all elements annihilated by $J$
in~$M$.
 So ${}_JM\simeq\Hom_R(R/J,M)$.
 According to the previous paragraph, an $R$\+module $\cM$ is
$I$\+torsion if and only if $\cM=\bigcup_{n\ge1}{}_{I^n}\cM$.

\begin{lem} \label{tensor-product-is-torsion}
 Let $R$ be a commutative ring, $I\subset R$ be an ideal, $\cM$ be
an $I$\+torsion $R$\+module, and $N$ be an $R$\+module.
 Then the $R$\+module\/ $\cM\ot_RN$ is $I$\+torsion.
\end{lem}

\begin{proof}
 The assertion is obvious; see~\cite[Lemma~6.1(a)]{Pcta} for
some details.
\end{proof}

 Now let $R$ be a commutative ring endowed with an \emph{adic
topology} in the sense of Section~\ref{adic-rings-subsecn}.
 Let $I$, $J\subset R$ be two ideals of definition of the ring~$R$.
 Then an $R$\+module $\cM$ is $I$\+torsion if and only if it is
$J$\+torsion.
 If this is the case, we say that $\cM$ is a \emph{torsion $R$\+module}
(with respect to the given adic topology of~$R$).

 An $R$\+module $\cM$ is torsion if and only if, for every element
$x\in\cM$, the annihilator of~$x$ in $R$ is an open ideal.
 In other words, an $R$\+module $\cM$ is torsion if and only if it is
\emph{discrete} (as a module over the topological ring~$R$) in
the sense of~\cite[Section~VI.4]{Sten}, \cite[Section~1.4]{Beil},
\cite[Section~D.5.1]{Psemi}, \cite[Section~E.2]{Pcosh},
\cite[Sections~1.4 and~2.1]{Prev}, \cite[Section~2.3]{Pcoun},
\cite[Section~2.4]{Pproperf}, \cite[Section~7.2]{PS1},
\cite[Section~1]{PS3}, etc.

 For any ring $R$, we denote by $R\Modl$ the abelian category of (left)
$R$\+modules.
 Given an ideal $I$ in a commutative ring $R$, the notation
$R\Modl_{I\dtors}\subset R\Modl$ stands for the full subcategory of
$I$\+torsion $R$\+modules.
 For an adic topological ring $R$, we will write simply $R\Modl_\tors$
instead of $R\Modl_{I\dtors}$, presuming any ideal of definition $I$ of
the topological ring~$R$.
 We will also write $R\Tors=R\Modl_\tors$.

 For any ideal $I$ in a commutative ring $R$, the full subcategory
$R\Modl_{I\dtors}$ is closed under submodules, quotients, extensions,
and infinite direct sums in $R\Modl$.
 In other words, $R\Modl_{I\dtors}$ is a \emph{localizing subcategory},
or in a different terminology, a \emph{hereditary torsion
class}~\cite[Sections~VI.2\+-3]{Sten} in $R\Modl$.

 It follows that the category $R\Modl_{I\dtors}$ is a Grothendieck
abelian category.
 The identity inclusion functor $R\Modl_{I\dtors}\rarrow R\Modl$ is
exact and preserves infinite direct sums (hence also all colimits),
but \emph{not} infinite products.
 The infinite direct product functors in $R\Modl_{I\dtors}$ are
usually \emph{not} exact.

\subsection{Adic completions} \label{prelim-adic-completions-subsecn}
 Let $R$ be a commutative ring and $I\subset R$ be a finitely
generated ideal.
 Then the construction of the \emph{$I$\+adic completion}
$$
 \Lambda_I(M)=\varprojlim\nolimits_{n\ge1} M/I^nM
$$
of an arbitrary $R$\+module $M$ is reasonably well-behaved~\cite{Yek}.
 The $I$\+adic completion $\Lambda_I(M)$ is a naturally a module
over the $I$\+adic completion $\Lambda_I(R)=\varprojlim_{n\ge1}R/I^n$
of the ring~$R$.

 There is a natural \emph{completion map} $\lambda_{I,M}\:
M\rarrow\Lambda_I(M)$.
 We say that an $R$\+module $M$ is \emph{$I$\+adically separated}
if the map~$\lambda_{I,M}$ is injective, and that $M$ is
\emph{$I$\+adically complete} if $\lambda_{I,M}$~is surjective.
 For any $R$\+module $M$, the $R$\+module $\Lambda_I(M)$ is
$I$\+adically separated and complete~\cite[Corollary~3.6]{Yek}.
 The latter assertion depends on the assumption that the ideal
$I\subset R$ is finitely generated~\cite[Example~1.8]{Yek}.
 The functor $\Lambda_I$ is left adjoint to the inclusion of
the full subcategory of $I$\+adically separated and complete
$R$\+modules into $R\Modl$, and the map~$\lambda_{I,M}$ is
the adjunction unit~\cite[Theorem~5.8]{Pcta}.

 Put $\fR=\Lambda_I(R)=\varprojlim_{n\ge1}R/I^n$.
 The extension $\fR I$ of the ideal $I\subset R$ to the ring $\fR$
coincides with the $I$\+adic completion $\fI=\varprojlim_{n\ge2}
I/I^n\subset\fR$ of the ideal~$I$.
 The commutative ring $\fR$ is naturally endowed with the projective
limit topology (of discrete rings $R/I^n$), which coincides with
the $I$\+adic topology and with the $\fI$\+adic topology on~$\fR$
\,\cite[proof of Theorem~5.8]{Pcta}.
 So $\fR$ is an adic topological ring, too; and $\fI$ is an ideal
of definition of~$\fR$.
 The topological ring $\fR$ is determined by the adic topological
ring $R$ and does \emph{not} depend on the choice of a particular ideal 
of definition $I\subset R$ (while the ideal $\fI\subset\fR$ depends
on the ideal $I\subset R$, of course).

 Given an adic topological ring $R$, the completion functor
$\Lambda_I\:R\Modl\rarrow\fR\Modl$ and the natural
map~$\lambda_{I,M}$ depend on the adic topology on the ring
$R$ only, and \emph{not} on the choice of a particular ideal of
definition $I\subset R$.
 Consequently, the properties of an $R$\+module to be $I$\+adically
separated or $I$\+adically complete are also determined by
the $I$\+adic topology on $R$ and do not depend on the choice
of an ideal of definition $I$ for the given adic topology.
 So we will write $\Lambda_R=\Lambda=\Lambda_I$ and
$\lambda_M=\lambda_{I,M}$, presuming any ideal of definition $I$
of the topological ring~$R$.

 Following the terminology above, an adic topological ring $R$ is
said to be \emph{separated} if the map $\lambda_R\:R\rarrow\fR$
is injective and \emph{complete} if the map~$\lambda_R$ is
surjective.
 For any adic topological ring $R$, the adic topological ring
$\fR=\Lambda(R)$ is separated and complete, and the completion
homomorphism of topological rings $\lambda_R\:R\rarrow\fR$ is
continuous, tight, and flat (in the sense of
Section~\ref{adic-rings-subsecn}).

 Let $R$ be an adic topological ring and $\fR=\Lambda(R)$ be its
completion.
 Then the rules $\fI=\fR\lambda_R(I)$ and $I=\lambda_R^{-1}(\fI)$
establish a natural bijection between the ideals of definition
$I\subset R$ and the ideals of definition $\fI\subset\fR$.
 More generally, the same rules provide a bijection between
open ideals $I\subset R$ and open ideals $\fI\subset\fR$.
 The map~$\lambda_R$ induces ring isomorphisms $R/I\simeq\fR/\fI$
for open ideals $I\subset R$ and $\fI\subset\fR$ corresponding
to each other under this bijection.

 Any continuous homomorphism of adic topological rings $f\:R\rarrow S$
induces a continuous homomorphism of the completions
$\Lambda(f)\:\Lambda(R)\rarrow\Lambda(S)$.

\begin{lem} \label{induced-map-of-complete-rings-lemma}
 Let $f\:R\rarrow S$ be a tight continuous homomorphism of adic
topological rings.
 Then \par
\textup{(a)} the induced continuous ring map of the completions
$\Lambda(f)\:\Lambda(R)\rarrow\Lambda(S)$ is also tight; \par
\textup{(b)} the tight continuous ring map~$f$ is flat if and only
if the tight continuous ring map $\Lambda(f)$ is flat.
\end{lem}

\begin{proof}
 Put $\fR=\Lambda(R)$, \,$\fS=\Lambda(S)$, and $\ff=\Lambda(f)$.

 Part~(a): let $I\subset R$ be an ideal of definition in~$R$.
 Then, by Lemma~\ref{continuous-and-tight-ring-map-lemma}(5),
the ideal $J=Sf(I)\subset S$ is an ideal of definition in~$S$.
 Following the discussion above, the ideal $\fI=\fR\lambda_R(I)
\subset\fR$ is an ideal of definition in $\fR$ and the ideal
$\fJ=\fS\lambda_S(J)$ is an ideal of definition in~$\fS$.
 Clearly, we have $\fJ=\fS\ff(\fI)$.
 By Lemma~\ref{continuous-and-tight-ring-map-lemma}(4), it follows
that the ring map~$\ff$ is tight.

 Part~(b): in the notation of the previous paragraph, we have
$R/I\simeq\fR/\fI$ and $S/J\simeq\fS/\fJ$.
 Hence $S/J$ is a flat module over $R/I$ if and only if $\fS/\fJ$ is
a flat module over $\fR/\fI$.
\end{proof}

\begin{rem} \label{tightness-completion-remark}
 The converse assertion to
Lemma~\ref{induced-map-of-complete-rings-lemma}(a) is \emph{not} true:
 There exist continuous maps of adic topological rings $f\:R\rarrow S$
such that the continuous ring map $\Lambda(f)$ is tight but the map~$f$
is not.
 For a counterexample, it suffices to take $R=S=k[x]$ and $f=\id_{k[x]}$
as in Examples~\ref{non-continuous-and-or-non-tight-examples}, with
the discrete ($(0)$\+adic) topology on $R$ and the indiscrete
($(1)$\+adic) topology on~$S$.
 Then one has $\Lambda(R)=R$ and $\Lambda(S)=0$; so the ring map
$\ff=\Lambda(f)$ is tight; but the ring map~$f$ is not.
 We will continue this discussion in
Section~\ref{prelim-tight-taut-subsecn}; see
Corollary~\ref{completion-preserves-reflects-tautness-cor} and
Proposition~\ref{complete-tightness=tautness-prop}.
\end{rem}

 Let $R$ be an adic topological ring and $\fR=\Lambda(R)$ be its
completion.
 The adic completion functor $\Lambda\:R\Modl\rarrow
R\Modl$ or $\Lambda\:R\Modl\rarrow\fR\Modl$ is neither left,
nor right exact~\cite[Example~3.20]{Yek}.
 This fact is closely related to the fact that the category of
separated and complete $R$\+modules is \emph{not}
abelian~\cite[Example~2.7(1)]{Pcta}.
 To overcome these difficulties, one considers contramodule
$R$\+modules instead of the separated and complete $R$\+modules;
see the next section.

\subsection{Contramodules} \label{prelim-contramodules-subsecn}
 Let $R$ be a commutative ring and $s\in R$ be an element.
 An $R$\+module $P$ is said to be an \emph{$s$\+contramodule} if
$\Hom_R(R[s^{-1}],P)=0=\Ext_R^1(R[s^{-1}],P)$.
 The ring $R[s^{-1}]$ is considered as an $R$\+module here.
 The definition above does not mention any $\Ext^i$ with $i\ge2$
because the projective dimension of the $R$\+module $R[s^{-1}]$
never exceeds~$1$ \,\cite[Section~1.1]{Pcosh}, \cite[proof of
Lemma~B.7.1]{Pweak}, \cite[proof of Lemma~2.1]{Pcta}.
 A detailed discussion of this definition, which goes back
to~\cite[Remark~A.1.1]{Psemi}, \cite[Theorem~B.1.1]{Pweak},
\cite[Section~D.1]{Pcosh}, and~\cite[Section~2]{Pmgm},
can be found in~\cite[Sections~1\+-3]{Pcta}.

 Let $I\subset R$ be an ideal.
 An $R$\+module $\fP$ is said to be an \emph{$I$\+contramodule} (or
an \emph{$I$\+contramodule $R$\+module}) if $\fP$ is
an $s$\+contramodule for every $s\in I$.
 It suffices to check this condition for the element~$s$ ranging over
any chosen set of generators of the ideal~$I$
\,\cite[Theorem~B.1.1]{Pweak}, \cite[Section~2]{Pmgm}.
 A detailed discussion with two proofs can be found
in~\cite[Theorem~5.1]{Pcta}.

\begin{lem} \label{Hom-is-a-contramodule}
 Let $R$ be a commutative ring and $I\subset R$ be an ideal.
 Then \par
\textup{(a)} for any $I$\+torsion $R$\+module\/ $\cM$ and any
$R$\+module $P$, the $R$\+module\/ $\Hom_R(\cM,P)$ is
an $I$\+contramodule; {\hfuzz=3pt\par}
\textup{(b)} for any $R$\+module $M$ and any $I$\+contramodule\/
$R$\+module\/ $\fP$, the $R$\+module\/ $\Hom_R(M,\fP)$ is
an $I$\+contramodule.
\end{lem}

\begin{proof}
 This is~\cite[Lemma~6.1(b)]{Pcta}.
\end{proof}

 Assume that the ideal $I\subset R$ is finitely generated.
 Then any $I$\+contramodule $R$\+module is $I$\+adically
complete~\cite[Theorem~5.6]{Pcta}, but it \emph{need not} be
$I$\+adically separated~\cite[Example~2.5]{Sim},
\cite[Example~3.20]{Yek}, \cite[Example~4.33]{PSY},
\cite[Section~A.1.1]{Psemi}, \cite[Section~1.5]{Prev},
\cite[Example~2.7(1)]{Pcta}.
 Any $I$\+adically separated and complete $R$\+module is
an $I$\+contramodule~\cite[Lemma~5.7]{Pcta}.

 Without reproducing the definition of a \emph{weakly proregular}
finitely generated ideal in a commutative ring (see, e.~g.,
\cite[Sections~4 and~6]{PSY}, \cite[Sections~1 and~2]{Pmgm},
\cite[Section~5]{Yek2}, \cite[Section~3]{Yek3}, \cite[Section~3]{Pdc},
or~\cite[Section~1]{Ppt}, as well as
Section~\ref{dedualizing-complexes-over-rings-subsecn} below),
we recall that \emph{all ideals in Noetherian commutative rings are
weakly proregular}~\cite[Theorem~4.34]{PSY}, \cite[Section~1]{Pmgm},
\cite[Theorem~3.3]{Yek3}.
 This reminder may help the reader to put into proper context
some assumptions and assertions in the references below.

 An $I$\+contramodule $R$\+module is said to be
\emph{quotseparated}~\cite[Section~5.5]{PSl1}, \cite[Section~1]{Pdc}
if it is a quotient $R$\+module of an $I$\+adically separated and
complete $R$\+module.
 The references in the previous paragraph provide examples of
quotseparated $I$\+contramodule $R$\+modules that are not
$I$\+adically separated.
 For an example of an $I$\+contramodule $R$\+module that is not
quotseparated, see~\cite[Example~2.6]{Pmgm},
\cite[Examples~1.8]{Pdc}.
 Over a Noetherian commutative ring $R$, all $I$\+contramodule
$R$\+modules are quotseparated~\cite[Theorem~B.1.1]{Pweak},
\cite[Corollary~3.7 and Remark~3.8]{Pdc}.
 A perhaps more detailed discussion can be found
in~\cite[Section~1]{Ppt}.

 Now let $R$ be a commutative ring endowed with an \emph{adic
topology} in the sense of Section~\ref{adic-rings-subsecn}.
 Let $I$, $J\subset R$ be two ideals of definition of the ring~$R$.
 Then an $R$\+module $\fP$ is an $I$\+contramodule if and only if it is
a $J$\+contramodule~\cite[Remark~5.5]{Pcta}.
 If this is the case, we will say that $\fP$ is a \emph{contramodule
$R$\+module} (with respect to the given adic topology of~$R$).
 Moreover, a contramodule $R$\+module is quotseparated as
an $I$\+contramodule $R$\+module if and only if it is quotseparated
as a $J$\+contramodule $R$\+module; so one can speak of
\emph{quotseparated contramodule $R$\+modules}.

 Given an ideal $I$ in a commutative ring $R$, we denote by
$R\Modl_{I\dctra}\subset R\Modl$ the full subcategory of
$I$\+contramodule $R$\+modules.
 When the ideal $I\subset R$ is finitely generated, the notation
$R\Modl_{I\dctra}^\qs\subset R\Modl_{I\dctra}$ stands for the full
subcategory of quotseparated $I$\+contramodule $R$\+modules.
 For an adic topological ring $R$, we will write simply
$R\Modl_\ctra^\qs\subset R\Modl_\ctra$ instead of
$R\Modl_{I\dctra}^\qs\subset R\Modl_{I\dctra}$, presuming any
ideal of definition $I$ of the topological ring~$R$.

 For any ideal $I$ in a commutative ring $R$, the full subcategory
$R\Modl_{I\dctra}$ is closed under kernels, cokernels, extensions,
and infinite direct products in $R\Modl$ \,\cite[Proposition~1.1]{GL},
\cite[Section~2]{Pmgm}, \cite[Theorem~1.2(a)]{Pcta}.
 It follows that $R\Modl_{I\dctra}$ is an abelian category.
 In fact, $R\Modl_{I\dctra}$ is a locally $\aleph_1$\+presentable
abelian category with enough projective
objects~\cite[Example~4.1(3)]{PR}
(see also~\cite[Example~1.3(4)]{Pper}).
 The identity inclusion functor $R\Modl_{I\dctra}\rarrow R\Modl$
is exact and preserves infinite products (hence also all limits) and
$\aleph_1$\+directed colimits, but \emph{not} infinite direct sums.
 The infinite direct sum functors in $R\Modl_{I\dctra}$ are usually
\emph{not} exact.

 In the case of a finitely generated ideal $I\subset R$, the full
subcategory $R\Modl_{I\dctra}^\qs$ is closed under kernels, cokernels,
and infinite direct products in $R\Modl_{I\dctra}$ and in $R\Modl$,
and also under subobjects and quotient objects in $R\Modl_{I\dctra}$
\,\cite[Lemma~1.2]{Pdc}.
 (However, the full subcategory of quotseparated contramodules is
\emph{not} closed under extensions; in fact, any $I$\+contramodule
$R$\+module is an extension of two quotseparated $I$\+contramodule
$R$\+modules~\cite[Proposition~1.6]{Pdc}.)
 It follows that $R\Modl_{I\dctra}^\qs$ is an abelian category.
 In fact, $R\Modl_{I\dctra}$ is also a locally $\aleph_1$\+presentable
abelian category with enough projective objects~\cite[end of
Section~1]{Pdc}.
 The identity inclusion functors $R\Modl_{I\dctra}^\qs\rarrow
R\Modl_{I\dctra}$ and $R\Modl_{I\dctra}^\qs\rarrow R\Modl$ are
exact and preserve infinite products (hence also all limits) and
$\aleph_1$\+directed colimits.
 The inclusion functor $R\Modl_{I\dctra}^\qs\rarrow R\Modl$ does
\emph{not} preserve infinite direct sums.
 The infinite direct sum functors in $R\Modl_{I\dctra}^\qs$ are
usually \emph{not} exact (the global dimension~$1$ case is one major
exception~\cite[Remark~1.2.1]{Pweak}).

 To any complete, separated topological ring $\fR$ where open right
ideals form a base of neighborhoods of zero, the abelian category
$\fR\Contra$ of \emph{left\/ $\fR$\+contramodules} is naturally
assigned~\cite[Remark~A.3]{Psemi}, \cite[Section~1.2]{Pweak},
\cite[Section~2.1]{Prev}, \cite[Sections~1.2 and~5]{PR},
\cite[Section~2.7]{Pcoun}, \cite[Sections~2.6\+-2.7]{Pproperf},
\cite[Section~6]{PS1}, \cite[Section~1]{PS3}.
 The category $\fR\Contra$ comes equipped with a natural exact,
faithful forgetful functor $\fR\Contra\rarrow\fR\Modl$, preserving
infinite products (hence all limits).

 In the case of an adic topological ring $R$, the composition of
forgetful functors $\fR\Contra\rarrow\fR\Modl\rarrow R\Modl$ is
fully faithful, and its essential image coincides with the full
subcategory of quotseparated contramodule
$R$\+modules~\cite[Proposition~1.5]{Pdc}.
 So there is a natural equivalence (in fact, isomorphism)
of categories $\fR\Contra\simeq R\Modl_\ctra^\qs$.

\begin{lem} \label{Hom-is-a-quotseparated-contramodule}
 Let $R$ be an adic topological ring.  Then \par
\textup{(a)} for any torsion $R$\+module\/ $\cM$ and any $R$\+module
$P$, the $R$\+module\/ $\Hom_R(\cM,P)$ is a quotseparated
contramodule $R$\+module (in fact, a separated and complete
$R$\+module); \par
\textup{(b)} for any $R$\+module $M$ and any complete, separated
$R$\+module\/ $\fP$, the $R$\+module\/ $\Hom_R(M,\fP)$ is complete
and separated; \par
\textup{(c)} for any $R$\+module $M$ and any quotseparated
contramodule $R$\+module\/ $\fP$, the $R$\+module\/ $\Hom_R(M,\fP)$
is a quotseparated contramodule $R$\+module.
\end{lem}

\begin{proof}
 Part~(a) is the assertion in parentheses in~\cite[Lemma~1.1(b)]{Ppt}.
 To prove parts~(b) and~(c), represent $M$ as a quotient $R$\+module
of a free $R$\+module~$F$.
 Then $\Hom_R(M,\fP)$ is an $R$\+submodule of the $R$\+module
$\Hom_R(F,\fP)$, which an infinite product of copies of~$\fP$.
 Furthermore, $\Hom_R(M,\fP)$ is a contramodule $R$\+module by
Lemma~\ref{Hom-is-a-contramodule}(b).
 Now assertion~(b) follows from the fact that the full subcategory
of separated $R$\+modules is closed under subobjects and infinite
products in $R\Modl$.
 Since the full subcategory $R\Modl_\ctra^\qs$ is closed under
subobjects, quotient objects, and infinite products in
$R\Modl_\ctra$, part~(c) follows as well.
 For an argument based on the notion of a contramodule over
a topological ring, see~\cite[Section~B.2]{Pweak}.
\end{proof}

\subsection{Injective torsion modules, projective and flat
contramodules} \label{prelim-inj-proj-flat-torsion-contra-subsecn}
 For any abelian category $\sA$, we denote by $\sA^\inj\subset\sA$
the full subcategory of injective objects in~$\sA$.
 More generally, the same notation applies to any exact category~$\sA$
(in the sense of Quillen).
 Dually, the notation $\sB_\prj\subset\sB$ stands for the full
subcategory of projective objects in an abelian or exact
category~$\sB$.

 Let $I$ be an ideal in a commutative ring~$R$.
 Given an $R$\+module $M$, we denote by $\Gamma_I(M)\subset M$
the submodule of all $I$\+torsion elements in~$M$.
 The functor $M\longmapsto\Gamma_I(M)\:R\Modl\rarrow R\Modl_{I\dtors}$
is right adjoint to the identity inclusion functor
$R\Modl_{I\dtors}\rarrow R\Modl$.

 As any Grothendieck category, the abelian category $R\Modl_{I\dtors}$
has enough injective objects.
 For every injective $R$\+module $K$, the torsion $R$\+module
$\Gamma_I(K)$ is injective as an object of $R\Modl_{I\dtors}$.
 Indeed, the maximal torsion submodule functor $\Gamma_I$ is right
adjoint to an exact functor, so it takes injective objects to
injective objects.
 There are enough injective objects of the form $\Gamma_I(K)$ in
$R\Modl_{I\dtors}$, where $K\in R\Modl^\inj$.
 Hence an $I$\+torsion $R$\+module $\cK$ is injective in
$R\Modl_{I\dtors}$ if and only if it is a direct summand of
the $I$\+torsion $R$\+module $\Gamma_I(K)$ for some injective
$R$\+module~$K$.

 Let $R$ be an adic topological ring with an ideal of definition
$I\subset R$.
 Clearly, the submodule $\Gamma(M)=\Gamma_I(M)\subset M$ only depends
on the adic topology on $R$, and not on the choice of a particular
ideal of definition~$I$.
 Using Baer's criterion of injectivity, one can easily check that
a torsion $R$\+module $\cK$ is an injective object of $R\Modl_\tors$
if and only if the $R/I^n$\+module ${}_{I^n}\cK$ is injective
for every $n\ge1$.
 Equivalently, a torsion $R$\+module $\cK$ is injective (as a torsion
$R$\+module) if and only if the $R/J$\+module ${}_J\cK$ is injective
for every ideal of definition $J\subset R$.
 If this is the case, then the $R/J$\+module ${}_J\cK$ is also
injective for every open ideal $J\subset R$.

 Let $\fR=\Lambda(R)$ be the adic completion of the ring~$R$.
 Then the $R$\+module structure of every torsion $R$\+module can be
extended to a torsion $\fR$\+module structure in a unique way.
 A map of torsion $\fR$\+modules is $\fR$\+linear if and only if it is
$R$\+linear.
 So the abelian categories $R\Modl_\tors$ and $\fR\Modl_\tors$ are
naturally equivalent (in fact, isomorphic).

 When the ring $R$ is Noetherian (as an abstract ring), one can
easily prove using the Artin--Rees lemma that any injective object
of $R\Modl_\tors$ is also injective in $R\Modl$.
 Thus one has $R\Modl_\tors^\inj=R\Modl_\tors\cap R\Modl^\inj$
in this case.
 The following lemma summarizes the observations above.

\begin{lem} \label{injective-torsion-modules-characterizations}
 Let $R$ be a Noetherian commutative ring endowed with an adic
topology, let $I\subset R$ be an ideal of definition of the topological
ring $R$, and let\/ $\cK$ be a torsion $R$\+module.
 Then the following six conditions are equivalent:
\begin{enumerate}
\item $\cK$ is an injective object of $R\Modl_\tors$;
\item there exists an injective $R$\+module $K$ such that\/ $\cK$ is
a direct summand of the torsion $R$\+module\/ $\Gamma_I(K)$;
\item for every $n\ge1$, the $R/I^n$\+module ${}_{I^n}\cK$ is
injective;
\item for every ideal of definition $J\subset R$, the $R/J$\+module
${}_J\cK$ is injective;
\item for every open ideal $J\subset R$, the $R/J$\+module
${}_J\cK$ is injective;
\item $\cK$ is injective as an object of $R\Modl$;
\item $\cK$ is injective as an object of\/ $\fR\Modl$. \qed
\end{enumerate}
\end{lem}

\begin{cor} \label{Noetherian-torsion-Ext-isomorphism-cor}
 Let $R$ be a Noetherian commutative ring endowed with an adic
topology.
 Then the fully faithful inclusion of abelian categories $R\Modl_\tors
\rarrow R\Modl$ induces isomorphisms on all the Ext groups/modules.
\end{cor}

\begin{proof}
 The assertion follows from
Lemma~\ref{injective-torsion-modules-characterizations}%
\,(1)\,$\Leftrightarrow$\,(6).
 For a stronger result concerning the full-and-faithfulness of
the induced functor between the unbounded derived categories,
see~\cite[Theorem~1.3]{Pmgm} or
Theorem~\ref{torsion-cohomology-modules-triang-equivalence} below.
 The Noetherianity assumption can be relaxed (by replacing it with
the weak proregularity assumption), but it \emph{cannot} be
simply dropped; see~\cite[Theorem~4.1]{Pdc}.
\end{proof}

 For any ideal $I$ in a commutative ring $R$, the identity inclusion
functor $R\Modl_{I\dctra}\allowbreak\rarrow R\Modl$ has a left adjoint
functor $\Delta_I\:R\Modl\rarrow R\Modl_{I\dctra}$
\,\cite[Examples~4.1(2\+-3)]{PR}, \cite[Example~1.3(4)]{Pper}.
 An explicit construction of the functor $\Delta_I$ in the case of
a finitely generated ideal $I\subset R$ can be found
in~\cite[Proposition~2.1]{Pmgm}; see~\cite[Sections~6\+-7]{Pcta}
for a more detailed discussion.

 For every projective $R$\+module $P$, the $I$\+contramodule
$R$\+module $\Delta_I(P)$ is projective as an object of
$R\Modl_{I\dctra}$.
 Indeed, the functor $\Delta_I$ is left adjoint to an exact functor,
so it takes projective objects to projective objects.
 The $R$\+modules $\Delta_I(F)$, where $F$ ranges over the free
$R$\+modules, are called the \emph{free $I$\+contramodule
$R$\+modules}.
 There are enough free $I$\+contramodule $R$\+modules in
$R\Modl_{I\dctra}$.
 Hence an $I$\+contramodule $R$\+module $\fP$ is projective in
$R\Modl_{I\dctra}$ if and only if it is a direct summand of
some free $I$\+contramodule $R$\+module.

 Let $R$ be an adic topological ring.
 Then the full subcategory $R\Modl_{I\dctra}\subset R\Modl$ does not
depend on the choice of an ideal of definition $I\subset R$ of
the topological ring~$R$.
 So we will write $\Delta=\Delta_I$, and speak of \emph{free
contramodule $R$\+modules} presuming any ideal of definition~$I$
of the adic topology on~$R$.

 The $0$\+th left derived functor $\boL_0\Lambda$ of the (neither left
nor right exact) adic completion functor $\Lambda$ is the left adjoint
functor $\boL_0\Lambda\:R\Modl\rarrow R\Modl_\ctra^\qs$ to
the identity inclusion functor $R\Modl_\ctra^\qs\rarrow
R\Modl$ \,\cite[Proposition~1.3]{Pdc}.
 Similarly to the discussion above, for every projective $R$\+module
$P$, the quotseparated contramodule $R$\+module $\boL_0\Lambda(P)$
is projective as an object of $R\Modl_\ctra^\qs$.
 Notice that, by the definition of the $0$\+th derived functor, one
has $\boL_0\Lambda(P)=\Lambda(P)$ for all projective $R$\+modules~$P$.

 The $R$\+modules $\Lambda(F)=\boL_0\Lambda(F)$, where $F$ ranges over
the free $R$\+modules, are called the \emph{free quotseparated
contramodule $R$\+modules}.
 There are enough free quotseparated contramodule $R$\+modules in
$R\Modl_\ctra^\qs$.
 Hence quotseparated contramodule $R$\+module $\fP$ is projective in
$R\Modl_\ctra^\qs$ if and only if it is a direct summand of some
free quotseparated contramodule $R$\+module.

 Recall some notation from the theory of contramodules over
topological rings~\cite[Section~1.2]{Pweak},
\cite[Sections~1.2 and~5]{PR}, \cite[Section~2.1]{Prev},
\cite[Section~2.7]{Pcoun}, \cite[Sections~2.6\+-2.7]{Pproperf},
\cite[Section~6]{PS1}, \cite[Section~1]{PS3}.
 For any abelian group $A$ and any set $X$, we denote by
$A[X]=A^{(X)}$ the direct sum of copies of $A$ indexed by~$X$.
 The elements of $A[X]$ are interpreted as finite formal linear
combinations of elements of $X$ with the coefficients in~$A$.
 For any complete, separated topological abelian group $\fA$ where
open subgroups form a base of neighborhoods of zero, we put
$\fA[[X]]=\varprojlim_{\fU\subset\fA}(\fA/\fU)[X]$, where
$\fU$ ranges over the open subgroups of~$\fA$.
 The elements of $\fA[[X]]$ are interpreted as infinite formal
linear combinations of elements of $X$ with the families of
coefficients converging to zero in~$\fA$.

 For any complete, separated topological ring $\fR$ where open
right ideals form a base of neighborhoods of zero, and for any
set $X$, the set/abelian group $\fR[[X]]$ has a natural structure
of left $\fR$\+contramodule, called the \emph{free $\fR$\+contramodule
spanned by~$X$}.
 The free $\fR$\+contramodules are projective objects in $\fR\Contra$,
there are enough of them, and all projective $\fR$\+contramodules
are direct summands of free ones.

 Given an adic topological ring $R$ with an ideal of definition
$I\subset R$, consider the adic completion $\fR=\Lambda(R)$ of
the ring~$R$.
 Then, for any set $X$, one has an obvious isomorphism of
$R$\+modules $\fR[[X]]=\varprojlim_{n\ge1}((R/I^n)[X])=\Lambda(R[X])$,
which is also an isomorphism in $\fR\Contra\simeq R\Modl_\ctra^\qs$.
 Thus the free quotseparated contramodule $R$\+modules are the same
things as free $\fR$\+contramodules in the sense of the general
theory of contramodules over topological rings, and our terminology
is consistent.

 When the ring $R$ is Noetherian (as an abstract ring), one has
$R\Modl_\ctra=R\Modl_\ctra^\qs$ (see
Section~\ref{prelim-contramodules-subsecn}).
 Hence a natural isomorphism $\Delta\simeq\boL_0\Lambda$ of functors
$R\Modl\rarrow R\Modl_\ctra$.

\begin{lem} \label{derived-completion-reduction-isomorphism}
 Let $R$ be an adic topological ring, let $I\subset R$ be an open
ideal, and let $P$ be an $R$\+module.
 Then the natural adjunction maps $P\rarrow\Delta(P)\rarrow
\boL_0\Lambda(P)\rarrow\Lambda(P)$ induce isomorphisms of $R/I$\+modules
$R/I\ot_RP\simeq R/I\ot_R\Delta(P)\simeq R/I\ot_R\boL_0\Lambda(P)
\simeq R/I\ot_R\Lambda(P)$.
\end{lem}

\begin{proof}
 Follows from the facts that the functor $P\longmapsto R/I\ot_RP$ is
left adjoint to the inclusion of full subcategory $R/I\Modl\rarrow
R\Modl$, and that all $R/I$\+modules are separated and complete
$R$\+modules.
\end{proof}

\begin{cor} \label{derived-completion-tensor-with-torsion-isom-cor}
 Let $R$ be an adic topological ring, let\/ $\cM$ be a torsion
$R$\+module, and let $P$ be an $R$\+module.
 Then the natural adjunction maps $P\rarrow\Delta(P)\rarrow
\boL_0\Lambda(P)\rarrow\Lambda(P)$ induce isomorphisms of $R$\+modules\/
$\cM\ot_RP\simeq\cM\ot_R\Delta(P)\simeq\cM\ot_R\boL_0\Lambda(P)
\simeq\cM\ot_R\Lambda(P)$.
\end{cor}

\begin{proof}
 Any torsion $R$\+module $\cM$ is a direct union of $R/I$\+modules,
where $I$ ranges over the ideals of definition of~$R$.
 Since the tensor product functor~$\ot_R$ preserves direct limits,
it suffices to consider the case of an $R/I$\+module $\cM=M$.
 In that case, the assertion follows from
Lemma~\ref{derived-completion-reduction-isomorphism} in view of
the natural isomorphism $M\ot_RQ\simeq M\ot_{R/I}(R/I\ot_RQ)$ for
any $R$\+module~$Q$.
\end{proof}

\begin{lem} \label{projective-contramodules-characterizations}
 Let $R$ be a Noetherian commutative ring endowed with an adic
topology, let $I\subset R$ be an ideal of definition of the topological
ring $R$, and let\/ $\fP$ be a contramodule $R$\+module.
 Then the following seven conditions are equivalent:
\begin{enumerate}
\item $\fP$ is a projective object of $R\Modl_\ctra=R\Modl_\ctra^\qs$;
\item there exists a free $R$\+module $F$ such that\/ $\fP$ is
a direct summand of the contramodule $R$\+module\/
$\Delta(F)=\boL_0\Lambda(F)=\Lambda(F)$;
\item for every $n\ge1$, the $R/I^n$\+module\/ $\fP/I^n\fP$
is projective;
\item for every ideal of definition $J\subset R$, the $R/J$\+module\/
$\fP/J\fP$ is projective;
\item for every open ideal $J\subset R$, the $R/J$\+module\/
$\fP/J\fP$ is projective;
\item for every $n\ge1$, the $R/I^n$\+module\/ $\fP/I^n\fP$
is flat, and the $R/I$\+module\/ $\fP/I\fP$ is projective;
\item $\fP$ is a flat $R$\+module and the $R/I$\+module $\fP/I\fP$
is projective.
\end{enumerate}
\end{lem}

\begin{proof}
 (1)~$\Longleftrightarrow$~(2) Explained above.
 No Noetherianity assumption on $R$ is needed here.
 (In the context of $R\Modl_{I\dctra}$ and $\Delta_I$, the assertion
even holds for infinitely generated ideals $I\subset R$.)

 (3)~$\Longleftrightarrow$~(4)~$\Longleftrightarrow$~(5) This is
obvious and holds for any adic topological ring~$R$ (does not need
the Noetherianity assumption).

 (3)~$\Longleftrightarrow$~(6) This is a particular case
of a quite general lemma about an associative ring with
a nilpotent ideal~\cite[Lemma~B.10.2]{Pweak}.

 (2)~$\Longrightarrow$~(3) Both in the contexts of $R\Modl_\ctra$
and $R\Modl_\ctra^\qs$, this holds for any adic topological ring $R$
and does not need the Noetherianity assumption.
 The fact that the $R/I$\+module $R/I\ot_R\fF$ is free whenever $\fF$
is either a free contramodule $R$\+module or a free quotseparated
contramodule $R$\+module follows from
Lemma~\ref{derived-completion-reduction-isomorphism} (applied to
a free $R$\+module~$M$).

 (1)~$\Longleftrightarrow$~(3) In the context of $R\Modl_\ctra^\qs$,
this also holds for any adic topological ring $R$ and does not need
the Noetherianity assumption.
 See~\cite[Corollary~E.1.10(a)]{Pcosh} for an even more general result.

 (7)~$\Longrightarrow$~(6) Obvious.
 
 (1), (3), or~(6)~$\Longrightarrow$~(7) The assumption that $R$ is
Noetherian is essential here.
 See~\cite[Corollary~B.8.2]{Pweak} or~\cite[Theorem~10.5]{Pcta} for
a complete proof of (1) $\Longleftrightarrow$ (3) $\Longleftrightarrow$
(6) $\Longleftrightarrow$~(7) for Noetherian rings~$R$.
\end{proof}

 A discussion of the functor of \emph{contratensor product} $\odot_\fR$
of discrete right $\fR$\+modules and left $\fR$\+contramodules, for
a complete, separated topological ring $\fR$ where open right ideals
form a base of neighborhoods of zero, can be found
in~\cite[Section~3.3]{Prev}, \cite[Section~5]{PR},
\cite[Section~2.8]{Pcoun}, \cite[Sections~2.8]{Pproperf},
\cite[Section~7.2]{PS1}, \cite[Section~1]{PS3}.
 In the context of an adic topological ring $R$ and its adic
completion $\fR$, the natural maps
$$
 \cM\ot_R\fP\lrarrow\cM\ot_\fR\fP\lrarrow\cM\odot_\fR\fP
$$
are isomorphisms for any torsion $R$\+module $\cM$ and any
quotseparated contramodule $R$\+module~$\fP$.
 This assertion, which is easy to check directly, follows also
from the fact that the forgetful functor $\fR\Contra\rarrow R\Modl$
is fully faithful~\cite[Lemma~7.11]{PS1}.
 For this reason, we will not dwell on the contratensor product
functor~$\odot_\fR$, working with the tensor product~$\ot_R$
or~$\ot_\fR$ instead.

 So, for an adic topological ring $R$, we are interested in the tensor
product functor
$$
 \ot_R\:R\Modl_\tors\times R\Modl_\ctra\lrarrow R\Modl_\tors,
$$
which is well-defined by Lemma~\ref{tensor-product-is-torsion}.
 A contramodule $R$\+module $\fF$ is said to be
(\emph{contra})\emph{flat} (as a contramodule $R$\+module) if
the tensor product functor ${-}\ot_R\fF\:\allowbreak R\Modl_\tors
\rarrow R\Modl_\tors$ is exact.
 Let $I\subset R$ be an ideal of definition.
 Since any short exact sequence of torsion $R$\+modules is
a direct limit of short exact sequences of $R/I^n$\+modules,
$n\ge1$, a contramodule $R$\+module $\fF$ is flat if and only if
the $R/I^n$\+module $\fF/I^n\fF$ is flat for all $n\ge1$.

 For any set $X$ and any torsion $R$\+module $\cM$, we have
a natural isomorphism of torsion $R$\+modules
$$
 \cM\ot_R\fR[[X]]\simeq\cM[X]=\cM^{(X)}.
$$
 Therefore, the free quotseparated contramodule $R$\+modules are
flat (as contramodule $R$\+modules), and it follows that
the projective quotseparated contramodule $R$\+modules are flat, too.

 The projective contramodule $R$\+modules, i.~e., the projective
objects of the abelian category $R\Modl_\ctra$ are flat contramodule
$R$\+modules as well~\cite[Section~3]{Ppt} (see also the proof of
Lemma~\ref{projective-contramodules-characterizations}%
\,(2)\,$\Rightarrow$\,(3) above).
 So flat contramodule $R$\+modules \emph{need not} be separated if
they are not quotseparated.
 However, the following assertion holds.

\begin{lem} \label{flat-quotseparated-contramodules-separated}
 For any adic topological ring $R$, any flat quotseparated
contramodule $R$\+module is (complete and) separated.
 In particular, if the ring $R$ is Noetherian (as an abstract ring),
then all flat contramodule $R$\+modules are separated.
\end{lem}

\begin{proof}
 This is a special case of the much more general result
of~\cite[Corollary~E.1.7]{Pcosh}, or of the even more
general~\cite[Corollary~6.15]{PR}.
 The Noetherian adic case can be found
in~\cite[proof of Lemma~B.9.2]{Pweak}
or~\cite[Corollary~10.3(b)]{Pcta}.
\end{proof}

 The following corollary is a partial generalization
of~\cite[Lemma~3.5 and Proposition~3.6]{PSY}.

\begin{cor} \label{flat-reductions-derived-Lambda-is-underived-cor}
 Let $R$ be an adic topological ring with an ideal of definition
$I\subset R$, and let $F$ be an $R$\+module.
 Then the following four conditions are equivalent:
\begin{enumerate}
\item the $R/I^n$\+module $F/I^nF$ is flat for every $n\ge1$;
\item the contramodule $R$\+module\/ $\Delta(F)$ is flat;
\item the contramodule $R$\+module\/ $\boL_0\Lambda(F)$ is flat;
\item the contramodule $R$\+module\/ $\Lambda(F)$ is flat.
\end{enumerate}
 If any one of the conditions~\textup{(1\+-4)} holds, then
the $R$\+module\/ $\boL_0\Lambda(F)$ is complete and separated,
and the natural $R$\+module map\/ $\boL_0\Lambda(F)
\rarrow\Lambda(F)$ is an isomorphism.
 Consequently, the completion map $F\rarrow\Lambda(F)$ induces
an isomorphism of $R$\+modules\/ $\Hom_R(\Lambda(F),\fP)\simeq
\Hom_R(F,\fP)$ for any quotseparated contramodule $R$\+module\/ $\fP$
in this case.
\end{cor}

\begin{proof}
 The equivalence of the four conditions~(1\+-4) follows from
Lemma~\ref{derived-completion-reduction-isomorphism}.
 If condition~(3) holds, then the $R$\+module $\boL_0\Lambda(F)$ is
complete and separated by
Lemma~\ref{flat-quotseparated-contramodules-separated}.
 Passing to the projective limit of the isomorphisms from
Lemma~\ref{derived-completion-reduction-isomorphism} over the integers
$n\ge1$, we conclude that $\boL_0\Lambda(F)\simeq\Lambda(F)$.
 Since the functor $\boL_0\Lambda\:R\Modl\rarrow R\Modl_\ctra^\qs$
is left adjoint to the inclusion functor $R\Modl_\ctra^\qs\rarrow
R\Modl$, the last assertion of the corollary follows.
\end{proof}

\begin{cor} \label{flat-completion-of-quotseparated-implies-separated}
 Let $R$ be an adic topological ring with an ideal of definition
$I\subset R$, and let\/ $\fF$ be a quotseparated contramodule
$R$\+module.
 Assume that the separated contramodule $R$\+module\/ $\Lambda(\fF)$
is flat.
 Then the $R$\+module\/ $\fF$ is separated, so\/
 $\fF\simeq\Lambda(\fF)$.
\end{cor}

\begin{proof}
 This can be deduced from the next-to-last assertion of
Corollary~\ref{flat-reductions-derived-Lambda-is-underived-cor}
by noticing that $\fF\simeq\boL_0\Lambda(\fF)$ for any quotseparated
contramodule $R$\+module~$\fF$.
 Alternatively, one can say that the quotseparated contramodule
$R$\+module $\fF$ is flat by
Corollary~\ref{flat-reductions-derived-Lambda-is-underived-cor}%
\,(4)\,$\Rightarrow$\,(1), and refer to
Lemma~\ref{flat-quotseparated-contramodules-separated}.
\end{proof}

\begin{lem} \label{flat-quotseparated-contramodules-well-behaved}
 Let $R$ be an adic topological ring and\/ $\cM$ be a torsion
$R$\+module.
 Then \par
\textup{(a)} the full subcategory of flat quotseparated contramodule
$R$\+modules is closed under extensions and kernels of epimorphisms
in $R\Modl_\ctra^\qs$; \par
\textup{(b)} for any short exact sequence of quotseparated contramodule
$R$\+modules\/ $0\rarrow\fH\rarrow\fG\rarrow\fF\rarrow0$ with
a flat quotseparated contramodule $R$\+module\/ $\fF$, the short
sequence of torsion $R$\+modules\/ $0\rarrow\cM\ot_R\fH\rarrow
\cM\ot_R\fG\rarrow\cM\ot_R\fF\rarrow0$ is exact.
\end{lem}

\begin{proof}
 The argument is based on the category equivalence/isomorphism
$R\Modl_\ctra^\qs\simeq\fR\Contra$.
 Part~(a) is a special case of~\cite[Lemmas~E.1.4 and~E.1.5]{Pcosh}.
 A further generalization can be found in~\cite[Corollaries~6.8
and~6.13]{PR}.
 Part~(b) is a special case of~\cite[Lemma~6.10]{PR}.
\end{proof}

\begin{cor} \label{completed-tensor-product-exact-on-flat-qs}
 Let $R$ be an adic topological ring and $M$ be an $R$\+module.
 Then, for any short exact sequence of quotseparated contramodule
$R$\+modules\/ $0\rarrow\fH\rarrow\fG\rarrow\fF\rarrow0$ with
a flat quotseparated contramodule $R$\+module\/ $\fF$, the short
sequence of quotseparated contramodule $R$\+modules\/
$0\rarrow\Lambda(M\ot_R\fH)\rarrow\Lambda(M\ot_R\fG)\rarrow
\Lambda(M\ot_R\fF)\rarrow0$ is exact.
\end{cor}

\begin{proof}
 Follows from
Lemma~\ref{flat-quotseparated-contramodules-well-behaved}(b)
together with the fact of exactness of projective limits of countable
sequences of surjective maps of abelian groups.
\end{proof}

 The definition of an \emph{adically Noetherian} adic topological
ring was given in Section~\ref{adic-rings-subsecn}.
 Following~\cite[Section~9]{Ppt}, we will say that an adic topological
ring $R$ is \emph{adically coherent} if the ring $R/I$ is coherent
for \emph{every} ideal of defintion $I\subset R$.
 Checking this condition for just one ideal of definition $I$ is
\emph{not} sufficient; see~\cite[Remark~9.2]{Ppt}.
 However, if the ring $R$ is coherent as an abstract ring,
then it is also adically coherent as an adic topological
ring~\cite[Lemma~9.1(b)]{Ppt}.

\begin{lem} \label{sums-of-inj-tors-prods-of-flat-contramods}
\textup{(a)} Let $R$ be an adically Noetherian adic topological ring.
 Then the full subcategory of injective objects $R\Modl_\tors^\inj$
is closed under filtered direct limits (and in particular, under
infinite direct sums) in the abelian category $R\Modl_\tors$. \par
\textup{(b)} Let $R$ be an adically coherent adic topological ring.
 Then the full subcategory of flat contramodule $R$\+modules is closed
under infinite direct products in the abelian category $R\Modl_\ctra$.
\end{lem}

\begin{proof}
 Part~(a): following the discussion in the beginning of this section,
a torsion $R$\+module $\cK$ is injective (as an object of
$R\Modl_\tors$) if and only if the $R/I$\+module ${}_I\cK$ is injective
for every ideal of definition $I\subset R$.
 Clearly, the reduction functor $\cK\longmapsto{}_I\cK\:R\Modl_\tors
\rarrow R/I\Modl$ preserves infinite direct sums.
 Since we assume all our ideals of definition to be finitely generated,
this reduction functor also preserves filtered direct limits.
 It remains to recall that filtered direct limits (and in particular,
infinite direct sums) of injective left modules over a left Noetherian
associative ring are injective.
 In particular, this applies to the ring~$R/I$.

 Part~(b): this is~\cite[Lemma~9.8]{Ppt}.
 Following the discussion above in this section, a contramodule
$R$\+module $\fF$ is flat (as a contramodule $R$\+module) if and only
if the $R/I$\+module $\fF/I\fF$ is flat for every ideal of definition $I\subset R$.
 Since we assume all our ideals of definition to be finitely generated,
the reduction functor $\fF\longmapsto\fF/I\fF\:R\Modl_\ctra\rarrow
R/I\Modl$ preserves infinite products.
 It remains to recall that infinite products of flat left modules over
a right coherent associative ring are flat.
\end{proof}

 When the ring $R$ is Noetherian (as an abstract ring), one can show
that a contramodule $R$\+module is flat as a contramodule
$R$\+module if and only if it is flat as $R$\+module.
 This is the main assertion of the following lemma.

\begin{lem} \label{flat-contramodules-characterizations}
 Let $R$ be a Noetherian commutative ring endowed with an adic
topology, let $I\subset R$ be an ideal of definition of the topological
ring $R$, let\/ $\fR=\Lambda(R)$ be the adic completion of $R$, and
let\/ $\fF$ be a contramodule $R$\+module.
 Then the following seven conditions are equivalent:
\begin{enumerate}
\item $\fF$~is a flat contramodule $R$\+module;
\item $\fF$~is a flat contramodule\/ $\fR$\+module;
\item for every $n\ge1$, the $R/I^n$\+module\/ $\fF/I^n\fF$
is flat;
\item for every ideal of definition $J\subset R$, the $R/J$\+module\/
$\fF/J\fF$ is flat;
\item for every open ideal $J\subset R$, the $R/J$\+module\/
$\fF/J\fF$ is flat;
\item $\fF$~is a flat $R$\+module;
\item $\fF$~is a flat\/ $\fR$\+module.
\end{enumerate}
\end{lem}

\begin{proof}
 (1)~$\Longleftrightarrow$~(3) Explained above.
 No Noetherianity assumption on $R$ is needed here.

 (3)~$\Longleftrightarrow$~(4)~$\Longleftrightarrow$~(5) This is
obvious and holds for any adic topological ring~$R$ (does not need
the Noetherianity assumption).

 (1)~$\Longleftrightarrow$~(2) Without the Noetherianity assumption,
one has to be careful.
 There is a natural isomorphism of the categories of 
\emph{quotseparated} contramodules $R\Modl_\ctra^\qs
\simeq\fR\Modl_\ctra^\qs$, as both the categories are isomorphic to
$\fR\Contra$.
 Simply put, any quotseparated contramodule $R$\+module admits
a unique extension of its $R$\+module structure to an $\fR$\+module
structure, making it a quotseparated contramodule $\fR$\+module;
and any $R$\+module morphism of quotseparated contramodule
$R$\+modules is an $\fR$\+module morphism.
 However, a nonquotseparated contramodule $R$\+module does \emph{not}
admit an extension of its $R$\+module structure to an $\fR$\+module
structure, in general.
 Instead, the nonquotseparated contramodule $R$\+modules have natural
module structures over the ring $\Delta(R)$
\,\cite[Example~5.2(3)]{Pper}.

 Nevertheless, in the context of quotseparated contramodules,
the equivalence (1)~$\Longleftrightarrow$~(2) holds for any adic
topological ring $R$ and follows immediately from the equivalence
(1)~$\Longleftrightarrow$~(3).

 (6)~$\Longrightarrow$~(5) Obvious.
 
 (3)~$\Longrightarrow$~(6) The assumption that $R$ is Noetherian
is essential here.
 A proof can be found in~\cite[Lemma~B.9.2]{Pweak} or
in~\cite[Corollary~10.3(a)]{Pcta}.

 (2)~$\Longleftrightarrow$~(7) This is a particular case of
(1)~$\Longleftrightarrow$~(6).
\end{proof}

\begin{cor} \label{Noetherian-adic-flat-ring-maps-characterized-cor}
 Let $R$ and $S$ be Noetherian commutative rings endowed with adic
topologies, and let $f\:R\rarrow S$ be a tight continuous ring map.
 Then the following conditions are equivalent:
\begin{enumerate}
\item $f$~is flat as a tight continuous map of adic
topological rings;
\item $\Lambda(S)$ is a flat $R$\+module;
\item $\Lambda(S)$ is a flat\/ $\Lambda(R)$\+module.
\end{enumerate}
\end{cor}

\begin{proof}
 Notice that $\fS=\Lambda(S)$ is a separated and complete module
over the adic topological ring~$R$ (as the adic/completion topology
on $\fS$ coincides with the topology induced on $\fS$ by
the adic topology on $R$, due to the assumptions that~$f$ is tight
and continuous).
 Therefore, $\fS$ is a contramodule $R$\+module.
 Furthermore, for every open ideal $I\subset R$, the ideal
$Sf(I)\subset S$ is open in $S$, hence $S/Sf(I)\simeq
\fS/\fS\lambda_S(f(I))$.
 The assertion of the Corollary can be obtained by applying
Lemma~\ref{flat-contramodules-characterizations}%
\,(5)\,$\Leftrightarrow$\,(6)\,$\Leftrightarrow$\,(7) to
the contramodule $R$\+module $\fF=\fS$.
\end{proof}

\begin{prop} \label{Noetherian-contramodule-Ext-isomorphism-prop}
 Let $R$ be a Noetherian commutative ring endowed with an adic
topology.
 Then the fully faithful inclusion of abelian categories
$R\Modl_\ctra^\qs=R\Modl_\ctra\rarrow R\Modl$ induces isomorphisms
on all the Ext groups/modules.
\end{prop}

\begin{proof}
 This is~\cite[Theorem~B.8.1]{Pweak}.
 For a stronger result concerning the full-and-faithfulness of
the induced functor between the unbounded derived categories,
see~\cite[Theorem~2.9]{Pmgm}.
 The Noetherianity assumption can be relaxed, but it \emph{cannot} be
simply dropped; see~\cite[Theorems~4.2 and~4.3]{Pdc}.
\end{proof}

\subsection{Change of scalars} \label{prelim-change-of-scalars-subsecn}
 Given a homomorphism of associative rings $f\:R\rarrow S$, we denote
by $f_*\:S\Modl\rarrow R\Modl$ the functor of \emph{restriction of
scalars}, taking every left $S$\+module $N$ to its underlying
left $R$\+module~$N$.
 The left adjoint functor to~$f_*$, taking every left $R$\+module $M$
to the left $S$\+module $S\ot_RM$, will be denoted by
$f^*\:R\Modl\rarrow S\Modl$; it is called the functor of
\emph{extension of scalars} with respect to~$f$.
 The right adjoint functor to~$f_*$, taking every left $R$\+module $P$
to the left $S$\+module $\Hom_R(S,P)$, will be denoted by
$f^!\:R\Modl\rarrow S\Modl$ and called the functor of
\emph{coextension of scalars}.

 Let $R$ and $S$ be adic topological rings, and let $f\:R\rarrow S$
be a ring homomorphism.
 Parts~(b) and~(d) of the following lemma are generalizations
of~\cite[Lemma~D.4.1(a)]{Pcosh}.

\begin{lem} \label{torsion-contra-restriction-of-scalars}
 Assume that the ring homomorphism~$f$ is continuous.
 Consider the functor of restriction of scalars
$f_*\:S\Modl\rarrow R\Modl$.
 Then \par
\textup{(a)} the forgetful functor~$f_*$ takes torsion $S$\+modules to
torsion $R$\+modules; \par
\textup{(b)} the forgetful functor~$f_*$ takes contramodule
$S$\+modules to contramodule $R$\+mod\-ules; \par
\textup{(c)} the forgetful functor~$f_*$ takes separated $S$\+modules
to separated $R$\+modules; \par
\textup{(d)} the forgetful functor~$f_*$ takes quotseparated
contramodule $S$\+modules to quotseparated contramodule $R$\+modules.
\end{lem}

\begin{proof}
 Recall that, by Lemma~\ref{continuous-ring-map-lemma}(2),
there exist an ideal of definition $I\subset R$ and an ideal of
definition $J\subset S$ such that $f(I)\subset J$.

 Part~(a): more generally, for any homomorphism of commutative rings
$f\:R\rarrow S$ and any two ideals $I\subset R$ and $J\subset S$ such
that $f(I)\subset J$, the functor $S\Modl\rarrow R\Modl$ takes
$J$\+torsion $S$\+modules to $I$\+torsion $R$\+modules.
 Indeed, for any $S$\+module $N$ and any fixed element $r\in R$,
an element $y\in N$ is $f(r)$\+torsion (as an element of
the $S$\+module~$N$) if and only if $y$~is $r$\+torsion (as
an element of the underlying $R$\+module of~$N$).

 Part~(b): once again, for any homomorphism of commutative rings
$f\:R\rarrow S$ and any two ideals $I\subset R$ and $J\subset S$ such
that $f(I)\subset J$, the functor $S\Modl\rarrow R\Modl$ takes
$J$\+contramodule $S$\+modules to $I$\+contramodule $R$\+modules.
 The point is that, for any $S$\+module $Q$ and any fixed element
$r\in R$, there is a natural isomorphism of the Ext groups
$\Ext^*_S(S[f(r)^{-1}],Q)\simeq\Ext^*_R(R[r^{-1}],Q)$, because
$R[r^{-1}]$ is a flat $R$\+module and $S\ot_RR[r^{-1}]\simeq
S[f(r)^{-1}]$.
 See, e.~g., \cite[Lemma~4.1(a)]{PSl1}; cf.~\cite[Lemma~2.1]{Pcta}.
 So the $S$\+module $Q$ is an $f(r)$\+contramodule if and only if
the underlying $R$\+module of $Q$ is an $r$\+contramodule.

 Part~(c): for every $S$\+module $C$ and every $n\ge1$, one has
$I^nC\subset J^nC\subset C$.
 Hence $\bigcap_{n\ge1}J^nC=0$ implies $\bigcap_{n\ge1}I^nC=0$.

 Part~(d) follows from parts~(b) and~(c).
 Recall that all contramodule $S$\+modules are complete, and
all separate complete $S$\+modules are contramodules, as
mentioned in Section~\ref{prelim-contramodules-subsecn}.
 So the class of separated and complete $S$\+modules coincides
with the class of separated contramodule $S$\+modules.

 Let $\fQ$ be a quotseparated contramodule $S$\+module.
 Then $\fQ$ is a contramodule $S$\+module and a quotient $S$\+module
of a separated contramodule $S$\+module~$\fC$.
 Then, by~(b) and~(c), the underlying $R$\+module of $\fQ$ is
a contramodule $R$\+module and a quotient $R$\+module of
a separated contramodule $R$\+module~$\fC$.
\end{proof}

\begin{cor} \label{finer-and-coarser-adic-topology-completeness-cor}
 Let $R$ be a commutative ring and $I\subset J\subset R$ be two
finitely generated ideals.
 Let $M$ be a $J$\+adically separated and complete $R$\+module.
 Then $M$ is also an $I$\+adically separated and complete $R$\+module.
\end{cor}

\begin{proof}
 Once again, an ``adically separated and complete module'' is
the same thing as a ``separated contramodule''.
 Denote by $R_I$ the ring $R$ endowed with the $I$\+adic topology and
by $R_J$ the ring $R$ endowed with the $J$\+adic topology.
 Then $f=\id_R\:R_I\rarrow R_J$ is a continuous homomorphism of
adic topological rings.
 Applying Lemma~\ref{torsion-contra-restriction-of-scalars}(b,c) to
this continuous ring homomorphism, we obtain the desired conclusion.
\end{proof}

 Let $f\:R\rarrow S$ be a continuous homomorphism of adic topological
rings.
 Then Lemma~\ref{torsion-contra-restriction-of-scalars}(a) provides
an exact, faithful functor of restriction of scalars
$$
 f_\diamond\:S\Modl_\tors\lrarrow R\Modl_\tors.
$$
 The functor~$f_\diamond$ has a right adjoint functor
$$
 f^\diamond\:R\Modl_\tors\lrarrow S\Modl_\tors,
$$
which can be computed by the formula $f^\diamond(\cM)=
\Gamma_J(\Hom_R(S,\cM))$ for any $\cM\in R\Modl_\tors$,
where $J$~is an ideal of definition in~$S$.
 Moreover, for any $R$\+module $M$ one has a natural isomorphism
of torsion $S$\+modules
$$
 f^\diamond(\Gamma_I(M))\simeq\Gamma_J(f^!(M)),
$$
where $I$~is an ideal of definition in~$R$.
 As a right adjoint functor to an exact functor,
the functor~$f^\diamond$ takes injective torsion $R$\+modules to
injective torsion $S$\+modules.
 The same assertions hold in the greater generality of a commutative
ring homomorphism $f\:R\rarrow S$ and any two ideals $I\subset R$,
\,$J\subset S$ such that $f(I)\subset J$ \,\cite[proof of
Lemma~15.1(a)]{Ppt}.
 See~\cite[Section~2.9]{Pproperf} and~\cite[Section~6]{Pcs} for
a generalization to the categories of discrete modules over topological
associative rings.

 We call~$f^\diamond$ the functor of \emph{coextension of scalars}
(for torsion modules, with respect to a continuous homomorphism of
adic topological rings).

 Furthermore, Lemma~\ref{torsion-contra-restriction-of-scalars}(b)
provides an exact, faithful functor of restriction of scalars
$$
 f_\#\:S\Modl_\ctra\lrarrow R\Modl_\ctra.
$$
 The functor~$f_\#$ has a left adjoint functor
$$
 f^\#\:R\Modl_\ctra\lrarrow S\Modl_\ctra,
$$
which can be computed by the formula $f^\#(\fP)=
\Delta_J(S\ot_R\fP)$ for any $\fP\in R\Modl_\ctra$,
where $J$~is an ideal of definition in~$S$.
 Moreover, for any $R$\+module $P$ one has a natural isomorphism
of contramodule $S$\+modules
$$
 f^\#(\Delta_I(P))\simeq\Delta_J(f^*(P)),
$$
where $I$~is an ideal of definition in~$R$.
 As a left adjoint functor to an exact functor,
the functor~$f^\#$ takes projective contramodule $R$\+modules to
projective contramodule $S$\+modules.
 Specifically, one has $f^\#(\Delta_I(R[X]))=\Delta_J(S[X])$
for any set~$X$ (so the functor~$f^\#$ also takes free
contramodule $R$\+modules to free contramodule $S$\+modules).
 The same assertions hold in the greater generality of a commutative
ring homomorphism $f\:R\rarrow S$ and any two ideals $I\subset R$,
\,$J\subset S$ such that $f(I)\subset J$ \,\cite[proof of
Lemma~15.1(b)]{Ppt}.

 We call~$f^\#$ the functor of \emph{contraextension of scalars}
(for contramodules, with respect to a continuous homomorphism of adic
topological rings).

 Finally, Lemma~\ref{torsion-contra-restriction-of-scalars}(d)
provides an exact, faithful functor of restriction of scalars
$$
 f_\sharp\:S\Modl_\ctra^\qs\lrarrow R\Modl_\ctra^\qs.
$$
 The functor~$f_\sharp$ has a left adjoint functor
$$
 f^\sharp\:R\Modl_\ctra^\qs\lrarrow S\Modl_\ctra^\qs,
$$
which can be computed by the formula $f^\sharp(\fP)=
\boL_0\Lambda_J(S\ot_R\fP)$ for any $\fP\in R\Modl_\ctra^\qs$,
where $J$~is an ideal of definition in~$S$ \,\cite[proof of
Lemma~15.1(c)]{Ppt}.
 Moreover, for any $R$\+module $P$ one has a natural isomorphism
of quotseparated contramodule $S$\+modules
$$
 f^\sharp(\boL_0\Lambda_I(P))\simeq\boL_0\Lambda_J(f^*(P)),
$$
where $I$~is an ideal of definition in~$R$.
 As a left adjoint functor to an exact functor,
the functor~$f^\sharp$ takes projective quotseparated contramodule
$R$\+modules to projective quotseparated contramodule $S$\+modules.
 Specifically, one has $f^\sharp(\Lambda_I(R[X]))=\Lambda_J(S[X])$
for any set~$X$ (so the functor~$f^\sharp$ also takes free
quotseparated contramodule $R$\+modules to free quotseparated
contramodule $S$\+modules).
 See~\cite[Section~2.9]{Pproperf} and~\cite[Section~6]{Pcs} for
a generalization to the categories of contramodules over topological
associative rings.

 By construction, the functors $f_\sharp$ and~$f_\#$ agree with
each other, in the obvious sense.
 The functors $f^\sharp$ and~$f^\#$ do \emph{not} agree, but compare as
follows: for any quotseparated contramodule $R$\+module $\fP$, there is
a natural isomorphism of quotseparated contramodule $S$\+modules
$$
 f^\sharp(\fP)\simeq\boL_0\Lambda_J(f^\#(\fP)),
$$
where $J$ is an ideal of definition in~$S$.

 The functor~$f^\sharp$ is called the functor of \emph{contraextension
of scalars} (for quotseparated contramodules, with respect to
a continuous homomorphism of adic topological rings).

\begin{lem} \label{map-of-adic-topol-rings-tensor-product-adjunction}
 Let $f\:R\rarrow S$ be a continuous homomorphism of adic topological
rings.
 Then \par
\textup{(a)} for any contramodule $R$\+module\/ $\fP$ and any torsion
$S$\+module\/ $\cN$, there is a natural isomorphism of $S$\+modules
$$
 f^\#(\fP)\ot_S\cN\simeq\fP\ot_Rf_\diamond(\cN);
$$ \par
\textup{(b)} for any quotseparated contramodule $R$\+module\/ $\fP$ and
any torsion $S$\+module\/ $\cN$, there is a natural isomorphism of
$S$\+modules
$$
 f^\sharp(\fP)\ot_S\cN\simeq\fP\ot_Rf_\diamond(\cN).
$$
\end{lem}

\begin{proof}
 Part~(a): both the left-hand side and the right-hand side of
the desired isomorphism are torsion $S$\+modules (cf.\
Lemma~\ref{tensor-product-is-torsion}).
 Now for any torsion $S$\+module $\cL$ we have natural isomorphisms
of abelian groups/$S$\+modules
\begin{multline*}
 \Hom_S(\fP\ot_Rf_\diamond(\cN),\>\cL)\simeq
 \Hom_R(\fP,\Hom_S(\cN,\cL)) \\ \simeq
 \Hom_S(f^\sharp(\fP),\Hom_S(\cN,\cL))\simeq
 \Hom_S(f^\sharp(\fP)\ot_S\cN,\>\cL),
\end{multline*}
where the second (middle, adjunction) isomorphism holds because
$\Hom_S(\cN,\cL)$ is a contramodule $S$\+module by
Lemma~\ref{Hom-is-a-contramodule}(a).

 The proof of part~(b) is similar and based on the fact that
$\Hom_S(\cN,\cL)$ is a quotseparated contramodule $S$\+module by
Lemma~\ref{Hom-is-a-quotseparated-contramodule}(a).
 Alternatively, part~(b) can be deduced from part~(a) using
the comparison of the functors $f^\sharp$ and~$f^\#$ above together
with Corollary~\ref{derived-completion-tensor-with-torsion-isom-cor}.
\end{proof}

\begin{lem} \label{into-discrete-ring-co-contra-extension-of-scalars}
 Let $R$ be an adic topological ring, $S$ be a discrete (i.~e.,
$(0)$\+adic) commutative ring, and $f\:R\rarrow S$ be a continuous
ring homomorphism.
 Then \par
\textup{(a)} for every torsion $R$\+module\/ $\cM$, there is
a natural isomorphism of $S$\+modules
$$
 f^\diamond(\cM)\simeq f^!(\cM)=\Hom_R(S,\cM);
$$ \par
\textup{(b)} for every contramodule $R$\+module\/ $\fP$, there is
a natural isomorphism of $S$\+modules
$$
 f^\#(\fP)\simeq f^*(\fP)=S\ot_R\fP;
$$ \par
\textup{(c)} for every quotseparated contramodule $R$\+module\/ $\fP$,
there is a natural isomorphism of $S$\+modules
$$
 f^\sharp(\fP)\simeq f^*(\fP)=S\ot_R\fP.
$$
\end{lem}

\begin{proof}
 Part~(a): more generally, given a continuous homomorphism of
adic topological rings $f\:R\rarrow S$ and a torsion $R$\+module
$\cM$, one has $f^\diamond(\cM)=f^!(\cM)$ whenever $f^!(\cM)$ is
a torsion $S$\+module.
 When the ring $S$ is discrete, all $S$\+modules are torsion.
 
 Part~(b): more generally, given a continuous homomorphism of
adic topological rings $f\:R\rarrow S$ and a contramodule $R$\+module
$\fP$, one has $f^\#(\fP)=f^*(\fP)$ whenever $f^*(\fP)$ is
a contramodule $S$\+module.
 When the ring $S$ is discrete, all $S$\+modules are contramodules.

 Part~(c): more generally, given a continuous homomorphism of
adic topological rings $f\:R\rarrow S$ and a quotseparated contramodule
$R$\+module $\fP$, one has $f^\sharp(\fP)=f^*(\fP)$ whenever $f^*(\fP)$
is a quotseparated contramodule $S$\+module.
 When the ring $S$ is discrete, all $S$\+modules are quotseparated
contramodules.
\end{proof}

\begin{lem} \label{reductions-of-co-contra-extension-of-scalars}
 Let $f\:R\rarrow S$ be a continuous homomorphism of adic topological
rings, and let $I\subset R$ and $J\subset S$ be two ideals of
definition such that $f(I)\subset J$.
 Then \par
\textup{(a)} for every torsion $R$\+module\/ $\cM$, there is a natural
isomorphism of $S/J$\+modules
$$
 {}_J f^\diamond(\cM)\simeq\Hom_{R/I}(S/J,{}_I\cM);
$$ \par
\textup{(b)} for every contramodule $R$\+module\/ $\fP$, there is
a natural isomorphism of $S/J$\+modules
$$
 f^\#(\fP)/Jf^\#(\fP)\simeq S/J\ot_{R/I}(\fP/I\fP);
$$ \par
\textup{(c)} for every quotseparated contramodule $R$\+module\/ $\fP$,
there is a natural isomorphism of $S/J$\+modules
$$
 f^\sharp(\fP)/Jf^\sharp(\fP)\simeq S/J\ot_{R/I}(\fP/I\fP).
$$
\end{lem}

\begin{proof}
 Let us view $R/I$ and $S/J$ as adic topological rings with
the discrete topologies.
 Denote the related continuous homomorphisms of adic topological rings
by $p\:R\rarrow R/I$, \ $p'\:S\rarrow S/J$, and $f'\:R/I\rarrow S/J$.
 Then we have a commutative square diagram of continuous homomorphisms
of adic topological rings $p'f=f'p$.
 Hence the commutative diagrams of the related functors of
coextension and contraextension of scalars $p'{}^\diamond f^\diamond
\simeq f'{}^\diamond p^\diamond$, $p'{}^\# f^\#\simeq f'{}^\# p^\#$,
and $p'{}^\sharp f^\sharp\simeq f'{}^\sharp p^\sharp$.
 It remains to use
Lemma~\ref{into-discrete-ring-co-contra-extension-of-scalars} for
the computation of the functors of co/contraextension of scalars
with respect to~$p$, $p'$, and~$f'$ in order to obtain
the desired natural isomorphisms.

 Alternatively, part~(a) is the adjunction isomorphism
$\Hom_S(S/J,f^\diamond(\cM))\simeq\Hom_R(S/J,\cM)$.
 Parts~(b) and~(c) are the special cases of the two parts of
Lemma~\ref{map-of-adic-topol-rings-tensor-product-adjunction}
for the torsion $S$\+module $\cN=S/J$.
\end{proof}

 The following corollary is a generalization
of~\cite[Lemma~D.4.3]{Pcosh}.

\begin{cor} \label{flat-contramodules-contraextension-of-scalars}
 Let $f\:R\rarrow S$ be a continuous homomorphism of adic topological
rings, and let $I\subset R$ and $J\subset S$ be two ideals of
definition such that $f(I)\subset J$.
 Then, for any flat quotseparated contramodule $R$\+module\/ $\fF$,
there is a natural isomorphism of quotseparated contramodule
$S$\+modules
$$
 f^\sharp(\fF)\simeq\varprojlim\nolimits_{n\ge1}
 \bigl((S/J^n)\ot_{R/I^n}(\fF/I^n\fF)\bigr).
$$
 In particular, the quoseparated contramodule $S$\+module
$f^\sharp(\fF)$ is flat.
 The functor~$f^\sharp$ takes short exact sequences of flat
quotseparated contramodule $R$\+modules to short exact sequences of
flat quotseparated contramodule $S$\+modules.

 Furthermore, the contramodule $S$\+module $f^\#(\fF)$ is flat for
any flat contramodule $S$\+module~$\fF$.
\end{cor}

\begin{proof}
 All the assertions are based on
Lemma~\ref{reductions-of-co-contra-extension-of-scalars}(b\+-c).
 The flatness of $f^\sharp(\fF)$ and $f^\#(\fF)$ can be established,
under the respective assumptions, using the facts that a contramodule
$R$\+module $\fF$ is flat if and only if the $R/I^n$\+module
$\fF/I^n\fF$ is flat for every $n\ge1$, and similarly for
contramodule $S$\+modules (see
Section~\ref{prelim-inj-proj-flat-torsion-contra-subsecn}).
 The formula for $f^\sharp(\fF)$ follows by virtue of
Lemma~\ref{flat-quotseparated-contramodules-separated}.
 The exactness assertion holds in view of
Lemma~\ref{flat-quotseparated-contramodules-well-behaved}(b),
as the functor of projective limit of sequences of surjective maps
of abelian groups is exact.
 Alternatively, one can also apply the argument
from~\cite[proof of Lemma~D.4.3]{Pcosh} using
Lemma~\ref{flat-quotseparated-contramodules-well-behaved}.
 The result of~\cite[Theorem~1.2 or~2.8]{Yek2}
or~\cite[Lemma~E.1.3]{Pcosh} is relevant here.
\end{proof}

\begin{cor} \label{flat-contramodules-adjusted-to-contraextension}
 Let $f\:R\rarrow S$ be a continuous homomorphism of adic topological
rings, and let\/ $0\rarrow\fQ\rarrow\fP\rarrow\fF\rarrow0$ be
a short exact sequence of quotseparated contramodule $R$\+modules with
a flat quotseparated contramodule $R$\+module\/~$\fF$.
 Then the short sequence of quotseparated contramodule $S$\+modules\/
$0\rarrow f^\sharp(\fQ)\rarrow f^\sharp(\fP)\rarrow f^\sharp(\fF)
\rarrow0$ is exact.
\end{cor}

\begin{proof}
 This is a special case of the second assertion
of~\cite[Lemma~6.2]{Pcs}.
 The claim in question follows purely formally from the exactness
assertion of
Corollary~\ref{flat-contramodules-contraextension-of-scalars}.
 One needs to use the facts that the functor~$f^\sharp$ is right exact
(as a left adjoint functor), there are enough projective quotseparated
contramodule $R$\+modules in $R\Modl_\ctra^\qs$, the projective
quotseparated contramodule $R$\+modules are flat (see
Section~\ref{prelim-inj-proj-flat-torsion-contra-subsecn}), and
the kernels of surjective morphisms of flat quotseparated contramodule
$R$\+modules are flat (see
Lemma~\ref{flat-quotseparated-contramodules-well-behaved}(a)).
\end{proof}

 Now we pass from the continuous maps of adic topological rings to
the tight ones.
 Part~(b) of the following lemma is a generalization
of~\cite[Lemma~D.4.1(b)]{Pcosh}.

\begin{lem} \label{tight-reflects-torsion-contra-lemma}
 Let $R$ and $S$ be adic topological rings, and let $f\:R\rarrow S$
be a tight ring map.
 In this context: \par
\textup{(a)} if\/ $\cN$ is an $S$\+module such that the $R$\+module
$f_*(\cN)$ is torsion, then the $S$\+module\/ $\cN$ is torsion; \par
\textup{(b)} if\/ $\fQ$ is an $S$\+module such that the $R$\+module
$f_*(\fQ)$ is a contramodule, then the $S$\+module\/ $\fQ$ is
a contramodule.
\end{lem}

\begin{proof}
 Recall that, by Lemma~\ref{tight-ring-map-lemma}(1), there exist
an ideal of definition $I\subset R$ and an ideal of definition
$J\subset S$ such that $J\subset Sf(I)$.
 Put $J'=Sf(I)$; so we have $J\subset J'$.

 Part~(a): as mentioned in Section~\ref{prelim-torsion-modules-subsecn},
an $S$\+module $\cN$ is $J'$\+torsion whenever it is $s$\+torsion for
all~$s$ ranging over some chosen set of generators of the ideal
$J'\subset S$.
 The ideal $J'$ is generated by its subset $f(I)\subset J'$.
 For any element $r\in I$, the underlying $R$\+module of $\cN$ is
$r$\+torsion by assumption.
 According to the proof of
Lemma~\ref{torsion-contra-restriction-of-scalars}(a), it follows
that the $S$\+module $\cN$ is $f(r)$\+torsion.
 Thus $\cN$ is $J'$\+torsion, and therefore also $J$\+torsion.
 
 Part~(b): as mentioned in Section~\ref{prelim-contramodules-subsecn},
an $S$\+module $\fQ$ is a $J'$\+contramodule whenever it is
an $s$\+contramodule for all~$s$ ranging over some chosen set of
generators of the ideal $J'\subset S$.
 The ideal $J'$ is generated by its subset $f(I)\subset J'$.
 For any element $r\in I$, the underlying $R$\+module of $\fQ$ is
an $r$\+contramodule by assumption.
 According to the proof of
Lemma~\ref{torsion-contra-restriction-of-scalars}(b), it follows
that the $S$\+module $\fQ$ is an $f(r)$\+contramodule.
 Thus $\fQ$ is a $J'$\+contramodule, and therefore also
a $J$\+contramodule.
\end{proof}

 A quotseparated version of
Lemma~\ref{tight-reflects-torsion-contra-lemma}(b) holds under
an additional flatness assumption on the morphism~$f$.
 See Theorem~\ref{flat-reflects-quotseparated-contramodules-thm}
and Remark~\ref{nonflat-not-reflects-quotseparated} below.

 Part~(b) of the following lemma is a generalization
of~\cite[Lemma~D.4.2(a)]{Pcosh}.

\begin{lem} \label{tight-co-extension-of-scalars}
 Let $R$ and $S$ be adic topological rings, and let $f\:R\rarrow S$
be a tight ring map.
 Then \par
\textup{(a)} for any torsion $R$\+module\/ $\cM$, the $S$\+module
$f^*(\cM)=S\ot_R\cM$ is also torsion; \par
\textup{(b)} for any contramodule $R$\+module\/ $\fP$,
the $S$\+module $f^!(\fP)=\Hom_R(S,\fP)$ is also a contramodule.
\end{lem}

\begin{proof}
 Part~(a): by Lemma~\ref{tensor-product-is-torsion}, the $R$\+module
$S\ot_R\cM$ is torsion.
 By Lemma~\ref{tight-reflects-torsion-contra-lemma}(a), it follows
that the $S$\+module $S\ot_R\cM$ is torsion.  {\hfuzz=3pt\par}

 Part~(b): by Lemma~\ref{Hom-is-a-contramodule}(b), the $R$\+module
$\Hom_R(S,\fP)$ is a contramodule.
 By Lemma~\ref{tight-reflects-torsion-contra-lemma}(b), it follows
that the $S$\+module $\Hom_R(S,\fP)$ is a contramodule.
\end{proof}

 Let $f\:R\rarrow S$ be a tight map of adic topological rings.
 Then Lemma~\ref{tight-co-extension-of-scalars}(a) claims that
the functor of extension of scalars $f^*\:R\Modl\rarrow S\Modl$
restricts to a functor
$$
 f^*\:R\Modl_\tors\lrarrow S\Modl_\tors,
$$
which we will also call the \emph{extension of scalars}.
 When the map~$f$ is both continuous and tight, the functor~$f^*$
is left adjoint to the functor of restriction of scalars
$f_\diamond\:S\Modl_\tors\rarrow R\Modl_\tors$.

 Dual-analogously, for any tight map of adic topological rings
$f\:R\rarrow S$, Lemma~\ref{tight-co-extension-of-scalars}(b) claims
that the functor of coextension of scalars $f^!\:R\Modl\rarrow S\Modl$
restricts to a functor
$$
 f^!\:R\Modl_\ctra\lrarrow S\Modl_\ctra,
$$
which we will also call the \emph{coextension of scalars}.
 When the map~$f$ is both continuous and tight, the functor~$f^!$
is right adjoint to the functor of restriction of scalars
$f_\#\:S\Modl_\ctra\rarrow R\Modl_\ctra$.

 For a quotseparated version of
Lemma~\ref{tight-co-extension-of-scalars}(b), which we can only prove
under a flatness assumption on the map~$f$, see
Proposition~\ref{flat-quotseparated-coextension-of-scalars-prop} below.

\begin{rem}
 For a continuous map of adic topological rings $f\:R\rarrow S$ that
is \emph{not} tight, the functor of restriction of scalars
$f_\diamond\:S\Modl_\tors\rarrow R\Modl_\tors$ does \emph{not} have
a left adjoint functor, generally speaking.
 It suffices to consider the case of the ring of polynomials in one
variable $R=S=k[x]$ over a field~$k$, with the identity map
$f=\id_{k[x]}\:R\rarrow S$, the $(0)$\+adic (i.~e., discrete)
topology on $R$, and the $(x)$\+adic topology on~$S$.
 Then the forgetful functor $f_\diamond\:S\Modl_\tors\rarrow
R\Modl_\tors$ is fully faithful, but it does \emph{not} preserve
infinite direct products (because the full subcategory of
$x$\+torsion modules is not closed under infinite products in
$k[x]\Modl$).
 Accordingly, the functor~$f_\diamond$ is \emph{not} a right adjoint.

 Similarly, for a continuous map of adic topological rings
$f\:R\rarrow S$ that is \emph{not} tight, the functors of restriction
of scalars $f_\#\:S\Modl_\ctra\rarrow R\Modl_\ctra$ and
$f_\sharp\:S\Modl_\ctra^\qs\rarrow R\Modl_\ctra^\qs$ do \emph{not} have
right adjoint functors, generally speaking.
 In the same example of adic topological rings $R$, $S$ and the map~$f$
as above, the forgetful functor $f_\#=f_\sharp\:S\Modl_\ctra\rarrow
R\Modl_\ctra$ is fully faithful, but it does \emph{not} preserve
infinite direct sums (because the full subcategory of $x$\+contramodule
modules is not closed under infinite direct sums in $k[x]\Modl$).
 Accordingly, the functor $f_\#=f_\sharp$ is \emph{not} a left
adjoint.
 (See the discussion in~\cite[Section~2.9]{Pproperf} and
the details in~\cite[Example~10.1]{Pcs}.)
\end{rem}

 The definition of a flat map of adic topological rings was
given at the end of Section~\ref{adic-rings-subsecn}.

\begin{lem} \label{flat-map-of-adic-rings-direct-image-lemma}
 Let $f\:R\rarrow S$ be a flat tight continuous map of adic
topological rings.
 Then \par
\textup{(a)} the functor $f^*\:R\Modl_\tors\rarrow S\Modl_\tors$
is exact; \par
\textup{(b)} the functor $f_\diamond\:S\Modl_\tors\rarrow
R\Modl_\tors$ takes injective torsion $S$\+modules to injective
torsion $R$\+modules; \par
\textup{(c)} the functor $f_\#\:S\Modl_\ctra\rarrow R\Modl_\ctra$
takes flat contramodule $S$\+modules to flat contramodule
$R$\+modules; \par
\textup{(d)} the functor $f_\sharp\:S\Modl_\ctra^\qs\rarrow
R\Modl_\ctra^\qs$ takes flat quotseparated contramodule $S$\+modules
to flat quotseparated contramodule $R$\+modules.
\end{lem}

\begin{proof}
 The assertions of parts~(a) and~(b) are equivalent restatements of
each other, because there are enough injective objects in
$S\Modl_\tors$.
 Nevertheless, we prove them separately below.
 Let $I\subset R$ be an ideal of definition in~$R$;
then, by Lemma~\ref{continuous-and-tight-ring-map-lemma}(5),
\,$Sf(I)\subset S$ is an ideal of definition in~$S$.

 Part~(a): every short exact sequence of $R$\+modules is a direct
limit of short exact sequences of $R/I^n$\+modules, $n\ge1$.
 For any $R/I^n$\+module $E$, one has $f^*(E)=S\ot_RE\simeq
S/Sf(I^n)\ot_{R/I^n}E$.
 Since $S/Sf(I^n)$ is a flat module over $R/I^n$, the functor~$f^*$
takes short exact sequences of $R/I^n$\+modules to short exact
sequences of $S/Sf(I^n)$\+modules.

 Part~(b): for any $S$\+module $K$, we have
$f_*({}_{Sf(I)}K)={}_If_*(K)$.
 For an injective torsion $S$\+module $\cK$, the $S/Sf(I)$\+module
${}_{Sf(I)}\cK$ is injective.
 Since $S/Sf(I)$ is a flat $R/I$\+module, it follows that
${}_If_\diamond(\cK)$ is injective as an $R/I$\+module, as desired.

 Part~(c): for any $S$\+module $P$, we have
$f_*(P/Sf(I)P)=f_*(P)/If_*(P)$.
 For a flat contramodule $S$\+module $\fF$, the $S/Sf(I)$\+module
$\fF/Sf(I)\fF$ is flat.
 Since $S/Sf(I)$ is a flat $R/I$\+module, it follows that
$f_\#(\fF)/If_\#(\fF)$ is flat as an $R/I$\+module, as desired.

 Part~(d) is the particular case of part~(c) for quotseparated
contramodules.
\end{proof}

 Let $R$ be an adic topological ring.
 Then the completion map $\lambda_R\:R\rarrow\fR$ is a flat
tight continuous map of adic topological rings
(as mentioned in Section~\ref{prelim-adic-completions-subsecn}).

 Clearly, any torsion $R$\+module admits a unique extension of
its $R$\+module structure to an $\fR$\+module structure.
 The functor $\lambda_R\hth{}_\diamond\:\fR\Modl_\tors\rarrow
R\Modl_\tors$ is an equivalence (in fact, isomorphism)
of abelian categories.
 Consequently, the functors $\lambda_R^*\:R\Modl_\tors\rarrow
\fR\Modl_\tors$ and $\lambda_R^\diamond\:R\Modl_\tors\rarrow
\fR\Modl_\tors$ adjoint to~$\lambda_R\hth{}_\diamond$ on the left and
on the right are also equivalences of categories, isomorphic to each
other and inverse to~$\lambda_R\hth{}_\diamond$.

 Furthermore, it is clear from the discussion of the category of
contramodules over the topological ring $\fR$ in
Section~\ref{prelim-contramodules-subsecn} that the functor
$\lambda_R\hth{}_\sharp\:\fR\Modl_\ctra^\qs\rarrow R\Modl_\ctra^\qs$
is an equivalence (in fact, isomorphism) of abelian categories.
 Consequently, the functors $\lambda_R^\sharp\:R\Modl_\ctra^\qs
\rarrow\fR\Modl_\ctra^\qs$ and $\lambda_R^!\:R\Modl_\ctra^\qs
\rarrow\fR\Modl_\ctra^\qs$ adjoint to~$\lambda_R\hth{}_\sharp$ on
the left and on the right are also equivalences of categories,
isomorphic to each other and inverse to~$\lambda_R\hth{}_\sharp$.

 Notice that the functor $\lambda_R\hth{}_\#\:\fR\Modl_\ctra\rarrow
R\Modl_\ctra$ is always fully faithful, but it \emph{need not} be
an equivalence of categories in general.
 See the discussion in the proof of
Lemma~\ref{flat-contramodules-characterizations}.

\subsection{Tensor products of adic topological rings}
\label{prelim-tensor-products-of-adic-topological-subsecn}
 Let $R$ be a commutative ring, let $S$ and $T$ be two adic topological
rings, and let $f\:R\rarrow S$ and $g\:R\rarrow T$ be two ring
homomorphisms.
 Consider the ring $W=T\ot_RS$, and denote by $f'\:T\rarrow T\ot_RS$
and $g'\:S\rarrow T\ot_RS$ the two induced ring homomorphisms.

\begin{lem} \label{tensor-product-topology-lemma}
\textup{(a)} There exists a unique finest ring topology on the ring $W$
with a base of neighborhoods of zero formed by open ideals such that
the ring maps $f'$ and~$g'$ are continuous. \par
\textup{(b)} There exists a unique coarsest ring topology on the ring
$W$ with a base of neighborhoods of zero formed by open ideals such
that, for every topological commutative ring $V$ with a base of
neighborhoods of zero formed by open ideals, and for every ring
homomorphism $h\:W\rarrow V$ for which the compositions $hf'\:T\rarrow
V$ and $hg'\:S\rarrow V$ are continuous, the ring homomorphism~$h$
is continuous as well. \par
\textup{(c)} The ring topologies on $W$ defined in parts~\textup{(a)}
and~\textup{(b)} coincide, and it is an adic topology. \par
\textup{(d)} Let $J\subset S$ be an ideal of definition in $S$ and
$K\subset T$ be an ideal of definition in~$T$.
 Then $L=W(f'(K)+g'(J))\subset W$ is an ideal of definition in
the adic topology on $W$ defined in parts~\textup{(a\+-c)}.
\end{lem}

\begin{proof}
 Consider the adic topology on $W$ with the ideal of definition
constructed in part~(d).
 Let us check that this topology on $W$ satisfies the conditions stated
in parts~(a) and~(b).

 Part~(a): the ring maps $f'$ and~$g'$ are continuous by
Lemma~\ref{continuous-ring-map-lemma}(2).
 Conversely, let $\tau$~be a ring topology on $W$, with a base of
neighborhoods of zero formed by open right ideals, such that
the maps $f'$ and $g'$ are continuous with respect to~$\tau$.
 To show that $\tau$~is equal to or coarser than the adic topology on
$W$ defined above, we need to check that all the ideals in $W$ that are
open in~$\tau$ are also open in our adic topology.
 In other words, we need to check that for every ideal $E\subset W$
that is open in~$\tau$ there exists an integer $n\ge1$ such
$L^n\subset E$.

 Indeed, since the map $f'\:T\rarrow W$ is continuous with respect
to~$\tau$, there exists an integer $n_1\ge1$ such that
$f'(K^{n_1})\subset E$.
 Similarly, since the map $g'\:S\rarrow W$ is continuous with respect
to~$\tau$, there exists an integer $n_2\ge1$ such that
$g'(J^{n_2})\subset E$.
 Put $n=n_1+n_2-1$.
 Then we have $(f'(K)+g'(J))^n\subset E$, hence $L^n\subset E$.

 Part~(b): firstly we assume that the compositions $hf'$ and~$hg'$
are continuous, and prove that the map~$h$ is continuous with respect
to our adic topology on~$W$.
 Let $E\subset V$ be an open ideal.
 We need to show that there exists an integer $n\ge1$ such that
$h(L^n)\subset E$.
 Indeed, since the map~$hf'$ is continuous, there exists an integer
$n_1\ge1$ such that $hf'(K^{n_1})\subset E$.
 Similarly, since the map~$hg'$ is continuous, there exists an integer
$n_2\ge1$ such that $hg'(J^{n_2})\subset E$.
 Put $n=n_1+n_2-1$.
 Then we have $(hf'(K)+hg'(J))^n\subset E$, hence $h(L^n)\subset E$.

 Now let $\tau$~be a ring topology on $W$ such that, for every
topological commutative ring $V$ with a base of neighborhoods of zero
formed by open ideals, and for every ring homomorphism $h\:W\rarrow V$
for which the compositions $hf'\:T\rarrow V$ and $hg'\:S\rarrow V$ are
continuous, the map~$h$ is continuous with respect to~$\tau$.
 To show that $\tau$~is equal to or finer than the adic topology on
$W$ defined above, we need to check that all the ideals in $W$ that
are open in our adic topology are also open in~$\tau$.

 Put $V=W$, and endow the ring $V$ with our adic topology on~$W$.
 Let $h\:W\rarrow V$ be the identity map.
 Then the maps $hf'$ and~$hg'$ are continuous, as we have shown in
the very beginning of the proof of part~(a) above.
 So the map~$h$ is continuous with respect to the topology~$\tau$
on $W$ and our adic topology on~$V$.
 This means precisely that $\tau$~is equal to or finer than
our adic topology on~$W$.
\end{proof}

 The adic topology on the ring $W=T\ot_RS$ constructed in
Lemma~\ref{tensor-product-topology-lemma} will be called
the \emph{tensor product topology}.

 In the sequel, we will sometimes abuse the notation and write
$L=K\ot_RS+T\ot_RJ\subset T\ot_RS=W$ instead of $L=W(f'(K)+g'(J))
\subset W$ in the context of
Lemma~\ref{tensor-product-topology-lemma}(d).
 Notice that ideals of the form $L=K\ot_RS+T\ot_RJ$, where $J$ ranges
over the ideals of definition in $S$ and $K$ ranges over the ideals of
definition in $T$, form a base of neighborhoods of zero in the tensor
product topology on $T\ot_RS$.

\begin{lem} \label{topol-tensor-product-monoidal-structure-axioms}
\textup{(a)} Let $f\:R\rarrow S$ be a continuous homomorphism of adic
topological rings.
 Endow the tensor product $R\ot_RS$ with the tensor product topology.
 Then the natural ring isomorphisms $S\rarrow R\ot_RS\rarrow S$ are
isomorphisms of topological rings. \par
\textup{(b)} Let $R$ and $S$ be commutative rings, $T$, $U$, and $V$
be adic topological rings, and $R\rarrow T$, \,$R\rarrow U$,
\,$S\rarrow U$, and $S\rarrow V$ be ring homomorphisms.
 Endow all the relevant tensor products with the tensor product 
topologies.
 Then the natural ring isomorphism $(T\ot_RU)\ot_SV\simeq
T\ot_R(U\ot_SV)$ is an isomorphism of topological rings.
\end{lem}

\begin{proof}
 Part~(a) Let $J\subset S$ be an ideal of definition.
 By Lemma~\ref{continuous-ring-map-lemma}(3), there exists an ideal
of definition $I\subset R$ such that $f(I)\subset J$.
 Now the isomorphism $R\ot_RS\simeq S$ transforms the ideal
$I\ot_RS+R\ot_RJ=R\ot_RJ\subset R\ot_RS$ into the ideal $J\subset S$.
 By Lemma~\ref{tensor-product-topology-lemma}(d), \,$I\ot_RS+R\ot_RJ$
is an ideal of definition in $R\ot_RS$.
 Thus $R\ot_RS\simeq S$ is an isomorphism of topological rings.
 Part~(b) is obvious.
\end{proof}

\begin{lem} \label{tight-flat-adic-ring-base-change-context}
 In the notation above, endow the ring $W=T\ot_RS$ with the tensor
product topology.
 Assume that the ring $R$ is endowed with an adic topology such that
$f\:R\rarrow S$ is a tight continuous ring map and $g\:R\rarrow T$ is
a continuous ring map.
 Then \par
\textup{(a)} the continuous ring map $f'\:T\rarrow W$ is also tight;
\par
\textup{(b)} if the tight continuous ring map~$f$ is flat, then
the tight continuous ring map~$f'$ is flat, too; \par
\textup{(c)} if the ring map\/ $\Lambda(f)\:\Lambda(R)\rarrow\Lambda(S)$
is a (topological) isomorphism, then the ring map\/ $\Lambda(f')\:
\Lambda(T)\rarrow\Lambda(W)$ is an isomorphism as well.
\end{lem}

\begin{proof}
 Part~(a): by Lemma~\ref{continuous-ring-map-lemma}(2), there exist
an ideal of definition $I\subset R$ and an ideal of definition
$K\subset T$ such that $g(I)\subset K$.
 By Lemma~\ref{continuous-and-tight-ring-map-lemma}(5), the ideal
$J=Sf(I)\subset S$ is an ideal of definition in~$S$.
 By Lemma~\ref{tensor-product-topology-lemma}(4), the ideal
$L=W(f'(K)+g'(J))\subset W$ is an ideal of definition in~$W$.
 Now we have $W(f'(K)+g'(J))=W(f'(K)+g'f(I))=W(f'(K)+f'g(I))=Wf'(K)$.
 By Lemma~\ref{tight-ring-map-lemma}(1)
or~\ref{continuous-and-tight-ring-map-lemma}(4), it follows that
the map~$f'$ is tight.

 Part~(b): now we apply Lemma~\ref{continuous-ring-map-lemma}(3)
to the effect that, for every ideal of definition $K\subset T$,
there exists an ideal of definition $I\subset R$ such that
$f(I)\subset K$.
 Then, in the notation of the previous paragraph, we have
$W/L\simeq T/K\ot_{R/I}S/J$.
 The ring $S/J$ is flat as an $R/I$\+module by
Lemma~\ref{flat-homomorphism-of-adic-rings}(2), hence
the ring $W/L$ is flat as a $T/K$\+module.

 Part~(c): first of all, for a tight continuous map~$f$, the map
$\Lambda(f)$ is also tight (and continuous) by
Lemma~\ref{induced-map-of-complete-rings-lemma}(a).
 Hence $\Lambda(f)$ is an isomorphism of abstract rings if and only
if $\Lambda(f)$ is an isomorphism of topological rings.

 Now let $K\subset T$ be an ideal of definition.
 Once again, in the notation above, we have
$W/L\simeq T/K\ot_{R/I}S/J$.
 If $\Lambda(f)$~is an isomorphism, then so is the map $R/I\rarrow
S/J$.
 Hence the map $T/K\rarrow W/L$ is an isomorphism, too.
 As this holds for all ideals of definition $K\subset T$, we can
pass to the projective limit and conclude that the map
$\Lambda(f')$ is an isomorphism.
\end{proof}

 Now consider the following commutative diagram of homomorphisms of
commutative rings:
\begin{equation} \label{tensor-product-functoriality-datum-diagram}
\begin{gathered}
 \xymatrix{
  \widetilde S & & \widetilde T \\
  S \ar[u]^{h_S} & \widetilde R
  \ar[ul]^{\tilde f} \ar[ur]_{\tilde g} & T \ar[u]_{h_T} \\
  & R \ar[ul]^f \ar[u]^-{h_R} \ar[ur]_g
 }
\end{gathered}
\end{equation}
 Assume that adic topologies are given on the rings $S$, $\widetilde S$,
$T$, and $\widetilde T$ such that $h_S$ and~$h_T$ are continuous ring
maps.
 Consider the tensor product rings $T\ot_RS$ and $\widetilde T
\ot_{\widetilde R}\widetilde S$, and endow them with the tensor product
topologies.
 Then the homomorphism of tensor products $h\:T\ot_RS\rarrow
\widetilde T\ot_{\widetilde R}\widetilde S$ induced by the ring
homomorphisms $h_S$, $h_T$, and~$h_R$ is a continuous ring map.

 The following lemma is a generalization of
Lemma~\ref{tight-flat-adic-ring-base-change-context}(c).

\begin{lem} \label{completed-tensor-product-of-completions-lemma}
 Let $S$, $T$, and $R$ be adic topological rings, and let
$f\:R\rarrow S$ and $g\:R\rarrow T$ be continuous ring maps.
 Consider the induced continuous maps of the completed rings\/
$\Lambda(f)\:\Lambda(R)\rarrow\Lambda(S)$ and\/
$\Lambda(g)\:\Lambda(R)\rarrow\Lambda(T)$.
 Furthermore, consider the induced continuous homomorphism of
the tensor product rings
$$
 h\:T\ot_RS\lrarrow\Lambda(T)\ot_{\Lambda(R)}\Lambda(S).
$$
 Then the induced homomorphism of the completions of the tensor
products
$$
 \Lambda(h)\:\Lambda(T\ot_RS)\lrarrow
 \Lambda(\Lambda(T)\ot_{\Lambda(R)}\Lambda(S)).
$$
is an isomorphism of topological rings.
\end{lem}

\begin{proof}
 Put $W=T\ot_RS$.
 Let $J\subset S$ be an ideal of definition and $K\subset T$ be
an ideal of definition.
 Since the product of any two ideals of definition in $R$ is
also an ideal of definition in $R$,
\,Lemma~\ref{continuous-ring-map-lemma}(3) implies existence of
an ideal of definition $I\subset R$ such that $f(I)\subset J$
and $g(I)\subset K$.
 Then, according to Lemma~\ref{tensor-product-topology-lemma}(d),
\,$L=W(f'(K)+g'(J))\subset W$ is an ideal of definition in~$W$.
 We have $W/L\simeq T/K\ot_{R/I}S/J$.

 Put $\fR=\Lambda(R)$, \,$\fS=\Lambda(S)$, \,$\fT=\Lambda(T)$,
and $\widetilde W=\fT\ot_\fR\fS$ (the latter notation is intended to
emphasize the fact that the topological ring $\fT\ot_\fR\fS$
\emph{need not} be separated or complete).
 Consider the ideals $\fI=\fR\lambda_R(I)\subset\fR$,
\ $\fJ=\fS\lambda_S(J)\subset\fS$, and $\fK=\fT\lambda_S(K)\subset\fT$
corresponding to the ideals $I$, $J$, and $K$ under the bijection
described in Section~\ref{prelim-adic-completions-subsecn}.
 Put $\ff=\Lambda(f)$ and $\fg=\Lambda(g)$, and denote by
$\ff'\:\fT\rarrow\widetilde W$ and $\fg'\:\fS\rarrow\widetilde W$
the induced homomorphisms of topological rings.
 Then $\widetilde L=\widetilde W(\ff'(\fK)+\fg'(\fJ))\subset
\widetilde W$ is an ideal of definition in $\widetilde W$, and
we have $\widetilde W/\widetilde L\simeq \fT/\fK\ot_{\fR/\fI}\fS/\fJ$.

 Now the ring homomorphisms $R/I\rarrow\fR/\fI$, \
$S/J\rarrow\fS/\fJ$, and $T/K\rarrow\fT/\fK$ are isomorphisms.
 Hence so is the ring homomorphism $W/L\rarrow
\widetilde W/\widetilde L$.
 It remains to point out that one has $\Lambda(W)=
\varprojlim_{J,K}W/L$ and $\Lambda(\widetilde W)=
\varprojlim_{J,K}\widetilde W/\widetilde L$, where $J$ and $K$
range over the directed posets of all ideals of definition in $S$
and~$T$.
 The point is that the ideals of definition of the form
$L=W(f'(K)+g'(J))$ form a base of neighborhoods of zero in $W$, and
similarly for the ideals $\widetilde L\subset\widetilde W$.
 Furthermore, the adic topology on $\Lambda(W)$ is the topology of
projective limt of discrete rings $W/L$, while the topology on
$\Lambda(\widetilde W)$ is the topology of projective limit of
discrete rings $\widetilde W/\widetilde L$.
 Thus the map $\Lambda(h)\:\Lambda(W)\rarrow\Lambda(\widetilde W)$
is an isomorphism of topological rings.
\end{proof}

\begin{cor} \label{tensor-product-isomorphism-after-completion-cor}
 In the context of the commutative diagram of commutative ring
homomorphisms~\eqref{tensor-product-functoriality-datum-diagram},
suppose that all the six rings are endowed with adic topologies
and all the seven maps are continuous ring homomorphisms.
 Assume further that the induced maps of the completions\/
$\Lambda(h_S)\:\Lambda(S)\rarrow\Lambda(\widetilde S)$,
\ $\Lambda(h_T)\:\Lambda(T)\rarrow\Lambda(\widetilde T)$,
and\/ $\Lambda(h_R)\:\Lambda(R)\rarrow\Lambda(\widetilde R)$ are
isomorphisms of topological rings.
 Then the induced map of the completions of the tensor products
$$
 \Lambda(h)\:\Lambda(T\ot_RS)\lrarrow
 \Lambda(\widetilde T\ot_{\widetilde R}\widetilde S)
$$
is also an isomorphism of topological rings.
\end{cor}

\begin{proof}
 Follows from
Lemma~\ref{completed-tensor-product-of-completions-lemma}.
\end{proof}

\subsection{Formal open immersions and formal open coverings}
\label{prelim-formal-open-immersions-coverings-subsecn}
 Let us say that a homomorphism of (discrete) commutative rings
$f\:R\rarrow S$ is an \emph{open immersion} if the induced map of
the spectra $\Spec S\rarrow\Spec R$ is an open immersion of affine
schemes.
 Equivalently, $f$~is an open immersion if and only if it is
a flat ring epimorphism of finite
presentation~\cite[Section~XI.2]{Sten},
\cite[Theoreme~IV.17.9.1]{EGAIV4}.
 In particular, we will say that a homomorphism of commutative
rings $R\rarrow S$ is a \emph{principal open immersion} if
there exists an element $r\in R$ such that $S$ is isomorphic to
$R[r^{-1}]$ as a commutative $R$\+algebra.

 Similarly, let us say that a collection of homomorphisms of
commutative rings $f_\alpha\:R\rarrow S_\alpha$ is a (\emph{principal})
\emph{open covering} if every map~$f_\alpha$ is a (principal) open
immersion and the collection of the induced maps of the spectra
$\Spec S_\alpha\rarrow\Spec R$ is an affine open covering of
the affine scheme $\Spec R$.
 Equivalently, the latter condition means that, for any field~$k$ and
any ring homomorphism $R\rarrow k$, there exists an index~$\alpha$
such that $S_\alpha\ot_Rk\ne0$.
 The following criterion is a basic result concerning affine open
subschemes in affine schemes.

\begin{lem} \label{open-immersions-in-terms-of-principal}
 A homomorphism of commutative rings $f\:R\rarrow S$ is an open
immersion if and only if there exists a finite collection of elements
$r_1$,~\dots, $r_m\in R$ such that the induced ring maps $R[r_j^{-1}]
\rarrow S[f(r_j)^{-1}]$ are isomorphisms for all\/ $1\le j\le m$
and the collection of principal open immersions $S\rarrow
S[f(r_j)^{-1}]$, \,$1\le j\le m$, is an open covering. \qed
\end{lem}

\begin{lem} \label{open-immersions-via-reductions}
 Let $f\:R\rarrow S$ be a homomorphism of commutative rings such that
$S$ is a flat $R$\+module, and let $I\subset R$ be a nilpotent ideal.
 Then \par
\textup{(a)} $f\:R\rarrow S$ is an isomorphism if and only if
$\bar f\:R/I\rarrow S/Sf(I)$ is an isomorphism; \par
\textup{(b)} $f\:R\rarrow S$ is a principal open immersion if and only
if $\bar f\:R/I\rarrow S/Sf(I)$ is a principal open immersion; \par
\textup{(c)} $f\:R\rarrow S$ is an open immersion if and only if
$\bar f\:R/I\rarrow S/Sf(I)$ is an open immersion.
\end{lem}

\begin{proof}
 Part~(a): more generally, for any nilpotent ideal $I$ in
a commutative ring $R$, a homomorphism of flat $R$\+modules
$f\:F\rarrow G$ is an isomorphism if and only if the induced map
$\bar f\:F/IF\rarrow G/IG$ is an isomorphism.
 This is provable using the Nakayama lemma for nilpotent ideals.

 Part~(b): the ``only if'' assertion is obvious.
 To prove the ``if'', let $\bar r\in R/I$ be an element such that
the $R/I$\+algebra $S/Sf(I)$ is isomorphic to $(R/I)[\bar r^{-1}]$.
 Let $r\in R$ be any element such that $\bar r = r+I$.
 Applying part~(a) to the $R$\+algebra morphisms $R[r^{-1}]
\rarrow S[f(r)^{-1}]\larrow S$ and the nilpotent ideals
$I[r^{-1}]\subset R[r^{-1}]$ and $SI\subset S$, one concludes
that both of these $R$\+algebra morphisms are isomorphisms.

 Part~(c): once again, the ``only if'' assertion is obvious.
 To prove the ``if'', let $\bar r_1$,~\dots, $\bar r_m\in R/I$ be
a finite collection of elements satisfying the conditions of
Lemma~\ref{open-immersions-in-terms-of-principal} for
the ring homomorphism $R/I\rarrow S/Sf(I)$.
 Let $r_j\in R$, \,$1\le j\le m$, be any element such that
$\bar r_j = r_j+I$.
 Applying part~(a) to the ring homomorphisms $R[r_j^{-1}]\rarrow
S[f(r_j)^{-1}]$ and the nilpotent ideals $I[r_j^{-1}]\subset
R[r_j^{-1}]$, one concludes that these ring homomorphisms
are isomorphisms.
 To prove that the collection of principal open immersions
$S\rarrow S[f(r_j)^{-1}]$, \,$1\le j\le m$, is an open covering,
it suffices to observe that any homomorphism from $S$ to a field~$k$
factorizes through $S/Sf(I)$.
\end{proof}

\begin{lem} \label{open-coverings-via-reductions}
 Let $f_\alpha\:R\rarrow S_\alpha$ be a collection of homomorphisms of
commutative rings such that the rings $S_\alpha$ are flat $R$\+modules,
and let $I\subset R$ be a nilpotent ideal.
 Consider the induced ring homomorphisms $\bar f_\alpha\:R/I\rarrow
S_\alpha/S_\alpha f_\alpha(I)$.
 Then the collection $(f_\alpha)$ is an open covering if and only if
the collection $(\bar f_\alpha)$ is an open covering.
\end{lem}

\begin{proof}
 One needs to use Lemma~\ref{open-immersions-via-reductions}(c).
 On top of that, the final argument from the proof of the same lemma
works:
 Any ring homomorphism from $R$ to a field~$k$ factorizes
through~$R/I$.
\end{proof}

\begin{lem} \label{open-coverings-covered-by-principal}
 A homomorphism of commutative rings $f\:R\rarrow S$ is an open
immersion if and only if there exists an open covering $g_\alpha\:
S\rarrow T_\alpha$ of the ring $S$ such that the compositions
$g_\alpha f\:R\rarrow T_\alpha$ are open immersions.
 If this is the case, then one can choose $(g_\alpha)$ to be
a principal open covering for which all the compositions
$g_\alpha f$ are also principal open immersions.
\end{lem}

\begin{proof}
 The first assertion expresses locality of the notion of an open
immersion of schemes.
 It is immediate from the definitions.
 The second assertion is a restatement of the fact that principal
affine open subschemes form a topology base in $\Spec R$.
\end{proof}

\begin{thm} \label{open-immersions-nilpotent-reduction-theorem}
 Let $R$ be a commutative ring and $I\subset R$ be a nilpotent ideal.
 Then the functor $S\longmapsto S/Sf(I)$ provides an equivalence
between the category of commutative $R$\+algebras $f\:R\rarrow S$
such that~$f$ is an open immersion and the category of commutative
$R/I$\+algebras $\bar f\:R/I\rarrow\overline S$ such that
$\bar f$~is an open immersion.
 In other words, the map\/ $\Spec S\longmapsto\Spec S/Sf(I)$ is
an isomorphism between the poset of affine open subschemes in\/
$\Spec R$ and the poset of affine open subschemes in\/ $\Spec R/I$
(with respect to the inclusion order).
 An open immersion~$f$ is principal if and only if the open
immersion~$\bar f$ is.
\end{thm}

\begin{proof}
 We are not aware of any argument within the realm of commutative
ring theory.
 The following proof is, basically, a collection of references.
 By~\cite[Theorem~IV.18.1.2]{EGAIV4}, the functor $X\longmapsto
\Spec R/I\times_{\Spec R}X$ is an equivalence between the category
of \'etale schemes over $\Spec R$ and the category of \'etale
schemes over $\Spec R/I$.
 By~\cite[Theorem~IV.17.9.1]{EGAIV4}, an \'etale morphism is an open
immersion if and only if it is radicial (in the sense
of~\cite[Section~I.3.7]{EGAI}).
 It follows easily that an \'etale morphism $X\rarrow\Spec R$ is
an open immersion if and only if the corresponding morphism
$\Spec R/I\times_{\Spec R}X\rarrow\Spec R/I$ is an open immersion.
 So the map $\Spec R/I\times_{\Spec R}{-}$ is an equivalence between
the category of open subschemes in $\Spec R$ and the category of
open subschemes in $\Spec R/I$.
 The same conclusion can be simply obtained from the facts that
open subschemes in schemes correspond bijectively to open subsets
in the underlying topology, and the underlying topological spaces
of $\Spec R$ and $\Spec R/I$ coincide.

 The key step is to prove that an open subscheme $Y\subset\Spec R$
is affine if and only if the scheme $\Spec R/I\times_{\Spec R}Y$ is
affine.
 The ``only if'' assertion is obvious.
 To prove the ``if'', notice that the underlying topological spaces
of $Y$ and $\Spec R/I\times_{\Spec R}Y$ coincide.
 Quasi-compactness and quasi-separatedness are properties of
the underlying topological space of a scheme; and all affine schemes
are quasi-compact and quasi-separated.
 So the scheme $Y$ is quasi-compact and quasi-separated (in fact,
separated, as an open subscheme of an affine scheme).
 By~\cite[Th\'eor\`eme~II.5.2.1(d)]{EGAII}, \cite[IV.1.7.17]{EGAIV1},
\cite[Lemmas Tags~01XB and~01XF or~01XG]{SP}, a quasi-compact
scheme $Y$ is affine if and only if $H^1(Y,\cF)=0$ for every
quasi-coherent sheaf $\cF$ on~$Y$.
 Here $H^*(Y,{-})$ denotes the cohomology groups of sheaves of
abelian groups on $Y$, which only depend on the underlying topological
space of~$Y$.
 It remains to point out that every quasi-coherent sheaf $\cF$ on $Y$
has a finite decreasing filtration by quasi-coherent subsheaves
$\cF\supset I\cF\supset I^2\cF\supset\dotsb\supset I^n\cF=0$, where
$n\ge1$ is an integer such that $I^n=0$ in $R$, and the successive
quotient sheaves $I^m\cF/I^{m+1}\cF$ are quasi-coherent sheaves
on $\Spec R/I\times_{\Spec R}Y$.

 Finally, if $f\:R\rarrow S$ is an open immersion and the $R/I$\+algebra
$S/Sf(I)$ is isomorphic to $(R/I)[\bar r^{-1}]$ for some element
$\bar r\in R/I$, then the $R$\+algebra $S$ is isomorphic to $R[r^{-1}]$
for any element $r\in R$ such that $\bar r=r+I$.
 This assertion follows immediately from the bijection estalished above.
\end{proof}

\begin{lem} \label{formal-open-immersions-lemma}
 Let $f\:R\rarrow S$ be a tight continuous map of adic
topological rings.
 Then the following three conditions are equivalent:
\begin{enumerate}
\item for every ideal of definition $I\subset R$, the induced ring map
$\bar f\:R/I\rarrow S/Sf(I)$ is an open immersion;
\item there exists an ideal of definition $I\subset R$ such that
the induced ring map $\bar f\:R/I^n\rarrow S/Sf(I^n)$ is an open
immersion for every $n\ge1$;
\item $f$~is flat (as a map of adic topological rings) and there
exists an ideal of definition $I\subset R$ such that the induced ring
map $\bar f\:R/I\rarrow S/Sf(I)$ is an open immersion.
\end{enumerate}
\end{lem}

\begin{proof}
 (1)~$\Longrightarrow$~(2) Holds because $I^n\subset R$ is an ideal
of definition for any ideal of definition $I\subset R$.
 
 (2)~$\Longrightarrow$~(1) The point is that if $I\subset I'\subset R$
are two ideals and the ring map $R/I\rarrow S/Sf(I)$ is an open
immersion, then the ring map $R/I'\rarrow S/Sf(I')$ is an open
immersion, too.

 (1)~$\Longrightarrow$~(3) Holds because, for all open immersions of
commutative rings $\bar f\:\overline R\rarrow\overline S$, the ring
$\overline S$ is a flat $\overline R$\+module.

 (3)~$\Longrightarrow$~(2) Follows from
Lemma~\ref{open-immersions-via-reductions}(c).
\end{proof}

 We will say that a tight continuous map of adic topological rings
$f\:R\rarrow S$ is a \emph{formal open immersion} if it satisfies
the equivalent conditions of Lemma~\ref{formal-open-immersions-lemma}.
 Using the criterion of
Lemma~\ref{formal-open-immersions-lemma}(1), one can easily check that
the composition of any two formal open immersions is a formal open
immersion.

\begin{lem} \label{formal-open-immersion-as-formal-ring-epimorphism}
 Let $f\:R\rarrow S$ be a continuous map of adic topological rings.
 Consider the ring $S\ot_RS$, and endow it with the tensor product
topology.
 In this context: \par
\textup{(a)} the induced map of topological rings
$S\ot_RS\rarrow S$ is continuous and tight; \par
\textup{(b)} the induced maps of topological rings
$S\rightrightarrows S\ot_RS$ are continuous; \par
\textup{(c)} if the map~$f$ is tight, then so are the maps
$S\rightrightarrows S\ot_RS$; \par
\textup{(d)} if the map~$f$ is a formal open immersion, then all
the three induced maps of the completions
$$
 \Lambda(S)\rightrightarrows\Lambda(S\ot_RS)\rarrow\Lambda(S)
$$
are isomorphisms of topological rings.
\end{lem}

\begin{proof}
 Parts~(a) and~(b) do not depend on the topology on~$R$ (one can use
the discrete topology on $R$, for example).
 Part~(a): denote the map in question by $p\:S\ot_RS\rarrow S$.
 Let $J\subset S$ be an ideal of definition.
 Then, by Lemma~\ref{tensor-product-topology-lemma}(d),
\,$J\ot_RS+S\ot_RJ$ is an ideal of definition in $S\ot_RS$.
 Now we have $J=p(J\ot_RS+S\ot_RJ)$, so $p$~is a tight continuous map
by Lemma~\ref{continuous-and-tight-ring-map-lemma}(4).

 Part~(b) follows from the discussion in
Section~\ref{prelim-tensor-products-of-adic-topological-subsecn}
(including
Lemma~\ref{topol-tensor-product-monoidal-structure-axioms}(a)).

 Part~(c): denote the maps in question by $i_1$ and $i_2\:S
\rightrightarrows S\ot_RS$.
 Let $I\subset R$ be an ideal of definition.
 By Lemma~\ref{continuous-and-tight-ring-map-lemma}(5), the ideal
$J=Sf(I)\subset S$ is an ideal of definition in~$S$.
 By Lemma~\ref{tensor-product-topology-lemma}(d), the ideal
$J\ot_RS+S\ot_RJ\subset S\ot_RS$ is an ideal of definition in $S\ot_RS$.
 Now we have $Sf(I)\ot_RS=S\ot_RSf(I)$, hence
$J\ot_RS+S\ot_RJ=J\ot_RS=S\ot_RJ\subset S\ot_RS$.
 By Lemma~\ref{continuous-and-tight-ring-map-lemma}(4), it follows that
$i_1$ and~$i_2$ are tight ring maps.

 Part~(d): in the notation of the previous paragraph, let $I$ range
over all the ideals of definition in~$R$.
 Then the ideals $J=Sf(I)$ form a base of neighborhoods of zero in $S$,
while the ideals $J\ot_RS=S\ot_RJ$ form a base of neighborhoods of
zero in $S\ot_RS$.
 Now we have $S/(J\ot_RS+S\ot_RJ)\simeq S/J\ot_{R/I}S/J$.
 The natural maps $S/J\rightrightarrows S/J\ot_{R/I}S/J\rarrow S/J$ are
ring isomorphisms, since $R/I\rarrow S/J$ is an open immersion
(of abstract commutative rings) by
Lemma~\ref{formal-open-immersions-lemma}(1).
 Passing to the projective limit over $I$, we obtain the desired
isomorphisms of complete adic topological rings.
\end{proof}

 For a discussion of the change-of-scalar functors with respect to
a formal open immersion of adic topological rings, see
Proposition~\ref{formal-open-immersion-change-of-scalars} below.

\begin{lem} \label{principal-formal-open-immersions-lemma}
 Let $f\:R\rarrow S$ be a tight continuous map of adic
topological rings.
 Then the following four conditions are equivalent:
\begin{enumerate}
\item for every ideal of definition $I\subset R$, the induced ring
map $\bar f\:R/I\rarrow S/Sf(I)$ is a principal open immersion;
\item there exists an ideal of definition $I\subset R$ such that
the induced ring map $\bar f\:R/I^n\rarrow S/Sf(I^n)$ is a principal 
open immersion for every $n\ge1$;
\item $f$~is flat (as a map of adic topological rings) and there
exists an ideal of definition $I\subset R$ such that the induced ring
map $\bar f\:R/I\rarrow S/Sf(I)$ is a principal open immersion;
\item there exists an element $r\in R$ such that, endowing the ring
$R[r^{-1}]$ with the adic topology for which the natural ring map
$l\:R\rarrow R[r^{-1}]$ is continuous and tight, there exists a unique
isomorphism of topological rings\/ $\Lambda(S)\simeq
\Lambda(R[r^{-1}])$ making the following diagram commutative:
\begin{equation} \label{principal-formal-open-immersion-diagram}
\begin{gathered}
 \xymatrix{
  R[r^{-1}] \ar[rr]^-{\lambda_{R[r^{-1}]}}
  && \Lambda(R[r^{-1}]) \ar@{=}[r]
  & \Lambda(S) && S \ar[ll]_-{\lambda_S} \\
  R \ar[rr]^-{\lambda_R} \ar[u]^-l && \Lambda(R) \ar[u]^-{\Lambda(l)}
  \ar@{=}[r]^{\id} \ar[u]
  & \Lambda(R) \ar[u]_-{\Lambda(f)} && R \ar[ll]_-{\lambda_R} \ar[u]_-f
 }
\end{gathered}
\end{equation}
\end{enumerate}
\end{lem}

\begin{proof}
 (1)~$\Longrightarrow$~(2) and (1)~$\Longrightarrow$~(3)
 Similar to the proof of Lemma~\ref{formal-open-immersions-lemma}.

 (2)~$\Longrightarrow$~(1) The point is that if $I\subset I'\subset R$
are two ideals and the ring map $R/I\rarrow S/Sf(I)$ is a principal
open immersion, then the ring map $R/I'\rarrow S/Sf(I')$ is
a principal open immersion, too.

 (3)~$\Longrightarrow$~(2) Follows from
Lemma~\ref{open-immersions-via-reductions}(b).

 (4)~$\Longrightarrow$~(1) Obvious.

 (4)~$\Longrightarrow$~(3) Follows from
Lemma~\ref{induced-map-of-complete-rings-lemma}(b).

 (3)~$\Longrightarrow$~(4) Let $\bar r\in R/I$ be an element such that
the $R/I$\+algebra $S/Sf(I)$ is isomorphic to $(R/I)[\bar r^{-1}]$.
 Pick any element $r\in R$ for which $\bar r=r+I$.
 Let $n\ge1$ be an integer; put $\tilde r=r+I^n\in R/I^n$.
 By the proof of Lemma~\ref{open-immersions-via-reductions}(b),
the $R/I^n$\+algebra $S/Sf(I^n)$ is isomorphic to
$(R/I^n)[\tilde r^{-1}]$.
 As the category of open subschemes in $\Spec R$ is a poset,
an isomorphism of $R/I^n$\+algebras $(R/I^n)[\tilde r^{-1}]\simeq
S/Sf(I^n)$ is unique.
 Hence, as $n$~varies, such isomorphisms form a projective system.
 Passing to the projective limit over $n\ge1$, we obtain
the desired isomorphism of adic topological rings and
$\Lambda(R)$\+algebras $\Lambda(R[r^{-1}])\simeq\Lambda(S)$.
\end{proof}

 We will say that a tight continuous map of adic topological rings
$f\:R\rarrow S$ is a \emph{principal formal open immersion} if
it satisfies the equivalent conditions of
Lemma~\ref{principal-formal-open-immersions-lemma}.
 Using the criterion of
Lemma~\ref{principal-formal-open-immersions-lemma}(1), one can easily
check that the composition of any two principal formal open immersions
is a principal formal open immersion.

\begin{cor} \label{formal-open-immersions-reduction-corollary}
 Let $R$ be an adic topological ring and $I\subset R$ be an ideal of
definition.
 Then the functor $S\longmapsto S/Sf(I)$ provides an equivalence
between the category of commutative $R$\+algebras $f\:R\rarrow S$
with an adic topology on $S$ such that $f$~is a formal open immersion
and the category of commutative $R/I$\+algebras $\bar f\:R/I\rarrow
\overline S$ such that $\bar f$~is an open immersion.
 A formal open immersion~$f$ is principal if and only if the open
immersion~$\bar f$ is principal.
\end{cor}

\begin{proof}
 This is a corollary of
Theorem~\ref{open-immersions-nilpotent-reduction-theorem}.
 Let us just explain how to construct the inverse functor.
 Suppose given an open immersion $\bar f\:R/I\rarrow
\overline S$.
 By Theorem~\ref{open-immersions-nilpotent-reduction-theorem},
for every $n\ge1$ there exists an open immersion $f_n\:R/I^n
\rarrow S_n$ together with a $R/I$\+algebra isomorphism
$R/I\ot_{R/I^n}S_n\simeq\overline S$.
 The uniqueness assertion of
Theorem~\ref{open-immersions-nilpotent-reduction-theorem}
implies existence and uniqueness of $R/I^m$\+algebra isomorphisms
$R/I^m\ot_{R/I^n}S_n\simeq S_m$ for all $n\ge m\ge1$.
 It remains to set $S=\varprojlim_{n\ge1}S_n$.
 The result of~\cite[Theorem~1.2 or~2.8]{Yek2}
or~\cite[Lemma~E.1.3]{Pcosh} is helpful here.
 To prove the final assertion of the corollary, one needs to use
Lemma~\ref{principal-formal-open-immersions-lemma}(1).
\end{proof}

\begin{lem} \label{open-immersion-adic-ring-base-change-context}
 Let $R$, $S$, and $T$ be adic topological rings, let $f\:R\rarrow S$
be a tight continuous ring map, and let $g\:R\rarrow T$ be
a continuous ring map.
 Put $W=T\ot_RS$, and denote by $f'\:T\rarrow W$ and $g'\:S\rarrow W$
the induced ring maps.
 Endow the commutative ring $W$ with the tensor product topology, as
in Section~\ref{prelim-tensor-products-of-adic-topological-subsecn}.
 In this setting: \par
\textup{(a)} if the tight continuous ring map~$f$ is a formal open
immersion, then the tight continuous ring map~$f'$ is a formal open
immersion, too; \par
\textup{(b)} if the tight continuous ring map~$f$ is a principal
formal open immersion, then the tight continuous ring map~$f'$ is
a principal formal open immersion, too.
\end{lem}

\begin{proof}
 Similar to the proof of
Lemma~\ref{tight-flat-adic-ring-base-change-context}(b).
\end{proof}

\begin{lem} \label{formal-open-coverings-lemma}
 Let $f_\alpha\:R\rarrow S_\alpha$ be a collection of tight continuous
maps of adic topological rings.
 Then the following three conditions are equivalent:
\begin{enumerate}
\item for every ideal of definition $I\subset R$, the collection of
induced ring maps $\bar f_\alpha\:R/I\rarrow
S_\alpha/S_\alpha f_\alpha(I)$ is an open covering;
\item there exists an ideal of definition $I\subset R$ such that
the collection of induced ring maps $\bar f_\alpha\:R/I^n\rarrow
S_\alpha/S_\alpha f_\alpha(I^n)$ is an open covering for every $n\ge1$;
\item the maps~$f_\alpha$ are flat (as maps of adic topological rings)
and there exists an ideal of definition $I\subset R$ such that
the collection of induced ring maps $\bar f_\alpha\:R/I\rarrow
S_\alpha/S_\alpha f_\alpha(I)$ is an open covering.
\end{enumerate}
\end{lem}

\begin{proof}
 Follows from Lemmas~\ref{open-coverings-via-reductions}
and~\ref{formal-open-immersions-lemma}.
\end{proof}

 We will say that a collection of tight continuous maps of adic
topological rings $f_\alpha\:R\rarrow S_\alpha$ is a \emph{formal open
covering} if it satisfies the equivalent conditions of
Lemma~\ref{formal-open-coverings-lemma}.
 It is clear from Lemma~\ref{formal-open-coverings-lemma}(3) that
any formal open covering of an adic topological ring has a finite 
subcovering.
 A \emph{principal formal open covering} is a formal open covering
$(f_\alpha)$ such that every map~$f_\alpha$ is a principal formal
open immersion.

\begin{lem} \label{formal-open-coverings-covered-by-principal}
 A tight continuous map of adic topological rings $f\:R\rarrow S$ is
an formal open immersion if and only if there exists a formal open
covering $g_\alpha\:S\rarrow T_\alpha$ of the adic topological ring $S$
such that the compositions $g_\alpha f\:R\rarrow T_\alpha$ are formal
open immersions.
 If this is the case, then one can choose $(g_\alpha)$ to be
a principal formal open covering for which all the compositions
$g_\alpha f$ are also principal formal open immersions.
\end{lem}

\begin{proof}
 First assertion: to prove the ``only if'', it suffices to take
the set of indices $\{\alpha\}$ to be the singleton $\{0\}$,
and the trivial formal open covering with $T_0=S$ and $g_0=\id_S$.
 To prove the ``if'', let $I\subset R$ be an ideal of definition.
 By Lemma~\ref{continuous-and-tight-ring-map-lemma}(5),
the ideal $Sf(I)\subset S$ is an ideal of definition in $S$,
while the ideals $T_\alpha g_\alpha f(I)\subset T_\alpha$ are
ideals of definition in~$T_\alpha$.
 By Lemma~\ref{formal-open-coverings-lemma}(1), the collection of
induced ring maps $\bar g_\alpha\:S/Sf(I)\rarrow
T_\alpha/T_\alpha g_\alpha f(I)$ is an open covering.
 By Lemma~\ref{formal-open-immersions-lemma}(1), the compositions
$\bar g_\alpha\bar f\:R/I\rarrow T_\alpha/T_\alpha g_\alpha f(I)$
are open immersions.
 By the first assertion of
Lemma~\ref{open-coverings-covered-by-principal}, it follows that
the induced ring map $\bar f\:R/I\rarrow S/Sf(I)$ is an open immersion.
 Using Lemma~\ref{formal-open-immersions-lemma}(1) again, we conclude
that $f\:R\rarrow S$ is a formal open immersion.

 Second assertion: let $I\subset R$ be an ideal of definition.
 By Lemma~\ref{formal-open-immersions-lemma}(1), the induced ring map
$\bar f\:R/I\rarrow S/Sf(I)$ is an open immersion.
 By the sectond assertion of
Lemma~\ref{open-coverings-covered-by-principal}, there exists
a principal open covering $\bar g_\alpha\:S/Sf(I)\rarrow\overline
T_\alpha$ such that the compositions $\bar g_\alpha\bar f\:R/I
\rarrow\overline T_\alpha$ are principal open immersions.
 By Corollary~\ref{formal-open-immersions-reduction-corollary},
there exist principal formal open immersions $g_\alpha\:S\rarrow
T_\alpha$ together with $S/Sf(I)$\+algebra isomorphisms
$T_\alpha/T_\alpha g_\alpha f(I)\simeq\overline T_\alpha$.
 By Lemma~\ref{formal-open-coverings-lemma}(1),
the collection~$(g_\alpha)$ is a formal open covering.
\end{proof}

\begin{lem} \label{open-covering-adic-ring-base-change-context}
 Let $R$ and $S$ be adic topological rings, let $f\:R\rarrow S$
be a continuous ring map, and let $g_\alpha\:R\rarrow T_\alpha$ be
a formal open covering.
 Put $W_\alpha=T_\alpha\ot_RS$, and denote by $f'_\alpha\:T_\alpha
\rarrow W_\alpha$ and $g'_\alpha\:S\rarrow W_\alpha$
the induced ring maps.
 Endow the commutative rings $W_\alpha$ with the tensor product
topologies.
 Then the collection of ring maps $g'_\alpha\:S\rarrow W_\alpha$
is a formal open covering, too.
\end{lem}

\begin{proof}
 The argument, based on the fact that any base change of an open
covering (of discrete commutative rings) is an open covering, is
straightforward.
\end{proof}

\begin{lem} \label{locality-of-contramodule-flatness}
 Let $f_\alpha\:R\rarrow S_\alpha$ be a formal open covering of
an adic topological ring~$R$.
 Then a quotseparated contramodule $R$\+module\/ $\fF$ is flat
if and only if the quotseparated contramodule $S_\alpha$\+module
$f_\alpha^\sharp(\fF)$ is flat for every index~$\alpha$.
 A contramodule $R$\+module\/ $\fF$ is flat if and only if
the contramodule $S_\alpha$\+module $f_\alpha^\#(\fF)$ is flat
for every index\/~$\alpha$.
\end{lem}

\begin{proof}
 Let us prove the quotseparated version.
 Let $I\subset R$ be an ideal of definition in~$R$.
 Then $J_\alpha=S_\alpha f_\alpha(I)\subset S_\alpha$ is an ideal
of definition in $S_\alpha$ for every index~$\alpha$.
 By Lemma~\ref{formal-open-coverings-lemma}(1), the collection of
induced ring maps $\bar f_\alpha\:R/I\rarrow S_\alpha/J_\alpha$ is
an open covering.
 By Lemma~\ref{reductions-of-co-contra-extension-of-scalars}(c),
we have $f_\alpha^\sharp(\fF)/Jf_\alpha^\sharp(\fF)\simeq
S_\alpha/J_\alpha\ot_{R/I}\fF/I\fF$.
 Now it remains to point out that the $R/I$\+module $\fF/I\fF$ is flat
if and only if the $S_\alpha/J_\alpha$\+module
$S_\alpha/J_\alpha\ot_{R/I}\fF/I\fF$ is flat for every index~$\alpha$.
 The proof of the nonquotseparated case is similar and based on
Lemma~\ref{reductions-of-co-contra-extension-of-scalars}(b).
\end{proof}

\subsection{Locality of Noetherianity}
\label{prelim-noetherianity-subsecn}
 The following lemma is well-known.

\begin{lem} \label{abstract-commutative-rings-noetherianity-local}
\textup{(a)} Let $f\:R\rarrow S$ be an open immersion of commutative
rings.
 Assume that the ring $R$ is Noetherian.
 Then the ring $S$ is Noetherian. \par
\textup{(b)} Let $f\:R\rarrow S_\alpha$ be an open covering of
a commutative ring~$R$.
 Assume that all the rings $S_\alpha$ are Noetherian.
 Then the ring $R$ is Noetherian.
\end{lem}

\begin{proof}
 See, e.~g., \cite[Lemma Tag~01OW]{SP}.
\end{proof}

 Recall from Section~\ref{adic-rings-subsecn} that an adic topological
ring $R$ with an ideal of definition $I\subset R$ is said to be
\emph{adically Noetherian} if the ring $R/I$ is Noetherian.

\begin{cor} \label{adic-topological-rings-adic-noetherianity-local-cor}
\textup{(a)} Let $f\:R\rarrow S$ be an formal open immersion of
adic topological rings.
 Assume that the ring $R$ is adically Noetherian.
 Then the ring $S$ is adically Noetherian. \par
\textup{(b)} Let $f\:R\rarrow S_\alpha$ be a formal open covering of
an adic topological ring~$R$.
 Assume that all the adic topological rings $S_\alpha$ are adically
Noetherian.
 Then the ring $R$ is adically Noetherian.
\end{cor}

\begin{proof}
 Follows from
Lemma~\ref{abstract-commutative-rings-noetherianity-local}.
\end{proof}

\begin{lem} \label{completion-noetherianity}
 Let $R$ be an adic topological ring.
 Then the following conditions are equivalent:
\begin{enumerate}
\item the adic topological ring $R$ is adically Noetherian;
\item the complete, separated adic topological ring\/ $\Lambda(R)$ is
adically Noetherian;
\item the ring\/ $\Lambda(R)$ is Noetherian (as an abstract
commutative ring).
\end{enumerate}
 In particular, a complete, separated adic topological ring\/ $\fR$ is
adically Noetherian if and only if it is Noetherian as an abstract
commutative ring.
\end{lem}

\begin{proof}
 (1)~$\Longleftrightarrow$~(2) Obvious.

 (2)~$\Longleftrightarrow$~(3) This is~\cite[Corollary~5 in
Section~III.2.10]{Bour}.
\end{proof}

\begin{cor} \label{adic-topological-rings-complet-noether-local-cor}
\textup{(a)} Let $f\:R\rarrow S$ be an formal open immersion of
adic topological rings.
 Assume that the ring\/ $\Lambda(R)$ is Noetherian.
 Then the ring\/ $\Lambda(S)$ is Noetherian. \par
\textup{(b)} Let $f\:R\rarrow S_\alpha$ be a formal open covering of
an adic topological ring~$R$.
 Assume that all the rings $\Lambda(S_\alpha)$ are Noetherian.
 Then the ring $\Lambda(R)$ is Noetherian.
\end{cor}

\begin{proof}
 Follows from
Corollary~\ref{adic-topological-rings-adic-noetherianity-local-cor}
and Lemma~\ref{completion-noetherianity}.
\end{proof}

\subsection{Tightness and tautness}  \label{prelim-tight-taut-subsecn}
 This section is largely based on~\cite[Section Tag~0GX1]{SP}
and~\cite[Sections~4\+-5]{Pcs}.
 The main result of the section is our version
of~\cite[Lemma Tag~0GBS]{SP} and~\cite[Lemmas~5.4 and~8.3]{Pcs}.

 In the following lemma, all the seven conditions are similar to
the respective conditions from Lemma~\ref{tight-ring-map-lemma},
while conditions~(5\+-6) are special cases of conditions~(1\+-2)
from~\cite[Lemma~4.2]{Pcs}.

\begin{lem} \label{taut-ring-map-lemma}
 Let $R$ and $S$ be adic topological rings, and let $f\:R\rarrow S$ be
a ring map.
 Then the following seven conditions are equivalent:
\begin{enumerate}
\item there exist an ideal of definition $I\subset R$ and an ideal of
definition $J\subset S$ such that $J$ is contained in the closure of
the ideal $Sf(I)\subset S$;
\item for every ideal of definition $J\subset S$ there exist an ideal
of definition $I\subset R$ and an integer $n\ge1$ such that the ideal
$J^n$ is contained in the closure of the ideal $Sf(I)\subset S$;
\item for any two ideals of definition $I\subset R$ and $J\subset S$
there exists an integer $n\ge1$ such that the ideal $J^n$ is contained
in the closure of the ideal $Sf(I)\subset S$;
\item for every ideal of definition $I\subset R$ there exists an ideal
of definition $J\subset S$ such that $J$ is contained in the closure of
the ideal $Sf(I)\subset S$;
\item for every open ideal $I\subset R$, the closure $J$ of the ideal
$Sf(I)\subset S$ is an open ideal in~$S$;
\item the open ideals $I\subset R$ for which the closure $J$ of
the ideal $Sf(I)\subset S$ is an open ideal in $S$ form a base of
neighborhoods of zero in~$R$;
\item there exists an ideal of definition $I\subset R$ such that
the closure $J$ of the ideal $Sf(I)\subset S$ is an open ideal in~$S$.
\end{enumerate}
\end{lem}

\begin{proof}
 (1)~$\Longrightarrow$~(5) Let $I'\subset R$ and $J'\subset S$ be
ideals of definition such that $J'$ is contained in the closure of
the ideal $Sf(I')\subset S$.
 Then there exists an integer $n\ge1$ such that $(I')^n\subset I$.
 Under any continuous map of topological spaces, the image of
the closure of a subset is contained in the closure of the image.
 Hence the closure of the ideal $Sf(I')^n\subset S$ contains $(J')^n$,
so the closure of $Sf(I')^n$ is an open ideal in~$S$.
 Since $Sf(I')^n\subset Sf(I)$, it follows that the closure of
$Sf(I)$ is an open ideal in $S$, too.

 (5)~$\Longrightarrow$~(3) $\Longrightarrow$~(2) $\Longrightarrow$~(1)
and (3)~$\Longrightarrow$~(4) $\Longrightarrow$~(1) Obvious.

 (5)~$\Longleftrightarrow$~(6) $\Longrightarrow$~(7)
$\Longleftrightarrow$~(1) Obvious.
\end{proof}

 We will say that a ring homomorphism $f\:R\rarrow S$ acting between
two adic topological rings is a \emph{taut} ring map if the equivalent
conditions of Lemma~\ref{taut-ring-map-lemma} hold
(cf.~\cite[Section Tag~0GX1]{SP} and~\cite[Section~4]{Pcs}).
 Notice that a taut ring map \emph{need not} be continuous in our
definition.
 Using the criterion of Lemma~\ref{taut-ring-map-lemma}(4) or~(5), one
can easily check that the composition of any two taut ring maps of
adic topological rings is taut.

 In the following lemma, we do \emph{not} say that ``$J$~is an ideal
of definition in~$S$'' because, in our terminology, an ideal of
definition must be finitely generated.

\begin{lem} \label{continuous-and-taut-ring-map-lemma}
 Let $R$ and $S$ be adic topological rings, and let $f\:R\rarrow S$ be
a ring map.
 Then the following four conditions are equivalent:
\begin{enumerate}
\item the ring homomorphism~$f$ is continuous and taut;
\item there exists an ideal of definition $I\subset R$ such that
the closure $J$ of the ideal $Sf(I)\subset S$ is an open ideal in $S$
and the powers $J^n$ of the ideal $J$ form a base of neighborhoods of
zero in~$S$;
\item for every ideal of definition $I\subset R$, the closure $J$ of
the ideal $Sf(I)\subset S$ is an open ideal in $S$ and the powers $J^n$
of the ideal $J$ form a base of neighborhoods of zero in~$S$;
\item for every ideal of definition $I\subset R$, the closure $J$ of
the ideal $Sf(I)\subset S$ is an open ideal in $S$, for every
integer $n\ge1$, the closure of the ideal $Sf(I^n)\subset S$ coincides
with the ideal $J^n\subset S$, and the ideals $J^n$ form a base of
neighborhoods of zero in~$S$.
\end{enumerate}
\end{lem}

\begin{proof}
 (2)~$\Longrightarrow$~(1) Assuming~(2), the ring map~$f$ is taut by
Lemma~\ref{taut-ring-map-lemma}(1).
 Furthermore, let $J'\subset S$ be an ideal of definition.
 Then there exists $n\ge1$ such that $J^n\subset J'$.
 Hence $f(I^n)\subset J^n\subset J'$, and the map~$f$ is continuous by
Lemma~\ref{continuous-ring-map-lemma}(2).

 (4)~$\Longrightarrow$~(3) $\Longrightarrow$~(2) Obvious.
 
 (1)~$\Longrightarrow$~(4) By Lemma~\ref{continuous-ring-map-lemma}(4),
there exists an ideal of definition $J'\subset S$ such that
$f(I)\subset J'$.
 Since $J'$ is an open ideal in $S$, we have $J\subset J'$.
 By Lemma~\ref{taut-ring-map-lemma}, \,$J$ is an open ideal in~$S$.
 Hence $J^n$ is an open ideal in $S$, too.
 In particular, $J^n$ is a closed ideal in $S$, and it follows that
$J^n$ is the closure of the ideal $Sf(I^n)\subset S$.
 Since $J\subset J'$, the ideals $J^n$ form a base of neighborhoods
of zero in~$S$.
\end{proof}

 Let $A$ be a topological abelian group where open subgroups
form a base of neighborhoods of zero.
 Then the \emph{completion} of $A$ is the topological abelian group
$\fA=\Lambda(A)=\varprojlim_{U\subset A}A/U$.
 Here the projective limit is taken over the direct poset of open
subgroups $U\subset A$, and the group $\fA$ is endowed with
the topology of projective limit of discrete groups~$A/U$.
 We denote by\/ $\lambda_A\:A\rarrow\Lambda(A)$ the natural (completion)
map.

\begin{lem} \label{open-closure-preserved-reflected-by-completion}
\textup{(a)} Let $A$ be a topological abelian group where open subgroups
form a base of neighborhoods of zero.
 Let $B\subset A$ be a subgroup.
 Then the closure of $B$ is open in $A$ if and only if the closure of
the subgroup\/ $\lambda_A(B)$ is open in\/~$\Lambda(A)$. \par
\textup{(b)} Let $R$ be an adic topological ring, and let $I\subset R$
be an ideal.
 Then the closure of $I$ is open in $R$ if and only if the closure
of the ideal $\Lambda(R)\lambda_R(I)$ is open in the adic topological
ring\/~$\Lambda(R)$.
\end{lem}

\begin{proof}
 Part~(a): the ``if'' assertion holds because the closure of $B$ in $A$
is the preimage of the closure of $\lambda_A(B)$ in $\Lambda(A)$ under
the natural (continuous) map~$\lambda_A$.
 To prove the ``only if'' assertion, one can observe that the preimage
map $\fV\longmapsto\lambda_A^{-1}(\fV)$ establishes a bijection
between the closed subgroups in $\Lambda(A)$ and the closed subgroups
in $A$, as well as between the open subgroups in $\Lambda(A)$ and
the open subgroups in~$A$.

 Part~(b) follows from part~(a), because the closure of the ideal
$\Lambda(R)\lambda_R(I)\subset\Lambda(R)$ coincides with the closure
of the subgroup $\lambda_R(I)\subset\Lambda(R)$.
 The point is that $I$ is an ideal in $R$ and $\Lambda(R)$ is
of its subring\/ $\lambda_R(R)\subset\Lambda(R)$.
\end{proof}

 Obviously, any tight ring map between adic topological rings is taut.
 The counterexample from
Remark~\ref{tightness-completion-remark} above shows that
a taut continuous ring map between adic topological rings 
\emph{need not} be tight.
 The following corollary should be also compared with
Remark~\ref{tightness-completion-remark}.

\begin{cor} \label{completion-preserves-reflects-tautness-cor}
 Let $R$ and $S$ be adic topological rings, and let $f\:R\rarrow S$
be a continuous ring map.
 Then the ring map~$f$ is taut if and only if the induced ring map\/
$\Lambda(f)\:\Lambda(R)\rarrow\Lambda(S)$ is taut.
\end{cor}

\begin{proof}
 Follows from Lemma~\ref{taut-ring-map-lemma}(5) or~(7) together with
Lemma~\ref{open-closure-preserved-reflected-by-completion}(b) (applied
to the adic topological ring $S$ with the ideal $Sf(I)\subset S$).
\end{proof}

 The following two lemmas form our version of~\cite[Lemma~4.6]{Pcs}.

\begin{lem} \label{taut-characterization}
 Let $R$ and $S$ be adic topological rings, and let $f\:R\rarrow S$ be
a ring map.
 Then the following five conditions are equivalent:
\begin{enumerate}
\item $f$~is a taut ring map;
\item there exist ideals of definition $I\subset R$ and $J\subset S$
such that $J^n\subset Sf(I^n)+J^{n+1}$ for every $n\ge1$;
\item for every ideal of definition $I\subset R$ there exists
an ideal of definition $J\subset S$ such that
$J^n\subset Sf(I^n)+J^{n+1}$ for every $n\ge1$;
\item there exist ideals of definition $I\subset R$ and $J\subset S$
such that $J\subset Sf(I)+J^2$;
\item for every ideal of definition $I\subset R$ there exists
an ideal of definition $J\subset S$ such that
$J\subset Sf(I)+J^2$.
\end{enumerate}
\end{lem}

\begin{proof}
 (3)~$\Longrightarrow$~(2) Obvious.

 (1)~$\Longrightarrow$~(3) By Lemma~\ref{taut-ring-map-lemma}(4),
there exists an ideal of definition $J\subset S$ such that $J$ is
contained in the closure of the ideal $Sf(I)\subset S$.
 It follows that the ideal $J^n$ is contained in the closure of
the ideal $Sf(I^n)\subset S$.
 Since $J^{n+1}$ is an open ideal in $S$, we can conclude that
$J^n\subset Sf(I^n)+J^{n+1}$.

 (2)~$\Longrightarrow$~(1)  We have
\begin{equation*}
\begin{split}
 J^n &\subset Sf(I^n) + J^{n+1}\\
       &\subset Sf(I^n) + Sf(I^{n+1}) + J^{n+2} \\
       &= Sf(I^n) + J^{n+2} \\
       &\subset Sf(I^n) + Sf(I^{n+2}) + J^{n+3} \\
       &= Sf(I^n) + J^{n+3} \\
       &\subset \dotsb \\
       &\subset Sf(I^n) + J^m
\end{split}
\end{equation*}
for every $m>n$.
 Since the ideals $J^m$ form a base of neighborhoods of zero in $S$,
it follows that $J^n$ is contained in the closure of $Sf(I^n)$ in~$S$.
 By Lemma~\ref{taut-ring-map-lemma}(1), this means that $f$~is a taut
ring map.

 (2)~$\Longrightarrow$~(4) and (3)~$\Longrightarrow$~(5) Obvious.
 
 (4)~$\Longrightarrow$~(2) and (5)~$\Longrightarrow$~(3) The inclusion
$J\subset Sf(I)+J^2$ implies $J^n\subset J^{n-1}f(I)+J^{n+1}\subset
Sf(I)+J^{n+1}$ for every $n\ge1$.
\end{proof}

\begin{lem} \label{continuous-taut-characterization}
 Let $R$ and $S$ be adic topological rings, and let $f\:R\rarrow S$ be
a ring map.
 Then the following three conditions are equivalent:
\begin{enumerate}
\item $f$~is a taut continuous ring map;
\item there exist an ideal of definition $I\subset R$ and
a descending chain of open ideals $J_1\supset J_2\supset J_3\supset
\dotsb$ forming a base of neighborhoods of zero in $S$ such that
$J_n=Sf(I^n)+J_{n+1}\subset S$ for every $n\ge1$;
\item for any ideal of definition $I\subset R$, there exists
a descending chain of open ideals $J_1\supset J_2\supset J_3\supset
\dotsb$ forming a base of neighborhoods of zero in $S$ such that
$J_n=Sf(I^n)+J_{n+1}\subset S$ for every $n\ge1$.
\end{enumerate}
\end{lem}

\begin{proof}
 (3)~$\Longrightarrow$~(2) Obvious.

 (1)~$\Longrightarrow$~(3)
 Denote by $J$ the closure of the ideal $Sf(I)\subset S$ in~$S$.
 By Lemma~\ref{continuous-and-taut-ring-map-lemma}(4), \,$J$ is
an open ideal in $S$, the ideals $J^n$ form a base of neighborhoods
of zero in~$S$, and $J^n$ is the closure of the ideal
$Sf(I^n)\subset S$ in~$S$.
 Since the ideal $J^{n+1}$ is open in $S$ and $J^{n+1}\subset J^n$,
we have $J^n=Sf(I^n)+J^{n+1}$.
 It remains to put $J_n=J^n$ for every $n\ge1$.

 (2)~$\Longrightarrow$~(1) We have
\begin{equation*}
\begin{split}
 J_n &= Sf(I^n) + J_{n+1}\\
       &= Sf(I^n) + Sf(I^{n+1}) + J_{n+2} \\
       &= Sf(I^n) + J_{n+2} \\
       &= Sf(I^n) + Sf(I^{n+2}) + J_{n+3} \\
       &= Sf(I^n) + J_{n+3} \\
       &= \dotsb \\
       &= Sf(I^n) + J_m
\end{split}
\end{equation*}
for every $m>n$.
 Since the ideals $J_m$ form a base of neighborhoods of zero in $S$,
it follows that $J_n$ is the closure of $Sf(I^n)$ in~$S$.
 Using Lemma~\ref{taut-ring-map-lemma}(6) or~(7), we can conclude
that the map~$f$ is taut.

 Now let $J\subset S$ be an open ideal.
 Then there exists $n\ge1$ such that $J_n\subset J$.
 Hence $I^n\subset f^{-1}(J_n)\subset f^{-1}(J)$.
 So the map~$f$ is continuous.
\end{proof}

 The following proposition is our version
of~\cite[Lemma Tag~0GBS]{SP}.

\begin{prop} \label{complete-tightness=tautness-prop}
 Let $R$ be an adic topological ring, $\fS$ be complete, separated adic
topological ring, and let $f\:R\rarrow\fS$ be a ring map.
 Then the map~$f$ is tight if and only if it is taut.
\end{prop}

\begin{proof}
 The ``only if'' assertion obviously holds for any ring map
$f\:R\rarrow S$ between adic topological rings $R$ and~$S$.
 To prove the ``if'', let $I\subset R$ be an ideal of definition,
and let $\fJ\subset\fS$ be a related ideal of definition from
Lemma~\ref{taut-characterization}(5).
 In order to prove that $\fS f(I)$ is an open ideal in $\fS$, it
suffices to show that $\fJ\subset\fS f(I)$.
 The following argument is our version of~\cite[proof of
Lemma~5.4]{Pcs}.

 Let $s\in\fJ$ be an element, and let $r_1$,~\dots, $r_m\in I$ be
a finite set of generators of the ideal $I\subset R$.
 Put $s_1=s$.
 By Lemma~\ref{taut-characterization}(5), we have
$\fJ\subset\fS f(I)+\fJ^2$.
 So there exists a sequence of elements $t_{0,1}$,~\dots,
$t_{0,m}\in\fS$ such that
$s_2=s-\sum_{j=1}^m t_{0,j}f(r_j)\in\fJ^2$.

 Furthermore, the inclusion $\fJ\subset\fS f(I)+\fJ^2$ implies
$\fJ^2\subset\fJ f(I)+\fJ^3$.
 Hence there exists a sequence of elements $t_{1,1}$,~\dots,
$t_{1,m}\in\fJ$ such that $s_3=s_2-\sum_{j=1}^mt_{1,j}f(r_j)
\in\fJ^3$.
 Similarly, the inclusion $\fJ\subset\fS f(I)+\fJ^2$ implies
$\fJ^3\subset\fJ^2 f(I)+\fJ^4$.
 Hence there exists a sequence of elements $t_{2,1}$,~\dots,
$t_{2,m}\in\fJ^2$ such that $s_4=s_3-\sum_{j=1}^mt_{2,j}f(r_j)
\in\fJ^4$.

 Proceeding in this way, we construct a sequence of elements
$s_n\in\fJ^n$, \,$n\ge1$, \,$s_1=s$, together with
$m$~sequences of elements $t_{n,j}\in\fJ^n$, \,$1\le j\le m$, such that
$s_{n+1}=s_n-\sum_{j=1}^m t_{n-1,j}f(r_j)$ for every $n\ge1$.
 Finally, in the topology of the complete, separated adic
topological ring $\fS$, we have
$$
 s=s_1=\sum\nolimits_{n=1}^\infty\sum\nolimits_{j=1}^m
 t_{n-1,j}f(r_j)
 = \sum\nolimits_{j=1}^m \left(f(r_j)\sum\nolimits_{n=0}^\infty
 t_{n,j}\right),
$$
where the sums $\sum_{n=0}^\infty t_{n,j}$ converge in $\fS$
for all $1\le j\le m$ since $t_{n,j}\in\fJ^n$.
 Thus $s\in f(I)\fS$.

 This proves the inclusion $\fJ\subset\fS f(I)$.
 By Lemma~\ref{tight-ring-map-lemma}(1), (4), or~(7), it follows that
the ring map~$f$ is tight.
\end{proof}

\subsection{Locality of tautness and tightness}
\label{prelim-locality-tautness-tightness-subsecn}
 The following technical lemma will be used in the proofs of
Propositions~\ref{locality-of-tautness-in-total-space-prop}
and~\ref{locality-of-tautness-in-the-base-prop}.

\begin{lem} \label{formal-open-covering-topology-base-lemma}
 Let $f_\alpha\:R\rarrow S_\alpha$, \,$1\le\alpha\le N$, be a finite
formal open covering of an adic topological ring~$R$.
 Suppose given, for every\/~$\alpha$, a descending chain of open
ideals $J_{\alpha,1}\supset J_{\alpha,2}\supset J_{\alpha,3}\supset
\dotsb$ forming a base of neighborhoods of zero in~$S_\alpha$.
 Then the open ideals $I_n=\bigcap_\alpha f_\alpha^{-1}(J_{\alpha,n})$
form a base of neighborhoods of zero in~$R$.
\end{lem}

\begin{proof}
 Let $I'\subset R$ be an open ideal.
 Then, for every index~$\alpha$, the ideal $J'_\alpha=
S_\alpha f_\alpha(I')$ is open in $S_\alpha$ (since the ring
map~$f_\alpha$ is tight).
 Hence there exists an integer $n\ge1$ such that $J_{\alpha,n}
\subset J'_\alpha$ for every~$\alpha$.
 We claim that $I_n\subset I'$.
 Indeed, consider the discrete quotient rings $\overline R=R/I'$
and $\overline S_\alpha=S_\alpha/J'_\alpha$, and let
$\bar f_\alpha\:\overline R\rarrow\overline S_\alpha$ be the induced
ring maps.
 Then the collection of ring maps~$\bar f_\alpha$ is an open covering.
 Therefore, the map $\overline R\rarrow\bigoplus_\alpha\overline
S_\alpha$ is injective, and it follows that $I'=\bigcap_\alpha
f_\alpha^{-1}(J'_\alpha)$.
\end{proof}

\begin{prop} \label{locality-of-tautness-in-total-space-prop}
 Let $f\:R\rarrow S$ be a continuous homomorphism of adic topological
rings, and let $g_\alpha\:S\rarrow T_\alpha$ be a formal open covering.
 Then the ring map~$f$ is taut if and only if the composition
$g_\alpha f\:R\rarrow T_\alpha$ is taut for every index\/~$\alpha$.
\end{prop}

\begin{proof}
 The ``only if'' assertion holds since the formal open immersions are
tight (hence taut) and the compositions of taut ring maps are taut.
 To prove the ``if'', pick an ideal of definition $I\subset R$.
 Without loss of generality we can assume that
the set of indices~$\alpha$ is finite.
 For every $n\ge1$ and every index~$\alpha$, denote by $K_{\alpha,n}$
the closure of the ideal $T_\alpha g_\alpha f(I^n)\subset T_\alpha$.
 By Lemma~\ref{taut-ring-map-lemma}(5), the ideal $K_{\alpha,n}\subset
T_\alpha$ is open.

 Furthermore, denote by $J_n\subset S$ the intersection of the ideals
$g_\alpha^{-1}(K_{\alpha,n})\subset S$.
 Since the maps~$g_\alpha$ are continuous and the set of
indices~$\alpha$ is finite, the ideal $J_n$ is open in $S$ as well.
 Clearly, we have $f(I^n)\subset J_n$, hence $K_{\alpha,n}$ is
the closure of the ideal $T_\alpha g_\alpha(J_n)\subset T_\alpha$.
 As the map~$g_\alpha$ is tight, we actually have $K_{\alpha,n}=
T_\alpha g_\alpha(J_n)$.

 Put $\overline R=R/I^{n+1}$, \,$\overline S=S/J_{n+1}$, and
$\overline T_\alpha=T_\alpha/K_{\alpha,n+1}$.
 Let $\overline I=I^n/I^{n+1}\subset\overline R$, \
$\overline J=J_n/J_{n+1}\subset\overline S$, and
$\overline K_\alpha=K_{\alpha,n}/K_{\alpha,n+1}\subset
\overline T_\alpha$ be the related notation for ideals.
 Then the collection of ring homomorphisms $\bar g_\alpha\:\overline S
\rarrow\overline T_\alpha$ is an open covering.
 Furthermore, we have $\overline T_\alpha \bar g_\alpha\bar f
(\overline I)=\overline K_\alpha=\overline T_\alpha
\bar g_\alpha(\overline J)$ (where $\bar f$ is the ring homomorphism
$\overline R\rarrow\overline S$).
 Since an ideal in $\overline S$ is determined by the collection of
its extensions to $\overline T_\alpha$, it follows that
$\overline J=\overline S\bar f(\overline I)$.

 Thus $J_n=Sf(I^n)+J_{n+1}\subset S$.
 It remains to point that, for every index~$\alpha$, the collection
of ideals $K_{\alpha,n}$ is a base of neighborhoods of zero
in $T_\alpha$, by Lemma~\ref{continuous-and-taut-ring-map-lemma}(4).
 Hence the ideals $J_n$ form a base of neighborhoods of zero in $S$,
by Lemma~\ref{formal-open-covering-topology-base-lemma}.
 Applying Lemma~\ref{continuous-taut-characterization}(2) or~(3),
we conclude that $f$~is a taut ring map.
\end{proof}

\begin{cor} \label{locality-of-tightness-in-total-space-cor}
 Let $f\:R\rarrow\fS$ be a continuous homomorphism of adic topological
rings, with a complete, separated adic topological ring\/~$\fS$.
 Let $g_\alpha\:\fS\rarrow T_\alpha$ be a formal open covering.
 Then the ring map~$f$ is tight if and only if the composition
$g_\alpha f\:R\rarrow T_\alpha$ is tight for every index\/~$\alpha$.
\end{cor}

\begin{proof}
 The ``only if'' implication holds since the formal open immersions are
tight and the compositions of tight ring maps are tight.
 To prove the ``if'' implication, notice that tightness of
the maps~$g_\alpha f$ obviously implies their tautness.
 By Propositions~\ref{locality-of-tautness-in-total-space-prop},
it follows that the ring map~$f$ is taut.
 Applying Proposition~\ref{complete-tightness=tautness-prop}, we
conclude that $f$~is tight.
\end{proof}

\begin{lem} \label{tautness-base-change-context}
 Let $R$, $S$, and $T$ be adic topological rings, let $f\:R\rarrow S$
be a taut continuous ring map, and let $g\:R\rarrow T$ be
a continuous ring map.
 Put $W=T\ot_RS$, and denote by $f'\:T\rarrow W$ and $g'\:S\rarrow W$
the induced ring maps.
 Endow the commutative ring $W$ with the tensor product topology, as
in Section~\ref{prelim-tensor-products-of-adic-topological-subsecn}.
 Then the continuous ring map~$f'$ is also taut.
\end{lem}

\begin{proof}
 By Corollary~\ref{completion-preserves-reflects-tautness-cor},
the continuous ring map $\Lambda(f)\:\Lambda(R)\rarrow\Lambda(S)$
is taut.
 By Proposition~\ref{complete-tightness=tautness-prop}, it follows
that the ring map $\Lambda(f)$ is tight.
 Hence, by Lemma~\ref{tight-flat-adic-ring-base-change-context}(a),
the ring map $\Lambda(T)\rarrow\Lambda(T)\ot_{\Lambda(R)}\Lambda(S)$
is tight, and consequently, taut.
 By Lemma~\ref{completed-tensor-product-of-completions-lemma},
we have $\Lambda(W)=\Lambda(\Lambda(T)\ot_{\Lambda(R)}\Lambda(S))$.
 Using Corollary~\ref{completion-preserves-reflects-tautness-cor} again,
we conclude that the ring map $\Lambda(f')\:\Lambda(T)\rarrow\Lambda(W)$
is taut.
 By the same corollary, it follows that the ring map~$f'$ is taut.
\end{proof}

\begin{prop} \label{locality-of-tautness-in-the-base-prop}
 Let $f\:R\rarrow S$ be a continuous homomorphism of adic topological
rings, and let $g_\alpha\:R\rarrow T_\alpha$ be a formal open covering.
 Put $W_\alpha=T_\alpha\ot_RS$, and endow the rings $W_\alpha$ with
the tensor product topologies.
 Denote by $f'_\alpha\:T_\alpha\rarrow W_\alpha$ and $g'_\alpha\:
S\rarrow W_\alpha$ the induced ring maps.
 Then the ring map~$f$ is taut if and only if the map~$f'_\alpha$ is
taut for every index\/~$\alpha$.
\end{prop}

\begin{proof}
 The ```only if'' assertion holds by
Lemma~\ref{tautness-base-change-context}.
 The proof of the ``if'' resembles that of the ``only if'' assertion
in Proposition~\ref{locality-of-tautness-in-total-space-prop}.

 Without loss of generality we can assume that the set of
indices~$\alpha$ is finite.
 Pick an ideal of definition $I\subset R$.
 For every $n\ge1$ and every index~$\alpha$, put $K_{\alpha,n}=
T_\alpha g_\alpha(I^n)\subset T_\alpha$, and denote by $L_{\alpha,n}$
the closure of the ideal $W_\alpha f'_\alpha(K_{\alpha,n})\subset
W_\alpha$.
 By Lemmas~\ref{tight-ring-map-lemma}(5)
and~\ref{taut-ring-map-lemma}(5), the ideals $K_{\alpha,n}\subset 
T_\alpha$ and $L_{\alpha,n}\subset W_\alpha$ are open (since
the maps~$g_\alpha$ are tight and the maps~$f'_\alpha$ are taut).
 Furthermore, by Lemmas~\ref{continuous-and-tight-ring-map-lemma}(5)
and~\ref{continuous-and-taut-ring-map-lemma}(4), the ideals
$K_{\alpha,n}$, \,$n\ge1$, form a base of neighborhoods of zero
in $T_\alpha$, while the ideals $L_{\alpha,n}$, \,$n\ge1$, form a base
of neighborhoods of zero in $W_\alpha$ for every~$\alpha$.

 Furthermore, denote by $J_n\subset S$ the intersection of the ideals
$(g'_\alpha)^{-1}(L_{\alpha,n})\subset S$.
 Since the maps~$g'_\alpha$ are continuous and the set of
indices~$\alpha$ is finite, the ideal $J_n$ is open in~$S$.
 Clearly, we have $f(I^n)\subset J_n$, hence $L_{\alpha,n}$ is
the closure of the ideal $W_\alpha g'_\alpha(J_n)\subset W_\alpha$.
 As the map~$g'_\alpha$ is tight by
Lemma~\ref{tight-flat-adic-ring-base-change-context}(a)
(and in fact a formal open immersion by
Lemma~\ref{open-immersion-adic-ring-base-change-context}(a)),
we actually have $L_{\alpha,n}=W_\alpha g'_\alpha(J_n)$.

 Put $\overline R=R/I^{n+1}$, \,$\overline S=S/J_{n+1}$,
\,$\overline T_\alpha=T_\alpha/K_{\alpha,n+1}$, and
$\overline W_\alpha=W_\alpha/L_{\alpha,n+1}$.
 Denote by $\bar f\:\overline R\rarrow\overline S$, \
$\bar g_\alpha\:\overline R\rarrow\overline T_\alpha$, \
$\bar f'_\alpha\:\overline T_\alpha\rarrow\overline W_\alpha$,
and $\bar g'_\alpha\:\overline S\rarrow\overline W_\alpha$
the induced ring maps.
 Then the collection of ring homomorphisms~$\bar g_\alpha$ is an open
covering.
 The induced ring maps $\overline T_\alpha\ot_{\overline R}
\overline S\rarrow\overline W_\alpha$ are isomorphisms, since
$L_{\alpha,n+1}=W_\alpha g'_\alpha(J_{n+1})$.
 Hence the collection of ring homomorphisms~$\bar g'_\alpha$ is
an open covering, too.

 Let $\overline I=I^n/I^{n+1}\subset\overline R$, \
$\overline J=J_n/J_{n+1}\subset\overline S$, \
$\overline K_\alpha=K_{\alpha,n}/K_{\alpha,n+1}\subset
\overline T_\alpha$, and $\overline L_\alpha=
L_{\alpha,n}/L_{\alpha,n+1}\subset\overline W_\alpha$
be the related notation for ideals.
 We have $\overline W_\alpha\bar g'_\alpha(\overline J)=
\overline L_\alpha=\overline W_\alpha\bar f'_\alpha(\overline K)
=\overline W_\alpha\bar f'_\alpha\bar g_\alpha(\overline I)=
\overline W_\alpha\bar g'_\alpha\bar f(\overline I)$.
 Since an ideal in $\overline S$ is determined by its extensions
to $\overline W_\alpha$, it follows that $\overline J=
\overline S\bar f(\overline I)$.

 Thus $J_n=Sf(I^n)+J_{n+1}\subset S$.
 Since the collection of ideals $L_{\alpha,n}$ is a base of
neighborhoods of zero in $W_\alpha$, the ideals $J_n$ form a base of
neighborhoods of zero in~$S$
(by Lemmas~\ref{open-covering-adic-ring-base-change-context}
and~\ref{formal-open-covering-topology-base-lemma}).
 Applying Lemma~\ref{continuous-taut-characterization}(2) or~(3),
we conclude that $f$~is a taut ring map.
\end{proof}

\begin{cor} \label{locality-of-tightness-in-the-base-cor}
 Let $f\:R\rarrow S$ be a continuous homomorphism of adic topological
rings, and let $g_\alpha\:R\rarrow T_\alpha$ be a formal open covering.
 Put $W_\alpha=T_\alpha\ot_RS$, and endow the rings $W_\alpha$ with
the tensor product topologies.
 Denote by $f'_\alpha\:T_\alpha\rarrow W_\alpha$ and $g'_\alpha\:
S\rarrow W_\alpha$ the induced ring maps.
 Then the ring map\/ $\Lambda(f)$ is tight if and only if the ring
map\/ $\Lambda(f'_\alpha)$ is tight for every index\/~$\alpha$.
\end{cor}

\begin{proof}
 Follows from Proposition~\ref{locality-of-tautness-in-the-base-prop}
in view of Corollary~\ref{completion-preserves-reflects-tautness-cor}
and Proposition~\ref{complete-tightness=tautness-prop}.
\end{proof}

\subsection{Very flat and contraadjusted contramodules}
\label{prelim-veryflat-and-contraadjusted-subsecn}
 In the sequel, we will use a shorthand notation $R\Contra=
R\Modl_\ctra^\qs$ for any adic topological ring~$R$.
 In the case of a complete, separated adic topological ring $\fR$,
a contramodule over $\fR$ as a topological ring is essentially
the same thing as a quotseparated contramodule $\fR$\+module
(see the end of Section~\ref{prelim-contramodules-subsecn}),
so this notation is unambiguous.
 Let us also recall a similar notation $R\Tors=R\Modl_\tors$ from
Section~\ref{prelim-torsion-modules-subsecn}.

 We start with a recollection of the case of an abstract
commutative ring $R$ (without any adic or other topology).
 An $R$\+module $P$ is said to be
\emph{contraadjusted}~\cite[Section~1.1]{Pcosh}, \cite[Section~5]{ST},
\cite[Sections~2 and~8]{Pcta}, \cite[Section~0.5]{PSl1},
\cite[Section~2 and Example~3.2]{Pal}, \cite[Section~4.3]{Pphil}
if $\Ext^1_R(R[s^{-1}],P)=0$ for all elements $s\in R$.
 An $R$\+module $F$ is said to be \emph{very
flat}~\cite[Section~1.1]{Pcosh}, \cite[Section~2]{ST},
\cite[Section~8]{Pcta}, \cite[Section~0.5]{PSl1},
\cite[Section~2 and Example~2.5]{Pal}, \cite[Section~4.3]{Pphil}
if $\Ext^1_R(F,C)=0$ for all contraadjusted $R$\+modules~$C$.

 The full subcategory of contraadjusted $R$\+modules is closed under
extensions, quotients, and infinite products in $R\Modl$.
 The full subcategory of very flat $R$\+modules is closed under
extensions, kernels of epimorphisms, infinite direct sums, and
tensor products in $R\Modl$.
 All very flat $R$\+modules have projective dimensions at most~$1$.

 A homomorphism of commutative rings $f\:R\rarrow S$ is said to be
\emph{very flat} if, for every element $s\in S$, the $R$\+module
$S[s^{-1}]$ is very flat.
 If this is the case, the ring $S$ is also said to be a \emph{very
flat commutative $R$\+algebra}.
 Very flatness of the $R$\+module $S$ is \emph{not} sufficient for
the ring homomorphism $R\rarrow S$ to be very
flat~\cite[Example~9.7]{PSl1}, \cite[Section~6.1]{Pphil}.
 The composition of any two very flat homomorphisms of commutative
rings is very flat~\cite[Lemma~1.2.3(b)]{Pcosh}.
 Every open immersion of commutative rings is very
flat~\cite[Lemma~1.2.4]{Pcosh}.

\begin{lem} \label{very-flat-contramodules-lemma}
 Let $R$ be an adic topological ring and\/ $\fF$ be
a contramodule $R$\+module.
 Then the following two conditions are equivalent:
\begin{enumerate}
\item $\fF$~is a flat contramodule $R$\+module and there exists
an ideal of definition $I\subset R$ for which\/ $\fF/I\fF$ is
a very flat $R/IR$\+module;
\item for every open ideal $I\subset R$, the $R/I$\+module\/
$\fF/I\fF$ is very flat.
\end{enumerate}
\end{lem}

\begin{proof}
 (2)~$\Longrightarrow$~(1) holds because all very flat $R/I$\+modules
are flat.

 (1)~$\Longrightarrow$~(2) Let $I\subset I'\subset R$ be two ideals
of definition in~$R$.
 Then $I'/I$ is a finitely generated nilpotent ideal in the ring
$R/I$.
 Therefore, assuming that the $R/I$\+module $\fF/I\fF$ is flat,
the $R/I$\+module $\fF/I\fF$ is very flat if and only if
the $R/I'$\+module $\fF/I'\fF$ is very
flat~\cite[Lemma~1.7.10(b)]{Pcosh}.
 Now if $I\subset I'\subset R$ are two ideals such that $I$ is
an ideal of definition, and the $R/I$\+module $\fF/I\fF$ is very
flat, then the $R/I'$\+module $\fF/I'\fF$ is very flat
by~\cite[Lemma~1.2.2(b)]{Pcosh}.
\end{proof}

 Let $R$ be an adic topological ring.
 A contramodule $R$\+module $\fF$ is said to be \emph{very flat} if
it satisfies the equivalent conditions of
Lemma~\ref{very-flat-contramodules-lemma}.
 All projective contramodule $R$\+modules and all projective
quotseparated contramodule $R$\+modules are very flat (see the proof
of Lemma~\ref{projective-contramodules-characterizations}%
\,(2)\,$\Rightarrow$\,(3)).
 It follows from
Lemma~\ref{flat-quotseparated-contramodules-well-behaved} that
the full subcategory of very flat quotseparated contramodule
$R$\+modules is closed under extensions and kernels of epimorphisms
in $R\Modl_\ctra^\qs$.
 Furthermore, in view of~\cite[Corollary~E.1.10(a)]{Pcosh},
all very flat quotseparated contramodule $R$\+modules have
projective dimensions at most~$1$ as objects of $R\Modl_\ctra^\qs$.
 These assertions are special cases of~\cite[Corollary~E.4.1]{Pcosh}.

\begin{lem} \label{flat-derived-completion-Ext1-adjunction}
 Let $R$ be an adic topological ring and $F$ be a flat $R$\+module.
 Then, for any quotseparated contramodule $R$\+module\/ $\fP$,
there is a natural\/ $\Ext^1$\+adjunction isomorphism of abelian
groups/$R$\+modules
$$
 \Ext^1_R(F,\fP)\simeq\Ext^1_{R\Contra}(\Lambda(F),\fP).
$$
\end{lem}

\begin{proof}
 By~\cite[Lemma~3.5 or Proposition~3.6]{PSY}, we have
$\boL_i\Lambda(F)=0$ for all $i>0$.
 Therefore, for any short exact sequence of $R$\+modules
$0\rarrow M\rarrow N\rarrow F\rarrow0$, the short sequence of
quotseparated contramodule $R$\+modules $0\rarrow\boL_0\Lambda(M)
\rarrow\boL_0\Lambda(N)\rarrow\boL_0\Lambda(F)\rarrow0$ is exact.

 Now we have a pair of adjoint functors between abelian categories
$$
 \xymatrix{
  R\Contra \ar@{>->}@<-2pt>[rr]&& R\Modl,
  \ar@<-4pt>[ll]_-{\boL_0\Lambda}
 }
$$
where the right adjoint functor $R\Contra\rarrow R\Modl$ is the identity
inclusion.
 The right adjoint functor is exact, while the object $F\in R\Modl$ is
adjusted to the left adjoint functor $\boL_0\Lambda$ in the sense
stated in the previous paragraph.
 By~\cite[Lemma~1.7(e)]{Pal}, the desired $\Ext^1$\+adjunction
isomorphism follows.
 Here we are also using the natural $R$\+module isomorphism
$\boL_0\Lambda(F)\simeq\Lambda(F)$, which holds by the same results
from~\cite{PSY} or by
Corollary~\ref{flat-reductions-derived-Lambda-is-underived-cor} above.
\end{proof}

\begin{lem} \label{contraadjusted-contramodules-lemma}
 Let $R$ be an adic topological ring, $\fR$ be the adic completion
of $R$, and\/ $\fP$ be a quotseparated contramodule $R$\+module.
 Then the following seven conditions are equivalent:
\begin{enumerate}
\item the functor\/ $\Hom_R({-},\fP)$ takes short exact sequences of
very flat quotseparated contramodule $R$\+modules to short exact
sequences of abelian groups/$R$\+mod\-ules/quotseparated
contramodule $R$\+modules;
\item one has\/ $\Ext^1_{R\Contra}(\fF,\fP)=0$ for all very flat
quotseparated contramodule $R$\+modules\/~$\fF$ (where the Ext group
is computed in the abelian category $R\Contra=R\Modl_\ctra^\qs$);
\item one has\/ $\Ext^1_{R\Contra}(\Lambda(R[r^{-1}]),\fP)=0$ for all
elements $r\in R$;
\item $\fP$~is a contraadjusted $R$\+module;
\item $\fP$~is a contraadjusted\/ $\fR$\+module;
\item for every ideal $I\subset R$, the $R/I$\+module\/ $\fP/I\fP$
is contraadjusted;
\item for every ideal $I\subset\fR$, the\/ $\fR/I$\+module\/ $\fP/I\fP$
is contraadjusted.
\end{enumerate}
 Furthermore, if\/ $\fP$ is a separated contramodule $R$\+module,
then conditions~\textup{(1\+-7)} are also equivalent to the following
two conditions:
\begin{enumerate}
\setcounter{enumi}{7}
\item there exists an ideal of definition $I\subset R$ for which\/
$\fP/I\fP$ is a contraadjusted $R/I$\+module;
\item for every open ideal $I\subset R$, the $R/I$\+module\/ $\fP/I\fP$
is contraadjusted.
\end{enumerate}
\end{lem}

\begin{proof}
 (1)~$\Longleftrightarrow$~(2) Follows from the facts that
there are enough very flat quotseparated contramodules in
$R\Modl_\ctra^\qs$ and the full subcategory of very flat
quotseparated contramodules is closed under kernels of epimorphisms.
 See Lemma~\ref{Hom-from-to-short-exact-sequences-lemma}(a) below
for the details.

 (3)~$\Longleftrightarrow$~(4) Follows immediately from
Lemma~\ref{flat-derived-completion-Ext1-adjunction} (applied to
the flat $R$\+module $F=R[r^{-1}]$).

 (1)~$\Longleftrightarrow$~(4)
 This is a special case
of~\cite[first assertion of Corollary~E.4.8]{Pcosh}.
 In the case of a Noetherian ring $R$, one can also refer
to~\cite[Corollary~D.3.5(c,e)]{Pcosh}.

 (1)~$\Longleftrightarrow$~(5)
 This is a particular case of (1)~$\Longleftrightarrow$~(4).
 
 (6)~$\Longrightarrow$~(4) or (7)~$\Longrightarrow$~(5)
 Take $I=0$.
 
 (4)~$\Longrightarrow$~(6) or (5)~$\Longrightarrow$~(7)
 This is~\cite[Lemma~1.7.7(b)]{Pcosh}.
 
 (6)~$\Longrightarrow$~(9) and (7)~$\Longrightarrow$ (9)
$\Longrightarrow$~(8) Obvious.

 (8)~$\Longrightarrow$~(9)
 This is~\cite[Lemma~1.7.10(a)]{Pcosh}.
 
 (9)~$\Longrightarrow$~(2)
 This is a special case of~\cite[Lemma~E.4.3]{Pcosh}.
\end{proof}

 A quotseparated contramodule $R$\+module $\fP$ is said to be
\emph{contraadjusted} if it satisfies any one of the equivalent
conditions~(1\+-7) of Lemma~\ref{contraadjusted-contramodules-lemma}.
 The full subcategory of contraadjusted quotseparated contramodule
$R$\+modules is closed under extensions, quotients, and infinite
products in $R\Modl_\ctra^\qs$.

\begin{prop} \label{very-flat-cotorsion-pair-in-quotseparated}
 Let $R$ be an adic topological ring and\/ $\fM$ be a quotseparated
contramodule $R$\+module.  Then \par
\textup{(a)} there exists a short exact sequence of quotseparated
contramodule $R$\+modules\/ $0\rarrow\fP\rarrow\fF\rarrow\fM\rarrow0$
with a very flat quotseparated contramodule $R$\+module\/ $\fF$ and
a contraadjusted quotseparated contramodule $R$\+module\/~$\fP$; \par
\textup{(b)} there exists a short exact sequence of quotseparated
contramodule $R$\+modules\/ $0\rarrow\fM\rarrow\fP\rarrow\fF\rarrow0$
with a contraadjusted quotseparated contramodule $R$\+module\/ $\fP$
and a very flat quotseparated contramodule $R$\+module\/~$\fF$.
\end{prop}

\begin{proof}
 Part~(a) is a special case of~\cite[Corollary~E.4.5]{Pcosh}.
 For a Noetherian ring $R$, one can also refer
to~\cite[Corollary~D.3.5(b,e)]{Pcosh}.
 Part~(b) is a special case of~\cite[Corollary~E.4.6]{Pcosh}.
 For a Noetherian ring $R$, one can also refer
to~\cite[Corollary~D.3.5(a,e)]{Pcosh}.
 Another reference is~\cite[Example~7.12(3)]{PR}.
\end{proof}

\begin{cor} \label{very-flat-contramodules-cor}
 Let $R$ be an adic topological ring and\/ $\fF$ be a quotseparated
contramodule $R$\+module.
 Then the following conditions are equivalent:
\begin{enumerate}
\item the functor\/ $\Hom_R(\fF,{-})$ takes short exact sequences of
contraadjusted quotseparated contramodule $R$\+modules to short
exact sequences of abelian groups/$R$\+mod\-ules/quotseparated
contramodule $R$\+modules; {\hbadness=1725\par}
\item one has\/ $\Ext^1_{R\Contra}(\fF,\fP)=0$ for all
contraadjusted quotseparated contramodule $R$\+modules\/~$\fP$;
\item $\fF$~is a very flat quotseparated contramodule $R$\+module.
\end{enumerate}
\end{cor}

\begin{proof}
 (1)~$\Longleftrightarrow$~(2) Follows from the facts that
there are enough contraadjusted quotseparated contramodules in
$R\Modl_\ctra^\qs$ (by
Proposition~\ref{very-flat-cotorsion-pair-in-quotseparated}(b))
and the full subcategory of contraadjusted quotseparated contramodules
is closed under cokernels of monomorphisms.
 See Lemma~\ref{Hom-from-to-short-exact-sequences-lemma}(b) below
for the details.

 (3)~$\Longrightarrow$~(2) Holds by
Lemma~\ref{contraadjusted-contramodules-lemma}(2).

 (2)~$\Longrightarrow$~(3) This a standard argument (see
Lemma~\ref{direct-summand-lemma} below for a general formulation).
 By Proposition~\ref{very-flat-cotorsion-pair-in-quotseparated}(a),
there exists a short exact sequence $0\rarrow\fP\rarrow\fG\rarrow\fF
\rarrow0$ in $R\Modl_\ctra^\qs$ such that $\fG$ is very flat and
$\fP$ is contraadjusted.
 By~(2), it follows that $\fF$ is a direct summand of~$\fG$.
 It remains to make the obvious observation that the class of
very flat quotseparated contramodule $R$\+modules is closed under
direct summands.
 (For a Noetherian ring $R$, the reference
to~\cite[Corollary~D.3.5(d,e)]{Pcosh} is also applicable.)
\end{proof}

 In the terminology of Section~\ref{appx-cotorsion-pairs-subsecn}
below, the results of
Proposition~\ref{very-flat-cotorsion-pair-in-quotseparated}
and Corollary~\ref{very-flat-contramodules-cor}\,%
(2)\,$\Leftrightarrow$\,(3), together with the observation above
that the class of very flat quotseparated contramodules is
closed under kernels of epimorphisms in $R\Contra$, can be
summarized by the following formulation.
 For any adic topological ring $R$, the pair of classes (very flat
quotseparated contramodule $R$\+modules, contraadjusted quotseparated
contramodule $R$\+modules) is a \emph{hereditary complete cotorsion
pair} in the abelian category of quotseparated contramodule
$R$\+modules $R\Contra$.

 The next lemma is a generalization of~\cite[Lemma~1.2.1]{Pcosh}.

\begin{lem} \label{vfl-cta-tensor-Hom}
 Let $R$ be an adic topological ring.  Then \par
\textup{(a)} for any two very flat quotseparated contramodule
$R$\+modules\/ $\fF$ and\/ $\fG$, the quotseparated contramodule
$R$\+module\/ $\Lambda(\fF\ot_R\fG)$ is very flat; \par
\textup{(b)} for any very flat quotseparated contramodule $R$\+module\/
$\fF$ and any contraadjusted quotseparated contramodule $R$\+module\/
$\fP$, the quotseparated contramodule $R$\+module\/ $\Hom_R(\fF,\fP)$
is contraadjusted.
\end{lem}

\begin{proof}
 Part~(a): let $I\subset R$ be an open ideal.
 Then, by Lemma~\ref{derived-completion-reduction-isomorphism}, we have
$R/I\ot_R\Lambda(\fF\ot_R\fG)\simeq R/I\ot_R\fF\ot_R\fG\simeq
(R/I\ot_R\fF)\ot_{R/I}(R/I\ot_R\fG)$.
 Since the $R/I$\+modules $R/I\ot_R\fF$ and $R/I\ot_R\fG$ are very
flat, their tensor product is a very flat $R/I$\+module
by~\cite[Lemma~1.2.1(a)]{Pcosh}.

 Part~(b): in view of
Lemma~\ref{contraadjusted-contramodules-lemma}(1), it suffices to prove
that the functor $\fG\longmapsto\Hom_R(\fG,\Hom_(\fF,\fP))$ takes
short exact sequences of very flat quotseparated contramodule
$R$\+modules to short exact sequences of abelian groups.
 Indeed, we have $\Hom_R(\fG,\Hom(\fF,\fP))\simeq
\Hom_R(\fF\ot_R\fG,\>\fP)\simeq\Hom_R(\boL_0\Lambda(\fF\ot_R\fG),
\>\fP)\simeq\Hom_R(\Lambda(\fF\ot_R\fG),\>\fP)$ by
Corollary~\ref{flat-reductions-derived-Lambda-is-underived-cor}.
 In view of part~(a), it remains to point that the functor
$\fG\longmapsto\Lambda(\fF\ot_R\fG)$ takes short exact sequences
of very flat quotseparated contramodule $R$\+modules to short exact
sequences of (quotseparated contramodule) $R$\+modules.
 The latter property is a special case of
Corollary~\ref{completed-tensor-product-exact-on-flat-qs}.
\end{proof}

 The following lemma is the very flat counterpart of
Corollary~\ref{flat-contramodules-contraextension-of-scalars}.

\begin{lem} \label{contraextension-of-scalars-very-flat}
 Let $f\:R\rarrow S$ be a continuous homomorphism of adic topological
rings.
 Then, for any very flat quotseparated contramodule $R$\+module\/ $\fF$,
the quotseparated contramodule $S$\+module $f^\sharp(\fF)$ is very flat.
 For any very flat contramodule $R$\+module\/ $\fF$, the contramodule
$S$\+module $f^\#(\fF)$ is very flat.
\end{lem}

\begin{proof}
 The argument is similar to the proof of
Corollary~\ref{flat-contramodules-contraextension-of-scalars} and based
on Lemma~\ref{reductions-of-co-contra-extension-of-scalars}(b\+-c).
 The point is that, for any homomorphism of (abstract,
nontopological) commutative rings $\bar f\:\overline R\rarrow
\overline S$, the functor $\bar f^*\:\overline R\Modl\rarrow
\overline S\Modl$ takes very flat $R$\+modules to very flat
$S$\+modules~\cite[Lemma~1.2.2(b)]{Pcosh}.
\end{proof}

\begin{lem} \label{restriction-of-scalars-contraadjusted}
 Let $f\:R\rarrow S$ be a continuous homomorphism of adic topological
rings.
 Then, for any contraadjusted quotseparated contramodule $S$\+module\/
$\fP$, the quotseparated contramodule $R$\+module $f_\sharp(\fP)$
is contraadjusted.
\end{lem}

\begin{proof}
 Using the criterion of
Lemma~\ref{contraadjusted-contramodules-lemma}(4), the question
reduces to the case of a homomorphism of abstract (nontopological)
commutative rings, where we can refer to~\cite[Lemma~1.2.2(a)]{Pcosh}.
 Alternatively, an argument proceeding entirely within the realm of
quotseparated contramodules is possible, using the criterion of
Lemma~\ref{contraadjusted-contramodules-lemma}(2) and the results of
Corollary~\ref{flat-contramodules-adjusted-to-contraextension} and
Lemma~\ref{contraextension-of-scalars-very-flat}, and similar to
the proof of Lemma~\ref{restriction-of-scalars-cotorsion}
sketched below.
\end{proof}

 The next lemma is a very flat complement to
Lemma~\ref{flat-homomorphism-of-adic-rings}.

\begin{lem} \label{very-flat-homomorphism-of-adic-rings}
 Let $f\:R\rarrow S$ be a tight continuous map of adic
topological rings.
 Then the following three conditions are equivalent:
\begin{enumerate}
\item for every open ideal $I\subset R$, the induced ring map
$\bar f\:R/I\rarrow S/Sf(I)$ is very flat;
\item there exists an ideal of definition $I\subset R$ such that
the induced ring map $\bar f\:R/I^n\rarrow S/Sf(I^n)$ is very flat
for every $n\ge1$;
\item $f$~is flat (as a map of adic topological rings) and there exists
an ideal of definition $I\subset R$ such that the induced ring
map $\bar f\:R/I\rarrow S/Sf(I)$ is very flat.
\end{enumerate}
\end{lem}

\begin{proof}
 All the assertions follow from~\cite[Lemmas~1.2.2(b)
and~1.7.10(b)]{Pcosh}.
\end{proof}

 We will say that a tight continuous homomorphism of adic topological
rings $f\:R\rarrow S$ is \emph{very flat} if it satisfies the equivalent
conditions of Lemma~\ref{very-flat-homomorphism-of-adic-rings}.
 Any formal open immersion of adic topological rings is very flat.
 Using the criterion of
Lemma~\ref{very-flat-homomorphism-of-adic-rings}(1), one can easily
check that the composition of any two very flat tight continuous ring
maps of adic topological rings is very flat.

 The following lemma is the very flat version of
Lemma~\ref{induced-map-of-complete-rings-lemma}(b).
 The proof is the same.

\begin{lem} \label{induced-map-of-complete-rings-very-flat-lemma}
 Let $f\:R\rarrow S$ be a tight continuous homomorphism of adic
topological rings.
 The the map~$f$ is very flat if and only if the map\/ $\Lambda(f)$ is
very flat.  \qed
\end{lem}

 The following lemma is the very flat version of
Lemma~\ref{flat-map-of-adic-rings-direct-image-lemma}(c\+-d).

\begin{lem} \label{very-flat-map-of-adic-rings-direct-image-lemma}
 Let $f\:R\rarrow S$ be a very flat tight continuous map of adic
topological rings.
 Then \par
\textup{(a)} the functor $f_\#\:S\Modl_\ctra\rarrow R\Modl_\ctra$
takes very flat contramodule $S$\+modules to very flat contramodule
$R$\+modules; \par
\textup{(b)} the functor $f_\sharp\:S\Modl_\ctra^\qs\rarrow
R\Modl_\ctra^\qs$ takes very flat quotseparated contramodule
$S$\+modules to very flat quotseparated contramodule $R$\+modules.
\end{lem}

\begin{proof}
 The argument is similar to the proof of
Lemma~\ref{flat-map-of-adic-rings-direct-image-lemma}(c\+-d).
 The point is that, for any very flat homomorphism of (abstract,
nontopological) commutative rings $\bar f\:\overline R\rarrow
\overline S$, the functor $\bar f_*\:\overline S\Modl\rarrow
\overline R\Modl$ takes very flat $S$\+modules to very flat
$R$\+modules~\cite[Lemma~1.2.3(b)]{Pcosh},
\cite[Lemma~9.3(a)]{PSl1}, \cite[Lemma~6.1]{Pphil}.
\end{proof}

 The next lemma is the very flat version of
Lemma~\ref{locality-of-contramodule-flatness}.

\begin{lem} \label{locality-of-contramodule-very-flatness}
 Let $f_\alpha\:R\rarrow S_\alpha$ be a formal open covering of
an adic topological ring~$R$ (in the sense of
Section~\ref{prelim-formal-open-immersions-coverings-subsecn}).
 Then a quotseparated contramodule $R$\+module\/ $\fF$ is very flat
if and only if the quotseparated contramodule $S_\alpha$\+module
$f_\alpha^\sharp(\fF)$ is very flat for every index~$\alpha$.
 A contramodule $R$\+module\/ $\fF$ is very flat if and only if
the contramodule $S_\alpha$\+module $f_\alpha^\#(\fF)$ is very flat
for every index~$\alpha$.
\end{lem}

\begin{proof}
 The proof is similar to that of
Lemma~\ref{locality-of-contramodule-flatness} and
based on~\cite[Lemma~1.2.6(a)]{Pcosh} together with
Lemma~\ref{reductions-of-co-contra-extension-of-scalars}(b\+-c) above.
\end{proof}

\begin{lem} \label{very-flat-reductions-Hom-exactness-lemma}
 Let $R$ be an adic topological ring with an ideal of definition
$I\subset R$, and let $F$ be an $R$\+module such that
the $R/I^n$\+module $F/I^nF$ is very flat for all $n\ge1$.
 Let\/ $0\rarrow\fP\rarrow\fM\rarrow\fN\rarrow0$ be a short exact
sequence of quotseparated contramodule $R$\+modules with
a contraadjusted quotseparated contramodule $R$\+module\/~$\fP$.
 Then the short sequence of quotseparated contramodule $R$\+modules\/
$0\rarrow\Hom_R(F,\fP)\rarrow\Hom_R(F,\fM)\rarrow\Hom_R(F,\fN)\rarrow0$
is exact.
\end{lem}

\begin{proof}
 First of all, the Hom $R$\+modules in question are quotseparated
contramodules by Lemma~\ref{Hom-is-a-quotseparated-contramodule}(c).
 In the special case when $F$ is a flat $R$\+module, the exactness
of the short sequence in question follows from the fact that
$\Ext_R^1(F,\fP)=0$, see~\cite[Corollary~E.4.8]{Pcosh}.
 In the general case, we notice that $\Hom_R(F,\fN)\simeq
\Hom_R(\Lambda(F),\fN)$ and similarly for $\fM$ and $\fP$, 
by Corollary~\ref{flat-reductions-derived-Lambda-is-underived-cor}.
 By Lemma~\ref{derived-completion-reduction-isomorphism},
\,$\Lambda(F)$ is a very flat (quot)separated contramodule $R$\+module
under our current assumptions.
 Now the exactness of the short sequence of $\Hom_R(F,{-})$ follows
from the fact that $\Ext_{R\Contra}^1(\Lambda(F),\fP)=0$ by
Lemma~\ref{contraadjusted-contramodules-lemma}(2).
\end{proof}

 The following lemma is a partial very flat version
of Lemma~\ref{tight-flat-adic-ring-base-change-context}(b).

\begin{lem} \label{very-flat-open-immersion-base-change-context}
 Let $R$, $S$, and $T$ be adic topological rings, let $f\:R\rarrow S$
be a very flat tight continuous ring map, and let $g\:R\rarrow T$ be
a formal open immersion.
 Put $W=T\ot_RS$, and denote by $f'\:T\rarrow W$ and $g'\:S\rarrow W$
the induced ring maps.
 Endow the commutative ring $W$ with the tensor product topology, as
in Section~\ref{prelim-tensor-products-of-adic-topological-subsecn}.
 Then the tight continuous ring map~$f'$ is very flat, too.
\end{lem}

\begin{proof}
 The point is that if $\overline R\rarrow\overline S$ is a very flat
homomorphism of commutative rings and $\overline R\rarrow\overline T$
is an open immersion, then $\overline T\rarrow\overline T
\ot_{\overline R}\overline S$ is a very flat ring homomorphism, too.
 This is a special case of~\cite[Lemma~1.10.1]{Pcosh}.
\end{proof}

\subsection{Cotorsion contramodules}
\label{prelim-cotorsion-contramodules-subsecn}
 Once again, we start with a reminder of the case of an abstract
associative ring~$R$.
 A left $R$\+module $P$ is said to be \emph{cotorsion} (in
the sense of Enochs~\cite{En2}) if $\Ext^1_R(F,P)=0$ for all
flat left $R$\+modules~$P$.
 The full subcategory of cotorsion $R$\+modules is closed under
extensions, cokernels of monomorphisms, and infinite products
in $R\Modl$.

\begin{lem} \label{cotorsion-contramodules-lemma}
 Let $R$ be an adic topological ring and\/ $\fP$ be a quotseparated
contramodule $R$\+module.
 Then the following three conditions are equivalent:
\begin{enumerate}
\item the functor\/ $\Hom_R({-},\fP)$ takes short exact sequences of
flat quotseparated contramodule $R$\+modules to short exact
sequences of abelian groups/$R$\+mod\-ules/quotseparated
contramodule $R$\+modules;
\item one has\/ $\Ext^1_{R\Contra}(\fF,\fP)=0$ for all flat
quotseparated contramodule $R$\+mod\-ules\/~$\fF$;
\item one has\/ $\Ext^n_{R\Contra}(\fF,\fP)=0$ for all flat
quotseparated contramodule $R$\+mod\-ules\/ $\fF$ and all integers
$n\ge1$.
\end{enumerate}
\end{lem}

\begin{proof}
 This is an application of
Lemma~\ref{Hom-from-to-short-exact-sequences-lemma}(a) below.
 All the equivalences follow from the facts that there are enough flat
quotseparated contramodules in $R\Modl_\ctra^\qs$ and the full
subcategory of flat quotseparated contramodules is closed under kernels
of epimorphisms (see
Lemma~\ref{flat-quotseparated-contramodules-well-behaved}(a)).
\end{proof}

 Let $R$ be an adic topological ring.
 A quotseparated contramodule $R$\+module $\fP$ is said to be
\emph{cotorsion} if it satisfies the equivalent conditions of
Lemma~\ref{cotorsion-contramodules-lemma}.
 The full subcategory of cotorsion quotseparated contramodule
$R$\+modules is closed under extensions, cokernels of monomorphisms,
and infinite products in $R\Modl_\ctra^\qs$.

\begin{prop} \label{flat-cotorsion-pair-in-quotseparated}
 Let $R$ be an adic topological ring and\/ $\fM$ be a quotseparated
contramodule $R$\+module.  Then \par
\textup{(a)} there exists a short exact sequence of quotseparated
contramodule $R$\+modules\/ $0\rarrow\fP\rarrow\fF\rarrow\fM\rarrow0$
with a flat quotseparated contramodule $R$\+module\/ $\fF$ and
a cotorsion quotseparated contramodule $R$\+module\/~$\fP$; \par
\textup{(b)} there exists a short exact sequence of quotseparated
contramodule $R$\+modules\/ $0\rarrow\fM\rarrow\fP\rarrow\fF\rarrow0$
with a cotorsion quotseparated contramodule $R$\+module\/ $\fP$
and a flat quotseparated contramodule $R$\+module\/~$\fF$.
\end{prop}

\begin{proof}
 This is a special case of~\cite[Corollary~7.8]{PR}.
 For a Noetherian ring $R$, one can also refer
to~\cite[Corollary~D.2.8(a\+-b)]{Pcosh}.
\end{proof}

\begin{cor} \label{flat-contramodules-cor}
 Let $R$ be an adic topological ring and\/ $\fF$ be a quotseparated
contramodule $R$\+module.
 Then the following conditions are equivalent:
\begin{enumerate}
\item the functor\/ $\Hom_R(\fF,{-})$ takes short exact sequences of
cotorsion quotseparated contramodule $R$\+modules to short exact
sequences of abelian groups/$R$\+mod\-ules/quotseparated contramodule
$R$\+modules;
\item one has\/ $\Ext^1_{R\Contra}(\fF,\fP)=0$ for all
cotorsion quotseparated contramodule $R$\+modules\/~$\fP$;
\item one has\/ $\Ext^n_{R\Contra}(\fF,\fP)=0$ for all cotorsion
quotseparated contramodule $R$\+modules\/ $\fF$ and all integers
$n\ge1$;
\item $\fF$~is a flat quotseparated contramodule $R$\+module.
\end{enumerate}
\end{cor}

\begin{proof}
 The argument is similar to the proof of
Corollary~\ref{very-flat-contramodules-cor} and based on
Proposition~\ref{flat-cotorsion-pair-in-quotseparated} together with
Lemma~\ref{cotorsion-contramodules-lemma}(2\+-3).
 Lemmas~\ref{Hom-from-to-short-exact-sequences-lemma}(b)
and~\ref{direct-summand-lemma} below are relevant here.
 For a Noetherian ring $R$, one can also refer
to~\cite[Corollary~D.2.8(d)]{Pcosh}.
\end{proof}

 In the terminology of Section~\ref{appx-cotorsion-pairs-subsecn} below,
the results of Proposition~\ref{flat-cotorsion-pair-in-quotseparated}
and Corollary~\ref{flat-contramodules-cor}\,%
(2)\,$\Leftrightarrow$\,(3)\,$\Leftrightarrow$\,(4)
can be summarized by the following formulation.
 For any adic topological ring $R$, the pair of classes (flat
quotseparated contramodule $R$\+modules, cotorsion quotseparated
contramodule $R$\+modules) is a \emph{hereditary complete cotorsion
pair} in the abelian category of quotseparated contramodule
$R$\+modules $R\Contra=R\Modl_\ctra^\qs$.

 Parts~(b\+-c) of the next lemma form a generalization
of~\cite[Lemma~1.3.2(a\+-b)]{Pcosh}.

\begin{lem} \label{flat-cotorsion-injective-tensor-Hom}
 Let $R$ be an adic topological ring.  Then \par
\textup{(a)} for any two flat quotseparated contramodule $R$\+modules\/
$\fF$ and\/ $\fG$, the quotseparated contramodule $R$\+module\/
$\Lambda(\fF\ot_R\fG)$ is flat; \par
\textup{(b)} for any flat quotseparated contramodule $R$\+module\/ $\fF$
and any cotorsion quotseparated contramodule $R$\+module\/ $\fP$,
the quotseparated contramodule $R$\+module\/ $\Hom_R(\fF,\fP)$ is
cotorsion; \par
\textup{(c)} for any torsion $R$\+module\/ $\cM$ and any injective
torsion $R$\+module\/ $\cJ\in R\Tors^\inj$, the quotseparated
contramodule $R$\+module\/ $\Hom_R(\cM,\cJ)$ is cotorsion.
\end{lem}

\begin{proof}
 The argument for parts~(a\+-b) is similar to the proof of
Lemma~\ref{vfl-cta-tensor-Hom}.
 The proof of part~(b) is based on
Lemma~\ref{cotorsion-contramodules-lemma}(1).
 In part~(c), one can use
Lemma~\ref{cotorsion-contramodules-lemma}(1) together with
Lemma~\ref{flat-quotseparated-contramodules-well-behaved}(b).
\end{proof}

\begin{lem} \label{restriction-of-scalars-cotorsion}
 Let $f\:R\rarrow S$ be a continuous homomorphism of adic topological
rings.
 Then, for any cotorsion quotseparated contramodule $S$\+module\/ $\fP$,
the quotseparated contramodule $R$\+module $f_\sharp(\fP)$ is cotorsion.
\end{lem}

\begin{proof}
 This is a special case of~\cite[Corollary~6.3]{Pcs}.
 The assertion is provable using the criterion
of Lemma~\ref{cotorsion-contramodules-lemma}(2), by an argument based
on the $\Ext^1$\+adjunction lemma from~\cite[Lemma~1.7(e)]{Pal}
together with the results of
Corollaries~\ref{flat-contramodules-contraextension-of-scalars}
and~\ref{flat-contramodules-adjusted-to-contraextension}.
 In the case of Noetherian rings $R$ and $S$, one can also refer to
Lemma~\ref{cotorsion-contramodules-over-Noetherian}(2) below together
with the result of~\cite[Lemma~1.3.4(a)]{Pcosh}
or~\cite[Lemma~1.6(a)]{Pdomc}.
\end{proof}

 The following lemma is a cotorsion version of
Lemma~\ref{very-flat-reductions-Hom-exactness-lemma}.

\begin{lem} \label{flat-reductions-Hom-exactness-lemma}
 Let $R$ be an adic topological ring with an ideal of definition
$I\subset R$, and let $F$ be an $R$\+module such that
the $R/I^n$\+module $F/I^nF$ is flat for all $n\ge1$.
 Let\/ $0\rarrow\fP\rarrow\fM\rarrow\fN\rarrow0$ be a short exact
sequence of quotseparated contramodule $R$\+modules with a cotorsion
quotseparated contramodule $R$\+module\/~$\fP$.
 Then the short sequence of quotseparated contramodule $R$\+modules\/
$0\rarrow\Hom_R(F,\fP)\rarrow\Hom_R(F,\fM)\rarrow\Hom_R(F,\fN)\rarrow0$
is exact.
\end{lem}

\begin{proof}
 The argument is similar to the proof of
Lemma~\ref{very-flat-reductions-Hom-exactness-lemma}.
 One notices that $\boL_0\Lambda(F)\simeq\Lambda(F)$ is a flat
(quot)separated contramodule $R$\+module and uses the fact that
$\Ext^1_{R\Contra}(\Lambda(F),\fP)=0$
by Lemma~\ref{cotorsion-contramodules-lemma}(2).
\end{proof}

\begin{lem} \label{cotorsion-contramodules-over-Noetherian}
 Let $R$ be a Noetherian commutative ring endowed with an adic
topology, $\fR$ be the adic completion of $R$, and\/ $\fP$ be
a contramodule $R$\+module.
 Then the following conditions are equivalent:
\begin{enumerate}
\item $\fP$~is a cotorsion (quotseparated) contramodule $R$\+module;
\item $\fP$~is a cotorsion $R$\+module;
\item $\fP$~is a cotorsion\/ $\fR$\+module.
\end{enumerate}
\end{lem}

\begin{proof}
 Recall that a contramodule $R$\+module $\fF$ is flat if and only if
it is flat as an $R$\+module (or as an $\fR$\+module), by
Lemma~\ref{flat-contramodules-characterizations}.
 Furthermore, the inclusion of abelian categories $R\Contra\rarrow
R\Modl$ induces isomorphisms of all the Ext groups/modules,
by Proposition~\ref{Noetherian-contramodule-Ext-isomorphism-prop}.
 In view of these observations, the equivalence
(1)~$\Longleftrightarrow$~(2) is provided
by~\cite[Corollary~D.2.8(c)]{Pcosh}.
 The equivalence (1)~$\Longleftrightarrow$~(3) is a particular
case of (1)~$\Longleftrightarrow$~(2).
\end{proof}

\begin{lem} \label{flat-cotorsion-contramods-in-terms-of-reductions}
 Let $R$ be a Noetherian commutative ring endowed with an adic
topology and\/ $\fF$ be a flat contramodule $R$\+module.
 Then the following conditions are equivalent:
\begin{enumerate}
\item $\fF$~is a cotorsion (quotseparated) contramodule $R$\+module;
\item for every ideal $I\subset R$, the $R/I$\+module\/ $\fF/I\fF$ is cotorsion;
\item for every open ideal $I\subset R$, the $R/I$\+module\/
$\fF/I\fF$ is cotorsion;
\item for every ideal of definition $I\subset R$, the $R/I$\+module\/
$\fF/I\fF$ is cotorsion;
\item there exists an ideal of definition $I\subset R$ for which\/
$\fF/I\fF$ is a cotorsion $R/I$\+module.
\end{enumerate}
\end{lem} 

\begin{proof}
 The implications (2)~$\Longrightarrow$~(1) and
(2)~$\Longrightarrow$~(3) $\Longrightarrow$~(4)
$\Longrightarrow$~(5) are obvious.

 (1)~$\Longrightarrow$~(2) By
Lemmas~\ref{flat-contramodules-characterizations}
and~\ref{cotorsion-contramodules-over-Noetherian}, any flat
cotorsion contramodule $R$\+module is a flat cotorsion $R$\+module.
 So the result of~\cite[Lemma~1.7.9(a)]{Pcosh} is applicable.

 (1)~$\Longrightarrow$~(3) This is a special case
of the ``only if'' assertion in~\cite[Corollary~E.5.5]{Pcosh}.

 (5)~$\Longrightarrow$~(4) This is~\cite[Lemma~1.7.10(d)]{Pcosh}.

 (4)~$\Longrightarrow$~(1) This is a special case
of~\cite[Corollary~E.5.2]{Pcosh}, which is applicable in view
of~\cite[Proposition~1.5]{Pdc}.
\end{proof}

 By the \emph{adic Krull dimension} of an adic topological ring $R$
we mean the Krull dimension of the quotient ring $R/I$, for any
ideal of definition $I\subset R$.
 Clearly, the Krull dimension of $R/I$ does not depend on the choice
of an ideal of definition $I$ of the given adic topology on~$R$.
 The definition of an \emph{adically Noetherian} adic topological ring
was given in Section~\ref{adic-rings-subsecn}.

 From some point on, we will need to use the language of
(\emph{co})\emph{resolving subcategories} and
(\emph{co})\emph{resolution dimensions} from~\cite[Section~2]{Sto}
and~\cite[Sections~A.2, A.3, and A.5]{Pcosh}.
 The reader can find the definitions and basic results collected
in Section~\ref{appx-resolving-subsecn} below.
 In particular, the \emph{cotorsion dimension} of a left module $M$
over an associative ring~$R$ \,\cite[Section~1.5]{Pcosh} is defined
as the minimal length~$d$ of a coresolution $0\rarrow M\rarrow
C^0\rarrow C^1\rarrow\dotsb\rarrow C^d\rarrow0$ of the $R$\+module
$M$ by cotorsion left $R$\+modules $C^i$, \,$0\le i\le d$,
in the abelian category $R\Modl$.

 Similarly, given an adic topological ring $R$, the \emph{cotorsion
dimension} of a quotseparated contramodule $R$\+module $\fM$ is
defined as the minimal length~$d$ of a coresolution $0\rarrow\fM\rarrow
\fC^0\rarrow\fC^1\rarrow\dotsb\rarrow\fC^d\rarrow0$ of the quotseparated
contramodule $R$\+module $\fM$ by cotorsion quotseparated contramodule
$R$\+modules $\fC^i$, \,$0\le i\le d$, in the abelian category
$R\Contra$.
 In view of
Corollary~\ref{intersection-co-resolution-dimension-cor}(b) below,
it follows from Lemma~\ref{cotorsion-contramodules-over-Noetherian}
that the cotorsion dimension of any contramodule $R$\+module $\fM$
coincides with its cotorsion dimension as an abstract $R$\+module
whenever the ring $R$ is Noetherian as an abstract ring.
 
\begin{prop} \label{adic-raynaud-gruson-prop}
 Let $R$ be an adically Noetherian adic topological ring of finite
adic Krull dimension~$D$.
 Then \par
\textup{(a)} the projective dimensions of flat quotseparated
contramodule $R$\+modules in the abelian category $R\Contra$
do not exceed~$D$; \par
\textup{(b)} the cotorsion dimensions of quotseparated contramodule
$R$\+modules do not exceed~$D$.
\end{prop}

\begin{proof}
 In other words, part~(b) claims that the coresolution dimension of any
quotseparated contramodule $R$\+module with respect to the coresolving
subcategory of cotorsion quotseparated contramodule $R$\+modules
$R\Contra^\cot\subset R\Contra$ does not exceed~$D$.
 In view of Proposition~\ref{flat-cotorsion-pair-in-quotseparated},
the assertions~(a) and~(b) are equivalent to each other
by~\cite[Lemma~B.1.9]{Pcosh} or~\cite[Lemma~6.4]{Pfltp}.
 Part~(a) follows from~\cite[Lemma~E.2.2(c)]{Pcosh} together with
the fact that the projective dimensions of flat $R/I^n$\+modules
(as objects of $R/I^n\Modl$) do not exceed~$D$, which is a celebrated
theorem of Raynaud and Gruson~\cite[Corollaire~II.3.2.7]{RG}.

 The setting of~\cite[Lemma~E.2.2(c)]{Pcosh} is way more general than
needed for the purposes of the present proposition.
 Lemmas~\ref{projective-contramodules-characterizations}\,%
(1)\,$\Leftrightarrow$\,(6)
and~\ref{flat-quotseparated-contramodules-well-behaved} are relevant
here.
 Basically, let $\fF$ be a flat quotseparated contramodule $R$\+module
and $\fP_\bu$ be its projective resolution in $R\Contra$.
 Consider the truncated fragment $0\rarrow\fG\rarrow\fP_{D-1}\rarrow
\fP_{D-2}\rarrow\dotsb\rarrow\fP_0\rarrow\fF\rarrow0$ of
the resolution~$\fP_\bu$.
 By Lemma~\ref{flat-quotseparated-contramodules-well-behaved}(a),
the quotseparated contramodule $R$\+module $\fG$ is flat (and so
are all the other quotsepatated contramodule $R$\+modules of cycles
in the complex~$\fP_\bu$).
 By Lemma~\ref{flat-quotseparated-contramodules-well-behaved}(b),
the complex of $R/I$\+modules $0\rarrow\fG/I\fG\rarrow
\fP_{D-1}/I\fP_{D-1}\rarrow\dotsb\rarrow\fP_0/I\fP_0\rarrow\fF/I\fF
\rarrow0$ is exact.
 As the $R/I$\+modules $\fP_i/I\fP_i$ are projective and
the $R/I$\+module $\fF/I\fF$ is flat, it follows by virtue of
the Raynaud--Gruson theorem that the $R/I$\+module $\fG/I\fG$ is
projective.
 It remains to apply
Lemma~\ref{projective-contramodules-characterizations}\,%
(6)\,$\Rightarrow$\,(1) (which does not require the Noetherianity
assumption on $R$, as per the proof) in order to conclude that
the quotseparated contramodule $R$\+module $\fG$ is projective.
\end{proof}

\subsection{Locality of flatness and very flatness of morphisms}
\label{prelim-morphism-very-flatness-locality-subsecn}
 This section is a continuation of
Sections~\ref{prelim-tight-taut-subsecn}\+-%
\ref{prelim-locality-tautness-tightness-subsecn}.

 The following lemma is to be compared with
Lemma~\ref{flat-homomorphism-of-adic-rings}.

\begin{lem} \label{taut-flat-homomorphism-of-adic-rings}
 Let $f\:R\rarrow S$ be a taut continuous map of adic topological rings.
 Then the following two conditions are equivalent:
\begin{enumerate}
\item there exists an ideal of definition $I\subset R$ such that,
denoting by $J_n$ the closure of the ideal $Sf(I^n)\subset S$,
the $R/I^n$\+module $S/J_n$ is flat for every $n\ge1$;
\item for every open ideal $I\subset R$, denoting by $J$ the closure of
the ideal $Sf(I)\subset S$, the $R/I$\+module $S/J$ is flat.
\end{enumerate}
\end{lem}

\begin{proof}
 (2)~$\Longrightarrow$~(1) Obvious.

 (1)~$\Longrightarrow$~(2) Let $I'\subset I$ be an open ideal
such that, denoting by $J'$ the closure of the ideal $Sf(I')\subset S$,
the $R/I'$\+module $S/J'$ is flat.
 We need to show that the $R/I$\+module $S/J$ is flat whenever
$I'\subset I$.
 Indeed, $J'$ is an open ideal in $S$ by
Lemma~\ref{taut-ring-map-lemma}(5).
 Since $J'\subset J$, it follows that $J=Sf(I)+J'$.
 Thus $S/J\simeq R/I\ot_{R/I'}S/J'$, and the assertion follows.
\end{proof}

 We will say that a taut continuous homomorphism of adic topological
rings $f\:R\rarrow S$ is \emph{taut-flat} if it satisfies the equivalent
conditions of Lemma~\ref{taut-flat-homomorphism-of-adic-rings}.
 Notice that a tight continuous map of adic topological rings is
taut-flat if and only if it is flat.
 Using the criterion of
Lemma~\ref{taut-flat-homomorphism-of-adic-rings}(2), one can easily
check that the composition of any two taut-flat continuous ring maps of
adic topological rings is taut-flat.

 The next lemma is to be compared with
Lemma~\ref{very-flat-homomorphism-of-adic-rings}.

\begin{lem} \label{taut-very-flat-homomorphism-of-adic-rings}
 Let $f\:R\rarrow S$ be a taut continuous map of adic topological rings.
 Then the following three conditions are equivalent:
\begin{enumerate}
\item for every open ideal $I\subset R$, denoting by $J$ the closure of
the ideal $Sf(I)\subset S$, the induced ring map $\bar f\:R/I\rarrow
S/J$ is very flat;
\item there exists an ideal of definition $I\subset R$ such that,
denoting by $J_n$ the closure of the ideal $Sf(I^n)\subset S$,
the induced ring map $\bar f\:R/I^n\rarrow S/J_n$ is very flat
for every $n\ge1$;
\item $f$~is taut-flat and there exists an ideal of definition
$I\subset R$ such that, denoting by $J$ the closure of the ideal
$Sf(I)\subset S$, the induced ring map $\bar f\:R/I\rarrow S/J$
is very flat.
\end{enumerate}
\end{lem}

\begin{proof}
 (1)~$\Longrightarrow$~(2) $\Longrightarrow$~(3) Obvious.

 (2)~$\Longrightarrow$~(1) Let $I'\subset I$ be an open ideal
such that, denoting by $J'$ the closure of the ideal $Sf(I')\subset S$,
the ring map $R/I'\rarrow S/J'$ is very flat.
 We need to show that the ring map $R/I\rarrow S/J$ is very flat
whenever $I'\subset I$.
 Following the proof of
Lemma~\ref{taut-flat-homomorphism-of-adic-rings}, we have
$S/J\simeq R/I\ot_{R/I'}S/J'$.
 It remains to point out that the ring maps $S/J'\rarrow S/J$ is
surjective (so every localization of $S/J$ at an element arises from
a localization of $S/J'$ at an element), and refer
to~\cite[Lemma~1.2.2(b)]{Pcosh}.

 (3)~$\Longrightarrow$~(2) We have $S/J\simeq R/I\ot_{R/I^n}S/J_n$,
so~\cite[Lemma~1.7.10(b)]{Pcosh} is applicable.
\end{proof}

 We will say that a taut continuous homomorphism of adic topological
rings $f\:R\rarrow S$ is \emph{taut-very flat} if it satisfies
the equivalent conditions of
Lemma~\ref{taut-very-flat-homomorphism-of-adic-rings}.
 Notice that a tight continuous map of adic topological rings is
taut-very flat if and only if it is very flat.

 Using the criterion of
Lemma~\ref{taut-very-flat-homomorphism-of-adic-rings}(1), one can
easily check that the composition of any two taut-very flat continuous
ring maps of adic topological rings is taut-very flat.
 The fact that the composition of any two very flat homomorphisms of
commutative rings is very flat (see
Section~\ref{prelim-veryflat-and-contraadjusted-subsecn}) needs to
be used here.

 The following lemma is a generalization of
Lemmas~\ref{induced-map-of-complete-rings-lemma}(b)
and~\ref{induced-map-of-complete-rings-very-flat-lemma}.

\begin{lem} \label{tight-taut-completion-very-flatness}
 Let $f\:R\rarrow S$ be a continuous map of adic topological rings.
 In this context: \par
\textup{(a)} the map~$f$ is taut-flat if and only if the map\/
$\Lambda(f)$ is flat; \par
\textup{(b)} the map~$f$ is taut-very flat if and only if the map\/
$\Lambda(f)$ is very flat.
\end{lem}

\begin{proof}
 First of all, the map~$f$ is taut if and only if the map $\Lambda(f)$
is tight, by Corollary~\ref{completion-preserves-reflects-tautness-cor}
and Proposition~\ref{complete-tightness=tautness-prop}.
 Assume that this is the case.
 Put $\fR=\Lambda(R)$ and $\fS=\Lambda(S)$, and denote by
$\ff=\Lambda(f)$ the induced ring map $\fR\rarrow\fS$.
 Let $I\subset R$ be an open ideal and $\fI\subset\fR$ be the related
open ideal, as per the discussion in
Section~\ref{prelim-adic-completions-subsecn}.
 Let $J\subset S$ be the closure of the ideal $Sf(I)\subset S$.
 Then $\fJ=\fS\lambda_S(J)$ is the open ideal in $\fS$ corresponding
to~$J$.
 Clearly, $\fJ$ is the closure of the ideal $\fS\ff(\fI)\subset\fS$.
 Since the map~$\ff$ is tight, it follows that $\fJ=\fS\ff(\fI)$.
 Now we have $R/I\simeq\fR/\fI$ and $S/J\simeq\fS/\fJ$.
 Both the assertions of parts~(a) and~(b) follow.
\end{proof}

\begin{lem} \label{locality-of-taut-very-flatness-in-total-space}
 Let $f\:R\rarrow S$ be a continuous homomorphism of adic topological
rings.
 Let $g_\alpha\:S\rarrow T_\alpha$ be a formal open covering.
 In this context: \par
\textup{(a)} the map~$f$ is taut-flat if and only if the composition
$g_\alpha f\:R\rarrow T_\alpha$ is taut-flat for every index\/~$\alpha$;
\par
\textup{(b)} the map~$f$ is taut-very flat if and only if
the composition $g_\alpha f\:R\rarrow T_\alpha$ is taut-very flat for
every index\/~$\alpha$.
\end{lem}

\begin{proof}
 The ``only if'' assertions hold, since the maps~$g_\alpha$ are
taut-very flat, the compositions of taut-flat morphisms are taut-flat,
and the compositions of taut-very flat morphisms are taut-very flat.
 To prove the ``if'', we use the notation from the proof of
Proposition~\ref{locality-of-tautness-in-total-space-prop}.
 According to that proof, we have $J_n=Sf(I^n)+J_{n+1}\subset S$, hence
$J_n$ is the closure of the ideal $Sf(I^n)\subset S$ by the proof of
Lemma~\ref{continuous-taut-characterization}\,(2)\,$\Rightarrow$\,(1).
 We also know that the ring homomorphisms~$\bar g_\alpha$ form an open
covering.

 Part~(a): we know that $\overline T_\alpha$ are flat
$\overline R$\+modules, and need to prove that $\overline S$ is
a flat $\overline R$\+module.
 This is a standard fact from commutative algebra; see, e.~g.,
\cite[Lemma Tag~01U5]{SP} or~\cite[Lemma~1.8.2]{Pcosh}.

 Part~(b): we know that the ring homomorphisms~$\bar g_\alpha\bar f$
are very flat, and need to prove that the ring homomorphism~$\bar f$
is very flat.
 This result can be found in~\cite[Lemma~1.8.6]{Pcosh}.
\end{proof}

\begin{cor} \label{locality-of-very-flatness-in-total-space-cor}
 Let $f\:R\rarrow S$ be a continuous homomorphism of adic topological
rings.
 Let $g_\alpha\:S\rarrow T_\alpha$ be a formal open covering.
 In this setting: \par
\textup{(a)} the map\/ $\Lambda(f)$ is flat (as a map of adic
topological rings) if and only if the map\/ $\Lambda(g_\alpha f)$
is flat for every index\/~$\alpha$; \par
\textup{(a)} the map\/ $\Lambda(f)$ is very flat (as a map of adic
topological rings) if and only if the map\/ $\Lambda(g_\alpha f)$
is very flat for every index\/~$\alpha$.
\end{cor}

\begin{proof}
 Follows from Lemmas~\ref{tight-taut-completion-very-flatness}
and~\ref{locality-of-taut-very-flatness-in-total-space}.
\end{proof}

\begin{lem} \label{taut-very-flatness-base-change-context}
 Let $R$, $S$, and $T$ be adic topological rings, and let
$f\:R\rarrow S$ and $g\:R\rarrow T$ be a continuous ring maps.
 Put $W=T\ot_RS$, and denote by $f'\:T\rarrow W$ and $g'\:S\rarrow W$
the induced ring maps.
 Endow the commutative ring $W$ with the tensor product topology, as
in Section~\ref{prelim-tensor-products-of-adic-topological-subsecn}.
 In this setting: \par
\textup{(a)} if the ring map~$f$ is taut-flat, then the continuous
ring map~$f'$ is also taut-flat; \par
\textup{(b)} if the ring map~$f$ is taut-very flat and the ring map~$g$
is a formal open immersion, then the continuous ring map~$f'$ is also
taut-very flat.
\end{lem}

\begin{proof}
 Part~(a): by Lemma~\ref{tight-taut-completion-very-flatness}(a),
the continuous ring map $\Lambda(f)\:\Lambda(R)\rarrow\Lambda(S)$ is
flat.
 Hence, by Lemma~\ref{tight-flat-adic-ring-base-change-context}(b),
the ring map $\Lambda(T)\rarrow\Lambda(T)\ot_{\Lambda(R)}\Lambda(S)$
is flat, and consequently, taut-flat.
 By Lemma~\ref{completed-tensor-product-of-completions-lemma},
we have $\Lambda(W)=\Lambda(\Lambda(T)\ot_{\Lambda(R)}\Lambda(S))$.
 Using Lemma~\ref{tight-taut-completion-very-flatness}(a) again,
we conclude that the ring map $\Lambda(f')\:\Lambda(T)\rarrow\Lambda(W)$
is flat.
 By the same lemma, it follows that the ring map~$f'$ is taut-flat.
 
 The proof of part~(b) is similar and based on
Lemmas~\ref{tight-taut-completion-very-flatness}(b)
and~\ref{very-flat-open-immersion-base-change-context}.
\end{proof}

\begin{lem} \label{locality-of-taut-very-flatness-in-the-base}
 Let $f\:R\rarrow S$ be a continuous homomorphism of adic topological
rings, and let $g_\alpha\:R\rarrow T_\alpha$ be a formal open covering.
 Put $W_\alpha=T_\alpha\ot_RS$, and endow the rings $W_\alpha$ with
the tensor product topologies.
 Denote by $f'_\alpha\:T_\alpha\rarrow W_\alpha$ and $g'_\alpha\:
S\rarrow W_\alpha$ the induced ring maps.
 In this context: \par
\textup{(a)} the map~$f$ is taut-flat if and only if
the map~$f'_\alpha$ is taut-flat for every index\/~$\alpha$; \par
\textup{(b)} the map~$f$ is taut-very flat if and only if
the map~$f'_\alpha$ is taut-very flat for every index\/~$\alpha$.
\end{lem}

\begin{proof}
 The ``only if'' assertions hold by
Lemma~\ref{taut-very-flatness-base-change-context}.
 To prove the ``if'', we use the notation from the proof of
Proposition~\ref{locality-of-tautness-in-the-base-prop}.
 According to that proof, we have $J_n=Sf(I^n)+J_{n+1}\subset S$, hence
$J_n$ is the closure of the ideal $Sf(I^n)\subset S$ by the proof of
Lemma~\ref{continuous-taut-characterization}\,(2)\,$\Rightarrow$\,(1).

 Part~(a): we know that $\overline W_\alpha$ are flat
$\overline T_\alpha$\+modules, and need to prove that $\overline S$ is
a flat $\overline R$\+module.
 This is a standard fact from commutative algebra; see, e.~g.,
\cite[Lemma Tag~01U5]{SP} or~\cite[Lemma~1.8.3]{Pcosh}.

 Part~(b): we know that the ring homomorphisms~$\bar f'_\alpha$
are very flat, and need to prove that the ring homomorphism~$\bar f$
is very flat.
 This follows from~\cite[Lemma~1.2.6(a)]{Pcosh}.
\end{proof}

\begin{cor} \label{locality-of-very-flatness-in-the-base-cor}
 Let $f\:R\rarrow S$ be a continuous homomorphism of adic topological
rings, and let $g_\alpha\:R\rarrow T_\alpha$ be a formal open covering.
 Put $W_\alpha=T_\alpha\ot_RS$, and endow the rings $W_\alpha$ with
the tensor product topologies.
 Denote by $f'_\alpha\:T_\alpha\rarrow W_\alpha$ and $g'_\alpha\:
S\rarrow W_\alpha$ the induced ring maps.
 Then the map\/ $\Lambda(f)$ is flat (as a map of adic topological
rings) if and only if the map\/ $\Lambda(f'_\alpha)$ is flat for every
index\/~$\alpha$.
\end{cor}

\begin{proof}
 Follows from Lemmas~\ref{tight-taut-completion-very-flatness}
and~\ref{locality-of-taut-very-flatness-in-the-base}.
\end{proof}

\subsection{The coextension of scalars is quotseparated}
\label{prelim-coextension-of-scalars-quotseparated-subsecn}
 The aim of this section is to formulate and prove quotseparated
versions of Lemmas~\ref{tight-reflects-torsion-contra-lemma}(b)
and~\ref{tight-co-extension-of-scalars}(b).
 Our exposition is based on the ideas and techniques from the new
preprint~\cite{Pcs}.

 For any adic topological ring $R$ with an ideal of definition
$I\subset R$ and any $R$\+module $M$, we will denote
the kernel of the $I$\+adic completion map $\lambda_{I,M}\:M\rarrow
\Lambda_I(M)$ by $\Omega_R(M)=\bigcap_{n\ge1}I^nM\subset M$.

\begin{lem} \label{cont-tight-reflects-separated-and-complete}
 Let $R$ and $S$ be adic topological rings, and let $f\:R\rarrow S$
be a tight continuous ring map.
 Then an $S$\+module\/ $\fQ$ is separated and complete if and only if
the $R$\+module $f_*(\fQ)$ is separated and complete.
 Moreover, for any $S$\+module $N$, one has\/ $\Omega_R(N)=\Omega_S(N)
\subset N$.
 Denoting by $I\subset R$ and $J\subset S$ any ideals of definition
in $R$ and $S$, we also have a natural isomorphism of $R$\+modules\/
$\Lambda_I(N)\simeq\Lambda_J(N)$ forming a commutative triangular
diagram with the completion maps\/ $\lambda_{I,N}\:N\rarrow\Lambda_I(N)$
and\/ $\lambda_{J,N}\:N\rarrow\Lambda_J(N)$.
\end{lem}

\begin{proof}
 All the assertions follow immediately from
Lemma~\ref{continuous-and-tight-ring-map-lemma}(4).
\end{proof}

 Let $f\:R\rarrow S$ be a tight continuous map of adic topological
rings.
 With Lemma~\ref{cont-tight-reflects-separated-and-complete} in
mind, we will write simply $\Omega(N)$ instead of $\Omega_R(N)=
\Omega_S(N)$, and $\Lambda(N)$ instead of $\Lambda_I(N)\simeq
\Lambda_J(N)$, for any $S$\+module~$N$.
 The completion map $N\rarrow\Lambda(N)$ will be denoted, as usual,
by~$\lambda_N$.

\begin{thm} \label{flat-reflects-quotseparated-contramodules-thm}
 Let $R$ and $S$ be adic topological rings, and let $f\:R\rarrow S$
be a flat tight continuous ring map.
 Then an $S$\+module\/ $\fQ$ is a quotseparated contramodule if and
only if the $R$\+module $f_*(\fQ)$ is a quotseparated contramodule.
\end{thm}

\begin{proof}
 This is our version of~\cite[Theorem~9.5]{Pcs}.
 The ``only if'' assertion holds by
Lemma~\ref{torsion-contra-restriction-of-scalars}(d); the nontrivial
implication is the ``if''.

 So let $\fQ$ be an $S$\+module whose underlying $R$\+module is
a quotseparated contramodule.
 Then $\fQ$ is a contramodule $S$\+module by
Lemma~\ref{tight-reflects-torsion-contra-lemma}(b).
 We need to pove that $\fQ$ is a quotseparated contramodule
$S$\+module.

 Any contramodule $R$\+module is complete, as mentioned in
Section~\ref{prelim-contramodules-subsecn}.
 Hence the completion map $\lambda_\fQ\:\fQ\rarrow\Lambda(\fQ)$ is
surjective, and we have a short exact sequence of $S$\+modules
$0\rarrow\Omega(\fQ)\rarrow\fQ\rarrow\Lambda(\fQ)\rarrow0$.

 Pick a flat quotseparated contramodule $S$\+module $\fF$ together
with a surjective $S$\+module map $\fF\rarrow\Lambda(\fQ)$.
 For example, it suffices to take a free $S$\+module $F$ together
with a surjective $S$\+module morphism $F\rarrow\Lambda(\fQ)$, and
put $\fF=\Lambda(F)$.
 Consider the pullback diagram in the category of $S$\+modules
\begin{equation} \label{Q-Q1-pullback-diagram}
\begin{gathered}
 \xymatrix{
  0 \ar[r] & \Omega(\fQ) \ar@{>->}[r] &
  \fQ \ar@{->>}[r] & \Lambda(\fQ) \ar[r] & 0 \\
  0 \ar[r] & \Omega(\fQ) \ar@{>->}[r] \ar@{=}[u] &
  \fQ_1 \ar@{->>}[u] \ar@{->>}[r] & \fF \ar@{->>}[u] \ar[r] & 0
 }
\end{gathered}
\end{equation}
 Here the $\fQ_1$ is the pullback of the pair of surjective maps
$\lambda_\fQ\:\fQ\rarrow\Lambda(\fQ)$ and $\fF\rarrow\Lambda(\fQ)$.
 Both the top and bottom rows in~\eqref{Q-Q1-pullback-diagram}
are short exact sequences.
 Since the underlying $R$\+modules of $\fQ$, \,$\Lambda(\fQ)$,
and $\fF$ are quotseparated contramodule $R$\+modules, it follows
that the underlying $R$\+module of $\fQ_1$ is a quotseparated
contramodule $R$\+module, too.

 By Lemma~\ref{flat-map-of-adic-rings-direct-image-lemma}(d),
the underlying $R$\+module of $\fF$ is a flat quotseparated
contramodule $R$\+module.
 By Lemma~\ref{flat-quotseparated-contramodules-separated},
\,$\fF$ is a separated contramodule $R$\+module (as well as
a separated contramodule $S$\+module).
 Applying the functor $\Lambda$ and the natural
transformation~$\lambda$ to the surjective $S$\+module map
$\fQ_1\rarrow\fF$, we obtain a factorization of that map as
the composition of surjective $S$\+module maps
$$
 \xymatrix{
  \fQ_1 \ar@{->>}[r] & \Lambda(\fQ_1) \ar@{->>}[r] & \fF.
 }
$$
 It follows, in view of~\eqref{Q-Q1-pullback-diagram}, that
the surjective $S$\+module map $\fQ_1\rarrow\fQ$ induces
an injective $S$\+module map $\Omega(\fQ_1)\rarrow\Omega(\fQ)$.
 
 Let us prove that the map $\Omega(\fQ_1)\rarrow\Omega(\fQ)$
\emph{cannot} be an isomorphism if $\Omega(\fQ)\ne0$.
 Indeed, if the map $\Omega(\fQ_1)\rarrow\Omega(\fQ)$ is
an isomorphism, then $\Lambda(\fQ_1)\simeq\fF$.
 So $\Lambda(\fQ_1)$ is a flat (quot)separated contramodule $R$\+module.
 Since $\fQ_1$ is a quotseparated contramodule $R$\+module,
Corollary~\ref{flat-completion-of-quotseparated-implies-separated}
is applicable, and it follows that $\Omega(\fQ_1)=0$.
 Hence $\Omega(\fQ)=0$.

 Now we proceed by transfinite induction, iterating the construction
above at the successor steps, and passing to the projective limit at
the limit steps.
 Let $\aleph$ be a cardinal greater than the cardinality of~$\fQ$.
 Viewing $\aleph$ as an ordinal, we will construct a projective system
of $S$\+modules $\fQ_i$, \,$0\le i<\aleph$, and surjective
$S$\+module morphisms $\fQ_j\rarrow\fQ_i$ for all ordinals
$0\le i\le j<\aleph$.
 The underlying $R$\+module of $\fQ_i$ will be a quotseparated
contramodule $R$\+module for every $0\le i<\aleph$.
 The induced maps $\Omega(\fQ_j)\rarrow\fQ_i$ will be injective
for all $0\le i\le j<\aleph$, and the map $\Omega(\fQ_{i+1})
\rarrow\Omega(\fQ_i)$ will be a proper embedding (i.~e., not
an isomorphism) whenever $\Omega(\fQ_i)\ne0$.

 Put $\fQ_0=\fQ$.
 For every successor ordinal $j=i+1<\aleph$, pick a flat quotseparated
contramodule $S$\+module $\fF_i$ together with a surjective
$S$\+module map $\fF_i\rarrow\Lambda(\fQ_i)$.
 Set $\fQ_j$ to be the pullback of the two surjective $S$\+module
maps $\fQ_i\rarrow\Lambda(\fQ_i)$ and $\fF_i\rarrow\Lambda(\fQ_i)$.
 For every limit ordinal $j<\aleph$, put $\fQ_j=\varprojlim_{i<j}\fQ_i$.
 Clearly, it follows that the maps $\fQ_j\rarrow\fQ_i$ are surjective
for all $0\le i\le j<\aleph$.

 The only property we still need to check is that, for a limit
ordinal~$j$, the map $\Omega(\fQ_j)\rarrow\Omega(\fQ_i)$ is injective
for all $i<j$.
 It suffices to show that the map $\Omega(\fQ_j)\rarrow
\varprojlim_{i<j}\Omega(\fQ_i)$ is injective.
 Indeed, $\Omega(\fQ_j)$ is a submodule of $\fQ_j$, while
$\varprojlim_{i<j}\Omega(\fQ_i)$ is a submodule of
$\varprojlim_{i<j}\fQ_i$.
 Commutativity of the square diagram $\Omega(\fQ_j)\rarrow
\varprojlim_{i<j}\Omega(\fQ_i)\rarrow\varprojlim_{i<j}\fQ_i$,
\ $\Omega(\fQ_j)\rarrow\fQ_j\simeq\varprojlim_{i<j}\fQ_i$
implies the desired injectivity of the map $\Omega(\fQ_j)\rarrow
\varprojlim_{i<j}\Omega(\fQ_i)$.

 After the construction is finished, we observe that a set $\fQ$
of the cardinality less than~$\aleph$ cannot have a descending
chain of subsets of length~$\aleph$.
 Hence there exists an ordinal $k<\aleph$ for which the map
$\Omega(\fQ_{k+1})\rarrow\Omega(\fQ_k)$ is an isomorphism.
 As we have shown above, it follows that $\Omega(\fQ_k)=0$.

 Thus $\fQ_k$ is an $S$\+module whose underlying $R$\+module is
separated and complete.
 By Lemma~\ref{cont-tight-reflects-separated-and-complete},
it follows that $\fQ_k$ is separated and complete as an $S$\+module.
 Thus the contramodule $S$\+module $\fQ$ is a quotient $S$\+module
of a separated and complete $S$\+module $\fQ_k$, so $\fQ$ is
a quotseparated contramodule $S$\+module.
\end{proof}

\begin{rem} \label{nonflat-not-reflects-quotseparated}
 Notice that the assertion of
Theorem~\ref{flat-reflects-quotseparated-contramodules-thm} with
the flatness condition on~$f$ replaced by the tightness condition
(i.~e., the obvious version of
Lemma~\ref{tight-reflects-torsion-contra-lemma}(b) for quotseparated
contramodules) is \emph{not} true.
 Indeed, let $S$ be an adic topological ring with an ideal of
definition $J$ such that there exists a nonquotseparated contramodule
$S$\+module~$\fQ$ (see Section~\ref{prelim-contramodules-subsecn}).
 Let $s_1$,~\dots, $s_m\in J$ be a finite set of generators of
the ideal~$J$.
 Consider the ring of polynomials $R=\boZ[x_1,\dotsc,x_m]$ in
$m$~variables over the ring of integers $\boZ$, and let
$f\:R\rarrow S$ be the ring homomorphism taking~$x_j$ to~$s_j$
for all $1\le j\le m$.
 Let $I=(x_1,\dotsc,x_m)\subset R$ be the ideal generated by
the elements $x_1$,~\dots,~$x_m$.
 Endow $R$ with the $I$\+adic topology.
 Then $f\:R\rarrow S$ is a tight continuous ring map by
Lemma~\ref{continuous-and-tight-ring-map-lemma}(4),
and $f_*(\fQ)$ is a contramodule $R$\+module by
Lemma~\ref{torsion-contra-restriction-of-scalars}(b).
 Since the ring $R$ is Noetherian, all contramodule $R$\+modules
are quotseparated; so $f_*(\fQ)$ is a quotseparated contramodule
$R$\+module.
 But $\fQ$ is \emph{not} a quotseparated contramodule $S$\+module.

 On the other hand, the conclusion of 
Theorem~\ref{flat-reflects-quotseparated-contramodules-thm} certainly
does \emph{not} imply flatness of the map~$f$, i.~e., the assertion
of the theorem holds true for \emph{some} nonflat (but still tight)
ring maps $f\:R\rarrow S$.
 For example, if both the rings $R$ and $S$ are Noetherian and
the map~$f$ is tight and continuous, then the assertion of the theorem
holds by Lemma~\ref{torsion-contra-restriction-of-scalars}(b) or~(d)
and Lemma~\ref{tight-reflects-torsion-contra-lemma}(b).
 Alternatively, the assertion of the theorem trivially holds for
any two discrete (i.~e., $(0)$\+adic) rings $R$ and $S$, and any
ring homomorphism~$f$ (since any module over a discrete commutative
ring is a separated, hence quotseparated, contramodule).
\end{rem}

\begin{lem} \label{separated-coextension-of-scalars-prop}
 Let $R$ and $S$ be adic topological rings, and let $f\:R\rarrow S$
be a tight continuous ring map.
 Then, for any separated contramodule $R$\+module\/ $\fP$,
the $S$\+module $f^!(\fP)=\Hom_R(S,\fP)$ is also a separated
contramodule.
\end{lem}

\begin{proof}[First proof]
 By Lemma~\ref{Hom-is-a-quotseparated-contramodule}(b),
the underlying $R$\+module of $\Hom_R(S,\fP)$ is separated and
complete.
 By Lemma~\ref{cont-tight-reflects-separated-and-complete}, it
follows that the $S$\+module $\Hom_R(S,\fP)$ is separated and
complete as well.
\end{proof}

\begin{proof}[Second proof]
 The following direct computational argument is a special case
of~\cite[proof of Proposition~10.3]{Pcs}.
 Let $I\subset R$ be an ideal of definition in~$R$.
 By Lemma~\ref{continuous-and-tight-ring-map-lemma}(5), the ideal
$Sf(I)\subset S$ is an ideal of definition in~$S$.
 We have $\fP\simeq\varprojlim_{n\ge1}\fP/I^n\fP$.
 Applying the functor $\Hom_R(S,{-})$, we obtain an $S$\+module
isomorphism $\Hom_R(S,\fP)\simeq\varprojlim_{n\ge1}
\Hom_R(S,\fP/I^n\fP)$.
 Now we compute
$$
 \Hom_R(S,\fP/I^n\fP)\simeq\Hom_{R/I^n}(S/Sf(I^n),\>\fP/I^n\fP)=
 \Hom_{R/I^n}(S/J^n,\fP/I^n\fP).
$$
 Finally, $\Hom_{R/I^n}(S/J^n,\fP/I^n\fP)$ is an $S/J^n$\+module.
 So the $S$\+module $\Hom_R(S,{-})$ is a projective limit of
$S/J^n$\+modules.
 It remains to point out that all $S/J^n$\+modules are separated
contramodule $S$\+modules, and the full subcategory of separated
and the full subcategory of separated contramodule $S$\+modules
is closed under all limits in $S\Modl$.
 The latter assertion holds since the full subcategory of contramodule
$S$\+modules is closed under limits in $S\Modl$, while the full
subcategory of separated $S$\+modules is closed under submodules
and infinite products.
\end{proof}

\begin{prop} \label{flat-quotseparated-coextension-of-scalars-prop}
 Let $R$ and $S$ be adic topological rings, and let $f\:R\rarrow S$
be a flat tight continuous ring map.
 Then, for any quotseparated contramodule $R$\+module\/ $\fP$,
the $S$\+module $f^!(\fP)=\Hom_R(S,\fP)$ is also a quotseparated
contramodule.
\end{prop}

\begin{proof}[First proof]
 By Lemma~\ref{Hom-is-a-quotseparated-contramodule}(c),
the underlying $R$\+module of $\Hom_R(S,\fP)$ is a quotseparated
contramodule $R$\+module.
 By Theorem~\ref{flat-reflects-quotseparated-contramodules-thm},
it follows that $\Hom_R(S,\fP)$ is a quotseparated contramodule
$S$\+module.
\end{proof}

\begin{proof}[Second proof]
 The following argument is a special case
of~\cite[proof of Theorem~10.4]{Pcs}.
 The $S$\+module $\Hom_R(S,\fP)$ is a contramodule by
Lemma~\ref{tight-co-extension-of-scalars}(b); we only have to prove
that it is a quotseparated contramodule.
 By Proposition~\ref{flat-cotorsion-pair-in-quotseparated}, there
exists a short exact sequence of quotseparated contramodule
$R$\+modules $0\rarrow\fC\rarrow\fF\rarrow\fP\rarrow0$ with a flat
quotseparated contramodule $R$\+module $\fF$ and a cotorsion
quotseparated contramodule $R$\+module~$\fC$.
 By Lemma~\ref{flat-reductions-Hom-exactness-lemma} applied to
the $R$\+module $F=S$, the short sequence of $S$\+modules
$0\rarrow\Hom_R(S,\fC)\rarrow\Hom_R(S,\fF)\rarrow\Hom_R(S,\fP)\rarrow0$
is exact.
 Now the contramodule $R$\+module $\fF$ is separated by
Lemma~\ref{flat-quotseparated-contramodules-separated}, hence
the contramodule $S$\+module $\Hom_R(S,\fP)$ is separated
by Lemma~\ref{separated-coextension-of-scalars-prop}.
 Thus the contramodule $S$\+module $\Hom_R(S,\fP)$ is a quotient
$S$\+module of a separated contramodule $S$\+module
$\Hom_R(S,\fF)$, as desired.
\end{proof}

 Now we can finish the discussion of right adjoint functors to
the contramodule restriction of scalars that was started in
Section~\ref{prelim-change-of-scalars-subsecn}.
 For any flat tight continuous map of adic topological rings
$f\:R\rarrow S$,
Proposition~\ref{flat-quotseparated-coextension-of-scalars-prop} claims
that the functor of coextension of scalars $f^!\:R\Modl\rarrow S\Modl$
restricts to a functor
$$
 f^!\:R\Modl_\ctra^\qs\lrarrow S\Modl_\ctra^\qs,
$$
which we will also call the \emph{coextension of scalars}.
 The functor~$f^!$ is right adjoint to the functor of restriction of
scalars $f_\sharp\:S\Modl_\ctra^\qs\rarrow R\Modl_\ctra^\qs$.

\subsection{Colocalization of contraadjustedness and colocality of
cotorsion} \label{prelim-colocalization-and-colocality-subsecn}
 We continue to use the notation $R\Contra=R\Modl_\ctra^\qs$.

\begin{cor} \label{colocalization-contraadjusted-exactness-cor}
 Let $f\:R\rarrow S$ be a very flat tight continuous map of adic
topological rings, and let\/ $0\rarrow\fP\rarrow\fM\rarrow\fN\rarrow0$
be a short exact sequence of quotseparated contramodule $R$\+modules
with a contraadjusted quotseparated contramodule $R$\+module\/~$\fP$.
 Then the short sequence of quotseparated contramodule $S$\+modules\/
$0\rarrow f^!\fP\rarrow f^!\fM\rarrow f^!\fN\rarrow0$ is exact.
\end{cor}

\begin{proof}
 This is the particular case of
Lemma~\ref{very-flat-reductions-Hom-exactness-lemma}
for the $R$\+module $F=S$.
\end{proof}

\begin{cor} \label{colocalization-cotorsion-exactness-cor}
 Let $f\:R\rarrow S$ be a flat tight continuous map of adic topological
rings, and let\/ $0\rarrow\fP\rarrow\fM\rarrow\fN\rarrow0$ be a short
exact sequence of quotseparated contramodule $R$\+modules with
a cotorsion quotseparated contramodule $R$\+module\/~$\fP$.
 Then the short sequence of quotseparated contramodule $S$\+modules\/
$0\rarrow f^!\fP\rarrow f^!\fM\rarrow f^!\fN\rarrow0$ is exact.
\end{cor}

\begin{proof}
 This is the particular case of
Lemma~\ref{flat-reductions-Hom-exactness-lemma}
for the $R$\+module $F=S$.
\end{proof}

\begin{prop} \label{colocalization-of-contraadjusted-contramodule-prop}
 Let $f\:R\rarrow S$ be a very flat tight continuous map of adic
topological rings, and let\/ $\fP$ be a contraadjusted quotseparated
contramodule $R$\+module.
 Then the quotseparated contramodule $S$\+module
$f^!(\fP)=\Hom_R(S,\fP)$ is contraadjusted.
\end{prop}

\begin{proof}
 We have a pair of adjoint functors between abelian categories
$$
 \xymatrix{
  R\Contra \ar@<-2pt>[rr]_{f^!} && S\Contra. \ar@<-2pt>[ll]_{f_\sharp}
 }
$$
 The left adjoint functor~$f_\sharp$ is exact.
 The object $\fP\in R\Contra$ is adjusted to the right adjoint
functor~$f^!$ in the following sense.
 By Lemma~\ref{very-flat-reductions-Hom-exactness-lemma} applied to
the $R$\+module $F=S$, the short sequence of contramodule $S$\+modules
$0\rarrow f^!\fP\rarrow f^!\fM\rarrow f^!\fN\rarrow0$ is exact for any
short exact sequence of quotseparated contramodule $R$\+modules
$0\rarrow\fP\rarrow\fM\rarrow\fN\rarrow0$.
 By~\cite[Lemma~1.7(b)]{Pal}, for any quotseparated contramodule
$R$\+module $\fP$ and any quotseparated contramodule
$S$\+module $\fG$, we have a natural injective map of abelian group
$$
 \Ext^1_{S\Contra}(\fG,f^!\fP)\lrarrow
 \Ext^1_{R\Contra}(f_\sharp\fG,\fP),
$$
which is an isomorphism when $\fP$ is contraadjusted
by~\cite[Lemma~1.7(e)]{Pal}.
 Now if the quotseparated contramodule $S$\+module $\fG$ is very flat,
then so is the quotseparated contramodule $R$\+module $f_\sharp\fG$,
by Lemma~\ref{very-flat-map-of-adic-rings-direct-image-lemma}(b).
 Using the criterion of
Lemma~\ref{contraadjusted-contramodules-lemma}(2), one concludes
that the functor~$f^!$ takes contraadjusted quotseparated contramodule
$R$\+modules to contraadjusted quotseparated contramodule $S$\+modules.
\end{proof}

\begin{prop} \label{colocalization-of-cotorsion-contramodule-prop}
 Let $f\:R\rarrow S$ be a flat tight continuous map of adic
topological rings, and let\/ $\fP$ be a cotorsion quotseparated
contramodule $R$\+module.
 Then the quotseparated contramodule $S$\+module
$f^!(\fP)=\Hom_R(S,\fP)$ is cotorsion.
\end{prop}

\begin{proof}
 This is the cotorsion version of
Proposition~\ref{colocalization-of-contraadjusted-contramodule-prop},
and the proof is very similar.
 One needs to use Lemma~\ref{flat-reductions-Hom-exactness-lemma},
Lemma~\ref{flat-map-of-adic-rings-direct-image-lemma}(d),
and the criterion of Lemma~\ref{cotorsion-contramodules-lemma}(2).
\end{proof}

 We recall the notation ${}_IM=\Hom_R(R/I,M)$ for the submodule of all
elements in an $R$\+module $M$ annihilated by an ideal $I\subset R$.

\begin{lem} \label{adic-contramodule-nakayama}
 Let $I$ be a finitely generated ideal in a commutative ring~$R$.
 In this context: \par
\textup{(a)} Let\/ $\cM$ be an $I$\+torsion $R$\+module.
 Assume that ${}_I\cM=0$.
 Then\/ $\cM=0$. \par
\textup{(b)} Let\/ $\fP$ be an $I$\+contramodule $R$\+module.
 Assume that\/ $\fP=I\fP$.
 Then\/ $\fP=0$.
\end{lem}

\begin{proof}
 Part~(a) is obvious.
 Part~(b) is a very basic example of an assertion from the class of
results known as ``contramodule Nakayama lemmas''.
 See~\cite[Lemma~4.2 and Remark~4.3]{Pcta}, which are applicable
by~\cite[Theorems~3.3 and~4.1]{Pcta}.
 A simple argument by induction on the number of generators of~$I$,
using the passage from $\fP$ to $\fP/s\fP$ (where $s\in I$ is one
of the generators), is also possible.
\end{proof}

\begin{cor} \label{contramodule-nakayama-surjectivity}
 Let $R$ be an adic topological ring with an ideal of definition
$I\subset R$.
 In this setting: \par
\textup{(a)} Let $g\:\cM\rarrow\cN$ be a morphism of torsion
$R$\+modules.
 Then the map~$g$ is injective if and only if the induced map
$\bar g\:{}_I\cM\rarrow{}_I\cN$ is injective. \par
\textup{(b)} Let $g\:\fP\rarrow\fQ$ be a morphism of contramodule
$R$\+modules.
 Then the map $g$~is surjective if and only if the induced map
$\bar g\:\fP/I\fP\rarrow\fQ/I\fQ$ is surjective.
\end{cor}

\begin{proof}
 The point is that, in part~(a), the kernel of the map~$g$ is a torsion
$R$\+module, while, in part~(b), the cokernel of the map~$g$ is
a contramodule $R$\+module.
 Hence the assertions follow from
Lemma~\ref{adic-contramodule-nakayama} applied to the respective kernel
or cokernel.
\end{proof}

\begin{lem} \label{flat-quotsep-contramodules-reduction-exactness}
 Let $R$ be an adic topological ring with an ideal of definition
$I\subset R$.
 In this context: \par
\textup{(a)} Let\/ $0\rarrow\cK^0\rarrow\cK^1\rarrow\cK^2\rarrow\dotsb$
be a bounded below complex of injective torsion $R$\+modules.
 Then the complex\/ $\cK^\bu$ is acyclic if and only if the complex
of $R/I$\+modules\/ $\Hom_R(R/I,\cK^\bu)$ is acyclic. \par
\textup{(b)} Let\/ $\dotsb\rarrow\fF_2\rarrow\fF_1\rarrow
\fF_0\rarrow0$ be a bounded above complex of flat quotseparated
contramodule $R$\+modules.
 Then the complex\/ $\fF_\bu$ is acyclic if and only if the complex
of $R/I$\+modules\/ $R/I\ot_R\fF_\bu$ is acyclic.
\end{lem}

\begin{proof}
 Part~(a): the ``only if'' assertion holds since any acyclic bounded
below complex of injective objects (in an abelian/exact category) is
contractible.
 To prove the ``if'', let $n\ge0$ be the minimal integer for which
the complex $\cK^\bu$ is \emph{not} exact at the cohomological
degree~$n$.
 Then we have a finite contractible complex of injective torsion
$R$\+modules $0\rarrow\cK^0\rarrow\cK^1\rarrow\dotsb\rarrow\cK^{n-1}
\rarrow\cK^n\rarrow\cJ^n\rarrow0$, where $\cJ^n$ is the cokernel of
the differential $\cK^{n-1}\rarrow\cK^n$.
 Furthermore, the natural map $\cJ^n\rarrow\cK^{n+1}$ is \emph{not}
injective.
 Now ${}_I\cJ^n$ is the cokernel of the map ${}_I\cK^{n-1}\rarrow
{}_I\cK^n$, while Corollary~\ref{contramodule-nakayama-surjectivity}(a)
implies that the map ${}_I\cJ^n\rarrow{}_I\cK^{n+1}$ cannot be
injective.

 Part~(b): the ``only if'' assertion holds by
Lemma~\ref{flat-quotseparated-contramodules-well-behaved}.
 To prove the ``if'', let $n\ge0$ be the minimal integer for which
the complex $\fF_\bu$ is \emph{not} exact at the homological degree~$n$.
 Then we have a finite exact complex of flat quotseparated
contramodule $R$\+modules $0\rarrow\fG_n\rarrow\fF_n\rarrow\fF_{n-1}
\rarrow\dotsb\rarrow\fF_1\rarrow\fF_0\rarrow0$, where $\fG_n$ is
the kernel of the differential $\fF_n\rarrow\fF_{n-1}$.
 Furthermore, the natural map $\fF_{n+1}\rarrow\fG_n$ is \emph{not}
surjective.
 Now it follows from
Lemma~\ref{flat-quotseparated-contramodules-well-behaved}
that $R/I\ot_R\fG_n$ is the kernel of the map $R/I\ot_R\fF_n\rarrow
R/I\ot_R\fF_{n-1}$, while
Corollary~\ref{contramodule-nakayama-surjectivity}(b) implies that
the map $R/I\ot_R\fF_{n+1}\rarrow R/I\ot_R\fG_n$ cannot be surjective.
\end{proof}

 In the rest of this section, we work in the following setting.
 Let $f_\alpha\:R\rarrow S_\alpha$, \,$1\le\alpha\le N$, be a finite
formal open covering of an adic topological ring~$R$.
 For every subsequence of indices $1\le\alpha_1<\dotsb<\alpha_k
\le N$, we endow the commutative ring $S_{\alpha_1,\dotsc,\alpha_k}=
S_{\alpha_1}\ot_R\dotsb\ot_R S_{\alpha_k}$ with the tensor product
topology, and denote by $f_{\alpha_1,\dotsc,\alpha_k}$ the induced ring
map $R\rarrow S_{\alpha_1,\dotsc,\alpha_k}$.

\begin{lem} \label{flat-quotsep-cohomol-Cech-sequence-lemma}
 Let\/ $\fF$ be a flat quotseparated contramodule $R$\+module.
 Then the cohomological \v Cech sequence of flat quotseparated
contramodule $R$\+modules
\begin{multline} \label{flat-quotsep-cohomol-Cech-sequence-eqn}
 0\lrarrow\fF\lrarrow\bigoplus\nolimits_{\alpha=1}^N
 f_\alpha{}_\sharp f_\alpha^\sharp\fF\lrarrow
 \bigoplus\nolimits_{1\le\alpha<\beta\le N}
 f_{\alpha,\beta}{}_\sharp f_{\alpha,\beta}^\sharp\fF \\ \lrarrow\dotsb
 \lrarrow f_{1,\dotsc,N}{}_\sharp f_{1,\dotsc,N}^\sharp\fF\lrarrow0
\end{multline}
is exact.
\end{lem}

\begin{proof}
 First of all, the differentials in 
the complex~\eqref{flat-quotsep-cohomol-Cech-sequence-eqn} are
constructed in terms of the adjunction units for the adjoint
functors~$f^\sharp$ and~$f_\sharp$ involved, as usual for
the \v Cech complexes.
 Furthermore, since the compositions of formal open immersions are
formal open immersions, it follows from
Lemma~\ref{open-immersion-adic-ring-base-change-context}(a) that
the maps $f_{\alpha_1,\dotsc,\alpha_k}$ are formal open immersions.
 In particular, $f_{\alpha_1,\dotsc,\alpha_k}$ are flat tight
continuous ring maps.
 Hence the functors $f_{\alpha_1,\dotsc,\alpha_k}{}_\sharp$ take flat
quotseparated contramodules to flat quotseparated contramodules
by Lemma~\ref{flat-map-of-adic-rings-direct-image-lemma}(d).
 The functors $f_{\alpha_1,\dotsc,\alpha_k}^\sharp$ also take flat
quotseparated contramodules to flat quotseparated contramodules by
Corollary~\ref{flat-contramodules-contraextension-of-scalars}.
 Thus \eqref{flat-quotsep-cohomol-Cech-sequence-eqn}~is a finite
complex of flat quotseparated contramodule $R$\+modules.

 Now let $I\subset R$ be an ideal of definition.
 Let $J_{\alpha_1,\dotsc,\alpha_k}\subset S_{\alpha_1,\dotsc,\alpha_k}$
be the ideals spanned by the $f_{\alpha_1,\dotsc,\alpha_k}(I)\subset
S_{\alpha_1,\dotsc,\alpha_k}$.
 Put $\overline R=R/I$ and $\overline S_{\alpha_1,\dotsc,\alpha_k}=
S_{\alpha_1,\dotsc,\alpha_k}/J_{\alpha_1,\dotsc,\alpha_k}$, and
let $\bar f_{\alpha_1,\dotsc,\alpha_k}\:\overline R\rarrow
\overline S_{\alpha_1,\dotsc,\alpha_k}$ be the induced maps of
the quotient rings.
 Denote also by $\overline F$ the $\overline R$\+module $R/I\ot_R\fF$.

 Applying the reduction functor $R/I\ot_R{-}$ to
the complex~\eqref{flat-quotsep-cohomol-Cech-sequence-eqn},
we obtain the complex of flat $\overline R$\+modules
\begin{multline} \label{flat-reduction-cohomol-Cech-sequence-eqn}
 0\lrarrow\overline F\lrarrow\bigoplus\nolimits_{\alpha=1}^N
 \bar f_\alpha{}_*\bar f_\alpha^*\overline F\lrarrow
 \bigoplus\nolimits_{1\le\alpha<\beta\le N}
 \bar f_{\alpha,\beta}{}_*\bar f_{\alpha,\beta}^*\overline F \\
 \lrarrow\dotsb\lrarrow
 \bar f_{1,\dotsc,N}{}_*\bar f_{1,\dotsc,N}^*\overline F\lrarrow0,
\end{multline}
as one can see from
Lemma~\ref{reductions-of-co-contra-extension-of-scalars}(c).
 Now \eqref{flat-reduction-cohomol-Cech-sequence-eqn}~is just
the \v Cech complex for the $\overline R$\+mod\-ule $\overline F$
with respect to the open covering $\bar f_\alpha\:\overline R
\rarrow\overline S_\alpha$ of the (abstract, nontopological)
commutative ring~$\overline R$.
 So the complex~\eqref{flat-reduction-cohomol-Cech-sequence-eqn}
is exact, and by
Lemma~\ref{flat-quotsep-contramodules-reduction-exactness}(b)
it follows that
the complex~\eqref{flat-quotsep-cohomol-Cech-sequence-eqn}
is exact as well.
\end{proof}

\begin{lem} \label{injective-torsion-homol-Cech-sequence-lemma}
 Let\/ $\cK$ be an injective torsion $R$\+module.
 Then the homological \v Cech complex of injective torsion
$R$\+modules
\begin{multline} \label{injective-torsion-homol-Cech-sequence-eqn}
 0\lrarrow f_{1,\dotsc,N}{}_\diamond f_{1,\dotsc,N}^\diamond\cK
 \lrarrow\dotsb \\ \lrarrow\bigoplus\nolimits_{1\le\alpha<\beta\le N}
 f_{\alpha,\beta}{}_\diamond f_{\alpha,\beta}^\diamond\cK\lrarrow
 \bigoplus\nolimits_{\alpha=1}^N
 f_\alpha{}_\diamond f_\alpha^\diamond\cK\lrarrow\cK\lrarrow0
\end{multline}
is exact (or equivalently, contractible).
\end{lem}

\begin{proof}
 This lemma is dual-analogous to
Lemma~\ref{flat-quotsep-cohomol-Cech-sequence-lemma}.
 First of all, the differentials in 
the complex~\eqref{injective-torsion-homol-Cech-sequence-eqn} are
constructed in terms of the adjunction counits for the adjoint
functors~$f_\diamond$ and~$f^\diamond$ involved.
 As mentioned in the proof of
Lemma~\ref{flat-quotsep-cohomol-Cech-sequence-lemma}, the maps
$f_{\alpha_1,\dotsc,\alpha_k}$ are flat tight continuous ring maps.
 Hence the functors $f_{\alpha_1,\dotsc,\alpha_k}{}_\diamond$ take
injective torsion modules to injective torsion modules
by Lemma~\ref{flat-map-of-adic-rings-direct-image-lemma}(b).
 The functors $f_{\alpha_1,\dotsc,\alpha_k}^\diamond$ also take
injective torsion modules to injective torsion modules,
as mentioned in Section~\ref{prelim-change-of-scalars-subsecn}.
 Thus \eqref{injective-torsion-homol-Cech-sequence-eqn}~is a finite
complex of injective torsion $R$\+modules.

 We use the notation from the proof of
Lemma~\ref{flat-quotsep-cohomol-Cech-sequence-lemma}.
 Denote also by $\overline K$ the $\overline R$\+module
$\Hom_R(R/I,\cK)$.
 Following the discussion in
Section~\ref{prelim-inj-proj-flat-torsion-contra-subsecn},
the $R/I$\+module $\overline K$ is injective.
 In particular, it follows that $\overline K$ is a contraadjusted
$\overline R$\+module.

 Applying the functor $\Hom_R(R/I,{-})$ to
the complex~\eqref{injective-torsion-homol-Cech-sequence-eqn},
we obtain the complex of injective $\overline R$\+modules
\begin{multline} \label{injective-reduction-homol-Cech-sequence-eqn}
 0\lrarrow \bar f_{1,\dotsc,N}{}_*\bar f_{1,\dotsc,N}^!\overline K
 \lrarrow\dotsb \\ \lrarrow\bigoplus\nolimits_{1\le\alpha<\beta\le N}
 f_{\alpha,\beta}{}_* f_{\alpha,\beta}^!\overline K\lrarrow
 \bigoplus\nolimits_{\alpha=1}^N
 f_\alpha{}_*f_\alpha^!\overline K\lrarrow\overline K\lrarrow0,
\end{multline}
as one can see from 
Lemma~\ref{reductions-of-co-contra-extension-of-scalars}(a).
 Now \eqref{injective-reduction-homol-Cech-sequence-eqn}~is
the homological \v Cech complex~\cite[formula~(1.3) from
Section~1.2]{Pcosh} for the contraadjusted $\overline R$\+module
$\overline K$ with respect to the open covering
$\bar f_\alpha\:\overline R\rarrow\overline S_\alpha$ of
the (abstract, nontopological) commutative ring~$\overline R$.
 By~\cite[Lemma~1.2.6(b)]{Pcosh},
the complex~\eqref{injective-reduction-homol-Cech-sequence-eqn}
is exact, and
by Lemma~\ref{flat-quotsep-contramodules-reduction-exactness}(a)
above it follows that
the complex~\eqref{injective-torsion-homol-Cech-sequence-eqn}
is exact as well.
\end{proof}

\begin{cor} \label{for-the-rings-cohomol-Cech-sequence-cor}
 The cohomological \v Cech sequence of very flat quotseparated
contramodule $R$\+modules
\begin{multline} \label{for-the-rings-cohomol-Cech-sequence-eqn}
 0\lrarrow\Lambda(R)\lrarrow\bigoplus\nolimits_{\alpha=1}^N
 \Lambda(S_\alpha)\lrarrow\bigoplus\nolimits_{1\le\alpha<\beta\le N}
 \Lambda(S_{\alpha,\beta}) \\ \lrarrow\dotsb\lrarrow
 \Lambda(S_{1,\dotsc,N})\lrarrow0
\end{multline}
is exact.
\end{cor}

\begin{proof}
 In the context of Lemma~\ref{flat-quotsep-cohomol-Cech-sequence-lemma},
put $\fF=\Lambda(R)$.
 Clearly, $\fF$ is a very flat (in fact, free) quotseparated
contramodule $R$\+module.
 As explained in the proof of
Lemma~\ref{flat-quotsep-cohomol-Cech-sequence-lemma}, the maps
$f_{\alpha_1,\dotsc,\alpha_k}$ are formal open immersions;
so they are very flat tight continuous ring maps.
 By Lemma~\ref{very-flat-map-of-adic-rings-direct-image-lemma}(b), it
follows that
the complex~\eqref{flat-quotsep-cohomol-Cech-sequence-eqn} for
$\fF=\Lambda(R)$ is a complex of very flat quotseparated contramodule
$R$\+modules.
 Furthermore, since $f_{\alpha_1,\dotsc,\alpha_k}$ are tight continuous
ring maps, the topology on $S_{\alpha_1,\dotsc,\alpha_k}$ coincides
with the $I$\+adic topology (for any ideal of definition $I\subset R$,
by Lemma~\ref{continuous-and-tight-ring-map-lemma}(3)); so
the notation $\Lambda(S_{\alpha_1,\dotsc,\alpha_k})$ is unambiguous.
 Finally, following the discussion in
Section~\ref{prelim-change-of-scalars-subsecn}, we have
$f_{\alpha_1,\dotsc,\alpha_k}^\sharp(\Lambda(R))\simeq
\Lambda(S_{\alpha_1,\dotsc,\alpha_k})$; so
the complex~\eqref{for-the-rings-cohomol-Cech-sequence-eqn} is
the special case of~\eqref{flat-quotsep-cohomol-Cech-sequence-eqn}
for $\fF=\Lambda(R)$.
\end{proof}

\begin{lem} \label{torsion-cohomol-Cech-sequence-lemma}
 Let\/ $\cM$ be a torsion $R$\+module.
 Then the cohomological \v Cech sequence of torsion $R$\+modules
\begin{multline} \label{torsion-cohomol-Cech-sequence-eqn}
 0\lrarrow\cM\lrarrow\bigoplus\nolimits_{\alpha=1}^N
 f_\alpha{}_\diamond f_\alpha^*\cM\lrarrow
 \bigoplus\nolimits_{1\le\alpha<\beta\le N}
 f_{\alpha,\beta}{}_\diamond f_{\alpha,\beta}^*\cM \\ \lrarrow\dotsb
 \lrarrow f_{1,\dotsc,N}{}_\diamond f_{1,\dotsc,N}^*\cM\lrarrow0
\end{multline}
is exact.
\end{lem}

\begin{proof}
 The differentials in
the complex~\eqref{torsion-cohomol-Cech-sequence-eqn} are constructed
in terms of the adjunction units for the adjoint functors~$f^*$
and~$f_\diamond$ involved.
 We observe that the same
complex~\eqref{torsion-cohomol-Cech-sequence-eqn} can be obtained
by applying the functor $\cM\ot_R{-}$ to
the complex~\eqref{for-the-rings-cohomol-Cech-sequence-eqn}.
 It follows from
Lemma~\ref{flat-quotseparated-contramodules-well-behaved} that applying
the functor $\cM\ot_R{-}$ takes
the acyclic complex~\eqref{for-the-rings-cohomol-Cech-sequence-eqn} to
an acyclic complex.
\end{proof}

\begin{lem} \label{contraadj-homol-Cech-sequence-lemma}
 Let\/ $\fP$ be a contraadjusted quotseparated contramodule $R$\+module.
 Then the homological \v Cech sequence of contraadjusted quotseparated
contramodule $R$\+modules
\begin{multline} \label{contraadj-homol-Cech-sequence-eqn}
 0\lrarrow f_{1,\dotsc,N}{}_\sharp f_{1,\dotsc,N}^!\fP\lrarrow\dotsb \\
 \lrarrow\bigoplus\nolimits_{1\le\alpha<\beta\le N}
 f_{\alpha,\beta}{}_\sharp f_{\alpha,\beta}^!\fP\lrarrow
 \bigoplus\nolimits_{\alpha=1}^N f_\alpha{}_\sharp f_\alpha^!\fP
 \lrarrow\fP\lrarrow0
\end{multline}
is exact.
\end{lem}

\begin{proof}
 This lemma is dual-analogous to
Lemma~\ref{torsion-cohomol-Cech-sequence-lemma}.
 First of all, the differentials in 
the complex~\eqref{contraadj-homol-Cech-sequence-eqn} are
constructed in terms of the adjunction counits for the adjoint
functors~$f_\sharp$ and~$f^!$ involved.
 Furthermore, the terms of
the complex~\eqref{contraadj-homol-Cech-sequence-eqn} are
contraadjusted quotseparated contramodule $R$\+modules by
Lemma~\ref{restriction-of-scalars-contraadjusted} and
Proposition~\ref{colocalization-of-contraadjusted-contramodule-prop}.
 Now we observe that the same
complex~\eqref{contraadj-homol-Cech-sequence-eqn} can be obtained
by applying the functor $\Hom_R({-},\fP)$ to
the complex~\eqref{for-the-rings-cohomol-Cech-sequence-eqn}.
 The class of very flat quotseparated contramodule $R$\+modules is
closed under kernels of epimorphisms in $R\Modl_\ctra^\qs$, so
the $R$\+modules of cocycles
in~\eqref{for-the-rings-cohomol-Cech-sequence-eqn} are very flat
quotseparated contramodule $R$\+modules.
 Therefore, it follows from
Lemma~\ref{contraadjusted-contramodules-lemma}(1) that applying
the functor $\Hom_R({-},\fP)$ takes
the acyclic complex~\eqref{for-the-rings-cohomol-Cech-sequence-eqn}
to an acyclic complex.
\end{proof}

\begin{cor} \label{colocality-of-cotorsion-presuming-contraadjusted}
 Let\/ $\fP$ be a contraadjusted quotseparated contramodule
$R$\+mod\-ule.
 Assume that, for every index\/~$\alpha$, the quotseparated
contramodule $S_\alpha$\+module $f_\alpha^!(\fP)$ is cotorsion.
 Then the quotseparated contramodule $R$\+module\/ $\fP$ is cotorsion.
\end{cor}

\begin{proof}
 By Lemma~\ref{restriction-of-scalars-cotorsion} and
Proposition~\ref{colocalization-of-cotorsion-contramodule-prop},
if all the quotseparated contramodule $S_\alpha$\+modules
$f_\alpha^!(\fP)$ are cotorsion, then all the quotseparated
contramodule $R$\+modules $f_{\alpha_1,\dotsc,\alpha_k}{}_\sharp
f_{\alpha_1,\dotsc,\alpha_k}^!\fP$ for $k\ge1$ are cotorsion, too.
 Recall that the class of cotorsion quotseparated contramodule
$R$\+modules is closed under cokernels of monomorphisms, as mentioned
in Section~\ref{prelim-cotorsion-contramodules-subsecn}.
 Going from the leftmost to the rightmost end of
the complex~\eqref{contraadj-homol-Cech-sequence-eqn} and arguing
by induction, one proves that the quotseparated contramodule
$R$\+module $\fP$ is cotorsion as well.
\end{proof}

 The definition of the \emph{cotorsion dimension} of a quotseparated
contramodule $R$\+module was given near the end of
Section~\ref{prelim-cotorsion-contramodules-subsecn}.
 The next corollary is our version of~\cite[Lemma~1.5.4(c)]{Pcosh}.

\begin{cor} \label{colocality-of-cotorsion-dimension}
 Let\/ $\fP$ be a contraadjusted quotseparated contramodule
$R$\+mod\-ule.
 Then the cotorsion dimension of the quotseparated contramodule
$R$\+module\/ $\fP$ is equal to the supremum of the cotorsion dimensions
of the quotseparated contramodule $S_\alpha$\+modules
$f_\alpha^!(\fP)$.
\end{cor}

\begin{proof}
 The assertion follows from
Corollary~\ref{colocalization-contraadjusted-exactness-cor},
Proposition~\ref{colocalization-of-cotorsion-contramodule-prop},
and Corollary~\ref{colocality-of-cotorsion-presuming-contraadjusted}.
\end{proof}

\subsection{Colocality of exactness}
\label{prelim-colocality-of-exactness-subsecn}
 Let $R$ be an adic topological ring.
 We will denote by $R\Contra^\cta\subset R\Contra=R\Modl_\ctra^\qs$
the full subcategory of contraadjusted quotseparated contramodule
$R$\+modules.
 Similarly, the notation $R\Contra^\cot\subset R\Contra$ stands for
the full subcategory of cotorsion quotseparated contramodule
$R$\+modules.
 
 As mentioned in
Sections~\ref{prelim-veryflat-and-contraadjusted-subsecn}
and~\ref{prelim-cotorsion-contramodules-subsecn}, the full subcategories
$R\Contra^\cta$ and $R\Contra^\cot$ are closed under extensions,
cokernels of monomorphisms, and infinite products in $R\Contra$.
 (The full subcategory $R\Contra^\cta$ is even closed under quotients.)
 Thus these full subcategories inherit exact category structures from
the abelian exact structure on $R\Contra$.
 We will always view $R\Contra^\cta$ and $R\Contra^\cot$ as exact
categories with the exact structures inherited from (the abelian
exact structure on) $R\Contra$.

 We will also denote by $R\Contra_\fl\subset R\Contra$ the full
subcategory of flat quotseparated contramodule $R$\+modules, and
by $R\Contra_\vfl\subset R\Contra$ the full subcategory of
very flat quotseparated contramodule $R$\+modules.
 According to
Lemma~\ref{flat-quotseparated-contramodules-well-behaved}(a) and
the discussion in 
Section~\ref{prelim-veryflat-and-contraadjusted-subsecn},
the full subcategories $R\Contra_\fl$ and $R\Contra_\vfl$ are
closed under extensions and kernels of epimorphisms in $R\Contra$.
 Thus these full subcategories also inherit exact category structures
from the abelian exact structure on $R\Contra$.
 We will always view $R\Contra_\fl$ and $R\Contra_\vfl$ as exact
categories with the exact structures inherited from $R\Contra$.

 For any abstract (nontopological) commutative ring $R$, the similar
notation $R\Modl^\cta\subset R\Modl$ and $R\Modl_\vfl\subset R\Modl$
stands for the exact categories of contraadjusted and very flat
$R$\+modules.
 For an associative ring $R$, we denote by $R\Modl^\cot\subset R\Modl$
and $R\Modl_\fl\subset R\Modl$ the exact subcategories of cotorsion
and flat left $R$\+modules in the abelian category $R\Modl$.

\begin{lem} \label{contractible-injective-complexes}
 Let $R$ be an adic topological ring and\/ $\cK^\bu$ be a complex of
injective torsion $R$\+modules.
 Then the following three conditions are equivalent:
\begin{enumerate}
\item the complex\/ $\cK^\bu$ is contractible;
\item there exists an ideal of definition $I\subset R$ such that, for
every $n\ge1$, the complex of injective $R/I^n$\+modules
${}_{I^n}\cK^\bu$ is contractible;
\item for every open ideal $I\subset R$, the complex of injective
$R/I$\+modules ${}_I\cK^\bu$ is contractible.
\end{enumerate}
\end{lem}

\begin{proof}
 (1)~$\Longrightarrow$~(3) $\Longrightarrow$~(2) Obvious.

 (2)~$\Longrightarrow$~(1) For every pair of integers $n\ge1$
and $i\in\boZ$, denote by $L_n^i$ the image of the differential
${}_{I^n}\cK^i\rarrow{}_{I^n}\cK^{i+1}$.
 Then condition~(3) means that the $R/I^n$\+modules $L_n^i$ are
injective and the short sequences of $R/I^n$\+modules $0\rarrow L_n^i
\rarrow{}_{I^n}\cK^i\rarrow L_n^{i+1}\rarrow0$ are (split) exact.
 As $n$~varies for $i$~fixed, the latter short exact sequences form
an inductive system.
 Hence we have the short exact sequences of inductive limits
$0\rarrow\varinjlim_{n\ge1}L_n^i\rarrow\cK^i\rarrow\varinjlim_{n\ge1}
L_n^{i+1}\rarrow0$.

 Furthermore, the $R/I^n$\+module $L_n^i$ can be also constructed as
the kernel of the differential ${}_{I_n}\cK^i\rarrow{}_{I^n}\cK^{i+1}$.
 As the functors $\Hom_R(R/I^n,{-})$ preserve kernels, it follows that
the natural maps $L_n^i\rarrow{}_{I^n}L_{n+1}^i$ are isomorphisms.
 Put $\cL^i=\varinjlim_{n\ge1}L_n^i$; then we have
$L_n^i={}_{I^n}\cL^i$.
 So the $R$\+modules $\cL^i$ are injective torsion $R$\+modules.
 Finally, the complex $\cK^\bu$ can be obtained by splicing
the (necessarily split) short exact sequences of injective torsion
$R$\+modules $0\rarrow\cL^i\rarrow\cK^i\rarrow\cL^{i+1}\rarrow0$;
so it is contractible.

 For a generalization, see~\cite[Lemma~4.21]{Psemten}, which is
applicable in view of~\cite[Section~2.4(6)]{Psemten}.
\end{proof}

\begin{lem} \label{acyclic-flat-quotseparated-complexes}
 Let $R$ be an adic topological ring and\/ $\fF^\bu$ be a complex of
flat quotseparated contramodule $R$\+modules.
 Then the following three conditions are equivalent:
\begin{enumerate}
\item the complex\/ $\fF^\bu$ is exact in the exact category
$R\Contra_\fl$;
\item there exists an ideal of definition $I\subset R$ such that, for
every $n\ge1$, the complex of flat $R/I^n$\+modules $R/I^n\ot_R\fF^\bu$
is exact in the exact category $R/I^n\Modl_\fl$;
\item for every open ideal $I\subset R$, the complex of flat
$R/I$\+modules $R/I\ot_R\fF^\bu$ is exact in the exact category
$R/I\Modl_\fl$.
\end{enumerate}
\end{lem}

\begin{proof}
 (1)~$\Longrightarrow$~(3) Follows from
Lemma~\ref{flat-quotseparated-contramodules-well-behaved}(b).

 (3)~$\Longrightarrow$~(2) Obvious.
 
 (2)~$\Longrightarrow$~(1) For every pair of integers $n\ge1$ and
$i\in\boZ$, denote by $G_n^i$ the image of the differential
$\fF^{i-1}/I^n\fF^{i-1}\rarrow\fF^i/I^n\fF^i$.
 Then condition~(2) means that the $R/I^n$\+modules $G_n^i$ are
flat and the short sequences of $R/I^n$\+modules $0\rarrow G_n^i
\rarrow\fF^i/I^n\fF^i\rarrow G_{n+1}^i\rarrow0$ are exact.
 As $n$~varies for $i$~fixed, the latter short exact sequences form
a projective system with termwise surjective transition maps.
 Hence we have the short exact sequences of projective limits
$0\rarrow\varprojlim_{n\ge1}G_n^i\rarrow\varprojlim_{n\ge1}
\fF^i/I^n\fF^i\rarrow\varprojlim_{n\ge1}G_n^{i+1}\rarrow0$.
 Following the discussion in Section~\ref{prelim-contramodules-subsecn}
and Lemma~\ref{flat-quotseparated-contramodules-separated},
the $R$\+modules $\fF^i$ are separated and complete; so
$\varprojlim_{n\ge1}\fF^i/I^n\fF^i\simeq\fF^i$.

 Furthermore, the $R/I^n$\+module $G_n^i$ can be also constructed as
the cokernel of the differential
$\fF^{i-2}/I^n\fF^{i-2}\rarrow\fF^{i-1}/I^n\fF^{i-1}$.
 As the tensor product functors preserve cokernels, it follows
that the natural maps $R/I^n\ot_{R/I^{n+1}}G_{n+1}^i\rarrow G_n^i$
are isomorphisms.
 Put $\fG^i=\varprojlim_{n\ge1}G_n^i$; then, by~\cite[Theorem~1.2
or~2.8]{Yek2} or~\cite[Lemma~E.1.3]{Pcosh}, we have
$G_n^i\simeq\fG^i/I^n\fG^i$.
 So the $R$\+modules $\fG_n^i$ are flat quotseparated contramodule
$R$\+modules.
 Finally, the complex $\fF^\bu$ can be obtained by splicing the short
exact sequences of flat quotseparated contramodule $R$\+modules
$0\rarrow\fG^i\rarrow\fF^i\rarrow\fG^{i+1}\rarrow0$, as desired.

 For a generalization, see~\cite[Lemma~4.13]{Psemten}, which is
applicable in view of~\cite[Example~3.8(2)]{Psemten}.
\end{proof}

\begin{lem} \label{acyclic-veryflat-quotseparated-complexes}
 Let $R$ be an adic topological ring and\/ $\fF^\bu$ be a complex of
very flat quotseparated contramodule $R$\+modules.
 Then the following three conditions are equivalent:
\begin{enumerate}
\item the complex\/ $\fF^\bu$ is exact in the exact category
$R\Contra_\vfl$;
\item there exists an ideal of definition $I\subset R$ such that,
for every $n\ge1$, the complex of very flat $R/I^n$\+modules
$R/I^n\ot_R\fF^\bu$ is exact in the exact category $R/I^n\Modl_\vfl$;
\item for every open ideal $I\subset R$, the complex of very flat
$R/I$\+modules $R/I\ot_R\fF^\bu$ is exact in the exact category
$R/I\Modl_\vfl$.
\end{enumerate}
\end{lem}

\begin{proof}
 Similar to the proof of
Lemma~\ref{acyclic-flat-quotseparated-complexes}.
\end{proof}

 Now let $f_\alpha\:R\rarrow S_\alpha$, \,$1\le\alpha\le N$, be
a finite formal open covering of an adic topological ring~$R$.
 As in Section~\ref{prelim-colocalization-and-colocality-subsecn},
we put $S_{\alpha_1,\dotsc,\alpha_k}=
S_{\alpha_1}\ot_R\dotsb\ot_R S_{\alpha_k}$ for all
$1\le\alpha_1<\dotsb<\alpha_k\le N$, endow the rings
$S_{\alpha_1,\dotsc,\alpha_k}$ with the tensor product topologies, and
denote the induced ring maps by $f_{\alpha_1,\dotsc,\alpha_k}\:
R\rarrow S_{\alpha_1,\dotsc,\alpha_k}$.

\begin{cor} \label{locality-of-surjectivity-for-contramodules}
\textup{(a)} A morphism of contramodule $R$\+modules\/ $\fP\rarrow\fQ$
is an epimorphism in $R\Modl_\ctra$ if and only if the induced morphism
of contramodule $S_\alpha$\+modules $f_\alpha^\#(\fP)\rarrow
f_\alpha^\#(\fQ)$ is an epimorphism in $S_\alpha\Modl_\ctra$
for every\/~$\alpha$. \par
\textup{(b)} A morphism of quotseparated contramodule $R$\+modules\/
$\fP\rarrow\fQ$ is an epimorphism in $R\Contra$ if and only if
the induced morphism of quotseparated contramodule $S_\alpha$\+modules
$f_\alpha^\sharp(\fP)\rarrow f_\alpha^\sharp(\fQ)$ is an epimorphism
in $S_\alpha\Contra$ for every\/~$\alpha$.
\end{cor}

\begin{proof}
 Part~(b): the ``only if'' assertion holds because the
functor~$f^\sharp$ is right exact (as a left adjoint functor) for any
continuous homomorphism of adic topological rings $f\:R\rarrow S$.
 To prove the ``if'', let $I\subset R$ be an ideal of definition in~$R$.
 By Corollary~\ref{contramodule-nakayama-surjectivity}(b), it suffices
to show that the map $\fP/I\fP\rarrow\fQ/I\fQ$ is surjective.
 Put $J_\alpha=S_\alpha f_\alpha(I)\subset S_\alpha$.
 Then the collection of ring homomorphisms $R/I\rarrow
S_\alpha/J_\alpha$ is an open covering.
 So it is sufficient to check that the map
$S_\alpha/J_\alpha\ot_{R/I}\fP/I\fP\rarrow
S_\alpha/J_\alpha\ot_{R/I}\fQ/I\fQ$ is surjective for every~$\alpha$.
 For this purpose, we notice that the map
$f_\alpha^\sharp(\fP)/J_\alpha f_\alpha^\sharp(\fP)\rarrow
f_\alpha^\sharp(\fQ)/J_\alpha f_\alpha^\sharp(\fQ)$ is surjective, and
refer to Lemma~\ref{reductions-of-co-contra-extension-of-scalars}(c).
 The proof of part~(a) is similar and based on
Lemma~\ref{reductions-of-co-contra-extension-of-scalars}(b).
\end{proof}

\begin{lem} \label{locality-of-exactness-for-flat-quotseparated}
\textup{(a)} A short sequence of flat quotseparated contramodule
$R$\+modules\/ $0\rarrow\fF\rarrow\fG\rarrow\fH\rarrow0$ is exact if
and only if the short sequence of flat quotseparated contramodule
$S_\alpha$\+modules\/ $0\rarrow f_\alpha^\sharp(\fF)\rarrow
f_\alpha^\sharp(\fG)\rarrow f_\alpha^\sharp(\fH)\rarrow0$ is exact
for every\/~$\alpha$. \par
\textup{(b)} A morphism of flat quotseparated contramodule
$R$\+modules\/ $\fF\rarrow\fG$ is an admissible monomorphism in
$R\Contra_\fl$ if and only if the induced morphism of flat
quotseparated contramodule $S_\alpha$\+modules $f_\alpha^\sharp(\fF)
\rarrow f_\alpha^\sharp(\fG)$ is an admissible monomorphism in
$S_\alpha\Contra_\fl$ for every\/~$\alpha$. \par
\textup{(c)} A morphism of flat quotseparated contramodule
$R$\+modules\/ $\fG\rarrow\fH$ is an admissible epimorphism in
$R\Contra_\fl$ if and only if the induced morphism of flat
quotseparated contramodule $S_\alpha$\+modules $f_\alpha^\sharp(\fG)
\rarrow f_\alpha^\sharp(\fH)$ is an admissible epimorphism in
$S_\alpha\Contra_\fl$ for every\/~$\alpha$. \par
\textup{(d)} A complex of flat quotseparated contramodule
$R$\+modules\/ $\fF^\bu$ is exact in $R\Contra_\fl$ if and only if
the complex of flat quotseparated contramodule $S_\alpha$\+modules
$f_\alpha^\sharp(\fF^\bu)$ is exact in $S_\alpha\Contra_\fl$ for
every\/~$\alpha$.
\end{lem}

\begin{proof}
 Part~(a): the ``only if'' assertion holds by
Corollary~\ref{flat-contramodules-contraextension-of-scalars}.
 To prove the ``if'', notice that, by the same corollary, exactness
of the short sequences $0\rarrow f_\alpha^\sharp(\fF)\rarrow
f_\alpha^\sharp(\fG)\rarrow f_\alpha^\sharp(\fH)\rarrow0$ implies
exactness of the short sequences $0\rarrow
f_{\alpha_1,\dotsc,\alpha_k}^\sharp(\fF)\rarrow
f_{\alpha_1,\dotsc,\alpha_k}^\sharp(\fG)\rarrow
f_{\alpha_1,\dotsc,\alpha_k}^\sharp(\fH)\rarrow0$
for all $k\ge1$.
 Now the three-term sequence of flat quotseparated contramodule
$R$\+modules $0\rarrow\fF\rarrow\fG\rarrow\fH\rarrow0$ induces
a three-term sequence of finite exact
sequences~\eqref{flat-quotsep-cohomol-Cech-sequence-eqn}
from Lemma~\ref{flat-quotsep-cohomol-Cech-sequence-lemma}.
 Since the latter three-term sequence is exact on all the terms
of the finite exact
sequences~\eqref{flat-quotsep-cohomol-Cech-sequence-eqn} except
perhaps the leftmost ones, it follows that the short sequence of
the leftmost terms $0\rarrow\fF\rarrow\fG\rarrow\fH\rarrow0$
is exact, too.

 Part~(b): the ``only if'' assertion follows from part~(a).
 To prove the ``if'', notice that the map $\fF\rarrow
\bigoplus_\alpha f_\alpha{}_\sharp f_\alpha^\sharp(\fF)$ is
an admissible monomorphism in $R\Contra_\fl$
by Lemmas~\ref{flat-quotsep-cohomol-Cech-sequence-lemma}
and~\ref{flat-quotseparated-contramodules-well-behaved}(a).
 If the map $f_\alpha^\sharp(\fF)\rarrow f_\alpha^\sharp(\fG)$ is
an admissible monomorphism in $S_\alpha\Contra_\fl$, then the map
$f_\alpha{}_\sharp f_\alpha^\sharp(\fF)\rarrow f_\alpha{}_\sharp
f_\alpha^\sharp(\fG)$ is an admissible monomorphism in $R\Contra_\fl$
by Lemma~\ref{flat-map-of-adic-rings-direct-image-lemma}(d) for
every~$\alpha$.
 Hence the composition $\fF\rarrow\bigoplus_\alpha f_\alpha{}_\sharp
f_\alpha^\sharp(\fF)\rarrow\bigoplus_\alpha f_\alpha{}_\sharp
f_\alpha^\sharp(\fG)$ is an admissible monomorphism in $R\Contra_\fl$.
 The latter composition is equal to the composition
$\fF\rarrow\fG\rarrow\bigoplus_\alpha f_\alpha{}_\sharp
f_\alpha^\sharp(\fG)$, and it follows that $\fF\rarrow\fG$ is
an admissible monomorphism in $R\Contra_\fl$ by the dual version
of~\cite[Proposition~7.6]{Bueh}.

 Part~(c): a morphism in $R\Contra_\fl$ is an admissible epimorphism
if and only if it is an epimorphism in $R\Contra$.
 So the assertion follows from
Corollary~\ref{locality-of-surjectivity-for-contramodules}(b).

 Part~(d): the ``only if'' assertion follows from part~(a).
 To prove the ``if'', let $I\subset R$ be an ideal of definition
in~$R$.
 By Lemma~\ref{acyclic-flat-quotseparated-complexes}(2), it suffices
to show that the complex $\fF^\bu/I\fF^\bu$ is exact in $R/I\Modl_\fl$.
 Continuing in the notation from the proof of
Corollary~\ref{locality-of-surjectivity-for-contramodules},
it is sufficient to check that the complex $S_\alpha/J_\alpha\ot_{R/I}
\fF^\bu/I\fF^\bu$ is exact in $S_\alpha/J_\alpha\Modl_\fl$
for every~$\alpha$.
 For this purpose, we notice that the complex
$f_\alpha^\sharp(\fF^\bu)/J_\alpha f_\alpha^\sharp(\fF^\bu)$ is
acyclic in $S_\alpha/J_\alpha\Modl_\fl$ and refer to
Lemma~\ref{reductions-of-co-contra-extension-of-scalars}(c).
 Alternatively, the argument similar to the proof of part~(a) and
dual-analogous to the proof of
Lemma~\ref{colocality-of-split-exactness-for-injective}(c) below
also works (see also the proof of
Lemma~\ref{colocality-of-exactness-for-contraadjusted}(c)).
\end{proof}

\begin{lem} \label{locality-of-exactness-for-very-flat-quotseparated}
\textup{(a)} A short sequence of very flat quotseparated contramodule
$R$\+modules\/ $0\rarrow\fF\rarrow\fG\rarrow\fH\rarrow0$ is exact if
and only if the short sequence of very flat quotseparated contramodule
$S_\alpha$\+modules\/ $0\rarrow f_\alpha^\sharp(\fF)\rarrow
f_\alpha^\sharp(\fG)\rarrow f_\alpha^\sharp(\fH)\rarrow0$ is exact
for every\/~$\alpha$. \par
\textup{(b)} A morphism of very flat quotseparated contramodule
$R$\+modules\/ $\fF\rarrow\fG$ is an admissible monomorphism in
$R\Contra_\vfl$ if and only if the induced morphism of very flat
quotseparated contramodule $S_\alpha$\+modules $f_\alpha^\sharp(\fF)
\rarrow f_\alpha^\sharp(\fG)$ is an admissible monomorphism in
$S_\alpha\Contra_\vfl$ for every\/~$\alpha$. \par
\textup{(c)} A morphism of very flat quotseparated contramodule
$R$\+modules\/ $\fG\rarrow\fH$ is an admissible epimorphism in
$R\Contra_\vfl$ if and only if the induced morphism of very flat
quotseparated contramodule $S_\alpha$\+modules $f_\alpha^\sharp(\fG)
\rarrow f_\alpha^\sharp(\fH)$ is an admissible epimorphism in
$S_\alpha\Contra_\vfl$ for every\/~$\alpha$. \par
\textup{(d)} A complex of very flat quotseparated contramodule
$R$\+modules\/ $\fF^\bu$ is exact in $R\Contra_\vfl$ if and only if
the complex of flat quotseparated contramodule $S_\alpha$\+modules
$f_\alpha^\sharp(\fF^\bu)$ is exact in $S_\alpha\Contra_\vfl$ for
every\/~$\alpha$.
\end{lem}

\begin{proof}
 Part~(a) is a particular case of
Lemma~\ref{locality-of-exactness-for-flat-quotseparated}(a)
(the result of Lemma~\ref{contraextension-of-scalars-very-flat}
is presumed).
 The proof of part~(b) is similar to that of
Lemma~\ref{locality-of-exactness-for-flat-quotseparated}(b).
 One needs to use Lemmas~\ref{contraextension-of-scalars-very-flat}
and~\ref{very-flat-map-of-adic-rings-direct-image-lemma}(b).
 Part~(c) is similar to
Lemma~\ref{locality-of-exactness-for-flat-quotseparated}(c).
 The proof of part~(d) is similar to that of
Lemma~\ref{locality-of-exactness-for-flat-quotseparated}(d).
 One needs to use the criterion of
Lemma~\ref{acyclic-veryflat-quotseparated-complexes}(2).
\end{proof}

\begin{lem} \label{locality-of-exactness-for-torsion}
\textup{(a)} A short sequence of torsion $R$\+modules\/ $0\rarrow\cL
\rarrow\cM\rarrow\cN\rarrow0$ is exact if and only if the short
sequence of torsion $S_\alpha$\+modules\/ $0\rarrow f_\alpha^*(\cL)
\rarrow f_\alpha^*(\cM)\rarrow f_\alpha^*(\cN)\rarrow0$ is exact
for every\/~$\alpha$. \par
\textup{(b)} A morphism of torsion $R$\+modules\/ $\cL\rarrow\cM$ is
injective if and only if the induced morphism of torsion
$S_\alpha$\+modules $f_\alpha^*(\cL)\rarrow f_\alpha^*(\cM)$ is
injective for every\/~$\alpha$. \par
\textup{(c)} A morphism of torsion $R$\+modules\/ $\cM\rarrow\cN$
is surjective if and only if the induced morphism of torsion
$S_\alpha$\+modules $f_\alpha^*(\cM)\rarrow f_\alpha^*(\cN)$ is
surjective for every\/~$\alpha$. \par
\textup{(d)} A complex of torsion $R$\+modules\/ $\cM^\bu$ is exact
if and only if the complex of torsion $S_\alpha$\+modules
$f_\alpha^*(\cM^\bu)$ is exact for every\/~$\alpha$.
\end{lem}

\begin{proof}
 The categories $R\Tors$ and $S_\alpha\Tors$ are abelian, and
the functors $f_\alpha^*$ are exact by
Lemma~\ref{flat-map-of-adic-rings-direct-image-lemma}(a).
 Hence, in order to prove all the assertions~(a\+-d), it suffices
to check that one has $\cM=0$ whenever $f_\alpha^*(\cM)=0$ for
all~$\alpha$ and some torsion $R$\+module~$\cM$.
 The latter claim follows immediately from
Lemma~\ref{torsion-cohomol-Cech-sequence-lemma}.
\end{proof}

\begin{lem} \label{colocality-of-split-exactness-for-injective}
\textup{(a)} A short sequence of injective torsion $R$\+modules\/
$0\rarrow\cM\rarrow\cL\rarrow\cK\rarrow0$ is (split) exact if and
only if the short sequence of injective torsion $S_\alpha$\+modules\/
$0\rarrow f_\alpha^\diamond(\cM)\rarrow f_\alpha^\diamond(\cL)
\rarrow f_\alpha^\diamond(\cK)\rarrow0$ is (split) exact
for every\/~$\alpha$. \par
\textup{(b)} A morphism of injective torsion $R$\+modules\/
$\cL\rarrow\cK$ is a split epimorphism in $R\Tors^\inj$ if and only if
the induced morphism of injective torsion $S_\alpha$\+modules
$f_\alpha^\diamond(\cL)\rarrow f_\alpha^\diamond(\cK)$ is a split 
epimorphism in $S_\alpha\Tors^\inj$ for every\/~$\alpha$. \par
\textup{(c)} A complex of injective torsion $R$\+modules\/ $\cK^\bu$
is contractible if and only if the complex of injective torsion
$S_\alpha$\+modules $f_\alpha^\diamond(\cK^\bu)$ is contractible
for every\/~$\alpha$.
\end{lem}

\begin{proof}
 This lemma is dual-analogous to
Lemmas~\ref{locality-of-exactness-for-flat-quotseparated}(a,b,d)
and~\ref{locality-of-exactness-for-very-flat-quotseparated}(a,b,d).
 The ``only if'' assertions in all parts~(a\+-c) are obvious.
 Let us prove the ``if''.

 Part~(a): notice that (split) exactness of the short sequences
$0\rarrow f_\alpha^\diamond(\cM)\rarrow f_\alpha^\diamond(\cL)
\rarrow f_\alpha^\diamond(\cK)\rarrow0$ implies (split) exactness
of the short sequences $0\rarrow
f_{\alpha_1,\dotsc,\alpha_k}^\diamond(\cM)\rarrow
f_{\alpha_1,\dotsc,\alpha_k}^\diamond(\cL)\rarrow
f_{\alpha_1,\dotsc,\alpha_k}^\diamond(\cK)\rarrow0$
for all $k\ge1$.
 Now the three-term sequence of injective torsion $R$\+modules
$0\rarrow\cM\rarrow\cL\rarrow\cK\rarrow0$ induces a three-term
sequence of finite split exact
complexes~\eqref{injective-torsion-homol-Cech-sequence-eqn}
from Lemma~\ref{injective-torsion-homol-Cech-sequence-lemma}.
 Since the latter three-term sequence is exact on all the terms
of the finite exact
sequences~\eqref{injective-torsion-homol-Cech-sequence-eqn} except
perhaps the rightmost ones, it follows that the short sequence of
the rightmost terms $0\rarrow\cM\rarrow\cL\rarrow\cK\rarrow0$
is exact, too.

 Part~(b): notice that the map $\bigoplus_\alpha f_\alpha{}_\diamond
f_\alpha^\diamond(\cK)\rarrow\cK$ is a split epimorphism
in $R\Tors^\inj$
by Lemma~\ref{injective-torsion-homol-Cech-sequence-lemma}.
 If the map $f_\alpha^\diamond(\cL)\rarrow f_\alpha^\diamond(\cK)$ is
a split epimorphism in $S_\alpha\Tors^\inj$, then the map
$f_\alpha{}_\diamond f_\alpha^\diamond(\cL)\rarrow
f_\alpha{}_\diamond f_\alpha^\diamond(\cK)$ is a split epimorphism
in $R\Tors^\inj$ for every~$\alpha$ (cf.\
Lemma~\ref{flat-map-of-adic-rings-direct-image-lemma}(b)).
 Hence the composition $\bigoplus_\alpha f_\alpha{}_\diamond
f_\alpha^\diamond(\cL)\rarrow\bigoplus_\alpha f_\alpha{}_\diamond
f_\alpha^\diamond(\cK)\rarrow\cK$ is a split epimorphism in
$R\Tors^\inj$.
 The latter composition is equal to the composition
$\bigoplus_\alpha f_\alpha{}_\diamond
f_\alpha^\diamond(\cL)\rarrow\cL\rarrow\cK$, and it follows that
$\cL\rarrow\cK$ is a split epimorphism in $R\Tors^\inj$, too.

 Part~(c): notice that if the complexes $f_\alpha^\diamond(\cK^\bu)$
are contractible, then so are the complexes
$f_{\alpha_1,\dotsc,\alpha_k}^\diamond(\cK^\bu)$ for all $k\ge1$.
 Applying the construction of the finite split exact
sequence~\eqref{injective-torsion-homol-Cech-sequence-eqn}
from Lemma~\ref{injective-torsion-homol-Cech-sequence-lemma}
to every term of the complex $\cK^\bu$, we obtain a finite
termwise split exact sequence of complexes in $R\Modl^\tors$.
 Since all the complexes appearing in the latter finite termwise
split exact sequence, except perhaps the rightmost one, are split
exact, it follows that the rightmost complex $\cK^\bu$ is split
exact, too.
\end{proof}

\begin{lem} \label{colocality-of-exactness-for-contraadjusted}
\textup{(a)} A short sequence of contraadjusted quotseparated
contramodule $R$\+modules\/ $0\rarrow\fL\rarrow\fM\rarrow\fN\rarrow0$
is exact if and only if the short sequence of contraadjusted 
quotseparated contramodule $S_\alpha$\+modules\/ $0\rarrow
f_\alpha^!(\fL)\rarrow f_\alpha^!(\fM)\rarrow f_\alpha^!(\fN)\rarrow0$
is exact for every\/~$\alpha$. \par
\textup{(b)} A morphism of contraadjusted quotseparated contramodule
$R$\+modules\/ $\fM\rarrow\fN$ is an admissible epimorphism in
$R\Contra^\cta$ if and only if the induced morphism of contraadjusted
quotseparated contramodule $S_\alpha$\+modules $f_\alpha^!(\fM)
\rarrow f_\alpha^!(\fN)$ is an admissible epimorphism in
$S_\alpha\Contra^\cta$ for every\/~$\alpha$. \par
\textup{(c)} A complex of contraadjusted quotseparated contramodule
$R$\+modules\/ $\fM^\bu$ is exact in $R\Contra^\cta$ if and only if
the complex of contraadjusted quotseparated contramodule
$S_\alpha$\+modules $f_\alpha^!(\fM^\bu)$ is exact in
$S_\alpha\Contra^\cta$ for every\/~$\alpha$.
\end{lem}

\begin{proof}
 Part~(a): the ``only if'' assertion holds by
Corollary~\ref{colocalization-contraadjusted-exactness-cor}
(the result of
Proposition~\ref{colocalization-of-contraadjusted-contramodule-prop}
is presumed).
 To prove the ``if'', notice that, by the same corollary,
exactness of the short sequences $0\rarrow f_\alpha^!(\fL)\rarrow
f_\alpha^!(\fM)\rarrow f_\alpha^!(\fN)\rarrow0$ implies exactness of
the short sequences $0\rarrow f_{\alpha_1,\dotsc,\alpha_k}^!(\fL)
\rarrow f_{\alpha_1,\dotsc,\alpha_k}^!(\fM)\rarrow
f_{\alpha_1,\dotsc,\alpha_k}^!(\fN)\rarrow0$ for all $k\ge1$.
 Now the three-term sequence of contraadjusted quotseparated
contramodule $R$\+modules $0\rarrow\fL\rarrow\fM\rarrow\fN\rarrow0$
induces a three-terms sequence of finite exact
sequences~\eqref{contraadj-homol-Cech-sequence-eqn}
from Lemma~\ref{contraadj-homol-Cech-sequence-lemma}.
 Since the latter three-term sequence is exact on all the terms of
the finite exact
sequences~\eqref{contraadj-homol-Cech-sequence-eqn} except perhaps
the rightmost ones, it follows that the short sequence of the rightmost
terms $0\rarrow\fL\rarrow\fM\rarrow\fN\rarrow0$ is exact, too.

 Part~(b): the ``only if'' assertion follow from part~(a).
 To prove the ``if'', notice that the map $\bigoplus_\alpha
f_\alpha{}_\sharp f_\alpha^!(\fN)\rarrow\fN$ is an admissible
epimorphism in $R\Contra^\cta$ by
Lemma~\ref{contraadj-homol-Cech-sequence-lemma}.
 Here we are using the fact that the full subcategory $R\Contra^\cta$
is closed under cokernels (in fact, even under quotients) in $R\Contra$.
 If the map $f_\alpha^!(\fM)\rarrow f_\alpha^!(\fN)$ is an admissible
epimorphism in $S_\alpha\Contra^\cta$, then the map
$f_\alpha{}_\sharp f_\alpha^!(\fM)\rarrow
f_\alpha{}_\sharp f_\alpha^!(\fN)$ is an admissible epimorphism in
$R\Contra^\cta$ by Lemma~\ref{restriction-of-scalars-contraadjusted}.
 Hence the composition $\bigoplus_\alpha
f_\alpha{}_\sharp f_\alpha^!(\fM)\rarrow\bigoplus_\alpha
f_\alpha{}_\sharp f_\alpha^!(\fN)\rarrow\fN$ is an admissible
epimorphism in $R\Contra^\cta$.
 The latter composition is equal to the composition
$\bigoplus_\alpha f_\alpha{}_\sharp f_\alpha^!(\fM)\rarrow\fM
\rarrow\fN$, and it follows that $\fM\rarrow\fN$ is an admissible
epimorphism in $R\Contra^\cta$ by~\cite[Proposition~7.6]{Bueh}.

 Part~(c): the ``only if'' assertion follow from part~(a).
 To prove the ``if'', notice that if the complexes $f_\alpha^!(\fM^\bu)$
are exact in $S_\alpha\Contra^\cta$, then the complexes
$f_{\alpha_1,\dotsc,\alpha_k}^!(\fM^\bu)$ are exact in
$S_{\alpha_1,\dotsc,\alpha_k}\Contra^\cta$ for all $k\ge1$
(by the ``only if'' assertion).
 Hence the complexes $f_{\alpha_1,\dotsc,\alpha_k}{}_\sharp
f_{\alpha_1,\dotsc,\alpha_k}^!(\fM^\bu)$ are exact in $R\Contra^\cta$
(by Lemma~\ref{restriction-of-scalars-contraadjusted}).
 Applying the construction of the finite exact
sequence~\eqref{contraadj-homol-Cech-sequence-eqn}
from Lemma~\ref{contraadj-homol-Cech-sequence-lemma} to every term
of the complex $\fM^\bu$, we obtain a finite termwise exact
sequence of complexes in $R\Contra^\cta$.
 Since all the complexes appearing in the latter finite termwise
exact sequence, except perhaps the rightmost one, are exact
in $R\Contra^\cta$, it follows that the rightmost complex $\fM^\bu$
is exact in $R\Contra^\cta$, too (e.~g.,
by~\cite[Corollary~3.6(i)]{Bueh}).
\end{proof}

\begin{lem} \label{colocality-of-exactness-for-cotorsion}
\textup{(a)} A short sequence of cotorsion quotseparated contramodule
$R$\+modules\/ $0\rarrow\fL\rarrow\fM\rarrow\fN\rarrow0$ is exact if
and only if the short sequence of cotorsion quotseparated contramodule
$S_\alpha$\+modules\/ $0\rarrow f_\alpha^!(\fL)\rarrow f_\alpha^!(\fM)
\rarrow f_\alpha^!(\fN)\rarrow0$ is exact for every\/~$\alpha$. \par
\textup{(b)} A morphism of cotorsion quotseparated contramodule
$R$\+modules\/ $\fM\rarrow\fN$ is an admissible epimorphism in
$R\Contra^\cot$ if and only if the induced morphism of cotorsion
quotseparated contramodule $S_\alpha$\+modules $f_\alpha^!(\fM)
\rarrow f_\alpha^!(\fN)$ is an admissible epimorphism in
$S_\alpha\Contra^\cot$ for every\/~$\alpha$. \par
\textup{(c)} A complex of cotorsion quotseparated contramodule
$R$\+modules\/ $\fM^\bu$ is exact in $R\Contra^\cot$ if and only if
the complex of cotorsion quotseparated contramodule $S_\alpha$\+modules
$f_\alpha^!(\fM^\bu)$ is exact in $S_\alpha\Contra^\cot$
for every\/~$\alpha$.
\end{lem}

\begin{proof}
 This is similar to
Lemma~\ref{colocality-of-exactness-for-contraadjusted}.
 The results of
Proposition~\ref{colocalization-of-cotorsion-contramodule-prop}
and Lemma~\ref{restriction-of-scalars-cotorsion} are relevant.
\end{proof}

 Notice that the obvious analogues of
Lemmas~\ref{colocality-of-split-exactness-for-injective}(b),
\ref{colocality-of-exactness-for-contraadjusted}(b),
and~\ref{colocality-of-exactness-for-cotorsion}(b) for admissible
(or split) monomorphisms are \emph{not} true (not even for
discrete regular Noetherian commutative rings $R$ and~$S_\alpha$
of Krull dimension~$1$).
 The counterexample can be discerned from~\cite[Example~3.2.1]{Pcosh}.

\begin{prop} \label{contraadjusted-and-cotorsion-periodicity-prop}
 Let $R$ be an adic topological ring.
 In this context: \par
\textup{(a)} a complex of contraadjusted quotseparated contramodule
$R$\+modules\/ $\fP^\bu$ is exact in $R\Contra^\cta$ if and only if
it is exact in $R\Contra$; \par
\textup{(b)} a complex of cotorsion quotseparated contramodule
$R$\+modules\/ $\fP^\bu$ is exact in $R\Contra^\cot$ if and only if
it is exact in $R\Contra$.
\end{prop}

\begin{proof}
 Part~(a) is easy: the assertion follows from the fact that the full
subcategory $R\Contra^\cta$ is closed under quotients in $R\Contra$.
 Part~(b) is the cotorsion periodicity theorem for contramodules;
see~\cite[Corollary~12.4]{Pflcc}.
\end{proof}

\subsection{Colocality of flatness of contramodules}
\label{colocality-of-flatness-of-contramods-subsecn}
 This section is our version of the beginning
of~\cite[Section~1.7]{Pcosh}.
 The following two lemmas form a partial generalization
of~\cite[Corollary~1.7.4]{Pcosh}.

\begin{lem} \label{tensor-with-finitely-presented-lemma}
 Let $R$ be an adic topological ring, $\fP$ be a contramodule
$R$\+module, $R$\+module, and $M$ be a finitely presented $R$\+module.
 Then \par
\textup{(a)} the tensor product\/ $\fP\ot_RM$ is a contramodule
$R$\+module; \par
\textup{(b)} if\/ $\fP$ is a quotseparated contramodule $R$\+module,
then the tensor product\/ $\fP\ot_RM$ is a quoseparated contramodule
$R$\+module; \par
\textup{(c)} if\/ $\fP$ is a contraadjusted quotseparated contramodule
$R$\+module, then the tensor product\/ $\fP\ot_RM$ is a contraadjusted
quotseparated contramodule $R$\+module.
\end{lem}

\begin{proof}
 Part~(a) holds because the full subcategory $R\Modl_\ctra$ is closed
under finite direct sums and cokernels in $R\Modl$.
 Part~(b) holds because the full subcategory $R\Contra=R\Modl_\ctra^\qs$
is closed under finite direct sums and cokernels in $R\Modl$.
 Part~(c) holds because the class of contraadjusted quotseparated
contramodule $R$\+modules is closed under finite direct sums and
quotients in $R\Contra$.
\end{proof}

\begin{lem} \label{Hom-from-very-flat-to-flat-contraadjusted}
 Let $R$ be a Noetherian commutative ring endowed with an adic
topology.
 Let\/ $\fF$ be a very flat contramodule $R$\+module, $\fP$ be a flat
contraadjusted contramodule $R$\+module, and $M$ be a finitely
generated $R$\+module.
 Then the natural map of abelian groups/contramodule $R$\+modules\/
$\Hom_R(\fF,\fP)\ot_RM\rarrow\Hom_R(\fF,\>\fP\ot_RM)$ is an isomorphism,
and the contramodule $R$\+module\/ $\Hom_R(\fF,\fP)$ is flat (and
contraadjusted).
\end{lem}

\begin{proof}
 Notice first of all that the $R$\+module $\Hom_R(F,\fQ)$ is
a contramodule for any $R$\+module $F$ and any contramodule
$R$\+module~$\fQ$ (Lemma~\ref{Hom-is-a-contramodule}(b)
or~\ref{Hom-is-a-quotseparated-contramodule}(c)).
 In view of Lemma~\ref{tensor-with-finitely-presented-lemma}(a) or~(b),
both $\Hom_R(\fF,\fP)\ot_RM$ and $\Hom_R(\fF,\>\fP\ot_RM)$ are
contramodule $R$\+modules.
 A natural map $\Hom_R(F,P)\ot_RM\rarrow\Hom_R(F,\>P\ot_RM)$ is
well-defined for any modules $F$, $P$, and $M$ over any commutative
ring~$R$.
 The contramodule $R$\+module $\Hom_R(\fF,\fQ)$ is contraadjusted for
any very flat contramodule $R$\+module $\fF$ and any contraadjusted
contramodule $R$\+module $\fQ$ by Lemma~\ref{vfl-cta-tensor-Hom}(b).

 The category $R\modl$ of finitely generated $R$\+modules $M$
is abelian, since the ring $R$ is Noetherian by assumption.
 By Lemma~\ref{flat-contramodules-characterizations}, any flat
contramodule $R$\+module is flat as an abstract $R$\+module; so
the $R$\+module $\fP$ is flat.
 In view of Lemma~\ref{tensor-with-finitely-presented-lemma}(a) or~(b),
we have an exact functor $\fP\ot_R{-}\,\:R\modl\rarrow R\Contra$.
 By Lemma~\ref{tensor-with-finitely-presented-lemma}(c), the image of
the latter functor is contained in the full subcategory of
contraadjusted contramodule $R$\+modules $R\Contra^\cta\subset
R\Contra$.
 In view of Lemma~\ref{contraadjusted-contramodules-lemma}(2), it
follows that the functor $M\longmapsto\Hom_R(\fF,\>\fP\ot_RM)\:
R\modl\rarrow R\Contra$ is exact.

 On the other hand, the functor $M\longmapsto\Hom_R(\fF,\fP)\ot_RM$
is right exact.
 The map $\Hom_R(\fF,\fP)\ot_RM\rarrow\Hom_R(\fF,\>\fP\ot_RM)$ is
clearly an isomorphism for all finitely generated free (or projective)
$R$\+modules~$M$.
 Since a right exact functor on an abelian category with enough
projective objects is determined by its action on projective objects,
it follows that the map $\Hom_R(\fF,\fP)\ot_RM\rarrow
\Hom_R(\fF,\>\fP\ot_RM)$ is an isomorphism for all finitely generated
$R$\+modules~$M$.
 Thus the functor $M\longmapsto\Hom_R(\fF,\fP)\ot_RM$ is exact on
$R\modl$, and we can conclude that the $R$\+module $\Hom_R(\fF,\fP)$
is flat.
 Hence it is also flat as a contramodule $R$\+module.
\end{proof}

 The next lemma is essentially a particular case
of~\cite[Corollary~1.7.5]{Pcosh}.

\begin{lem} \label{Hom-from-flat-to-flat-cotorsion}
 Let $R$ be a Noetherian commutative ring endowed with an adic
topology.
 Let\/ $\fF$ be a flat contramodule $R$\+module, $\fP$ be a flat
cotorsion contramodule $R$\+module, and $M$ be a finitely generated
$R$\+module.
 Then the contramodule $R$\+module $\fP\ot_RM$ is cotorsion,
the natural map of abelian groups/contramodule $R$\+modules\/
$\Hom_R(\fF,\fP)\ot_RM\rarrow\Hom_R(\fF,\>\fP\ot_RM)$ is an isomorphism,
and the contramodule $R$\+module\/ $\Hom_R(\fF,\fP)$ is flat (and
cotorsion).
\end{lem}

\begin{proof}
 The proof is similar to that of
Lemma~\ref{Hom-from-very-flat-to-flat-contraadjusted}.
 The contramodule $R$\+module $\Hom_R(\fF,\fQ)$ is cotorsion for any
flat contramodule $R$\+module $\fF$ and any cotorsion contramodule
$R$\+module $\fQ$ by Lemma~\ref{flat-cotorsion-injective-tensor-Hom}(b).

 The key step is to prove that the contramodule $R$\+module $\fP\ot_RM$
is cotorsion.
 Indeed, any flat contramodule $R$\+module is flat as an abstract
$R$\+module by Lemma~\ref{flat-contramodules-characterizations}.
 One can also notice that a contramodule $R$\+module is cotorsion if
and only if it is cotorsion as an abstract $R$\+module, by
Lemma~\ref{cotorsion-contramodules-over-Noetherian}; so it remains to
refer to~\cite[Corollary~1.7.5(a)]{Pcosh}.
 Alternatively, one can argue in the way similar to the proof
of~\cite[Lemma~1.7.3 and Corollary~1.7.5(a)]{Pcosh}, staying mostly
within the realm of contramodule $R$\+modules.
 Then one needs to use
Proposition~\ref{contraadjusted-and-cotorsion-periodicity-prop}(b)
together with the fact that any contramodule $R$\+module is
a sub(contra)module of a cotorsion contramodule $R$\+module (a weak
version of Proposition~\ref{flat-cotorsion-pair-in-quotseparated}(b)).

 Then one needs to use
Lemma~\ref{cotorsion-contramodules-lemma}(2), and the rest of
the argument proceeds similarly to the proof of
Lemma~\ref{Hom-from-very-flat-to-flat-contraadjusted}.
\end{proof}

 The next two lemmas form a dual-analogous version of the previous
three ones.

\begin{lem} \label{Hom-from-finitely-generated-lemma}
 Let $R$ be an adic topological ring, $\cN$ be a torsion $R$\+module,
and $M$ be a finitely generated $R$\+module.
 Then the $R$\+module\/ $\Hom_R(M,\cN)$ is torsion.
\end{lem}

\begin{proof}
 The $R$\+module $\Hom_R(M,\cN)$ is a submodule of a finite direct
sum of copies of $\cN$, so the assertion follows from the fact that
the full subcategory $R\Modl_\tors$ is closed under finite direct
sums and subobjects in $R\Modl$.
\end{proof}

\begin{lem} \label{Hom-into-tensor-and-injective-tensored-with-flat}
\textup{(a)} Let $R$ be an adic topological ring, $\fF$ be a flat
contramodule $R$\+module, $\cN$ be a torsion $R$\+module, and $M$ be
a finitely presented $R$\+module.
 Then the natural map of torsion $R$\+modules\/ $\Hom_R(M,\cN)\ot_R\fF
\rarrow\Hom_R(M,\>\cN\ot_R\fF)$ is an isomorphism. \par
\textup{(b)} Let $R$ be a Noetherian commutative ring endowed with
an adic topology, $\fF$ be a flat contramodule $R$\+module, and\/ $\cJ$
be an injective torsion $R$\+module.
 Then the torsion $R$\+module\/ $\cJ\ot_R\fF$ is injective.
\end{lem}

\begin{proof}
 Part~(a): let $Q\rarrow P\rarrow M\rarrow0$ be a finite presentation 
of~$M$ (so $P$ and $Q$ are finitely generated free $R$\+modules).
 Then we have a left exact sequence of torsion $R$\+modules
$0\rarrow\Hom_R(M,\cN)\rarrow\Hom_R(P,\cN)\rarrow\Hom_R(Q,\cN)$,
which remains left exact after the functor ${-}\ot_R\fF$ is applied.

 Part~(b): in view of
Lemmas~\ref{injective-torsion-modules-characterizations}
and~\ref{flat-contramodules-characterizations}, it suffices to
consider the case of the discrete topology on $R$, which is well-known.
\end{proof}

 The following lemma and proposition shed some light on the concept of
a formal open immersion of adic topological rings defined in
Section~\ref{prelim-formal-open-immersions-coverings-subsecn}.

\begin{lem} \label{contramodulation-of-multiplication-lemma}
 Let $f\:R\rarrow S$ be a formal open immersion of adic topological
rings.
 Consider the induced $S$\+$S$\+bimodule (multiplication) map
$S\ot_RS\rarrow S$; it can be viewed an $S$\+module map with
the left action of $S$, or as an $R$\+module map.
 Then the following induced maps are isomorphisms: \par
\textup{(a)} $\Lambda_R(S\ot_RS)\simeq\Lambda_R(S)$; \par
\textup{(b)} $\boL_0\Lambda_R(S\ot_RS)\simeq\boL_0\Lambda_R(S)$; \par
\textup{(c)} $\Delta_R(S\ot_RS)\simeq\Delta_R(S)$; \par
\textup{(d)} $\Lambda_S(S\ot_RS)\simeq\Lambda_S(S)$; \par
\textup{(e)} $\boL_0\Lambda_S(S\ot_RS)\simeq\boL_0\Lambda_S(S)$; \par
\textup{(f)} $\Delta_S(S\ot_RS)\simeq\Delta_S(S)$.
\end{lem}

\begin{proof}
 Part~(c) implies (a) and~(b), in view of the isomorphisms of functors
$\Lambda_R\circ\Delta_R\simeq\Lambda_R$ and $\boL_0\Lambda_R\circ
\Delta_R\simeq\boL_0\Lambda_R$.
 Similarly, part~(f) implies (d) and~(e).
 Let us prove part~(f); the proof of part~(c) is similar.

 The (left) $S$\+module map $S\ot_RS\rarrow S$ is a split epimorphism,
with the map $S\rarrow S\ot_RS$ taking every element $s\in S$ to
the element $s\ot1\in S\ot_RS$ providing the splitting.
 The reduction functor $M\longmapsto M/IM$ transforms the map
$S\ot_RS\rarrow S$ into the map $S/IS\ot_{R/I}S/IS\rarrow S/IS$
for any ideal of definition $I\subset R$.
 Since the ring homomorphism $f$ is a formal open immersion, the ring
homomorphism $R/I\rarrow S/IS$ is an open immersion, hence the map
$S/IS\ot_{R/I}S/IS\rarrow S/IS$ is an isomorphism.

 Denote by $K$ the kernel of the map $S\ot_RS$.
 It follows from the previous paragraph that $IK=K$ for any ideal of
definition $I\subset R$.
 Since the ring homomorphism~$f$ is continuous, one also has
$JK=K$ for any ideal of definition $J\subset S$.

 It remains to point out that $\Delta_S(M)=0$ for any $S$\+module $M$
such that $JM=M$ for some ideal of definition $J\subset S$.
 Indeed, the image of the adjunction map $M\rarrow\Delta_S(M)$ vanishes
by~\cite[Lemma~4.2 and Remark~4.3]{Pcta}, and it follows easily that
$\Delta_S(M)=0$.
 (Cf.\ Lemma~\ref{adic-contramodule-nakayama}(b).) 
\end{proof}

\begin{prop} \label{formal-open-immersion-change-of-scalars}
 Let $f\:R\rarrow S$ be a formal open immersion of adic topological
rings.
 Then \par
\textup{(a)} the restriction-of-scalars functor
$f_\diamond\:S\Modl_\tors\rarrow R\Modl_\tors$ from
Section~\ref{prelim-change-of-scalars-subsecn} is fully faithful; \par
\textup{(b)} the adjunction morphism\/ $\cN\rarrow f^\diamond
f_\diamond(\cN)$ for the pair of adjoint functors
$(f_\diamond,f^\diamond)$ between the abelian
categories $R\Modl_\tors$ and $S\Modl_\tors$, as per
Section~\ref{prelim-change-of-scalars-subsecn}, is an isomorphism
for any torsion $S$\+module\/~$\cN$; \par
\textup{(c)} the adjunction morphism\/ $f^*f_\diamond(\cN)\rarrow\cN$
for the pair of adjoint functors $(f^*,f_\diamond)$ between
the abelian categories $R\Modl_\tors$ and $S\Modl_\tors$, as per
Section~\ref{prelim-change-of-scalars-subsecn}, is an isomorphism
for any torsion $S$\+module\/~$\cN$; \par
\textup{(d)} the restriction-of-scalars functor
$f_\#\:S\Modl_\ctra\rarrow R\Modl_\ctra$ from
Section~\ref{prelim-change-of-scalars-subsecn} is fully faithful; \par
\textup{(e)} the adjunction morphism\/ $f^\#f_\#(\fQ)\rarrow\fQ$
for the pair of adjoint functors $(f^\#,f_\#)$ between the abelian
categories $R\Modl_\ctra$ and $S\Modl_\ctra$, as per
Section~\ref{prelim-change-of-scalars-subsecn}, is an isomorphism
for any contramodule $S$\+module\/~$\fQ$; \par
\textup{(f)} the adjunction morphism\/ $\fQ\rarrow f^!f_\#(\fQ)$
for the pair of adjoint functors $(f_\#,f^!)$ between the abelian
categories $R\Modl_\ctra$ and $S\Modl_\ctra$, as per
Section~\ref{prelim-change-of-scalars-subsecn}, is an isomorphism
for any contramodule $S$\+module\/~$\fQ$; \par
\textup{(g)} the restriction-of-scalars functor
$f_\sharp\:S\Modl_\ctra^\qs\rarrow R\Modl_\ctra^\qs$ from
Section~\ref{prelim-change-of-scalars-subsecn} is fully faithful; \par
\textup{(h)} the adjunction morphism\/ $f^\sharp f_\sharp(\fQ)\rarrow
\fQ$ for the pair of adjoint functors $(f^\sharp,f_\sharp)$ between
the abelian categories $R\Modl_\ctra^\qs$ and $S\Modl_\ctra^\qs$, as per
Section~\ref{prelim-change-of-scalars-subsecn}, is an isomorphism
for any quotseparated contramodule $S$\+module\/~$\fQ$; \par
\textup{(i)} the adjunction morphism\/ $\fQ\rarrow f^!f_\sharp(\fQ)$
for the pair of adjoint functors $(f_\sharp,f^!)$ between the abelian
categories $R\Modl_\ctra^\qs$ and $S\Modl_\ctra^\qs$, as per
Section~\ref{prelim-coextension-of-scalars-quotseparated-subsecn}, is
an isomorphism for any quotseparated contramodule $S$\+module\/~$\fQ$.
\end{prop}

\begin{proof}
 The triple equivalences of properties
(a)~$\Longleftrightarrow$~(b) $\Longleftrightarrow$~(c), \
(d)~$\Longleftrightarrow$~(e) $\Longleftrightarrow$~(f),
and (g)~$\Longleftrightarrow$~(h) $\Longleftrightarrow$~(i) hold
for any flat tight continuous homomorphism of adic topological rings
$f\:R\rarrow S$ by~\cite[Proposition~I.1.3]{GZ}.
 The implication (d)~$\Longrightarrow$~(g) is obvious (since
the functors~$f_\#$ and~$f_\sharp$ agree).
 So it suffices to prove parts~(c) and~(f).

 Part~(c): one has $\cN\simeq S\ot_S\cN$ and $f^*f_\diamond(\cN)=
S\ot_R\cN\simeq S\ot_RS\ot_S\cN$, and the adjunction map
$f^*f_\diamond(\cN)\rarrow\cN$ is induced by the multiplication
map $S\ot_RS\rarrow S$.
 In view of
Corollary~\ref{derived-completion-tensor-with-torsion-isom-cor}, it
follows from Lemma~\ref{contramodulation-of-multiplication-lemma}(d),
(e), or~(f) that the map $S\ot_S\cN\rarrow S\ot_RS\ot_S\cN$ is
an isomorphism.

 Part~(f): one has $\fQ=\Hom_S(S,\fQ)$ and $f^!f_\#(\fQ)=
\Hom_R(S,\fQ)\simeq\Hom_S(S\ot_RS,\>\fQ)$.
 The adjunction map $\fQ\rarrow f^!f_\#(\fQ)$ is induced by
the multiplication map $S\ot_RS\rarrow S$.
 It follows from
Lemma~\ref{contramodulation-of-multiplication-lemma}(f) that the map
$\Hom_S(S,\fQ)\rarrow\Hom_S(S\ot_RS,\>\fQ)$ is an isomorphism.
\end{proof}

 The next corollary is our version
of~\cite[Corollary~1.2.5(b,d)]{Pcosh}.

\begin{cor} \label{formal-open-immersion-reflects-properties}
 Let $f\:R\rarrow S$ be a formal open immersion of adic topological
rings.
 Then \par
\textup{(a)} a torsion $S$\+module\/ $\cJ$ is injective as
a torsion $S$\+module if and only if it is injective as a torsion
$R$\+module; \par
\textup{(b)} a contramodule $S$\+module\/ $\fF$ is flat as
a contramodule $S$\+module if and only if it is flat as
a contramodule $R$\+module; \par
\textup{(c)} a contramodule $S$\+module\/ $\fF$ is very flat as
a contramodule $S$\+module if and only if it is very flat as
a contramodule $R$\+module; \par
\textup{(d)} a quotseparated contramodule $S$\+module\/ $\fQ$ is
contraadjusted as a quotseparated contramodule $S$\+module if and
only if it is contraadjusted as a quotseparated contramodule
$R$\+module; \par
\textup{(e)} a quotseparated contramodule $S$\+module\/ $\fQ$ is
cotorsion as a quotseparated contramodule $S$\+module if and
only if it is cotorsion as a quotseparated contramodule $R$\+module.
\end{cor}

\begin{proof}
 Part~(a): any injective torsion $S$\+module is injective as a torsion
$R$\+module by Lemma~\ref{flat-map-of-adic-rings-direct-image-lemma}(b).
 Conversely, by
Proposition~\ref{formal-open-immersion-change-of-scalars}(b) we have
$\cJ\rarrow f^\diamond f_\diamond(\cJ)$.
 If the torsion $R$\+module $f_\diamond(\cJ)$ is injective, then
the torsion $S$\+module $\cJ$ is injective, since
the functor~$f^\diamond$ takes injectives to injectives (see
Section~\ref{prelim-change-of-scalars-subsecn}).

 Part~(b): any flat contramodule $S$\+module is flat as a contramodule
$R$\+module by Lemma~\ref{flat-map-of-adic-rings-direct-image-lemma}(c).
 Conversely, by
Proposition~\ref{formal-open-immersion-change-of-scalars}(e) we have
$\fF\simeq f^\#f_\#(\fF)$.
 If the contramodule $R$\+module $f_\#(\fF)$ is flat, then
the contramodule $S$\+module $\fF$ is flat, since the functor~$f^\#$
takes flat contramodule $R$\+modules to flat contramodule $S$\+modules
by Corollary~\ref{flat-contramodules-contraextension-of-scalars}.

 Part~(c): any very flat contramodule $S$\+module is very flat as
a contramodule $R$\+module by
Lemma~\ref{very-flat-map-of-adic-rings-direct-image-lemma}(a).
 Conversely, by
Proposition~\ref{formal-open-immersion-change-of-scalars}(e) we have
$\fF\simeq f^\#f_\#(\fF)$.
 If the contramodule $R$\+module $f_\#(\fF)$ is very flat,
then the contramodule $S$\+module $\fF$ is very flat, since
the functor~$f^\#$ takes very flat contramodule $R$\+modules
to very flat contramodule $S$\+modules by
Lemma~\ref{contraextension-of-scalars-very-flat}.

 Part~(d): any contraadjusted quotseparated contramodule $S$\+module is
contraadjusted as a quotseparated contramodule $R$\+module by
Lemma~\ref{restriction-of-scalars-contraadjusted}.
 Conversely, by
Proposition~\ref{formal-open-immersion-change-of-scalars}(i) we have
$\fQ\simeq f^!f_\sharp(\fQ)$.
 If the quotseparated contramodule $R$\+module $f_\sharp(\fQ)$ is
contraadjusted, then the quotseparated contramodule $S$\+module $\fQ$
is contraadjusted, since the functor~$f^!$ takes contraadjusted
quotseparated contramodule $R$\+modules to contraadjusted quotseparated 
contramodule $S$\+modules by
Proposition~\ref{colocalization-of-contraadjusted-contramodule-prop}.

 Part~(e): any cotorsion quotseparated contramodule $S$\+module is
cotorsion as a quotseparated contramodule $R$\+module by
Lemma~\ref{restriction-of-scalars-cotorsion}.
 Conversely, by
Proposition~\ref{formal-open-immersion-change-of-scalars}(i) we have
$\fQ\simeq f^!f_\sharp(\fQ)$.
 If the quotseparated contramodule $R$\+module $f_\sharp(\fQ)$ is
cotorsion, then the quotseparated contramodule $S$\+module $\fQ$ is
cotorsion, since the functor~$f^!$ takes cotorsion quotseparated
contramodule $R$\+modules to cotorsion quotseparated  contramodule
$S$\+modules by
Proposition~\ref{colocalization-of-cotorsion-contramodule-prop}.
\end{proof}

 The following corollary is our version
of~\cite[Corollary~1.7.6(a)]{Pcosh}.

\begin{cor} \label{colocalization-of-flat-contramodules-cor}
 Let $R$ be a Noetherian commutative ring endowed with an adic
topology, and let $f\:R\rarrow S$ be a formal open immersion of adic
topological rings.
 Then the contramodule $S$\+module $f^!(\fP)=\Hom_R(S,\fP)$ is flat
and contraadjusted for any flat contraadjusted contramodule
$R$\+module\/~$\fP$.
\end{cor}

\begin{proof}
 The contramodule $S$\+module $f^!(\fP)$ is contraadjusted by
Proposition~\ref{colocalization-of-contraadjusted-contramodule-prop}.
 Furthermore, we have $\Lambda_S(S)=\Lambda_R(S)$ since the ring
homomorphism~$f$ is tight and continuous, and
$\Lambda_R(S)=\boL_0\Lambda_R(S)$ by
Corollary~\ref{flat-reductions-derived-Lambda-is-underived-cor}
applied to the $R$\+module $F=S$.
 Hence the contramodule $R$\+module $f_\sharp f^!(\fP)=
\Hom_R(S,\fP)\simeq\Hom_R(\Lambda(S),\fP)$ is flat by
Lemma~\ref{Hom-from-very-flat-to-flat-contraadjusted} (applied
to the very flat contramodule $R$\+module $\fF=\Lambda(S)$).
 By Corollary~\ref{formal-open-immersion-reflects-properties}(b),
it follows that the contramodule $S$\+module $f^!(\fP)$ is flat.
\end{proof}

 The next corollary is the dual-analogous version of the previous one.

\begin{cor} \label{localization-of-injective-torsion-modules-cor}
 Let $R$ be a Noetherian commutative ring endowed with an adic
topology, and let $f\:R\rarrow S$ be a formal open immersion of adic
topological rings.
 Then the torsion $S$\+module $f^*(\cJ)=S\ot_R\cJ$ is injective
for any injective torsion $R$\+module\/~$\cJ$.
\end{cor}

\begin{proof}
 We have $\Lambda_S(S)=\Lambda_R(S)$ according to the previous proof,
and $S\ot_R\cJ\simeq\Lambda(S)\ot_R\cJ$ by
Corollary~\ref{derived-completion-tensor-with-torsion-isom-cor}.
 Hence the torsion $R$\+module $f_\diamond f^*(\cJ)=S\ot_R\cJ$ is
injective by
Lemma~\ref{Hom-into-tensor-and-injective-tensored-with-flat}(b)
(applied to the flat contramodule $R$\+module $\fF=\Lambda(S)$).
 By Corollary~\ref{formal-open-immersion-reflects-properties}(a),
it follows that the torsion $S$\+module $f^*(\cJ)$ is injective.
\end{proof}

 The final three theorems of this section are results from adic
commutative algebra whose proofs use the algebraic geometry of
formal schemes.
 They are the adic (or dual-analogous adic) versions
of~\cite[Corollary~1.7.6(b)]{Pcosh}.
 We will prove these theorems in
Section~\ref{injective-projective-secn}.

\begin{thm} \label{locality-of-injective-torsion-modules-thm}
 Let $f_\alpha\:R\rarrow S_\alpha$ be a formal open covering of
an adic topological ring $R$ whose underlying commutative ring $R$
is Noetherian.
 Then a torsion $R$\+module\/ $\cK$ is injective if and only if
the torsion $S_\alpha$\+module $S_\alpha\ot_R\cK$ is injective
for every\/~$\alpha$.
\end{thm}

 The proof of
Theorem~\ref{locality-of-injective-torsion-modules-thm} will be
given below in Section~\ref{injective-torsion-sheaves-subsecn}.

\begin{thm} \label{colocality-of-flat-contramods-if-cotorsion-thm}
 Let $f_\alpha\:R\rarrow S_\alpha$ be a formal open covering of
an adic topological ring $R$ whose underlying commutative ring $R$
is Noetherian.
 Then a cotorsion contramodule $R$\+module\/ $\fP$ is flat if and only
if the cotorsion contramodule $S_\alpha$\+module\/
$\Hom_R(S_\alpha,\fP)$ is flat for every\/~$\alpha$.
\end{thm}

 The proof of
Theorem~\ref{colocality-of-flat-contramods-if-cotorsion-thm} will be
given below in Section~\ref{projective-lct-cosheaves-subsecn}.

\begin{thm} \label{colocality-of-flat-contramods-if-finite-Krull-thm}
 Let $f_\alpha\:R\rarrow S_\alpha$ be a formal open covering of
an adic topological ring $R$ whose underlying commutative ring $R$
is Noetherian.
 Assume that for some (or equivalently, for every) ideal of definition
$I\subset R$, the quotient ring $R/I$ has finite Krull dimension.
 Then a contraadjusted contramodule $R$\+module\/ $\fP$ is flat if and
only if the contraadjusted contramodule $S_\alpha$\+module\/
$\Hom_R(S_\alpha,\fP)$ is flat for every\/~$\alpha$.
\end{thm}

 The proof of
Theorem~\ref{colocality-of-flat-contramods-if-finite-Krull-thm} will be
given below in Section~\ref{flat-cosheaves-of-contramods-subsecn}.

\subsection{Flat cotorsion contramodules over adic Noetherian rings}
\label{flat-cotorsion-contramods-over-adic-Noetherian-subsecn}
 The discussion of flat cotorsion contramodules over adic Noetherian
rings was started in
Lemma~\ref{flat-cotorsion-contramods-in-terms-of-reductions}.
 In this section we formulate the adic versions of Matlis'
classification of injective modules~\cite{Mat} and Enochs'
classification of flat cotorsion modules~\cite{En2} over Noetherian
commutative rings.

 Given an associative ring $R$ and an $R$\+module $M$, we denote
(as usual) by $E_R(M)$ an injective envelope of the $R$\+module~$M$.
 Now let $R$ be a commutative ring and $\fp\subset R$ be a prime ideal.
 Then we denote (also as usual) by $R_\fp=(R\setminus\fp)^{-1}R$
the localization of $R$ at the prime ideal~$\fp$.
 So $R_\fp$ is a local ring and $R_\fp\fp=(R\setminus\fp)^{-1}\fp$
is its maximal ideal.
 Notice that the injective envelope $E_R(R/\fp)$ of the cyclic
$R$\+module $R/\fp$ is a module over~$R_\fp$, and in fact,
an $R_\fp\fp$\+torsion $R_\fp$\+module.

 Moreover, denote the completion of the local ring $R_\fp$ by
$\widehat R_\fp=\Lambda_{R_\fp\fp}(R_\fp)$.
 Then $E_R(R/\fp)$ is also the injective envelope of
the $\widehat R_\fp$\+module $R_\fp/R_\fp\fp\simeq\widehat R_\fp/
\widehat R_\fp\fp$ in the category of $R_\fp$\+modules and in
the category of $\widehat R_\fp$\+modules; so we have
$E_R(R/\fp)\simeq E_{R_\fp}(R_\fp/R_\fp\fp)\simeq
E_{\widehat R_\fp}(\widehat R_\fp/\widehat R_\fp\fp)$.

\begin{prop} \label{injective-torsion-modules-classification-prop}
 Let $R$ be a Noetherian commutative ring endowed with an adic
topology.
 Then the injective objects of the abelian category of torsion
$R$\+modules $R\Tors$ are precisely the (possibly infinite) direct sums
of copies of the $R$\+modules $E_R(R/\fp)$, where\/ $\fp$ ranges over
the prime ideals of $R$ containing an ideal of definition of
the adic topology on~$R$.
\end{prop}

\begin{proof}
 By~\cite[Section~3]{Mat}, the injective $R$\+modules are precisely
the direct sums of copies of the $R$\+modules $E_R(R/\fp)$, where
$\fp$~ranges over the prime ideals of~$R$.
 Let $I$ be an ideal of definition of the adic topology on~$R$.
 By Lemma~\ref{injective-torsion-modules-characterizations}, all
injective $I$\+torsion $R$\+modules are injective as $R$\+modules.
 It remains to find out which of the injective $R$\+modules are
$I$\+torsion.
 Notice that the full subcategory of torsion $R$\+modules $R\Tors$
is closed under direct summands and infinite direct sums in $R\Modl$.
 The final and key observation is that the indecomposable injective
$R$\+module $E_R(R/\fp)$ is $I$\+torsion if and only if $I\subset\fp$.
\end{proof}

\begin{thm} \label{all-contramodules-are-relatively-cotorsion-thm}
 Let $R$ be a commutative ring and $I\subset R$ be a finitely
generated ideal.
 Let $F$ be a flat $R$\+module such that the $R/I$\+module $F/IF$ is
projective.
 Then, for any $I$\+contramodule $R$\+module\/ $\fP$, one has
$\Ext_R^i(F,\fP)=0$ for all $i\ge1$.
\end{thm}

\begin{proof}
 This is~\cite[Lemma~9.4(a) and Theorem~9.5]{Pcta}; see
also~\cite[Proposition~B.10.1]{Pweak}.
\end{proof}

\begin{cor} \label{maximal-ideal-in-Noetherian-ring-cor}
 Let $R$ be a commutative Noetherian ring and\/ $\fm\subset R$ be
a maximal ideal.
 Then \par
\textup{(a)} all\/ $\fm$\+contramodule $R$\+modules are cotorsion
$R$\+modules; \par
\textup{(b)} an\/ $\fm$\+contramodule $R$\+module is flat if and only
if it is projective as an object of the abelian category
$R\Modl_{\fm\dctra}$.
\end{cor}

\begin{proof}
 Notice first of all that an $\fm$\+contramodule $R$\+module is flat
as an $R$\+module if and only if it is flat as an $\fm$\+contramodule
$R$\+module (by Lemma~\ref{flat-contramodules-characterizations}).
 Furthermore, an $\fm$\+contramodule $R$\+module is cotorsion as
an $R$\+module if and only if it is cotorsion as an $\fm$\+contramodule
$R$\+module (by Lemma~\ref{cotorsion-contramodules-over-Noetherian}).
 Finally, all projective contramodules are flat (as contramodules),
according to Section~\ref{prelim-inj-proj-flat-torsion-contra-subsecn}.

 In view of these observations, the desired assertions~(a) and~(b)
are equivalent to each other by
Lemma~\ref{cotorsion-contramodules-lemma}(2)
and Corollary~\ref{flat-contramodules-cor}.
 Part~(a) follows immediately from
Theorem~\ref{all-contramodules-are-relatively-cotorsion-thm};
see~\cite[Theorem~9.3]{Pcta} and/or~\cite[Proposition~1.3.7(a)]{Pcosh}.
 Part~(b) follows immediately from
Lemma~\ref{projective-contramodules-characterizations}(7).
\end{proof}

\begin{lem} \label{projective-contramods-are-free-lemma}
 Let $R$ be a commutative Noetherian ring, $\fm\subset R$ be
a maximal ideal, and\/ $\fR=\Lambda_\fm(R)$ be the\/ $\fm$\+adic
completion of~$R$.
 Then all the projective objects\/ $\fP$ of the abelian category
$R\Modl_{\fm\dctra}^\qs=R\Modl_{\fm\dctra}$ are free (quotseparated)
contramodule $R$\+modules, i.~e., all of them have the form\/
$\fP\simeq\fR[[X]]$ for some sets~$X$.
 In other words, the classes of free and projective contramodule
$R$\+modules coincide.
\end{lem}

\begin{proof}
 This is~\cite[Corollary~10.7]{Pcta}.
\end{proof}

\begin{prop} \label{flat-cotorsion-contramods-classification-prop}
 Let $R$ be a Noetherian commutative ring endowed with an adic
topology with an ideal of definition $I\subset R$.
 Then the flat cotorsion contramodule $R$\+modules are precisely
the $R$\+modules\/ $\fF$ of the form\/ $\fF\simeq\prod_{\fp\supset I}
\widehat R_\fp[[X_\fp]]$.
 Here\/ $\fp$~ranges over the prime ideals\/ $\fp\subset R$ such that
$I\subset\fp$, and $X_\fp$ are arbitrary sets
(so $\widehat R_\fp[[X_\fp]]$ are the free/projective/flat
contramodule $\widehat R_\fp$\+modules).
\end{prop}

\begin{proof}
 The classes of flat, projective, and free contramodule
$\widehat R_\fp$\+modules coincide by
Corollary~\ref{maximal-ideal-in-Noetherian-ring-cor}(b)
and Lemma~\ref{projective-contramods-are-free-lemma}.
 By~\cite{En2} and~\cite[Theorem~1.3.8]{Pcosh}, the flat cotorsion
$R$\+modules are precisely the $R$\+modules of the form
$\prod\widehat R_\fp[[X_\fp]]$, where $\fp$~ranges over the prime
ideals of~$R$.
 By Lemmas~\ref{flat-contramodules-characterizations})
and~\ref{cotorsion-contramodules-over-Noetherian}, a contramodule
$R$\+module is flat (respectively, cotorsion) if and only if it is
flat (resp., cotorsion) as an $R$\+module.
 It remains to find out which of the flat cotorsion $R$\+modules
are $I$\+contramodules.

 Notice that the full subcategory of contramodule $R$\+modules
$R\Contra$ is closed under direct summands and infinite products
in $R\Modl$.
 Now if $I\subset\fp$, then (for any set $X_\fp$) the $R$\+module
$\widehat R_\fp[[X_\fp]]$ is a $\fp$\+contramodule, and any
$\fp$\+contramodule is an $I$\+contramodule.
 Conversely, if $I\not\subset\fp$, then the $R$\+module
$\widehat R_\fp$ is a direct summand of $\widehat R_\fp[[X_\fp]]$
for any nonempty set $X_\fp$, and it remains to show that
$\widehat R_\fp$ is \emph{not} an $I$\+contramodule.
 Indeed, let $s\in I$ be an element such that $s\notin\fp$.
 Then $s$~acts by an automorphism of $\widehat R_\fp$,
hence $\Hom_R(R[s^{-1}],\widehat R_\fp)\simeq\widehat R_\fp\ne0$.
\end{proof}

\subsection{Finite tight morphisms of adic topological rings}
\label{prelim-finite-ring-morphisms-subsecn}
 A ring homomorphism $f\:R\rarrow S$ is said to be \emph{finite} if
$S$ is a finitely generated $R$\+module.
 The following lemma and corollary are simple but illuminating.
 
\begin{lem} \label{finite-maps-of-adic-topological-rings-lemma}
 Let $f\:\fR\rarrow\fS$ be a tight continuous homomorphism of
complete, separated adic topological rings.
 Let\/ $\fI\subset\fR$ be an ideal of definition in\/ $\fR$ and\/
$\fJ=\fS f(\fI)$ be the corresponding ideal of definition in\/~$\fS$.
 Let $s_1$,~\dots, $s_m\in\fS$ be a finite collection of elements
and $\bar s_1$,~\dots, $\bar s_m\in\fS/\fJ\fS$ be their images under
the natural surjective ring homomorphism\/ $\fS\rarrow\fS/\fJ$.
 Then the\/ $\fR$\+module\/ $\fS$ is generated by the elements
$s_1$,~\dots, $s_m$ if and only if the\/ $\fR/\fI$\+module\/
$\fS/\fJ$ is generated by the elements $\bar s_1$,~\dots,~$\bar s_m$.
\end{lem}

\begin{proof}
 Notice that both $\fS$ and the direct sum of $m$~copies of $\fR$
are complete, separated $\fR$\+modules.
 All complete, separated $\fR$\+modules are contramodule $\fR$\+modules
as per Section~\ref{prelim-contramodules-subsecn}.
 Hence it remains to apply
Corollary~\ref{contramodule-nakayama-surjectivity}(b) to
the $\fR$\+module map $g\:\bigoplus_{j=1}^m\fR\rarrow\fS$ given by
the collection of elements $s_j\in\fS$.
\end{proof}

\begin{cor} \label{finite-maps-of-adic-topological-rings-cor}
 Let $f\:\fR\rarrow\fS$ be a tight continuous homomorphism of
complete, separated adic topological rings.
 Then the following conditions are equivalent:
\begin{enumerate}
\item the ring map~$f$ (viewed as a map of abstract commutative rings)
is a finite ring homomorphism;
\item there exists an ideal of definition\/ $\fI\subset\fR$ for which
the induced ring map $\bar f\:\fR/\fI\rarrow\fS/\fS f(\fI)$ is
a finite ring homomorphism;
\item for every ideal of definition\/ $\fI\subset\fR$, the induced
ring map $\bar f\:\fR/\fI\rarrow\fS/\fS f(\fI)$ is a finite ring
homomorphism.
\end{enumerate}
\end{cor}

\begin{proof}
 Follows immediately from
Lemma~\ref{finite-maps-of-adic-topological-rings-lemma}.
\end{proof}

 We will say that a tight continuous map of complete, separated adic
topological rings is \emph{finite} if it satisfies any one of
the equivalent conditions of
Corollary~\ref{finite-maps-of-adic-topological-rings-cor}.

\begin{cor} \label{surjective-maps-of-adic-topological-rings-cor}
 Let $f\:\fR\rarrow\fS$ be a tight continuous homomorphism of
complete, separated adic topological rings.
 Then the following conditions are equivalent:
\begin{enumerate}
\item the ring map~$f$ is surjective;
\item there exists an ideal of definition\/ $\fI\subset\fR$ for which
the induced ring map $\bar f\:\fR/\fI\rarrow\fS/\fS f(\fI)$ is
surjective;
\item for every ideal of definition\/ $\fI\subset\fR$, the induced
ring map $\bar f\:\fR/\fI\rarrow\fS/\fS f(\fI)$ is surjective.
\end{enumerate}
\end{cor}

\begin{proof}
 Apply Lemma~\ref{finite-maps-of-adic-topological-rings-lemma}
to $m=1$ and $s_1=1\in\fS$.
\end{proof}

 The following lemma may help the reader feel more comfortable.

\begin{lem} \label{finitely-generated-modules-are-complete-separated}
\textup{(a)} Let\/ $\fR$ be a complete, separated adic topological ring.
 Then any finitely presented\/ $\fR$\+module is a quotseparated
contramodule\/ $\fR$\+module. \par
\textup{(b)} Let\/ $\fR$ be a complete, separated Noetherian commutative
ring with an adic topology.
 Then any finitely generated\/ $\fR$\+module is separated and
complete.
\end{lem}

\begin{proof}
 Part~(a): the assertion follows from the the fact that the full
subcategory of quotseparated contramodule $\fR$\+modules $\fR\Contra$
contains the free $\fR$\+module with one generator $\fR$ and is closed
under finite direct sums and cokernels in $\fR\Modl$.

 Part~(b): more generally, for any Noetherian commutative ring $R$ with
an adic topology and any finitely generated $R$\+module $M$,
one has $\Lambda_R(M)\simeq\Lambda_R(R)\ot_RM$.
 This is a well-known corollary of the Artin--Rees lemma;
see, e.~g., \cite[Theorem~3(ii) in Section~III.3.4]{Bour}.
\end{proof}

 The next lemma is to be compared with
Lemma~\ref{into-discrete-ring-co-contra-extension-of-scalars}.

\begin{lem} \label{finpres-topol-ring-map-co-contra-extensn-of-scalars}
 Let $f\:R\rarrow S$ be a tight continuous homomorphism of adic
topological rings.
 In this context: \par
\textup{(a)} Assume that $S$ is a finitely generated $R$\+module.
 Then, for every torsion $R$\+module\/ $\cM$, there is a natural
isomorphism of $S$\+modules
$$
 f^\diamond(\cM)\simeq f^!(\cM)=\Hom_R(S,\cM).
$$ \par
\textup{(b)} Assume that $S$ is a finitely presented $R$\+module.
 Then, for every contramodule $R$\+module\/ $\fP$, there is a natural
isomorphism of $S$\+modules
$$
 f^\#(\fP)\simeq f^*(\fP)=S\ot_R\fP.
$$
\end{lem}

\begin{proof}
 Part~(a): according to the proof of
Lemma~\ref{into-discrete-ring-co-contra-extension-of-scalars}(a),
it suffices to show that $\Hom_R(S,\cM)$ is a torsion $S$\+module.
 Indeed, $\Hom_R(S,\cM)$ is a torsion $R$\+module by
Lemma~\ref{Hom-from-finitely-generated-lemma}.
 By Lemma~\ref{tight-reflects-torsion-contra-lemma}(a), it follows
that $\Hom_R(S,\cM)$ is a torsion $S$\+module.

 Part~(b): according to the proof of
Lemma~\ref{into-discrete-ring-co-contra-extension-of-scalars}(b),
it suffices to show that $S\ot_R\fP$ is a contramodule $S$\+module.
 Indeed, $S\ot_R\fP$ is a contramodule $R$\+module by
Lemma~\ref{tensor-with-finitely-presented-lemma}(a).
 By Lemma~\ref{tight-reflects-torsion-contra-lemma}(b), it follows
that $S\ot_R\fP$ is a contramodule $S$\+module.
\end{proof}

\begin{cor}
 Let $f\:R\rarrow S$ be a tight continuous homomorphism of Noetherian
commutative rings with adic topologies.
 Assume that $f$~is a finite ring map.
 For any set $X$, consider the free (quotseparated) contramodule
$R$\+module\/ $\fR[[X]]=\Lambda_R(R[X])$ and the free (quotseparated)
contramodule $S$\+module\/ $\fS[[X]]=\Lambda_S(S[X])$.
 Then there is a natural isomorphism of $S$\+modules\/
$S\ot_R\fR[[X]]\simeq\fS[[X]]$.
\end{cor}

\begin{proof}
 For any continuous homomorphism of adic topological rings
$f\:R\rarrow S$, one has $f^\sharp(\Lambda_R(R[X])\simeq
\Lambda_S(S[X])$, as mentioned in
Section~\ref{prelim-change-of-scalars-subsecn}.
 In the Noetherian case, all contramodule modules are quotseparated,
so there is no difference between~$f^\sharp$ and~$f^\#$.
 It remains to refer to
Lemma~\ref{finpres-topol-ring-map-co-contra-extensn-of-scalars}(b).
\end{proof}

 The following lemma is an adic version of~\cite[Lemma~1.7.7]{Pcosh}.

\begin{lem} \label{contraadjusted-extension-of-scalars-lemma}
 Let $f\:R\rarrow S$ be a surjective tight continuous homomorphism of
Noetherian commutative rings with adic topologies.
 Then \par
\textup{(a)} a (quotseparated) contramodule $S$\+module\/ $\fQ$ is
contraadjusted if and only if $f_\#(\fQ)=f_\sharp(\fQ)$ is
a contraadjusted contramodule $R$\+module; \par
\textup{(b)} the contramodule $S$\+module $f^\#(\fQ)=f^\sharp(\fQ)$ is
contraadjusted for any contraadjusted (quotseparated) contramodule
$R$\+module\/~$\fQ$; \par
\textup{(c)} for any very flat contramodule $R$\+module\/ $\fF$ and
any flat contraadjusted contramodule $R$\+module\/ $\fP$, the natural
morphism of contramodule $S$\+modules $f^\#(\Hom_R(\fF,\fP))\rarrow
\Hom_S(f^\#(\fF),f^\#(\fP))$ is an isomorphism. {\hbadness=1550\par}
\end{lem}

\begin{proof}
 Part~(a): this part does not depend on the Noetherianity assumptions
yet.
 Lemma~\ref{contraadjusted-contramodules-lemma}(4) reduces the desired
assertion to the case of discrete rings $R$ and $S$, which is covered
by~\cite[Lemma~1.7.7(a)]{Pcosh}.

 Part~(b): Lemmas~\ref{contraadjusted-contramodules-lemma}(4)
and~\ref{finpres-topol-ring-map-co-contra-extensn-of-scalars}(b)
reduce the desired assertion to the case of discrete rings $R$ and $S$,
which is covered by~\cite[Lemma~1.7.7(b)]{Pcosh}.

 Part~(c): compare Lemma~\ref{Hom-from-very-flat-to-flat-contraadjusted}
with Lemma~\ref{finpres-topol-ring-map-co-contra-extensn-of-scalars}(b).
\end{proof}

 The next lemma is an adic generalization of~\cite[Lemma~1.7.9]{Pcosh}.
 Part~(b) of this lemma, as well as its original version
in~\cite[Lemma~1.7.9(b)]{Pcosh}, provides an example for
the discussion of \emph{relatively perfect ring homomorphisms}
in the paper~\cite[Section~3]{PBas}.

\begin{lem} \label{cotorsion-extension-of-scalars-lemma}
 Let $f\:R\rarrow S$ be a tight continuous homomorphism of Noetherian
commutative rings with adic topologies.
 Assume that $f$~is a finite ring map.
 Then \par
\textup{(a)} the contramodule $S$\+module $f^\#(\fP)=f^\sharp(\fP)$ is
flat and cotorsion for any flat cotorsion contramodule
$R$\+module\/~$\fP$; \par
\textup{(b)} assuming that the adic topological ring $S$ has finite
adic Krull dimension, a contramodule $S$\+module\/ $\fQ$ is cotorsion
if and only if the contramodule $R$\+module $f_\#(\fQ)=f_\sharp(\fQ)$
is cotorsion; \par
\textup{(c)} for any flat contramodule $R$\+module\/ $\fF$ and any flat
cotorsion contramodule $R$\+module\/ $\fP$, the natural morphism of
contramodule $S$\+modules $f^\#(\Hom_R(\fF,\fP))\rarrow
\Hom_S(f^\#(\fF),f^\#(\fP))$ is an isomorphism.
\end{lem}

\begin{proof}
 Part~(a): Lemmas~\ref{flat-contramodules-characterizations},
\ref{cotorsion-contramodules-over-Noetherian},
and~\ref{finpres-topol-ring-map-co-contra-extensn-of-scalars}(b)
reduce the desired assertion to the case of discrete rings $R$ and $S$,
which is covered by~\cite[Lemma~1.7.9(a)]{Pcosh}.
 Alternatively, one can use
Corollary~\ref{flat-contramodules-contraextension-of-scalars},
Lemma~\ref{reductions-of-co-contra-extension-of-scalars}(b) or~(c),
and Lemma~\ref{flat-cotorsion-contramods-in-terms-of-reductions}.

 Part~(b): the key observation, based on
Proposition~\ref{flat-cotorsion-contramods-classification-prop} above
and~\cite[proof of Lemma~1.7.9(a)]{Pcosh}, is that every flat
cotorsion contramodule $S$\+module $\fG$ is a direct summand of
a flat cotorsion contramodule $S$\+module $f^\#(\fF)=S\ot_R\fF$
for some flat cotorsion contramodule $R$\+module~$\fF$.
 Here it is helpful to notice that, given two ideals of definition
$I\subset R$ and $J=Sf(I)\subset S$, a prime ideal $\fqq\subset S$
contains $J$ if and only if the prime ideal $\fp=f^{-1}(\fqq)\subset R$
contains~$I$.

 The ``only if'' implication in part~(b) holds by
Lemma~\ref{restriction-of-scalars-cotorsion}.
 To prove the ``if'', it suffices to show that
$\Ext^1_{S\Contra}(\fG,\fQ)=0$ for every flat contramodule
$S$\+module~$\fG$ (by Lemma~\ref{cotorsion-contramodules-lemma}(2)).
 By Proposition~\ref{Noetherian-contramodule-Ext-isomorphism-prop},
the $\Ext$ groups computed in the categories $S\Contra$ and $S\Modl$
agree (and so do the $\Ext$ groups computed in the categories
$R\Contra$ and $R\Modl$).

 Furthermore, recall that the cotorsion dimensions of contramodule
$S$\+modules are finite by
Proposition~\ref{adic-raynaud-gruson-prop}(b).
 By Proposition~\ref{flat-cotorsion-pair-in-quotseparated}(b),
it follows that any flat contramodule $S$\+module has a finite
coresolution by flat cotorsion contramodule $S$\+modules (it is helpful
to keep the invariance property of the coresolution dimension in mind,
see Lemma~\ref{co-resol-dim-independent-of-co-resolution-lem}(b)
or Corollary~\ref{intersection-co-resolution-dimension-cor}(b)).
 Hence it suffices to prove that $\Ext_S^n(\fG,\fQ)=0$ for all
flat cotorsion contramodule $S$\+modules $\fG$ and all $n\ge1$.
 Now the $S$\+module $\fG$ is a direct summand of $S\ot_R\fF$ for some
flat (cotorsion contramodule) $R$\+module $\fF$, and
$\Ext_S^n(S\ot_R\fF,\>\fQ)\simeq\Ext_R^n(\fF,\fQ)=0$ by
Lemma~\ref{cotorsion-contramodules-lemma}(3).
 
 Part~(c): compare Lemma~\ref{Hom-from-flat-to-flat-cotorsion} with
Lemma~\ref{finpres-topol-ring-map-co-contra-extensn-of-scalars}(b).
\end{proof}

\begin{lem} \label{co-extension-of-scalars-tensor-product-lemma}
 Let $f\:R\rarrow S$ be a tight continuous homomorphism of adic
topological rings.
 Assume that $S$ is a finitely presented $R$\+module.
 Then \par
\textup{(a)} for any flat contramodule $R$\+module\/ $\fF$ and any
torsion $R$\+module\/ $\cL$, the natural morphism of torsion
$S$\+modules $f^\#(\fF)\ot_S f^\diamond(\cL)\rarrow
f^\diamond(\fF\ot_R\cL)$ is an isomorphism; \par
\textup{(b)} for any flat quotseparated contramodule $R$\+module\/ $\fF$
and any torsion $R$\+module\/ $\cL$, the natural morphism of torsion
$S$\+modules $f^\sharp(\fF)\ot_S f^\diamond(\cL)\rarrow
f^\diamond(\fF\ot_R\cL)$ is an isomorphism.
\end{lem}

\begin{proof}
 The natural morphisms in parts~(a) and~(b) exist without the flatness
assumption on $\fF$, without the tightness assumption on~$f$, and
without the finite presentability assumption on~$S$.
 Such morphisms can be constructed using the natural isomorphisms of
tensor products from
Lemma~\ref{map-of-adic-topol-rings-tensor-product-adjunction}(a\+-b)
and the adjunction of the functors~$f_\diamond$ and~$f^\diamond$.
 The assertion of part~(a) is provable by comparing
Lemma~\ref{Hom-into-tensor-and-injective-tensored-with-flat}(a) with
Lemma~\ref{finpres-topol-ring-map-co-contra-extensn-of-scalars}(a\+-b).
 To deduce part~(b) from part~(a), it suffices to recall the natural
isomorphisms $f^\sharp(\fP)\simeq\boL_0\Lambda_S(f^\#(\fP))$ for
any quotseparated contramodule $R$\+module~$\fP$ (as per
Section~\ref{prelim-change-of-scalars-subsecn}) and
$\boL_0\Lambda_S(Q)\ot_S\cN\simeq Q\ot_S\cN$ for any $S$\+module $Q$
and any torsion $S$\+module~$\cN$ (as per
Corollary~\ref{derived-completion-tensor-with-torsion-isom-cor}).
\end{proof}

\Section{Contraherent Cosheaves of Contramodules}
\label{contraherent-cosheaves-of-contramods-secn}

\subsection{The basics of formal schemes}
\label{basics-of-formal-schemes-subsecn}
 Let $\fX$ be a locally Noetherian formal scheme and $\fO=\fO_\fX$
be its structure sheaf.
 For every affine open formal subscheme $\fU\subset\fX$, the ring
of functions $\fR=\fO_\fX(\fU)$ is a Noetherian commutative ring
endowed with a complete, separated adic topology (in the sense
of Section~\ref{adic-rings-subsecn}).
 The results collected in Section~\ref{prelim-noetherianity-subsecn}
are relevant to the fact that the notion of a (locally) Noetherian
formal scheme is reasonably well-behaved.

 From this point on, by an \emph{adic Noetherian ring} we will mean
a Noetherian commutative ring endowed with adic topology.
 So the adic Noetherian rings are the adic topological rings $R$ whose
underlying abstract ring $R$ is Noetherian.

 An ideal $I\subset R$ is said to be an \emph{ideal of definition}
of a topological commutative ring $R$ if the topology of $R$ is
the $I$\+adic topology.
 This terminology was discussed already in
Section~\ref{adic-rings-subsecn}.
 The ideals of definition $I\subset R$ are characterized by
two properties: for every integer $n\ge1$, the ideal $I^n$~is open 
in~$R$; and the descending sequence $I^n$, \,$n\ge1$, of powers
of the ideal $I$ is a base of open neighborhoods of zero in~$R$.

 \emph{Closed subschemes of definition} $\overline X\subset\fX$
correspond bijectively to coherent sheaves of ideals of definition.
 The latter term means coherent sheaves of ideals $\fI$ on $\fX$
such that, for every affine open formal subscheme $\fU\subset\fX$,
the ideal $\fI(\fU)\subset\fO_\fX(\fU)$ is an ideal of definition
in the adic Noetherian ring $\fR=\fO_\fX(\fU)$.
 In particular, the \emph{minimal closed subscheme of definition}
$\overline X_\mn\subset\fX$ corresponds to the coherent sheaf formed
by the maximal ideals of definition $\fI_\mx\subset\fO_\fX(\fU)$.

 Let $\fV\subset\fX$ be another affine open formal subscheme in $\fX$
such that $\fV\subset\fU$.
 Put $\fS=\fO_\fX(\fV)$.
 Then we have the restriction homomorphism of topological rings
$f\:\fR\rarrow\fS$.
 For every ideal of definition $\fI\subset\fR$, the extension
$\fJ=\fS f(\fI)$ of the ideal $\fI$ to the ring $\fS$ is
an ideal of definition $\fJ\subset\fS$.
 So $f\:\fR\rarrow\fS$ is a tight continuous ring map in the sense
of Section~\ref{adic-rings-subsecn}.
 The discrete quotient ring $\overline S=\fS/\fJ$ is the ring of
regular functions on an affine open subscheme $\Spec\overline S$ in
the affine scheme $\Spec\overline R$, where $\overline R=\fR/\fI$.
 In other words, $f\:\fR\rarrow\fS$ is a formal open immersion of
adic topological rings in the sense of
Section~\ref{prelim-formal-open-immersions-coverings-subsecn}.

 If $\fU=\bigcup_\alpha\fV_\alpha$ is an affine open covering of
an affine formal scheme $\fU$ and $\fI\subset\fR$ is an ideal of
definition in $\fR=\fO(\fU)$, then the collection of affine open
subschemes $\Spec\overline S_\alpha$ is an affine open covering of
the affine scheme $\Spec\overline R$.
 Here $\overline R=\fR/\fI$ and $\overline S_\alpha=
\fS_\alpha/\fJ_\alpha$, where $\fJ_\alpha=\fS_\alpha f_\alpha(\fI)$
denotes the extension of the ideal $\fI$ to the ring
$\fS_\alpha=\fO_\fU(\fV_\alpha)$.
 The notation~$f_\alpha$ stands for the restriction homomorphism
of topological rings $f_\alpha\:\fR\rarrow\fS_\alpha$.
 So the collection of ring maps $f_\alpha\:\fR\rarrow\fS_\alpha$ is
a formal open covering of the adic topological ring $\fR$ in the sense
of Section~\ref{prelim-formal-open-immersions-coverings-subsecn}.

 Conversely, to any complete, separated adic Noetherian ring $\fR$,
the affine Noetherian formal scheme $\Spf\fR$ (the \emph{formal
spectrum} of~$\fR$) is naturally assigned.
 The contravariant functor $\Spf$ establishes an anti-equivalence
between the category of complete, separated adic Noetherian rings
(and continuous ring homomorphisms) and the category of affine
Noetherian formal schemes.

 Under the anti-equivalence of categories from the previous paragraph,
open immersions of affine Noetherian formal schemes $\Spf\fS\rarrow
\Spf\fR$ correspond precisely to formal open immersions of complete,
separated adic Noetherian rings $\fR\rarrow\fS$.
 Affine open coverings of affine Noetherian formal schemes
$\Spf\fS_\alpha\rarrow\Spf\fR$ correspond precisely to formal
open coverings of complete, separated adic Noetherian rings
(by complete, separated adic Noetherian rings).

 The abelian category of coherent sheaves on $\fU=\Spf\fR$ is naturally
equivalent to the abelian category of finitely generated
$\fR$\+modules.
 The closed subschemes of definition $\overline U\subset\fU$ correspond
bijectively to the ideals of definition $\fI\subset\fR$; the rule
is $\overline U=\Spec\fR/\fI$.
 Arbitrary closed subschemes $\overline U'\subset\fU$ correspond
bijectively to the open ideals $\fI'\subset\fR$; once again, the rule
is $\overline U'=\Spec\fR/\fI'$.

 A locally Noetherian formal scheme $\fX$ is said to be Noetherian if
$\fX$ is quasi-compact (i.~e., the underlying topological space of
$\fX$ is quasi-compact).
 A locally Noetherian formal scheme $\fX$ is said to be
\emph{semi-separated} if the intersection of any two affine open
formal subschemes in $\fX$ is an affine open formal subscheme.
 Notice that \emph{every} locally Noetherian formal subschemes $\fX$ is
\emph{quasi-separated} in the sense that the intersection of any two
quasi-compact open formal subschemes in $\fX$ is quasi-compact.

 Let $\fX$ be a locally Noetherian formal scheme and $\overline X
\subset\fX$ be a closed subscheme of definition in~$\fX$.
 Then the formal scheme $\fX$ is affine if and only if the scheme
$\overline X$ is affine.
 Indeed, the underlying topological spaces of $\fX$ and $\overline X$
coincide; hence without loss of generality we can assume $\fX$ and
$\overline X$ to be quasi-compact (i.~e., Noetherian).
 Then the assertion follows from the fact that, for any Noetherian
scheme $X$ with a closed subscheme $\overline X$ such that
the underlying sets of $X$ and $\overline X$ coincide,
the scheme $X$ is affine if and only if $\overline X$ is
affine~\cite[Exercise~III.3.1]{HarAG}, \cite[Lemma Tag~01YQ]{SP}.

 Consequently, a locally Noetherian formal scheme $\fX$ is
semi-separated if and only if it has a semi-separated closed subscheme
of definition $\overline X$, or equivalently, if and only if every
closed subscheme of definition in $\fX$ is semi-separated.

 Given an open covering $\bW$ of a locally Noetherian formal
scheme~$\fX$, an open formal subscheme $\fY\subset\fX$ is said to be
\emph{subordinate to\/~$\bW$} if there exists an open formal subscheme
$\fW\in\bW$ for which $\fY\subset\fW$.
 An open covering $\fX=\bigcup_\alpha\fY_\alpha$ is said to be
\emph{subordinate to\/~$\bW$} if $\fY_\alpha$ is subordinate to $\bW$
for every~$\alpha$.
 More generally, the same terminology is applied to open coverings
of arbitrary topological spaces.

\subsection{Cosheaves of modules over a ringed space}
\label{co-sheaves-of-modules-subsecn}
 Let $(\fX,\fO)$ be a ringed space and $\bB$ be a base of open subsets
in~$\fX$.
 Recall the definitions of a \emph{copresheaf of\/ $\fO$\+modules
on\/~$\bB$} and a \emph{cosheaf of\/ $\fO$\+modules on\/~$\bB$}
from~\cite[Section~2.1]{Pcosh} and~\cite[Section~2.1]{Pdomc}.

 Let us start with (pre)sheaves.
 The set $\bB$ is partially ordered by inclusion, and every poset
can be viewed as a category.
 A \emph{presheaf of abelian groups} $\cM$ on $\bB$ is a contravariant
functor $\cM\:\bB^\sop\rarrow\Ab$ (where $\Ab=\boZ\Modl$ is
the category of abelian groups).
 A \emph{presheaf of\/ $\fO$\+modules} $\cM$ on $\bB$ is a presheaf
of abelian groups such that, for every open subset $\fU\subset\fX$,
\,$\fU\in\bB$, an $\fO(\fU)$\+module structure is defined on
the abelian group~$\cM(\fU)$.
 It is required that the restriction map $\cM(\fU)\rarrow\cM(\fV)$
be an $\fO(\fU)$\+module morphism for every pair of open subsets
$\fV\subset\fU\subset\fX$, \,$\fV$, $\fU\in\bB$.

 A presheaf of abelian groups $\cM$ on $\bB$ is called a \emph{sheaf}
if the following \emph{sheaf axiom} is satisfied.
 Let $\fU\subset\fX$, \,$\fU\in\bB$ be an open subset,
$\fU=\bigcup_\alpha\fV_\alpha$ be an open covering of $\fU$ by open
subsets $\fV_\alpha\in\bB$, and let
$\fV_\alpha\cap\fV_\beta=\bigcup_\gamma\fW_{\alpha\beta\gamma}$ be open
coverings of the intersections $\fV_\alpha\cap\fV_\beta$ by open subsets
$\fW_{\alpha\beta\gamma}\in\bB$.
 Then the short sequence
\begin{equation} \label{topology-base-sheaf-axiom}
 0\lrarrow\cM(\fU)\lrarrow\prod\nolimits_\alpha\cM(\fV_\alpha)
 \lrarrow\prod\nolimits_{\alpha,\beta,\gamma}
 \cM(\fW_{\alpha\beta\gamma})
\end{equation}
should be left exact.
 A presheaf of $\fO$\+modules $\cM$ on $\bB$ is called a \emph{sheaf
of\/ $\fO$\+modules} if the underlying presheaf of abelian groups
of $\cM$ is a sheaf of abelian groups on~$\bB$.

 We denote the additive category of sheaves of $\fO$\+modules on
$\bB$ by $(\bB,\fO)\Sh$.
 In the case of sheaves of $\fO$\+modules defined an \emph{all}
the open subsets of~$\fX$ (i.~e., when the base of open subsets
under consideration consists of all the open subsets in~$\fX$),
we denote the corresponding category of sheaves by $(\fX,\fO)\Sh$.

 The main result in connection with the definition of a sheaf of
$\fO$\+modules on $\fX$ is the following theorem
from~\cite[Section~0.3.2]{EGAI}, \cite[Lemma Tag~009U]{SP},
\cite[Proposition~2.1.3]{Pcosh}, \cite[Theorem~2.1(a)]{Pdomc}:
\emph{the natural restriction functor $(\fX,\fO)\Sh\rarrow
(\bB,\fO)\Sh$ is a category equivalence}.
 So any sheaf of $\fO$\+modules on $\bB$ can be extended to a sheaf
of $\fO$\+modules on $\fX$ in a unique way.

 Now we turn to co(pre)sheaves.
 A \emph{copresheaf of abelian groups} $\fP$ on $\bB$ is
a covariant functor $\fP\:\bB\rarrow\Ab$.
 We denote by $\fP[\fU]$ the abelian group that the copresheaf $\fP$
assigns to an open subset $\fU\subset\fX$, \,$\fU\in\bB$.
 The group $\fP[\fU]$ is called the group of \emph{cosections} of
the copresheaf $\fP$ over an open subset~$\fU$.
 The abelian group map $\fP[\fV]\rarrow\fP[\fU]$ that the copresheaf
$\fP$ assigns to an identity inclusion $\fV\subset\fU$ of open
subsets $\fV$, $\fU\in\bB$ is called the \emph{corestriction} map.

 A \emph{copresheaf of\/ $\fO$\+modules} $\fP$ on $\bB$ is
a copresheaf of abelian groups such that, for every open subset
$\fU\subset\fX$, \,$\fU\in\bB$, an $\fO(\fU)$\+module structure is
defined on the abelian group $\fP[\fU]$.
 It is required that the corestriction map $\fP[\fV]\rarrow\fP[\fU]$
be an $\fO(\fU)$\+module morphism for every pair of open subsets
$\fV\subset\fU\subset\fX$, \,$\fV$, $\fU\in\bB$.

 A copresheaf of abelian groups $\fP$ on $\bB$ is called
a \emph{cosheaf} if the following \emph{cosheaf axiom} is satisfied.
 As above, let $\fU\subset\fX$, \,$\fU\in\bB$ be an open subset,
$\fU=\bigcup_\alpha\fV_\alpha$ be an open covering of $\fU$ by open
subsets $\fV_\alpha\in\bB$, and let
$\fV_\alpha\cap\fV_\beta=\bigcup_\gamma\fW_{\alpha\beta\gamma}$ be open
coverings of the intersections $\fV_\alpha\cap\fV_\beta$ by open subsets
$\fW_{\alpha\beta\gamma}\in\bB$.
 Then the short sequence
\begin{equation} \label{topology-base-cosheaf-axiom}
 \bigoplus\nolimits_{\alpha,\beta,\gamma}\fP[\fW_{\alpha\beta\gamma}]
 \lrarrow\bigoplus\nolimits_\alpha\fP[\fV_\alpha]\lrarrow\fP[\fU]
 \lrarrow0
\end{equation}
should be right exact.
 A copresheaf of $\fO$\+modules $\fP$ on $\bB$ is called a \emph{cosheaf
of\/ $\fO$\+modules} if the underlying copresheaf of abelian groups
of $\fP$ is a cosheaf of abelian groups on~$\bB$.

 We denote the additive category of cosheaves of $\fO$\+modules on
$\bB$ by $(\bB,\fO)\Cosh$.
 In the case of cosheaves of $\fO$\+modules defined an \emph{all}
the open subsets of~$\fX$ (i.~e., when the base of open subsets
under consideration consists of all the open subsets in~$\fX$),
we denote the corresponding category of cosheaves by $(\fX,\fO)\Cosh$.

 The main result in connection with the definition of a cosheaf of
$\fO$\+modules on $\fX$ is the following theorem
from~\cite[Theorem~2.1.2]{Pcosh}, \cite[Theorem~2.1(b)]{Pdomc}:
\emph{the natural restriction functor $(\fX,\fO)\Cosh\rarrow
(\bB,\fO)\Cosh$ is a category equivalence}.
 So any cosheaf of $\fO$\+modules on $\bB$ can be extended to a cosheaf
of $\fO$\+modules on $\fX$ in a unique way.

 It is clear immediately from the definitions that infinite products
exist in the category of sheaves $(\fX,\fO)\Sh$, while infinite
direct sums exist in the category of cosheaves $(\fX,\fO)\Cosh$.
 The functors of sections over all open subsets preserve the infinite
products in $(\fX,\fO)\Sh$, while the functors of cosections over
all open subsets preserve the infinite direct sums in $(\fX,\fO)\Cosh$.

\begin{rems} \label{on-co-sheaf-axioms-remark}
 (1)~The sheaf axiom~\eqref{topology-base-sheaf-axiom} implies,
in particular, injectivity of the maps $\cM(\fU)\rarrow
\prod_\alpha\cM(\fV_\alpha)$.
 It follows from the latter property that one need not check
axiom~\eqref{topology-base-sheaf-axiom} for all the choices of
the open coverings $\fV_\alpha\cap\fV_\beta=
\bigcup_\gamma\fW_{\alpha\beta\gamma}$ of the intersections
$\fV_\alpha\cap\fV_\beta$.
 It suffices to check~\eqref{topology-base-sheaf-axiom} for any one
given choice of the open coverings of the pairwise intersections of
the open subsets~$\fV_\alpha$.

 Dual-analogously, the cosheaf
axiom~\eqref{topology-base-cosheaf-axiom} implies, in particular,
surjectivity of the maps $\bigoplus_\alpha\fP[\fV_\alpha]\rarrow
\fP[\fU]$.
 It follows from the latter property that it suffices to check
axiom~\eqref{topology-base-cosheaf-axiom} for any one given choice
of the open coverings $\fV_\alpha\cap\fV_\beta=
\bigcup_\gamma\fW_{\alpha\beta\gamma}$ of the intersections
$\fV_\alpha\cap\fV_\beta$.

\smallskip
 (2)~Assume that the base of open subsets $\bB$ in $\fX$ consists of
quasi-compact open subsets, and furthermore, the intersection of
any two open subsets from $\bB$ that are contained in a third open
subset from $\bB$ is also quasi-compact.
 Then, following~\cite[Remark~2.1.4]{Pcosh}, it suffices to check
axioms~\eqref{topology-base-sheaf-axiom}
and~\eqref{topology-base-cosheaf-axiom} for \emph{finite} open
coverings $\fU=\bigcup_\alpha\fV_\alpha$ and finite open coverings
of the intersections $\fV_\alpha\cap\fV_\beta=
\bigcup_\gamma\fW_{\alpha\beta\gamma}$.

\smallskip
 (3)~As a corollary of~(2), under the assumptions of~(2) one can
immediately see that the category of sheaves $(\fX,\fO)\Sh\simeq
(\bB,\fO)\Sh$ has infinite direct sums, while the category of
cosheaves $(\fX,\fO)\Cosh\simeq(\bB,\fO)\Cosh$ has infinite products.
 The functors of sections over open subsets belonging to $\bB$
preserve the infinite direct sums in $(\fX,\fO)\Sh$, while
the functors of cosections over open subsets belonging to $\bB$
preserve the infinite products in $(\fX,\fO)\Cosh$.
\end{rems}

 Let $(\fX,\fO)$ be a ringed space, $\fW\subset\fX$ be an open
subset, $\cM$ be a presheaf of $\fO$\+modules on $\fX$, and
$\fP$ be a copresheaf of $\fO$\+modules on~$\fX$.
 Then the notation $\cM|_\fW$ and $\fP|_\fW$ stands for
the restrictions of $\cM$ and $\fP$ to open subsets $\fU\subset\fX$
contained in $\fW$, i.~e., $\fU\subset\fW$.
 So $\cM$ is a presheaf of $\fO|_\fW$\+modules on $\fW$ and
$\fP|_\fW$ is a copresheaf of $\fO|_\fW$\+modules on~$\fW$.
 The restriction of a (co)sheaf of $\fO$\+modules on $\fX$ to $\fW$
is a (co)sheaf of $\fO|_\fW$\+modules on~$\fW$.

 Given a presheaf of $\fO$\+modules $\cM$ on a ringed space $(\fU,\fO)$
and a finite open covering $\fU=\bigcup_{\alpha=1}^N\fV_\alpha$ of
$\fU$, we denote by $C^\bu(\{\fV_\alpha\},\cM)$ the finite cohomological
\v Cech complex of abelian groups (in fact, $\fO(\fU)$\+modules)
\begin{multline} \label{Cech-complex-of-abelian-groups-for-presheaf}
 0\lrarrow
 \bigoplus\nolimits_{\alpha=1}^N\cM(\fV_\alpha)\lrarrow
 \bigoplus\nolimits_{1\le\alpha<\beta\le N}
 \cM(\fV_\alpha\cap\fV_\beta) \\ \lrarrow\dotsb\lrarrow
 \cM(\fV_1\cap\dotsb\cap\fV_N)\lrarrow0.
\end{multline}
 The sheaf axiom~\eqref{topology-base-sheaf-axiom} (for finite
open coverings $\fU=\bigcup_\alpha\fV_\alpha$) essentially says
that the natural map of $\fO(\fU)$\+modules
$\cM(\fU)\rarrow H^0C^\bu(\{\fV_\alpha\},\cM)$ is an isomorphism.

 Dual-analogously, given a copresheaf of $\fO$\+modules $\fP$ on
$(\fU,\fO)$, we denote by $C_\bu(\{\fV_\alpha\},\fP)$ the finite
homological \v Cech complex of abelian groups/$\fO(\fU)$\+modules
\begin{multline} \label{Cech-complex-of-abelian-groups-for-copresheaf}
 0\lrarrow\fP[\fV_1\cap\dotsb\cap\fV_N]\lrarrow\dotsb \\
 \lrarrow\bigoplus\nolimits_{1\le\alpha<\beta\le N}
 \fP[\fV_\alpha\cap\fV_\beta]\lrarrow
 \bigoplus\nolimits_{\alpha=1}^N\fP[\fV_\alpha]\lrarrow0.
\end{multline}
 The cosheaf axiom~\eqref{topology-base-cosheaf-axiom} (for finite
open coverings $\fU=\bigcup_\alpha\fV_\alpha$) essentially says
that the natural map of $\fO(\fU)$\+modules
$H_0C^\bu(\{\fV_\alpha\},\fP)\rarrow\fP[\fU]$ is an isomorphism.

\subsection{Cosheaves of contramodules over a formal scheme}
\label{cosheaves-of-contramodules-subsecn}
 For any Noetherian commutative ring $R$ with an adic topology,
we denote by $R\Tors=R\Modl_\tors\subset R\Modl$ the full subcategory
of torsion $R$\+modules and by $R\Contra=R\Modl_\ctra^\qs=R\Modl_\ctra
\subset R\Modl$ the full subcategory of contramodule $R$\+modules.
 All contramodule $R$\+modules are quotseparated due to the assumption
that the ring $R$ is Noetherian.
 See the notation and discussion in
Sections~\ref{prelim-torsion-modules-subsecn},
\ref{prelim-contramodules-subsecn},
and~\ref{prelim-veryflat-and-contraadjusted-subsecn}.

 Let $\fX$ be a locally Noetherian formal scheme.
 A sheaf of $\fO_\fX$\+modules $\cM$ on $\fX$ is said to be
a \emph{sheaf of torsion\/ $\fO_\fX$\+modules} if, for every affine
open formal subscheme $\fU\subset\fX$, the $\fO_\fX(\fU)$\+module
$\cM(\fU)$ is torsion/discrete as a module over a topological
commutative ring with an adic topology.
 This means that the annihilator of any element of $\cM(\fU)$ is
an open ideal in $\fO_\fX(\fU)$, or equivalently, for every ideal of
definition $\fI\subset\fO_\fX(\fU)$ and every element $x\in\cM(\fU)$
there exists an integer $n\ge1$ such that $\fI^nx=0$ in $\cM(\fU)$.

 The sheaf axiom for $\cM$ implies that it suffices to check
the torsion condition for affine open formal subschemes $\fU\subset\fX$
subordinate to a given open covering of~$\fX$.
 Indeed, if $\fU=\bigcup_\alpha\fV_\alpha$ is a finite affine open
covering of an affine formal scheme $\fU$ and $\cN$ is a sheaf of
$\fO_\fU$\+modules on $\fU$, then $\cN(\fU)$ is a torsion module
over $\fR=\fO_\fU(\fU)$ whenever $\cN(\fV_\alpha)$ is a torsion module
over $\fS_\alpha=\fO_\fU(\fV_\alpha)$ for all indices~$\alpha$.
 The latter assertion follows from the facts that the map
$\cN(\fU)\rarrow\bigoplus_\alpha\cN(\fV_\alpha)$ is injective,
the functors of restrictions of scalars $f_\alpha{}_*\:\fS_\alpha\Modl
\rarrow\fR\Modl$ take $\fS_\alpha\Tors$ into $\fR\Tors$ (see
Section~\ref{prelim-change-of-scalars-subsecn}), and any submodule
of a torsion $\fR$\+module is torsion.

 A cosheaf of $\fO_\fX$\+modules $\fP$ on $\fX$ is said to be
a \emph{cosheaf of contramodule\/ $\fO_\fX$\+mod\-ules} if, for every
affine open formal subscheme $\fU\subset\fX$, the $\fO_\fX(\fU)$\+module
$\fP[\fU]$ is a contramodule over the adic topological ring
$\fO_\fX(\fU)$.
 We refer to Section~\ref{prelim-contramodules-subsecn} for
the definition of contramodule $R$\+modules for an adic topological
ring~$R$.
 Let us only point out that, since the commutative ring $\fO_\fX(\fU)$
is Noetherian by assumption, all the contramodule
$\fO_\fX(\fU)$\+modules are quotseparated.

 The cosheaf axiom for $\fP$ implies that is suffices to check
the contramoduleness condition for affine open formal subschemes
$\fU\subset\fX$ subordinate to a given open covering of~$\fX$.
 Indeed, let $\fU=\bigcup_\alpha\fV_\alpha$ be a finite affine open
covering of an affine formal scheme~$\fU$.
 Let $\fQ$ be a cosheaf of $\fO_\fU$\+modules on~$\fU$.
 Then the $\fO_\fU(\fU)$\+module $\fQ[\fU]$ is a contramodule over
the adic topological ring $\fR=\fO_\fU(\fU)$ whenever
the $\fO_\fU(\fV_\alpha)$\+module $\fQ[\fV_\alpha]$ is a contramodule
over the adic topological ring $\fS_\alpha=\fO_\fU(\fV_\alpha)$
\emph{and} the $\fO_\fU(\fV_{\alpha,\beta})$\+module
$\fQ[\fV_{\alpha,\beta}]$ is a contramodule over the adic topological
ring $\fS_{\alpha,\beta}=\fO_\fU(\fV_{\alpha,\beta})$ for all indices
$\alpha$ and~$\beta$.
 Here we put $\fV_{\alpha,\beta}=\fV_\alpha\cap\fV_\beta\subset\fU$.
 
 The latter assertion follows from right exactness of the short
sequence~\eqref{topology-base-cosheaf-axiom},
$$
 \bigoplus\nolimits_{\alpha,\beta}\fQ[\fV_{\alpha,\beta}]
 \lrarrow\bigoplus\nolimits_\alpha\fQ[\fV_\alpha]\lrarrow\fQ[\fU]
 \lrarrow0,
$$
together with the facts that the functors of restriction of scalars
$\fS_\alpha\Modl\rarrow\fR\Modl$ and $\fS_{\alpha,\beta}\Modl\rarrow
\fR\Modl$ take $\fS_\alpha\Contra$ and $\fS_{\alpha,\beta}\Contra$
into $\fR\Contra$ (see Section~\ref{prelim-change-of-scalars-subsecn}),
and that the full subcategory of contramodule $\fR$\+modules
$\fR\Contra$ is closed under cokernels in $\fR\Modl$.

\subsection{Quasi-coherent torsion sheaves}
\label{qcoh-torsion-sheaves-subsecn}
 Let $\fX$ be a locally Noetherian formal scheme with an open
covering~$\bW$.
 Denote by $\bB$ the base of open subsets in $\fX$ consisting of
all the affine open formal subschemes $\fU\subset\fX$ subordinate
to~$\bW$.

 Let $\cM$ be a presheaf of $\fO_\fX$\+modules on~$\bB$.
 We will say that $\cM$ is a \emph{quasi-coherent torsion presheaf}
if the following conditions hold:
\begin{enumerate}
\renewcommand{\theenumi}{\roman{enumi}}
\item For every affine open formal subscheme $\fU\subset\fX$
subordinate to $\bW$, the $\fO_\fX(\fU)$\+module $\cM(\fU)$ is torsion
(as a module over the adic topological ring~$\fO_\fX(\fU)$).
\item For every pair of affine open formal subschemes
$\fV\subset\fU\subset\fX$ subordinate to $\bW$,
the $\fO_\fX(\fV)$\+module map
$$
 \fO_\fX(\fV)\ot_{\fO_\fX(\fU)}\cM(\fU)\lrarrow\cM(\fV)
$$
corresponding by adjunction to the restriction map
$\cM(\fU)\rarrow\cM(\fV)$ is an isomorphism.
\end{enumerate}
 We will call condition~(i) is the \emph{torsion axiom}, while
condition~(ii) will be called the \emph{quasi-coherence axiom}.
 This definition is to be compared with the Enochs--Estrada
description of quasi-coherent sheaves on schemes~\cite[Section~2]{EE},
\cite[Section~2.2]{Pdomc}.

\begin{lem} \label{quasi-coherence-implies-sheaf-axiom}
 The torsion and quasi-coherence axioms imply the sheaf axiom:
 Any quasi-coherent torsion presheaf on\/ $\bB$ is a sheaf of\/
$\fO_\fX$\+modules on\/~$\bB$.
\end{lem}

\begin{proof}
 The assertion follows from exactness of the leftmost fragment of
sequence~\eqref{torsion-cohomol-Cech-sequence-eqn}
in Lemma~\ref{torsion-cohomol-Cech-sequence-lemma} together with
Remarks~\ref{on-co-sheaf-axioms-remark}(1\+-2).

 In fact, under our assumption of local Noetherianity of $\fX$,
one need not use the torsion axiom for this argument, and
the quasi-coherence axiom alone is sufficient.
 Inceed, the finite exact
sequence~\eqref{for-the-rings-cohomol-Cech-sequence-eqn} from
Corollary~\ref{for-the-rings-cohomol-Cech-sequence-cor} is
an exact sequence of flat $R$\+modules by
Lemma~\ref{flat-contramodules-characterizations}(6) or~(7),
so it remains exact after taking the tensor product with
an arbitrary $R$\+module $\cM(\fU)$ for $R=\fO_\fX(\fU)$.
\end{proof}

 So any quasi-coherent torsion presheaf $\cM$ on $\bB$ is, in fact,
a sheaf on $\bB$, and we will speak of $\cM$ as a \emph{quasi-coherent
torsion sheaf} on~$\bB$.
 By Lemma~\ref{quasi-coherence-implies-sheaf-axiom} and in view of
the discussion in Section~\ref{co-sheaves-of-modules-subsecn}, any
quasi-coherent torsion sheaf on $\bB$ can be extended to
a sheaf of $\fO_\fX$\+modules on the whole formal scheme $\fX$ in
a unique way.
 Moreover, according to the discussion in
Section~\ref{cosheaves-of-contramodules-subsecn}, the resulting
sheaf is a sheaf of torsion $\fO_\fX$\+modules on~$\fX$.

 The following proposition establishes the fact that the notion of
a quasi-coherent torsion sheaf is local (i.~e., does not depend on
the open covering~$\bW$).

\begin{prop} \label{torsion-quasi-coherence-is-local-prop}
 Denote by\/ $\bW_\fX$ the trivial open covering\/ $\{\fX\}$ of
the formal scheme\/ $\fX$, and let\/ $\bB_\fX$ be the base of open
subsets in\/ $\fX$ consisting of all affine open formal subschemes
of\/~$\fX$.
 Let\/ $\cM$ be a sheaf of\/ $\fO_\fX$\+modules on\/~$\bB_\fX$.
 Assume that the restriction of\/ $\cM$ to\/ $\bB$ is a quasi-coherent
torsion sheaf on\/~$\bB$.
 Then\/ $\cM$ is a quasi-coherent torsion sheaf on\/~$\bB_\fX$.
\end{prop}

\begin{proof}
 Let $\fV\subset\fU\subset\fX$ be two affine open formal subschemes,
and let $\fU=\bigcup_{\alpha=1}^N\fW_\alpha$ be a finite open
covering of $\fU$ by affine open formal subschemes $\fW_\alpha$
subordinate to~$\bW$.
 The argument from Section~\ref{cosheaves-of-contramodules-subsecn}
explains that $\cM(\fU)$ is a torsion module over $\fO_\fX(\fU)$
whenever $\cM(\fW_\alpha)$ are torsion modules over
$\fO_\fX(\fW_\alpha)$ for all~$\alpha$.
 So the torsion axiom~(i) is satisfied for the sheaf of
$\fO_\fX$\+modules $\cM$ on~$\bB_\fX$.
 We need to check the quasi-coherence axiom~(ii) for $\cM$ with
respect to $\fV\subset\fU$.

 Notice that $\fV=\bigcup_{\alpha=1}^N(\fV\cap\fW_\alpha)$ is
an finite open covering of $\fV$ by affine open formal subschemes
subordinate to~$\bW$.
 By the sheaf axiom~\eqref{topology-base-sheaf-axiom}, we have
a left exact sequence of torsion $\fO_\fX(\fU)$\+modules
\begin{equation} \label{sheaf-axiom-covering-of-U-by-W-alpha}
 0\lrarrow\cM(\fU)\lrarrow\bigoplus\nolimits_{\alpha=1}^N\cM(\fW_\alpha)
 \lrarrow\bigoplus\nolimits_{1\le\alpha<\beta\le N}
 \cM(\fW_\alpha\cap\fW_\beta)
\end{equation}
and a left exact sequence of torsion $\fO_\fX(\fV)$\+modules
\begin{equation} \label{sheaf-axiom-covering-of-V-by-V-cap-W-alpha}
 0\lrarrow\cM(\fV)\lrarrow\bigoplus\nolimits_{\alpha=1}^N
 \cM(\fV\cap\fW_\alpha)\lrarrow
 \bigoplus\nolimits_{1\le\alpha<\beta\le N}
 \cM(\fV\cap\fW_\alpha\cap\fW_\beta).
\end{equation}

 Recall that $\fO_\fX(\fV)$ is a flat contramodule
$\fO_\fX(\fU)$\+module by Lemma~\ref{formal-open-immersions-lemma}(3)
(and a flat $\fO_\fX(\fU)$\+module by
Lemma~\ref{flat-contramodules-characterizations}(6) or~(7)).
 Hence the tensor product functor $\fO_\fX(\fV)\ot_{\fO_\fX(\fU)}{-}$
preserves exactness of the left exact
sequence~\eqref{sheaf-axiom-covering-of-U-by-W-alpha}, and we
obtain a left exact sequence of $\fO_\fX(\fV)$\+modules
\begin{multline} \label{sheaf-axiom-tensored-with-O(V)-sequence}
 0\lrarrow\fO_\fX(\fV)\ot_{\fO_\fX(\fU)}\cM(\fU)
 \lrarrow\bigoplus\nolimits_{\alpha=1}^N
 \fO_\fX(\fV)\ot_{\fO_\fX(\fU)}\cM(\fW_\alpha) \\
 \lrarrow\bigoplus\nolimits_{1\le\alpha<\beta\le N}
 \fO_\fX(\fV)\ot_{\fO_\fX(\fU)}\cM(\fW_\alpha\cap\fW_\beta).
\end{multline}

 The restriction maps in the (pre)sheaf $\cM$ induce a natural
morphism from the left exact
sequence~\eqref{sheaf-axiom-tensored-with-O(V)-sequence}
to the left exact
sequence~\eqref{sheaf-axiom-covering-of-V-by-V-cap-W-alpha}.
 It remains to refer to the next
Lemma~\ref{intersection-of-affine-open-tensor-with-discrete}
to the effect that the map
from~\eqref{sheaf-axiom-tensored-with-O(V)-sequence}
to~\eqref{sheaf-axiom-covering-of-V-by-V-cap-W-alpha} is an isomorphism
on the middle and rightmost terms.
 Hence the natural map $\fO_\fX(\fV)\ot_{\fO_\fX(\fU)}\cM(\fU)\rarrow
\cM(\fV)$ of the leftmost terms in an isomorphism, too.
\end{proof}

\begin{lem} \label{intersection-of-affine-open-tensor-with-discrete}
 Let\/ $\cM$ be a quasi-coherent torsion sheaf on\/ $\bB$, let
$\fV\subset\fU\subset\fX$ be affine open formal subschemes, and let\/
$\fW\subset\fU$ be an affine open formal subscheme subordinate
to\/~$\bW$.
 Then the restriction map of\/ $\fO_\fX(\fW)$\+modules\/
$\cM(\fW)\rarrow\cM(\fV\cap\fW)$ induces an isomorphism of\/
$\fO_\fX(\fV)$\+modules\/ $\fO_\fX(\fV)\ot_{\fO_\fX(\fU)}\cM(\fW)
\rarrow\cM(\fV\cap\fW)$.
\end{lem}

\begin{proof}
 Both $\fW$ and $\fV\cap\fW$ are affine open formal subschemes in $\fX$
subordinate to $\bW$, so we have $\fW\in\bB$ and $\fV\cap\fW\in\bB$.
 By the quasi-coherence axiom~(ii), the restriction map
$\cM(\fW)\rarrow\cM(\fV\cap\fW)$ induces an isomorphism of
$\fO_\fX(\fV\cap\fW)$\+modules
$$
 \fO_\fX(\fV\cap\fW)\ot_{\fO_\fX(\fW)}\cM(\fW)\simeq\cM(\fV\cap\fW).
$$
 There is also an obvious isomorphism
$$
 \fO_\fX(\fV)\ot_{\fO_\fX(\fU)}\cM(\fW)\simeq
 \bigl(\fO_\fX(\fV)\ot_{\fO_\fX(\fU)}\fO_\fX(\fW)\bigr)
 \ot_{\fO_\fX(\fW)}\cM(\fW).
$$
 It remains to show that the map
$$
 \bigl(\fO_\fX(\fV)\ot_{\fO_\fX(\fU)}\fO_\fX(\fW)\bigr)
 \ot_{\fO_\fX(\fW)}\cM(\fW)
 \lrarrow\fO_\fX(\fV\cap\fW)\ot_{\fO_\fX(\fW)}\cM(\fW)
$$
induced by the multiplication map of rings
\begin{equation} \label{affine-open-formal-intersect-multiplicat-map}
 \fO_\fX(\fV)\ot_{\fO_\fX(\fU)}\fO_\fX(\fW)\lrarrow
 \fO_\fX(\fV\cap\fW)
\end{equation}
is an isomorphism.
 Now the map~\eqref{affine-open-formal-intersect-multiplicat-map}
becomes an isomorphism after the adic completion functor $\Lambda$
is applied, while the $\fO_\fX(\fW)$\+module $\cM(\fW)$ is torsion,
and we can refer to
Corollary~\ref{derived-completion-tensor-with-torsion-isom-cor}.
\end{proof}

 Finally, let $\cM$ be a sheaf of torsion $\fO_\fX$\+modules on~$\fX$.
 We will say that $\cM$ is a \emph{quasi-coherent torsion sheaf}
on $\fX$ if the quasi-coherence axiom~(ii) holds for all affine open
formal subschemes $\fV\subset\fU\subset\fX$ (or equivalently, for all
affine open formal subschemes $\fV\subset\fU\subset\fX$ subordinate
to~$\bW$).
 According to the preceding discussion in this section, the category
of quasi-coherent torsion sheaves on $\fX$ is equivalent to
the category of quasi-coherent torsion (pre)sheaves on the topology
base $\bB$ consisting of all affine open formal subschemes of $\fX$
subordinate to~$\bW$ (for any given open covering $\bW$ of the formal
scheme~$\fX$).

 We will denote the category of quasi-coherent torsion sheaves on $\fX$
by $\fX\Tors$.
 One can see that $\fX\Tors$ is a Grothendieck abelian category
and, for every affine open formal subscheme $\fU\subset\fX$,
the functor of sections ${-}(\fU)\:\fX\Tors\rarrow\fO_\fX(\fU)\Tors
\subset\fO_\fX(\fU)\Modl$ is an exact functor preserving infinite
direct sums.
 Lemma~\ref{flat-map-of-adic-rings-direct-image-lemma}(a) is
relevant here.
 The infinite direct sums in $\fX\Tors$ agree with the ones in
$(\fX,\fO_\fX)\Sh$ (cf.\ Remark~\ref{on-co-sheaf-axioms-remark}(3)).
 So the functor of sections ${-}(\fU)\:\fX\Tors\rarrow
\fO_\fX(\fU)\Modl$ preserves infinite direct sums for any
quasi-compact open formal subscheme $\fU\subset\fX$.

 A short sequence $0\rarrow\cL\rarrow\cM\rarrow\cN\rarrow0$ of
quasi-coherent torsion sheaves on $\fX$ is exact in $\fX\Tors$ if
and only if, for every affine open formal subscheme $\fU\subset\fX$,
the short sequence of torsion $\fO_\fX(\fU)$\+modules
$0\rarrow\cL(\fU)\rarrow\cM(\fU)\rarrow\cN(\fU)\rarrow0$ is exact.
 By Lemma~\ref{locality-of-exactness-for-torsion}(a), it suffices
to check this condition for affine open formal subschemes $\fU$
belonging to any chosen affine open covering of~$\fX$.

 Moreover, a complex $\cM^\bu$ in $\fX\Tors$ is exact if and only if,
for every affine open formal subscheme $\fU\subset\fX$, the complex of
torsion $\fO_\fX(\fU)$\+modules is exact in $\fO_\fX(\fU)\Tors$.
 By Lemma~\ref{locality-of-exactness-for-torsion}(d), it suffices
to check this condition for affine open formal subschemes $\fU$
belonging to any chosen affine open covering of~$\fX$.

 For any affine Noetherian formal scheme $\fU$, the functor of
global sections $\cM\longmapsto\cM(\fU)$ establishes an equivalence
between the abelian category of quasi-coherent torsion sheaves
$\fU\Tors$ on $\fU$ and the abelian category $\fO_\fU(\fU)\Tors$ of
torsion modules over the adic Noetherian ring $\fO_\fU(\fU)$.

\subsection{Contraherent cosheaves of contramodules}
\label{contraherent-cosheaves-of-contramods-subsecn}
 We keep the notation of Section~\ref{qcoh-torsion-sheaves-subsecn}.
 So $\fX$ is a locally Noetherian formal scheme, $\bW$ is an open
covering of $\fX$, and $\bB$ is the base of open subsets in $\fX$
consisting of all the affine open formal subschemes $\fU\subset\fX$
subordinate to~$\bW$.

 Let $\fP$ be a copresheaf of $\fO_\fX$\+modules on~$\bB$.
 We will say that $\fP$ is a \emph{contraherent copresheaf of
contramodules} (or a \emph{contraherent contramodule copresheaf}) if
the following conditions hold:
\begin{enumerate}
\setcounter{enumi}{2}
\renewcommand{\theenumi}{\roman{enumi}}
\item For every affine open formal subscheme $\fU\subset\fX$
subordinate to $\bW$, the $\fO_\fX(\fU)$\+module $\fP[\fU]$ is
a contramodule (as a module over the adic topological
ring~$\fO_\fX(\fU)$).
\item For every pair of affine open formal subschemes
$\fV\subset\fU\subset\fX$ subordinate to $\bW$,
the $\fO_\fX(\fV)$\+module map
$$
 \fP[\fV]\lrarrow\Hom_{\fO_\fX(\fU)}(\fO_\fX(\fV),\fP[\fU])
$$
corresponding by adjunction to the corestriction map
$\fP[\fV]\rarrow\fP[\fU]$ is an isomorphism.
\item For every pair of affine open formal subschemes
$\fV\subset\fU\subset\fX$ subordinate to $\bW$, one has
$$
 \Ext^1_{\fO_\fX(\fU)\Contra}(\fO_\fX(\fV),\fP[\fU])=0.
$$
\end{enumerate}
 We will call condition~(iii) is the \emph{contramoduleness axiom},
while condition~(iv) will be called the \emph{contraherence axiom},
and condition~(v) will be called the \emph{contraadjustedness axiom}.
 Cf.\ the definitions of a (locally) cohereherent cosheaf of modules
over a scheme in~\cite[Sections~2.2 and~3.1]{Pcosh}
and~\cite[Section~2.2]{Pdomc}.

 Comparing the conditions of
Lemma~\ref{contraadjusted-contramodules-lemma}(2) and~(3), one can
see that condition~(v) holds for a \emph{fixed} affine open formal
subscheme $\fU\subset\fX$ subordinate to $\bW$ and \emph{all} affine
open formal subschemes $\fV\subset\fU$ if and only if
the contramodule $\fO_\fX(\fU)$\+module $\fP[\fU]$ is contraadjusted.
 So our terminology is consistent.

\begin{lem} \label{contraherence-implies-cosheaf-axiom}
 The contramoduleness, contraherence, and contraadjustedness axioms
imply the cosheaf axiom:
 Any contraherent copresheaf of contramodules on\/ $\bB$ is a cosheaf
of\/ $\fO_\fX$\+modules on\/~$\bB$.
\end{lem}

\begin{proof}
 The assertion follows from exactness of the rightmost fragment of
sequence~\eqref{contraadj-homol-Cech-sequence-eqn}
from Lemma~\ref{contraadj-homol-Cech-sequence-lemma} together with
Remarks~\ref{on-co-sheaf-axioms-remark}(1\+-2).
\end{proof}

 Thus any contraherent contramodule copresheaf $\fP$ on $\bB$ is, in
fact, a cosheaf on $\bB$, and we will speak of $\fP$ as
a \emph{contraherent contramodule cosheaf} (or a \emph{contraherent
cosheaf of contramodules}) on~$\bB$.
 By Lemma~\ref{contraherence-implies-cosheaf-axiom} and in view of
the discussion in Section~\ref{co-sheaves-of-modules-subsecn}, any
contraherent cosheaf of contramodules on $\bB$ can be extended to
a cosheaf of $\fO_\fX$\+modules on the whole formal scheme $\fX$ in
a unique way.
 Moreover, according to the discussion in
Section~\ref{cosheaves-of-contramodules-subsecn}, the resulting
cosheaf is a cosheaf of contramodule $\fO_\fX$\+modules on~$\fX$.

 The notion of a contraherent cosheaf over a scheme (hence also of
a contraherent cosheaf of contramodules over a formal scheme) is
\emph{not} local.
 See~\cite[Example~3.2.1]{Pcosh} for a counterexample.
 
 More precisely, in the context of a formal scheme $\fX$, for any
contraherent cosheaf of contramodules $\fP$ on $\bB$ extended to
a cosheaf of $\fO_\fX$\+modules $\fP$ on the whole of $\fX$,
\emph{both} the contramoduleness axiom~(iii) and the contraadjustedness
axiom~(v) hold for all affine open formal subschemes $\fU\subset\fX$.
 For the contramoduleness axiom, this was mentioned in the previous
paragraph.
 For the contraadjustedness axiom, this follows from the facts that
the restrictions of scalars preserve contraadjustedness (see
Lemma~\ref{restriction-of-scalars-contraadjusted}), finite direct
sums of contraadjusted contramodules are contraadjusted, and any
quotient contramodule of a contraadjusted contramodule is
contraadjusted (as mentioned in
Section~\ref{prelim-veryflat-and-contraadjusted-subsecn}).
 It is the contraherence axiom~(iv) that is not local.

\subsection{Locally contraherent cosheaves of contramodules}
\label{locally-contraherent-cosheaves-of-contramods-subsecn}
 We still keep the notation of
Sections~\ref{qcoh-torsion-sheaves-subsecn}
and~\ref{contraherent-cosheaves-of-contramods-subsecn}.

 Let\/ $\fP$ be a cosheaf of contramodule $\fO_\fX$\+modules on~$\fX$.
 We will say that $\fP$ is a \emph{$\bW$\+locally contraherent
cosheaf of contramodules} (or a \emph{$\bW$\+locally contraherent
contramodule cosheaf}) on $\fX$ if the contraherence axiom~(iv) and
the contraadjustedness axiom~(v) hold for all affine open formal
subschemes $\fV\subset\fU\subset\fX$ subordinate to~$\bW$.
 If this is the case, then the contraadjustedness axiom actually
holds for all affine open formal subschemes $\fV\subset\fU\subset\fX$,
as explained in
Section~\ref{contraherent-cosheaves-of-contramods-subsecn}.
 According to the discussion in 
Section~\ref{contraherent-cosheaves-of-contramods-subsecn},
the category of $\bW$\+locally contraherent cosheaves of contramodules
on $\fX$ is equivalent to the category of contraherent co(pre)sheaves
on the topology base $\bB$ consisting of all affine open formal
subschemes of $\fX$ subordinate to~$\bW$.

 A cosheaf of contramodule $\fO_\fX$\+modules $\fP$ is said to be
a \emph{contraherent cosheaf of contramodules} (or a \emph{contraherent
contramodule cosheaf}) on $\fX$ if the contraherence axiom~(iv) and
the contraadjustedness axiom~(v) hold for all affine open formal
subschemes $\fV\subset\fU\subset\fX$.
 So a contraherent cosheaf of contramodules is the same thing as
a $\bW_\fX$\+locally contraherent cosheaf of contramodules (where
$\bW_\fX=\{\fX\}$ is the trivial open covering of~$\fX$).

 The next proposition is the formal scheme version of the homological
criterion of contraherence from~\cite[Lemma~3.2.2]{Pcosh}.
 The definition of the homological \v Cech complex
$C_\bu(\{\fV_\alpha\},\fP)$ was given in
Section~\ref{co-sheaves-of-modules-subsecn}.

\begin{prop} \label{homological-criterion-of-contraherence-prop}
 Let\/ $\fU$ be an affine Noetherian formal scheme with an open
covering\/ $\bW$, and let\/ $\fU=\bigcup_\alpha\fW_\alpha$ be a finite
affine open covering of\/ $\fU$ subordinate to\/~$\bW$.
 Let\/ $\fP$ be a\/ $\bW$\+locally contraherent cosheaf of
contramodules on\/~$\fU$.
 Then the following three conditions are equivalent:
\begin{enumerate}
\item $H_1C_\bu(\{\fW_\alpha\},\fP)=0$;
\item $H_iC_\bu(\{\fW_\alpha\},\fP)=0$ for all $i>0$;
\item $\fP$~is a (globally) contraherent cosheaf of contramodules
on\/~$\fU$.
\end{enumerate}
\end{prop}

\begin{proof}
 (2)~$\Longrightarrow$~(1) Obvious.

 (3)~$\Longrightarrow$~(2) The coagumented homological \v Cech
complex $C_\bu(\{\fW_\alpha\},\fP)\rarrow\fP[\fU]$ is acyclic by
Lemma~\ref{contraadj-homol-Cech-sequence-lemma} (applied to
the contraadjusted contramodule $\fP[\fU]$ over the adic topological
ring $R=\fO_\fU(\fU)$).

 (1)~$\Longrightarrow$~(3) This is dual-analogous to (but more
complicated than) the proof of
Proposition~\ref{torsion-quasi-coherence-is-local-prop}.
 Let $\fV\subset\fU\subset\fU$ be two affine open formal subschemes.
 We need to check the contraherence axiom~(iv) for $\fP$ with
respect to $\fV\subset\fU$.

 Notice that $\fV=\bigcup_{\alpha=1}^N(\fV\cap\fW_\alpha)$ is
an finite open covering of $\fV$ by affine open formal subschemes
subordinate to~$\bW$.
 By the cosheaf axiom~\eqref{topology-base-cosheaf-axiom} and
the assumption~(1), we have a four-term right exact sequence
of contramodule $\fO_\fU(\fU)$\+modules
\begin{multline} \label{Cech-H1-vanishing-four-term-right-exact-seq}
 \bigoplus\nolimits_{1\le\alpha<\beta<\gamma\le N}
 \fP[\fW_\alpha\cap\fW_\beta\cap\fW_\gamma]
 \lrarrow\bigoplus\nolimits_{1\le\alpha<\beta\le N}
 \fP[\fW_\alpha\cap\fW_\beta] \\ \lrarrow
 \bigoplus\nolimits_{\alpha=1}^N\fP[\fW_\alpha]
 \lrarrow\fP[\fU]\lrarrow0.
\end{multline}
 By the cosheaf axiom~\eqref{topology-base-cosheaf-axiom}, we also
have a right exact sequence of contramodule $\fO_\fU(\fV)$\+modules
\begin{equation} \label{cosheaf-axiom-covering-of-V-by-V-cap-W-alpha}
 \bigoplus\nolimits_{1\le\alpha<\beta\le N}
 \fP[\fV\cap\fW_\alpha\cap\fW_\beta] \lrarrow
 \bigoplus\nolimits_{\alpha=1}^N\fP[\fV\cap\fW_\alpha]
 \lrarrow\fP[\fV]\lrarrow0.
\end{equation}

 All the terms
of~\eqref{Cech-H1-vanishing-four-term-right-exact-seq}
are contraadjusted contramodule $\fO_\fU(\fU)$\+modules
by Lemma~\ref{restriction-of-scalars-contraadjusted}.
 Since the class of contraadjusted contramodule $R$\+modules is
closed under quotients in $R\Contra$ for any adic topological ring $R$,
we obtain a four-term exact sequence of contramodule {\hbadness=1325
$\fO_\fU(\fU)$\+modules
\begin{multline} \label{contraadjusted-four-term-exact-sequence}
 0\lrarrow\fQ\lrarrow\bigoplus\nolimits_{1\le\alpha<\beta\le N}
 \fP[\fW_\alpha\cap\fW_\beta] \\ \lrarrow
 \bigoplus\nolimits_{\alpha=1}^N\fP[\fW_\alpha]
 \lrarrow\fP[\fU]\lrarrow0.
\end{multline}
with} a contraadjusted contramodule $\fO_\fU(\fU)$\+module~$\fQ$.
 The contramodule $\fO_\fU(\fU)$\+module $\fO_\fU(\fV)$ is very flat,
so applying the functor $\Hom_{\fO_\fU(\fU)}(\fO_\fU(\fV),{-})$
to the sequence~\eqref{contraadjusted-four-term-exact-sequence}, we
obtain a four-term exact sequence of $\fO_\fU(\fV)$\+modules
{\hbadness=2625
\begin{multline*}
 0\lrarrow\Hom_{\fO_\fU(\fU)}(\fO_\fU(\fV),\fQ)\lrarrow
 \bigoplus\nolimits_{1\le\alpha<\beta\le N}
 \Hom_{\fO_\fU(\fU)}(\fO_\fU(\fV),\>\fP[\fW_\alpha\cap\fW_\beta]) \\
 \lrarrow\bigoplus\nolimits_{\alpha=1}^N
 \Hom_{\fO_\fU(\fU)}(\fO_\fU(\fV),\fP[\fW_\alpha])
 \lrarrow\Hom_{\fO_\fU(\fU)}(\fO_\fU(\fV),\fP[\fU])\lrarrow0.
\end{multline*}
 So, in} particular, the sequence
\begin{multline} \label{Hom-into-contraadjusted-four-term-exact-seq}
 \lrarrow\bigoplus\nolimits_{1\le\alpha<\beta\le N}
 \Hom_{\fO_\fU(\fU)}(\fO_\fU(\fV),\>\fP[\fW_\alpha\cap\fW_\beta]) \\
 \lrarrow\bigoplus\nolimits_{\alpha=1}^N
 \Hom_{\fO_\fU(\fU)}(\fO_\fU(\fV),\fP[\fW_\alpha])
 \lrarrow\Hom_{\fO_\fU(\fU)}(\fO_\fU(\fV),\fP[\fU])\lrarrow0.
\end{multline}
is right exact.

 The corestriction maps in the cosheaf $\fP$ induce a natural morphism
from the right exact
sequence~\eqref{cosheaf-axiom-covering-of-V-by-V-cap-W-alpha} to
the right exact
sequence~\eqref{Hom-into-contraadjusted-four-term-exact-seq}.
 It remains to refer to the next
Lemma~\ref{intersection-of-affine-open-Hom-into-contramod}
to the effect that the map
from~\eqref{cosheaf-axiom-covering-of-V-by-V-cap-W-alpha}
to~\eqref{Hom-into-contraadjusted-four-term-exact-seq} is
an isomorphism on the middle and leftmost terms.
 Hence the natural map $\fP[\fV]\rarrow
\Hom_{\fO_\fU(\fU)}(\fO_\fU(\fV),\fP[\fU])$ on the rightmost
terms is an isomorphism, too.
\end{proof}

\begin{lem} \label{intersection-of-affine-open-Hom-into-contramod}
 Let\/ $\fU$ be an affine Noetherian formal scheme with an open
covering\/ $\bW$, let\/ $\fP$ be a\/ $\bW$\+locally contraherent
cosheaf of contramodules on\/ $\fU$, let\/ $\fV\subset\fU$ be
an affine open formal subscheme, and let\/ $\fW\subset\fU$ be
an affine open formal subscheme subordinate to\/~$\bW$.
 Then the corestriction map of\/ $\fO_\fU(\fW)$\+modules\/
$\fP[\fV\cap\fW]\rarrow\fP[\fW]$ induces an isomorphism of\/
$\fO_\fU(\fV)$\+modules\/ $\fP[\fV\cap\fW]\rarrow
\Hom_{\fO_\fU(\fU)}(\fO_\fU(\fV),\fP[\fW])$.
\end{lem}

\begin{proof}
 This is the dual-analogous version of
Lemma~\ref{intersection-of-affine-open-tensor-with-discrete}.
 By the contraherence axiom~(iv), the restriction map
$\fP[\fV\cap\fW]\rarrow\fP[\fW]$ induces an isomorphism of
$\fO_\fU(\fV\cap\fW)$\+modules
$$
 \fP[\fV\cap\fW]\simeq
 \Hom_{\fO_\fU(\fW)}(\fO_\fU(\fV\cap\fW),\>\fP[\fW]).
$$
 There is also an obvious isomorphism
$$
 \Hom_{\fO_\fU(\fU)}(\fO_\fU(\fV),\fP[\fW])
 \simeq\Hom_{\fO_\fU(\fW)}
 (\fO_\fU(\fV)\ot_{\fO_\fU(\fU)}\fO_\fU(\fW),\>\fP[\fW]).
$$
 It remains to show that the map
$$
 \Hom_{\fO_\fU(\fW)}(\fO_\fU(\fV\cap\fW),\fP[\fW])
 \lrarrow\Hom_{\fO_\fU(\fW)}
 (\fO_\fU(\fV)\ot_{\fO_\fU(\fU)}\fO_\fU(\fW),\>\fP[\fW])
$$
induced by the multiplication map of rings
\begin{equation} \label{affine-open-formal-intersect-multiplicat-map2}
 \fO_\fU(\fV)\ot_{\fO_\fU(\fU)}\fO_\fU(\fW)\lrarrow
 \fO_\fU(\fV\cap\fW)
\end{equation}
is an isomorphism.
 Finally,
the map~\eqref{affine-open-formal-intersect-multiplicat-map2} becomes
an isomorphism after the adic completion functor $\Lambda$ is applied,
while the $\fO_\fU(\fW)$\+module $\fP[\fW]$ is a (quotseparated)
contramodule, and we can refer to
Corollary~\ref{flat-reductions-derived-Lambda-is-underived-cor}
(whose flatness assumption is obviously satisfied for
the $\fO_\fU(\fW)$\+module
$F=\fO_\fU(\fV)\ot_{\fO_\fU(\fU)}\fO_\fU(\fW)$).
\end{proof}

 We will denote the category of $\bW$\+locally contraherent cosheaves
of contramodules on $\fX$ by $\fX\Lcth_\bW$.
 The category of contraherent cosheaves of contramodules on $\fX$ will
be denoted by $\fX\Ctrh=\fX\Lcth_{\{\fX\}}$.
 A cosheaf of contramodule $\fO_\fX$\+modules on $\fX$ is said
to be \emph{locally contraherent} if it is $\bW$\+locally contraherent
with respect to some open covering $\bW$ of the formal scheme~$\fX$.
 The category of locally contraherent cosheaves of contramodules on
$\fX$ will be denoted by $\fX\Lcth=\bigcup_\bW\fX\Lcth_\bW$.

 All infinite direct products exists in the category $\fX\Lcth_\bW$
(in particular, in $\fX\Ctrh$).
 The functors of cosections $\fP\longmapsto\fP[\fU]$ over affine open
formal subschemes $\fU\subset\fX$ subordinate to $\bW$ preserve
infinite products.
 In fact, the condition that $\fU$ is subordinate to $\bW$ is not
needed here; and the same assertion holds for all quasi-compact
open formal subschemes $\fU\subset\fX$.
 In other words, the infinite products in $\fX\Lcth_\bW$ agree with
the ones in $(\fX,\fO_\fX)\Cosh$ (see
Remark~\ref{on-co-sheaf-axioms-remark}(3)).

 A short sequence of $\bW$\+locally contraherent cosheaves of
contramodules $0\rarrow\fL\rarrow\fM\rarrow\fN\rarrow0$ on $\fX$
is said to be (\emph{admissible}) exact if the short sequence of
(contramodule) $\fO_\fX(\fU)$\+modules $0\rarrow\fL[\fU]\rarrow
\fM[\fU]\rarrow\fN[\fU]\rarrow0$ is exact for every affine open formal
subscheme $\fU\subset\fX$ subordinate to~$\bW$.
 By Lemma~\ref{colocality-of-exactness-for-contraadjusted}(a),
it suffices to check this condition for affine open formal subschemes
$\fU$ belonging to any chosen affine open covering of $\fX$ subordinate
to~$\bW$.

 Hence a short sequence of contraherent cosheaves of contramodules on
$\fX$ is exact in $\fX\Ctrh$ if and only if it is exact in
$\fX\Lcth_\bW$.
 More generally, let $\bW'$ be an open covering of $\fX$ subordinate
to~$\bW$.
 Then it follows from
Lemma~\ref{colocality-of-exactness-for-contraadjusted}(a) that
a short sequence in $\fX\Lcth_\bW$ is exact in $\fX\Lcth_\bW$ if and
only if it is exact in $\fX\Lcth_{\bW'}$.
 A short sequence is said to be exact in $\fX\Lcth$ if it is exact
in $\fX\Lcth_\bW$ for some open covering $\bW$ of~$\fX$.

\begin{lem} \label{contrah-cosheaves-of-contramods-lemma}
\textup{(a)} Let\/ $\fX$ be a locally Noetherian formal scheme with
an open covering\/~$\bW$.
 Then the category\/ $\fX\Lcth_\bW$ of\/ $\bW$\+locally contraherent
cosheaves of contramodules on\/ $\fX$, endowed with the class of
(admissible) short exact sequences specified above, is an exact
category (in the sense of Quillen~\cite{Bueh}). \par
\textup{(b)} Let\/ $\fU$ be an affine Noetherian formal scheme.
 Then the functor of global cosections\/ $\fP\longmapsto\fP[\fU]$
establishes an equivalence between the exact category\/
$\fU\Ctrh$ of contraherent cosheaves of contramodules on\/ $\fU$
and the exact category\/ $\fO_\fU(\fU)\Contra^\cta$ of contraadjusted
contramodules over the adic Noetherian ring\/~$\fO_\fU(\fU)$.
\end{lem}

\begin{proof}
 Both the assertions follow from
Corollary~\ref{colocalization-contraadjusted-exactness-cor} and
Proposition~\ref{colocalization-of-contraadjusted-contramodule-prop}.
\end{proof}

 It follows from
Proposition~\ref{homological-criterion-of-contraherence-prop}
that the full subcategory $\fX\Ctrh$ is closed under extensions
in the exact category $\fX\Lcth_\bW$.
 More generally, the full subcategory $\fX\Lcth_\bW$ is closed
under extensions in the exact category $\fX\Lcth_{\bW'}$.
 The full subcategory $\fX\Lcth_\bW$ is also closed under infinite
products in $\fX\Lcth_{\bW'}$.

 Furthermore, it follows from
Lemma~\ref{colocality-of-exactness-for-contraadjusted}(b) that
the full subcategory $\fX\Ctrh$ is closed under kernels of
admissible epimorphisms in $\fX\Lcth_\bW$.
 More generally, the full subcategory $\fX\Lcth_\bW$ is closed
under kernels of admissible epimorphisms in $\fX\Lcth_{\bW'}$.
 So a morphism in $\fX\Lcth_\bW$ is an admissible epimorphism in
$\fX\Lcth_\bW$ if and only if it is an admissible epimorphism in
$\fX\Lcth_{\bW'}$.
 The similar assersions for admissible monomorphisms are \emph{not}
true~\cite[Example~3.2.1]{Pcosh}.

\subsection{Locally cotorsion locally contraherent cosheaves}
\label{locally-contraherent-lct-cosheaves-subsecn}
 A $\bW$\+locally contraherent cosheaf of contramodules $\fP$ on $\fX$
is said to be \emph{locally cotorsion} if the contramodule
$\fO_\fX(\fU)$\+module $\fP[\fU]$ is cotorsion for every affine open
formal subscheme $\fU\subset\fX$ subordinate to~$\bW$.
 By Proposition~\ref{colocalization-of-cotorsion-contramodule-prop}
and Corollary~\ref{colocality-of-cotorsion-presuming-contraadjusted},
it suffices to check that condition for affine open formal
subschemes $\fU$ belonging to any chosen affine open covering of $\fX$
subordinate to~$\bW$.

 We will denote the full subcategory of locally cotorsion $\bW$\+locally
contraherent cosheaves of contramodules on $\fX$ by $\fX\Lcth_\bW^\lct
\subset\fX\Lcth_\bW$.
 The category of locally cotorsion contraherent cosheaves of
contramodules on $\fX$ will be denoted by
$\fX\Ctrh^\lct=\fX\Lcth_{\{\fX\}}^\lct$.
 It follows from the previous paragraph that
$\fX\Ctrh^\lct=\fX\Ctrh\cap\fX\Lcth_\bW^\lct$ and
$\fX\Lcth_\bW^\lct=\fX\Lcth_\bW\cap\fX\Lcth_{\bW'}^\lct$ (in
the notation from the end of the previous section).
 The category of locally cotorsion locally contraherent cosheaves of
contramodules on $\fX$ will be denoted by
$\fX\Lcth^\lct=\bigcup_\bW\fX\Lcth_\bW^\lct$.

 Clearly, the full subcategory $\fX\Lcth_\bW^\lct$ is closed under
extensions, cokernels of admissible monomorphisms, and infinite
direct products in $\fX\Lcth_\bW$.
 So the full subcategories $\fX\Ctrh^\lct$, \,$\fX\Lcth_\bW^\lct$,
and $\fX\Lcth^\lct$ inherit exact category structures from
the ambient exact categories $\fX\Ctrh$, \,$\fX\Lcth_\bW$, and
$\fX\Lcth$.
 A morphism in $\fX\Lcth_\bW^\lct$ is an admissible monomorphism in
$\fX\Lcth_\bW^\lct$ if and only if it is an admissible monomorphism
in $\fX\Lcth_\bW$.

\begin{lem} \label{lct-contrah-cosheaves-of-contramods-on-affine-lemma}
 Let\/ $\fU$ be an affine Noetherian formal scheme.
 Then the equivalence of exact categories\/
$\fU\Ctrh\simeq\fO_\fU(\fU)\Contra^\cta$ from
Lemma~\ref{contrah-cosheaves-of-contramods-lemma}(b) restricts to
an equivalence\/ $\fU\Ctrh^\lct\simeq\fO_\fU(\fU)\Contra^\cot$
between the exact category of locally cotorsion contraherent cosheaves
of contramodules on\/ $\fU$ and the exact category of cotorsion
contramodules over the adic Noetherian ring\/~$\fO_\fU(\fU)$.
\end{lem}

\begin{proof}
 Follows from
Proposition~\ref{colocalization-of-cotorsion-contramodule-prop}.
\end{proof}

\begin{lem} \label{complex-of-cosheaves-exactness-criterion}
 Let\/ $\fX$ be a locally Noetherian formal scheme with an open
covering\/~$\bW$.
 In this context: \par
\textup{(a)} A complex\/ $\fP^\bu$ in the exact category of\/
$\bW$\+locally contraherent cosheaves of contramodules\/ $\fX\Lcth_\bW$
is exact if and only if, for every affine open formal subscheme\/
$\fU\subset\fX$ subordinate to\/~$\bW$, the complex of\/
$\fO_\fX(\fU)$\+modules\/ $\fP^\bu[\fU]$ is exact.
 It suffices to check this condition for affine open formal
subschemes\/ $\fU\subset\fX$ belonging to any chosen open covering
of\/ $\fX$ subordinate to\/~$\bW$. \par
\textup{(a)} A complex\/ $\fP^\bu$ in the exact category of locally
cotorsion\/ $\bW$\+locally contraherent cosheaves of contramodules\/
$\fX\Lcth_\bW^\lct$ is exact if and only if, for every affine open
formal subscheme\/ $\fU\subset\fX$ subordinate to\/~$\bW$, the complex
of\/ $\fO_\fX(\fU)$\+modules\/ $\fP^\bu[\fU]$ is exact.
 It suffices to check this condition for affine open formal
subschemes\/ $\fU\subset\fX$ belonging to any chosen open covering
of\/ $\fX$ subordinate to\/~$\bW$.
\end{lem}

\begin{proof}
 This is similar to~\cite[Lemma~3.1.2]{Pcosh}.
 The results of
Lemmas~\ref{colocality-of-exactness-for-contraadjusted}(c)
and~\ref{colocality-of-exactness-for-cotorsion}(c) are relevant,
as well as the results of
Proposition~\ref{contraadjusted-and-cotorsion-periodicity-prop}.
\end{proof}

 To distinguish them from the objects of the exact category
$\fX\Lcth^\lct$ of locally cotorsion locally contraherent cosheaves
of contramodules, the objects of the ambient exact category of
locally contraherent cosheaves of contramodules $\fX\Lcth$ are
sometimes called \emph{locally contraadjusted} locally contraherent
cosheaves of contramodules.
 So all our locally contraherent cosheaves of contramodules are always,
at least, locally contraadjusted (while some of them are even
locally cotorsion).

\subsection{Morphisms of formal schemes}
\label{morphisms-of-formal-schemes-subsecn}
 Let $\ff\:\fY\rarrow\fX$ be a morphism of locally Noetherian formal
schemes.
 Then, for every open formal subscheme $\fU\subset\fX$, the preimage
$\fV=\ff^{-1}(\fU)$ is an open formal subscheme $\fV\subset\fY$.

 Let $\fU\subset\fX$ and $\fV\subset\fY$ be a pair of affine open
formal subschemes such that $\ff(\fV)\subset\fU$.
 Then the induced morphism of affine formal schemes $\ff|_\fV^\fU\:
\fV\rarrow\fU$ corresponds to a continuous homomorphism of adic
Noetherian rings $f_{\fU,\fV}\:\fO_\fX(\fU)\rarrow\fO_\fY(\fV)$.
 Given an affine open formal subscheme $\fU'\subset\fU$, the affine
open formal subscheme $\ff^{-1}(\fU')\cap\fV\subset\fY$ can be
described as
$$
 \ff^{-1}(\fU')\cap\fV=
 \Spf\Lambda\bigl(\fO_\fX(\fU')\ot_{\fO_\fX(\fU)}\fO_\fY(\fV)\bigr),
$$
where the tensor product of adic topological rings
$\fO_\fX(\fU')\ot_{\fO_\fX(\fU)}\fO_\fY(\fV)$ is endowed with
the tensor product topology, as in
Section~\ref{prelim-tensor-products-of-adic-topological-subsecn},
and $\Lambda$ is the adic completion functor
from Section~\ref{prelim-adic-completions-subsecn}.
 Notice that the ring $\Lambda\bigl(\fO_\fX(\fU')\ot_{\fO_\fX(\fU)}
\fO_\fY(\fV)\bigr)$ is Noetherian by
Lemma~\ref{open-immersion-adic-ring-base-change-context}(a) and
Corollary~\ref{adic-topological-rings-complet-noether-local-cor}(a).

 We will say that a morphism of formal schemes $\ff\:\fY\rarrow\fX$ is
\emph{tight} (or, in a more traditional terminology, \emph{adic}) if,
for every pair of affine open formal subschemes $\fU\subset\fX$ and
$\fV\subset\fY$ such that $\ff(\fV)\subset\fU$, the continuous map of
adic Noetherian rings $f_{\fU,\fV}$ is tight (in the sense of
Section~\ref{adic-rings-subsecn}).
 By Proposition~\ref{complete-tightness=tautness-prop}, this is
equivalent to the map $f_{\fU,\fV}$ being taut in the sense of
Section~\ref{prelim-tight-taut-subsecn}.

 The tightness property of morphisms of formal schemes is local in
the following sense.
 A morphism~$\ff$ is tight if and only if there exist two affine open
coverings $\fX=\bigcup_\alpha\fU_\alpha$ and $\fY=\bigcup_\alpha
\fV_\alpha$, indexed by the same set of indices~$\alpha$, such that
$\ff(\fV_\alpha)\subset\fU_\alpha$ for every index~$\alpha$ and
the continuous ring map~$f_{\fU_\alpha,\fV_\alpha}$ is tight
for every~$\alpha$.
 This follows from
Corollaries~\ref{locality-of-tightness-in-total-space-cor}
and~\ref{locality-of-tightness-in-the-base-cor}.

 We will say that a morphism of formal schemes $\ff\:\fY\rarrow\fX$ is
\emph{flat} if, for every pair of affine open formal subschemes
$\fU\subset\fX$ and $\fV\subset\fY$ such that $\ff(\fV)\subset\fU$,
the continuous map of adic Noetherian rings $f_{\fU,\fV}$ is flat
(in the sense of Section~\ref{adic-rings-subsecn}).
 By Lemma~\ref{tight-taut-completion-very-flatness}(a), this is
equivalent to the map $f_{\fU,\fV}$ being taut-flat in the sense of
Section~\ref{prelim-morphism-very-flatness-locality-subsecn}.
 Notice that any flat morphism of formal schemes is tight by
the definition.

 The flatness property of morphisms of formal schemes is local.
 In other words, a morphism~$\ff$ is flat if and only if there exist
two affine open coverings $\fX=\bigcup_\alpha\fU_\alpha$ and
$\fY=\bigcup_\alpha\fV_\alpha$, indexed by the same set of
indices~$\alpha$, such that $\ff(\fV_\alpha)\subset\fU_\alpha$ for
every index~$\alpha$ and the continuous ring
map~$f_{\fU_\alpha,\fV_\alpha}$ is flat for every~$\alpha$.
 This follows from
Corollaries~\ref{locality-of-very-flatness-in-total-space-cor}(a)
and~\ref{locality-of-very-flatness-in-the-base-cor}(a).

 We will say that a morphism of formal schemes $\ff\:\fY\rarrow\fX$ is
\emph{very flat} if, for every pair of affine open formal subschemes
$\fU\subset\fX$ and $\fV\subset\fY$ such that $\ff(\fV)\subset\fU$,
the continuous map of adic Noetherian rings $f_{\fU,\fV}$ is very flat
(in the sense of
Section~\ref{prelim-veryflat-and-contraadjusted-subsecn}).
 By Lemma~\ref{tight-taut-completion-very-flatness}(b), this is
equivalent to the map $f_{\fU,\fV}$ being taut-very flat in the sense
of Section~\ref{prelim-morphism-very-flatness-locality-subsecn}.
 Any very flat morphism of formal schemes is flat (hence tight) by
the definition.

 The very flatness property of morphisms of formal schemes is local.
 In other words, a morphism~$\ff$ is very flat if and only if there
exist two affine open coverings $\fX=\bigcup_\alpha\fU_\alpha$ and
$\fY=\bigcup_\alpha\fV_\alpha$, indexed by the same set of
indices~$\alpha$, such that $\ff(\fV_\alpha)\subset\fU_\alpha$ for
every index~$\alpha$ and the continuous ring
map~$f_{\fU_\alpha,\fV_\alpha}$ is very flat for every~$\alpha$.
 This follows from
Corollaries~\ref{locality-of-very-flatness-in-total-space-cor}(b)
and~\ref{locality-of-very-flatness-in-the-base-cor}(b).

 A morphism of locally Noetherian formal schemes $\ff\:\fY\rarrow\fX$
is said to be \emph{quasi-compact} if, for every quasi-compact open
formal subscheme $\fU\subset\fX$, the open formal subscheme
$\ff^{-1}(\fU)\subset\fY$ is quasi-compact.
 It suffices to check this condition for affine open formal subschemes
$\fU\subset\fX$, or for affine open formal subschemes belonging to
any chosen affine open covering of $\fX$, or even for quasi-compact
open formal subschemes belonging to any chosen quasi-compact open
covering of~$\fX$.
 (The fact that all locally Noetherian formal schemes are
quasi-separated is relevant to the validity of the latter formulation.)

 A morphism of locally Noetherian formal schemes $\ff\:\fY\rarrow\fX$
is said to be \emph{semi-separated} if, for every semi-separated open
formal subscheme $\fU\subset\fX$, the open formal subscheme
$\ff^{-1}(\fU)\subset\fY$ is semi-separated.
 It suffices to check this condition for affine open formal subschemes
$\fU\subset\fX$, or for affine open formal subschemes belonging to
any chosen affine open covering of $\fX$, or for semi-separated
open formal subschemes belonging to any chosen semi-separated open
covering of~$\fX$.

 A morphism of locally Noetherian formal schemes $\ff\:\fY\rarrow\fX$
is said to be \emph{affine} if, for every affine open formal
subscheme $\fU\subset\fX$, the open formal subscheme
$\ff^{-1}(\fU)\subset\fY$ is affine.
 It suffices to check this condition for affine open formal subschemes
$\fU$ belonging to any chosen affine open covering of~$\fX$.
 In view of the discussion in
Section~\ref{basics-of-formal-schemes-subsecn}, it suffices to prove
the latter assertion for (conventional, not formal) schemes, for
which it can be found in~\cite[D\'efinition~II.1.2.1
and Corollaire~II.1.3.2]{EGAII} or~\cite[Lemma Tag~01S8]{SP}.
 Clearly, any affine morphism of formal schemes is both quasi-compact
and semi-separated.

\subsection{Direct images of quasi-coherent torsion sheaves}
\label{direct-images-of-qcoh-tors-subsecn}
 Let $\ff\:(\fY,\fO_\fY)\rarrow(\fX,\fO_\fX)$ be a morphism of ringed
spaces.
 For any sheaf of $\fO_\fY$\+modules $\cN$ on $\fY$, the sheaf of
$\fO_\fX$\+modules $\ff_*\cN$ on $\fX$ is defined by the rule
$$
 (\ff_*\cN)(\fU)=\cN(\ff^{-1}(\fU))
$$
for all open subsets $\fU\subset\fX$.
 The $\fO_\fY(\ff^{-1}(\fU))$\+module $\cN(\ff^{-1}(\fU))$ is endowed
with an $\fO_\fX(\fU)$\+module structure using the restriction of
scalars with respect to the ring homomorphism
$\fO_\fX(\fU)\rarrow\fO_\fY(\ff^{-1}(\fU))$.
 The sheaf axiom for $\ff_*\cN$ follows easily from the sheaf axiom
for~$\cN$.

 Let $(\fX,\fO_\fX)$ be a ringed space and $\fY\subset\fX$ be an open
subset.
 Consider $\fY$ as a ringed space $(\fY,\fO_\fY)$ whose structure sheaf
$\fO_\fY=\fO_\fX|_\fY$ is the restriction of $\fO_\fX$ to~$\fY$.
 Then the identity inclusion $\fj\:\fY\rarrow\fX$ is naturally
a morphism of ringed spaces.
 For any sheaf of $\fO_\fX$\+modules $\cM$ on $\fX$, the restriction
$\cM|_\fY$ is a sheaf of $\fO_\fY$\+modules on~$\fY$.
 In this context, the direct image functor
$\fj_*\:(\fY,\fO_\fY)\Sh\rarrow(\fX,\fO_\fX)\Sh$ is right adjoint to
the restriction functor $\cM\longmapsto\cM|_\fY\:(\fX,\fO_\fX)\Sh
\rarrow(\fY,\fO_\fY)\Sh$, i.~e., there is a natural isomorphism of
abelian groups
\begin{equation} \label{sheaves-restriction-direct-image-adjunction}
 \Hom_{\fO_\fY}(\cM|_\fY,\cN)\simeq\Hom_{\fO_\fX}(\cM,\fj_*\cN)
\end{equation}
for any sheaf of $\fO_\fX$\+modules $\cM$ and any sheaf of
$\fO_\fY$\+modules~$\cN$.
 Here $\Hom_{\fO_\fY}({-},{-})$ and $\Hom_{\fO_\fX}({-},{-})$ denote
the groups of morphisms in the categories $(\fY,\fO_\fY)\Sh$ and
$(\fX,\fO_\fX)\Sh$, respectively.

\begin{prop} \label{torsion-sheaves-qcoh-direct-images-prop}
 Let\/ $\ff\:\fY\rarrow\fX$ be a quasi-compact morphism of locally
Noetherian formal schemes.
 Then \par
\textup{(a)} for any sheaf of torsion\/ $\fO_\fY$\+modules\/ $\cN$ on\/
$\fY$, the sheaf of\/ $\fO_\fX$\+modules\/ $\ff_*\cN$ on\/ $\fX$ is
a sheaf of torsion\/ $\fO_\fX$\+modules; \par
\textup{(b)} for any quasi-coherent sheaf of torsion\/
$\fO_\fY$\+modules\/ $\cN$ on\/ $\fY$, the sheaf of torsion\/
$\fO_\fX$\+modules\/ $\ff_*\cN$ on\/ $\fX$ is quasi-coherent.
\end{prop}

\begin{proof}
 Part~(a): this is similar to the argument from
Section~\ref{cosheaves-of-contramodules-subsecn}.
 Let $\fU\subset\fX$ be an affine open formal subscheme, and
let $\ff^{-1}(\fU)=\bigcup_\alpha\fW_\alpha$ be a finite affine open
covering of the open formal subscheme $\ff^{-1}(\fU)\subset\fY$.
 Then the $\fO_\fX(\fU)$\+module $(\ff_*\cN)(\fU)$ is a submodule of
the direct sum $\bigoplus_\alpha\cN(\fW_\alpha)$.
 Since any torsion $\fO_\fY(\fW_\alpha)$\+module is also a torsion
$\fO_\fX(\fU)$\+module by
Lemma~\ref{torsion-contra-restriction-of-scalars}(a), and the class of
torsion modules is closed under subobjects and direct sums in $R\Modl$
for any adic topological ring $R$, the assertion follows.

 Part~(b): this is similar to the proof of
Proposition~\ref{torsion-quasi-coherence-is-local-prop}.
 Let $\fV\subset\fU\subset\fX$ be two affine open formal subschemes,
and let $\ff^{-1}(\fU)=\bigcup_{\alpha=1}^N\fW_\alpha$ be a finite
affine open covering of the open formal subscheme
$\ff^{-1}(\fU)\subset\fY$.
 For every pair of indices $1\le\alpha<\beta\le N$, pick a finite
affine open covering $\fW_\alpha\cap\fW_\beta=\bigcup_\gamma
\fT_{\alpha\beta\gamma}$ of the intersection
$\fW_\alpha\cap\fW_\beta\subset\fY$.
 Notice that $\ff^{-1}(\fV)=\bigcup_{\alpha=1}^N
(\ff^{-1}(\fV)\cap\fW_\alpha)$ is a finite affine open covering of
the open formal subscheme $\ff^{-1}(\fV)\subset\fY$, and
$\ff^{-1}(\fV)\cap\fW_\alpha\cap\fW_\beta=\bigcup_\gamma
(\ff^{-1}(\fV)\cap\fT_{\alpha\beta\gamma})$ are also finite affine
open coverings.

 By the sheaf axiom~\eqref{topology-base-sheaf-axiom} and part~(a),
we have  a left exact sequence of torsion $\fO_\fX(\fU)$\+modules
\begin{equation} \label{sheaf-axiom-covering-of-preimage-of-U}
 0\lrarrow(\ff_*\cN)(\fU)\lrarrow
 \bigoplus\nolimits_{\alpha=1}^N\cN(\fW_\alpha)
 \lrarrow\bigoplus\nolimits_{\alpha,\beta,\gamma}
 \cN(\fT_{\alpha\beta\gamma})
\end{equation}
and a left exact sequence of torsion $\fO_\fX(\fV)$\+modules
\begin{equation} \label{sheaf-axiom-covering-of-preimage-of-V}
 0\lrarrow(\ff_*\cN)(\fV)\lrarrow\bigoplus\nolimits_{\alpha=1}^N
 \cN(\ff^{-1}(\fV)\cap\fW_\alpha)\lrarrow
 \bigoplus\nolimits_{\alpha,\beta,\gamma}
 \cN(\ff^{-1}(\fV)\cap\fT_{\alpha\beta\gamma}).
\end{equation}
 As explained in the proof of
Proposition~\ref{torsion-quasi-coherence-is-local-prop},
the tensor product functor $\fO_\fX(\fV)\ot_{\fO_\fX(\fU)}{-}$
preserves exactness of sequences of torsion $\fO_\fX(\fU)$\+modules.
 So, applying this functor
to~\eqref{sheaf-axiom-covering-of-preimage-of-U}, we obtain a left
exact sequence of $\fO_\fX(\fV)$\+modules {\hbadness=1500
\begin{multline} \label{sheaf-axiom-tensored-with-O(V)-sequence-II}
 0\lrarrow\fO_\fX(\fV)\ot_{\fO_\fX(\fU)}(\ff_*\cN)(\fU)
 \lrarrow\bigoplus\nolimits_{\alpha=1}^N
 \fO_\fX(\fV)\ot_{\fO_\fX(\fU)}\cN(\fW_\alpha) \\
 \lrarrow\bigoplus\nolimits_{\alpha,\beta,\gamma}
 \fO_\fX(\fV)\ot_{\fO_\fX(\fU)}\cN(\fT_{\alpha\beta\gamma}).
\end{multline}

 The} restriction maps in the sheaf $\cN$ induce a natural morphism
from the left exact
sequence~\eqref{sheaf-axiom-tensored-with-O(V)-sequence-II}
to the left exact
sequence~\eqref{sheaf-axiom-covering-of-preimage-of-V}.
 It remains to refer to the next
Lemma~\ref{preimage-of-affine-open-tensor-with-discrete}
to the effect that the map
from~\eqref{sheaf-axiom-tensored-with-O(V)-sequence-II}
to~\eqref{sheaf-axiom-covering-of-preimage-of-V} is an isomorphism
on the middle and rightmost terms.
 Hence the natural map $\fO_\fX(\fV)\ot_{\fO_\fX(\fU)}(\ff_*\cN)(\fU)
\rarrow(\ff_*\cN)(\fV)$ of the leftmost terms in an isomorphism, too,
as desired.
\end{proof}

\begin{lem} \label{preimage-of-affine-open-tensor-with-discrete}
 Let\/ $\fg\:\fW\rarrow\fU$ be a morphism of affine Noetherian formal
schemes, and let\/ $\fV\subset\fU$ be an affine open formal subscheme.
 Let\/ $\cN$ be a quasi-coherent torsion sheaf on\/~$\fW$.
 Then the restriction map of\/ $\fO_\fW(\fW)$\+modules\/
$\cN(\fW)\rarrow\cN(\fg^{-1}(\fV))$ induces an isomorphism of\/
$\fO_\fU(\fV)$\+modules\/ $\fO_\fU(\fV)\ot_{\fO_\fU(\fU)}\cN(\fW)
\rarrow\cN(\fg^{-1}(\fV))$.
\end{lem}

\begin{proof}
 This is similar to
Lemma~\ref{intersection-of-affine-open-tensor-with-discrete}.
 By the quasi-coherence axiom, the restriction map
$\cN(\fW)\rarrow\cN(\fg^{-1}(\fV))$ induces an isomorphism of
$\fO_\fW(\fg^{-1}(\fV))$\+modules
$$
 \fO_\fW(\fg^{-1}(\fV))\ot_{\fO_\fW(\fW)}\cN(\fW)
 \simeq\cN(\fg^{-1}(\fV)).
$$
 There is also an obvious isomorphism
$$
 \fO_\fU(\fV)\ot_{\fO_\fU(\fU)}\cN(\fW)\simeq
 \bigl(\fO_\fU(\fV)\ot_{\fO_\fU(\fU)}\fO_\fW(\fW)\bigr)
 \ot_{\fO_\fW(\fW)}\cN(\fW).
$$
 It remains to show that the map
$$
 \bigl(\fO_\fU(\fV)\ot_{\fO_\fU(\fU)}\fO_\fW(\fW)\bigr)
 \ot_{\fO_\fW(\fW)}\cN(\fW)
 \lrarrow\fO_\fW(\fg^{-1}(\fV))\ot_{\fO_\fW(\fW)}\cN(\fW)
$$
induced by the multiplication map of rings
\begin{equation} \label{affine-open-formal-preimage-multiplicat-map}
 \fO_\fU(\fV)\ot_{\fO_\fU(\fU)}\fO_\fW(\fW)\lrarrow
 \fO_\fW(\fg^{-1}(\fV))
\end{equation}
is an isomorphism.
 Now the map~\eqref{affine-open-formal-preimage-multiplicat-map} becomes
an isomorphism after the adic completion functor $\Lambda$ is applied
(as mentioned in Section~\ref{morphisms-of-formal-schemes-subsecn}),
while the $\fO_\fW(\fW)$\+module $\cN(\fW)$ is torsion, and we can refer
to Corollary~\ref{derived-completion-tensor-with-torsion-isom-cor}.
\end{proof}

 For any quasi-compact morphism of locally Noetherian formal schemes
$\ff\:\fY\rarrow\fX$, we have constructed the functor of direct image
of quasi-coherent torsion sheaves
\begin{equation} \label{direct-image-qcoh-tors}
 \ff_*\:\fY\Tors\lrarrow\fX\Tors.
\end{equation}
 The direct image functor~$\ff_*$ is left exact and preserves infinite
direct sums (as well as filtered direct limits).
 In the case of an affine morphism~$\ff$, the functor~$\ff_*$ is
exact.

\subsection{Direct images of contraherent cosheaves of contramodules}
\label{direct-images-of-ctrh-contramod-subsecn}
 Let $\ff\:(\fY,\fO_\fY)\rarrow(\fX,\fO_\fX)$ be a morphism of ringed
spaces.
 For any cosheaf of $\fO_\fY$\+modules $\fQ$ on $\fY$, the cosheaf of
$\fO_\fX$\+modules $\ff_!\fQ$ on $\fX$ is defined by the rule
{\hbadness=1550
$$
 (\ff_!\fQ)[\fU]=\fQ[\ff^{-1}(\fU)]
$$
for} all open subsets $\fU\subset\fX$.
 The $\fO_\fY(\ff^{-1}(\fU))$\+module $\fQ[\ff^{-1}(\fU)]$ is endowed
with an $\fO_\fX(\fU)$\+module structure using the restriction of
scalars with respect to the ring homomorphism
$\fO_\fX(\fU)\rarrow\fO_\fY(\ff^{-1}(\fU))$.
 The cosheaf axiom for $\ff_!\fQ$ follows easily from the cosheaf axiom
for~$\fQ$.

 Let $(\fX,\fO_\fX)$ be a ringed space and $\fY\subset\fX$ be an open
subset.
 Consider $\fY$ as a ringed space $(\fY,\fO_\fY)$ with the structure
sheaf $\fO_\fY=\fO_\fX|_\fY$, and denote by $\fj\:\fY\rarrow\fX$
the natural open immersion morphism of ringed spaces, as in
Section~\ref{direct-images-of-qcoh-tors-subsecn}.
 For any cosheaf of $\fO_\fX$\+modules $\fP$ on $\fX$, the restriction
$\fP|_\fY$, defined by the rule $\fP|_\fY[\fV]=\fP[\fV]$ for all open
subsets $\fV\subset\fY$, is a cosheaf of $\fO_\fY$\+modules on~$\fY$.

 In this context, the direct image functor
$\fj_!\:(\fY,\fO_\fY)\Cosh\rarrow(\fX,\fO_\fX)\Cosh$ is left adjoint to
the restriction functor $\fP\longmapsto\fP|_\fY\:(\fX,\fO_\fX)\Cosh
\rarrow(\fY,\fO_\fY)\Cosh$, i.~e., there is a natural isomorphism of
abelian groups~\cite[formula~(2.7) in Section~2.3]{Pcosh},
\cite[formula~(19) in Section~2.4]{Pdomc}
\begin{equation} \label{cosheaves-direct-image-restriction-adjunction}
 \Hom^{\fO_\fY}(\fQ,\fP|_\fY)\simeq\Hom^{\fO_\fX}(\fj_!\fQ,\fP)
\end{equation}
for any cosheaf of $\fO_\fY$\+modules $\fQ$ and any cosheaf of
$\fO_\fX$\+modules~$\fP$.
 Here $\Hom^{\fO_\fY}({-},{-})$ and $\Hom^{\fO_\fX}({-},{-})$ denote
the groups of morphisms in the categories $(\fY,\fO_\fY)\Cosh$ and
$(\fX,\fO_\fX)\Cosh$, respectively.

\begin{lem} \label{contramodule-cosheaves-direct-images-lemma}
 Let\/ $\ff\:\fY\rarrow\fX$ be a quasi-compact morphism of locally
Noetherian formal schemes.
 Then, for any cosheaf of contramodule\/ $\fO_\fY$\+modules\/ $\fQ$
on\/ $\fY$, the cosheaf of\/ $\fO_\fX$\+modules\/ $\ff_!\fQ$ on\/ $\fX$
is a cosheaf of contramodule\/ $\fO_\fX$\+modules.
\end{lem}

\begin{proof}
 This is similar to the argument from
Section~\ref{cosheaves-of-contramodules-subsecn}.
 Let $\fU\subset\fX$ be an affine open formal subscheme, and
let $\ff^{-1}(\fU)=\bigcup_\alpha\fW_\alpha$ be a finite affine open
covering of the open formal subscheme $\ff^{-1}(\fU)\subset\fY$.
 For every pair of indices $1\le\alpha<\beta\le N$, pick a finite
affine open covering $\fW_\alpha\cap\fW_\beta=\bigcup_\gamma
\fT_{\alpha\beta\gamma}$ of the intersection
$\fW_\alpha\cap\fW_\beta\subset\fY$.
 By the cosheaf axiom~\eqref{topology-base-cosheaf-axiom},
we have a right exact sequence of $\fO_\fX(\fU)$\+modules
$$
 \bigoplus\nolimits_{\alpha,\beta,\gamma}\fQ[\fT_{\alpha\beta\gamma}]
 \lrarrow\bigoplus\nolimits_{\alpha=1}^N\fQ[\fW_\alpha]\lrarrow
 (\ff_!\fQ)[\fU]\lrarrow0.
$$
 Since any contramodule $\fO_\fY(\fW_\alpha)$\+module and any
contramodule $\fO_\fY(\fT_{\alpha\beta\gamma})$\+module are also
contramodule $\fO_\fX(\fU)$\+modules by
Lemma~\ref{torsion-contra-restriction-of-scalars}(b) or~(d), and
the class of contramodules (or quotseparated contramodules) is closed
under cokernels and direct products in $R\Modl$ for any adic
topological ring $R$, the assertion follows.
\end{proof}

 Let $\ff\:\fY\rarrow\fX$ be a morphism of locally Noetherian formal
schemes.
 Let $\bW$ be an open covering of $\fX$ and $\bT$ be an open covering
of~$\fY$.
 Similarly to~\cite[Section~3.3]{Pcosh} and~\cite[Section~2.4]{Pdomc},
we will say that the morphism~$\ff$ is \emph{$(\bW,\bT)$\+affine} if,
for every affine open formal subscheme $\fU\subset\fX$ subordinate to
$\bW$, the open formal subscheme $\ff^{-1}(\fU)\subset\fY$ is affine
\emph{and} subordinate to~$\bT$.
 It is clear from the discussion of affine morphisms of formal
schemes in
Section~\ref{morphisms-of-formal-schemes-subsecn} that
any $(\bW,\bT)$\+affine morphism is affine.

\begin{lem} \label{lcth-contram-W-T-affine-direct-image-lemma}
 Let\/ $\fX$ be a locally Noetherian formal scheme with an open
covering\/ $\bW$, let\/ $\fY$ be a locally Noetherian formal scheme
with an open covering\/ $\bT$, and let\/ $\ff\:\fY\rarrow\fX$ be
a $(\bW,\bT)$\+affine morphism of formal schemes.
 In this setting: \par
\textup{(a)} for any\/ $\bW$\+locally contraherent cosheaf of
contramodules\/ $\fQ$ on\/ $\fY$, the direct image\/ $\ff_!\fQ$ is
a\/ $\bT$\+locally contraherent cosheaf of contramodules on\/~$\fX$;
\par
\textup{(b)} for any locally cotorsion\/ $\bW$\+locally contraherent
cosheaf of contramodules\/ $\fQ$ on\/ $\fY$, the direct image\/
$\ff_!\fQ$ is a locally cotorsion\/ $\bT$\+locally contraherent cosheaf
of contramodules on\/~$\fX$.
\end{lem}

\begin{proof}
 Part~(a): the contramoduleness axiom~(iii) from
Section~\ref{contraherent-cosheaves-of-contramods-subsecn} holds
by Lemma~\ref{contramodule-cosheaves-direct-images-lemma}
(or more directly, by
Lemma~\ref{torsion-contra-restriction-of-scalars}(b) or~(d)).
 The contraadjustedness axiom~(v) holds by
Lemma~\ref{restriction-of-scalars-contraadjusted}.
 It remains to check the contraherence axiom~(iv) for any pair of affine
open formal subschemes $\fV\subset\fU\subset\fX$ subordinate to~$\bW$.
 For this purpose, we refer to the next
Lemma~\ref{preimage-of-affine-open-Hom-into-contramod}.
 Part~(b) follows from part~(a) by virtue of
Lemma~\ref{restriction-of-scalars-cotorsion}.
\end{proof}

\begin{lem} \label{preimage-of-affine-open-Hom-into-contramod}
 Let\/ $\fg\:\fT\rarrow\fU$ be a morphism of affine Noetherian formal
schemes, and let\/ $\fV\subset\fU$ be an affine open formal subscheme.
 Let\/ $\fQ$ be a contraherent cosheaf of contramodules on\/~$\fT$.
 Then the corestriction map of\/ $\fO_\fT(\fT)$\+modules\/
$\fQ[\fg^{-1}(\fV)]\rarrow\fQ[\fT]$ induces an isomorphism of\/
$\fO_\fU(\fV)$\+modules\/ $\fQ[\fg^{-1}(\fV)]\rarrow
\Hom_{\fO_\fU(\fU)}(\fO_\fU(\fV),\fQ[\fT])$.
\end{lem}

\begin{proof}
 This is similar to
Lemma~\ref{intersection-of-affine-open-Hom-into-contramod}
and dual-analogous to
Lemma~\ref{preimage-of-affine-open-tensor-with-discrete}.
 By the contraherence axiom, the restriction map
$\fQ[\fg^{-1}(\fV)]\rarrow\fQ[\fT]$ induces an isomorphism of
$\fO_\fT(\fg^{-1}(\fV))$\+modules
$$
 \fQ[\fg^{-1}(\fV)]\simeq
 \Hom_{\fO_\fT(\fT)}(\fO_\fT(\fg^{-1}(\fV)),\>\fQ[\fT]).
$$
 There is also an obvious isomorphism
$$
 \Hom_{\fO_\fU(\fU)}(\fO_\fU(\fV),\fQ[\fT])
 \simeq\Hom_{\fO_\fT(\fT)}
 (\fO_\fU(\fV)\ot_{\fO_\fU(\fU)}\fO_\fT(\fT),\>\fQ[\fT]).
$$
 It remains to show that the map
$$
 \Hom_{\fO_\fT(\fT)}(\fO_\fT(\fg^{-1}(\fV)),\fQ[\fT])
 \lrarrow\Hom_{\fO_\fT(\fT)}
 (\fO_\fU(\fV)\ot_{\fO_\fU(\fU)}\fO_\fT(\fT),\>\fQ[\fT])
$$
induced by the multiplication map of rings
\begin{equation} \label{affine-open-formal-preimage-multiplicat-map2}
 \fO_\fU(\fV)\ot_{\fO_\fU(\fU)}\fO_\fT(\fT)\lrarrow
 \fO_\fT(\fg^{-1}(\fV))
\end{equation}
is an isomorphism.
 Finally,
the map~\eqref{affine-open-formal-preimage-multiplicat-map2} becomes
an isomorphism after the adic completion functor $\Lambda$ is applied,
while the $\fO_\fT(\fT)$\+module $\fQ[\fT]$ is a (quotseparated)
contramodule, and we can refer to
Corollary~\ref{flat-reductions-derived-Lambda-is-underived-cor}.
 The flatness assumption of the latter corollary is satisfied for
the $\fO_\fT(\fT)$\+module
$F=\fO_\fU(\fV)\ot_{\fO_\fU(\fU)}\fO_\fT(\fT)$ by
Lemma~\ref{tight-flat-adic-ring-base-change-context}(b)
or Lemma~\ref{open-immersion-adic-ring-base-change-context}(a).
\end{proof}

 For any $(\bW,\bT)$\+affine morphism of locally Noetherian formal
schemes $\ff\:\fY\rarrow\fX$, we have constructed the functor of
direct image of $\bT$\+locally contraherent cosheaves of contramodules
\begin{equation} \label{loc-ctrh-contramod-W-T-affine-direct-image}
 \ff_!\:\fY\Lcth_\bT\lrarrow\fX\Lcth_\bW,
\end{equation}
as well as the functor of direct image of locally cotorsion
$\bT$\+locally contraherent cosheaves of contramodules
\begin{equation} \label{loc-ctrh-lct-contramod-W-T-affine-direct-image}
 \ff_!\:\fY\Lcth^\lct_\bT\lrarrow\fX\Lcth^\lct_\bW.
\end{equation}
 Both the direct image functors~$\ff_!$
\,\eqref{loc-ctrh-contramod-W-T-affine-direct-image}
and~\eqref{loc-ctrh-lct-contramod-W-T-affine-direct-image} are exact
(as functors of exact categories) and preserve infinite products.

 In particular, specializing to the case of globally contraherent
cosheaves (i.~e., choosing the trivial open coverings $\bW_\fX=\{\fX\}$
and $\bT_\fY=\{\fY\}$), one obtains, for any affine morphism of locally
Noetherian formal schemes $\ff\:\fY\rarrow\fX$, the functor of direct
image of contraherent cosheaves of contramodules
\begin{equation} \label{ctrh-contramod-affine-direct-image}
 \ff_!\:\fY\Ctrh\lrarrow\fX\Ctrh,
\end{equation}
as well as the functor of direct image of locally cotorsion
contraherent cosheaves of contramodules
\begin{equation} \label{ctrh-lct-contramod-affine-direct-image}
 \ff_!\:\fY\Ctrh^\lct\lrarrow\fX\Ctrh^\lct.
\end{equation}

 Notice that, for nonaffine morphisms of Noetherian (even not formal)
schemes, the functor of direct image of cosheaves of modules does
\emph{not} preserve either contraherence or local contraherence,
generally speaking.
 A counterexample can be found in~\cite[Remark~2.3.1 and
Example~2.3.2]{Pcosh}.
 Additional conditions stronger than contraherence need to be imposed
on a cosheaf of contramodules to make sure that its direct image under
a nonaffine morphism of formal schemes remains contraherent.
 For some results of this kind, see
Corollaries~\ref{ssep-flat-direct-image-of-antilocally-flat-cor},
\ref{co-flasque-direct-image-cor}(b\+-c),
and~\ref{codilute-direct-images} below, as well as their special cases
such as Corollaries~\ref{ssep-proj-contrah-direct-image-cor},
\ref{proj-lct-direct-images}, and~\ref{flat-contrah-direct-images}.

 On the other hand, one can observe that, for every quasi-compact
morphism of locally Noetherian formal schemes $\ff\:\fY\rarrow\fX$,
the functor of direct image of cosheaves of $\fO$\+modules
$\ff_!\:(\fY,\fO_\fY)\Cosh\rarrow(\fX,\fO_\fX)\Cosh$ preserves
infinite products.
 This is clear from Remark~\ref{on-co-sheaf-axioms-remark}(3).

\subsection{Inverse images of quasi-coherent torsion sheaves}
\label{inverse-images-of-qcoh-tors-subsecn}
 We skip the conventional discussion of the inverse images of
sheaves of modules over ringed spaces, based as it is on
the sheafification functor.
 The cosheafification of copresheaves of modules is generally not
well-behaved, because the filtered projective limit functors in
module categories are not exact.
 So, being primarily interested in contraherent cosheaves, we prefer
to avoid the (co)sheafification altogether.

 Let $\fX$ be a locally Noetherian formal scheme with an open covering
$\bW$, let $\fY$ be a locally Noetherian formal scheme with an open
covering $\bT$, and let $\ff\:\fY\rarrow\fX$ be a morphism of formal
schemes.
 We will say that the morphism~$\ff$ is
\emph{$(\bW,\bT)$\+coaffine}~\cite[Sections~2.3 and~3.3]{Pcosh},
\cite[Section~2.5]{Pdomc} if for every affine open formal subscheme
$\fV\subset\fY$ subordinate to $\bT$ there exists an affine open
formal subscheme $\fU\subset\fX$ subordinate to $\bW$ such that
$\ff(\fV)\subset\fU$.

 Given a morphism of formal schemes $\ff\:\fY\rarrow\fX$ and an open
covering $\bW$ of the scheme $\fX$, consider the open covering $\bT$
of the scheme $\fY$ defined by the rule
$$
 \bT=\{\,\ff^{-1}(\fU)\subset\fY\mid\fU\subset\fX \text{ is an affine
 open formal subscheme subordinate to~$\bW$}\,\}.
$$
 Then the morphism~$\ff$ is $(\bW,\bT)$\+coaffine.
 Moreover, if the morphism~$\ff$ happens to be affine, then it is also
$(\bW,\bT)$\+affine for the open covering $\bT$ defined above.

 Let $\ff\:\fY\rarrow\fX$ be a tight morphism of locally Noetherian
formal schemes (in the sense of the definition
in Section~\ref{morphisms-of-formal-schemes-subsecn}).
 Then the functor of inverse image of quasi-coherent torsion sheaves
\begin{equation} \label{inverse-image-qcoh-tors}
 \ff^*\:\fX\Tors\lrarrow\fY\Tors
\end{equation}
is constructed as follows.

 Pick open coverings $\bW$ of the formal scheme $\fX$ and $\bT$ of
the formal scheme $\fY$ such that the morphism~$\ff$ is
$(\bW,\bT)$\+coaffine.
 It is convenient to think of the objects of $\fX\Tors$ as
quasi-coherent torsion presheaves (in the sense of
Section~\ref{qcoh-torsion-sheaves-subsecn}) on the topology base of
the formal scheme $\fX$ consisting of all the affine open formal
subschemes $\fU\subset\fX$ subordinate to~$\bW$.
 Similarly, let us think of the objects of $\fY\Tors$ as quasi-coherent
torsion presheaves on the topology base of the formal scheme $\fY$
consisting of all the affine open formal subschemes $\fV\subset\fY$
subordinate to~$\bT$.
 This point of view is based on the discussion in
Section~\ref{qcoh-torsion-sheaves-subsecn}.

 Let $\cM$ be a quasi-coherent torsion sheaf on $\fX$ and
$\fV\subset\fY$ be an affine open formal subscheme subordinate to~$\bT$.
 Pick an affine open formal subscheme $\fU\subset\fX$ subordinate to
$\bW$ such that $\ff(\fV)\subset\fU$.
 Denote by $f_{\fU,\fV}\:\fO_\fX(\fU)\rarrow\fO_\fY(\fV)$ the map of
adic topological rings corresponding to the morphism of formal
schemes $\ff|_\fV^\fU\:\fV\rarrow\fU$, as in
Section~\ref{morphisms-of-formal-schemes-subsecn}.
 Put
$$
 (\ff^*\cM)(\fV)=f_{\fU,\fV}^*(\cM(\fU))=
 \fO_\fY(\fV)\ot_{\fO_\fX(\fU)}\cM(\fU),
$$
where the notation~$f^*$ comes from
Section~\ref{prelim-change-of-scalars-subsecn}.
 By Lemma~\ref{tight-co-extension-of-scalars}(a),
the $\fO_\fY(\fV)$\+mod\-ule $f_{\fU,\fV}^*(\cM(\fU))$ is torsion.

 Let $\fU'\subset\fU$ be another affine open formal subscheme in $\fX$
such that $\ff(\fV)\subset\fU'$.
 Let $f_{\fU',\fV}\:\fO_\fX(\fU')\rarrow\fO_\fY(\fV)$ be
the corresponding map of adic topological rings.
 Then we have
\begin{multline*}
 f_{\fU',\fV}^*(\cM(\fU'))=\fO_\fY(\fV)\ot_{\fO_\fX(\fU')}\cM(\fU')
 \simeq \fO_\fY(\fV)\ot_{\fO_\fX(\fU')}
 (\fO_\fX(\fU')\ot_{\fO_\fX(\fU)}\cM(\fU)) \\
 \simeq\fO_\fY(\fV)\ot_{\fO_\fX(\fU)}\cM(\fU)=
 f_{\fU,\fV}^*(\cM(\fU)),
\end{multline*}
since $\cM(\fU')\simeq\fO_\fX(\fU')\ot_{\fO_\fX(\fU)}\cM(\fU)$ by
the quasi-coherence axiom~(ii) from
Section~\ref{qcoh-torsion-sheaves-subsecn}.
 In this sense at least, the torsion $\fO_\fY(\fV)$\+module
$f_{\fU,\fV}^*(\cM(\fU))$ does not depend on the choice of an affine
open formal subscheme $\fU\subset\fX$ subordinate to $\bW$ such that
$\ff(\fV)\subset\fU$.

 Let us first assume that the formal scheme $\fX$ is semi-separated.
 Then, for any pair of affine open formal subschemes $\fU_1$ and
$\fU_2\subset\fX$ subordinate to $\bW$ such that $\ff(\fV)\subset\fU_1$
and $\ff(\fV)\subset\fU_2$, we have an affine open formal subscheme
$\fU=\fU_1\cap\fU_2$ subordinate to $\bW$ such that $\fU\subset\fU_1$,
\,$\fU\subset\fU_2$, and $\ff(\fV)\subset\fU$.
 In view of the computation in the previous paragraph, it follows that
the torsion $\fO_\fY(\fV)$\+module $(\ff^*\cM)(\fV)$ is well-defined.

 Let us check that $\ff^*\cM$ is a quasi-coherent torsion presheaf on
the topology base of $\fY$ consisting of all the affine open formal
subschemes $\fV\subset\fY$ subordinate to~$\bT$.
 Indeed, we have explained already that the torsion axiom~(i) holds.
 To establish the quasi-coherence axiom, let $\fV'\subset\fV$ be
a pair of affine open formal subschems in $\fY$ subordinate to $\bT$,
and let $\fU\subset\fX$ be an affine open formal subscheme
subordinate to $\bW$ such that $\ff(\fV)\subset\fU$.
 Let $f_{\fU,\fV}\:\fO_\fX(\fU)\rarrow\fO_\fY(\fV)$ and
$f_{\fU,\fV'}\:\fO_\fX(\fU)\rarrow\fO_\fY(\fV')$ be the corresponding
maps of adic topological rings.
 Then we have
\begin{multline*}
 f_{\fU,\fV'}^*(\cM(\fU))=\fO_\fY(\fV')\ot_{\fO_\fX(\fU)}\cM(\fU) \\
 \simeq \fO_\fY(\fV')\ot_{\fO_\fY(\fV)}
 (\fO_\fY(\fV)\ot_{\fO_\fX(\fU)}\cM(\fU))
 =\fO_\fY(\fV')\ot_{\fO_\fY(\fV)}f_{\fU,\fV}^*(\cM(\fU)),
\end{multline*}
as desired.

 In the case of a semi-separated formal scheme $\fX$, the inverse
image functor~$\ff^*$ is constructed.
 One can easily check (by refining the coverings) that
the functor~$\ff^*$ does not depend on the choice of the open
coverings~$\bW$ and~$\bT$.

 Now we are ready to discuss the case of a non-semi-separated formal
scheme~$\fX$.
 Let $\fV\subset\fY$ be an affine open formal subscheme subordinate
to $\bT$, and let $\fU_1$ and $\fU_2\subset\fX$ be two affine open
formal subschemes subordinate to $\bW$ such that $\ff(\fV)\subset\fU_1$
and $\ff(\fV)\subset\fU_2$.
 Let $f_{\fU_1,\fV}\:\fO_\fX(\fU_1)\rarrow\fO_\fY(\fV)$ and
$f_{\fU_2,\fV}\:\fO_\fX(\fU_2)\rarrow\fO_\fY(\fV)$ be the corresponding
maps of adic topological rings. 
 Put $(\ff^*\cM)_1(\fV)=f_{\fU_1,\fV}^*(\cM(\fU_1))$
and $(\ff^*\cM)_2(\fV)=f_{\fU_2,\fV}^*(\cM(\fU_1))$.
 We need to construct a natural isomorphism of torsion
$\fO_\fY(\fV)$\+modules $(\ff^*\cM)_1(\fV)\simeq(\ff^*\cM)_2(\fV)$.

 The open formal subscheme $\fU_1\cap\fU_2\subset\fX$ need not be
affine, but it is always semi-separated (as an open formal subscheme
of an affine, hence semi-separated, formal scheme $\fU_1$ or~$\fU_2$).
 We follow the exposition in~\cite[Section~3.3]{Pcosh} with some
details added.

 Put $\fX'=\fU_1\cap\fU_2$ and $\fY'=\fV$, and denote by
$\ff'\:\fY'\rarrow\fX'$ the induced morphism of formal schemes.
 Put $\cM'=\cM|_{\fX'}$.
 Furthermore, consider the trivial open covering
$\bW'=\bW_{\fX'}=\{\fX'\}$ of the formal scheme~$\fX'$.
 Pick any open covering $\bT'$ of the formal scheme $\fY'$ such that
the morphism~$\ff'$ is $(\bW',\bT')$\+coaffine.
 Applying the construction above to the morphism~$\ff'$ and the open
coverings $\bW'$ and $\bT'$, we obtain a quasi-coherent torsion sheaf
$\ff'{}^*(\cM')$ on~$\fY'$.

 For each index $i=1$ or~$2$ and any affine open formal subscheme
$\fV'\subset\fY'$, let $f_{\fU_i,\fV'}\:\fO_\fX(\fU_i)\rarrow
\fO_\fX(\fV')$ be the map of adic topological rings corresponding to
the morphism of affine formal schemes $\ff|_{\fV'}^{\fU_i}\:\fV'
\rarrow\fU_i$.
 Put $\ff^*(\cM)'_i(\fV')=f_{\fU_i,\fV'}^*(\cM(\fU_i))=
\fO_\fY(\fV')\ot_{\fO_\fX(\fU_i)}\cM(\fU_i)$.
 In other words, consider the morphisms of affine formal schemes
$\ff'_i=\ff|_\fV^{\fU_i}\:\fV\rarrow\fU_i$, and put
$\cM'_i=\cM|_{\fU_i}$.
 Consider the trivial open coverings $\bW_{\fU_i}=\{\fU_i\}$ and
$\bT_\fV=\{\fV\}$ of the affine formal schemes $\fU_i$ and~$\fV$.
 Then the morphisms~$\ff_i$ are $(\bW_{\fU_i},\bT_\fV)$\+coaffine,
and we have $\ff^*(\cM)'_i=\ff'_i{}^*(\cM'_i)$.
 So $\ff^*(\cM)'_1$ and $\ff^*(\cM)'_2$ are quasi-coherent torsion
presheaves on the topology base of all affine open formal subschemes
in $\fY'$, as it is clear from the discussion above.
 Thus both $\ff^*(\cM)'_1$ and $\ff^*(\cM)'_2$ can be viewed as
quasi-coherent torsion sheaves on~$\fY'$.

 Furthermore, the restrictions of all the three quasi-coherent
torsion sheaves $\ff^*(\cM)'_1$, \,$\ff'{}^*(\cM')$, and
$\ff^*(\cM)'_2$ on the topology base of all affine open formal
subschemes subordinate to $\bT'$ in $\fY'$ are naturally isomorphic to
each other, as it is also clear from the preceding discussion.
 Since a sheaf is determined by its restriction to any topology base
(according to Section~\ref{co-sheaves-of-modules-subsecn}), we have
natural isomorphisms of quasi-coherent torsion sheaves
$\ff^*(\cM)'_1\simeq\ff'{}^*(\cM')\simeq\ff^*(\cM)'_2$ on~$\fY'$.
 Finally, we have
$$
 \ff^*(\cM)_1(\fV)=\ff^*(\cM)'_1(\fV)\simeq\ff'{}^*(\cM')(\fV)
 \simeq\ff^*(\cM)'_2(\fV)=\ff^*(\cM)_2(\fV).
$$
 Hence the promised natural isomorphism of $\fO_\fY(\fV)$\+modules
$\ff^*(\cM)_1(\fV)\simeq\ff^*(\cM)_2(\fV)$.

 The inverse image functor $\ff^*$ \,\eqref{inverse-image-qcoh-tors}
is now constructed for all tight morphisms of locally Noetherian
formal schemes.
 Once again, one can easily check that the functor~$\ff^*$ does not
depend on the choice of the open coverings~$\bW$ and~$\bT$.
 The functor $\ff^*$ \,\eqref{inverse-image-qcoh-tors} is right
exact and preserves infinite direct sums.
 For a flat tight morphism of locally Noetherian formal schemes~$\ff$,
the inverse image functor~$\ff^*$ is exact (in view of
Lemma~\ref{flat-map-of-adic-rings-direct-image-lemma}(a)).

 Any open immersion morphism of locally Noetherian formal schemes
$\fj\:\fY\rarrow\fX$ is tight.
 In this setting, the inverse image functor~$\fj^*$ agrees with
the restriction functor $\cM\longmapsto\cM|_\fY$ discussed in
the beginning of Section~\ref{direct-images-of-qcoh-tors-subsecn}.
 So one has a natural isomorphism $\fj^*\cM\simeq\cM|_\fY$ of
quasi-coherent torsion sheaves on $\fY$ for any quasi-coherent
torsion sheaf $\cM$ on~$\fX$.

 For any quasi-compact tight morphism of locally Noetherian formal
schemes $\ff\:\fY\rarrow\fX$, the direct image functor $\ff_*$
\,\eqref{direct-image-qcoh-tors} is right adjoint to the inverse
image functor $\ff^*$~\,\eqref{inverse-image-qcoh-tors}.
 In other words, for any quasi-coherent torsion sheaf $\cM$ on $\fX$
and any quasi-coherent torsion sheaf $\cN$ on $\fY$, there is a natural
isomorphism of abelian groups
\begin{equation} \label{qcoh-tors-inverse-direct-adjunction}
 \Hom_\fY(\ff^*\cM,\cN)\simeq\Hom_\fX(\cM,\ff_*\cN),
\end{equation}
where $\Hom_\fY({-},{-})$ and $\Hom_\fX({-},{-})$ denote the groups
of morphisms in the categories $\fY\Tors$ and $\fX\Tors$.

 More generally, for any tight morphism of locally Noetherian formal
schemes $\ff\:\fY\rarrow\fX$, the adjunction isomorphism
\begin{equation} \label{sheaves-inverse-direct-adjunction}
 \Hom_{\fO_\fY}(\ff^*\cM,\cN)\simeq\Hom_{\fO_\fX}(\cM,\ff_*\cN)
\end{equation}
of the groups of morphisms in the categories $(\fY,\fO_\fY)\Sh$ and
$(\fX,\fO_\fX)\Sh$ holds for every quasi-coherent torsion sheaf $\cM$
on $\fX$ and every sheaf of $\fO_\fY$\+modules $\cN$ on~$\fY$.
 Indeed, both the abelian groups
in~\eqref{sheaves-inverse-direct-adjunction} are naturally isomorphic
to the group of all compatible systems of $\fO_\fX(\fU)$\+module
maps $\cM(\fU)\rarrow\cN(\fV)$ defined for all affine open formal
subschemes $\fU\subset\fX$ subordinate to $\bW$ and all affine open
formal subschemes $\fV\subset\fX$ subordinate to $\bT$ such that
$\ff(\fV)\subset\fU$.

\subsection{Inverse images of contraherent cosheaves of contramodules}
\label{inverse-images-of-ctrh-contramod-subsecn}
 Let $\fX$ be a locally Noetherian formal scheme with an open covering
$\bW$, let $\fY$ be a locally Noetherian formal scheme with an open
covering $\bT$, and let $\ff\:\fY\rarrow\fX$ be a $(\bW,\bT)$\+coaffine
very flat tight morphism of formal schemes (in the sense of
the definitions in Sections~\ref{morphisms-of-formal-schemes-subsecn}
and~\ref{inverse-images-of-qcoh-tors-subsecn}).
 Then the functor of inverse image of locally contraherent cosheaves
of contramodules
\begin{equation} \label{inverse-image-lcth-W-T-contramod}
 \ff^!\:\fX\Lcth_\bW\lrarrow\fY\Lcth_\bT
\end{equation}
is constructed as follows.

 Let $\fP$ be a $\bW$\+locally contraherent cosheaf of contramodules
on $\fX$ and $\fV\subset\fY$ be an affine open formal subscheme
subordinate to~$\bT$.
 Pick an affine open formal subscheme $\fU\subset\fX$ subordinate to
$\bW$ such that $\ff(\fV)\subset\fU$.
 Denote by $f_{\fU,\fV}\:\fO_\fX(\fU)\rarrow\fO_\fY(\fV)$ the map of
adic topological rings corresponding to the morphism of formal
schemes $\ff|_\fV^\fU\:\fV\rarrow\fU$, as in
Section~\ref{morphisms-of-formal-schemes-subsecn}.
 Put
$$
 (\ff^!\fP)[\fV]=f_{\fU,\fV}^!(\fP[\fU])=
 \Hom_{\fO_\fX(\fU)}(\fO_\fY(\fV),\fP[\fU]),
$$
where the notation~$f^!$ comes from
Sections~\ref{prelim-change-of-scalars-subsecn}
and~\ref{prelim-coextension-of-scalars-quotseparated-subsecn}.
 By Lemma~\ref{tight-co-extension-of-scalars}(b) or
Proposition~\ref{flat-quotseparated-coextension-of-scalars-prop},
the $\fO_\fY(\fV)$\+module $f_{\fU,\fV}^!(\fP[\fU])$ is a contramodule.

 Let $\fU'\subset\fU$ be another affine open formal subscheme in $\fX$
such that $\ff(\fV)\subset\fU'$.
 Let $f_{\fU',\fV}\:\fO_\fX(\fU')\rarrow\fO_\fY(\fV)$ be
the corresponding map of adic topological rings.
 Then we have
\begin{multline*}
 f_{\fU',\fV}^!(\fP[\fU'])=\Hom_{\fO_\fX(\fU')}(\fO_\fY(\fV),\fP[\fU'])
 \\ \simeq \Hom_{\fO_\fX(\fU')}\bigl(\fO_\fY(\fV),
 \Hom_{\fO_\fX(\fU)}(\fO_\fX(\fU'),\fP[\fU])\bigr) \\
 \simeq \Hom_{\fO_\fX(\fU)}(\fO_\fY(\fV),\fP[\fU])
 = f_{\fU,\fV}^!(\fP[\fU]),
\end{multline*}
since $\fP[\fU']\simeq\Hom_{\fO_\fX(\fU)}(\fO_\fX(\fU'),\fP[\fU])$ by
the contraherence axiom~(iv) from
Section~\ref{contraherent-cosheaves-of-contramods-subsecn}.
 In this sense at least, the contramodule $\fO_\fY(\fV)$\+module
$f_{\fU,\fV}^!(\fP[\fU])$ does not depend on the choice of an affine
open formal subscheme $\fU\subset\fX$ subordinate to $\bW$ such that
$\ff(\fV)\subset\fU$.
 Similarly to the discussion in
Section~\ref{inverse-images-of-qcoh-tors-subsecn}, this proves that
the contramodule $\fO_\fY(\fV)$\+module $(\ff^!\fP)[\fV]$ is
well-defined if the formal scheme $\fX$ is semi-separated.

 Assuming first that $\fX$ is semi-separated, let us check that
$\ff^!\fP$ is a contraherent copresheaf of contramodules on
the topology base of $\fY$ consisting of all the affine open
formal subschemes $\fV\subset\fY$ subordinate to~$\bT$.
 Indeed, we have explained already that the contramoduleness
axiom~(iii) holds.
 The contraadjustedness axiom~(v) holds by
Proposition~\ref{colocalization-of-contraadjusted-contramodule-prop}.
 To establish the contraherence axiom, let $\fV'\subset\fV$ be
a pair of affine open formal subschems in $\fY$ subordinate to $\bT$,
and let $\fU\subset\fX$ be an affine open formal subscheme
subordinate to $\bW$ such that $\ff(\fV)\subset\fU$.
 Let $f_{\fU,\fV}\:\fO_\fX(\fU)\rarrow\fO_\fY(\fV)$ and
$f_{\fU,\fV'}\:\fO_\fX(\fU)\rarrow\fO_\fY(\fV')$ be the corresponding
maps of adic topological rings.
 Then we have
\begin{multline*}
 f_{\fU,\fV'}^!(\fP[\fU])=\Hom_{\fO_\fX(\fU)}(\fO_\fY(\fV'),\fP[\fU])
 \\ \simeq \Hom_{\fO_\fY(\fV)}\bigl(\fO_\fY(\fV'),
 \Hom_{\fO_\fX(\fU)}(\fO_\fY(\fV),\fP[\fU])\bigr) \\ =
 \Hom_{\fO_\fY(\fV)}\bigl(\fO_\fY(\fV'),f_{\fU,\fV}^!(\fP[\fU])\bigr),
\end{multline*}
as desired.
 Following the discussion in
Sections~\ref{contraherent-cosheaves-of-contramods-subsecn}\+-%
\ref{locally-contraherent-cosheaves-of-contramods-subsecn},
the copresheaf $\ff^!\fP$ can be uniquely extended to
a $\bT$\+locally contraherent cosheaf of contramodules on~$\fY$.

 In the case of a semi-separated formal scheme $\fX$, the inverse
image functor~$\ff^!$ is constructed.
 One can easily check that, as the coverings $\bW$ and $\bT$ vary,
the functors~$\ff^!$ agree with each other, in the sense that for
a given locally contraherent cosheaf of contramodules $\fP$ on $\fX$,
the locally contraherent cosheaf of contramodules $\ff^!\fP$ on $\fY$
does not depend on the choice of a pair of open coverings $\bW$ and
$\bT$ for which the cosheaf $\fP$ is $\bW$\+locally contraherent.

 Following~\cite[Section~3.3]{Pcosh} and
Section~\ref{inverse-images-of-qcoh-tors-subsecn} above, let us discuss
the case of a non-semi-separated formal scheme~$\fX$.
 We keep the notation of
Section~\ref{inverse-images-of-qcoh-tors-subsecn} concerning
the open formal subschemes, the open coverings, and the ring maps.
 In that notation, put
$(\ff^!\fP)_1[\fV]=f_{\fU_1,\fV}^!(\fP[\fU_1])$
and $(\ff^!\fP)_2[\fV]=f_{\fU_2,\fV}^!(\fP[\fU_2])$.
 We need to construct a natural isomorphism of contramodule
$\fO_\fY(\fV)$\+modules $(\ff^!\fP)_1[\fV]\simeq(\ff^!\fP)_2[\fV]$.

 Put $\fP=\fP|_{\fX'}$.
 Applying the construction above to the morphism $\ff'\:\fY'\rarrow
\fX'$ and the open coverings $\bW'$ and $\bT'$, we obtain
a $\bT$\+locally contraherent cosheaf of contramodules
$\ff'{}^!(\fP')$ on~$\fY'$.

 For each index $i=1$ or~$2$ and any affine open formal subscheme
$\fV'\subset\fY'$, put $\ff^!(\fP)'_i[\fV']=
f_{\fU_i,\fV'}^!(\fP[\fU_i])=
\Hom_{\fO_\fX(\fU_i)}(\fO_\fY(\fV'),\fP[\fU_i])$.
 In other words, setting $\fP'_i=\fP|_{\fU_i}$, we put
$\ff^!(\fP)'_i=\ff'_i{}^!(\fP'_i)$.
 Then $\ff^!(\fP)'_1$ and $\ff^!(\fP)'_2$ are contraherent copresheaves
of contramodules on the topology base of all affine open formal
subschemes in $\fY'$, as it is clear from the discussion above.
 So both $\ff^!(\fP)'_1$ and $\ff^!(\fP)'_2$ can be viewed as
(globally) contraherent cosheaves of contramodules on~$\fY'$.

 Furthermore, the restrictions of all the three locally contraherent
cosheaves of contramodules $\ff^!(\fP)'_1$, \,$\ff'{}^!(\fP')$, and
$\ff^!(\fP)'_2$ on the topology base of all affine open formal
subschemes subordinate to $\bT'$ in $\fY'$ are naturally isomorphic to
each other, as it is also clear from the preceding discussion.
 Since a cosheaf is determined by its restriction to any topology base
(according to Section~\ref{co-sheaves-of-modules-subsecn}), we have
natural isomorphisms of (locally) contraherent cosheaves of
contramodules $\ff^!(\fP)'_1\simeq\ff'{}^!(\fP')\simeq\ff^!(\fP)'_2$
on~$\fY'$.
 Finally, we have
$$
 \ff^!(\fP)_1[\fV]=\ff^!(\fP)'_1[\fV]\simeq\ff'{}^!(\fP')[\fV]
 \simeq\ff^!(\fP)'_2[\fV]=\ff^!(\fP)_2[\fV].
$$
 Hence the promised natural isomorphism of $\fO_\fY(\fV)$\+modules
$\ff^!(\fP)_1[\fV]\simeq\ff^!(\fP)_2[\fV]$.

 The inverse image functor $\ff^!$
\,\eqref{inverse-image-lcth-W-T-contramod} is now constructed for all
very flat tight morphisms of locally Noetherian formal schemes.
 Once again, one can easily check that, as the open coverings $\bW$
and $\bT$ vary, the related functors~$\ff^!$ agree.
 Passing to the inductive limit of categories with respect to
the refinements of the open coverings $\bW$ and $\bT$, we obtain
the inverse image functor
\begin{equation} \label{inverse-image-lcth-contramod-nocovering}
 \ff^!\:\fX\Lcth\lrarrow\fY\Lcth.
\end{equation}

 Now let $\fX$ be a locally Noetherian formal scheme with an open
covering $\bW$, let $\fY$ be a locally Noetherian formal scheme with
an open covering $\bT$, and let $\ff\:\fY\rarrow\fX$ be
a $(\bW,\bT)$\+coaffine flat tight morphism of formal schemes
(in the sense of the definitions in
Sections~\ref{morphisms-of-formal-schemes-subsecn}
and~\ref{inverse-images-of-qcoh-tors-subsecn}).
 Then the functor of inverse image of locally cotorsion
locally contraherent cosheaves of contramodules
\begin{equation} \label{inverse-image-lcth-lct-W-T-contramod}
 \ff^!\:\fX\Lcth^\lct_\bW\lrarrow\fY\Lcth^\lct_\bT
\end{equation}
is constructed in the way very similar to the discussion above.
 The reference to
Proposition~\ref{colocalization-of-cotorsion-contramodule-prop}
is relevant.
 As the open coverings $\bW$ and $\bT$ vary, the related
functors $\ff^!$ \,\eqref{inverse-image-lcth-lct-W-T-contramod} agree.
 Passing to the inductive limit of categories with respect to
the refinements of the open coverings $\bW$ and $\bT$, we obtain
the inverse image functor
\begin{equation} \label{inverse-image-lcth-lct-contramod-nocovering}
 \ff^!\:\fX\Lcth^\lct\lrarrow\fY\Lcth^\lct.
\end{equation}
 The functors of inverse image of locally contraadjusted locally
contraherent cosheaves of
contramodules~\eqref{inverse-image-lcth-contramod-nocovering}
and of locally cotorsion locally contraherent cosheaves of
contramodules~\eqref{inverse-image-lcth-lct-contramod-nocovering}
agree with each other wherever both of them are defined.

 The inverse image functors~$\ff^!$
\,\eqref{inverse-image-lcth-W-T-contramod},
\eqref{inverse-image-lcth-contramod-nocovering},
\eqref{inverse-image-lcth-lct-W-T-contramod},
and~\eqref{inverse-image-lcth-lct-contramod-nocovering} are exact
(as functors of exact categories).
 The references to
Corollaries~\ref{colocalization-contraadjusted-exactness-cor}
and~\ref{colocalization-cotorsion-exactness-cor} are relevant here.
 The inverse image functors~\eqref{inverse-image-lcth-W-T-contramod}
and~\eqref{inverse-image-lcth-lct-W-T-contramod} also preserve
infinite products.

 Any open immersion morphism of locally Noetherian formal schemes
$\fj\:\fY\rarrow\fX$ is tight and very flat.
 Given an open covering $\bW$ of $\fX$, the restriction
$\bT=\bW|_\fY=\{\fY\cap\fW\mid\fW\in\bW\}$ of the open covering $\bW$
to $\fY$ has the property that the morphism~$\fj$ is
$(\bW,\bT)$\+coaffine.
 In this setting, the inverse image functors~$\fj^!$
\,\eqref{inverse-image-lcth-contramod-nocovering}
and~\eqref{inverse-image-lcth-lct-contramod-nocovering} agree with
the restriction functor $\fP\longmapsto\fP|_\fY$ defined and discussed
in the beginning of
Section~\ref{direct-images-of-ctrh-contramod-subsecn}.
 So one has a natural isomorphism $\fj^!\fP\simeq\fP|_\fY$ of
$\bW|_\fY$\+locally contraherent cosheaves of contramodules on $\fY$
for any $\bW$\+locally contraherent cosheaf of contramodules $\fP$
on~$\fX$.

 For any $(\bW,\bT)$\+affine $(\bW,\bT)$\+coaffine very flat tight
morphism of locally Noetherian formal schemes $\ff\:\fY\rarrow\fX$,
the direct image functor $\ff_!$
\,\eqref{loc-ctrh-contramod-W-T-affine-direct-image} is left adjoint
to the inverse image functor
$\ff^!$~\,\eqref{inverse-image-lcth-W-T-contramod}.
 In other words, for any $\bW$\+locally contraherent cosheaf of
contramodules $\fP$ on $\fX$ and any $\bT$\+locally contraherent
cosheaf of contramodules $\fQ$ on $\fY$, there is a natural
isomorphism of abelian groups
\begin{equation} \label{lcth-contramod-direct-inverse-adjunction}
 \Hom^\fY(\fQ,\ff^!\fP)\simeq\Hom^\fX(\ff_!\fQ,\fP),
\end{equation}
where $\Hom^\fY({-},{-})$ and $\Hom^\fX({-},{-})$ denote the groups
of morphisms in the categories $\fY\Lcth$ and $\fX\Lcth$.

 Similarly, for any $(\bW,\bT)$\+affine $(\bW,\bT)$\+coaffine flat
tight morphism of locally Noetherian formal schemes
$\ff\:\fY\rarrow\fX$, the direct image functor $\ff_!$
\,\eqref{loc-ctrh-lct-contramod-W-T-affine-direct-image}
is left adjoint to the inverse image functor
$\ff^!$~\,\eqref{inverse-image-lcth-lct-W-T-contramod}.
 So the natural isomorphism of abelian
groups~\eqref{lcth-contramod-direct-inverse-adjunction}
holds for any locally cotorsion $\bW$\+locally contraherent cosheaf of
contramodules $\fP$ on $\fX$ and any locally cotorsion $\bT$\+locally
contraherent cosheaf of contramodules $\fQ$ on $\fY$ in this context.

 More generally, for any very flat tight morphism of locally Noetherian
formal schemes $\ff\:\fY\rarrow\fX$, the adjunction isomorphism
\begin{equation} \label{cosheaves-direct-inverse-adjunction}
 \Hom^{\fO_\fY}(\fQ,\ff^!\fP)\simeq\Hom^{\fO_\fX}(\ff_!\fQ,\fP)
\end{equation}
of the groups of morphisms in the categories $(\fY,\fO_\fY)\Cosh$ and
$(\fX,\fO_\fX)\Cosh$ holds for every locally contraherent cosheaf of
contramodules $\fP$ on $\fX$ and every cosheaf of $\fO_\fY$\+modules
$\fQ$ on~$\fY$.
 Similarly, for any flat tight morphism of locally Noetherian
formal schemes $\ff\:\fY\rarrow\fX$, the adjunction
isomorphism~\eqref{cosheaves-direct-inverse-adjunction} holds
for every locally cotorsion locally contraherent cosheaf of
contramodules $\fP$ on $\fX$ and every cosheaf of $\fO_\fY$\+modules
$\fQ$ on~$\fY$.

 Indeed, both the abelian groups
in~\eqref{cosheaves-direct-inverse-adjunction} are naturally isomorphic
to the group of all compatible systems of $\fO_\fX(\fU)$\+module
maps $\fQ[\fV]\rarrow\fP[\fU]$ defined for all affine open formal
subschemes $\fU\subset\fX$ subordinate to $\bW$ and all affine open
formal subschemes $\fV\subset\fX$ subordinate to $\bT$ such that
$\ff(\fV)\subset\fU$.
 Here we presume the cosheaf $\fP$ on $\fX$ to be
$\bW$\+locally contraherent.

\subsection{Base change}  \label{base-change-subsecn}
 In this paper, we only consider locally Noetherian formal schemes.
 The fibered products of (locally) Noetherian (formal) schemes need
not be (locally) Noetherian, in general.
 For this reason, we restrict ourselves to fibered products where
one of the two morphisms is an open immersion in this section.

\begin{lem} \label{tors-openimmers-dir-tight-inv-basechange-lemma}
 Let\/ $\ff\:\fY\rarrow\fX$ be a tight morphism of locally Noetherian
formal schemes and\/ $\fW\subset\fX$ be an open formal subscheme with
the open immersion morphism\/ $\fj\:\fW\rarrow\fX$.
 Consider the open formal subscheme\/ $\fT=\ff^{-1}(\fW)\subset\fY$,
and denote by\/ $\fj'\:\fT\rarrow\fY$ and\/ $\ff'\:\fT\rarrow\fW$
the natural morphisms.
 Assume that the morphism\/~$\fj$ is quasi-compact (then
the morphism\/~$\fj'$ is quasi-compact as well).
 Let\/ $\cL$ be a quasi-coherent torsion sheaf on\/~$\fW$.
 Then there is a natural morphism
\begin{equation} \label{tors-openimmers-dir-tight-inv-basechange-eqn}
 \ff^*\fj_*\cL\lrarrow\fj'_*\ff'{}^*\cL
\end{equation}
of quasi-coherent torsion sheaves on\/~$\fY$.
 The morphism of quasi-coherent torsion
sheaves~\eqref{tors-openimmers-dir-tight-inv-basechange-eqn}
is an isomorphism whenever either \par
\textup{(a)} the morphism of formal schemes\/~$\fj$ is affine; or \par
\textup{(b)} the morphism of formal schemes\/~$\ff$ is flat.
\end{lem}

\begin{proof}
 The direct image functor $\fj_*\:\fW\Tors\rarrow\fX\Tors$ is right
adjoint to the inverse image/restriction functor $\fj^*\:\fX\Tors\rarrow
\fX\Tors$ (by formula~\eqref{qcoh-tors-inverse-direct-adjunction}),
so we have the adjunction (counit) morphism $\fj^*\fj_*\cL\rarrow\cL$
in $\fW\Tors$.
 Applying the inverse image functor~$\ff'{}^*$, we obtain a morphism
$\ff'{}^*\fj^*\fj_*\cL\rarrow\ff'{}^*\cL$ in $\fT\Tors$.
 The two morphisms of formal schemes $\fj\ff'$ and $\ff\fj'\:\fT
\rarrow\fX$ are equal to each other, hence a natural isomorphism of
the compositions of inverse image functors $\ff'{}^*\fj^*\simeq
\fj'{}^*\ff^*$.
 So we arrive at a morphism $\fj'{}^*\ff^*\fj_*\cL\rarrow\ff'{}^*\cL$
in $\fT\Tors$, which corresponds to the desired
morphism~\eqref{tors-openimmers-dir-tight-inv-basechange-eqn}
in $\fY\Tors$ under the adjunction
isomorphism~\eqref{qcoh-tors-inverse-direct-adjunction}.

 Alternatively, one can start from the adjunction (unit) morphism of
sheaves of $\fO_\fW$\+modules $\cL\rarrow\ff'_*\ff'{}^*\cL$ provided
by formula~\eqref{sheaves-inverse-direct-adjunction}.
 Applying the direct image functor $\fj_*\:(\fW,\fO_\fW)\Sh\rarrow
(\fX,\fO_\fX)\Sh$, we obtain a morphism $\fj_*\cL\rarrow
\fj_*\ff'_*\ff'{}^*\cL$ in $(\fX,\fO_\fX)\Sh$.
 The two morphisms of formal schemes $\fj\ff'$ and $\ff\fj'\:\fT
\rarrow\fX$ are equal to each other, hence a natural isomorphism of
the compositions of direct image functors $\fj_*\ff'_*\simeq
\ff_*\fj'_*$.
 So we arrive at a morphism $\fj_*\cL\rarrow\ff_*\fj'_*\ff'{}^*\cL$
in $(\fX,\fO_\fX)\Sh$, which corresponds to the desired
morphism~\eqref{tors-openimmers-dir-tight-inv-basechange-eqn}
in $\fY\Tors$ under the adjunction
isomorphism~\eqref{sheaves-inverse-direct-adjunction}.

 It remains to prove the base change isomorphism assertions~(a)
and~(b).
 The question is local in both $\fX$ and $\fY$, so without loss
of generality we can assume $\fX$ and $\fY$ to be affine
formal schemes.
 Then we have complete, separated adic Noetherian rings $\fR$ and $\fS$
such that $\fX=\Spf\fR$ and $\fY=\Spf\fS$.
 The morphism $\ff\:\fY\rarrow\fX$ corresponds to a tight morphism
of adic topological rings $f\:\fR\rarrow\fS$.

 Under the assumption of part~(a), the formal schemes $\fW$ and $\fT$
are affine as well; so we have $\fW=\Spf\widetilde\fR$ and
$\fT=\Spf\widetilde\fS$ for some complete, separated adic
Noetherian rings $\widetilde\fR$ and~$\widetilde\fS$ (cf.\
Corollary~\ref{adic-topological-rings-complet-noether-local-cor}(a)).
 Let $j\:\fR\rarrow\widetilde\fR$, \ $j'\:\fS\rarrow\widetilde\fS$,
and $f'\:\widetilde\fR\rarrow\widetilde\fS$ be the continuous
homomorphisms of adic topological rings corresponding to~$\fj$,
$\fj'$, and~$\ff'$.
 We have $\widetilde\fS\simeq\Lambda(\widetilde\fR\ot_\fR\fS)$, where
the tensor product ring $\widetilde\fR\ot_\fR\fS$ is endowed with
the tensor product topology.
 The ring map~$f'$ is tight, while the ring maps~$j$ and~$j'$ are
formal open immersions (by
Lemmas~\ref{induced-map-of-complete-rings-lemma}(a),
\ref{tight-flat-adic-ring-base-change-context}(a),
and~\ref{open-immersion-adic-ring-base-change-context}(a)).

 The quasi-coherent torsion sheaf $\cL$ on $\fW=\Spf\widetilde R$
corresponds to a torsion $\widetilde\fR$\+module $\cM=\cL(\fW)$.
 The assertion of part~(a) claims a natural isomorphism of
torsion $\fS$\+modules
$$
 \fS\ot_\fR\cM\simeq \widetilde\fS\ot_{\widetilde\fR}\cM.
$$
 Obviously, $\fS\ot_\fR\cM\simeq(\widetilde\fR\ot_\fR\fS)
\ot_{\widetilde\fR}\cM$, and it remains to refer to
Corollary~\ref{derived-completion-tensor-with-torsion-isom-cor}.
 Notice that $\widetilde\fR\rarrow\widetilde\fR\ot_\fR\fS$ is
a tight continuous ring map by
Lemma~\ref{tight-flat-adic-ring-base-change-context}(a),
so the tensor product topology on
$\widetilde\fR\ot_\fR\fS$ coincides with the topology induced
by the adic topology of $\fR$ by
Lemma~\ref{continuous-and-tight-ring-map-lemma}(2) or~(3).

 Under the assumption of part~(b), the map of adic Noetherian rings
$f\:\fR\rarrow\fS$ is tight and flat, while the formal schemes $\fW$
and $\fT$ are quasi-compact and semi-separated.
 Let $\fW=\bigcup_{\alpha=1}^N\fU_\alpha$ be a finite affine open
covering of~$\fW$.
 Put $\fV_\alpha=\ff^{-1}(\fU_\alpha)\subset\fT\subset\fY$; then
$\fT=\bigcup_{\alpha=1}^N\fV_\alpha$ is a finite affine open
covering of~$\fT$.
 The sheaf axiom~\eqref{topology-base-sheaf-axiom} provides a left
exact sequence of torsion $\fR$\+modules
\begin{equation} \label{flat-base-change-sheaf-axiom-downstairs}
 0\lrarrow\cL(\fW)\lrarrow\bigoplus\nolimits_{\alpha=1}^N\cL(\fU_\alpha)
 \lrarrow\bigoplus\nolimits_{1\le\alpha<\beta\le N}
 \cL(\fU_\alpha\cap\fU_\beta)
\end{equation}
as well as a left exact sequence of torsion $\fS$\+modules
\begin{equation} \label{flat-base-change-sheaf-axiom-upstairs}
 0\lrarrow(\ff'{}^*\cL)(\fT)\lrarrow
 \bigoplus\nolimits_{\alpha=1}^N(\ff'{}^*\cL)(\fV_\alpha)
 \lrarrow\bigoplus\nolimits_{1\le\alpha<\beta\le N}
 (\ff'{}^*\cL)(\fV_\alpha\cap\fV_\beta).
\end{equation}
 Since $\fS$ is a flat contramodule $\fR$\+module, the tensor product
functor ${-}\ot_\fR\fS$ takes the left exact
sequence~\eqref{flat-base-change-sheaf-axiom-downstairs} to a left
exact sequence of $\fS$\+modules
\begin{multline} \label{flat-base-change-sheaf-axiom-tensored}
 0\lrarrow\cL(\fW)\ot_\fR\fS\lrarrow
 \bigoplus\nolimits_{\alpha=1}^N\bigl(\cL(\fU_\alpha)\ot_\fR\fS\bigr)
 \\ \lrarrow\bigoplus\nolimits_{1\le\alpha<\beta\le N}
 \bigl(\cL(\fU_\alpha\cap\fU_\beta)\ot_\fR\fS\bigr).
\end{multline}
 The natural $\fR$\+linear map
from~\eqref{flat-base-change-sheaf-axiom-downstairs}
to~\eqref{flat-base-change-sheaf-axiom-upstairs} induces
a natural $\fS$\+linear map
from~\eqref{flat-base-change-sheaf-axiom-tensored}
to~\eqref{flat-base-change-sheaf-axiom-upstairs}.
 The latter map is a morphism of left exact sequences of $\fS$\+modules
that is an isomorphism on the middle and rightmost terms according to
(the proof of) part~(a).
 It follows that the natural map
$$
 (\ff^*\fj_*\cL)(\fY)=(\fj_*\cL)(\fX)\ot_\fR\fS=\cL(\fW)\ot_\fR\fS
 \lrarrow(\ff'{}^*\cL)(\fT)=(\fj'_*\ff'{}^*\cL)(\fY)
$$
is an isomorphism, which is the desired assertion of part~(b).
\end{proof}

\begin{lem} \label{cosheaf-arb-dir-openimmers-inv-basechange-lemma}
 Let\/ $\ff\:(\fY,\fO_\fY)\rarrow(\fX,\fO_\fX)$ be a morphism of
ringed spaces.
 Let\/ $\fW\subset\fX$ be an open subset, viewed as a ringed space
$(\fW,\fO_\fW)$ with\/ $\fO_\fW=\fO_\fX|_\fW$.
 Consider the open subset\/ $\fT=\ff^{-1}(\fW)\subset\fY$, and let
us view it as a ringed space $(\fT,\fO_\fT)$ with
$\fO_\fT=\fO_\fY|_\fT$.
 Denote the induced morphism of ringed spaces by\/
$\ff'\:(\fT,\fO_\fT)\rarrow(\fW,\fO_\fW)$.
 Let\/ $\fQ$ be a cosheaf of\/ $\fO_\fY$\+modules on\/~$\fY$.
 Then there is a natural isomorphism
$$
 (\ff_!\fQ)|_\fW\simeq\ff'_!(\fQ|_\fT)
$$
of cosheaves of\/ $\fO_\fW$\+modules on\/~$\fW$.
\end{lem}

\begin{proof}
 Follows immediately from the constructions (see
Section~\ref{direct-images-of-ctrh-contramod-subsecn}).
\end{proof}

\subsection{\v Cech (co)resolutions of sheaves and cosheaves}
 Let $\fX$ be a Noetherian formal scheme with a finite affine open
covering $\fX=\bigcup_{\alpha=1}^N\fU_\alpha$.
 Denote by $\fj_{\alpha_1,\dotsc,\alpha_i}$ the open immersion morphisms
$\fU_{\alpha_1}\cap\dotsb\cap\fU_{\alpha_i}\rarrow\fX$.

 Let $\cM$ be a quasi-coherent torsion sheaf on~$\fX$.
 Then there is a natural \v Cech exact sequence of quasi-coherent
torsion sheaves on~$\fX$
\begin{multline} \label{Cech-sequence-of-torsion-sheaves}
 0\lrarrow\cN\lrarrow\bigoplus\nolimits_{\alpha=1}^N
 \fj_\alpha{}_*\fj_\alpha^*\cM\lrarrow
 \bigoplus\nolimits_{1\le\alpha<\beta\le N}
 \fj_{\alpha,\beta}{}_*\fj_{\alpha,\beta}^*\cM \\
 \lrarrow\dotsb\lrarrow\fj_{1,\dotsc,N}{}_*\fj_{1,\dotsc,N}^*\cM
 \lrarrow0
\end{multline}
with the differentials defined in terms of
the adjunction~\eqref{sheaves-restriction-direct-image-adjunction}
or~\eqref{qcoh-tors-inverse-direct-adjunction}.
 To check exactness, it suffices to notice that the restrictions
of the sequence of quasi-coherent torsion
sheaves~\eqref{Cech-sequence-of-torsion-sheaves} to the open formal
subschemes $\fU_\alpha\subset\fX$ are contractible.
 In the case when $\fX$ is semi-separated, one can also say that
the sequence of sections of~\eqref{Cech-sequence-of-torsion-sheaves}
over any affine open formal subscheme $\fW\subset\fX$ is exact by
Lemma~\ref{torsion-cohomol-Cech-sequence-lemma}.

 Now assume that $\fX$ is semi-separated, and let $\bW$ be an open
covering of $\fX$ such that the open covering $\{\fU_\alpha\}$ is
subordinate to~$\bW$.
 Let $\fP$ be a $\bW$\+locally contraherent cosheaf of contramodules
on~$\fX$.
 Then there is a natural \v Cech exact sequence of $\bW$\+locally
contraherent cosheaves of contramodules on~$\fX$
\begin{multline} \label{Cech-sequence-of-cosheaves-of-contramods}
 0\lrarrow\fj_{1,\dotsc,N}{}_!\fj_{1,\dotsc,N}^!\fP
 \lrarrow\dotsb \\
 \lrarrow\bigoplus\nolimits_{1\le\alpha<\beta\le N}
 \fj_{\alpha,\beta}{}_!\fj_{\alpha,\beta}^!\fP\lrarrow
 \bigoplus\nolimits_{\alpha=1}^N\fj_\alpha{}_!\fj_\alpha^!\fP
 \lrarrow\fP\lrarrow0
\end{multline}
with the differentials defined in terms of
the adjunction~\eqref{cosheaves-direct-image-restriction-adjunction}
or~\eqref{lcth-contramod-direct-inverse-adjunction}.
 To show that the finite
sequence~\eqref{Cech-sequence-of-cosheaves-of-contramods} is exact
in the exact category $\fX\Lcth_\bW$, it suffices to point out that
the sequence of cosections
of~\eqref{Cech-sequence-of-cosheaves-of-contramods}
over any affine open formal subscheme $\fW\subset\fX$ subordinate to
$\bW$ is exact by Lemma~\ref{contraadj-homol-Cech-sequence-lemma}.

 Notice that \eqref{Cech-sequence-of-cosheaves-of-contramods}~is
a finite resolution of a $\bW$\+locally contraherent cosheaf of
contramodules $\fP$ on $\fX$ by (globally) contraherent cosheaves
of contramodules on~$\fX$.

 If $\fP$ is a locally cotorsion $\bW$\+locally contraherent cosheaf
of contramodules on $\fX$, then the finite
sequence~\eqref{Cech-sequence-of-cosheaves-of-contramods} is exact in
the exact category $\fX\Lcth_\bW^\lct$ by
Lemma~\ref{complex-of-cosheaves-exactness-criterion} or by the proof
of Corollary~\ref{colocality-of-cotorsion-presuming-contraadjusted}.
 So \eqref{Cech-sequence-of-cosheaves-of-contramods}~is
a finite resolution of a locally cotorsio $\bW$\+locally contraherent
cosheaf of contramodules $\fP$ on $\fX$ by locally cotorsion (globally)
contraherent cosheaves of contramodules on $\fX$ in this case.

\subsection{Contraherent $\fHom$}  \label{contraherent-Hom-subsecn}
 This section is a formal scheme counterpart
of~\cite[Section~2.5]{Pcosh} and~\cite[Sections~5.5\+-5.6]{Pdomc}.

 Let $\fX$ be a locally Noetherian formal scheme.
 Given an associative ring $R$ and a quasi-coherent torsion sheaf
$\cM$ on $\fX$, a \emph{right action of $R$ on\/~$\cM$} can be simply
defined as a ring homomorphism $R^\rop\rarrow\Hom_\fX(\cM,\cM)$
(where $R^\rop$ denotes the opposite ring to~$R$).
 This is a special case of the general concept of a ring acting on
an object of an additive category.
 Applying any covariant additive functor to an object with an action of
the ring $R$ produces an object endowed with a natural action of $R$
again.
 In particular, the direct and inverse images of quasi-coherent
torsion sheaves with a right action of $R$ also carry natural right
actions of~$R$.
 Applying a contravariant additive functor to an object with a right
action of $R$ produces an object with a left action of~$R$.

 Let $\cM$ be a quasi-coherent torsion sheaf on $\fX$ with a right
action of $R$ on $\cM$, and let $L$ be a left $R$\+module.
 Then the quasi-coherent torsion sheaf $\cM\ot_RL$ on $\fX$ is defined
by the rule
$$
 (\cM\ot_RL)(\fU)=\cM(\fU)\ot_RL
$$
for all affine open formal subschemes $\fU\subset\fX$.
 Clearly, the torsion axiom~(i) and the quasi-coherence axiom~(ii)
from Section~\ref{qcoh-torsion-sheaves-subsecn} are satisfied for
the presheaf of $\fO_\fX$\+modules $\cM\ot_RL$ on the topology base
of all affine open formal subschemes $\fU\subset\fX$ defined by
the rule above; so the quasi-coherent torsion sheaf $\cM\ot_RL$
on $\fX$ is well-defined.
 This sheaf-theoretic construction is a special case of a general
category-theoretic construction of the tensor product of an object with
a module, which is applicable to objects with a ring action in any
additive category with colimits.

 For an additive category $\sB$ with limits, a dual-analogous
category-theoretic construction assigns an object $\Hom_R(L,B)\in\sB$
to every object $B\in\sB$ with a left action of $R$ and every left
$R$\+module~$L$.
 Both the mentioned category-theoretic tensor product and
category-theoretic Hom object are easily characterized by universal
properties and constructed in terms of any presentation of
the $R$\+module $L$ as the cokernel of a morphism of free $R$\+modules.

\begin{lem} \label{torsion-sheaf-tensored-with-module-adjunction}
 Let\/ $\cM$ be a quasi-coherent torsion sheaf on\/ $\fX$ with a right
action of $R$ on\/ $\cM$, let $L$ be a left $R$\+module, and let\/ $\cN$
be a quasi-coherent torsion sheaf on\/~$\fX$.
 Then there is a natural adjunction isomorphism of abelian groups
$$
 \Hom_\fX(\cM\ot_RL,\>\cN)\simeq\Hom_R(L,\Hom_\fX(\cM,\cN)).
$$
\end{lem}

\begin{proof}
 First of all, the contravariant functor $\Hom_\fX({-},\cN)\:\fX\Tors
\rarrow\Ab$ takes an object $\cM$ with a right action of $R$ to
an abelian group with a left action of~$R$.
 So $\Hom_\fX(\cM,\cN)$ is a left $R$\+module.
 Furthermore, for any contravariant functor $H\:\sA\rarrow\sB$ from
an additive category $\sA$ with colimits to an additive category $\sB$
with limits such that $H$ transforms colimits into limits, and for
any object $M\in\sA$ with a right action of $R$ and any left $R$\+module 
$L$, one has a natural isomorphism $H(M\ot_RL)\simeq\Hom_R(L,H(M))$
in~$\sB$.
\end{proof}

 Now let $R$ be an adic topological ring and $\cM$ be a quasi-coherent
torsion sheaf on $\fX$ with an action of the ring~$R$.
 We will say that the sheaf $\cM$ is \emph{$R$\+torsion} if, for every
affine open formal subscheme $\fU\subset\fX$, the $R$\+module
$\cM(\fU)$ is torsion.
 It suffices to check this condition for affine open formal subschemes
$\fU$ belonging to any given affine open covering of the formal
scheme~$\fX$.
 It is clear from the constructions of the direct and inverse images
of quasi-coherent torsion sheaves on locally Noetherian formal schemes
in Sections~\ref{direct-images-of-qcoh-tors-subsecn}
and~\ref{inverse-images-of-qcoh-tors-subsecn} that the functors of
direct image under quasi-compact morphisms and the functors or inverse
image under tight morphisms take $R$\+torsion quasi-coherent torsion
sheaves to $R$\+torsion quasi-coherent torsion sheaves.

\begin{lem} \label{sheaf-tensored-with-module-sections-over-nonaffine}
\textup{(a)} Let\/ $\cM$ be a quasi-coherent torsion sheaf on\/ $\fX$
with a right action of a ring $R$ on\/ $\cM$, and let $L$ be
a left $R$\+module.
 Then, for any open formal subscheme\/ $\fY\subset\fX$, there is
a natural map
\begin{equation} \label{sheaf-module-tensor-sections-comparison-map}
 \cM(\fY)\ot_RL\lrarrow(\cM\ot_RL)(\fY).
\end{equation} \par
\textup{(b)} Let $R$ be an adic topological ring and\/ $\cM$ be
an $R$\+torsion quasi-coherent torsion sheaf on\/~$\fX$.
 Let\/ $\fF$ be a flat contramodule $R$\+module.
 Then, for any quasi-compact open formal subscheme\/ $\fY\subset\fX$,
the map\/ $\cM(\fY)\ot_R\fF\lrarrow(\cM\ot_R\fF)(\fY)$
\,\eqref{sheaf-module-tensor-sections-comparison-map} is an isomorphism.
\end{lem}

\begin{proof}
 Part~(a): quite generally, let $\sA$ and $\sB$ be additive categories
with colimits, and let $G\:\sA\rarrow\sB$ be an additive functor.
 Let $M\in\sA$ be an object with a right action of a ring $R$ and
$L$ be a left $R$\+module.
 Then there is a natural morphism $G(M)\ot_RL\rarrow G(M\ot_RL)$
in~$\sB$.

 Part~(b): let $\fY=\bigcup_\alpha\fU_\alpha$ be a finite affine open
covering of $\fY$, and let $\fU_\alpha\cap\fU_\beta=\bigcup_\gamma
\fV_{\alpha\beta\gamma}$ be finite affine open coverings of
the intersections.
 Then, according to the sheaf axiom~\eqref{topology-base-sheaf-axiom},
we have left exact sequences of $\fO_\fX(\fY)$\+modules
\begin{equation} \label{M-of-Y-sheaf-axiom}
 0\lrarrow\cM(\fY)\lrarrow\bigoplus\nolimits_\alpha\cM(\fU_\alpha)
 \lrarrow\bigoplus\nolimits_{\alpha,\beta,\gamma}
 \cM(\fV_{\alpha\beta\gamma})
\end{equation}
and
\begin{equation} \label{M-tensor-L-of-Y-sheaf-axiom}
 0\lrarrow(\cM\ot_RL)(\fY)\lrarrow\bigoplus\nolimits_\alpha
 (\cM\ot_RL)(\fU_\alpha)\lrarrow\bigoplus
 \nolimits_{\alpha,\beta,\gamma}(\cM\ot_RL)(\fV_{\alpha\beta\gamma}).
\end{equation}

 Now~\eqref{M-of-Y-sheaf-axiom} is actually a short exact sequence
of $\fO_\fX(\fY)$\+$R$\+bimodules, and all the three terms
of~\eqref{M-of-Y-sheaf-axiom} are torsion $R$\+modules.
 Since $\fF$ is flat contramodule $R$\+module, the functor ${-}\ot_R\fF$
takes~\eqref{M-of-Y-sheaf-axiom} to a left exact sequence of
$\fO_\fX(\fY)$\+modules
\begin{equation} \label{M-of-Y-sheaf-axiom-tensor-with-F}
 0\lrarrow\cM(\fY)\ot_R\fF\lrarrow\bigoplus\nolimits_\alpha
 \cM(\fU_\alpha)\ot_R\fF
 \lrarrow\bigoplus\nolimits_{\alpha,\beta,\gamma}
 \cM(\fV_{\alpha\beta\gamma})\ot_R\fF.
\end{equation}
 Putting $L=\fF$, we have a natural
morphism~\eqref{sheaf-module-tensor-sections-comparison-map}
from the left exact sequence~\eqref{M-of-Y-sheaf-axiom-tensor-with-F}
to the left exact sequence~\eqref{M-tensor-L-of-Y-sheaf-axiom}.
 This morphism of left exact sequences is an isomorphism on the middle
and rightmost terms (by the definition of the sheaf $\cM\ot_RL$), and
consequently on the leftmost terms as well.
\end{proof}

 Let $\ff\:\fY\rarrow\fX$ be a quasi-compact morphism of locally
Noetherian formal schemes, let $\cN$ be a quasi-coherent torsion
sheaf on $\fY$ with a right action of a ring $R$, and let $L$ be
a left $R$\+module.
 Then, according to the argument in the proof of
Lemma~\ref{sheaf-tensored-with-module-sections-over-nonaffine}(a),
there is a natural morphism of quasi-coherent torsion sheaves
\begin{equation} \label{sheaf-module-tensor-direct-image}
 \ff_*(\cN)\ot_RL\lrarrow\ff_*(\cN\ot_RL)
\end{equation}
on~$\fX$.
 Now let us assume that $R$ is an adic topological ring, and let
$\fF$ be a flat contramodule $R$\+module.
 Then it is clear from
Lemma~\ref{sheaf-tensored-with-module-sections-over-nonaffine}(b)
that the morphism $\ff_*(\cN)\ot_R\fF\rarrow\ff_*(\cN\ot_R\fF)$
\,\eqref{sheaf-module-tensor-direct-image} is an isomorphism
in $\fX\Tors$.

 When the morphism~$\ff$ is affine, the functor~$\ff_*$ is exact
and preserves infinite direct sums, so it preserves all colimits.
 It follows that the morphism~\eqref{sheaf-module-tensor-direct-image}
is an isomorphism for any ring $R$ and any $R$\+module $L$ in
this case.

 Let $\cM$ be a quasi-coherent torsion sheaf on $\fX$, and let $\cJ$
be an injective quasi-coherent torsion sheaf on~$\fX$.
 In this setting, the locally cotorsion contraherent cosheaf of
contramodules $\fHom_\fX(\cM,\cJ)$ is constructed as follows.

 Let $\fU$ be an affine open formal subscheme in~$\fX$.
 Denote by $\fj\:\fU\rarrow\fX$ the open immersion morphism.
 Then the quasi-coherent torsion sheaf $\fj^*\cM=\cM|_\fU$ on $\fU$
comes endowed with a natural action of the ring $R=\fO_\fX(\fU)$.
 Hence the ring $\fO_\fX(\fU)$ also acts on the sheaf $\fj_*\fj^*\cM$.
 Notice that $\fj$~is a quasi-compact morphism of locally Noetherian
formal schemes; so $\fj_*\fj^*\cM$ is a quasi-coherent torsion sheaf
on $\fX$ by Proposition~\ref{torsion-sheaves-qcoh-direct-images-prop}.

 Now we put $\fHom_\fX(\cM,\cJ)[\fU]=\Hom_\fX(\fj_*\fj^*\cM,\cJ)$.
 Following the discussion above, the abelian group
$\Hom_\fX(\fj_*\fj^*\cM,\cJ)$ has a natural $\fO_\fX(\fU)$\+module
structure.

\begin{lem} \label{fHom-cosection-cotorsion-contramodule}
 The\/ $\fO_\fX(\fU)$\+module\/ $\Hom_\fX(\fj_*\fj^*\cM,\cJ)$ is
a cotorsion (quotseparated) contramodule module over the adic
topological ring\/ $\fO_\fX(\fU)$.
\end{lem}

\begin{proof}
 Put $\fR=\fO_\fX(\fU)$.
 To prove that $\Hom_\fX(\fj_*\fj^*\cM,\cJ)$ is a contramodule
$\fR$\+module, pick an ideal of definition $\fI\subset\fR$.
 Then we have $\cM(\fU)=\bigcup_{n\ge1}{}_{\fI^n}\cM(\fU)$,
where ${}_{\fI^n}\cM(\fU)$ is the submodule of all elements
annihilated by $\fI^n$ in $\cM(\fU)$.
 Denoting by $(\fj^*\cM)_n\subset\fj^*\cM$ the quasi-coherent
torsion subsheaf corresponding to the submodule ${}_{\fI^n}\cM(\fU)
\subset\cM(\fU)$, we have $\fj^*\cM=\varinjlim_{n\ge1}(\fj^*\cM)_n$,
hence $\fj_*\fj^*\cM=\varinjlim_{n\ge1}\fj_*((\fj^*\cM)_n)$.
 Thus
$$
 \Hom_\fX(\fj_*\fj^*\cM,\cJ)=\varprojlim\nolimits_{n\ge1}
 \Hom_\fX(\fj_*((\fj^*\cM)_n),\cJ).
$$
 It remains to point out that $(\fj^*\cM)_n$ is a quasi-coherent
torsion sheaf with an action of the ring $\fR/\fI^n$, hence
$\Hom_\fX(\fj_*((\fj^*\cM)_n),\cJ)$ is an $\fR$\+module annihilated
by~$\fI^n$.
 Any such module is a (separated) contramodule $\fR$\+module, and
any projective limit of contramodule $\fR$\+modules in $\fR\Modl$ is
a contramodule $\fR$\+module.
 In fact, this argument proves that $\Hom_\fX(\fj_*\fj^*\cM,\cJ)$
is a separated contramodule $\fO_\fX(\fU)$\+module.

 To prove that $\fC=\Hom_\fX(\fj_*\fj^*\cM,\cJ)$ is a cotorsion
contramodule $\fR$\+module, we use the criterion of
Lemma~\ref{cotorsion-contramodules-lemma}(1).
 Let $0\rarrow\fH\rarrow\fG\rarrow\fF\rarrow0$ be a short exact
sequence of flat contramodule $\fR$\+modules.
 We need to check exactness of the short sequence of abelian groups
$0\rarrow\Hom_\fR(\fF,\fC)\rarrow\Hom_\fR(\fG,\fC)\rarrow
\Hom_\fR(\fH,\fC)\rarrow0$.
 For any $\fR$\+module $L$, we have
$$
 \Hom_\fR(L,\fC)=\Hom_\fR(L,\Hom_\fX(\fj_*\fj^*\cM,\cJ))\simeq
 \Hom_\fX((\fj_*\fj^*\cM)\ot_\fR L,\>\cJ)
$$
by Lemma~\ref{torsion-sheaf-tensored-with-module-adjunction}.

 Now $\fj^*\cM$ is an $\fR$\+torsion quasi-coherent torsion sheaf
on $\fU$, hence $\fj_*\fj^*\cM$ is an $\fR$\+torsion quasi-coherent
torsion sheaf on~$\fX$ (according to the discussion above in
this section).
 By Lemma~\ref{flat-quotseparated-contramodules-well-behaved}(b),
it follows that $0\rarrow(\fj_*\fj^*\cM)\ot_\fR\fH\rarrow
(\fj_*\fj^*\cM)\ot_\fR\fG\rarrow(\fj_*\fj^*\cM)\ot_\fR\fF\rarrow0$
is a short exact sequence of quasi-coherent torsion sheaves on~$\fX$.
 Since $\cJ$ is an injective object of $\fX\Tors$, the desired
exactness of the short sequence of abelian groups
$0\rarrow\Hom_\fX((\fj_*\fj^*\cM)\ot_\fR\fF,\>\cJ)\rarrow
\Hom_\fX((\fj_*\fj^*\cM)\ot_\fR\fG,\>\cJ)\rarrow
\Hom_\fX((\fj_*\fj^*\cM)\ot_\fR\fH,\>\cJ)\rarrow0$ follows.
\end{proof}

\begin{lem} \label{fHom-contraherence}
 Let\/ $\fV\subset\fU\subset\fX$ be a pair of affine open formal
subschemes in\/ $\fX$ with the open immersion morphisms\/
$\fj\:\fU\rarrow\fX$ and\/ $\fk\:\fV\rarrow\fX$.
 Then there is a natural isomorphism of\/ $\fO_\fX(\fV)$\+modules
$$
 \Hom_\fX(\fk_*\fk^*\cM,\cJ)\simeq
 \Hom_{\fO_\fX(\fU)}\bigl(\fO_\fX(\fV),\Hom_\fX(\fj_*\fj^*\cM,\cJ)\bigr).
$$
\end{lem}

\begin{proof}
 Denote by\/ $\fh\:\fV\rarrow\fU$ the natural open immersion morphism;
so $\fk=\fj\fh$.
 Then we have
$$
 \fk_*\fk^*\cM\simeq\fj_*\fh_*\fh^*\fj^*\cM\simeq
 \fj_*\bigl((\fj^*\cM)\ot_{\fO_\fX(\fU)}\fO_\fX(\fV)\bigr)
$$
by the quasi-coherence axiom~(ii) from
Section~\ref{qcoh-torsion-sheaves-subsecn}.
 According to the discussion of the natural
morphism~\eqref{sheaf-module-tensor-direct-image} applied to the flat
contramodule $\fO_\fX(\fU)$\+module $L=\fF=\fO_\fX(\fV)$, we have
$$
 \fj_*\bigl((\fj^*\cM)\ot_{\fO_\fX(\fU)}\fO_\fX(\fV)\bigr)\simeq
 (\fj_*\fj^*\cM)\ot_{\fO_\fX(\fU)}\fO_\fX(\fV).
$$
 Now it remains to refer to
Lemma~\ref{torsion-sheaf-tensored-with-module-adjunction} for
the isomorphism
$$
 \Hom_\fX\bigl((\fj_*\fj^*\cM)\ot_{\fO_\fX(\fU)}\fO_\fX(\fV),\>\cJ\bigr)
 \simeq\Hom_{\fO_\fX(\fU)}
 \bigl(\fO_\fX(\fV),\Hom_\fX(\fj_*\fj^*\cM,\cJ)\bigr).
$$
\end{proof}

 Lemmas~\ref{fHom-cosection-cotorsion-contramodule}
and~\ref{fHom-contraherence} say that $\fHom_\fX(\cM,\cJ)$ is
a locally cotorsion contraherent copresheaf of contramodules on
the topology base of affine open formal subschemes $\fU\subset\fX$.
 So the locally cotorsion (globally) contraherent cosheaf of
contramodules $\fHom_\fX(\cM,\cJ)$ on $\fX$ is well-defined.

\begin{prop} \label{cosections-of-fHom-prop}
 Let\/ $\fX$ be a locally Noetherian formal scheme, $\cM$ be
a quasi-coherent torsion sheaf on\/ $\fX$, and\/ $\cJ$ be an injective
quasi-coherent torsion sheaf on\/~$\fX$.
 Let\/ $\fY$ be a quasi-compact open formal subscheme in\/ $\fX$ with
the open immersion morphism\/ $\ff\:\fY\rarrow\fX$.
 Then there is a natural isomorphism of\/ $\fO_\fX(\fY)$\+modules
$$
 \fHom_\fX(\cM,\cJ)[\fY]\simeq
 \Hom_\fX(\ff_*\ff^*\cM,\cJ).
$$
\end{prop}

\begin{proof}
 This is our version of~\cite[Lemma~2.5.2(c)]{Pcosh}
and~\cite[Lemma~5.18(d)]{Pdomc}.
 Let $\fY=\bigcup_{\alpha=1}^N\fV_\alpha$ be a finite affine open
covering of a Noetherian formal scheme $\fY$ and
$\fV_\alpha\cap\fV_\beta=\bigcup_\gamma\fW_{\alpha\beta\gamma}$
be finite affine open coverings of the intersections
$\fV_\alpha\cap\fV_\beta$.
 Denote by $\fj_{\alpha_1,\dotsc,\alpha_i}$ the open immersion morphisms
$\fV_{\alpha_1}\cap\dotsb\cap\fV_{\alpha_i}\rarrow\fY$.
 Then, for any quasi-coherent torsion sheaf $\cN$ on $\fY$, there is
a finite \v Cech exact
sequence~\eqref{Cech-sequence-of-torsion-sheaves}
\begin{multline} \label{Cech-sequence-of-torsion-sheaves-on-Y}
 0\lrarrow\cN\lrarrow\bigoplus\nolimits_{\alpha=1}^N
 \fj_\alpha{}_*\fj_\alpha^*\cN\lrarrow
 \bigoplus\nolimits_{1\le\alpha<\beta\le N}
 \fj_{\alpha,\beta}{}_*\fj_{\alpha,\beta}^*\cN \\
 \lrarrow\dotsb\lrarrow\fj_{1,\dotsc,N}{}_*\fj_{1,\dotsc,N}^*\cN
 \lrarrow0
\end{multline}
of quasi-coherent torsion sheaves on~$\fY$.

 Denoting by $\fh_{\alpha\beta\gamma}$ the open immersion morphisms
$\fW_{\alpha\beta\gamma}\rarrow\fV_\alpha\cap\fV_\beta$ and
applying~\eqref{Cech-sequence-of-torsion-sheaves-on-Y} to
the quasi-coherent torsion sheaves $\fj_{\alpha,\beta}^*\cN$ on
the Noetherian formal schemes $\fV_\alpha\cap\fV_\beta$ with their
open coverings by $\fW_{\alpha\beta\gamma}$, we see that the morphism
\begin{equation} \label{intersection-covered-Cech-monomorphism}
 \fj_{\alpha,\beta}^*\cN\lrarrow\bigoplus\nolimits_\gamma
 \fh_{\alpha\beta\gamma}{}_*\fh_{\alpha\beta\gamma}^*
 \fj_{\alpha,\beta}^*\cN
\end{equation}
of quasi-coherent torsion sheaves on $\fV_\alpha\cap\fV_\beta$ is
injective.
 Denote by $\fk_{\alpha\beta\gamma}=\fj_{\alpha,\beta}
\fh_{\alpha\beta\gamma}$ the open immersion morphisms
$\fW_{\alpha\beta\gamma}\rarrow\fY$.
 Applying the left exact functor of direct image
$\fj_{\alpha,\beta}{}_*$
to~\eqref{intersection-covered-Cech-monomorphism}, we obtain
a monomorphism
\begin{equation} \label{intersection-covered-Cech-monomorphism2}
 \fj_{\alpha,\beta}{}_*\fj_{\alpha,\beta}^*\cN\lrarrow
 \bigoplus\nolimits_{\alpha\beta\gamma}
 \fk_{\alpha\beta\gamma}{}_*\fk_{\alpha\beta\gamma}^*\cN
\end{equation}
of quasi-coherent torsion sheaves on~$\fY$.
 Comparing~\eqref{Cech-sequence-of-torsion-sheaves-on-Y}
with~\eqref{intersection-covered-Cech-monomorphism2}, we arrive at
a left exact sequence
\begin{equation} \label{Cech-left-exact-sequence-of-torsion-sheaves}
 0\lrarrow\cN\lrarrow\bigoplus\nolimits_{\alpha=1}^N
 \fj_\alpha{}_*\fj_\alpha^*\cN\lrarrow
 \bigoplus\nolimits_{\alpha,\beta,\gamma}
 \fk_{\alpha\beta\gamma}{}_*\fk_{\alpha\beta\gamma}^*\cN
\end{equation}
of quasi-coherent torsion sheaves on~$\fY$.
 Finally, setting $\cN=\ff^*\cN=\cN|_\fY$ and applying the left exact
functor $\ff_*\:\fY\Tors\rarrow\fX\Tors$
to~\eqref{Cech-left-exact-sequence-of-torsion-sheaves},
we obtain a left exact sequence
\begin{equation} \label{Cech-left-exact-sequence-of-torsion-sheaves2}
 0\lrarrow\ff_*\ff^*\cM\lrarrow\bigoplus\nolimits_{\alpha=1}^N
 \ff_*\fj_\alpha{}_*\fj_\alpha^*\ff^*\cM\lrarrow
 \bigoplus\nolimits_{\alpha,\beta,\gamma}
 \ff_*\fk_{\alpha\beta\gamma}{}_*\fk_{\alpha\beta\gamma}^*\ff^*\cM
\end{equation}
of quasi-coherent torsion sheaves on~$\fX$.

 It remains to apply the exact functor $\Hom_\fX({-},\cJ)$ to
the left exact
sequence~\eqref{Cech-left-exact-sequence-of-torsion-sheaves2},
producing a right exact sequence of $\fO_\fX(\fY)$\+modules
\begin{multline} \label{cosections-of-fHom-computed}
 \bigoplus\nolimits_{\alpha,\beta,\gamma}\Hom_\fX
 (\ff_*\fk_{\alpha\beta\gamma}{}_*\fk_{\alpha\beta\gamma}^*\ff^*\cM,
 \cJ) \\ \lrarrow
 \bigoplus\nolimits_{\alpha=1}^N
 \Hom_\fX(\ff_*\fj_\alpha{}_*\fj_\alpha^*\ff^*\cM,\cJ)
 \lrarrow\Hom_\fX(\ff_*\ff^*\cM,\cJ)\rarrow0.
\end{multline}
 By construction, the leftmost morphism
in~\eqref{cosections-of-fHom-computed} is naturally isomorphic to
the leftmost morphism in the right exact sequence
\begin{multline} \label{fHom-cosheaf-axiom}
 \bigoplus\nolimits_{\alpha,\beta,\gamma}
 \fHom_\fX(\cM,\cJ)[\fW_{\alpha\beta\gamma}] \\
 \lrarrow\bigoplus\nolimits_{\alpha=1}^N
 \fHom_\fX(\cM,\cJ)[\fV_\alpha]\lrarrow\fHom_\fX(\cM,\cJ)[\fY]
 \lrarrow0
\end{multline}
from the cosheaf axiom~\eqref{topology-base-cosheaf-axiom} for
the contraherent cosheaf of contramodules $\fHom_\fX(\cM,\cJ)$
on $\fX$, the open subset $\fY\subset\fX$, and the open covering
$\fY=\bigcup_\alpha\fV_\alpha$.
 It follows that the cokernels of the two leftmost morphisms, i.~e.,
the rightmost terms of~\eqref{cosections-of-fHom-computed}
and~\eqref{fHom-cosheaf-axiom}, are also naturally isomorphic,
as desired. \hbadness=1550
\end{proof}

\subsection{Contratensor product}  \label{contratensor-product-subsecn}
 This section is a formal scheme version of~\cite[Section~2.6]{Pcosh}
and~\cite[Section~5.7]{Pdomc}.
 Let $\fX$ be a locally Noetherian formal scheme.
 By a \emph{diagram} (\emph{of affine open formal subschemes}) in
$\fX$ we mean a poset $\bD$ endowed with an order-preserving map
$a\longmapsto\fU_a$ from $\bD$ into the poset of all affine open formal
subschemes of $\fX$ (with respect to the inclusion).
 The order-preservation conditition means that for every pair of
elements $b\le a$ in $\bD$ one has $\fU_b\subset\fU_a$.

 A diagram $\bD$ in $\fX$ is said to be \emph{exhaustive} if
the following two conditions hold:
\begin{enumerate}
\renewcommand{\theenumi}{\roman{enumi}}
\item $\fX=\bigcup_{a\in\bD}\fU_a$;
\item for every pair of elements $a$, $b\in\bD$ one has
$\fU_a\cap\fU_b=\bigcup_{c\in\bD}^{c\le a,\,c\le b}\fU_c$.
 Here the notation means that the union is taken over all
the indices $c\in\bD$ such that $c\le a$ and $c\le b$.
\end{enumerate}
 For example, given an open covering $\bW$ of $\fX$, the poset of
all affine open formal subschemes $\fU\subset\fX$ subordinate to $\bW$,
ordered by inclusion, is an exhaustive diagram of affine open formal
subschemes in~$\fX$.
 More generally, any topology base $\bB$ consisting of affine open
formal subschemes in $\fX$, ordered by inclusion, is an exhaustive
diagram.

 Let $\cM$ be a quasi-coherent torsion sheaf and $\fP$ be
a cosheaf of $\fO_\fX$\+modules on~$\fX$.
 Let $\fV\subset\fU\subset\fX$ be a pair of affine open formal
subschemes with open immersion morphisms $\fj\:\fU\rarrow\fX$,
\ $\fk\:\fV\rarrow\fX$, and $\fh\:\fV\rarrow\fU$.
 Following the proof of Lemma~\ref{fHom-contraherence},
we have natural isomorphisms
$$
 \fk_*\fk^*\cM\simeq\fj_*\fh_*\fh^*\fj^*\cM\simeq
 (\fj_*\fj^*\cM)\ot_{\fO_\fX(\fU)}\fO_\fX(\fV)
$$
of quasi-coherent torsion sheaves on~$\fX$.
 Consequently, the corestriction map $\fP[\fV]\rarrow\fP[\fU]$ of
the cosheaf $\fP$ induces a natural morphism
\begin{equation} \label{contratensor-transition-map}
 (\fk_*\fk^*\cM)\ot_{\fO_\fX(\fV)}\fP[\fV]\lrarrow
 (\fj_*\fj^*\cM)\ot_{\fO_\fX(\fU)}\fP[\fU]
\end{equation}
of quasi-coherent torsion sheaves on~$\fX$.

 The \emph{contratensor product} of the quasi-coherent torsion sheaf
$\cM$ and the cosheaf of $\fO_\fX$\+modules $\fP$, computed on
a diagram of affine formal open subschemes $\bD$ in $\fX$, is defined
as the quasi-coherent torsion sheaf on $\fX$ constructed as
the inductive limit
$$
 \cM\ocn_{\fX,\bD}\fP = \varinjlim\nolimits_{a\in\bD}
 \bigl((\fj_a{}_*\fj_a^*\cM)\ot_{\fO_\fX(\fU_a)}\fP[\fU_a]\bigr).
$$
 Here the inductive limit indexed by the poset $\bD$ is computed in
the Grothendieck abelian category $\fX\Tors$ with respect to
the transition maps~\eqref{contratensor-transition-map}.
 The notation~$\fj_a$ stands for the open immersion morphism
$\fU_a\rarrow\fX$.
 Notice that the relevant diagrams $\bD$ are usually \emph{not} directed
(i.~e., the poset $\bD$ is not directed), so the inductive limit functor
involved is \emph{not} exact.

\begin{lem} \label{fHom-contratensor-adjunction-lemma}
 Let\/ $\fX$ be a locally Noetherian formal scheme, $\cM$ be
a quasi-coherent torsion sheaf, $\cJ$ be an injective quasi-coherent
torsion sheaf, and\/ $\fP$ be a cosheaf of\/ $\fO_\fX$\+modules
on\/~$\fX$.
 Let\/ $\bD$ be an exhaustive diagram of affine open formal
subschemes in\/~$\fX$.
 Then there is a natural adjunction isomorphism of abelian groups
\begin{equation} \label{fHom-contratensor-adjunction-eqn}
 \Hom_\fX(\cM\ocn_{\fX,\bD}\fP,\>\cJ)\simeq
 \Hom^{\fO_\fX}(\fP,\fHom_\fX(\cM,\cJ)).
\end{equation}
\end{lem}

\begin{proof}
 The point is that, for any cosheaf of $\fO_\fX$\+modules $\fP$ and
any contraherent cosheaf of contramodules $\fQ$ on $\fX$, the group of
morphisms $\Hom^{\fO_\fX}(\fP,\fQ)$ in the category $(\fX,\fO_\fX)\Cosh$
can be computed as the group of all families of $\fO_\fX(\fU_a)$\+module
maps $\fP[\fU_a]\rarrow\fQ[\fU_a]$, defined for all $a\in\bD$ and
satisfying the obvious compatibility condition for all $b\le a\in\bD$.
 Using this fact, the right-hand side
of~\eqref{fHom-contratensor-adjunction-eqn} is identified with the group
of all compatible collections of morphisms of $\fO_\fX(\fU_a)$\+modules
$\fP[\fU]\rarrow\Hom_\fX(\fj_a{}_*\fj_a^*\cM,\cJ)$.
 On the other hand, by the definition, the left-hand side
of~\eqref{fHom-contratensor-adjunction-eqn} is the group of all
compatible collections of morphisms of quasi-coherent torsion sheaves
$(\fj_a{}_*\fj_a^*\cM)\ot_{\fO_\fX(\fU_a)}\fP[\fU_a]\rarrow\cJ$.
 It remains to refer to
Lemma~\ref{torsion-sheaf-tensored-with-module-adjunction}.
\end{proof}

 It follows immediately from
Lemma~\ref{fHom-contratensor-adjunction-lemma} that the quasi-coherent
torsion sheaf of contratensor product $\cM\ocn_{\fX,\bD}\fP$ \emph{does
not depend} on the choice of a diagram $\bD$ \emph{provided that
the diagram\/ $\bD$ is exhaustive}.
 So we put
$$
 \cM\ocn_\fX\fP=\cM\ocn_{\fX,\bD}\fP
$$
for any exhaustive diagram $\bD$ of affine open formal subschemes
in~$\fX$.
 For example, one can use the poset of all affine open formal subschemes
$\fU\subset\fX$, ordered by inclusion, in the role of~$\bD$.

\subsection{Projection formulas}  \label{projection-formulas-subsecn}
 In this section we present the formal scheme counterparts of some
natural isomorphisms from~\cite[Section~3.8]{Pcosh}
and~\cite[Section~5.8]{Pdomc}.

\begin{prop} \label{fHom-projection-formula-prop}
 Let\/ $\ff\:\fY\rarrow\fX$ be a flat tight quasi-compact morphism of
locally Noetherian formal schemes.
 Let\/ $\cM$ be a quasi-coherent torsion sheaf on\/ $\fX$ and\/ $\cJ$
be an injective quasi-coherent torsion sheaf on\/~$\fY$.
 Then there is a natural isomorphism
\begin{equation} \label{fHom-projection-formula-eqn}
 \fHom_\fX(\cM,\ff_*\cJ)\simeq\ff_!\fHom_\fY(\ff^*\cM,\cJ).
\end{equation}
of locally cotorsion contraherent cosheaves of contramodules on\/~$\fX$.
\end{prop}

\begin{proof}
 This is our version of~\cite[formula~(3.24) in Section~3.8]{Pcosh}
and~\cite[formula~(49) in Section~5.8]{Pdomc}.
 First of all, the direct image functor~$\ff_*$ is right adjoint to
the exact functor of inverse image $\ff^*\:\fX\Tors\rarrow\fY\Tors$
(see the discussion in
Section~\ref{inverse-images-of-qcoh-tors-subsecn} and
formula~\eqref{qcoh-tors-inverse-direct-adjunction}).
 Hence the functor $\ff_*\:\fY\Tors\rarrow\fX\Tors$ takes injective
objects to injective objects
(Lemma~\ref{injective-torsion-sheaves-direct-image-lem} below),
the quasi-coherent torsion sheaf $\ff_*\cJ$ on $\fX$ is injective,
and the left-hand side of~\eqref{fHom-projection-formula-eqn}
is well-defined.
 By construction, the right-hand side
of~\eqref{fHom-projection-formula-eqn} is a locally cotorsion
(globally) contraherent cosheaf of contramodules on $\fX$ when
the morphism~$\ff$ is affine (by
Lemma~\ref{lcth-contram-W-T-affine-direct-image-lemma}(b)) and
a cosheaf of contramodule $\fO_\fX$\+modules on $\fX$ in the general
case (by Lemma~\ref{contramodule-cosheaves-direct-images-lemma}).
 After the isomorphism~\eqref{fHom-projection-formula-eqn} is proved,
it will follow that the right-hand side
of~\eqref{fHom-projection-formula-eqn} is actually a locally cotorsion
contraherent cosheaf of contramodules on $\fX$ in the general case.

 Indeed, let $\fU\subset\fX$ be an affine open formal subscheme with
the open immersion morphism $\fj\:\fU\rarrow\fX$.
 Denote by $\fV\subset\fY$ the open formal subscheme $\fV=\ff^{-1}(\fU)
\subset\fY$.
 Let $\fj'\:\fV\rarrow\fY$ and $\ff'\:\fV\rarrow\fU$ be the natural
morphisms.
 Then we have natural isomorphisms of $\fO_\fX(\fU)$\+modules
\begin{multline} \label{fHom-projection-formula-computation}
 \fHom_\fX(\cM,\ff_*\cJ)[\fU]=\Hom_\fX(\fj_*\fj^*\cM,\ff_*\cJ)
 \simeq\Hom_\fY(\ff^*\fj_*\fj^*\cM,\cJ) \\
 \simeq\Hom_\fY(\fj'_*\ff'{}^*\fj^*\cM,\cJ)\simeq
 \Hom_\fY(\fj'_*\fj'{}^*\ff^*\cM,\cJ)\simeq\fHom_\fY(\ff^*\cM,\cJ)[\fV],
\end{multline}
as desired.

 Here the first isomorphism (equality)
in~\eqref{fHom-projection-formula-computation} holds by the definition
and the second one by
the adjunction~\eqref{qcoh-tors-inverse-direct-adjunction}.
 The third (middle) isomorphism
in~\eqref{fHom-projection-formula-computation} is provided by
Lemma~\ref{tors-openimmers-dir-tight-inv-basechange-lemma}(b),
and the fourth one holds due to the isomorphism of the compositions
of inverse image functors $\ff'{}^*\fj^*\simeq\fj'{}^*\ff^*$,
which follows from the equality of the morphisms of formal schemes
$\fj\ff'=\ff\fj'$.
 The final fifth isomorphism
in~\eqref{fHom-projection-formula-computation} is provided by
Proposition~\ref{cosections-of-fHom-prop} (for the open immersion
morphism $\fj'\:\fV\rarrow\fY$).
\end{proof}

\begin{lem} \label{cosheaf-cosections-as-inductive-image-lemma}
 Let\/ $\fX$ be a topological space with a base of open subsets\/ $\bB$,
and let\/ $\fP$ be a cosheaf of abelian groups on\/~$\fX$.
 Then, for any open subset\/ $\fY\subset\fX$, the natural map
$$
 \varinjlim\nolimits_{\fU\subset\fY}^{\fU\in\bB}\fP[\fU]
 \lrarrow\fP[\fY]
$$
is an isomorphism.
 There the inductive limit is taken over the (nondirected) poset of
all open subsets\/ $\fU\subset\fX$ such that\/ $\fU\subset\fY$ and\/
$\fU\in\bB$, which is ordered by inclusion.
\end{lem}

\begin{proof}
 It may be more intuitive to pass from a cosheaf $\fP$ to the sheaf
$\cM$ assigning to every open subset $\fV\subset\fX$ the abelian group
$\cM(\fV)=\Hom_\boZ(\fP[\fV],I)$, for any chosen abelian group~$I$
(or just $I=\boQ/\boZ$), as in~\cite[proof of Theorem~2.1.2]{Pcosh}.
 Then the question reduces to the essentially obvious assertion
that the natural map
$$
 \cM(\fY)\lrarrow
 \varprojlim\nolimits_{\fU\subset\fY}^{\fU\in\bB}\cM(\fU)
$$
is an isomorphism for any sheaf of abelian groups $\cM$ on~$\fX$.
\end{proof}

\begin{prop} \label{contratensor-projection-formula-prop}
\textup{(a)} Let\/ $\ff\:\fY\rarrow\fX$ be a tight quasi-compact
morphism of locally Noetherian formal schemes.
 Let\/ $\cM$ be a quasi-coherent torsion sheaf on\/ $\fX$ and\/ $\fQ$
be a cosheaf of\/ $\fO_\fY$\+modules on\/~$\fY$.
 Then there is a natural morphism
\begin{equation} \label{contratensor-projection-formula-map-eqn}
 \cM\ocn_\fX\ff_!\fQ\lrarrow\ff_*(\ff^*\cM\ocn_\fY\fQ)
\end{equation}
of quasi-coherent torsion sheaves on~$\fX$. \par
\textup{(b)} Assuming that the morphism of formal schemes\/ $\ff$ is
flat and affine, the natural morphism of quasi-coherent torsion
sheaves~\eqref{contratensor-projection-formula-map-eqn}
is an isomorphism. \par
\textup{(c)} For any open immersion\/ $\fh\:\fY\rarrow\fX$ of
an \emph{affine} formal scheme\/ $\fY$ into a locally Noetherian formal
scheme\/ $\fX$, any quasi-coherent torsion sheaf\/ $\cM$ on\/ $\fX$,
and any cosheaf of contramodule\/ $\fO_\fY$\+modules\/ $\fF$ on\/ $\fY$
such that the contramodule\/ $\fO_\fY(\fY)$\+module\/ $\fF[\fY]$ is
flat, the morphism~\eqref{contratensor-projection-formula-map-eqn} is
an isomorphism, that is
\begin{equation} \label{contratensor-projection-formula-iso-eqn}
 \cM\ocn_\fX\fh_!\fF\simeq\fh_*(\fh^*\cM\ocn_\fY\fF).
\end{equation}
\end{prop}

\begin{proof}
 This is our version of~\cite[formulas~(3.25) and~(3.26) in
Section~3.8]{Pcosh} and~\cite[formula~(50) in Section~5.8]{Pdomc}.

 Part~(a): by the definition and by
Lemma~\ref{cosheaf-cosections-as-inductive-image-lemma}, we have
natural isomorphisms of quasi-coherent torsion sheaves on~$\fX$
\begin{multline} \label{contratensor-projection-morphism-constructed1}
 \cM\ocn_\fX\ff_!\fQ=\varinjlim\nolimits_{\fU\subset\fX}
 \bigl((\fj_*\fj^*\cM)\ot_{\fO_\fX(\fU)}\fQ[\ff^{-1}(\fU)]\bigr) \\
 \simeq\varinjlim\nolimits_{\fU\subset\fX}
 \bigl((\fj_*\fj^*\cM)\ot_{\fO_\fX(\fU)}
 \varinjlim\nolimits_{\fV\subset\ff^{-1}(\fU)}\fQ[\fV]\bigr) \\
 \simeq\varinjlim\nolimits_{\fU\subset\fX}
 \varinjlim\nolimits_{\fV\subset\ff^{-1}(\fU)}
 \bigl((\fj_*\fj^*\cM)\ot_{\fO_\fX(\fU)}\fQ[\fV]\bigr)
 \simeq\varinjlim\nolimits_{\fV\subset\ff^{-1}(\fU)}
 \bigl((\fj_*\fj^*\cM)\ot_{\fO_\fX(\fU)}\fQ[\fV]\bigr) \\
 \simeq\varinjlim\nolimits_{\fV\subset\ff^{-1}(\fU)}
 \Bigl(\bigl((\fj_*\fj^*\cM)\ot_{\fO_\fX(\fU)}\fO_\fY(\fV)\bigr)
 \ot_{\fO_\fY(\fV)}\fQ[\fV]\Bigr).
\end{multline}
 Here the inductive limit in the second formula
in~\eqref{contratensor-projection-morphism-constructed1} is taken
over all the affine open formal subschemes $\fU\subset\fX$,
the inductive limits in the third and fourth formulas are taken over all
the affine open formal subschemes $\fU\subset\fX$ and over all
the affine open formal subschemes $\fV\subset\ff^{-1}(\fU)\subset\fY$,
and the inductive limits in the last two (fifth and sixth) formulas
in~\eqref{contratensor-projection-morphism-constructed1} are taken
over all pairs $(\fU,\fV)$, where $\fU\subset\fX$ and $\fV\subset\fY$
are affine open formal subschemes such that $\fV\subset\ff^{-1}(\fU)$.
 The posets of affine open formal subschemes or pairs of affine open
formal subschemes are ordered by inclusion.
 The notation~$\fj$ stands for the open immersion morphism
$\fj\:\fU\rarrow\fX$.

 The first isomorphism (equality)
in~\eqref{contratensor-projection-morphism-constructed1} holds by
the definitions of the contratensor product and the direct image,
and the second one by
Lemma~\ref{cosheaf-cosections-as-inductive-image-lemma}.
 The third isomorphism holds because the functor of tensor product
of additive category objects and modules preserves colimits quite
generally.
 The fourth isomorphism
in~\eqref{contratensor-projection-morphism-constructed1} is
an obvious general property of inductive limits, and the last (fifth)
isomorphism is an obvious general property of the functors of
tensor product of additive category objects and modules.

 Furthermore, we have a natural morphism in $\fX\Tors$
\begin{multline} \label{contratensor-projection-morphism-constructed2}
 \varinjlim\nolimits_{\fV\subset\ff^{-1}(\fU)}
 \Bigl(\bigl((\fj_*\fj^*\cM)\ot_{\fO_\fX(\fU)}\fO_\fY(\fV)\bigr)
 \ot_{\fO_\fY(\fV)}\fQ[\fV]\Bigr) \\
 \lrarrow\varinjlim\nolimits_{\fV\subset\ff^{-1}(\fU)}
 \Bigl(\fj_*\bigl((\fj^*\cM)\ot_{\fO_\fX(\fU)}\fO_\fY(\fV)\bigr)
 \ot_{\fO_\fY(\fV)}\fQ[\fV]\Bigr)
\end{multline}
induced by the natural
transformation~\eqref{sheaf-module-tensor-direct-image} (for
the quasi-compact morphism of locally Noetherian formal schemes
$\fj\:\fU\rarrow\fX$).

 Introducing further notation $\fj'\:\fV\rarrow\fY$ and
$\ff'\:\fV\rarrow\fU$ for the natural morphisms of formal schemes,
we obviously have $(\fj^*\cM)\ot_{\fO_\fX(\fU)}\fO_\fY(\fV)\simeq
\ff'_*\ff'{}^*\fj^*\cM$, hence
$$
 \fj_*\bigl((\fj^*\cM)\ot_{\fO_\fX(\fU)}\fO_\fY(\fV)\bigr)\simeq
 \fj_*\ff'_*\ff'{}^*\fj^*\cM\simeq\ff_*\fj'_*\fj'{}^*\ff^*\cM,
$$
since $\fj_*\ff'_*\simeq\ff_*\fj'_*$ and
$\ff'{}^*\fj^*\simeq\fj'{}^*\ff^*$.
 So there is a natural isomorphism
\begin{multline} \label{contratensor-projection-morphism-constructed3}
 \varinjlim\nolimits_{\fV\subset\ff^{-1}(\fU)}
 \Bigl(\fj_*\bigl((\fj^*\cM)\ot_{\fO_\fX(\fU)}\fO_\fY(\fV)\bigr)
 \ot_{\fO_\fY(\fV)}\fQ[\fV]\Bigr) \\
 \simeq\varinjlim\nolimits_{\fV\subset\ff^{-1}(\fU)}
 \bigl((\ff_*\fj'_*\fj'{}^*\ff^*\cM)\ot_{\fO_\fY(\fV)}\fQ[\fV]\bigr)
\end{multline}
in $\fX\Tors$.

 Finally, we have natural morphisms in $\fX\Tors$
\begin{multline} \label{contratensor-projection-morphism-constructed4}
 \varinjlim\nolimits_{\fV\subset\ff^{-1}(\fU)}
 \bigl((\ff_*\fj'_*\fj'{}^*\ff^*\cM)\ot_{\fO_\fY(\fV)}\fQ[\fV]\bigr)
 \\ \lrarrow \varinjlim\nolimits_{\fV\subset\ff^{-1}(\fU)}
 \Bigl(\ff_*\bigl((\fj'_*\fj'{}^*\ff^*\cM)
 \ot_{\fO_\fY(\fV)}\fQ[\fV]\bigr)\Bigr) \\
 \lrarrow\ff_*\varinjlim\nolimits_{\fV\subset\ff^{-1}(\fU)}
 \bigl((\fj'_*\fj'{}^*\ff^*\cM)\ot_{\fO_\fY(\fV)}\fQ[\fV]\bigr),
\end{multline}
the former of which is induced by the natural
transformation~\eqref{sheaf-module-tensor-direct-image} (for
the quasi-compact morphism of locally Noetherian formal schemes
$\ff\:\fY\rarrow\fX$), while the latter one is the standard
category-theoretic constrution of a natural transformation arising
when one applies a functor to a colimit.

 It remains to mention the natural isomorphism
\begin{equation} \label{contratensor-projection-morphism-constructed5}
 \ff_*\varinjlim\nolimits_{\fV\subset\ff^{-1}(\fU)}
 \bigl((\fj'_*\fj'{}^*\ff^*\cM)\ot_{\fO_\fY(\fV)}\fQ[\fV]\bigr)
 \simeq\ff_*(\ff^*\cM\ocn_\fY\fQ)
\end{equation}
in $\fX\Tors$ following from the natural isomorphism
$$
 \varinjlim\nolimits_{\fV\subset\ff^{-1}(\fU)}
 \bigl((\fj'_*\fj'{}^*\ff^*\cM)\ot_{\fO_\fY(\fV)}\fQ[\fV]\bigr)
 \simeq\ff^*\cM\ocn_\fY\fQ
$$
in $\fY\Tors$.
 The isomorphism in $\fY\Tors$ can be obtained by computing
the contratensor product $\ff^*\cM\ocn_\fY\fQ$ on the exhaustive diagram
of affine open formal subschemes $\fV\subset\fY$ indexed by the poset
$\bD$ of all pairs $(\fU,\fV)$ of affine open formal subschemes
$\fU\subset\fX$ and $\fV\subset\fY$ such that $\fV\subset\ff^{-1}(\fU)$.
 See the discussion in Section~\ref{contratensor-product-subsecn}.

 The desired natural
transformation~\eqref{contratensor-projection-formula-map-eqn}
is now constructed as the composition of the natural morphisms and
isomorphisms~\eqref{contratensor-projection-morphism-constructed1},
\eqref{contratensor-projection-morphism-constructed2},
\eqref{contratensor-projection-morphism-constructed3},
\eqref{contratensor-projection-morphism-constructed4},
and~\eqref{contratensor-projection-morphism-constructed5}.

 Part~(b): the argument is based on the proof of part~(a).
 It suffices to show that, under the assumptions of~(b), all
the morphisms
in~\eqref{contratensor-projection-morphism-constructed2}
and~\eqref{contratensor-projection-morphism-constructed4}
are isomorphisms.

 Indeed, whenever the morphism~$\ff$ is flat, the contramodule
$\fO_\fX(\fU)$\+module $\fO_\fY(\fV)$ is flat, too.
 According to the discussion of the natural
transformation~\eqref{sheaf-module-tensor-direct-image}
in Section~\ref{contraherent-Hom-subsecn}, it follows that
the morphism~\ref{contratensor-projection-morphism-constructed2}
is an isomorphism.

 Whenever the morphism~$\ff$ is affine, the direct image
functor~$\ff_*$ is exact and preserves infinite direct sums,
so it preserves all colimits. 
 Hence both the morphisms
in~\eqref{contratensor-projection-morphism-constructed4}
are isomorphisms in this case.

 Notice that, in the case of an affine morphism~$\ff$, in
the context of the proof of part~(a), the open formal subscheme
$\ff^{-1}(\fU)\subset\fY$ is affine.
 For this reason, the argument in part~(a) can be simplified by
setting $\fV=\ff^{-1}(\fU)$ once and for all, rather than letting
$\fV$ vary over the poset of all affine open formal subschemes
in $\ff^{-1}(\fU)$.
 All the inductive limits involved in the proof of part~(a) remain
unchanged by such a specialization.
 The diagram of affine open formal subschemes $\ff^{-1}(\fU)\subset\fY$
indexed by the poset $\bD$ of all affine open formal subschemes
$\fU\subset\fX$ is exhaustive in $\fY$, so it can be used for
computation of contratensor products over~$\fY$.
 This produces a somewhat simpler argument for the proof of part~(b).

 To prove part~(c), we also follow the construction in the proof
of part~(a).
 The morphism $\ff=\fh$ is an open immersion now, so it is flat,
and the morphism~\eqref{contratensor-projection-morphism-constructed2}
is an isomorphism according to the proof of part~(b).

 Concerning~\eqref{contratensor-projection-morphism-constructed4},
under the assumptions of~(c) all the inductive limits involved
agree with the similar inductive limits taken over the poset of all
affine open formal subschemes $\fV\subset\fY$ (since one can always
take $\fU=\fh(\fV)$).
 The latter inductive limits can be simply computed by specializing
to $\fV=\fY$.
 This proves that the second morphism
in~\eqref{contratensor-projection-morphism-constructed4} is
an isomorphism.
 Concerning the first morphism (for $\ff=\fh$, \ $\fj'=\id_\fY$,
and $\fQ=\fF$), it remains to observe that the morphism
$$
 (\fh_*\fh^*\cM)\ot_{\fO_\fY(\fY)}\fF[\fY]
 \lrarrow\fh_*\bigl((\fh^*\cM)\ot_{\fO_\fY(\fY)}\fF[\fY]\bigl)
$$
is an isomorphism according to the discussion of the natural
transformation~\eqref{sheaf-module-tensor-direct-image}
in Section~\ref{contraherent-Hom-subsecn} (as the contramodule
$\fO_\fY(\fY)$\+module $\fF[\fY]$ is flat by assumption).
\end{proof}

\Section{Injective Torsion Sheaves and Projective Cosheaves of
Contramodules}  \label{injective-projective-secn}

\subsection{Semi-separated Noetherian formal schemes}
\label{semi-separated-inj-proj-subsecn}
 In this section we describe injective quasi-coherent torsion sheaves
on Noetherian formal schemes, as well as projective contraherent
cosheaves of contramodules and projective locally cotorsion contraherent
cosheaves of contramodules on semi-separated Noetherian formal schemes.
 The descriptions provided in this section are formal scheme versions of
results applicable to arbitrary (not necessarily Noetherian)
quasi-compact and quasi-separated or semi-separated schemes,
as in~\cite[Section~4.5]{Pcosh}.

\begin{lem} \label{injective-torsion-sheaves-on-qcomp}
 Let\/ $\fX$ be a Noetherian formal scheme with a finite affine open
covering\/ $\fX=\bigcup_{\alpha=1}^N\fU_\alpha$.
 Then the abelian category\/ $\fX\Tors$ of quasi-coherent torsion
sheaves on\/ $\fX$ has enough injective objects.
 A quasi-coherent torsion sheaf on\/ $\fX$ is injective if and only if
it is a direct summand of a finite direct sum over\/~$\alpha$ of
direct images of injective quasi-coherent torsion sheaves
from\/~$\fU_\alpha$.
\end{lem}

\begin{proof}
 For any locally Noetherian formal scheme $\fX$, the category of
quasi-coherent torsion sheaves $\fX\Tors$ is a Grothendieck abelian
category, hence it has enough injective objects (cf.\
Theorem~\ref{loc-Noetherian-injective-qcoh-torsion} below).
 Alternatively, in the case of a Noetherian formal scheme $\fX$,
the existence of enough injective quasi-coherent torsion sheaves on
$\fX$ will be established by the construction below in this proof.

 Denote by $\fj_\alpha$ the open immersion morphism
$\fU_\alpha\rarrow\fX$.
 Then the direct image functor $\fj_\alpha{}_*\:\fU_\alpha\Tors\rarrow
\fX\Tors$ is right adjoint to the inverse image/restriction functor
$\fj_\alpha^*\:\fX\Tors\rarrow\fU_\alpha\Tors$ (see
Proposition~\ref{torsion-sheaves-qcoh-direct-images-prop} and
formula~\eqref{direct-image-qcoh-tors}, together with
formula~\eqref{sheaves-restriction-direct-image-adjunction}
or~\eqref{qcoh-tors-inverse-direct-adjunction}).
 Since the inverse image functor~$\fj_\alpha^*$ is exact, it follows
that the direct image functor~$\fj_\alpha{}_*$ takes injective objects
to injective objects.
 This proves the ``if'' assertion of the lemma.

 To prove the ``only if'', it suffices to show that every
quasi-coherent torsion sheaf $\cM$ on $\fX$ can be embedded into
a direct sum over~$\alpha$ of the direct images of injective
quasi-coherent torsion sheaves from~$\fU_\alpha$.
 Indeed, we have category equivalences $\fU_\alpha\Tors\simeq
\fO(\fU_\alpha)\Tors$ as per Section~\ref{qcoh-torsion-sheaves-subsecn},
and the abelian categories $\fO(\fU_\alpha)\Tors$ have enough injective
objects according to
Section~\ref{prelim-inj-proj-flat-torsion-contra-subsecn}.
 Pick injective quasi-coherent torsion sheaves $\cJ_\alpha$ on
$\fU_\alpha$ together with injective morphisms of quasi-coherent
torsion sheaves $\fj_\alpha^*\cM\rarrow\cJ_\alpha$ on~$\fU_\alpha$.

 By adjunction, we have the induced morphisms of injective
quasi-coherent torsion sheaves $\cM\rarrow \fj_\alpha{}_*\cJ_\alpha$
on~$\fX$.
 It remains to point out that the resulting morphism
$\cM\rarrow\bigoplus_\alpha\fj_\alpha{}_*\cJ_\alpha$ is injective
in $\fX\Tors$, since the restriction functor~$\fj_\alpha^*$ takes
it into an injective morphism in $\fU_\alpha\Tors$ for every~$\alpha$.
\end{proof}

\begin{lem} \label{projective-contrah-on-qcomp-ssep}
 Let\/ $\fX$ be a semi-separated Noetherian formal scheme with an open
covering\/ $\bW$ and a finite affine open covering\/
$\fX=\bigcup_{\alpha=1}^N\fU_\alpha$ subordinate to\/~$\bW$.
 Then the exact category\/ $\fX\Lcth_\bW$ of\/ $\bW$\+locally
contraherent cosheaves of contramodules on\/ $\fX$ has enough
projective objects.
 A\/ $\bW$\+locally contraherent cosheaf of contramodules on\/ $\fX$
is projective if and only if it is a direct summand of a finite direct
sum over\/~$\alpha$ of direct images of projective contraherent
cosheaves of contramodules from\/~$\fU_\alpha$.
\end{lem}

\begin{proof}
 This is our version of~\cite[Lemma~4.5.1]{Pcosh}.
 We keep the notation~$\fj_\alpha$ for the open immersion morphisms
$\fU_\alpha\rarrow\fX$.
 The direct image functor $\fj_\alpha{}_!\:\fU_\alpha\Ctrh\rarrow
\fX\Lcth_\bW$ is left adjoint to the inverse image/restriction functor
$\fj_\alpha^!\:\fX\Lcth_\bW\rarrow\fU_\alpha\Ctrh$ (see
Lemma~\ref{lcth-contram-W-T-affine-direct-image-lemma}(a)
and formula~\eqref{loc-ctrh-contramod-W-T-affine-direct-image}, together
with formula~\eqref{cosheaves-direct-image-restriction-adjunction}
or~\eqref{lcth-contramod-direct-inverse-adjunction}).
 Since the inverse image functor~$\fj_\alpha^!$ is exact (as a functor
of exact categories), it follows that the direct image
functor~$\fj_\alpha{}_!$ takes projective objects to projective objects.
 This proves the ``if'' assertion.

 In order to prove both the ``only if'' and the existence of enough
projectives in $\fX\Lcth_\bW$, it remains to show that every
$\bW$\+locally contraherent cosheaf of contramodules $\fM$ on $\fX$ is
the target of an admissible epimorphism from a direct sum over~$\alpha$
of the direct images of projective contraherent cosheaves of
contramodules from~$\fU_\alpha$.
 Indeed, we have exact category equivalences $\fU_\alpha\Ctrh\simeq
\fO(\fU_\alpha)\Contra^\cta$ as per
Lemma~\ref{contrah-cosheaves-of-contramods-lemma}(b).
 The exact categories $\fO(\fU_\alpha)\Contra^\cta$ have enough
projective objects by
Proposition~\ref{very-flat-cotorsion-pair-in-quotseparated}(a);
specifically, the very flat contraadjusted contramodule
$\fO(\fU_\alpha)$\+modules are the projective objects of
$\fO(\fU_\alpha)\Contra^\cta$ (see
Lemma~\ref{cotorsion-pair-kernel-as-injectives-projectives}(b) below).
 Pick projective contraherent cosheaves of contramodules $\fF_\alpha$
on $\fU_\alpha$ together with admissible epimorphisms $\fF_\alpha
\rarrow \fj_\alpha^!\fM$ of contraherent cosheaves of contramodules
on~$\fU_\alpha$.

 By adjunction, we have the induced morphisms of $\bW$\+locally
contraherent cosheaves of contramodules $\fj_\alpha{}_!\fF_\alpha
\rarrow\fM$ on~$\fX$.
 It remains to check that the resulting morphism $\bigoplus_\alpha
\fj_\alpha{}_!\fF_\alpha\rarrow\fM$ is an admissible epimorphism in
$\fX\Lcth_\bW$.

 For this purpose, one can say that the map $\bigoplus_\alpha
\fj_\alpha{}_!\fF_\alpha\rarrow\fM$ is the composition
$\bigoplus_\alpha\fj_\alpha{}_!\fF_\alpha\rarrow
\bigoplus_\alpha\fj_\alpha{}_!\fj_\alpha^!\fM\rarrow\fM$.
 The morphism $\bigoplus_\alpha\fj_\alpha{}_!\fF_\alpha\rarrow
\bigoplus_\alpha\fj_\alpha{}_!\fj_\alpha^!\fM$ is an admissible
epimorphism since the functor~$\fj_\alpha{}_!$ is exact
(see Section~\ref{direct-images-of-ctrh-contramod-subsecn}), while
the morphism $\bigoplus_\alpha\fj_\alpha{}_!\fj_\alpha^!\fM\rarrow\fM$
is an admissible epimorphism by
formula~\eqref{Cech-sequence-of-cosheaves-of-contramods}.
 Alternatively, it suffices to point out that the morphism
$\bigoplus_\alpha\fj_\alpha{}_!\fF_\alpha\rarrow\fM$ is an admissible
epimorphism in restriction to every $\fU_\alpha$ and refer to
the locality of admissible epimorphisms in $\fX\Lcth_\bW$ (see
Section~\ref{locally-contraherent-cosheaves-of-contramods-subsecn}).
\end{proof}

\begin{cor} \label{projective-contrah-on-qcomp-ssep-corollary}
 Let\/ $\fX$ be a semi-separated Noetherian formal scheme with an open
covering\/~$\bW$.
 Then all the three exact categories\/ $\fX\Ctrh\subset\fX\Lcth_\bW
\subset\fX\Lcth$ have enough projective objects, and the classes of
projective objects in all the three exact categories coincide.
 In particular, all the projective objects in\/ $\fX\Lcth_\bW$ are
(globally) contraherent cosheaves of contramodules on\/~$\fX$.
\end{cor}

\begin{proof}
 Follows from Lemma~\ref{projective-contrah-on-qcomp-ssep}.
\end{proof}

\begin{lem} \label{projective-contrah-lct-on-qcomp-ssep}
 Let\/ $\fX$ be a semi-separated Noetherian formal scheme with an open
covering\/ $\bW$ and a finite affine open covering\/
$\fX=\bigcup_{\alpha=1}^N\fU_\alpha$ subordinate to\/~$\bW$.
 Then the exact category\/ $\fX\Lcth^\lct_\bW$ of\/ locally cotorsion
$\bW$\+locally contraherent cosheaves of contramodules on\/ $\fX$ has
enough projective objects.
 A locally cotorsion\/ $\bW$\+locally contraherent cosheaf of
contramodules on\/ $\fX$ is projective if and only if it is a direct
summand of a finite direct sum over\/~$\alpha$ of direct images of
projective locally cotorsion contraherent cosheaves of contramodules
from\/~$\fU_\alpha$.
\end{lem}

\begin{proof}
 This is our version of~\cite[Lemma~4.5.3]{Pcosh}.
 The argument is similar to the proof of
Lemma~\ref{projective-contrah-on-qcomp-ssep}.
 One needs to use
Lemma~\ref{lcth-contram-W-T-affine-direct-image-lemma}(b),
formula~\eqref{loc-ctrh-lct-contramod-W-T-affine-direct-image},
Lemma~\ref{lct-contrah-cosheaves-of-contramods-on-affine-lemma},
Proposition~\ref{flat-cotorsion-pair-in-quotseparated}(a), and
the discussion of exact categories of locally contraherent cosheaves
in Section~\ref{locally-contraherent-lct-cosheaves-subsecn}.
 The flat cotorsion contramodule $\fO(\fU_\alpha)$\+modules are
the projective objects of $\fO(\fU_\alpha)\Contra^\cot\simeq
\fU_\alpha\Ctrh^\lct$, see
Lemma~\ref{cotorsion-pair-kernel-as-injectives-projectives}(b) below.
\end{proof}

\begin{cor} \label{projective-contrah-lct-on-qcomp-ssep-corollary}
 Let\/ $\fX$ be a semi-separated Noetherian formal scheme with an open
covering\/~$\bW$.
 Then all the three exact categories\/ $\fX\Ctrh^\lct\subset
\fX\Lcth^\lct_\bW\subset\fX\Lcth^\lct$ have enough projective objects,
and the classes of projective objects in all the three exact categories
coincide.
 In particular, all the projective objects in\/ $\fX\Lcth^\lct_\bW$ are
locally cotorsion (globally) contraherent cosheaves of contramodules
on\/~$\fX$.
\end{cor}

\begin{proof}
 Follows from Lemma~\ref{projective-contrah-lct-on-qcomp-ssep}.
\end{proof}

 Following the general notational convention of
Section~\ref{prelim-inj-proj-flat-torsion-contra-subsecn}, we will
denote the additive category of injective quasi-coherent torsion
sheaves on a locally Noetherian formal scheme $\fX$ by
$\fX\Tors^\inj\subset\fX\Tors$.
 Similarly, the additive category of projective contraherent cosheaves
of contramodules on $\fX$ will be denoted by $\fX\Ctrh_\prj\subset
\fX\Ctrh$, and the additive category of projective locally cotorsion
contraherent cosheaves of contramodules on $\fX$ will be denoted by
$\fX\Ctrh^\lct_\prj\subset\fX\Ctrh^\lct$.
 Let us \emph{warn} the reader that projective locally cotorsion
contraherent cosheaves of contramodules are usually \emph{not}
projective as contraherent cosheaves of contramodules, that is
$\fX\Ctrh^\lct_\prj\not\subset\fX\Ctrh_\prj$
\,\cite[Section~4.5]{Pcosh}.

\begin{lem} \label{injective-torsion-sheaves-direct-image-lem}
 Let\/ $\ff\:\fY\rarrow\fX$ be a flat quasi-compact morphism of
locally Noetherian formal schemes.
 Then the functor of direct image of quasi-coherent torsion sheaves\/
$\ff_*\:\fY\Tors\rarrow\fX\Tors$ takes injective quasi-coherent torsion
sheaves on\/ $\fY$ to injective quasi-coherent torsion sheaves
on\/~$\fX$.
\end{lem}

\begin{proof}
 This essentially well known fact was already mentioned in
the proof of Proposition~\ref{fHom-projection-formula-prop}.
 According to formula~\eqref{qcoh-tors-inverse-direct-adjunction},
the direct image functor~$\ff_*$ is right adjoint to the inverse image
functor $\ff^*\:\fX\Tors\rarrow\fY\Tors$ (recall that any flat
morphism of formal schemes is tight by the definition).
 Since the morphism~$\ff$ is flat, the functor~$\ff^*$ is exact
(as per Section~\ref{inverse-images-of-qcoh-tors-subsecn}).
 Hence the right adjoint functor to~$\ff^*$ takes injective objects
to injective objects.
\end{proof}

\begin{cor} \label{ssep-proj-contrah-direct-image-cor}
 Let\/ $\ff\:\fY\rarrow\fX$ be a morphism of semi-separated Noetherian
formal schemes.
 In this setting: \par
\textup{(a)} if the morphism\/~$\ff$ is very flat, then the functor of
direct image of cosheaves of\/ $\fO$\+modules\/
$\ff_!\:(\fY,\fO_\fY)\Cosh\rarrow(\fX,\fO_\fX)\Cosh$ takes projective
contraherent cosheaves of contramodules on\/ $\fY$ to projective
contraherent cosheaves of contramodules on\/~$\fX$; \par
\textup{(b)} if the morphism\/~$\ff$ is flat, then the functor of direct
image of cosheaves of\/ $\fO$\+modules\/
$\ff_!\:(\fY,\fO_\fY)\Cosh\rarrow(\fX,\fO_\fX)\Cosh$ takes projective
locally cotorsion contraherent cosheaves of contramodules on\/ $\fY$ to
projective locally cotorsion contraherent cosheaves of contramodules
on\/~$\fX$.
\end{cor}

\begin{proof}
 This is our version of~\cite[Corollary~4.6.4]{Pcosh}.
 Let us first consider the case when the morphism~$\ff$ is affine.
 Then we have the direct image functors $\ff_!\:\fY\Ctrh\rarrow\fX\Ctrh$
\,\eqref{ctrh-contramod-affine-direct-image} and $\ff_!\:\fY\Ctrh^\lct
\rarrow\fX\Ctrh^\lct$ \,\eqref{ctrh-lct-contramod-affine-direct-image}
acting between the exact categories of contraherent cosheaves of
contramodules.

 Part~(a): by formula~\eqref{lcth-contramod-direct-inverse-adjunction},
the direct image functor $\ff_!\:\fY\Ctrh\rarrow\fX\Ctrh$ is
``partially left adjoint'' to the inverse image functor
$\ff^!\:\fX\Lcth\rarrow\fY\Lcth$, in the sense that the adjunction
isomorphism of abelian groups $\Hom^\fY(\fQ,\ff^!\fP)\simeq
\Hom^\fX(\ff_!\fQ,\fP)$ holds for all contraherent cosheaves of
contramodules $\fQ$ on $\fY$ and all locally contraherent cosheaves
of contramodules $\fP$ on~$\fX$.
 Furthermore, according to
Section~\ref{inverse-images-of-ctrh-contramod-subsecn}, the functor
$\ff^!\:\fX\Lcth\rarrow\fY\Lcth$
\,\eqref{inverse-image-lcth-contramod-nocovering} is exact.
 Finally, it is important for our argument that all the projective
objects of the exact category $\fY\Ctrh$ are also projective in
the exact category $\fY\Lcth$ (by
Corollary~\ref{projective-contrah-on-qcomp-ssep-corollary}).
 The desired assertion follows from these observations.
 The proof of part~(b) for affine morphisms~$\ff$ is similar and based
on exactness of
the functor~\eqref{inverse-image-lcth-lct-contramod-nocovering} and
Corollary~\ref{projective-contrah-lct-on-qcomp-ssep-corollary}.

 To prove the assertions of the corollary in the general case, notice 
that the functors of direct image of cosheaves of $\fO$\+modules are
compatible with the compositions of morphisms of ringed spaces.
 For any pair of composable morphisms of ringed spaces
$\fg\:(\fW,\fO_\fW)\rarrow(\fY,\fO_\fY)$ and
$\ff\:(\fY,\fO_\fY)\rarrow(\fX,\fO_\fX)$, one has
$(\ff\circ\fg)_!=\ff_!\circ\fg_!\:(\fW,\fO_\fW)\Cosh\rarrow
(\fX,\fO_\fX)\Cosh$.

 Now part~(a) follows from its affine morphism case together with
Lemma~\ref{projective-contrah-on-qcomp-ssep}, and part~(b)
follows from its affine morphism case together with
Lemma~\ref{projective-contrah-lct-on-qcomp-ssep}.
 The point is that if $\fY=\bigcup_\alpha\fV_\alpha$ is a finite
affine open covering of $\fY$ and $\fj_\alpha\:\fV_\alpha\rarrow
\fY_\alpha$ are the open immersion morphisms, then the compositions
$\ff\circ\fj_\alpha\:\fV_\alpha\rarrow\fX$ are affine morphisms of
formal schemes (since the formal scheme $\fX$ is semi-separated,
and any morphism from an affine scheme to a semi-separated one is
affine).
 The morphisms $\ff\circ\fj_\alpha$ are also very flat in the context
of part~(a) and flat in the context of part~(b).

 For a version of part~(a) with the semi-separatedness assumption
replaced by the assumption of finite Krull dimension, see
Corollary~\ref{flat-contrah-direct-images}(c) below.
 A far-reaching generalization of part~(b) will be presented in
Corollary~\ref{proj-lct-direct-images}(b).
\end{proof}

\subsection{Antilocally flat contraherent cosheaves of contramodules}
\label{antilocally-flat-subsecn}
 This section is a formal scheme version of~\cite[Section~4.4]{Pcosh}.

 Let $\fX$ be a locally Noetherian formal scheme with an open
covering~$\bW$.
 A $\bW$\+locally contraherent cosheaf of contramodules $\fF$ on $\fX$
is said to be \emph{antilocally flat} if, for any short exact sequence
$0\rarrow\fK\rarrow\fL\rarrow\fM\rarrow0$ of locally cotorsion
$\bW$\+locally contraherent cosheaves of contramodules on $\fX$,
the sequence of abelian groups $0\rarrow\Hom^\fX(\fF,\fK)\rarrow
\Hom^\fX(\fF,\fL)\rarrow\Hom^\fX(\fF,\fM)\rarrow0$ is exact.

 Clearly, all projective objects of the exact category $\fX\Lcth_\bW$,
as well as all projective objects of the exact category
$\fX\Lcth^\lct_\bW$, are antilocally flat as $\bW$\+locally
contraherent cosheaves of contramodules on~$\fX$.

 A contraherent cosheaf of contramodules $\fF$ on an affine formal
scheme $\fU$ with the open covering $\bW_\fU=\{\fU\}$ is antilocally
flat if and only if the contraadjusted contramodule
$\fO_\fU(\fU)$\+module $\fF[\fU]$ is flat (the ``if'' implication
follows from the definition of cotorsion contramodule
$\fO_\fU(\fU)$\+modules in
Lemma~\ref{cotorsion-contramodules-lemma}(2), and the ``only if''
holds by Corollary~\ref{flat-contramodules-cor}).
 It follows from the adjunction
isomorphism~\eqref{cosheaves-direct-inverse-adjunction}
and the formula~\eqref{inverse-image-lcth-lct-W-T-contramod}
that, for any flat tight $(\bW,\bT)$\+affine $(\bW,\bT)$\+coaffine
morphism of locally Noetherian formal schemes $\ff\:\fY\rarrow\fX$,
the functor of direct image of $\bT$\+locally contraherent cosheaves
of contramodules $\ff_!$
\,\eqref{loc-ctrh-contramod-W-T-affine-direct-image} takes antilocally
flat $\bT$\+locally contraherent cosheaves of contramodules on $\fY$
to antilocally flat $\bW$\+locally contraherent cosheaves of
contramodules on~$\fX$.

 Let $\fX$ be a semi-separated Noetherian formal scheme with an open
covering~$\bW$.
 We will prove in this section that any antilocally flat
$\bW$\+locally contraherent cosheaf of contramodules on $\fX$ is
(globally) contraherent.
 Moreover, the class of antilocally flat $\bW$\+locally contraherent
cosheaves of contramodules on $\fX$ does not depend on an open
covering~$\bW$; so it coincides with the class of antilocally flat
contraherent cosheaves of contramodules on~$\fX$.

\begin{lem} \label{locally-cotorsion-preenvelope-lemma}
 Let\/ $\fX$ be a semi-separated Noetherian formal scheme with an open
covering\/ $\bW$ and a finite affine open covering\/
$\fX=\bigcup_{\alpha=1}^N\fU_\alpha$ subordinate to\/~$\bW$.
 Then any\/ $\bW$\+locally contraherent cosheaf of contramodules\/
$\fM$ on\/ $\fX$ can be included in an (admissible) short exact
sequence\/ $0\rarrow\fM\rarrow\fP\rarrow\fF\rarrow0$ in\/
the exact category $\fX\Lcth_\bW$, where\/ $\fP$ is a locally
cotorsion\/ $\bW$\+locally contraherent cosheaf of contramodules on\/
$\fX$ and\/ $\fF$ is a finitely iterated extension of the direct images
of contraherent cosheaves of contramodules on\/ $\fU_\alpha$
corresponding to flat contraadjusted contramodule\/
$\fO_\fX(\fU_\alpha)$\+modules.
\end{lem}

\begin{proof}
 This is our version of~\cite[Lemma~4.4.1]{Pcosh}.
 The proof is similar to that of~\cite[Lemma~4.3.2]{Pcosh}.
 Let us spell out some details.

 Arguing by induction on $1\le\beta\le N$, we consider the open formal
subscheme $\fV=\bigcup_{\alpha<\beta}\fU_\alpha$ with the induced
open covering $\bW|_\fV=\{\fV\cap\fW\mid\fW\in\bW\}$ and the identity
open immersion morphism $\fh\:\fV\rarrow\fX$.
 Assume that we have constructed a short exact sequence
$0\rarrow\fM\rarrow\fQ\rarrow\fG\rarrow0$ of $\bW$\+locally contraherent
cosheaves of contramodules on $\fX$ such that the restriction
$\fh^!\fQ$ of the $\bW$\+locally contraherent cosheaf of contramodules
$\fQ$ to the open formal subscheme $\fV\subset\fX$ is locally cotorsion,
while the cosheaf $\fG$ on $\fX$ is a finitely iterated extension of
the direct images of antilocally flat contraherent cosheaves of
contramodules from the affine open formal subschemes $\fU_\alpha\subset
\fX$, \,$\alpha<\beta$.
 When $\beta=1$, it suffices to take $\fQ=\fM$ and $\fG=0$ for
the induction base.
 Set $\fU=\fU_\beta$ and denote by $\fj\:\fU\rarrow\fX$ the identity
open immersion morphism.

 Let $0\rarrow\fj^!\fQ\rarrow\fK\rarrow\fH\rarrow0$ be a short exact
sequence of contraherent cosheaves of contramodules on the affine
formal scheme $\fU$ such that the contraherent cosheaf of contramodules
$\fK$ is locally cotorsion, while the contraherent cosheaf of
contramodules $\fH$ is antilocally flat.
 The existence of such a short exact sequence in $\fU\Ctrh$ follows
from Lemmas~\ref{contrah-cosheaves-of-contramods-lemma}(b)
and~\ref{lct-contrah-cosheaves-of-contramods-on-affine-lemma} together
with Proposition~\ref{flat-cotorsion-pair-in-quotseparated}(b) and
the fact that the full subcategory $\fO(\fU)\Contra^\cta$ is closed
under cokernels (in fact, under quotients) in $\fO(\fU)\Contra$
(see Section~\ref{prelim-veryflat-and-contraadjusted-subsecn}).

 Consider the direct image $0\rarrow\fj_!\fj^!\fQ\rarrow\fj_!\fK
\rarrow\fj_!\fH\rarrow0$ of the short exact sequence $0\rarrow\fj^!\fQ
\rarrow\fK\rarrow\fH\rarrow0$ with respect to the affine open
immersion~$\fj$, and take its pushout with respect to the adjunction
morphism $\fj_!\fj^!\fQ\rarrow\fQ$ in $\fX\Lcth_\bW$.
 Let us show that, in the resulting short exact sequence
$0\rarrow\fQ\rarrow\fL\rarrow\fj_!\fH\rarrow0$, the $\bW$\+locally
contraherent cosheaf of contramodules $\fL$ on $\fX$ is locally
cotorsion in restriction to $\fU\cup\fV$.
 By Corollary~\ref{colocality-of-cotorsion-presuming-contraadjusted},
it suffices to show that the restrictions of $\fL$ to $\fU$ and $\fV$
are locally cotorsion.

 Indeed, in restriction to $\fU$ we have $\fj^!\fj_!\fj^!\fQ\simeq
\fj^!\fQ$, so the morphism $\fj_!\fj^!\fQ\rarrow\fQ$ becomes
an isomorphism after applying~$\fj^!$, hence $\fj^!\fL\simeq
\fj^!\fj_!\fK\simeq\fK$ is a locally cotorsion contraherent cosheaf
of contramodules.
 On the other hand, denote by $\fj'\:\fU\cap\fV\rarrow\fV$ and
$\fh'\:\fU\cap\fV\rarrow\fU$ the open immersions of $\fU\cap\fV$.
 Then we have a natural isomorphism $\fh^!\fj_!\fH\simeq\fj'_!
\fh'{}^!\fH$ of contraherent cosheaves of contramodules on~$\fV$
(by Lemma~\ref{cosheaf-arb-dir-openimmers-inv-basechange-lemma}).

 Furthermore, we have an isomorphism $\fh'{}^!\fj^!\fQ\simeq
\fj'{}^!\fh^!\fQ$ of contraherent cosheaves of contramodules
on $\fU\cap\fV$, since the morphisms of formal schemes $\fj\fh'$
and $\fh\fj'\:\fU\cap\fV\rarrow\fX$ coincide.
 Notice that the contraherent cosheaf of contramodules $\fj'{}^!
\fh^!\fQ$ is locally cotorsion, since the $\bW|_\fV$\+locally
contraherent cosheaf of contramodules $\fh^!\fQ$ on $\fV$ is locally
cotorsion by assumption.
 Hence the contraherent cosheaf of contramodules $\fh'{}^!\fH$
is locally cotorsion as the cokernel of the admissible monomorphism
of locally cotorsion contraherent cosheaves of contramodules
$\fh'{}^!\fj^!\fQ\rarrow\fh'{}^!\fK$ on $\fU\cap\fV$.

 Since the direct images of locally cotorsion $\bT$\+locally
contraherent cosheaves of contramodules with respect to
$(\bW,\bT)$\+affine morphisms are locally cotorsion (by
formula~\eqref{loc-ctrh-lct-contramod-W-T-affine-direct-image}),
the contraherent cosheaf of contramodules $\fj'_!\fh'{}^!\fH$
on $\fV$ is locally cotorsion, too.
 Now in the short exact sequence $0\rarrow\fh^!\fQ\rarrow\fh^!\fL
\rarrow\fh^!\fj_!\fH\rarrow0$ of $\bW|_\fV$\+locally contraherent
cosheaves of contramodules on $\fV$, the middle term is locally
cotorsion because so are the leftmost and the rightmost terms.

 Finally, the composition of the admissible monomorphisms of
$\bW$\+locally contraherent cosheaves of contramodules
$\fM\rarrow\fQ\rarrow\fL$ on $\fX$ is an admissible monomorphism
whose cokernel is an extension of the contraherent cosheaves
of contramodules $\fj_!\fH$ and $\fG$ on $\fX$, hence also a finitely
iterated extension of the direct images of antilocally flat contraherent
cosheaves of contramodules from the affine open formal subschemes
$\fU_\alpha\subset\fX$, \,$\alpha\le\beta$.
 This finishes the argument for the induction step, proving the lemma.
\end{proof}

 Let\/ $\fX$ be a semi-separated Noetherian formal scheme with an open
covering\/~$\bW$.
 We denote by $\Ext^{\fX,*}({-},{-})$ the $\Ext$ groups in the exact
category of $\bW$\+locally contraherent cosheaves of contramodules
on~$\fX$.
 Similarly to~\cite[Section~4.3]{Pcosh}, we notice that such $\Ext$
groups do not depend on the open covering $\bW$ and coincide with
the $\Ext$ groups computed in the whole category of locally
contraherent cosheaves of contramodules $\fX\Lcth$.
 Indeed, the full exact subcategory $\fX\Lcth_\bW$ is closed under
extensions and kernels of admissible epimorphisms in $\fX\Lcth$ (see
Section~\ref{locally-contraherent-cosheaves-of-contramods-subsecn}),
and for any object $\fP\in\fX\Lcth_\bW$ there exists an admissible
epimorphism onto $\fP$ from an object of $\fX\Ctrh\subset\fX\Lcth_\bW$
(see formula~\eqref{Cech-sequence-of-cosheaves-of-contramods}).
 So the exact categories $\fX\Ctrh\subset\fX\Lcth_\bW\subset\fX\Lcth$
are \emph{resolving} as full subcategories in each other, in
the sense of~\cite[Section~2]{Sto}, \cite[Section~A.3]{Pcosh},
or Section~\ref{appx-resolving-subsecn} below, and the result
of Proposition~\ref{infinite-co-resolutions-minus-plus-derived}(a)
below is applicable.

 For the same reasons (up to inverting the arrows), the $\Ext$ groups
computed in the exact subcategory of locally cotorsion $\bW$\+locally
contraherent cosheaves of contramodules $\fX\Lcth^\lct_\bW$ agree with
those in $\fX\Lcth_\bW$ (and also in $\fX\Lcth^\lct$).
 Indeed, the full subcategory $\fX\Lcth^\lct_\bW$ is closed under
extensions and cokernels of admissible monomorphisms in $\fX\Lcth_\bW$
(see Section~\ref{locally-contraherent-lct-cosheaves-subsecn}),
and Lemma~\ref{locally-cotorsion-preenvelope-lemma} provides
an admissible monomorphism from any $\bW$\+locally contraherent cosheaf
of contramodules on $\fX$ into a locally cotorsion $\bW$\+locally
contraherent cosheaf of contramodules.
 In other words, the full subcategory $\fX\Lcth^\lct_\bW$ is
\emph{coresolving} in $\fX\Lcth_\bW$.
 We refer to
Proposition~\ref{infinite-co-resolutions-minus-plus-derived}(b) below
and~\cite[Lemma~6.3]{Pfltp} for further details.

\begin{cor} \label{antilocally-flat-characterized-cor}
 Let\/ $\fX$ be a semi-separated Noetherian formal scheme with an open
covering\/~$\bW$. \par
\textup{(a)} A\/ $\bW$\+locally contraherent cosheaf of contramodules\/
$\fF$ on\/ $\fX$ is antilocally flat if and only if\/
$\Ext^{\fX,1}(\fF,\fP)=0$ for all locally cotorsion\/ $\bW$\+locally
contraherent cosheaves of contramodules\/ $\fP$ on\/ $\fX$, and if
and only if\/ $\Ext^{\fX,n}(\fF,\fP)=0$ for all locally cotorsion\/
$\bW$\+locally contraherent cosheaves of contramodules\/ $\fP$ on\/
$\fX$ and all $n\ge1$. \par
\textup{(b)} The class of antilocally flat\/ $\bW$\+locally contraherent
cosheaves of contramodules on\/ $\fX$ is closed under extensions and
kernels of admissible epimorphisms in the exact category\/
$\fX\Lcth_\bW$.
\end{cor}

\begin{proof}
 This is our version of~\cite[Corollary~4.4.2]{Pcosh}.
 Part~(a) follows from the existence of an admissible monomorphism from
any $\bW$\+locally contraherent cosheaf of contramodules on $\fX$ to
a locally cotorsion $\bW$\+locally contraherent cosheaf of contramodules
(a weak form of Lemma~\ref{locally-cotorsion-preenvelope-lemma})
by virtue of Lemma~\ref{Hom-from-to-short-exact-sequences-lemma}(b)
below.
 Part~(b) follows from part~(a).
\end{proof}

\begin{lem} \label{antilocally-flat-precover-lemma}
 Let\/ $\fX$ be a semi-separated Noetherian formal scheme with an open
covering\/ $\bW$ and a finite affine open covering\/
$\fX=\bigcup_{\alpha=1}^N\fU_\alpha$ subordinate to\/~$\bW$.
 Then any\/ $\bW$\+locally contraherent cosheaf of contramodules\/
$\fM$ on\/ $\fX$ can be included in an (admissible) short exact
sequence\/ $0\rarrow\fP\rarrow\fF\rarrow\fM\rarrow0$ in\/
the exact category $\fX\Lcth_\bW$, where\/ $\fP$ is a locally
cotorsion\/ $\bW$\+locally contraherent cosheaf of contramodules on\/
$\fX$ and\/ $\fF$ is a finitely iterated extension of the direct images
of contraherent cosheaves of contramodules on\/ $\fU_\alpha$
corresponding to flat contraadjusted contramodule\/
$\fO_\fX(\fU_\alpha)$\+modules.
\end{lem}

\begin{proof}
 This is our version of~\cite[Lemma~4.4.3]{Pcosh}.
 The assertion follows from
Lemmas~\ref{locally-cotorsion-preenvelope-lemma}
and~\ref{projective-contrah-on-qcomp-ssep} by virtue of an argument
known as the \emph{Salce lemma}~\cite{Sal},
\cite[Lemma~B.1.1(a)]{Pcosh}.
 See Section~\ref{appx-cotorsion-pairs-subsecn} below
for a brief discussion.

 Specifically, by Lemma~\ref{projective-contrah-on-qcomp-ssep},
there exist projective contraherent cosheaves of contramodules
$\fG_\alpha$ on $\fU_\alpha$ together with an admissible epimorphism
$\bigoplus_\alpha\fj_\alpha{}_!\fG_\alpha\rarrow\fM$ in $\fX\Lcth_\bW$.
 Denote by $\fK$ the kernel of this admissible epimorphism.
 By Lemma~\ref{locally-cotorsion-preenvelope-lemma}, there exists
a short exact sequence $0\rarrow\fK\rarrow\fP\rarrow\fH\rarrow0$
in $\fX\Lcth_\bW$ with a locally cotorsion $\bW$\+locally contrahernt
cosheaf of contramodules $\fP$ and a cosheaf $\fH$ that is a finitely
iterated extension of the direct images of antilocally flat
contraherent cosheaves of contramodules from~$\fU_\alpha$.

 Denote by $\fF$ the pushout of the two admissible monomorphisms
$\fK\rarrow\bigoplus_\alpha j_\alpha{}_!\fG_\alpha$ and
$\fK\rarrow\fP$ in $\fX\Lcth_\bW$.
 Then we have short exact sequences $0\rarrow\fP\rarrow\fF\rarrow
\fM\rarrow0$ and $0\rarrow\bigoplus_\alpha\fj_\alpha{}_!\fG_\alpha
\rarrow\fF\rarrow\fH\rarrow0$ in $\fX\Lcth_\bW$.
 The latter sequence shows that $\fF$ is a finitely iterated extension
of the direct images of antilocally flat contraherent cosheaves of
contramodules from~$\fU_\alpha$, as desired.
\end{proof}

\begin{cor} \label{antilocally-flat-corollary}
 Let\/ $\fX$ be a semi-separated Noetherian formal scheme with an open
covering\/ $\bW$ and a finite affine open covering\/
$\fX=\bigcup_{\alpha=1}^N\fU_\alpha$ subordinate to\/~$\bW$. \par
\textup{(a)} For any\/ $\bW$\+locally contraherent cosheaf of
contramodules\/ $\fM$ on\/ $\fX$ there exists an admissible monomorphism
from\/ $\fM$ into a locally cotorsion\/ $\bW$\+locally contraherent
cosheaf of contramodules\/ $\fP$ on\/ $\fX$ such that the cokernel\/
$\fF$ is an antilocally flat\/ $\bW$\+locally contraherent cosheaf of
contramodules. \par
\textup{(b)} For any\/ $\bW$\+locally contraherent cosheaf of
contramodules\/ $\fM$ on\/ $\fX$ there exists an admissible epimorphism
onto\/ $\fM$ from an antilocally flat\/ $\bW$\+locally contraherent
cosheaf of contramodules on\/ $\fX$ such that the kernel\/ $\fP$ is
a locally cotorsion\/ $\bW$\+locally contraherent cosheaf of 
contramodules. \par
\textup{(c)} A\/ $\bW$\+locally contraherent cosheaf of contramodules
on\/ $\fX$ is antilocally flat if and only if it is (a contraherent
cosheaf of contramodules and) a direct summand of a finitely iterated
extension of the direct images of contraherent cosheaves of
contramodules on\/ $\fU_\alpha$ corresponding to flat contraadjusted
contramodule\/ $\fO(\fU_\alpha)$\+modules.
\end{cor}

\begin{proof}
 This is our version of~\cite[Corollary~4.4.4]{Pcosh}.
 The ``if'' implication in part~(c) follows from
Corollary~\ref{antilocally-flat-characterized-cor}(b) together with
the remarks in the beginning of this section.
 After the ``if'' assertion of part~(c) is shown, part~(a) follows
from Lemma~\ref{locally-cotorsion-preenvelope-lemma} and part~(b)
follows from Lemma~\ref{antilocally-flat-precover-lemma}.

 The ``only if'' implication in part~(c) follows from
Lemma~\ref{antilocally-flat-precover-lemma} and
Corollary~\ref{antilocally-flat-characterized-cor}(a) by virtue of
the direct summand lemma (essentially, Lemma~\ref{direct-summand-lemma}
below).
 Specifically, let $\fF$ be an antilocally flat $\bW$\+locally
contraherent cosheaf of contramodules on~$\fX$.
 By Lemma~\ref{antilocally-flat-precover-lemma}, there exists a short
exact sequence $0\rarrow\fP\rarrow\fG\rarrow\fF\rarrow0$ in
$\fX\Lcth_\bW$ such that $\fP$ is a locally cotorsion $\bW$\+locally
contraherent cosheaf of contramodules on $\fX$, while $\fG$ is
a direct summand of a finitely iterated extension of the direct images
of antilocally flat contraherent cosheaves of contramodules
from~$\fU_\alpha$.
 By Corollary~\ref{antilocally-flat-characterized-cor}(a), we have
$\Ext^{\fX,1}(\fF,\fP)=0$, hence $\fF$ is a direct summand of~$\fG$.

 It remains to point out that the direct images with respect to
the open immersions $\fU_\alpha\rarrow\fX$ take contraherent cosheaves
of contramodules to (globally) contraherent cosheaves of contramodules
(by formula~\eqref{loc-ctrh-contramod-W-T-affine-direct-image}), and
the full subcategory of contraherent cosheaves of contramodules
$\fX\Ctrh$ is closed under extensions in $\fX\Lcth_\bW$ (as per
Section~\ref{locally-contraherent-cosheaves-of-contramods-subsecn}).
\end{proof}

 The result of Corollary~\ref{antilocally-flat-corollary}(c) can be
rephrased by saying that the class of antilocally flat contraherent
cosheaves of contramodules on semi-separated Noetherian formal schemes
is \emph{antilocal} in the sense similar to~\cite[Section~4]{Pal}.
 See Remark~\ref{flat-and-antilocally-flat-remark} below for
a further discussion.

\begin{cor} \label{antilocally-flat-independent-on-covering}
 Let\/ $\fX$ be a semi-separated Noetherian formal scheme with an open
covering\/~$\bW$.
 Then all antilocally flat\/ $\bW$\+locally contraherent cosheaves of
contramodules on\/ $\fX$ are (globally) contraherent.
 The full subcategory of antilocally flat\/ $\bW$\+locally
contraherent cosheaves of contramodules in the exact category\/
$\fX\Lcth$ does not depend on the open covering\/~$\bW$.
\end{cor}

\begin{proof}
 This is our version of~\cite[Corollary~4.4.5]{Pcosh}.
 The first assertion is explicitly a part of
Corollary~\ref{antilocally-flat-corollary}(c).
 To prove the second one, notice that, for any two open coverings\/
$\bW'$ and\/ $\bW''$ of the formal scheme $\fX$, there exists a finite
affine open covering $\fX=\bigcup_{\alpha=1}^N\fU_\alpha$ subordinate
to both $\bW'$ and~$\bW''$.
 Then it remains to apply Corollary~\ref{antilocally-flat-corollary}(c).
\end{proof}

\begin{cor} \label{locally-cotorsion-characterized-cor}
 Let\/ $\fX$ be a semi-separated Noetherian formal scheme with an open
covering\/~$\bW$.
 Then, for any\/ $\bW$\+locally contraherent cosheaf of contramodules\/
$\fP$ on\/ $\fX$, the following conditions are equivalent:
\begin{enumerate}
\item the functor\/ $\Hom^\fX({-},\fP)$ takes short exact sequences of
antilocally flat contraherent cosheaves of contramodules on\/ $\fX$ to
short exact sequences of abelian groups;
\item one has\/ $\Ext^{\fX,1}(\fF,\fP)=0$ for all antilocally
flat contraherent cosheaves of contramodules\/ $\fF$ on\/~$\fX$;
\item one has\/ $\Ext^{\fX,n}(\fF,\fP)=0$ for all antilocally
flat contraherent cosheaves of contramodules\/ $\fF$ on\/ $\fX$ and
all integers $n\ge1$;
\item $\fP$~is a locally cotorsion\/ $\bW$\+locally contraherent
cosheaf of contramodules on\/~$\fX$.
\end{enumerate}
\end{cor}

\begin{proof}
 The argument is dual to the proofs of
Corollaries~\ref{very-flat-contramodules-cor}
and~\ref{flat-contramodules-cor}.
 Let us spell out some specifics.

 (1)~$\Longleftrightarrow$~(2)~$\Longleftrightarrow$~(3)
 This is an application of
Lemma~\ref{Hom-from-to-short-exact-sequences-lemma}(a) below.
 The assertions follow from the facts that there are enough antilocally
flat contraherent cosheaves of contramodules in $\fX\Lcth_\bW$
(a weak version of Lemma~\ref{projective-contrah-on-qcomp-ssep}
or Corollary~\ref{antilocally-flat-corollary}(b)) and the class of
antilocally flat contraherent cosheaves of contramodules is closed
under kernels of admissible epimorphisms
(Corollary~\ref{antilocally-flat-characterized-cor}(b)).

 (4)~$\Longrightarrow$~(3) holds by
Corollary~\ref{antilocally-flat-characterized-cor}(a).

 (2)~$\Longrightarrow$~(4) follows from
Corollary~\ref{antilocally-flat-corollary}(a) by virtue of
the argument proving Lemma~\ref{direct-summand-lemma}.
\end{proof}

 In the terminology of Section~\ref{appx-cotorsion-pairs-subsecn} below,
the results of Corollaries~\ref{antilocally-flat-corollary}(a\+-b),
\ref{antilocally-flat-characterized-cor},
and~\ref{locally-cotorsion-characterized-cor}
can be summarized by the following formulation.
 For any semi-separated Noetherian formal scheme $\fX$ with an open
covering $\bW$, the pair of classes (antilocally flat contraherent
cosheaves of contramodules on~$\fX$, locally cotorsion $\bW$\+locally
contraherent cosheaves of contramodules on~$\fX$) is a \emph{hereditary
complete cotorsion pair} in the exact category of $\bW$\+locally
contraherent cosheaves of contramodules $\fX\Lcth_\bW$.

 We will denote the category of antilocally flat contraherent cosheaves
of contramodules on $\fX$ by $\fX\Ctrh_\alf$.
 As a full subcategory closed under extensions and the kernels of
admissible epimorphisms in $\fX\Ctrh$, \,$\fX\Lcth_\bW$, or $\fX\Lcth$,
the category $\fX\Ctrh_\alf$ inherits an exact category structure.

 Clearly, one has $\fX\Ctrh^\lct_\prj=\fX\Ctrh_\alf\cap
\fX\Ctrh^\lct=\fX\Ctrh_\alf\cap\fX\Lcth^\lct_\bW$.
 There are \emph{both} enough projective objects and enough injective
objects in the exact category $\fX\Ctrh_\alf$.
 The full subcategory of projective objects in $\fX\Ctrh_\alf$ is
$\fX\Ctrh_\prj$, while the full subcategory of injective objects in
$\fX\Ctrh_\alf$ is $\fX\Ctrh^\lct_\prj$.
 Here the assertions concerning projectives in $\fX\Ctrh_\alf$ follow
from Lemma~\ref{projective-contrah-on-qcomp-ssep} and
Corollary~\ref{antilocally-flat-characterized-cor}(b), while
the assertions concerning injectives follow from
Corollaries~\ref{antilocally-flat-corollary}(a)
and~\ref{antilocally-flat-characterized-cor}(b)
(cf.\ Lemmas~\ref{cotorsion-pair-kernel-as-injectives-projectives}
and~\ref{cotorsion-pair-the-same-projectives-injectives}(a) below).

\begin{cor} \label{antilocally-flat-products}
 Let\/ $\fX$ be a semi-separated Noetherian formal scheme with an open
covering\/~$\bW$.
 Then the full subcategory\/ $\fX\Ctrh_\alf$ is closed under infinite
direct products in\/ $\fX\Ctrh$ and\/ $\fX\Lcth_\bW$.
\end{cor}

\begin{proof}
 This is our version of~\cite[Corollary~4.4.8]{Pcosh}.
 The assertion follows from
Corollary~\ref{antilocally-flat-corollary}(c) together with the facts
that the infinite products in the exact category $\fX\Lcth_\bW$ are
exact, the functors of direct image of contraherent cosheaves of
contramodules with respect to affine morphisms of formal schemes
preserve infinite products, and the class of antilocally flat
contraherent cosheaves of contramodules on an affine Noetherian
formal scheme $\fU$ is closed under infinite products.
 The latter property holds because the class of flat contramodule
modules over an adically coherent adic topological ring $R$ is
preserved by infinite products in $R\Contra$
(see Lemma~\ref{sums-of-inj-tors-prods-of-flat-contramods}(b)).
\end{proof}

\begin{cor} \label{antilocally-flat-inverse-image-affine-open-immers}
 Let\/ $\fX$ be a semi-separated Noetherian formal scheme and\/
$\fY\subset\fX$ be an open formal subscheme such that the open
immersion morphism\/ $\fj\:\fY\rarrow\fX$ is affine.
 Then the inverse image functor\/~$\fj^!$ takes antilocally flat
contraherent cosheaves of contramodules on\/ $\fX$ to antilocally flat
contraherent cosheaves of contramodules on\/~$\fY$.
\end{cor}

\begin{proof}
 This is our version of~\cite[Corollary~4.4.9]{Pcosh}.
 The point is that any finite affine open covering $\{\fU_\alpha\}$
of $\fX$ restricts to a finite affine open covering
$\{\fY\cap\fU_\alpha\}$ of~$\fY$.
 So one can use Corollary~\ref{antilocally-flat-corollary}(c) together
with Lemma~\ref{cosheaf-arb-dir-openimmers-inv-basechange-lemma}
in order to reduce the question to the case of an affine Noetherian
formal scheme $\fX$ and its affine open formal subscheme~$\fY$.
 In the latter case, one needs to use
Corollary~\ref{colocalization-of-flat-contramodules-cor}.
\end{proof}

 Before we finish this section, let us have a brief introductory
discussion of \emph{flat cosheaves of contramodules}.

 Let $\fX$ be a locally Noetherian formal scheme with an open
covering~$\bW$.
 A cosheaf of contramodule $\fO_\fX$\+modules $\fF$ on $\fX$ is said
to be \emph{$\bW$\+flat} (cf.~\cite[Section~3.7]{Pcosh}) if
the contramodule $\fO_\fX(\fU)$\+module $\fF[\fU]$ is flat for
every affine open subscheme $\fU\subset\fX$ subordinate to~$\bW$.
 A cosheaf of contramodule $\fO_\fX$\+modules $\fF$ is said to be
\emph{flat} if it is $\bW_\fX$\+flat for the trivial open covering
$\bW_\fX=\{\fX\}$.
 It is clear from
Lemma~\ref{flat-map-of-adic-rings-direct-image-lemma}(c) or~(d) that
the direct image of a $\bT$\+flat cosheaf of contramodule
$\fO_\fY$\+modules with respect to a flat tight $(\bW,\bT)$\+affine
morphism of formal schemes $\ff\:\fY\rarrow\fX$ is $\bW$\+flat.
 By Corollary~\ref{colocalization-of-flat-contramodules-cor},
a contraherent cosheaf of contramodules $\fF$ on an affine Noetherian
formal scheme $\fU$ is flat if and only if the contraadjusted
contramodule $\fO_\fU(\fU)$\+module $\fF[\fU]$ is flat.

 We will denote the full subcategory of $\bW$\+flat $\bW$\+locally
contraherent cosheaves of contramodules on $\fX$ by
$\fX\Lcth_\bW^\fl\subset\fX\Lcth_\bW$, and the full subcategory of
flat contraherent cosheaves of contramodules on $\fX$ by
$\fX\Ctrh^\fl\subset\fX\Ctrh$.
 The full subcategory $\fX\Lcth_\bW^\fl$ is closed under extensions,
kernels of admissible epimorphisms, and infinite products in
$\fX\Lcth_\bW$ (cf.\ the second paragraph of the proof of
Corollary~\ref{antilocally-flat-products}).
 So we have the exact category structure on $\fX\Lcth_\bW^\fl$
inherited from $\fX\Lcth_\bW$.

 The next corollary is our version of~\cite[Corollary~4.4.7]{Pcosh}.

\begin{cor} \label{antilocally-flat-are-flat-cor}
 Any antilocally flat contraherent cosheaf of contramodules on
a semi-separated Noetherian formal scheme is flat.
\end{cor}

\begin{proof}
 Follows from Corollary~\ref{antilocally-flat-corollary}(c) in view
of the remarks in the paragraphs preceding the present corollary.
\end{proof}

\subsection{Homology of contraherent cosheaves of contramodules~I}
\label{homology-of-cosheaves-I-subsecn}
 As usual, we denote by $\Gamma(\fU,\cM)=\cM(\fU)$ the group/module
of sections of a sheaf $\cM$ on a topological space/formal scheme
$\fX$ over an open subset $\fU\subset\fX$.
 Similarly, following the notation in~\cite[Sections~3.2, 3.4,
and~4.6]{Pcosh}, we denote by $\Delta(\fU,\fP)=\fP[\fU]$
the group/module of cosections of a cosheaf $\fP$ over~$\fU$.

 The functor $\Gamma(\fU,{-})$ is left exact on the abelian category
of sheaves of abelian groups on a topological space $\fX$, hence also
on the abelian category $\fX\Tors$ of quasi-coherent torsion sheaves
on a locally Noetherian formal scheme $\fX$, for any open subset/formal
subscheme $\fU\subset\fX$.
 Dual-analogously, the functor $\Delta(\fU,{-})$ is right exact on
the exact category $\fX\Lcth$ of locally contraherent cosheaves of
contramodules on a locally Noetherian formal scheme $\fX$, for any
open formal subscheme $\fU\subset\fX$.
 The latter assertion means that, for any (admissible) short exact
sequence $0\rarrow\fK\rarrow\fL\rarrow\fM\rarrow0$ in $\fX\Lcth$,
the sequence of abelian groups/$\fO_\fX(\fU)$\+modules
$$
 \Delta(\fU,\fK)\lrarrow\Delta(\fU,\fL)\lrarrow\Delta(\fU,\fM)
 \lrarrow0
$$
is right exact.
 Indeed, the formula~\eqref{topology-base-cosheaf-axiom} computing
the group $\fP[\fU]$ for a $\bW$\+locally contraherent cosheaf of
contramodules $\fP$ on $\fX$ in terms of the groups $\fP[\fW]$ for
affine open formal subschemes $\fW\subset\fX$ subordinate to $\bW$
expresses the group $\fP[\fU]$ as the cokernel of a morphism of
(possibly infinite) direct sums of the groups $\fP[\fW]$; so this
is a right exact functor.

 Let $\fX$ be a semi-separated Noetherian formal scheme.
 Then the left derived functor $\boL_*\Delta(\fX,{-})\:\fX\Lcth\rarrow
\fO_\fX(\fX)\Modl$ of the functor of global cosections $\Delta(\fX,{-})$
is defined as usual projective resolutions in the exact category
$\fX\Lcth$ (which exist by
Corollary~\ref{projective-contrah-on-qcomp-ssep-corollary}).
 It follows from the same corollary that the derived functors of
$\Delta(\fX,{-})$ computed in the exact category $\fX\Lcth_\bW$ for
a fixed open covering $\bW$ and in the whole exact category $\fX\Lcth$
agree.

 Similarly to~\cite[Section~4.6]{Pcosh}, the groups/modules
$\boL_*\Delta(\fX,\fM)$ are called the \emph{homology groups/modules}
of a locally contraherent cosheaf of contramodules $\fM$ on
a semi-separated Noetherian formal scheme~$\fX$.

 In the special case of an affine formal scheme $\fU$ and the trivial
open covering $\bW_\fX=\{\fX\}$, the functor $\Delta(\fU,{-})\:
\fU\Ctrh\rarrow\fO_\fU(\fU)\Modl$ is exact.
 Hence the higher homology modules $\boL_n\Delta(\fU,\fN)$ vanish
for all (globally) contraherent cosheaves of contramodules $\fN$
on $\fU$ and all $n\ge1$.

 For a very flat $(\bW,\bT)$\+affine morphism of semi-separated
Noetherian formal schemes $\ff\:\fY\rarrow\fX$, we have an exact
functor of direct image $\ff_!\:\fY\Lcth_\bT\rarrow\fX\Lcth_\bW$
\,\eqref{loc-ctrh-contramod-W-T-affine-direct-image}.
 By the definition of~$\ff_!$, there is a natural isomorphism
of $\fO_\fX(\fX)$\+mod\-ules $\Delta(\fX,\ff_!\fN)\simeq\Delta(\fY,\fN)$
for any $\bT$\+locally contraherent cosheaf of contramodules $\fN$
on~$\fY$.
 Furthermore, by
Corollaries~\ref{projective-contrah-on-qcomp-ssep-corollary}
and~\ref{ssep-proj-contrah-direct-image-cor}(a), the functor~$\ff_!$
takes projective objects of $\fY\Lcth_\bT$ to projective objects
of $\fX\Lcth_\bW$.
 It follows that, for any $\bT$\+locally contraherent cosheaf of
contramodules $\fN$ of $\fY$, there are natural isomorphisms of
$\fO_\fX(\fX)$\+modules
\begin{equation} \label{veryflat-affine-dirim-cosheaf-homology-iso}
 \boL_n\Delta(\fX,\ff_!\fN)\simeq\boL_n\Delta(\fY,\fN),
 \qquad n\ge0.
\end{equation}
 As a particular case of the latter assertion, one obtains the higher
homology vanishing $\boL_n\Delta(\fX,\fj_!\fN)=0$ for any affine open
formal subscheme $\fU\subset\fX$ with the open immersion morphism
$\fj\:\fU\rarrow\fX$, any contraherent cosheaf of contramodules
$\fN$ on~$\fU$, and all $n\ge1$.

 The full subcategory in $\fX\Lcth$ formed by all the cosheaves $\fM$
such that $\boL_n\Delta(\fX,\fM)=0$ for all $n\ge1$ is obviously
closed under extensions and direct summands.
 Thus it follows from the previous paragraph together with
Corollary~\ref{antilocally-flat-corollary}(c) that one has
$\boL_n\Delta(\fX,\fF)=0$ for all antilocally flat contraherent
cosheaves $\fF$ on $\fX$ and all $n\ge1$.
 So one can compute the derived functor $\boL_*\Delta(\fX,{-})$ using
antilocally flat resolutions in $\fX\Lcth$.

 In particular, it follows that the homology of locally cotorsion
locally contraherent cosheaves of contramodules on $\fX$ can be
computed using projective locally cotorsion resolutions.
 Hence the derived functors $\boL_*\Delta(\fX,{-})$ computed in
the exact categories $\fX\Lcth_\bW$ and $\fX\Lcth_\bW^\lct$ agree,
and so do the exact functors $\boL_*\Delta(\fX,{-})$ computed in
the exact categories $\fX\Lcth$ and $\fX\Lcth^\lct$.

\begin{cor} \label{homological-criterion-of-contraherence-cor}
 Let\/ $\fU$ be an affine Noetherian formal scheme and\/ $\fP$ be
a locally contraherent cosheaf of contramodules on\/~$\fU$.
 Then the following three conditions are equivalent:
\begin{enumerate}
\item $\boL_1\Delta(\fU,\fP)=0$;
\item $\boL_n\Delta(\fU,\fP)=0$ for all $n\ge1$;
\item $\fP$ is a (globally) contraherent cosheaf of contramodules
on\/~$\fU$.
\end{enumerate}
\end{cor}

\begin{proof}
 This is our version of~\cite[Corollary~4.6.1]{Pcosh}.
 Let $\bW$ be an open covering of $\fU$ such that the cosheaf of
contramodules $\fP$ is $\bW$\+locally contraherent, and let
$\fU=\bigcup_\alpha\fW_\alpha$ be a finite affine open covering of
$\fU$ subordinate to~$\bW$.
 Then the \v Cech exact
sequence~\eqref{Cech-sequence-of-cosheaves-of-contramods} is 
a resolution of $\fP$ by finite direct sums of the direct images
of contraherent cosheaves of contramodules from affine open formal
subschemes of~$\fU$.
 So all the terms of~\eqref{Cech-sequence-of-cosheaves-of-contramods},
except perhaps the rightmost one, are annihilated by the higher
derived cosections functors $\boL\Delta_n(\fU,{-})$, \,$n\ge1$,
according to the discussion above.
 On the other hand, applying the functor $\Delta(\fU,{-})$ to
the truncated form of
the resolution~\eqref{Cech-sequence-of-cosheaves-of-contramods} produces
the \v Cech complex of abelian groups $C_\bu(\{\fW_\alpha\},\fP)$
\,\eqref{Cech-complex-of-abelian-groups-for-copresheaf}.
 Thus the assertion of the corollary is a restatement of
Proposition~\ref{homological-criterion-of-contraherence-prop}.
\end{proof}

\begin{cor} \label{ssep-flat-direct-image-of-antilocally-flat-cor}
 For any flat tight morphism of semi-separated Noetherian formal
schemes\/ $\ff\:\fY\rarrow\fX$, the functor of direct image of
cosheaves of\/ $\fO$\+modules\/ $\ff_!\:(\fY,\fO_\fY)\Cosh\rarrow
(\fX,\fO_\fX)\Cosh$ takes antilocally flat contraherent cosheaves of
contramodules on\/ $\fY$ to antilocally flat contraherent cosheaves of
contramodules on\/ $\fX$, and induces an exact functor\/
$\ff_!\:\fY\Ctrh_\alf\rarrow\fX\Ctrh_\alf$.
\end{cor}

\begin{proof}
 This is our version of~\cite[Corollary~4.6.3(c)]{Pcosh};
see also~\cite[Proposition~6.15]{Pns}.
 Given an affine open formal subscheme $\fU\subset\fX$, consider
the open formal subscheme $\fV=\ff^{-1}(\fU)\subset\fY$.
 Then the open immersion morphism $\fk\:\fV\rarrow\fY$ is affine
(since so is the open immersion morphism $\fj\:\fU\rarrow\fX$).
 By Corollary~\ref{antilocally-flat-inverse-image-affine-open-immers},
for any antilocally flat contraherent cosheaf of contramodules $\fG$
on $\fY$, the contraherent cosheaf of contramodules $\fk^!\fG$ on
$\fV$ is antilocally flat.
 According to the discussion above, it follows that
$\boL_n\Delta(\fV,\fk^!\fG)=0$ for all $n\ge1$.
 Hence the functor $\fG\longmapsto\fG[\fV]=(\ff_!\fG)[\fU]$ is
exact on the exact category $\fY\Ctrh_\alf$.

 By Corollary~\ref{antilocally-flat-corollary}(c), \,$\fG$~is
a direct summand of a finitely iterated extension of the direct images
of antilocally flat contraherent cosheaves of contramodules from affine
open formal subschemes of~$\fY$.
 Let $\fT\subset\fY$ be an affine open formal subscheme with
the open immersion morphism $\fh\:\fT\rarrow\fY$.
 Then the composition $\ff\circ\fh\:\fT\rarrow\fX$ is a flat affine
morphism of formal schemes.
 One can choose open coverings $\bW$ of the formal scheme $\fX$ and
$\bT$ of the formal scheme $\fT$ such that the morphism $\ff\circ\fh$
is $(\bW,\bT)$\+affine and $(\bW,\bT)$\+coaffine.
 Following the discussion in the beginning of
Section~\ref{antilocally-flat-subsecn}, the functor
$\ff_!\circ\fh_!=(\ff\circ\fh)_!\:\fT\Lcth_\bT\rarrow \fX\Lcth_\bW$
takes antilocally flat contraherent cosheaves of contramodules on $\fT$
to antilocally flat contraherent cosheaves of contramodules on~$\fX$.

 Passing to a finitely iterated extension and arguing by induction,
one can now easily check both the contraadjustedness axiom~(v)
and the contraherence axiom~(iv) from
Section~\ref{contraherent-cosheaves-of-contramods-subsecn} for
the cosheaf of $\fO_\fX$\+modules $\ff_!\fG$ on~$\fX$
(while the contramoduleness axiom~(iii) always holds by
Lemma~\ref{contramodule-cosheaves-direct-images-lemma}).
 The property of Corollary~\ref{very-flat-contramodules-cor}(1) is
relevant here for verifying the contraherence axiom.
 So we have $\ff_!\fG\in\fX\Ctrh$.
 Finally, the class of antilocally flat contraherent cosheaves of
contramodules on $\fX$ is closed under extensions in $\fX\Ctrh$,
so it follows that the contraherent cosheaf of contramodules
$\ff_!\fG$ on $\fX$ is antilocally flat.

 For a more general result, see
Corollary~\ref{codilute-direct-images}(a) below.
\end{proof}

\subsection{Coflasque contraherent cosheaves of contramodules}
\label{coflasque-subsecn}
 Let $\fX$ be a topological space.
 There is a classical notion of a \emph{flasque sheaf of abelian groups}
on~$\fX$ \,\cite[Section~II.3.1]{God}.
 Dual-analogously, a cosheaf of abelian groups $\fE$ on $\fX$ is said
to be \emph{coflasque} if the corestriction maps $\fE[\fV]\rarrow
\fE[\fU]$ are injective for all open subsets $\fV\subset\fU\subset\fX$
\,\cite[Section~3.4]{Pcosh}.
 By~\cite[Section~II.3.1]{God} and~\cite[Lemma~3.4.1(a)]{Pcosh},
the property of a (co)sheaf of abelian groups on $\fX$ to be (co)flasque
is local, i.~e., it can be checked in restriction to any open covering
of~$\fX$.

\begin{lem} \label{coflasque-lcth-are-globally-contraherent}
 Let\/ $\fX$ be a locally Noetherian formal scheme with an open
covering\/ $\bW$ and\/ $\fE$ be a\/ $\bW$\+locally contraherent cosheaf
of contramodules on\/~$\fX$.
 Suppose that the cosheaf\/ $\fE$ is coflasque.
 Then\/ $\fE$ is a (globally) contraherent cosheaf of contramodules
on\/~$\fX$.
\end{lem}

\begin{proof}
 This is our version of~\cite[Corollary~3.4.2]{Pcosh}.
 The assertion follows from the homological criterion of contraherence
(Proposition~\ref{homological-criterion-of-contraherence-prop}) and
the fact that the higher \v Cech homology of a coflasque cosheaf with
respect to any open covering vanish~\cite[Lemma~3.4.1(b)]{Pcosh}.
\end{proof}

 Given a locally Noetherian formal scheme $\fX$, we will denote by
$\fX\Tors^\fq\subset\fX\Tors$ the full subcategory of flasque
quasi-coherent torsion sheaves and by $\fX\Ctrh_\cfq\subset\fX\Lcth_\bW$
the full subcategory of coflasque contraherent cosheaves of
contramodules on~$\fX$.
 We also put $\fX\Ctrh^\lct_\cfq=\fX\Ctrh_\cfq\cap\fX\Lcth_\bW^\lct
\subset\fX\Lcth_\bW$; so $\fX\Ctrh^\lct_\cfq$ is the full subcategory
of coflasque locally cotorsion contraherent cosheaves of contramodules.
 The following lemma summarizes the dual-analogous versions of
some well-known properties of flasque sheaves.

\begin{lem} \label{coflasque-adjusted-to-cosections}
 Let\/ $\fX$ be a locally Noetherian formal scheme with an open
covering\/~$\bW$.
 Then \par
\textup{(a)} the full subcategory\/ $\fX\Ctrh_\cfq$ is closed under
extensions in\/ $\fX\Lcth_\bW$; \par
\textup{(b)} the full subcategory\/ $\fX\Ctrh_\cfq$ is closed under
kernels of admissible epimorphisms in\/ $\fX\Lcth_\bW$; \par
\textup{(c)} for any (admissible) short exact sequence\/ $0\rarrow\fM
\rarrow\fN\rarrow\fE\rarrow0$ in\/ $\fX\Lcth_\bW$ with\/
$\fE\in\fX\Ctrh_\cfq$ and any open formal subscheme\/ $\fY\subset\fX$,
the short sequence of abelian groups\/ $0\rarrow\fM[\fY]\rarrow
\fN[\fY]\rarrow\fE[\fY]\rarrow0$ is exact.
\end{lem}

\begin{proof}
 Similar to~\cite[Corollary~3.4.4]{Pcosh}.
\end{proof}

 The full subcategory $\fX\Ctrh_\cfq$ is also closed under infinite
products in $\fX\Lcth_\bW$, while the full subcategory $\fX\Tors^\fq$
is closed under infinite direct sums in $\fX\Tors$.
 Indeed, it suffices to check the (co)flasqueness in restriction to
affine or quasi-compact open formal subschemes of $\fX$, and all open
formal subschemes of Noetherian formal schemes are Noetherian
(i.~e., quasi-compact).
 Cf.~\cite[the paragraphs preceding Example~3.4.3]{Pcosh}.

 The following lemma is our version of~\cite[Lemma~3.4.6(c,d)]{Pcosh}.
 For a partial generalization, see
Lemma~\ref{fHom-from-flasque-contratensor-with-alf-co-flasque} below.

\begin{lem} \label{fHom-contratensor-co-flasque}
 Let\/ $\fU$ be an affine Noetherian formal scheme.
 In this setting: \par
\textup{(a)} Let\/ $\cE$ be a flasque quasi-coherent torsion sheaf and\/
$\cJ$ be an injective quasi-coherent torsion sheaf on\/~$\fU$.
 Then the locally cotorsion contraherent cosheaf of contramodules\/
$\fHom_\fU(\cE,\cJ)$ on\/ $\fU$ is coflasque. \par
\textup{(b)} Let\/ $\cE$ be a flasque quasi-coherent torsion sheaf and\/
$\fG$ be a flat cosheaf of contramodule\/ $\fO_\fU$\+modules on\/~$\fU$.
 Then the quasi-coherent torsion sheaf\/ $\cE\ocn_\fU\fG$ on\/ $\fU$ is
flasque.
\end{lem}

\begin{proof}
 Part~(b): notice that $\cE\ocn_\fU\fG$ is a quasi-coherent torsion
sheaf on $\fU$ corresponding to the torsion $\fO_\fU(\fU)$\+module
$(\cE\ocn_\fU\fG)(\fU)=\cE(\fU)\ot_{\fO_\fU(\fU)}\fG[\fU]$.
 Furthermore, $\fG[\fU]$ is a flat $\fO_\fU(\fU)$\+module, and all open
formal subschemes in $\fU$ are quasi-compact.
 The key observation is that the rule $\fV\longmapsto
\cE(\fV)\ot_{\fO_\fU(\fU)}\fG[\fU]$ for all (not necessarily affine!)
open formal subschemes $\fV\subset\fU$ defines a quasi-coherent
torsion sheaf on~$\fU$.
 This quasi-coherent torsion sheaf is clearly flasque and isomorphic to
$\cE\ocn_\fU\fG$.

 Part~(a): notice that $\fHom_\fU(\cE,\cJ)$ is a contraherent cosheaf
of contramodules on $\fU$ corresponding to the contramodule
$\fO_\fU(\fU)$\+module $\fHom_\fU(\cE,\cJ)[\fU]=
\Hom_{\fO_\fU(\fU)}(\cE(\fU),\allowbreak\cJ(\fU))$.
 Furthemore, $\cJ(\fU)$ is an injective $\fO_\fU(\fU)$\+module.
 The key observation is that the rule $\fV\longmapsto
\Hom_{\fO_\fU(\fU)}(\cE(\fV),\cJ(\fU))$ for all open formal subschemes
$\fV\subset\fU$ defines a contraherent cosheaf of contramodules
on~$\fU$.
 This contraherent cosheaf of contramodules is clearly coflasque and
isomorphic to $\fHom_\fU(\cE,\cJ)$.
\end{proof}

\begin{lem} \label{finite-Krull-dim-co-flasque-co-resol-dim}
 Let\/ $\fX$ be a Noetherian formal scheme with an open
covering\/~$\bW$.
 Assume that the underlying topological space of\/ $\fX$ has finite
Krull dimension\/~$\le D$.
 In this setting: \par
\textup{(a)} Let\/ $0\rarrow\cM\rarrow\cE^0\rarrow\dotsb\rarrow
\cE^{D-1}\rarrow\cE\rarrow0$ be an exact sequence in the abelian
category\/ $\fX\Tors$.
 Assume that the sheaves\/ $\cE^i$ are flasque for all\/
$0\le i\le D-1$.
 Then the sheaf\/ $\cE$ is flasque. \par
\textup{(b)} Let\/ $0\rarrow\fE\rarrow\fE_{D-1}\rarrow\dotsb\rarrow
\fE_0\rarrow\fP\rarrow0$ be an exact sequence in the exact category\/
$\fX\Lcth_\bW$.
 Assume that the cosheaves\/ $\fE_i$ are coflasque for all\/
$0\le i\le D-1$.
 Then the cosheaf\/ $\fE$ is coflasque.
\end{lem}

\begin{proof}
 Part~(a) is a special case of~\cite[Lemma~3.4.7(a)]{Pcosh}.
 Part~(b) is similar to~\cite[Lemma~3.4.7(b)]{Pcosh}.
\end{proof}

 It follows from Lemma~\ref{coflasque-adjusted-to-cosections}(a,b)
that the full subcategories of coflasque contraherent cosheaves
$\fX\Ctrh_\cfq\subset\fX\Lcth_\bW$ and
$\fX\Ctrh^\lct_\cfq\subset\fX\Lcth^\lct_\bW$ inherit exact category
structures from $\fX\Lcth_\bW$ and $\fX\Lcth^\lct_\bW$.
 Dual-analogously, the full subcategory of flasque quasi-coherent
torsion sheaves $\fX\Tors^\fq$ is closed under extensions and cokernels
of monomorphisms in $\fX\Tors$; so it inherits an exact category
structure from the abelian exact structure of $\fX\Tors$.

 The next corollary is our version of~\cite[Corollary~3.4.8]{Pcosh}.
 It is to be compared with Corollary~\ref{dilute-direct-images} below.

\begin{cor} \label{co-flasque-direct-image-cor}
 Let\/ $\ff\:\fY\rarrow\fX$ be a quasi-compact morphism of locally
Noetherian formal schemes.  Then \par
\textup{(a)} the functor of direct image of quasi-coherent torsion
sheaves\/ $\ff_*\:\fY\Tors\rarrow\fX\Tors$ takes flasque quasi-coherent
torsion sheaves on\/ $\fY$ to flasque quasi-coherent torsion sheaves
on\/ $\fX$, and induces an exact functor\/ $\ff_*\:\fY\Tors^\fq\rarrow
\fX\Tors^\fq$ between these exact categories; \par
\textup{(b)} the functor of direct image of cosheaves of\/
$\fO$\+modules\/
$\ff_!\:(\fY,\fO_\fY)\Cosh\rarrow(\fX,\fO_\fX)\Cosh$ takes coflasque
contraherent cosheaves of contramodules on\/ $\fY$ to coflasque
contraherent cosheaves of contramodules on\/ $\fX$, and induces
an exact functor\/ $\ff_!\:\fY\Ctrh_\cfq\rarrow\fX\Ctrh_\cfq$ between
these exact categories; \par
\textup{(c)} the functor of direct image of cosheaves of\/
$\fO$\+modules\/
$\ff_!\:(\fY,\fO_\fY)\Cosh\rarrow(\fX,\fO_\fX)\Cosh$ takes coflasque
locally cotorsion contraherent cosheaves of contramodules on\/ $\fY$
to coflasque locally cotorsion contraherent cosheaves of contramodules
on\/ $\fX$, and induces an exact functor\/ $\ff_!\:\fY\Ctrh^\lct_\cfq
\rarrow\fX\Ctrh^\lct_\cfq$ between these exact categories.
\end{cor}

\begin{proof}
 The argument is similar to the one
in~\cite[proof of Corollary~3.4.8]{Pcosh}.
 All the assumptions and properties claimed are local in $\fX$, so one
can assume the formal scheme $\fX$ to be affine.
 Then the formal scheme $\fY$ is quasi-compact (i.~e., Noetherian).
 We skip the discussion of part~(a), which is easier.

 Let $\fE$ be a coflasque contraherent cosheaf of contramodules
on $\fY$, consider a finite affine open covering
$\fY=\bigcup_\alpha\fV_\alpha$, and use the fact that the \v Cech
complex $C_\bu(\{\fV_\alpha\},\fE)$ is acyclic,
by~\cite[Lemma~3.4.1(b)]{Pcosh} (since the cosheaf $\fE$ is coflasque);
so $C_\bu(\{\fV_\alpha\},\fE)$ is a finite resolution of
the contramodule $\fO_\fX(\fX)$\+module~$\fE[\fY]$.
 Notice first of all that $\ff_!(\fE)$ is a cosheaf of contramodules
on $\fX$ by Lemma~\ref{contramodule-cosheaves-direct-images-lemma}.

 The case when the formal scheme $\fY$ is semi-separated needs to
be considered first, before passing to the general case.
 Then $C_\bu(\{\fV_\alpha\},\fE)$ is a finite resolution of $\fE[\fY]$
by contraadjusted contramodule $\fO_\fX(\fX)$\+modules (in the context
of~(b), by Lemma~\ref{restriction-of-scalars-contraadjusted}) or even
by cotorsion contramodule $\fO_\fX(\fX)$\+modules (in the context
of~(c), by Lemma~\ref{restriction-of-scalars-cotorsion}).
 It follows that the contramodule $\fO_\fX(\fX)$\+module
$(\ff_!\fE)[\fX]=\fE[\fY]$ is contraadjusted in the case~(b) and
cotorsion in the case~(c).

 Now let $\fU\subset\fX$ be an affine open formal subscheme.
 Using Lemma~\ref{preimage-of-affine-open-Hom-into-contramod}, one
shows that the complex $C_\bu(\{\ff^{-1}(\fU)\cap\fV_\alpha\},
\fE|_{\ff^{-1}(\fU)})$ can be obtained by applying the functor
$\Hom_{\fO_\fX(\fX)}(\fO_\fX(\fU),{-})$ to the complex
$C_\bu(\{\fV_\alpha\},\fE)$.
 In view of Corollary~\ref{colocalization-contraadjusted-exactness-cor},
it follows that $(\ff_!\fE)[\fU]=\fE[\ff^{-1}(\fU)]\simeq
\Hom_{\fO_\fX(\fX)}(\fO_\fX(\fU),\fE[\fY])$, so the contraherence
axiom holds for the cosheaf of contramodules $\ff_!\fE$ on~$\fX$.

 Finally, in the general case of a non-semi-separated formal scheme
$\fY$, one needs to use the fact that the intersections of nonempty
subsets of the set $\{\fV_\alpha\}$ are semi-separated (in fact,
separated) open formal subschemes of $\fY$, and for the restrictions
of the morphism~$\ff$ to such open formal subschemes we already know
the desired contraadjustedness/cotorsion and contraherence properties
of the direct images of coflasque contraherent cosheaves of
contramodules (such as the restrictions of $\fE$ to the respective
open formal subschemes of~$\fY$).
\end{proof}

 The next lemma is our version of~\cite[Lemma~6.10.2]{Pcosh}.

\begin{lem} \label{co-flasque-open-subscheme-adjunction-morphism-lemma}
 Let\/ $\fX$ be a locally Noetherian scheme and\/ $\fY\subset\fX$ be
an open formal subscheme with the open immersion morphism\/
$\fj\:\fY\rarrow\fX$.
 Then \par
\textup{(a)} for any flasque quasi-coherent torsion sheaf\/ $\cE$ on\/
$\fX$, the adjunction morphism\/ $\cE\rarrow\fj_*\fj^*\cE$ is
an admissible epimorphism in the exact category\/ $\fX\Tors^\fq$
(i.~e., an epimorphism in\/ $\fX\Tors$ between two flasque
quasi-coherent torsion sheaves with a flasque kernel); \par
\textup{(b)} for any coflasque contraherent cosheaf of contramodules\/
$\fE$ on $\fX$, the adjunction morphism\/ $\fj_!\fj^!\fE\rarrow\fE$
is an admissible monomorphism in the exact category\/ $\fX\Ctrh_\cfq$
(i.~e., an admissible monomorphism in\/ $\fX\Ctrh$ between two coflasque
contraherent cosheaves of contramodules with a coflasque cokernel); \par
\textup{(c)} for any coflasque locally cotorsion contraherent cosheaf
of contramodules\/ $\fE$ on $\fX$, the adjunction morphism\/
$\fj_!\fj^!\fE\rarrow\fE$ is an admissible monomorphism in the exact
category\/ $\fX\Ctrh^\lct_\cfq$ (i.~e., an admissible monomorphism in\/
$\fX\Ctrh^\lct$ between two coflasque locally cotorsion contraherent
cosheaves of contramodules with a coflasque cokernel).
\end{lem}

\begin{proof}
 Part~(a) is a special case of~\cite[Lemma~6.10.2]{Pcosh}.
 In part~(b), \,$\fj_!\fj^!\fE$ is a coflasque contraherent cosheaf of
contramodules on $\fX$ by
Corollary~\ref{co-flasque-direct-image-cor}(b).
 Dual-analogously to~\cite[proof of Lemma~6.10.2]{Pcosh}, one shows
that the map $(\fj_!\fj^!\fE)[\fU]\rarrow\fE[\fU]$ is injective
for every open formal subscheme $\fU\subset\fX$ (not even necessarily
affine), and the cosheaf formed by the cokernels of such maps
is coflasque.
 In part~(c), \,$\fj_!\fj^!\fE$ is a coflasque locally cotorsion
contraherent cosheaf of contramodules on $\fX$ by
Corollary~\ref{co-flasque-direct-image-cor}(c).
 With this observation in mind, part~(c) follows from part~(b).
\end{proof}

 We finish the section with our version
of~\cite[Corollary~3.4.9]{Pcosh}.

\begin{cor} \label{co-flasque-direct-image-of-acyclic-complexes-cor}
 Let\/ $\ff\:\fY\rarrow\fX$ be a quasi-compact morphism of locally
Noetherian formal schemes and\/ $\bT$ be an open covering of\/~$\fY$.
 Then \par
\textup{(a)} for any complex\/ $\cE^\bu$ of flasque quasi-coherent
torsion sheaves on\/ $\fY$ that is acyclic as a complex in the abelian
category\/ $\fY\Tors$, the complex\/ $\ff_*\cE^\bu$ of flasque
quasi-coherent torsion sheaves on\/ $\fX$ is acyclic as a complex in
the abelian category\/ $\fX\Tors$; \par
\textup{(b)} for any complex\/ $\fE^\bu$ of coflasque contraherent
cosheaves of contramodules on\/ $\fY$ that is acyclic as a complex in
the exact category\/ $\fY\Lcth_\bT$, the complex\/ $\ff_!\fE^\bu$ of
coflasque contraherent cosheaves of contramodules on\/ $\fX$ is
acyclic as a complex in the exact category\/ $\fX\Ctrh$; \par
\textup{(c)} or any complex\/ $\fE^\bu$ of coflasque locally cotorsion
contraherent cosheaves of contramodules on\/ $\fY$ that is acyclic as
a complex in the exact category\/ $\fY\Lcth_\bT^\lct$, the complex\/
$\ff_!\fE^\bu$ of coflasque locally cotorsion contraherent cosheaves of
contramodules on\/ $\fX$ is acyclic as a complex in the exact category\/
$\fX\Ctrh^\lct$.
\end{cor}

\begin{proof}
 The argument is similar to the proof of~\cite[Corollary~3.4.9]{Pcosh}.
 Lemma~\ref{complex-of-cosheaves-exactness-criterion} is relevant,
but one does not need to use it in any nontrivial way.
 The point is that even though an assumption of finite Krull dimension
is \emph{not} made in this corollary and
Lemma~\ref{finite-Krull-dim-co-flasque-co-resol-dim} is \emph{not}
applicable, the functors of (co)sections over Noetherian formal
schemes still can be computed by finite \v Cech complexes of abelian
groups, which are acyclic due to the (co)flasqueness condition on
the terms of the complex $\cE^\bu$ or~$\fE^\bu$.
 Finiteness of the \v Cech complexes plays a crucial role.
\end{proof}

\subsection{Injective quasi-coherent torsion sheaves}
\label{injective-torsion-sheaves-subsecn}
 In this section we spell out the classification of injective
quasi-coherent torsion sheaves over locally Noetherian formal schemes,
based on the classification of injective torsion modules over adic
Noetherian rings presented in
Proposition~\ref{injective-torsion-modules-classification-prop}.
 The latter, in turn, is based on Matlis' classification of injective
modules over Noetherian commutative rings~\cite{Mat}.

 The classification of injective quasi-coherent sheaves over locally
Noetherian schemes, based on Matlis' work, is due to
Hartshorne~\cite[Proposition~II.7.17]{HarRD}.
 Another approach, providing a part the results in a more general
context, can be found in~\cite[Theorem~A.3]{EP}.
 The locally Noetherian formal scheme case is not that much
different from the locally Noetherian scheme case; we spell it out
here for the sake of completeness of the exposition and as
a preparation for the discussion of projective contraherent
cosheaves of contramodules over locally Noetherian formal schemes
in the next section.
 Our exposition follows the approach of~\cite[Section~6.1]{Pcosh}.

 A result comparable to the main theorem of this section
(Theorem~\ref{loc-Noetherian-injective-qcoh-torsion} below)
can be found in Yekutieli's paper~\cite[Proposition~4.2]{Yek0}.

 Let $\fX$ be a (not necessarily semi-separated) locally Noetherian
formal scheme.
 Recall that the points $x\in\fX$ of the (underlying topological
space of) the formal scheme $\fX$ are the same things as the points
of any closed subscheme of definition $X\subset\fX$.
 The stalk $\fO_{x,\fX}$ of the structure sheaf $\fO_\fX$ at a point
$x\in\fX$ is a Noetherian commutative local ring.
 Following the notation of
Section~\ref{flat-cotorsion-contramods-over-adic-Noetherian-subsecn},
we denote by $\widehat\fO_{x,\fX}$ the completion of the local
ring~$\fO_{x,\fX}$.
 Denote by $\fm_{x,\fX}$ the maximal ideal of the local ring
$\fO_{x,\fX}$ and by $\widehat\fm_{x,\fX}$ the maximal ideal of
the complete local ring $\widehat\fO_{x,\fX}$.

\begin{rem} \label{tricky-adic-topology-on-completed-stalk-ring}
 \emph{Let us emphasize a tricky aspect that can be potentially
confusing.}
 The Noetherian local ring $\widehat\fO_{x,\fX}$ is separated and
complete in its $\widehat\fm_{x,\fX}$\+adic topology, but
\emph{by default we do not endow the ring\/ $\widehat\fO_{x,\fX}$
with the\/ $\widehat\fm_{x,\fX}$\+adic topology}.
 By default, we view $\widehat\fO_{x,\fX}$ as an adic Noetherian
ring with \emph{a different (finer) adic topology than
the\/ $\widehat\fm_{x,\fX}$\+adic one}.

 Let $\fU\subset\fX$ be an affine Noetherial open formal subscheme
such that $x\in\fU$.
 Then we have natural ring homomorphisms $\fO_\fX(\fU)\rarrow
\fO_{x,\fX}\rarrow\widehat\fO_{x,\fX}$.
 The commutative Noetherian ring $\fO_\fX(\fU)$ carries an adic
topology that is a part of the formal scheme structure of~$\fX$.
 We endow the ring $\widehat\fO_{x,\fX}$ with \emph{the unique
adic topology for which the ring homomorphism $\fO_\fX(\fU)\rarrow
\widehat\fO_{x,\fX}$ is continuous and tight}.
 By Corollary~\ref{finer-and-coarser-adic-topology-completeness-cor},
the ring $\widehat\fO_{x,\fX}$ is separated and complete in its adic
topology (defined in the previous sentence).

 The notation $\Spf\widehat\fO_{x,\fX}$ will always stand for
the formal spectrum of the adic Noetherian ring $\widehat\fO_{x,\fX}$
with the adic topology defined in the previous paragraph.
 Notice that, while the underlying topological space of the formal
spectrum of the complete Noetherian local ring $\widehat\fO_{x,\fX}$
with its $\widehat\fm_{x,\fX}$\+adic topology is a singleton set,
the underlying topological space of the affine Noetherian formal
scheme $\Spf\widehat\fO_{x,\fX}$ as defined above usually has
many points.

 One can readily check that our adic topology on the ring
$\widehat\fO_{x,\fX}$ does not depend on the choice of an affine
Noetherian open formal subscheme $\fU\subset\fX$ such that $x\in\fU$.
 The composition of morphisms of formal schemes $\Spf\widehat\fO_{x,\fX}
\rarrow\fU\rarrow\fX$ (where $\Spf\widehat\fO_{x,\fX}\rarrow\fU$ is
the morphism of affine formal schemes corresponding to the continuous
homomorphism of adic topological rings $\fO_\fX(\fU)\rarrow
\widehat\fO_{x,\fX}$) provides a natural morphism of formal schemes
$\Spf\widehat\fO_{x,\fX}\rarrow\fX$, which does not depend depend on
the choice of $\fU$, either.

 Nevertheless, in Theorem~\ref{loc-Noetherian-injective-qcoh-torsion}
below, we will be interested in $\widehat\fm_{x,\fX}$\+torsion
modules over the ring~$\widehat\fO_{x,\fX}$ (rather than merely
torsion ~$\widehat\fO_{x,\fX}$\+modules, which means torsion over
an ideal of definition of $\fO_\fX(\fU)$).

 To give an example, if $\fX=X$ is a (nonformal) locally Noetherian
scheme, then our adic topology on the complete Noetherian local ring
$\widehat\fO_{x,\fX}=\widehat{\mathcal O}_{x,X}$ is
the \emph{discrete} topology.
 So $\Spf\widehat\fO_{x,\fX}=\Spec\widehat{\mathcal O}_{x,X}$ is
merely a scheme (rather than a formal scheme) in this case.
 Nevertheless, we will be interested in $\widehat\fO_{x,\fX}$\+modules
that are torsion with respect to the (generally speaking, quite nonzero)
maximal ideal of the complete Noetherian local
ring~$\widehat\fO_{x,\fX}$.
\end{rem}

\begin{lem} \label{stalk-ring-map-tight-and-flat}
 Let\/ $\fU$ be an affine Noetherian formal scheme and $x\in\fU$ be
a point.
 Consider the complete Noetherian local ring\/ $\widehat\fO_{x,\fU}$
as above, and endow it with the adic topology defined in
Remark~\ref{tricky-adic-topology-on-completed-stalk-ring}.
 Then the natural map of adic Noetherian rings\/
$\fO_\fU(\fU)\rarrow\widehat\fO_{x,\fU}$ is continuous, tight,
and flat.
\end{lem}

\begin{proof}
 The map of adic topological rings $\fO_\fU(\fU)\rarrow
\widehat\fO_{x,\fU}$ is continuous and tight by construction.
 In view of
Corollary~\ref{Noetherian-adic-flat-ring-maps-characterized-cor},
it remains to show that the ring $\widehat\fO_{x,\fU}$ is flat as
an abstract (nontopological) module over the abstract
ring~$\fO_\fU(\fU)$.

 By the same corollary, for any affine open formal subscheme
$\fV\subset\fU$, the ring $\fO_\fU(\fV)$ is a flat module over
the ring $\fO_\fU(\fU)$.
 Hence the ring $\fO_{x,\fU}$ is a flat module over $\fO_\fU(\fU)$
as a direct limit of flat modules.
 Finally, it is well-known (or once again follows from the same
corollary) that the $I$\+adic completion of a Noetherian
commutative ring $R$ is a flat $R$\+module for any ideal $I\subset R$.
 In particular, the ring $\widehat\fO_{x,\fU}$ is a flat module
over the ring~$\fO_{x,\fU}$.
 By transitivity of flat homomorphisms of commutative rings,
$\widehat\fO_{x,\fU}$ is a flat module over~$\fO_\fU(\fU)$.
\end{proof}

\begin{lem} \label{stalk-completion-of-multiplication-lemma}
 Let\/ $\fU$ be an affine Noetherian formal scheme and $x\in\fU$ be
a point.
 Consider the multiplication map\/
$\widehat\fO_{x,\fU}\ot_{\fO_\fU(\fU)}\widehat\fO_{x,\fU}\rarrow
\widehat\fO_{x,\fU}$; let us view it an a morphism of\/
$\widehat\fO_{x,\fU}$\+modules with respect to the left action
of\/~$\widehat\fO_{x,\fU}$.
 Then the following induced maps are isomorphisms: \par
\textup{(a)} $\Lambda_{\widehat\fm_{x,\fU}}
(\widehat\fO_{x,\fU}\ot_{\fO_\fU(\fU)}\widehat\fO_{x,\fU})\rarrow
\Lambda_{\widehat\fm_{x,\fU}}(\widehat\fO_{x,\fU})=
\widehat\fO_{x,\fU}$; \par
\textup{(b)} $\boL_0\Lambda_{\widehat\fm_{x,\fU}}
(\widehat\fO_{x,\fU}\ot_{\fO_\fU(\fU)}\widehat\fO_{x,\fU})\rarrow
\boL_0\Lambda_{\widehat\fm_{x,\fU}}(\widehat\fO_{x,\fU})=
\widehat\fO_{x,\fU}$.
\end{lem}

\begin{proof}
 In the context of the present lemma,
the $\widehat\fm_{x,\fU}$\+adic topology on the ring
$\widehat\fO_{x,\fU}$ is relevant.
 The ring $\widehat\fO_{x,\fU}$ with its $\widehat\fm_{x,\fU}$\+adic 
topology is the completion of the ring $\fO_{x,\fU}$ with its
$\fm_{x,\fU}$\+adic topology, $\widehat\fO_{x,\fU}=
\Lambda_{\fm_{x,\fU}}(\fO_{x,\fU})$.

 Moreover, denote by $\fm_x(\fU)\subset\fO_\fU(\fU)$ the preimage
of the ideal $\fm_{x,\fU}\subset\fO_{x,\fU}$ under the ring
homomorphism $\fO_\fU(\fU)\rarrow\fO_{x,\fU}$.
 Then the ring $\widehat\fO_{x,\fU}$ with its
$\widehat\fm_{x,\fU}$\+adic  topology is the completion of the ring
$\fO_\fU(\fU)$ with its $\fm_x(\fU)$\+adic topology.
 The ideal $\widehat\fm_{x,\fU}\subset\widehat\fO_{x,\fU}$ is
the extension of the ideal $\fm_x(\fU)\subset\fO_\fU(\fU)$ in
the ring $\widehat\fO_{x,\fU}$, that is,
$\widehat\fm_{x,\fU}=\widehat\fO_{x,\fU}\cdot\fm_x(\fU)\subset
\widehat\fO_{x,\fU}$.

 Notice that the ring $\widehat\fO_{x,\fU}$ is flat as a module over
the ring $\fO_\fU(\fU)$, hence the tensor product
$\widehat\fO_{x,\fU}\ot_{\fO_\fU(\fU)}\widehat\fO_{x,\fU}$ is flat
as a left module over $\widehat\fO_{x,\fU}$.
 By Corollary~\ref{flat-reductions-derived-Lambda-is-underived-cor},
parts~(a) and~(b) of the lemma are equivalent restatements of each
other.
 So it suffices to prove part~(a).

 Similarly to the proof of
Lemma~\ref{contramodulation-of-multiplication-lemma},
the multiplication map $\widehat\fO_{x,\fU}\ot_{\fO_\fU(\fU)}
\widehat\fO_{x,\fU}\rarrow\widehat\fO_{x,\fU}$ is a split epimorphism
of left modules over~$\widehat\fO_{x,\fU}$.
 For any commutative ring $S$ with a finitely generated ideal
$J\subset S$, any $S$\+module annihilated by the functor $S/J\ot_S{-}$
is also annihilated by the functor~$\Lambda_J$.
 So it suffices to show that the functor of tensor product with
the residue field $k_{x,\fU}=\widehat\fO_{x,\fU}/\widehat\fm_{x,\fU}$
takes the multiplication map $\widehat\fO_{x,\fU}\ot_{\fO_\fU(\fU)}
\widehat\fO_{x,\fU}\rarrow\widehat\fO_{x,\fU}$ to an isomorphism.
 Now we have
\begin{multline*}
 \bigl(\widehat\fO_{x,\fU}/\widehat\fm_{x,\fU}\bigr)\ot_{\fO_\fU(\fU)}
 \widehat\fO_{x,\fU}\simeq\bigl(\fO_\fU(\fU)/\fm_x(\fU)\bigr)
 \ot_{\fO_\fU(\fU)}\widehat\fO_{x,\fU} \\ \simeq
 \widehat\fO_{x,\fU}/\bigl(\fm_x(\fU)\cdot\widehat\fO_{x,\fU}\bigr)=
 \widehat\fO_{x,\fU}/\widehat\fm_{x,\fU}.
\end{multline*}
\end{proof}

 Given a complete, separated adic Noetherian ring $\fR$ and
a torsion $\fR$\+module $M\in\fR\Tors$, we denote by
$\widetilde M\in\Spf\fR\Tors$ the quasi-coherent torsion sheaf on
the affine Noetherian formal scheme $\Spf\fR$ corresponding to
the torsion $\fR$\+module $M$ under the equivalence of abelian
categories $\fR\Tors\simeq\Spf\fR\Tors$ from
Section~\ref{qcoh-torsion-sheaves-subsecn}.

\begin{thm} \label{loc-Noetherian-injective-qcoh-torsion}
\textup{(a)} There are enough injective objects in the abelian
category of quasi-coherent torsion sheaves\/ $\fX\Tors$. \par
\textup{(b)} For any point $x\in\fX$, denote by $\iota_x\:
\Spf\widehat\fO_{x,\fX}\rarrow\fX$ the natural morphism of formal
schemes.
 Then a quasi-coherent torsion sheaf\/ $\cJ$ on\/ $\fX$ is injective
if and only if it is isomorphic to an infinite direct sum\/
$\bigoplus_{x\in\fX}\iota_x{}_*\widetilde J_x$ of the direct images
$\iota_x{}_*\widetilde J_x$ of quasi-coherent torsion sheaves
$\widetilde J_x$ on\/ $\Spf\widehat\fO_{x,\fX}$ corresponding to some
injective\/ $\widehat\fm_{x,\fX}$\+torsion\/
$\widehat\fO_{x,\fX}$\+modules~$J_x$
(i.~e., $J_x$ is an injective object of the abelian category\/
$\widehat\fO_{x,\fX}\Modl_{\widehat\fm_{x,\fX}\dtors}$).
\end{thm}

\begin{proof}
 To prove part~(a), it suffices to point out that $\fX\Tors$ is
a Grothendieck abelian category; hence it has enough injective
objects (cf.\ the first paragraph of the proof of
Lemma~\ref{injective-torsion-sheaves-on-qcomp}).
 Alternatively, the argument below, based on
Proposition~\ref{injective-torsion-modules-classification-prop},
proves parts~(a) and~(b) simultaneously.

 Notice that, when the formal scheme $\fX$ is not semi-separated,
the morphisms of formal schemes~$\iota_x$ \emph{need not} be affine.
 However, any locally Noetherian formal scheme is quasi-separated;
and any morphism from an affine to a quasi-separated (formal)
scheme is quasi-compact.
 Therefore, by Proposition~\ref{torsion-sheaves-qcoh-direct-images-prop}
and formula~\eqref{direct-image-qcoh-tors}, we have well-defined
direct image functors $\iota_x{}_*\:\Spf\widehat\fO_{x,\fX}\Tors
\rarrow\fX\Tors$.
 Hence the quasi-coherent torsion sheaves on $\fX$ of the form
$\cJ=\bigoplus_{x\in\fX}\iota_x{}_*\widetilde J_x$ are well-defined,
too.

 Let us show that sheaves of the form
$\bigoplus_{x\in\fX}\iota_x{}_*\widetilde J_x$, where
$J_x\in\widehat\fO_{x,\fX}\Modl_{\widehat\fm_{x,\fX}\dtors}^\inj$,
are injective objects of $\fX\Tors$.
 As in~\cite[proof of Theorem~6.1.1]{Pcosh}, we pick an affine open
covering $\fX=\bigcup_\alpha\fU_\alpha$ of the formal scheme $\fX$,
and choose a well-ordering of the set of indices~$\{\alpha\}$.
 Put $\fS_\alpha=\fU_\alpha\setminus\bigcup_{\gamma<\alpha}\fU_\gamma$
(so $\fS_\alpha$ is a subset of the underlying topological space
of~$\fX$).
 For any quasi-compact open formal subscheme $\fU\subset\fX$, one
has $\fU\cap\fS_\alpha=\varnothing$ for all but a finite subset
of indices~$\alpha$ (since the underlying topological space of
$\fU$ is Noetherian, so any ascending chain of open subsets in it
terminates).

 Next let us show that the restriction of the quasi-coherent torsion
sheaf $\iota_x{}_*\widetilde M_x$ on $\fX$ to an affine open formal
subscheme $\fU\subset\fX$ vanishes for every
$\widehat\fm_{x,\fX}$\+torsion $\widehat\fO_{x,\fX}$\+module $M_x$
whenever $x\notin\fU$.
 Indeed, it suffices to check that the module of sections
$\widetilde M_x(\fV)$ vanishes for every principal affine open formal
subscheme $\fV\subset\Spf\widehat\fO_{x,\fX}$ that does not contain
the closed point of $\Spf\widehat\fO_{x,\fX}$ (as $x\in\fX$ is
the image of the unique closed point of $\Spf\widehat\fO_{x,\fX}$ under
the morphism~$\iota_x$).
 So we need to check that $\widehat\fO_{x,\fX}[s^{-1}]
\ot_{\widehat\fO_{x,\fX}}M_x=0$ for all $s\in\widehat\fm_{x,\fX}$.
 This holds because $M_x$ is an $\widehat\fm_{x,\fX}$\+torsion
$\widehat\fO_{x,\fX}$\+module.

 Put $\cJ_\alpha=\bigoplus_{z\in\fS_\alpha}\iota_z{}_*\widetilde J_z
\in\fX\Tors$.
 Clearly, $\cJ\simeq\bigoplus_\alpha\cJ_\alpha$.
 It follows from the previous two paragraphs that, for any affine open
formal subscheme $\fU\subset\fX$, one has $\cJ_\alpha(\fU)=0$
for all but a finite subset of indices~$\alpha$.
 Thus the sheaf $\cJ$ is also the product of the sheaves $\cJ_\alpha$
(in the category $(\fX,\fO_\fX)\Sh$ of sheaves of $\fO_\fX$\+modules
on~$\fX$).
 Therefore, it suffices to check that each sheaf $\cJ_\alpha$ is
an injective object in $\fX\Tors$.

 Denote by~$\fj_\alpha$ the open immersion $\fU_\alpha\rarrow\fX$ and
by~$\kappa_z$ the natural morphisms $\Spf\widehat\fO_{z,\fX}\rarrow
\fU_\alpha$ for $z\in\fS_\alpha$.
 Then we have $\cJ_\alpha=\fj_\alpha{}_*\cK_\alpha$, where
$\cK_\alpha=\bigoplus_{z\in\fS_\alpha}\kappa_z{}_*\widetilde J_z$.
 The observation that the direct images with respect to quasi-compact
morphisms of locally Noetherian formal schemes preserve infinite
direct sums of quasi-coherent torsion sheaves (see
Section~\ref{direct-images-of-qcoh-tors-subsecn}) is used here.
 By Proposition~\ref{injective-torsion-modules-classification-prop},
\,$\cK_\alpha$ is an injective quasi-coherent torsion sheaf
on~$\fU_\alpha$.
 By the adjunction~\eqref{qcoh-tors-inverse-direct-adjunction}, it
follows that $\cJ_\alpha$ is an injective quasi-coherent torsion sheaf
on~$\fX$.

 Now let us construct for any quasi-coherent torsion sheaf $\cM$ on
$\fX$ an injective morphism from $\cM$ into a sheaf
$\cJ=\bigoplus_{x\in\fX}\iota_x{}_*\widetilde J_x$ as in part~(b) of
the theorem.
 For this purpose, we choose an affine open covering $\fU_\alpha$ as
above, and proceed by transfinite induction on~$\alpha$, constructing
quasi-coherent torsion sheaves $\cJ_\alpha=\bigoplus_{z\in\fS_\alpha}
\iota_z{}_*\widetilde J_z$ and morphisms of quasi-coherent torsion
sheaves $\cM\rarrow\cJ_\alpha$.

 Suppose such sheaves and morphisms have been constructed for all
indices $\alpha<\beta$, with some fixed index~$\beta$.
 Then, as explained above, there is the induced morphism of
quasi-coherent torsion sheaves $\cM\rarrow\bigoplus_{\alpha<\beta}
\cJ_\alpha=\prod_{\alpha<\beta}\cJ_\alpha$.
 Assume that the related morphism of the modules of sections
$\cM(\fU)\rarrow\bigoplus_{\alpha<\beta}\cJ_\alpha(\fU)=
\prod_{\alpha<\beta}\cJ_\alpha(\fU)$ is injective for every affine
open formal subscheme $\fU\subset\bigcup_{\alpha<\beta}\fU_\alpha
\subset\fX$.
 Our aim is to construct a quasi-coherent torsion sheaf
$\cJ_\beta=\bigoplus_{z\in\fS_\beta}\iota_z{}_*\widetilde J_z$
and a morphism of quasi-coherent torsion sheaves $\cM\rarrow\cJ_\beta$
such that the resulting morphism $\cM\rarrow\bigoplus_{\alpha\le\beta}
\cJ_\alpha$ satisfies the same condition as above for every affine
open formal subscheme $\fU\subset\bigcup_{\alpha\le\beta}\fU_\alpha
\subset\fX$.

 Choose a monomorphism in $\fU_\beta\Tors$ from the quasi-coherent
torsion sheaf $\fj_\beta^*\cM$ into an injective quasi-coherent
torsion sheaf $\cK$ on~$\fU_\beta$.
 By Proposition~\ref{injective-torsion-modules-classification-prop},
the sheaf $\cK$ decomposes into a direct sum $\cK\simeq
\bigoplus_{x\in\fU_\beta}\lambda_x{}_*\widetilde K_x$, where $K_x$
are some injective $\widehat\fm_{x,\fX}$\+torsion
$\widehat\fO_{x,\fX}$\+modules and $\lambda_x\:\Spf\widehat\fO_{x,\fX}
\rarrow\fU_\beta$ are the natural morphisms of formal schemes.
 Let us split this infinite direct sum into two summands as
$\cK=\bigoplus_{z\in\fS_\beta}\lambda_z{}_*\widetilde K_z\oplus
\bigoplus_{y\in\fU_\beta\setminus\fS_\beta}\lambda_y{}_*
\widetilde K_y$.
 Denote the former direct summand by $\cK_\beta$ and the latter one
by~$\cL$.
 Set $\cI=\bigoplus_{\alpha<\beta}\cJ_\alpha\in\fX\Tors$
and $\cJ_\beta=\fj_\beta{}_*\cK_\beta\in\fX\Tors$.

 The property of a morphism of quasi-coherent torsion sheaves to be
a monomorphism is local; the morphism $\cM\rarrow\cI$ becomes
a monomorphism when restricted to $\bigcup_{\alpha<\beta}\fU_\alpha$
by assumption; so we only need to show that the natural morphism
$\cM\rarrow\cI\oplus\cJ_\beta$ becomes a monomorphism when restricted
to~$\fU_\beta$.
 This means that the morphism $\fj_\beta^*\cM\rarrow\fj_\beta^*\cI
\oplus\cK_\beta$ should be a monomorphism in $\fU_\beta\Tors$.
 It suffices to check that the monomorphism $\fj_\beta^*\cM\rarrow
\cK$ factorizes through the morphism $\fj_\beta^*\cM\rarrow
\fj_\beta^*\cI\oplus\cK_\beta$, or that the morphism $\fj_\beta^*\cM
\rarrow\cL$ factorizes through the morphism $\fj_\beta^*\cM\rarrow
\fj_\beta^*\cI$.

 Put $\fV_\beta=\fU_\beta\setminus\fS_\beta=\fU_\beta\cap
\bigcup_{\alpha<\beta}\fU_\alpha$.
 So $\fV_\beta$ is a Noetherian formal scheme (not necessarily affine)
and an open formal subscheme in~$\fX$.
 Denote by $\fj$ the identity open immersion $\fV_\alpha\rarrow
\fU_\beta$.
 Then we have $\cL=\fj_*\cE$, where $\cE$ is an injective object in
$\fV_\beta\Tors$ (as we have already proved above).
 Hence the morphism $\fj^*\fj_\beta^*\cM\rarrow\cE$ factorizes through
the monomorphism $\fj^*\fj_\beta^*\cM\rarrow\fj^*\fj_\beta^*\cI$.
 Therefore, the morphism $\fj_\beta^*\cM\rarrow\fj_*\cE$ factorizes
through the morphism $\fj_\beta^*\cM\rarrow\fj_\beta^*\cI$, as desired.
 Here we are using the assumption that the morphism $\cM\rarrow\cI$
is a monomorphism in restriction to $\bigcup_{\alpha<\beta}\fU_\alpha
\supset\fV_\beta$, as well as the adjunction
isomorphism~\eqref{sheaves-restriction-direct-image-adjunction}.

 We have proved part~(a) and almost proved part~(b).
 In order to finish the proof of part~(b), it remains to show that
the class of quasi-coherent torsion sheaves of the form
$\cJ=\bigoplus_{x\in\fX}\iota_x{}_*\widetilde J_x$ on $\fX$ is closed
under direct summands.
 The argument is based on the following two lemmas.

\begin{lem} \label{sheaves-lower-upper-sets-fully-faithful-lemma}
\textup{(a)} Suppose the set of all points (of the underlying
topological space) of\/ $\fX$ is presented as a union of two
nonintersecting subsets\/ $\fX=\fS\sqcup\fT$ such that for any points
$z\in\fS$ and $y\in\fT$ the closure of~$z$ in\/ $\fX$ does not
contain~$y$.
 Then, for any sheaves of\/ $\fO$\+modules\/ $\cM_y$ over\/
$\Spf\widehat\fO_{y,\fX}$ and any\/ $\widehat\fm_{z,\fX}$\+torsion\/
$\widehat\fO_{z,\fX}$\+modules $M_z$ one has\/
$\Hom_{\fO_\fX}(\bigoplus_{z\in\fS}\iota_z{}_*\widetilde M_z,\>
\bigoplus_{y\in\fT}\iota_y{}_*\cM_y)=0$. \par
\textup{(b)} For any point $x\in\fX$, the functor assigning to
an\/ $\widehat\fm_{x,\fX}$\+torsion\/ $\widehat\fO_{x,\fX}$\+module
$M_x$ the quasi-coherent torsion sheaf\/ $\iota_x{}_*\widetilde M_x$
on\/ $\fX$ is fully faithful.
\end{lem}

\begin{proof}
 Part~(a): notice that the infinite direct sums in $\fX\Tors$ agree
with the ones in $(\fX,\fO_\fX)\Sh$ (see
Section~\ref{qcoh-torsion-sheaves-subsecn}).
 By the definition of the infinite direct sum, it suffices to check
that $\Hom_{\fO_\fX}(\iota_z{}_*\widetilde M_z,\>
\bigoplus_{y\in\fT}\iota_y{}_*\cM_y)=0$ for every $z\in\fS$.
 Let $\fZ$ be the closure of~$z$ in $\fX$ and $\fY=\fX\setminus\fZ$
be its open complement; then $\fT\subset\fY$ and $\fY$ is an open
formal subscheme in~$\fX$.
 Denote the related open immersion morphism by $\fj\:\fY\rarrow\fX$.
 For every point $y\in\fT$, denote the natural morphism of formal
schemes $\Spf\widehat\fO_{y,\fX}\rarrow\fY$ by~$\kappa_y$; so we have
$\iota_y=\fj\circ\kappa_y$.

 Now we have $\bigoplus_{y\in\fT}\iota_y{}_*\cM_y\simeq
\fj_*\bigoplus_{y\in\fT}\kappa_y{}_*\cM_y$ (recall that any open
immersion of locally Noetherian formal schemes is a quasi-compact
morphism).
 According to the adjunction~\eqref{sheaves-inverse-direct-adjunction},
we have
$$
 \Hom_{\fO_\fX}\bigl(\iota_z{}_*\widetilde M_z,\>
 \fj_*\bigoplus\nolimits_{y\in\fT}\kappa_y{}_*\cM_y\bigr)\simeq
 \Hom_{\fO_\fY}\bigl(\fj^*\iota_z{}_*\widetilde M_z,\>
 \bigoplus\nolimits_{y\in\fT}\kappa_y{}_*\cM_y\bigr).
$$
 It was shown above in the proof of the theorem that
$\fj^*\iota_z{}_*\widetilde M_z=0$, so the assertion of part~(a)
follows.

 Part~(b): the morphism of formal schemes~$\iota_x$ is tight (and flat)
by Lemma~\ref{stalk-ring-map-tight-and-flat}.
 So the adjunction~\eqref{qcoh-tors-inverse-direct-adjunction} is
applicable, and it suffices to show that the adjunction morphism
$\iota_x^*\iota_x{}_*\widetilde M_x\rarrow\widetilde M_x$ is
an isomorphism in $\Spf\widehat\fO_{x,\fX}\Tors$.
 One can replace the formal scheme $\fX$ by any affine open formal
subscheme $\fU\subset\fX$ containing~$x$; then we need to show that
the natural map of $\widehat\fO_{x,\fX}$\+modules
$\widehat\fO_{x,\fX}\ot_{\fO_\fU(\fU)}M_x\rarrow M_x$
is an isomorphism.
 Indeed, we have
$$
 \widehat\fO_{x,\fX}\ot_{\fO_\fU(\fU)}M_x\simeq
 \bigl(\widehat\fO_{x,\fX}\ot_{\fO_\fU(\fU)}\ot\widehat\fO_{x,\fX}\bigr)
 \ot_{\widehat\fO_{x,\fX}}M_x
$$
while $M_x\simeq\widehat\fO_{x,\fX}\ot_{\widehat\fO_{x,\fX}}M_x$,
and it remains to apply (the left-right opposite version of)
Lemma~\ref{stalk-completion-of-multiplication-lemma} together with
Corollary~\ref{derived-completion-tensor-with-torsion-isom-cor}.
\end{proof}

 The next category-theoretic lemma is just the opposite version
of~\cite[Lemma~6.1.3]{Pcosh}.
 We present the full formulation here for the reader's convenience.

\begin{lem} \label{infinite-coproduct-direct-summand-lemma}
 Let\/ $\sA$ be an additive category and $M_\alpha$ be a family of
objects of\/ $\sA$ indexed by a well-ordered set of
indices\/~$\{\alpha\}$.
 Assume that, for every index~$\beta$, the coproducts\/
$\bigoplus_{\alpha<\beta}M_\alpha$ and\/ $\bigoplus_{\alpha\ge\beta}
M_\alpha$ exist in $\sA$, and\/
$\Hom_\sA(\bigoplus_{\alpha<\beta}M_\alpha,\>
\bigoplus_{\alpha\ge\beta}M_\alpha)=0$.
 Assume further that the images of all idempotent endomorphisms of
the objects $M_\alpha$ exist in\/~$\sA$.
 Then the images of all idempotent endomorphisms of the object\/
$\bigoplus_\alpha M_\alpha$ exist in\/ $\sA$, too.
 Moreover, those images are isomorphic to coproducts over\/~$\alpha$
of images of idempotent endomorphisms of~$M_\alpha$.  \qed
\end{lem}

 Now we are ready to finish the proof of the theorem.
 Firstly, consider the case of a Noetherian (i.~e., quasi-compact)
formal scheme~$\fX$.
 Following the proof of~\cite[Theorem~6.1.1]{Pcosh}, there exists
a well-ordering relation~$<$ on the set of points of $\fX$ such that
$x<y$ whenever $x$ belongs to the closure of~$y$.
 Now we apply
Lemma~\ref{infinite-coproduct-direct-summand-lemma}
to the full additive subcategory $\sA$ of objects of the form
$\cJ=\bigoplus_{x\in\fX}\iota_x{}_*\widetilde J_x$ in $\fX\Tors$, and
claim that the images of all idempotent endomorphisms exist in~$\sA$.
 Indeed, the first assumption of
Lemma~\ref{infinite-coproduct-direct-summand-lemma} holds by
Lemma~\ref{sheaves-lower-upper-sets-fully-faithful-lemma}(a),
while the second assumption follows from
Lemma~\ref{sheaves-lower-upper-sets-fully-faithful-lemma}(b).

 For an arbitrary locally Noetherian formal scheme $\fX$, we consider
the \emph{direct product} decomposition $\cJ=\bigoplus_\alpha\cJ_\alpha
=\prod_\alpha\cJ_\alpha$ of the object $\cJ=\bigoplus_{x\in\fX}
\iota_x{}_*\widetilde J_x$ in the category of sheaves of
$\fO_\fX$\+modules, as per the discussion above in the beginning of
the proof of the theorem.
 Here $\fX=\bigcup_\alpha\fU_\alpha$ is an affine open covering of
$\fX$ indexed by a well-ordered set of indices~$\{\alpha\}$.
 We apply Lemma~\ref{infinite-coproduct-direct-summand-lemma} to
the \emph{opposite category} $\sA$ to the full additive subcategory
of objects of the form $\bigoplus_{x\in\fX}\iota_x{}_*\widetilde J_x$
in the category $(\fX,\fO_\fX)\Sh$, and claim that the images of
idempotent endomorphisms exist in~$\sA$.
 Once again, the first assumption of
Lemma~\ref{infinite-coproduct-direct-summand-lemma} holds by
Lemma~\ref{sheaves-lower-upper-sets-fully-faithful-lemma}(a).
 To check the second assumption, recall that $\cJ_\alpha\simeq
\fj_\alpha{}_*\cK_\alpha$, where $\cK_\alpha=\bigoplus_{z\in\fS_\alpha}
\kappa_z{}_*\widetilde J_z$.
 The functor~$\fj_\alpha{}_*$ is fully faithful (as one can see from
the adjunction~\eqref{sheaves-restriction-direct-image-adjunction}
or~\eqref{qcoh-tors-inverse-direct-adjunction}), so the images of
all the idempotent endomorphisms of the sheaves $\cJ_\alpha$ have
the desired form by the argument in the previous paragraph applied to
the Noetherian affine formal scheme $\fU_\alpha$ in the role of~$\fX$.
\end{proof}

 Recall the notation $\fX\Tors^\inj\subset\fX\Tors$ from
Section~\ref{semi-separated-inj-proj-subsecn} for the full subcategory
of injective quasi-coherent torsion sheaves on~$\fX$.
 According to the following corollary, being an injective
quasi-coherent torsion sheaf on a locally Noetherian formal scheme
is a local property.
 Furthermore, all such sheaves are flasque.
 For a locally Noetherian (nonformal) scheme, these results are
essentially due to Hartshorne~\cite[Theorem~II.7.18]{HarRD}.

\begin{cor} \label{injective-torsion-local-flasque}
 Let\/ $\fX$ be a locally Noetherian formal scheme. \par
\textup{(a)} Let\/ $\fY\subset\fX$ be an open formal subscheme.
 Then for any sheaf\/ $\cJ\in\fX\Tors^\inj$ the sheaf\/
$\cJ|_\fY$ belongs to\/ $\fY\Tors^\inj$. \par
\textup{(b)} In the context of part~\textup{(a)}, the restriction map\/
$\cJ(\fX)\rarrow\cJ(\fY)$ is surjective (so the sheaf\/ $\cJ$ is
flasque). \par
\textup{(c)} For an affine Noetherian formal scheme\/ $\fX$,
a (not necessarily affine) open formal subscheme\/ $\fY\subset\fX$,
and a sheaf\/ $\cJ\in\fX\Tors^\inj$, the restriction sequence of
sections\/ $0\rarrow\ker\bigl(\cJ(\fX)\to\cJ(\fY)\bigr)\rarrow\cJ(\fX)
\rarrow\cJ(\fY)\rarrow0$ is a split short exact sequence of injective
torsion\/ $\fO_\fX(\fX)$\+modules. \par
\textup{(d)} Let\/ $\fX=\bigcup_\alpha\fY_\alpha$ be an open covering
of a locally Noetherian formal scheme.
 Then a quasi-coherent torsion sheaf on $\fX$ belongs to\/
$\fX\Tors^\inj$ if and only if its restrictions to\/ $\fY_\alpha$
belong to\/ $\fY_\alpha\Tors^\inj$ for all\/~$\alpha$.
\end{cor}

\begin{proof}
 Part~(a): denote the open immersion morphism by $\fj\:\fY\rarrow\fX$.
 Then by Theorem~\ref{loc-Noetherian-injective-qcoh-torsion}(b) we have
$\fj^*\cJ\simeq\fj^*\bigl(\bigoplus_{x\in\fX}\iota_x{}_*\widetilde J_x
\bigr)\simeq\bigoplus_{x\in\fX}\fj^*\iota_x{}_*\widetilde J_x$.
 Furthermore, it was explained in the proof of the same theorem that
$\fj^*\iota_z{}_*\widetilde J_z=0$ for all $z\in\fX\setminus\fY$.
 Given a point $y\in\fY$, denote by $\kappa_y\:\Spf\widehat\fO_{y,\fX}
\rarrow\fY$ the natural morphism of formal schemes.
 Then we have $\fj^*\iota_y{}_*\widetilde J_y\simeq
\kappa_y{}_*\widetilde J_y$.
 Hence the isomorphism $\fj^*\cJ\simeq\bigoplus_{y\in\fY}
\kappa_y{}_*\widetilde J_y$, implying that $\fj^*\cJ\in\fY\Tors$
according to the same theorem.

 Part~(b): by the definition, we have
$(\iota_x{}_*\widetilde J_x)(\fX)\simeq J_x$.
 Since the sections over a quasi-compact open formal subscheme
preserve infinite direct sums of quasi-coherent torsion sheaves,
we conclude from the computation above that $\cJ(\fU)\simeq
\bigoplus_{x\in\fU}J_x$ for any quasi-compact open formal
subscheme $\fU\subset\fX$.
 Hence, once the isomorphism of quasi-coherent tosion sheaves
$\cJ\simeq\bigoplus_{x\in\fX}\iota_x{}_*\widetilde J_x$ has been
chosen, the $\fO_\fX(\fU)$\+module $\cJ(\fY\cap\fU)$ becomes naturally
a direct summand of the $\fO_\fX(\fU)$\+module $\cJ(\fU)$.
 On the other hand, the sheaf axiom~\eqref{topology-base-sheaf-axiom}
implies the isomorphism $\cJ(\fX)\simeq\varprojlim_{\fU\subset\fX}
\cJ(\fU)$, with the filtered projective limit taked over all
quasi-compact open formal subschemes $\fU\subset\fX$; and similarly
for~$\cJ(\fY)$.
 Thus the restriction map of $\fO_\fX(\fX)$\+modules $\cJ(\fX)\rarrow
\cJ(\fY)$ is the projection onto a direct summand.
 We have also proved part~(c) along the way.

 Part~(d): the ``only if'' assertion holds by part~(a).
 To prove the ``if'', we proceed by transfinite induction on
any well-ordering of the set of indices~$\alpha$.
 This reduces the question to the following assertion.
 Suppose $\fX$ is represented as the union of two open formal
subschemes, $\fX=\fW\cup\fY$.
 Suppose further that the restriction of a quasi-coherent torsion
sheaf $\cJ$ on $\fX$ onto $\fW$ is identified with the direct sum
$\bigoplus_{w\in\fW}\kappa_w{}_*\widetilde J_w$, where $J_w$ are
some injective $\widehat\fm_{w,\fX}$\+torsion
$\widehat\fO_{w,\fX}$\+modules, while~$\kappa_w$ are the natural
morphisms of formal schemes $\Spf\widehat\fO_{w,\fX}\rarrow\fW$.
 In addition, suppose that the restriction of $\cJ$ onto $\fY$ is
an injective quasi-coherent torsion sheaf.
 Then there exist some injective $\widehat\fm_{z,\fX}$\+torsion
$\widehat\fO_{z,\fX}$\+modules $J_z$ defined for all the points
$z\in\fX\setminus\fW$ and an isomorphism of quasi-coherent torsion
sheaves $\cJ\simeq\bigoplus_{x\in\fX}\iota_x{}_*\widetilde J_x$
whose restriction to $\fW$ coincides with the given isomorphism
$\cJ|_\fW\simeq\bigoplus_{w\in\fW}\kappa_w{}_*\widetilde J_w$.

 Indeed, by Theorem~\ref{loc-Noetherian-injective-qcoh-torsion} we have
$\cJ|_\fY\simeq\bigoplus_{y\in\fY}\iota'_y{}_*\widetilde J'_y$,
where $\iota_y'\:\Spf\widehat\fO_{y,\fX}\rarrow\fY$ are the natural
morphisms of formal schemes and $J'_y$ are some injective
$\widehat\fm_{y,\fX}$\+torsion $\widehat\fO_{y,\fX}$\+modules.
 Restricting to the intersection $\fV=\fW\cap\fY$, we obtain
isomorphism of quasi-coherent torsion sheaves
$\bigoplus_{v\in\fV}\kappa'_v{}_*\widetilde J_v\simeq
\bigoplus_{v\in\fV}\kappa'_v{}_*\widetilde J'_v$, where $\kappa'_v\:
\Spf\widehat\fO_{v,\fX}\rarrow\fV$ are the natural morphisms.
 It is clear from the arguments in the second half of the proof of
Theorem~\ref{loc-Noetherian-injective-qcoh-torsion}(b) that such
an isomorphism of infinite direct sums induces an isomorphism of
$\widehat\fm_{y,\fX}$\+torsion $\widehat\fO_{y,\fX}$\+modules
$J_v\simeq J'_v$.
 (Basically, if an uppertriangular $3\times 3$ matrix of morphisms
is invertible, then its middle diagonal morphism is invertible.)
 Let us identify $J'_v$ with $J_v$ using this isomorphism, and
set $J_y=J_y'$ for all $y\in\fX\setminus\fW\subset\fY$.
 Then our isomorphism of infinite direct sums can be viewed as
an automorphism~$\phi$ of the quasi-coherent torsion sheaf
$\bigoplus_{v\in\fV}\kappa'_v{}_*\widetilde J_v$ on~$\fV$.
 We would like to show that the sheaf $\cJ$ on $\fX$ is isomorphic to
$\bigoplus_{x\in\fX}\iota_x{}_*\widetilde J_x$ in a way compatible
with the given isomorphism $\cJ|_\fW\simeq\bigoplus_{w\in\fW}
\kappa_w{}_*\widetilde J_w$ over~$\fW$.

 Since a sheaf of $\fO_\fX$\+modules is determined by its restrictions
to $\fW$ and $\fY$ together with the induced isomorphism between
the restrictions of these to the intersection $\fV=\fW\cap\fY$, it
suffices to check that the automorphism~$\phi$ can be extended to
an automorphism of the sheaf $\bigoplus_{y\in\fY}\iota'_y{}_*
\widetilde J_y$ on~$\fY$.
 The latter sheaf decomposes into the direct sum of two direct sums:
the one taken over the points $v\in\fV$ and the one ranging over
the points $y\in\fY\setminus\fV$.
 So it remains to show that the automorphism~$\phi$ can be lifted to
an automorphism of the sheaf $\bigoplus_{v\in\fV}\iota'_v{}_*
\widetilde J_v$ on~$\fY$.
 Actually, the rings of endomorphisms of the quasi-coherent torsion
sheaves $\bigoplus_{v\in\fV}\kappa'_v{}_*\widetilde J_v$ on $\fV$
and $\bigoplus_{v\in\fV}\iota'_v{}_*\widetilde J_v$ on $\fY$ are
isomorphic.
 Indeed, the functor of direct image of quasi-coherent torsion
sheaves~$\fj_*$ with respect to the open immersion $\fj\:\fV\rarrow\fY$
is fully faithful (in view of
the adjunction~\eqref{sheaves-restriction-direct-image-adjunction}
or~\eqref{qcoh-tors-inverse-direct-adjunction}) and preserves
infinite direct sums.
\end{proof}

\begin{cor} \label{local-section-criterion-of-injectivity-cor}
 Let\/ $\fX$ be a locally Noetherian formal scheme with an open
covering\/ $\bW$, and let\/ $\cJ$ be a quasi-coherent torsion sheaf
on\/~$\fX$.
 Then one has\/ $\cJ\in\fX\Tors^\inj$ if and only if, for every
affine open formal subscheme\/ $\fU\subset\fX$ subordinate to\/ $\bW$,
the torsion\/ $\fO_\fX(\fU)$\+module\/ $\cJ(\fU)$ is injective.
\end{cor}

\begin{proof}
 If $\cJ$ is an injective quasi-coherent torsion sheaf on $\fX$,
then the torsion $\fO_\fX(\fU)$\+module $\cJ(\fU)$ is injective
for every affine open formal subscheme $\fU\subset\fX$ by
Corollary~\ref{injective-torsion-local-flasque}(a).
 Conversely, let $\fX=\bigcup_\alpha\fU_\alpha$ be an affine open
covering of $\fX$ subordinate to~$\bW$.
 If the torsion $\fO_\fX(\fU_\alpha)$\+module $\cJ(\fU_\alpha)$ is
injective for every~$\alpha$, then one has $\cJ\in\fX\Tors^\inj$ by
Corollary~\ref{injective-torsion-local-flasque}(d).
\end{proof}

 Now we can prove a theorem promised near the end of
Section~\ref{colocality-of-flatness-of-contramods-subsecn}.

\begin{proof}[Proof of
Theorem~\ref{locality-of-injective-torsion-modules-thm}]
 Let $R$ be an adic Noetherian ring and $f_\alpha\:R\rarrow S_\alpha$
be a formal open covering.
 Let $K$ be a torsion $R$\+module.
 We need to prove that the torsion $R$\+module $K$ is injective
if and only if the torsion $S_\alpha$\+module $S_\alpha\ot_RK$
is injective for every~$\alpha$.

 The ``only if'' implication holds by
Corollary~\ref{localization-of-injective-torsion-modules-cor}.
 To prove the ``if'', put $\fX=\Spf\Lambda(R)$ and $\fU_\alpha=
\Spf\Lambda(S_\alpha)$.
 Then $\fX=\bigcup_\alpha\fU_\alpha$ is an affine open covering
of the affine Noetherian formal scheme~$\fX$.
 Consider the quasi-coherent torsion sheaf $\cJ=\widetilde K$ on~$\fX$.
 We have $\Lambda_{S_\alpha}(S_\alpha)=\Lambda_R(S_\alpha)$ and
$\Lambda_R(S_\alpha)\ot_RK\simeq S_\alpha\ot_RK$, which is
an injective $S_\alpha$\+module by assumption.
 In view of Lemma~\ref{injective-torsion-modules-characterizations},
\,$\Lambda_{S_\alpha}(S_\alpha)\ot_RK$ is also an injective module
over $\Lambda_{S_\alpha}(S_\alpha)$.
 So the restriction of $\cJ$ to $\fU_\alpha$ belongs to
$\fU_\alpha\Tors^\inj$ for every~$\alpha$.
 By Corollary~\ref{injective-torsion-local-flasque}(d), it follows
that $\cJ\in\fX\Tors^\inj$.
 Thus $K=\cJ(\fX)$ is injective as a $\Lambda(R)$\+module, hence
also as an $R$\+module.
\end{proof}

\begin{cor} \label{qcoh-tors-direct-sums-of-injectives-cor}
 For any locally Noetherian formal scheme\/ $\fX$, the full
subcategory of injective objects\/ $\fX\Tors^\inj$ is closed under
infinite direct sums in\/ $\fX\Tors$.
\end{cor}

\begin{proof}
 Easily deduced either from
Corollary~\ref{local-section-criterion-of-injectivity-cor} together
with Lemma~\ref{sums-of-inj-tors-prods-of-flat-contramods}(a), or
directly from Theorem~\ref{loc-Noetherian-injective-qcoh-torsion}(b).
\end{proof}

\subsection{Projective locally cotorsion contraherent cosheaves of
contramodules} \label{projective-lct-cosheaves-subsecn}
 In this section we spell out the classification of projective
locally cotorsion contraherent cosheaves of contramodules over locally
Noetherian formal schemes, based on the classification of flat cotorsion
contramodules over adic Noetherian rings presented in
Proposition~\ref{flat-cotorsion-contramods-classification-prop}.
 The latter, in turn, is based on Enochs' classification of flat
cotorsion modules over Noetherian commutative rings~\cite{En2}.

 The locally Noetherian formal scheme case is not very different from
the locally Noetherian (nonformal) scheme case that was worked out
in~\cite[Section~6.1]{Pcosh}.
 Our exposition follows that in~\cite[Section~6.1]{Pcosh}.
 
 Let $\fX$ be a (not necessarily semi-separated) locally Noetherian
formal scheme with an open covering~$\bW$.
 We keep the notation and setting of the previous
Section~\ref{injective-torsion-sheaves-subsecn}, including
the definition of the adic topologies on the rings $\widehat\fO_{x,\fX}$
as per Remark~\ref{tricky-adic-topology-on-completed-stalk-ring}.
 In the following theorem, we are interested in
$\widehat\fm_{x,\fX}$\+contramodule $\widehat\fO_{x,\fX}$\+modules
(which are not to be confused with the merely contramodule
$\widehat\fO_{x,\fX}$\+modules, which means contramodules with
respect to the ideal of definition of $\fO_\fX(\fU)$ for $x\in\fU
\subset\fX$).

 Given a complete, separated adic Noetherian ring $\fR$ and
a contraadjusted contramodule $\fR$\+module $P\in\fR\Contra^\cta$,
we denote by $\widecheck P\in\Spf\fR\Ctrh$ the contraherent cosheaf
of contramodules on the affine Noetherian formal scheme $\Spf\fR$
corresponding to the contraadjusted contramodule $\fR$\+module $P$
under the equivalence of categories $\fR\Contra^\cta\simeq\Spf\fR\Ctrh$
from Lemma~\ref{contrah-cosheaves-of-contramods-lemma}(b).

\begin{thm} \label{loc-Noetherian-projective-lct-cosh-contramods}
\textup{(a)} There are enough projective objects in the exact categories
of locally cotorsion locally contraherent cosheaves of contramodules\/
$\fX\Lcth^\lct_\bW$ and\/ $\fX\Lcth^\lct$.
 All of these projective objects belong to the full subcategory of
locally cotorsion contraherent cosheaves of contramodules\/
$\fX\Ctrh^\lct$.
 The full subcategories of projective objects in the three exact
categories\/ $\fX\Ctrh^\lct\subset\fX\Lcth^\lct_\bW\subset\fX\Lcth^\lct$
coincide. \par
\textup{(b)} For any point $x\in\fX$, denote by $\iota_x\:
\Spf\widehat\fO_{x,\fX}\rarrow\fX$ the natural morphism of formal
schemes.
 Then a locally cotorsion contraherent cosheaf of contramodules\/ $\fF$
on\/ $\fX$ is projective if and only if it is isomorphic to an infinite
direct product\/ $\prod_{x\in\fX}\iota_x{}_!\widecheck F_x$ of
the direct images $\iota_x{}_!\widecheck F_x$ of contraherent cosheaves
of contramodules $\widecheck F_x$ on\/ $\Spf\widehat\fO_{x,\fX}$
corresponding to some free/projective/flat\/
$\widehat\fm_{x,\fX}$\+contramodule\/
$\widehat\fO_{x,\fX}$\+modules~$F_x$
(i.~e., $F_x$ is a projective object of the abelian category\/
$\widehat\fO_{x,\fX}\Modl_{\widehat\fm_{x,\fX}\dctra}$).
\end{thm}

\begin{proof}
 For a semi-separated Noetherian formal scheme $\fX$, all
the assertions of part~(a) were proved in
Section~\ref{semi-separated-inj-proj-subsecn};
see Lemma~\ref{projective-contrah-lct-on-qcomp-ssep}
and Corollary~\ref{projective-contrah-lct-on-qcomp-ssep-corollary}.
 The argument below, based on
Proposition~\ref{flat-cotorsion-contramods-classification-prop},
proves parts~(a) and~(b) simultaneously in the general case of
a locally Noetherian formal scheme~$\fX$.

 We recall that the classes of free, projective, and flat
$\widehat\fm_{x,\fX}$\+contramodule
$\widehat\fO_{x,\fX}$\+mod\-ules $F_x$ coincide by
Corollary~\ref{maximal-ideal-in-Noetherian-ring-cor}(b)
and Lemma~\ref{projective-contramods-are-free-lemma}.
 All $\widehat\fm_{x,\fX}$\+con\-tra\-mod\-ule
$\widehat\fO_{x,\fX}$\+mod\-ules are cotorsion
$\widehat\fO_{x,\fX}$\+modules by
Corollary~\ref{maximal-ideal-in-Noetherian-ring-cor}(a).
 By Lemma~\ref{lct-contrah-cosheaves-of-contramods-on-affine-lemma},
it follows that $\widecheck P_x$ is a locally cotorsion contraherent
cosheaf of contramodules on $\Spf\widehat\fO_{x,\fX}$ for any
$\widehat\fm_{x,\fX}$\+contramodule
$\widehat\fO_{x,\fX}$\+module~$P_x$.

 However, following the beginning of the proof of
Theorem~\ref{loc-Noetherian-injective-qcoh-torsion},
the morphism~$\iota_x$ is quasi-compact but \emph{not} necessarily
affine.
 So Lemma~\ref{lcth-contram-W-T-affine-direct-image-lemma}(b)
and formula~\eqref{ctrh-lct-contramod-affine-direct-image} are
\emph{not} applicable to~$\iota_x$, and the assertion that
$\iota_x{}_!\widecheck F_x$ is a locally cotorsion contraherent
cosheaf of contramodules on $\fX$ requires a separate proof.

 Let us prove a slightly more general statement: the cosheaf of
$\fO_\fX$\+modules $\iota_x{}_!\widecheck P_x$ is a locally cotorsion
contraherent cosheaf of contramodules on $\fX$ for any
$\widehat\fm_{x,\fX}$\+contramodule
$\widehat\fO_{x,\fX}$\+module~$P_x$.
 It suffices to check that the restriction of
$\iota_x{}_!\widecheck P_x$ to every affine open formal subscheme
$\fU\subset\fX$ is a locally cotorsion contraherent cosheaf of
contramodules on~$\fU$.
 If the point $x\in\fX$ belongs to $\fU$, then the morphism~$\iota_x$
factorizes through $\fU$, so we have a natural morphism of formal
schemes $\kappa_x\:\Spf\widehat\fO_{x,\fX}\rarrow\fU$.
 Clearly, $(\iota_x{}_!\widecheck P_x)|_\fU\simeq
\kappa_x{}_!\widecheck P_x$.
 The morphism~$\kappa_x$ is affine, so $\kappa_x{}_!\widecheck P_x$
is a locally cotorsion contraherent cosheaf on $\fU$ by
Lemma~\ref{lcth-contram-W-T-affine-direct-image-lemma}(b)
and formula~\eqref{ctrh-lct-contramod-affine-direct-image}.

 Let us show that the restriction of the cosheaf
$\iota_x{}_!\widecheck P_x$ on $\fX$ to an affine open formal subscheme
$\fU\subset\fX$ vanishes for every $\widehat\fm_{x,\fX}$\+contramodule
$\widehat\fO_{x,\fX}$\+module $P_x$ whenever $x\notin\fU$.
 Indeed, it suffices to check that the module of cosections
$\widecheck P_x[\fV]$ vanishes for every principal affine open formal
subscheme $\fV\subset\Spf\widehat\fO_{x,\fX}$ that does not contain
the closed point of $\Spf\widehat\fO_{x,\fX}$.
 So we need to check that $\Hom_{\widehat\fO_{x,\fX}}
(\widehat\fO_{x,\fX}[s^{-1}],P_x)=0$ for all $s\in\widehat\fm_{x,\fX}$.
 This holds because $P_x$ is
an $\widehat\fm_{x,\fX}$\+contramodule $\widehat\fO_{x,\fX}$\+module.

 It follows that cosheaves of contramodule $\fO_X$\+modules of
the form $\prod_{x\in\fX}\iota_x{}_!\widecheck P_x$ on $\fX$,
where $P_x$ are any $\widehat\fm_{x,\fX}$\+contramodule
$\widehat\fO_{x,\fX}$\+modules, are contraherent and locally cotorsion.
 Let us show that cosheaves
$\fF=\prod_{x\in\fX}\iota_x{}_!\widecheck F_x$, where $F_x$ denote some
free $\widehat\fm_{x,\fX}$\+contramodule $\widehat\fO_{x,\fX}$\+modules,
are projective objects in the exact category $\fX\Lcth^\lct_\bW$.

 Pick an affine open covering $\fX=\bigcup_\alpha\fU_\alpha$ of
the formal scheme $\fX$ subordinate to $\bW$, and choose a well-ordering
of the set of indices~$\{\alpha\}$.
 As in the proof of Theorem~\ref{loc-Noetherian-injective-qcoh-torsion},
we put $\fS_\alpha=\fU_\alpha\setminus\bigcup_{\gamma<\alpha}\fU_\gamma$
(so $\fS_\alpha$ is a subset of the underlying topological space
of~$\fX$).
 For any quasi-compact open formal subscheme $\fU\subset\fX$, one
has $\fU\cap\fS_\alpha=\varnothing$ for all but a finite subset
of indices~$\alpha$.
 
 Put $\fF_\alpha=\prod_{z\in\fS_\alpha}\iota_z{}_!\widecheck F_z\in
\fX\Ctrh^\lct$.
 Clearly, $\fF\simeq\prod_\alpha\fF_\alpha$.
 It follows from the previous two paragraphs that, for any affine open
formal subscheme $\fU\subset\fX$, one has $\fF_\alpha[\fU]=0$ for all
but a finite subset of indices~$\alpha$.
 Thus the cosheaf $\fF$ is also the coproduct of
the cosheaves~$\fF_\alpha$ (in the category $(\fX,\fO_X)\Cosh$ of
cosheaves of $\fO_\fX$\+modules on~$\fX$).
 Therefore, it suffices to check that each cosheaf $\fF_\alpha$ is
a projective object in $\fX\Lcth^\lct_\bW$.

 Denote by~$\fj_\alpha$ the open immersion $\fU_\alpha\rarrow\fX$ and
by~$\kappa_z$ the natural morphisms $\Spf\widehat\fO_{z,\fX}\rarrow
\fU_\alpha$ for $z\in\fS_\alpha$.
 Then we have $\fF_\alpha=\fj_\alpha{}_!\fG_\alpha$, where
$\fG_\alpha=\prod_{z\in\fS_\alpha}\kappa_z{}_!\widecheck F_z$.
 The observation that the direct images with respect to quasi-compact
morphisms of locally Noetherian formal schemes preserve infinite
products of cosheaves of $\fO$\+modules (see
Section~\ref{direct-images-of-ctrh-contramod-subsecn}) is used here.
 By Proposition~\ref{flat-cotorsion-contramods-classification-prop},
\,$\fG_\alpha[\fU_\alpha]$ is a flat cotorsion contramodule
$\fO_\fX(\fU_\alpha)$\+module, hence $\fG_\alpha$ is projective
locally cotorsion contraherent cosheaf of contramodules on~$\fU_\alpha$
(see the proof of Lemma~\ref{projective-contrah-lct-on-qcomp-ssep} and
Lemma~\ref{cotorsion-pair-kernel-as-injectives-projectives}(b) below).
 By the adjunction~\eqref{cosheaves-direct-inverse-adjunction}, it
follows that $\fF_\alpha$ is a projective locally cotorsion
$\bW$\+locally contraherent cosheaf of contramodules on~$\fX$.

 Now let us construct for any locally cotorsion $\bW$\+locally
contraherent cosheaf of contramodules $\fQ$ on $\fX$ an admissible
epimorphism in $\fX\Lcth^\lct_\bW$ onto $\fQ$ from a cosheaf
$\fF=\prod_{x\in\fX}\iota_x{}_!\widecheck F_x$ as in part~(b) of
the theorem.
 For this purpose, we choose an affine open covering $\fU_\alpha$ as
above, and proceed by transfinite induction on~$\alpha$, constructing
locally cotorsion contraherent cosheaves of contramodules
$\fF_\alpha=\prod_{z\in\fS_\alpha}\iota_z{}_!\widecheck F_z$
and morphisms of locally cotorsion locally contraherent cosheaves of
contramodules $\fF_\alpha\rarrow\fQ$.

 Suppose such cosheaves and morphisms have been constructed for all
indices $\alpha<\beta$, with some fixed index~$\beta$.
 Then, as explained above, there is the induced morphism of locally
cotorsion locally contraherent cosheaves of contramodules
$\prod_{\alpha<\beta}\fF_\alpha=\coprod_{\alpha<\beta}\fF_\alpha
\rarrow\fQ$.
 Assume that the related morphism of the modules of cosections
$\prod_{\alpha<\beta}\fF_\alpha[\fU]=
\bigoplus_{\alpha<\beta}\fF_\alpha[\fU]\rarrow\fQ[\fU]$ is an admissible
epimorphism of in the exact category of cotorsion contramodule
$\fO_\fX(\fU)$\+modules for every affine open formal subscheme
$\fU\subset\bigcup_{\alpha<\beta}\fU_\alpha\subset\fX$ subordinate
to~$\bW$.
 Our aim is to construct a contraherent cosheaf of contramodules
$\fF_\beta=\bigoplus_{z\in\fS_\beta}\iota_z{}_!\widecheck F_z$ and
a morphism of locally contraherent cosheaves of contramodules
$\fF_\beta\rarrow\fQ$ such that the resulting morphism
$\prod_{\alpha\le\beta}\fF_\alpha\rarrow\fQ$ satisfies the same
condition as above for every affine open formal subscheme
$\fU\subset\bigcup_{\alpha\le\beta}\fU_\alpha\subset\fX$
subordinate to~$\bW$.

 Choose an admissible epimorphism in $\fU_\beta\Ctrh^\lct$ onto
the locally cotorsion contraherent cosheaf of contramodules
$\fj_\beta^!\fQ$ from a projective locally cotorsion contraherent
cosheaf of contramodules $\fG$ on~$\fU_\beta$.
 Then $\fG[\fU_\beta]$ is a flat cotorsion contramodule
$\fO_{\fU_\beta}(\fX)$\+module.
 By Proposition~\ref{flat-cotorsion-contramods-classification-prop},
the cosheaf $\fG$ decomposes into a direct product
$\fG\simeq\prod_{x\in\fU_\beta}\lambda_x{}_!\widecheck G_x$,
where $G_x$ are some free $\widehat\fm_{x,\fX}$\+contramodule
$\widehat\fO_{x,\fX}$\+modules and $\lambda_x\:\Spf\widehat\fO_{x,\fX}
\rarrow\fU_\beta$ are the natural morphisms of formal schemes.
 Let us split this infinite product into two summands as
$\fG=\prod_{z\in\fS_\beta}\lambda_z{}_!\widecheck G_z\oplus
\prod_{y\in\fU_\beta\setminus\fS_\beta}\lambda_y{}_!\widecheck G_y$.
 Denote the former direct summand by $\fG_\beta$ and the latter one
by~$\fE$.
 Set $\fH=\prod_{\alpha<\beta}\fF_\alpha\in\fX\Ctrh^\lct$ and
$\fF_\beta=\fj_\beta{}_!\fG_\beta\in\fX\Ctrh^\lct$.

 The property of a morphism of locally cotorsion $\bW$\+locally
contraherent cosheaves of contramodules to be an admissible epimorphism
is local; the morphism $\fH\rarrow\fQ$ becomes a an admissible
epimorphism when restricted to $\bigcup_{\alpha<\beta}\fU_\alpha$ by
assumption; so we only need to show that the natural morphism
$\fH\oplus\fF_\beta\rarrow\fQ$ becomes an admissible epimorphism
in $\fU_\beta\Ctrh^\lct$ when restricted to~$\fU_\beta$.
 This means that the morphism $\fj_\beta^!\fH\oplus\fG_\beta\rarrow
\fj_\beta^!\fQ$ should be an admissible epimorphism in
$\fU_\beta\Ctrh^\lct$.
 It suffices to check that the admissible epimorphism $\fG\rarrow
\fj_\beta^!\fQ$ factorizes through the morphism $\fj_\beta^!\fH\oplus
\fG_\beta\rarrow\fj_\beta^!\fQ$ (see~\cite[Proposition~7.6]{Bueh}),
or that the morphism $\fE\rarrow\fj_\beta^!\fQ$ factorizes through
the morphism $\fj_\beta^!\fH\rarrow\fj_\beta^!\fQ$.

 Put $\fV_\beta=\fU_\beta\setminus\fS_\beta=\fU_\beta\cap
\bigcup_{\alpha<\beta}\fU_\alpha$.
 So $\fV_\beta$ is a Noetherian formal scheme and an open formal
subscheme in~$\fX$.
 Denote by $\fj$ the identity open immersion $\fV_\alpha\rarrow
\fU_\beta$.
 Then we have $\fE=\fj_!\fL$, where $\fL$ is a projective object in
$\fV_\beta\Ctrh^\lct$ (as we have already proved above).
 Hence the morphism $\fL\rarrow\fj^!\fj_\beta^!\fQ$ factorizes through
the admissible epimorphism of locally cotorsion contraherent cosheaves
of contramodules $\fj^!\fj_\beta^!\fH\rarrow\fj^!\fj_\beta^!\fQ$.
 Therefore, the morphism $\fj_!\fL\rarrow\fj_\beta^!\fQ$ factorizes
through the morphism $\fj_\beta^!\fH\rarrow\fj_\beta^!\fQ$, as desired.
 Here we are using the assumption that the morphism $\fH\rarrow\fQ$
is an admissible epimorphism of locally cotorsion locally contraherent
cosheaves of contramodules in restriction to
$\bigcup_{\alpha<\beta}\fU_\alpha\supset\fV_\beta$, as well as
the adjunction
isomorphism~\eqref{cosheaves-direct-image-restriction-adjunction}.

 We have proved part~(a) and almost proved part~(b).
 In order to finish the proof of part~(b), it remains to show that
the class of contraherent cosheaves of contramodules of the form
$\fF=\prod_{x\in\fX}\iota_x{}_!\widecheck F_x$ on $\fX$ is closed
under direct summands.
 The argument is based on
Lemma~\ref{infinite-coproduct-direct-summand-lemma}
together with the following lemma.

\begin{lem} \label{cosheaves-lower-upper-sets-fully-faithful-lemma}
\textup{(a)} Suppose the set of all points (of the underlying
topological space) of\/ $\fX$ is presented as a union of two
nonintersecting subsets\/ $\fX=\fS\sqcup\fT$ such that for any points
$z\in\fS$ and $y\in\fT$ the closure of~$z$ in\/ $\fX$ does not
contain~$y$.
 Then, for any cosheaves of\/ $\fO$\+modules\/ $\fP_y$ over\/
$\Spf\widehat\fO_{y,\fX}$ and any\/
$\widehat\fm_{z,\fX}$\+contramodule\/ $\widehat\fO_{z,\fX}$\+modules\/
$P_z$ one has\/ $\Hom^{\fO_\fX}(\prod_{y\in\fT}\iota_y{}_!\fP_y,\>
\prod_{z\in\fS}\iota_z{}_!\widecheck P_z)=0$. \par
\textup{(b)} For any point $x\in\fX$, the functor assigning to
an\/ $\widehat\fm_{x,\fX}$\+contramodule\/ $\widehat\fO_{x,\fX}$\+module
$P_x$ the locally cotorsion contraherent cosheaf of contramodules\/
$\iota_x{}_!\widecheck P_x$ on\/ $\fX$ is fully faithful.
\end{lem}

\begin{proof}
 Part~(a): notice that the infinite products in $\fX\Ctrh^\lct$ agree
with the ones in $\fX\Ctrh$ and $(\fX,\fO_\fX)\Cosh$ (see
Sections~\ref{locally-contraherent-cosheaves-of-contramods-subsecn}
and~\ref{locally-contraherent-lct-cosheaves-subsecn}).
 By the definition of the infinite product, it suffices to check
that $\Hom^{\fO_\fX}(\prod_{y\in\fT}\iota_y{}_!\fP_y,\>
\iota_z{}_!\widecheck P_z)=0$ for every $z\in\fS$.
 As in the proof of
Lemma~\ref{sheaves-lower-upper-sets-fully-faithful-lemma}(a),
we consider the closure $\fZ$ of the point~$z$ in $\fX$ and let
$\fY=\fX\setminus\fZ$ be its open complement.
 Denote the related open immersion morphism of formal schemes
by $\fj\:\fY\rarrow\fX$.
 For every point $y\in\fT$, denote the natural morphism of
formal schemes by $\kappa_y\:\Spf\widehat\fO_{y,\fX}\rarrow\fY$;
so $\iota_y=\fj\circ\kappa_y$.

 Now we have $\prod_{y\in\fT}\iota_y{}_!\fP_y\simeq
\fj_!\prod_{y\in\fT}\kappa_y{}_!\fP_y$.
 According to
the adjunction~\eqref{cosheaves-direct-inverse-adjunction},
there is an isomorphism of abelian groups
$$
 \Hom^{\fO_\fX}\Bigl(\fj_!\prod\nolimits_{y\in\fT}\kappa_y{}_!\fP_y,\>
 \iota_z{}_!\widecheck P_z\Bigr)\simeq
 \Hom^{\fO_\fY}\Bigl(\prod\nolimits_{y\in\fT}\kappa_y{}_!\fP_y,\>
 \fj^!\iota_z{}_!\widecheck P_z\Bigr).
$$
 It was shown above in the proof of the theorem that
$\fj^!\iota_z{}_!\widecheck P_z=0$, so the assertion of part~(a)
follows.

 Part~(b): the morphism of formal schemes~$\iota_x$ is tight and flat
(and coaffine for the trivial open coverings
$\{\fX\}$ and $\{\Spf\widehat\fO_{x,\fX}\}$
of $\fX$ and $\Spf\widehat\fO_{x,\fX}$)
by Lemma~\ref{stalk-ring-map-tight-and-flat}, so
the adjunction~\eqref{cosheaves-direct-inverse-adjunction} is
applicable, and it suffices to show that the adjunction morphism
$\widecheck P_x\rarrow\iota_x^!\iota_x{}_!\widecheck P_x$ is
an isomorphism in $\Spf\widehat\fO_{x,\fX}\Ctrh^\lct$.
 One can replace the formal scheme $\fX$ by any affine open formal
subscheme $\fU\subset\fX$ containing~$x$; then we need to show
that the natural map of $\widehat\fO_{x,\fX}$\+modules
$P_x\rarrow\Hom_{\fO_\fU(\fU)}(\widehat\fO_{x,\fX},P_x)$ is
an isomorphism.
 Indeed, we have
$$
 \Hom_{\fO_\fU(\fU)}(\widehat\fO_{x,\fX},P_x)\simeq
 \Hom_{\widehat\fO_{x,\fX}}
 \bigl(\widehat\fO_{x,\fX}\ot_{\fO_\fU(\fU)}\widehat\fO_{x,\fX},
 \>P_x\bigr),
$$
while $P_x\simeq\Hom_{\widehat\fO_{x,\fX}}(\widehat\fO_{x,\fX},P_x)$,
and it remains to use
Lemma~\ref{stalk-completion-of-multiplication-lemma}(b).
\end{proof}

 Now we are ready to finish the proof of the theorem.
 Firstly, consider the case of a Noetherian (i.~e., quasi-compact)
formal scheme~$\fX$.
 Following the proofs of
Theorem~\ref{loc-Noetherian-injective-qcoh-torsion}
and~\cite[Theorem~6.1.1]{Pcosh}, there exists
a well-ordering relation~$<$ on the set of points of $\fX$ such that
$x<y$ whenever $x$ belongs to the closure of~$y$.
 Now we apply
Lemma~\ref{infinite-coproduct-direct-summand-lemma} to
the opposite category $\sA$ to the full additive subcategory of objects
of the form $\fF=\prod_{x\in\fX}\iota_x{}_!\widecheck F_x$ in
$\fX\Ctrh^\lct$, and claim that the images of all idempotent
endomorphisms exist in~$\sA$.
 Indeed, the first assumption of
Lemma~\ref{infinite-coproduct-direct-summand-lemma} holds by
Lemma~\ref{cosheaves-lower-upper-sets-fully-faithful-lemma}(a),
while the second assumption follows from
Lemma~\ref{cosheaves-lower-upper-sets-fully-faithful-lemma}(b).

 For an arbitrary locally Noetherian formal scheme $\fX$, we consider
the \emph{direct sum} decomposition $\fF=\prod_\alpha\fF_\alpha
=\bigoplus_\alpha\fF_\alpha$ of the object $\fF=\prod_{x\in\fX}
\iota_x{}_!\widecheck F_x$ in the category of cosheaves of
$\fO_\fX$\+modules, as per the discussion above in the beginning of
the proof of the theorem.
 Here $\fX=\bigcup_\alpha\fU_\alpha$ is an affine open covering of
$\fX$ indexed by a well-ordered set of indices~$\{\alpha\}$.
 We apply Lemma~\ref{infinite-coproduct-direct-summand-lemma} to
the to the full additive subcategory $\sA$ of objects of the form
$\prod_{x\in\fX}\iota_x{}_!\widecheck F_x$
in the category $(\fX,\fO_\fX)\Cosh$, and claim that the images of
idempotent endomorphisms exist in~$\sA$.
 Once again, the first assumption of
Lemma~\ref{infinite-coproduct-direct-summand-lemma} holds by
Lemma~\ref{cosheaves-lower-upper-sets-fully-faithful-lemma}(a).
 To check the second assumption, recall that $\fF_\alpha\simeq
\fj_\alpha{}_!\fG_\alpha$, where $\fG_\alpha=\prod_{z\in\fS_\alpha}
\kappa_z{}_!\widecheck F_z$.
 The functor~$\fj_\alpha{}_!$ is fully faithful (as one can see from
the adjunction~\eqref{cosheaves-direct-image-restriction-adjunction}
or~\eqref{cosheaves-direct-inverse-adjunction}), so the images of
all the idempotent endomorphisms of the cosheaves $\fF_\alpha$ have
the desired form by the argument in the previous paragraph applied to
the Noetherian affine formal scheme $\fU_\alpha$ in the role of~$\fX$.
\end{proof}

 Recall the notation $\fX\Ctrh^\lct_\prj\subset\fX\Ctrh^\lct$ from
Section~\ref{semi-separated-inj-proj-subsecn} for the full subcategory
of projective locally cotorsion contraherent cosheaves of contramodules
on~$\fX$.
 According to the following corollary, which is our generalization
of~\cite[Corollary~6.1.4]{Pcosh} and a dual-analogous version of
Corollary~\ref{injective-torsion-local-flasque}, being a projective
locally cotorsion contraherent cosheaf of contramodules on a locally
Noetherian formal scheme is a local property.
 Furthermore, all such cosheaves are coflasque.

\begin{cor} \label{projective-lct-of-contramods-local-coflasque}
 Let\/ $\fX$ be a locally Noetherian formal scheme. \par
\textup{(a)} Let\/ $\fY\subset\fX$ be an open formal subscheme.
 Then for any cosheaf\/ $\fF\in\fX\Ctrh^\lct_\prj$ the cosheaf\/
$\fF|_\fY$ belongs to\/ $\fY\Ctrh^\lct_\prj$. \par
\textup{(b)} In the context of part~\textup{(a)}, the corestriction
map\/ $\fF[\fY]\rarrow\fF[\fX]$ is injective (so the cosheaf\/ $\fF$
is coflasque). \par
\textup{(c)} For an affine Noetherian formal scheme\/ $\fX$,
a (not necessarily affine) open formal subscheme\/ $\fY\subset\fX$,
and a cosheaf\/ $\fF\in\fX\Ctrh^\lct_\prj$, the corestriction sequence
of cosections\/ $0\rarrow\fF[\fY]\rarrow\fF[\fX]\rarrow\fF[\fX]/\fF[\fY]
\rarrow0$ is a split short exact sequence of flat cotorsion
contramodule\/ $\fO_\fX(\fX)$\+modules. \par
\textup{(d)} Let\/ $\fX=\bigcup_\alpha\fY_\alpha$ be an open covering
of a locally Noetherian formal scheme.
 Then a locally contraherent cosheaf of contramodules on\/ $\fX$
belongs to\/ $\fX\Ctrh^\lct_\prj$ if and only if its restrictions
to\/ $\fY_\alpha$ belong to\/ $\fY_\alpha\Ctrh^\lct_\prj$
for all\/~$\alpha$.
\end{cor}

\begin{proof}
 Part~(a): denote the open immersion morphism by $\fj\:\fY\rarrow\fX$.
 Then by Theorem~\ref{loc-Noetherian-projective-lct-cosh-contramods}(b)
we have $\fj^!\fF\simeq\fj^!\bigl(\prod_{x\in\fX}\iota_x{}_!
\widecheck F_x\bigr)\simeq
\prod_{x\in\fX}\fj^!\iota_x{}_!\widecheck F_x$.
 Furthermore, it was explained in the proof of the same theorem that
$\fj^!\iota_z{}_!\widecheck F_z=0$ for all $z\in\fX\setminus\fY$.
 Given a point $y\in\fY$, denote by $\kappa_y\:\Spf\widehat\fO_{y,\fX}
\rarrow\fY$ the natural morphism of formal schemes.
 Then we have $\fj^!\iota_y{}_!\widecheck F_y\simeq
\kappa_y{}_!\widecheck F_y$.
 Hence the isomorphism $\fj^!\fF\simeq\prod_{y\in\fY}
\kappa_y{}_!\widecheck F_y$, implying that
$\fj^!\fF\in\fY\Ctrh^\lct_\prj$ according to the same theorem.

 Part~(b): by the definition, we have
$(\iota_x{}_!\widecheck F_x)[\fX]\simeq F_x$.
 Since the cosections over a quasi-compact open formal subscheme
preserve infinite products of contraherent cosheaves of contramodules,
we conclude from the computation above that $\fF[\fU]\simeq
\prod_{x\in\fU}F_x$ for any quasi-compact open formal subscheme
$\fU\subset\fX$.
 Hence, once the isomorphism of contraherent cosheaves of contramodules
$\fF\simeq\prod_{x\in\fX}\iota_x{}_!\widecheck F_x$ has been chosen,
the $\fO_\fX(\fU)$\+module $\fF[\fY\cap\fU]$ becomes naturally
a direct summand of the $\fO_\fX(\fU)$\+module $\fF[\fU]$.
 On the other hand, the cosheaf
axiom~\eqref{topology-base-cosheaf-axiom} implies the isomorphism
$\fF[\fX]\simeq\varinjlim_{\fU\subset\fX}\fF[\fU]$, with the filtered
inductive limit taked over all quasi-compact open formal subschemes
$\fU\subset\fX$; and similarly for~$\fF[\fY]$.
 Thus the corestriction map of $\fO_\fX(\fX)$\+modules $\fF[\fY]\rarrow
\fF[\fX]$ is the inclusion of a direct summand.
 We have also proved part~(c) along the way.

 Part~(d): the ``only if'' assertion holds by part~(a).
 To prove the ``if'', we proceed by transfinite induction on
any well-ordering of the set of indices~$\alpha$.
 This reduces the question to the following assertion.
 Suppose $\fX$ is represented as the union of two open formal
subschemes, $\fX=\fW\cup\fY$.
 Suppose further that the restriction of a locally contraherent cosheaf
of contramodules $\fF$ on $\fX$ onto $\fW$ is identified with the direct
product $\prod_{w\in\fW}\kappa_w{}_!\widecheck F_w$, where $F_w$ are
some free $\widehat\fm_{w,\fX}$\+contramodule
$\widehat\fO_{w,\fX}$\+modules, while~$\kappa_w$ are the natural
morphisms of formal schemes $\Spf\widehat\fO_{w,\fX}\rarrow\fW$.
 In addition, suppose that the restriction of $\fF$ onto $\fY$ is
a projective locally cotorsion contraherent cosheaf of contramodules.
 Then there exist some free $\widehat\fm_{z,\fX}$\+contramodule
$\widehat\fO_{z,\fX}$\+modules $F_z$ defined for all the points
$z\in\fX\setminus\fW$ and an isomorphism of contraherent cosheaves
of contramodules $\fF\simeq\prod_{x\in\fX}\iota_x{}_!\widecheck F_x$
whose restriction to $\fW$ coincides with the given isomorphism
$\fF|_\fW\simeq\prod_{w\in\fW}\kappa_w{}_!\widecheck F_w$.

 Indeed, by Theorem~\ref{loc-Noetherian-projective-lct-cosh-contramods}
we have $\fF|_\fY\simeq\prod_{y\in\fY}\iota'_y{}_!\widecheck F'_y$,
where $\iota_y'\:\Spf\widehat\fO_{y,\fX}\rarrow\fY$ are the natural
morphisms of formal schemes and $F'_y$ are some free
$\widehat\fm_{y,\fX}$\+contramodule $\widehat\fO_{y,\fX}$\+modules.
 Restricting to the intersection $\fV=\fW\cap\fY$, we obtain
isomorphism of contraherent cosheaves of contramodules
$\prod_{v\in\fV}\kappa'_v{}_!\widecheck F_v\simeq
\prod_{v\in\fV}\kappa'_v{}_!\widecheck F'_v$, where $\kappa'_v\:
\Spf\widehat\fO_{v,\fX}\rarrow\fV$ are the natural morphisms.
 It is clear from the arguments in the second half of the proof of
Theorem~\ref{loc-Noetherian-projective-lct-cosh-contramods}(b) that
such an isomorphism of infinite products induces an isomorphism of
$\widehat\fm_{y,\fX}$\+contramodule $\widehat\fO_{y,\fX}$\+modules
$F_v\simeq F'_v$.
 (Cf.\ the proof of Corollary~\ref{injective-torsion-local-flasque}(d).)
 Let us identify $F'_v$ with $F_v$ using this isomorphism, and
set $F_y=F_y'$ for all $y\in\fX\setminus\fW\subset\fY$.
 Then our isomorphism of infinite direct products can be viewed as
an automorphism~$\phi$ of the contraherent cosheaf of contramodules
$\prod_{v\in\fV}\kappa'_v{}_!\widecheck F_v$ on~$\fV$.
 We would like to show that the cosheaf $\fF$ on $\fX$ is isomorphic to
$\prod_{x\in\fX}\iota_x{}_!\widecheck F_x$ in a way compatible
with the given isomorphism $\fF|_\fW\simeq\prod_{w\in\fW}
\kappa_w{}_!\widecheck F_w$ over~$\fW$.

 Since a cosheaf of $\fO_\fX$\+modules is determined by its restrictions
to $\fW$ and $\fY$ together with the induced isomorphism between
the restrictions of these to the intersection $\fV=\fW\cap\fY$, it
suffices to check that the automorphism~$\phi$ can be extended to
an automorphism of the cosheaf $\prod_{y\in\fY}\iota'_y{}_!
\widecheck F_y$ on~$\fY$.
 The latter cosheaf decomposes into the direct sum of two direct
products: the one taken over the points $v\in\fV$ and the one ranging
over the points $y\in\fY\setminus\fV$.
 So it remains to show that the automorphism~$\phi$ can be lifted to
an automorphism of the cosheaf $\prod_{v\in\fV}\iota'_v{}_!
\widecheck F_v$ on~$\fY$.
 Actually, the rings of endomorphisms of the contraherent cosheaves of
contramodules $\prod_{v\in\fV}\kappa'_v{}_!\widecheck F_v$ on $\fV$
and $\prod_{v\in\fV}\iota'_v{}_!\widecheck F_v$ on $\fY$ are
isomorphic.
 Indeed, the functor of direct image of cosheaves of
$\fO$\+modules~$\fj_!$ with respect to the open immersion
$\fj\:\fV\rarrow\fY$ is fully faithful (in view of
the adjunction~\eqref{cosheaves-direct-image-restriction-adjunction})
and preserves infinite products.
\end{proof}

 We refer to Section~\ref{antilocally-flat-subsecn} for the definitions
of flat and $\bW$\+flat cosheaves of contramodule $\fO_\fX$\+modules,
as well as antilocally flat $\bW$\+locally contraherent cosheaves of
contramodules.

\begin{cor} \label{W-flat-lct-coincide-with-projective-lct}
 Let\/ $\fX$ be a locally Noetherian formal scheme with an open
covering\/~$\bW$.
 Then the classes of projective locally cotorsion contraherent cosheaves
of contramodules and\/ $\bW$\+flat locally cotorsion\/ $\bW$\+locally
contraherent cosheaves of contramodules on\/ $\fX$ coincide.
 In particular, any\/ $\bW$\+flat locally cotorsion\/ $\bW$\+locally
contraherent cosheaf of contramodules on\/ $\fX$ is flat,
contraherent, and antilocally flat.
\end{cor}

\begin{proof}
 This is the formal scheme generalization
of~\cite[Corollary~6.1.5]{Pcosh}.
 All projective locally cotorsion contraherent cosheaves of
contramodules on\/ $\fX$ are flat by
Corollary~\ref{projective-lct-of-contramods-local-coflasque}(a).
 To prove that every $\bW$\+flat locally cotorsion\/ $\bW$\+locally
contraherent cosheaf of contramodules $\fF$ on $\fX$ belongs to
$\fX\Ctrh^\lct_\prj$, choose an affine open covering $\fU_\alpha$
of the formal scheme $\fX$ subordinate to~$\bW$.
 Then the $\fO_\fX(\fU_\alpha)$\+modules $\fF[\fU_\alpha]$ are,
by the assumptions, flat and cotorsion.
 Hence the locally cotorsion contraherent cosheaves of contramodules
$\fF|_{\fU_\alpha}$ are projective, and it remains to apply
Corollary~\ref{projective-lct-of-contramods-local-coflasque}(d).
 Finally, by
Theorem~\ref{loc-Noetherian-projective-lct-cosh-contramods}(a),
any projective object of the exact category $\fX\Ctrh^\lct$ is also
projective in $\fX\Lcth^\lct_\bW$, and any projective object of
the exact category $\fX\Lcth^\lct_\bW$ is antilocally flat as
a $\bW$\+locally contraherent cosheaf of contramodules on $\fX$
by the definition.
\end{proof}

 Now we can prove another theorem promised near the end of
Section~\ref{colocality-of-flatness-of-contramods-subsecn}.

\begin{proof}[Proof of
Theorem~\ref{colocality-of-flat-contramods-if-cotorsion-thm}]
 Let $R$ be an adic Noetherian ring and $f_\alpha\:R\rarrow S_\alpha$
be a formal open covering.
 Let $P$ be a cotorsion contramodule $R$\+module.
 We need to prove that the contramodule $R$\+module $P$ is flat
if and only if the contramodule $S_\alpha$\+module $\Hom_R(S_\alpha,P)$
is flat for every~$\alpha$.

 The ``only if'' implication holds by
Corollary~\ref{colocalization-of-flat-contramodules-cor}.
 To prove the ``if'', put $\fX=\Spf\Lambda(R)$ and $\fU_\alpha=
\Spf\Lambda(S_\alpha)$.
 Then $\fX=\bigcup_\alpha\fU_\alpha$ is an affine open covering
of the affine Noetherian formal scheme~$\fX$.
 Consider the locally cotorsion contraherent cosheaf of contramodules
$\fF=\widecheck P$ on~$\fX$.
 We have $\Lambda_{S_\alpha}(S_\alpha)=\Lambda_R(S_\alpha)$ (since
the ring map~$f_\alpha$ is continuous and tight) and
$\boL_0\Lambda_R(S_\alpha)\simeq\Lambda_R(S_\alpha)$
by Corollary~\ref{flat-reductions-derived-Lambda-is-underived-cor}
(since $f_\alpha$~is a flat map of adic topological rings).
 Hence $\Hom_R(\Lambda_{S_\alpha}(S_\alpha),P)\simeq
\Hom_R(S_\alpha,P)$, which is a flat $S_\alpha$\+module by assumption.
 In view of Lemma~\ref{flat-contramodules-characterizations},
\,$\Hom_R(\Lambda_{S_\alpha}(S_\alpha),P)$ is also a flat
module over $\Lambda_{S_\alpha}(S_\alpha)$.
 So the restriction of the cosheaf $\fF$ to $\fU_\alpha$ belongs
to $\fU_\alpha\Ctrh^\lct_\prj$ for every~$\alpha$.
 By Corollary~\ref{projective-lct-of-contramods-local-coflasque}(d),
it follows that $\fF\in\fX\Ctrh^\lct_\prj$.
 Thus $P=\fF[\fX]$ is flat as a $\Lambda(R)$\+module, hence also
as an $R$\+module.
 (Once again, we refer to the proof of
Lemma~\ref{projective-contrah-lct-on-qcomp-ssep} and
Lemma~\ref{cotorsion-pair-kernel-as-injectives-projectives}(b) below
for the assertion that the projective objects of $\fX\Ctrh^\lct$
correspond precisely to flat cotorsion contramodule
$\Lambda(R)$\+modules under the equivalence of categories from
Lemma~\ref{lct-contrah-cosheaves-of-contramods-on-affine-lemma}.)
\end{proof}

\begin{cor} \label{proj-lct-closed-under-products}
 For any locally Noetherian formal scheme\/ $\fX$, the full subcategory
of projective locally cotorsion contraherent cosheaves of
contramodules\/ $\fX\Ctrh^\lct_\prj$ is closed under infinite
products in\/ $\fX\Ctrh^\lct$ as well as in\/ $\fX\Ctrh$.
\end{cor}

\begin{proof}
 This is the formal scheme generalization
of~\cite[Corollary~6.1.6]{Pcosh}.
 The assertion follows from
Corollary~\ref{W-flat-lct-coincide-with-projective-lct} together with
Lemma~\ref{sums-of-inj-tors-prods-of-flat-contramods}(b).

 Alternatively, the assertion can be deduced directly from
Theorem~\ref{loc-Noetherian-projective-lct-cosh-contramods}(b).
 Then one needs to use the fact that the class of free
$\fm$\+contramodule $R$\+modules with respect to a maximal ideal~$\fm$
of a Noetherian commutative ring $R$ is closed under infinite products
in $R\Modl$.
 A direct reference for a much more general version of this fact
is~\cite[Lemma~1.3.7]{Pweak}.
 Alternatively, an $\fm$\+contramodule $R$\+module is free
as an $\fm$\+contramodule $R$\+module if and only if it is flat
as an $R$\+module (see
Corollary~\ref{maximal-ideal-in-Noetherian-ring-cor}(b)
and Lemma~\ref{projective-contramods-are-free-lemma}), and then
one needs to use the fact that infinite products of flat $R$\+modules
are flat.
\end{proof}

\begin{cor} \label{proj-lct-direct-images}
 Let\/ $\ff\:\fY\rarrow\fX$ be a quasi-compact morphism of locally
Noetherian formal schemes.
 Then \par
\textup{(a)} the functor of direct image of cosheaves of contramodule\/
$\fO$\+modules\/~$\ff_!$ takes projective locally cotorsion contraherent
cosheaves of contramodules on\/ $\fY$ to locally cotorsion contraherent
cosheaves of contramodules on\/~$\fX$; \par
\textup{(b)} if the morphism of formal schemes\/~$\ff$ is tight and
flat, then the functor of direct image of cosheaves of contramodule\/
$\fO$\+modules\/~$\ff_!$ takes projective locally cotorsion contraherent
cosheaves of contramodules on\/ $\fY$ to projective locally cotorsion
contraherent cosheaves of contramodules on\/~$\fX$.
\end{cor}

\begin{proof}
 This is the formal scheme generalization
of~\cite[Corollary~6.1.7]{Pcosh}.

 Part~(a): keep the notation~$\iota_x$ for the natural morphisms of
formal schemes $\Spf\widehat\fO_{x,\fX}\rarrow\fX$, and denote
by~$\kappa_y$ the similar natural morphisms of formal schemes
$\Spf\widehat\fO_{y,\fY}\rarrow\fY$ (where $x\in\fX$ and $y\in\fY$).
 Similarly to the argument in~\cite{Pcosh}, one shows that, for
every locally cotorsion contraherent cosheaf of contramodules of
the form $\fP=\prod_{y\in\fY}\kappa_y{}_!\widecheck P_y$ on $\fY$,
where $P_y$ are some $\widehat\fm_{y,\fY}$\+contramodule
$\widehat\fO_{y,\fY}$\+modules, one has $\ff_!\fP\simeq
\prod_{x\in\fX}\iota_x{}_!\widecheck P_x$, where
$P_x=\prod_{\ff(y)=x}f_{x,y}{}_\sharp P_y$.
 Here $f_{x,y}{}_\sharp=f_{x,y}{}_\#$ are the functors of
contramodule restriction of scalars from
Section~\ref{prelim-change-of-scalars-subsecn} corresponding to
the continuous homomorphisms of complete Noetherian local rings
$f_{x,y}\:\widehat\fO_{x,\fX}\rarrow\widehat\fO_{y,\fY}$
(with the $\widehat\fm_{x,\fX}$\+adic topology on
$\widehat\fO_{x,\fX}$ and the $\widehat\fm_{y,\fY}$\+adic topology
on $\widehat\fO_{y,\fY}$).
 Alternatively, the assertion of part~(a) can be obtained as a special
case of Corollary~\ref{co-flasque-direct-image-cor}(c).

 In part~(b), one has to be careful in that the ring homomorphism
$f_{x,y}\:\widehat\fO_{x,\fX}\rarrow\widehat\fO_{y,\fY}$ is usually
\emph{not} tight with respect to the $\widehat\fm_{x,\fX}$\+adic
topology on $\widehat\fO_{x,\fX}$ and the $\widehat\fm_{y,\fY}$\+adic
topology on $\widehat\fO_{y,\fY}$; so
Lemma~\ref{flat-map-of-adic-rings-direct-image-lemma}(c\+-d) is
\emph{not} applicable.
 Nevertheless, under assumptions of~(b), the ring $\widehat\fO_{y,\fY}$
is a flat module over the ring $\widehat\fO_{x,\fX}$, and one can use
the fact that a $\widehat\fm_{y,\fY}$\+contramodule
$\widehat\fO_{y,\fY}$\+module is a free
$\widehat\fm_{y,\fY}$\+contramodule $\widehat\fO_{y,\fY}$\+module
if and only if it is flat as a $\widehat\fO_{y,\fY}$\+module, and
similarly for~$\fX$
(see Lemma~\ref{projective-contramods-are-free-lemma}
and Corollary~\ref{maximal-ideal-in-Noetherian-ring-cor}(b)
together with Lemma~\ref{flat-contramodules-characterizations}).
 Alternatively, the assertion of part~(b) follows from part~(a) and
the adjunction~\eqref{cosheaves-direct-inverse-adjunction} together
with exactness of the inverse image functor $\ff^!\:\fX\Lcth^\lct\rarrow
\fY\Lcth^\lct$~\,\eqref{inverse-image-lcth-lct-contramod-nocovering}.
\end{proof}

\subsection{Flat contraherent cosheaves of contramodules}
\label{flat-cosheaves-of-contramods-subsecn}
 The notion of a \emph{flat contraherent cosheaf of contramodules}
was defined at the end of Section~\ref{antilocally-flat-subsecn}.
 In the present section, which is the formal scheme generalization
of~\cite[Section~6.2]{Pcosh}, we develop the theory of flat
contraherent cosheaves of contramodules on (locally) Noetherian
formal schemes.
 However, the following additional assumption is needed.

 By the \emph{Krull dimension} of a locally Noetherian formal scheme
$\fX$ we mean the Krull dimension of any closed subscheme of
definition $X\subset\fX$.
 Equivalently, the Krull dimension of an affine Noetherian formal
scheme $\fU$ is defined as the Krull dimension of the quotient
ring of $\fO_\fU(\fU)/\fI$ of the adic Noetherian ring $\fO_\fU(\fU)$
by any its ideal of definition $\fI\subset\fO_\fU(\fU)$.
 Then the Krull dimension of $\fX$ is the supremum of the Krull
dimensions of its affine open formal subschemes $\fU\subset\fX$.

 Most results of this section will be proved under the assumption that
the Krull dimension of $\fX$ is finite.
 On the other hand, \emph{unlike} in
Sections~\ref{semi-separated-inj-proj-subsecn}\+-%
\ref{homology-of-cosheaves-I-subsecn}, the main results of
the present section will be applicable to \emph{non-semi-separated}
Noetherian formal schemes~$\fX$.

 Even in the case of (nonformal) schemes, we obtain a significant
improvement upon the exposition in~\cite[Section~6.2]{Pcosh} in
the present section.
 See Proposition~\ref{finite-krulldim-locally-cotorsion-preenvelope}
below and some of the subsequent corollaries.

 The definition of the \emph{cotorsion dimension} of a contramodule
$R$\+module for an adic Noetherian ring $R$ was given near the end
of Section~\ref{prelim-cotorsion-contramodules-subsecn}.
 Let $\fX$ be a locally Noetherian formal scheme with an open
covering~$\bW$.
 A $\bW$\+locally contraherent cosheaf of contramodules $\fM$ on $\fX$
is said to have \emph{locally cotorsion dimension\/~$\le\nobreak d$}
if the cotorsion dimension of the contramodule $\fO_\fX(\fU)$\+module
$\fM[\fU]$ does not exceed~$d$ for every affine open formal subscheme
$\fU\subset\fX$ subordinate to~$\bW$.

 In view of Corollary~\ref{colocality-of-cotorsion-dimension},
the locally cotorsion dimension of a $\bW$\+locally contraherent
cosheaf of contramodules on $\fX$ does not change when the open
covering $\bW$ is replaced by its refinement.
 If a $\bW$\+locally contraherent cosheaf of contramodules $\fM$ has
a coresolution by locally cotorsion $\bW$\+locally contraherent
cosheaves of contramodules in the exact category $\fX\Lcth_\bW$, then
the locally cotorsion dimension of $\fM$ is equal to the minimal
length of such a coresolution. 

\begin{lem} \label{locally-cotorsion-coresolution-dimension}
 Let\/ $\fX$ be a semi-separated Noetherian formal scheme of Krull
dimension~$D$ with an open covering\/~$\bW$.
 In this context: \par
\textup{(a)} The full subcategory of locally cotorsion\/
$\bW$\+locally contraherent cosheaves of contramodules\/
$\fX\Lcth^\lct_\bW$ is coresolving in the exact category of\/
$\bW$\+locally contraherent cosheaves of contramodules\/ $\fX\Lcth_\bW$.
 The coresolution dimension of any object of\/ $\fX\Lcth_\bW$ with
respect to the coresolving subcategory\/ $\fX\Lcth^\lct_\bW\subset
\fX\Lcth_\bW$ does not exceed~$D$. \par
\textup{(b)} The full subcategory of projective locally cotorsion
contraherent cosheaves of contramodules\/ $\fX\Ctrh^\lct_\prj$ is
coresolving in the exact category of\/ $\bW$\+flat\/ $\bW$\+locally
contraherent cosheaves of contramodules\/ $\fX\Lcth_\bW^\fl$.
 The coresolution dimension of any object of\/ $\fX\Lcth_\bW^\fl$
with respect to the coresolving subcategory\/ $\fX\Ctrh^\lct_\prj
\subset\fX\Lcth_\bW^\fl$ does not exceed~$D$.
\end{lem}

\begin{proof}
 This is the formal scheme generalization of~\cite[Lemma~6.2.1]{Pcosh}.
 Part~(a): the locally cotorsion dimension of any locally contraherent
cosheaf of contramodules on a locally Noetherian formal scheme of Krull
dimension $D$ does not exceed $D$ by
Proposition~\ref{adic-raynaud-gruson-prop}(b).
 For a semi-separated Noetherian formal scheme $\fX$, the full
subcategory $\fX\Lcth^\lct_\bW$ is coresolving in $\fX\Lcth_\bW$
by virtue of Lemma~\ref{locally-cotorsion-preenvelope-lemma}
or Corollary~\ref{antilocally-flat-corollary}(a) together with
the results of Section~\ref{locally-contraherent-lct-cosheaves-subsecn}.
 The assertion now follows from the remark in the paragraph preceding
the lemma.

 Part~(b): we have
$\fX\Ctrh^\lct_\prj=\fX\Lcth^\lct_\bW\cap\fX\Lcth_\bW^\fl$
by Corollary~\ref{W-flat-lct-coincide-with-projective-lct}.
 The full subcategory $\fX\Ctrh^\lct_\prj$ is cogenerating in
$\fX\Lcth_\bW^\fl$ by Corollaries~\ref{antilocally-flat-corollary}(a),
\ref{antilocally-flat-independent-on-covering},
and~\ref{antilocally-flat-are-flat-cor}.
 So the assertions of~(b) follows from those of part~(a) in view of
Lemma~\ref{intersection-co-resolving-lemma}(b) and
Corollary~\ref{intersection-co-resolution-dimension-cor}(b) below.
\end{proof}

 Part~(a) of the following corollary is essentially a special case of
Corollary~\ref{finite-krulldim-flat-cosheaves-characterized}(b) below.
 See Remark~\ref{flat-and-antilocally-flat-remark} for a discussion.
 
\begin{cor} \label{W-flat-implies-flat-and-coflasque-cor}
\textup{(a)} On a semi-separated Noetherian formal scheme\/ $\fX$ of
finite Krull dimension, the classes of\/ $\bW$\+flat\/ $\bW$\+locally
contraherent cosheaves of contramodules and antilocally flat
contraherent cosheaves of contramodules coincide. \par
\textup{(b)} On a locally Noetherian formal scheme\/ $\fX$ of finite
Krull dimension, any\/ $\bW$\+flat\/ $\bW$\+locally contraherent
cosheaf of contramodules is (globally) flat and contraherent.
 So the categories\/ $\fX\Lcth_\bW^\fl$ and\/ $\fX\Ctrh^\fl$ coincide.
\par
\textup{(c)} On a locally Noetherian formal scheme\/ $\fX$ of finite
Krull dimension, any flat contraherent cosheaf of contramodules is
coflasque; so one has\/ $\fX\Ctrh^\fl\subset\fX\Ctrh_\cfq$. \par
\textup{(d)} For an affine Noetherian formal scheme\/ $\fX$ of finite
Krull dimension, a (not necessarily afine) open formal subscheme\/
$\fY\subset\fX$, and a cosheaf\/ $\fF\in\fX\Ctrh^\fl$, the corestriction
sequence of cosections\/ $0\rarrow\fF[\fY]\rarrow\fF[\fX]\rarrow
\fF[\fX]/\fF[\fY]\rarrow0$ is a short exact sequence of flat 
contraadjusted contramodule\/ $\fO_\fX(\fX)$\+modules.
\end{cor}

\begin{proof}
 This is the formal scheme generalization
of~\cite[Corollary~6.2.2]{Pcosh}.

 Part~(a): any antilocally flat contraherent cosheaf of contramodules
on $\fX$ is flat by Corollary~\ref{antilocally-flat-are-flat-cor}
(see also Corollary~\ref{antilocally-flat-independent-on-covering}).
 Conversely, by Lemma~\ref{locally-cotorsion-coresolution-dimension}(b),
any $\bW$\+flat $\bW$\+locally contraherent cosheaf of contramodules
$\fF$ on $\fX$ has a finite coresolution in $\fX\Lcth_\bW$ by
projective locally cotorsion contraherent cosheaves, which are
antilocally flat by the definition.
 Using Corollary~\ref{antilocally-flat-characterized-cor}(b), we
conclude that the cosheaf $\fF$ is antilocally flat.

 Part~(b): by the definition, the flatness and contraherence of locally
contraherent cosheaves of contramodules are checked over affine open
formal subschemes, i.~e., a locally contraherent cosheaf of
contramodules $\fF$ on $\fX$ is flat and contraherent if and only if
the cosheaf $\fF|_\fU$ is flat and contraherent for every affine open
formal subscheme $\fU\subset\fX$.
 On the other hand, for a $\bW$\+flat $\bW$\+locally contraherent
cosheaf of contramodules $\fF$ on $\fX$, the cosheaf of contramodules
$\fF|_\fU$ is $\bW|_\fU$\+flat and $\bW|_\fU$\+locally contraherent,
where $\bW|_\fU=\{\fU\cap\fW\mid\fW\in\bW\}$.
 This reduces the desired assertion to the case of affine Noetherian
formal schemes $\fX$, which is covered by part~(a).

 Part~(c): coflasqueness of cosheaves is a local property
by~\cite[Lemma~3.4.1(a)]{Pcosh} (as mentioned in the beginning of
Section~\ref{coflasque-subsecn}), and flatness of cosheaves of
contramodules is preserved by restictions to open formal subschemes
by the definition.
 So once again it suffices to consider the case of an affine Noetherian
formal scheme~$\fX$.
 Then we know that $\fF$ has a finite coresolution in $\fX\Lcth_\bW$
by cosheaves from $\fX\Ctrh^\lct_\prj$, which are coflasque by
Corollary~\ref{projective-lct-of-contramods-local-coflasque}(b).
 By Lemma~\ref{coflasque-adjusted-to-cosections}(b), it follows that
our coresolution is an exact sequence in $\fX\Ctrh_\cfq$ and
the cosheaf $\fF$ is coflasque.

 Part~(d): continuing with the coresolution from the previous paragraph,
by Lemma~\ref{coflasque-adjusted-to-cosections}(c) it remains exact
after applying the functor of cosections over any open formal subscheme
$\fY\subset\fX$.
 Considering the corestriction maps related to the pair of open formal
subschemes $\fY\subset\fX$, we obtain an injective morphism of exact
sequences of $\fO_\fX(\fX)$\+modules; so the induced sequence of
cokernels is also exact.
 So the sequence of cokernels is a finite exact sequence of
$\fO_\fX(\fX)$\+modules, and all of its terms, except perhaps
the leftmost one, are flat contramodule $\fO_\fX(\fX)$\+modules by
Corollary~\ref{projective-lct-of-contramods-local-coflasque}(c).
 It follows that the leftmost term $\fF[\fX]/\fF[\fY]$ is a flat
contramodule $\fO_\fX(\fX)$\+module, too.
 The $\fO_\fX(\fX)$\+module $\fF[\fY]$ is a contraadjusted contramodule
$\fO_\fX(\fX)$\+module by Corollary~\ref{co-flasque-direct-image-cor}(b)
applied to the open immersion morphism $\fY\rarrow\fX$.
\end{proof}

 Now we can prove the third theorem promised near the end of
Section~\ref{colocality-of-flatness-of-contramods-subsecn}.

\begin{proof}[Proof of
Theorem~\ref{colocality-of-flat-contramods-if-finite-Krull-thm}]
 Let $R$ be an adic Noetherian ring with an ideal of definition
$I\subset R$ such that the quotient ring $R/I$ has finite Krull
dimension.
 Let $f_\alpha\:R\rarrow S_\alpha$ be a formal open covering,
and let $P$ be a contraadjusted contramodule $R$\+module.
 We need to prove that the contramodule $R$\+module $P$ is flat
if and only if the contramodule $S_\alpha$\+module $\Hom_R(S_\alpha,P)$
is flat for every~$\alpha$.

 The ``only if'' implication holds by
Corollary~\ref{colocalization-of-flat-contramodules-cor}.
 To prove the ``if'', put $\fX=\Spf\Lambda(R)$ and $\fU_\alpha=
\Spf\Lambda(S_\alpha)$.
 Then $\fX=\bigcup_\alpha\fU_\alpha$ is an affine open covering
of the affine Noetherian formal scheme~$\fX$.
 Denote this open covering by~$\bW$.
 Consider the contraherent cosheaf of contramodules $\fF=\widecheck P$
on~$\fX$.
 Similarly to the proof of
Theorem~\ref{colocality-of-flat-contramods-if-cotorsion-thm}
given in Section~\ref{projective-lct-cosheaves-subsecn}, one shows
that $\Hom_R(\Lambda_{S_\alpha}(S_\alpha),P)$ is a flat module
over $\Lambda_{S_\alpha}(S_\alpha)$.
 So the contraherent cosheaf of contramodules $\fF$ on $\fX$ is
$\bW$\+flat (we use
Corollary~\ref{colocalization-of-flat-contramodules-cor} again here;
cf.\ the discussion at the end of
Section~\ref{antilocally-flat-subsecn}).
 By Corollary~\ref{W-flat-implies-flat-and-coflasque-cor}(b), it
follows that $\fF$ is a flat contraherent cosheaf of contramodules
on~$\fX$.
 Thus $P=\fF[\fX]$ is flat as a $\Lambda(R)$\+module, hence also
as an $R$\+module.
\end{proof}

 The next lemma is a special case of the more general
Corollary~\ref{flat-contrah-direct-images}(b) below.

\begin{lem} \label{from-semi-separated-flat-direct-image}
 Let\/ $\fX$ be a locally Noetherian formal scheme, $\fY$ be
a semi-separated locally Noetherian formal scheme of finite Krull
dimension, and\/ $\ff\:\fY\rarrow\fX$ be a flat quasi-compact
morphism of formal schemes.
 Then the functor of direct image of cosheaves of\/ $\fO$\+modules\/
$\ff_!\:(\fY,\fO_\fY)\Cosh\rarrow(\fX,\fO_\fX)\Cosh$ takes flat
contraherent cosheaves of contramodules on\/ $\fY$ to flat
contraherent cosheaves of contramodules on\/ $\fX$, and induces
an exact functor\/ $\fY\Ctrh^\fl\rarrow\fX\Ctrh^\fl$.
\end{lem}

\begin{proof}
 This is our version of~\cite[Lemma~6.2.3]{Pcosh}.
 It suffices to consider the case of an affine Noetherian formal
scheme~$\fX$.
 Then the formal scheme $\fY$ is Noetherian (i.~e., quasi-compact)
of finite Krull dimension and semi-separated.
 By Corollary~\ref{W-flat-implies-flat-and-coflasque-cor}(a), any flat
contraherent cosheaf of contramodules on $\fY$ is antilocally flat.
 By Corollary~\ref{ssep-flat-direct-image-of-antilocally-flat-cor},
we have an exact functor of direct image
$\ff_!\:\fY\Ctrh_\alf\rarrow\fX\Ctrh_\alf$.
 By Corollary~\ref{antilocally-flat-are-flat-cor}, any antilocally
flat contraherent cosheaf on $\fX$ is flat.
\end{proof}

\begin{cor} \label{projective-contrah-on-finite-Krull-dim}
 Let\/ $\fX$ be a Noetherian formal scheme of finite Krull dimension
with an open covering\/ $\bW$ and a finite affine open covering\/
$\fX=\bigcup_\alpha\fU_\alpha$ subordinate to\/~$\bW$. 
 In this context: \par
\textup{(a)} There are enough projective objects in the exact categories
of locally contraherent cosheaves of contramodules\/ $\fX\Lcth_\bW$
and\/ $\fX\Lcth$.
 All of these projective objects belong to the full subcategory of
contraherent cosheaves of contramodules\/ $\fX\Ctrh$.
 The full subcategories of projective objects in the three exact
categories\/ $\fX\Ctrh\subset\fX\Lcth_\bW\subset\fX\Lcth$ coincide. \par
\textup{(b)} A contraherent cosheaf of contramodules on\/ $\fX$ is
projective if and only if it is isomorphic to a direct summand of
a finite direct sum of the direct images of projective contraherent
cosheaves of contramodules from\/~$\fU_\alpha$. \par
\textup{(c)} Any projective contraherent cosheaf of contramodules
on\/ $\fX$ is flat (and consequently, coflasque).
\end{cor}

\begin{proof}
 This is the formal scheme generalization
of~\cite[Corollary~6.2.4]{Pcosh}.
 The argument is similar to the proofs of
Lemmas~\ref{projective-contrah-on-qcomp-ssep}
and~\ref{projective-contrah-lct-on-qcomp-ssep}.
 The main difference is that the formal scheme $\fX$ is not necessarily
semi-separated now, so one has to use also
Lemma~\ref{from-semi-separated-flat-direct-image}.
 Let $\fj_\alpha\:\fU_\alpha\rarrow\fX$ denote the open immersion
morphisms, and let $\fF_\alpha$ be projective contraherent cosheaves
of contramodules on~$\fU_\alpha$.
 By Corollary~\ref{antilocally-flat-are-flat-cor} and
Lemma~\ref{from-semi-separated-flat-direct-image}, the cosheaf of
contramodules $\fj_\alpha{}_!\fF_\alpha$ on $\fX$ is flat and
contraherent; and any flat contraherent cosheaf of contramodules
on $\fX$ is coflasque by
Corollary~\ref{W-flat-implies-flat-and-coflasque-cor}(c).
 The adjunction~\eqref{cosheaves-direct-image-restriction-adjunction}
or~\eqref{cosheaves-direct-inverse-adjunction} shows that
$\fj_\alpha{}_!\fF_\alpha$ is a projective object in $\fX\Lcth_\bW$.

 It remains to show that for every object $\fQ\in\fX\Lcth_\bW$ there
exists an admissible epimorphism onto $\fQ$ from an object of
the form $\bigoplus_\alpha\fj_\alpha{}_!\fF_\alpha$.
 For every~$\alpha$, choose an admissible epimorphism onto
the contraherent cosheaf of contramodules $\fj_\alpha^!\fQ$ from
a projective contraherent cosheaf of contramodules $\fF_\alpha$
on~$\fU_\alpha$.
 Then the same adjunction provides a natural morphism of contraherent
cosheaves of contramodules $\bigoplus_\alpha\fj_\alpha{}_!\fF_\alpha
\rarrow\fQ$.
 The latter morphism is an admissible epimorphism of $\bW$\+locally
contraherent cosheaves of contramodules, because it is so in
restriction to each open formal subscheme $\fU_\alpha\subset\fX$
(see the discussion at the end of
Section~\ref{locally-contraherent-cosheaves-of-contramods-subsecn}).
\end{proof}

 Let $\fX$ be a Noetherian formal scheme of finite Krull dimension.
 Recall the notation $\fX\Ctrh_\prj\subset\fX\Ctrh$ from
Section~\ref{semi-separated-inj-proj-subsecn} for the full subcategory
of projective contraherent cosheaves of contramodules on~$\fX$.
 Clearly, $\fX\Ctrh_\prj$ is also the full subcategory
of projective objects in $\fX\Ctrh^\fl$, and there are enough such
projective objects.
 We will see below in Corollary~\ref{flat-proj-lct-preenvelope-cor}(b)
that $\fX\Ctrh^\lct_\prj$ is the full subcategory of \emph{injective}
objects in $\fX\Ctrh^\fl$, and there are enough such injective objects.

 A more general version of the next lemma will be proved below as
Proposition~\ref{finite-krulldim-locally-cotorsion-preenvelope}.

\begin{lem} \label{coflasque-locally-cotorsion-preenvelope}
 Let\/ $\fX$ be a Noetherian formal scheme of finite Krull dimension $D$
with a finite affine open covering\/ $\fX=\bigcup_\alpha\fU_\alpha$.
\par
\textup{(a)} Any coflasque contraherent cosheaf of contramodules\/ $\fE$
on\/ $\fX$ can be included in an (admissible) short exact sequence\/
$0\rarrow\fE\rarrow\fP\rarrow\fF\rarrow0$ in the exact category\/
$\fX\Ctrh_\cfq$, where\/ $\fP$ is a coflasque locally cotorsion
contraherent cosheaf of contramodules on\/ $\fX$ and\/ $\fF$ is
a finitely iterated extension of the direct images of flat contraherent
cosheaves of contramodules from\/~$\fU_\alpha$. \par
\textup{(b)} The full subcategory of coflasque locally cotorsion
contraherent cosheaves of contramodules\/ $\fX\Ctrh^\lct_\cfq$ is
coresolving in the exact category of coflasque contraherent cosheaves
of contramodules\/ $\fX\Ctrh_\cfq$.
 The coresolution dimension of any object of\/ $\fX\Ctrh_\cfq$ with
respect to the coresolving subcategory\/ $\fX\Ctrh^\lct_\cfq$ does not
exceed~$D$.
\end{lem}

\begin{proof}
 This is the formal scheme generalization of~\cite[Lemma~6.2.5]{Pcosh}.

 Part~(a): the proof is similar to that
of Lemma~\ref{locally-cotorsion-preenvelope-lemma}.
 The main difference is that the formal scheme $\fX$ is not necessarily
semi-separated now, so one has to use the coflasqueness assumption and
Corollary~\ref{co-flasque-direct-image-cor}(b\+-c).
 Notice that flat contraherent cosheaves of contramodules
on $\fU_\alpha$ are coflasque by
Corollary~\ref{W-flat-implies-flat-and-coflasque-cor}(c), and
finitely iterated extensions of coflasque contraherent cosheaves of
contramodules in $\fX\Ctrh$ are coflasque by
Lemma~\ref{coflasque-adjusted-to-cosections}(a).

 One proceeds by induction on a linear ordering of the finite set of
indices~$\alpha$, considering the open formal subscheme
$\fV=\bigcup_{\alpha<\beta}\fU_\alpha\subset\fX$.
 Assume that we already have a short exact sequence $0\rarrow\fE
\rarrow\fK\rarrow\fL\rarrow0$ of coflasque contraherent cosheaves of
contramodules on $\fX$ such that the contraherent cosheaf of
contramodules $\fK|_\fV$ is locally cotorsion, while the contraherent
cosheaf of contramodules $\fL$ on $\fX$ is a finitely iterated
extension of the direct images of flat contraherent cosheaves of
contramodules from the affine open formal subschemes $\fU_\alpha$
with $\alpha<\beta$.
 Put $\fU=\fU_\beta$, and let $\fj\:\fU\rarrow\fX$ be the identity
open immersion.

 Choose a short exact sequence $0\rarrow\fj^!\fK\rarrow\fQ\rarrow\fG
\rarrow0$ of contraherent cosheaves of contramodules on $\fU$ such that
$\fQ$ is a locally cotorsion contraherent cosheaf of contramodules
and $\fG$ is a flat contraherent cosheaf of contramodules
(see Proposition~\ref{flat-cotorsion-pair-in-quotseparated}(a)
and Corollary~\ref{colocalization-of-flat-contramodules-cor}).
 Then the cosheaf $\fG$ is coflasque by
Corollary~\ref{W-flat-implies-flat-and-coflasque-cor}(c) and
the cosheaf $\fQ$ is coflasque by 
Lemma~\ref{coflasque-adjusted-to-cosections}(a).

 By Corollary~\ref{co-flasque-direct-image-cor}(b), the related
sequence of direct images $0\rarrow\fj_!\fj^!\fK\rarrow\fj_!\fQ
\rarrow\fj_!\fG\rarrow0$ is a short exact sequence of coflasque
contraherent cosheaves of contramodules on~$\fX$.
 By part~(c) of the latter corollary, the contraherent cosheaf of
contramodules $\fj_!\fQ$ is locally cotorsion.

 Let $0\rarrow\fK\rarrow\fT\rarrow\fj_!\fG\rarrow0$ be the pushout
of the short exact sequence $0\rarrow\fj_!\fj^!\fK\rarrow\fj_!\fQ
\rarrow\fj_!\fG\rarrow0$ with respect to the adjunction morphism
$\fj_!\fj^!\fK\rarrow\fK$ in the exact category $\fX\Ctrh_\cfq$.
 We are going to show that the coflasque contraherent cosheaf of
contramodules $\fT$ on $\fX$ is locally cotorsion in restriction to
$\fU\cup\fV$.

 It suffices to check that both the contraherent cosheaves of
contramodules $\fT|_\fU$ and $\fT|_\fV$ are locally cotorsion.
 In restriction to $\fU$, one has $\fj^!\fT\simeq\fQ$.
 In order to compute the restriction of $\fT$ to $\fV$, denote
by~$\fj'$ the open immersion morphism $\fU\cap\fV\rarrow\fV$.
 Then we have $(\fj_!\fG)|_\fV\simeq\fj'_!(\fG|_{\fU\cap\fV})$.
 The contraherent cosheaf of contramodules $\fK|_{\fU\cap\fV}$ is
locally cotorsion, hence so is the cokernel $\fG|_{\fU\cap\fV}$ of
the admissible monomorphism of locally cotorsion contraherent cosheaves
of contramodules $\fK|_{\fU\cap\fV}\rarrow\fQ|_{\fU\cap\fV}$.

 By Corollary~\ref{co-flasque-direct-image-cor}(c), the direct
image $\fj'_!(\fG|_{\fU\cap\fV})$ is a coflasque locally cotorsion
contraherent cosheaf of contramodules on~$\fV$.
 It follows that in the short exact sequence of contraherent cosheaves
of contramodules $0\rarrow\fK|_\fV\rarrow\fT|_\fV\rarrow(\fj_!\fG)|_\fV
\rarrow0$ the middle term is locally cotorsion, since the two side
terms are.

 Finally, it remains to point out that the composition
$\fE\rarrow\fK\rarrow\fT$ of admissible monomorphisms in $\fX\Ctrh_\cfq$
is again and admissible monomorphism, and its cokernel is an extension
of the flat contraherent cosheaves $\fj_!\fG$ and~$\fL$.

 Part~(b): the full subcategory $\fX\Ctrh^\lct_\cfq$ is cogenerating
in $\fX\Ctrh_\cfq$ by part~(a), so it is also coresolving in view of
the discussion in
Section~\ref{locally-contraherent-lct-cosheaves-subsecn}.
 The bound on the coresolution dimension holds by virtue of
Proposition~\ref{adic-raynaud-gruson-prop}(b).
\end{proof}

\begin{cor} \label{flat-proj-lct-preenvelope-cor}
 Let\/ $\fX$ be a Noetherian formal scheme of finite Krull dimension $D$
with a finite affine open covering\/ $\fX=\bigcup_\alpha\fU_\alpha$.
\par
\textup{(a)} Any flat contraherent cosheaf of contramodules\/ $\fG$
on\/ $\fX$ can be included in a short exact sequence\/
$0\rarrow\fG\rarrow\fP\rarrow\fF\rarrow0$ in the exact category\/
$\fX\Ctrh^\fl$, where\/ $\fP$ is a projective locally cotorsion
contraherent cosheaf of contramodules on\/ $\fX$ and\/ $\fF$ is
a finitely iterated extension of the direct images of flat contraherent
cosheaves of contramodules from\/~$\fU_\alpha$. \par
\textup{(b)} The homological dimension of the exact category\/
$\fX\Ctrh^\fl$ does not exceed~$D$.
 There are both enough projective and enough injective objects in\/
$\fX\Ctrh^\fl$; the full subcategory of projective objects is\/
$\fX\Ctrh_\prj$, and the full subcategory of injective objects\/
is\/ $\fX\Ctrh^\lct_\prj$.
\end{cor}

 In particular, part~(b) of the corollary implies that the resolution
dimensions of all the objects of $\fX\Ctrh^\fl$ with respect to
the resolving subcategory $\fX\Ctrh_\prj$ do not exceed $D$, and
the coresolution dimensions of all the objects of $\fX\Ctrh^\fl$ with
respect to the coresolving subcategory $\fX\Ctrh^\lct_\prj$ do not
exceed~$D$.

\begin{proof}[Proof of
Corollary~\ref{flat-proj-lct-preenvelope-cor}]
 This is the formal scheme generalization
of~\cite[Corollary~6.2.6]{Pcosh}.
 Part~(a): by Corollary~\ref{W-flat-implies-flat-and-coflasque-cor}(c),
the cosheaf $\fG$ is coflasque.
 Lemma~\ref{coflasque-locally-cotorsion-preenvelope}(a)
provides a short exact sequence $0\rarrow\fG\rarrow\fP\rarrow
\fF\rarrow0$ with $\fP\in\fX\Ctrh^\lct_\cfq$ and $\fF$ from the desired
class of finitely iterated extensions of the direct images.
 Now the cosheaf of contramodules $\fF$ is flat by
Lemma~\ref{from-semi-separated-flat-direct-image}, and it follows
that the cosheaf $\fP$ is flat, too.
 Using Corollary~\ref{W-flat-lct-coincide-with-projective-lct},
we conclude that $\fP\in\fX\Ctrh^\lct_\prj$.

 Part~(b): the full subcategory $\fX\Ctrh^\lct_\prj$ is cogenerating
in $\fX\Ctrh^\fl$ by part~(a).
 In view of Corollary~\ref{W-flat-lct-coincide-with-projective-lct}
and the discussion in
Section~\ref{locally-contraherent-lct-cosheaves-subsecn}, it follows
that the full subcategory $\fX\Ctrh^\lct_\prj\subset\fX\Ctrh^\fl$
is coresolving.
 Since $\fX\Ctrh^\lct_\prj$ is a split exact category and it is
closed under direct summands as a full subcategory in $\fX\Ctrh^\fl$,
it follows that $\fX\Ctrh^\lct_\prj$ is the full subcategory of
injective objects in $\fX\Ctrh^\fl$ and there are enough such
injective objects.

 By Corollary~\ref{W-flat-lct-coincide-with-projective-lct}
and Proposition~\ref{adic-raynaud-gruson-prop}(b), the coresolution
dimensions of all the objects of $\fX\Ctrh^\fl$ with respect to
the coresolving subcategory $\fX\Ctrh^\lct_\prj$ do not exceed~$D$.
 The coresolution dimension with respect to the coresolving subcategory
of injective objects is the homological dimension of an exact category
with enough injective objects.
 Finally, the assertions that $\fX\Ctrh_\prj$ is the full subcategory
of projective objects in $\fX\Ctrh^\fl$ and there are enough such
projective objects follow from
Corollary~\ref{projective-contrah-on-finite-Krull-dim}(a,c).
\end{proof}

 A more general version of the next lemma can be found below in
Corollary~\ref{finite-krulldim-flat-precover}.

\begin{lem} \label{coflasque-flat-precover}
 Let\/ $\fX$ be a Noetherian formal scheme of finite Krull dimension
with a finite affine open covering\/ $\fX=\bigcup_\alpha\fU_\alpha$.
\par
\textup{(a)} Any coflasque contraherent cosheaf of contramodules\/ $\fE$
on\/ $\fX$ can be included in a short exact sequence\/
$0\rarrow\fP\rarrow\fF\rarrow\fE\rarrow0$ in the exact category\/
$\fX\Ctrh_\cfq$, where\/ $\fP$ is a coflasque locally cotorsion
contraherent cosheaf of contramodules on\/ $\fX$ and\/ $\fF$ is
a finitely iterated extension of the direct images of flat contraherent
cosheaves of contramodules from\/~$\fU_\alpha$. \par
\textup{(b)} Any flat contraherent cosheaf of contramodules\/ $\fG$
on\/ $\fX$ can be included in a short exact sequence\/
$0\rarrow\fP\rarrow\fF\rarrow\fG\rarrow0$ in the exact category\/
$\fX\Ctrh^\fl$, where\/ $\fP$ is a projective locally cotorsion
contraherent cosheaf of contramodules on\/ $\fX$ and\/ $\fF$ is
a finitely iterated extension of the direct images of flat contraherent
cosheaves of contramodules from\/~$\fU_\alpha$.
\end{lem}

\begin{proof}
 This is the formal scheme generalization of~\cite[Lemma~6.2.7]{Pcosh}
and another application of a suitably general version of
the Salce lemma, see~\cite[Lemma~B.1.1(a)]{Pcosh}.

 Part~(a): according to
Corollary~\ref{projective-contrah-on-finite-Krull-dim}, there exists
an admissible epimorphism in $\fX\Ctrh$ onto any given contraherent
cosheaf of contramodules $\fE$ from a finite direct sum of
the direct images of projective contraherent cosheaves of contamodules
from~$\fU_\alpha$.
 Moreover, any such finite direct sum of the direct images is a flat
and coflasque contraherent cosheaf of contramodules on~$\fX$.
 If the cosheaf $\fE$ is coflasque, then such an admissible epimorphism
in $\fX\Ctrh$ is also an admissible epimorphism in $\fX\Ctrh_\cfq$ 
(by Lemma~\ref{coflasque-adjusted-to-cosections}(b)).
 It remains to apply
Lemma~\ref{coflasque-locally-cotorsion-preenvelope}(a) to the kernel
of our admissible epimorphism and use the construction from
the Salce lemma.

 Part~(b) can be deduced from part~(a) using
Corollaries~\ref{W-flat-implies-flat-and-coflasque-cor}(c)
and~\ref{W-flat-lct-coincide-with-projective-lct}.
 Alternatively, one can prove part~(b) directly in the same way as
part~(a), using Corollary~\ref{flat-proj-lct-preenvelope-cor}(a).
\end{proof}

 The next corollary uses the definitions of the exotic derived
categories, as well as the notation for the homotopy and derived
categories that can be found in
Sections~\ref{appx-exotic-positselski-subsecn}\+-%
\ref{appx-becker-subsecn} below.
 This corollary plays an auxiliary role, to be used in the subsequent
discussion of the groups $\Ext^{\fX,*}({-},{-})$ and in the proof
of Corollary~\ref{flat-contraherent-cosheaves-cor}(a).
 For a more detailed discussion, see
Section~\ref{background-equivalences-subsecn} below.

\begin{cor} \label{minus-ctr-bctr-derived-categories-agree}
\textup{(a)} Let\/ $\fX$ be a Noetherian formal scheme of finite Krull
dimension, and\/ $\st=\bb$, $+$, $-$, $\varnothing$, $\abs+$, $\abs-$,
$\abs$, $\co$, $\ctr$, $\bco$, or\/~$\bctr$ be a conventional or exotic
derived category symbol.
 Then the inclusions of additive/exact catgories\/ $\fX\Ctrh_\prj
\rarrow\fX\Ctrh^\fl$ and\/ $\fX\Ctrh^\lct_\prj\rarrow\fX\Ctrh^\fl$
induce equivalences of triangulated categories\/
$\Hot^\st(\fX\Ctrh_\prj)\simeq\sD^\st(\fX\Ctrh^\fl)$ and\/
$\Hot^\st(\fX\Ctrh^\lct_\prj)\simeq\sD^\st(\fX\Ctrh^\fl)$. \par
\textup{(b)} Let\/ $\fX$ be a locally Noetherian formal scheme with
an open covering\/ $\bW$, and\/ $\st=-$, $\ctr$, or\/~$\bctr$ be
a derived category symbol.
 Then the inclusions of additive/exact categories\/ $\fX\Ctrh^\lct_\prj
\rarrow\fX\Ctrh^\lct_\cfq\rarrow\fX\Ctrh^\lct\rarrow\fX\Lcth^\lct_\bW
\rarrow\fX\Lcth^\lct$ induce triangulated equivalences\/
$\Hot^\st(\fX\Ctrh^\lct_\prj)\simeq\sD^\st(\fX\Ctrh^\lct_\cfq)\simeq
\sD^\st(\fX\Ctrh^\lct)\simeq\sD^\st(\fX\Lcth^\lct_\bW)$, as well as
a triangulated equivalence\/ $\sD^-(\fX\Lcth^\lct_\bW)\simeq
\sD^-(\fX\Lcth^\lct)$. \par
\textup{(c)} Let\/ $\fX$ be a Noetherian formal scheme of finite
Krull dimension with an open covering\/ $\bW$, and\/ $\st=-$, $\ctr$,
or\/~$\bctr$ be a derived category symbol.
 Then the inclusions of additive/exact categories\/ $\fX\Ctrh_\prj)
\rarrow\fX\Ctrh^\fl\rarrow\fX\Ctrh_\cfq\rarrow\fX\Ctrh\rarrow
\fX\Lcth_\bW\rarrow\fX\Lcth$ induce triangulated equivalences\/
$\Hot^\st(\fX\Ctrh_\prj)\simeq\sD^\st(\fX\Ctrh^\fl)\simeq
\sD^\st(\fX\Ctrh_\cfq)\simeq\sD^\st(\fX\Ctrh)\simeq
\sD^\st(\fX\Lcth_\bW)$, as well as a triangulated equivalence\/
$\sD^-(\fX\Lcth_\bW)\simeq\sD^-(\fX\Lcth)$. \par
\textup{(d)} Let\/ $\fX$ be a Noetherian formal scheme of finite
Krull dimension with an open covering\/ $\bW$, and\/ $\st=-$, $\ctr$,
or\/~$\bctr$ be a derived category symbol.
 Then the inclusions of exact categories\/ $\fX\Lcth^\lct_\bW\rarrow
\fX\Lcth_\bW$ and\/ $\fX\Lcth^\lct\rarrow\fX\Lcth$ induce triangulated
equivalences\/ $\sD^\st(\fX\Lcth^\lct_\bW)\simeq\sD^\st(\fX\Lcth_\bW)$
and\/ $\sD^-(\fX\Lcth^\lct)\simeq\sD^-(\fX\Lcth)$.
\end{cor}

\begin{proof}
 This is our version of~\cite[Corollary~6.2.8]{Pcosh}.
 Part~(a) holds for any exact category of finite homological dimension
with enough projective and injective objects in the role of
$\fX\Ctrh^\fl$.
 See Corollary~\ref{flat-proj-lct-preenvelope-cor}(b)
and Theorem~\ref{fin-homol-dim-co-contra-derived-theorem}.

 In the contexts of parts~(b\+-c), the inclusions
$\fX\Ctrh^\lct_\prj\subset\fX\Ctrh^\lct_\cfq$
and $\fX\Ctrh^\fl\subset\fX\Ctrh_\cfq$ hold by
Corollaries~\ref{projective-lct-of-contramods-local-coflasque}(b)
and~\ref{W-flat-implies-flat-and-coflasque-cor}(c).
 All the respective full subcategories are closed under extensions
and kernels of admissible epimorphisms in each other (see
Sections~\ref{locally-contraherent-cosheaves-of-contramods-subsecn}\+-%
\ref{locally-contraherent-lct-cosheaves-subsecn} and
Lemma~\ref{coflasque-adjusted-to-cosections}(a\+-b)).
 All of them are also generating subcategories in each other
by Theorem~\ref{loc-Noetherian-projective-lct-cosh-contramods}(a)
and Corollary~\ref{projective-contrah-on-finite-Krull-dim}(a);
hence all of these subcategories are resolving in each other.
 All of these categories and full subcategories, with the exception of
$\fX\Ctrh_\prj$, $\fX\Lcth$, and $\fX\Lcth^\lct$, also have infinite
products and are closed under infinite products in each other
(see, in particular, Corollary~\ref{proj-lct-closed-under-products}).
 Hence the desired triangulated equivalences hold by
Propositions~\ref{infinite-co-resolutions-minus-plus-derived}(a),
\ref{infinite-co-resolutions-contra-co-derived}(a),
and~\ref{infinite-co-resolutions-becker-derived}(a).

 Finally, part~(d) is deduced by comparing parts~(a\+-c).
 We refer to Section~\ref{background-equivalences-subsecn} for
further details.
\end{proof}

 Let\/ $\fX$ be a Noetherian formal scheme of finite Krull dimension.
 It follows from Corollary~\ref{minus-ctr-bctr-derived-categories-agree}
that the $\Ext$ groups computed in the exact categories $\fX\Lcth_\bW$
and $\fX\Lcth_\bW^\lct$ agree with each other and with the $\Ext$
groups computed in the exact categories $\fX\Lcth$ and $\fX\Lcth^\lct$.
 Furthermore, they aslo agree with the $\Ext$ gropus computed in
the exact categories $\fX\Ctrh^\fl$, \,$\fX\Ctrh_\cfq$,
\,$\fX\Ctrh_\cfq^\lct$, etc.
 As in Section~\ref{antilocally-flat-subsecn}, we denote these $\Ext$
groups by $\Ext^{\fX,*}({-},{-})$.

\begin{cor} \label{flat-contraherent-cosheaves-cor}
 Let\/ $\fX$ be a Noetherian formal scheme of finite Krull dimension
with an open covering\/ $\bW$ and a finite affine open covering\/
$\fX=\bigcup_\alpha\fU_\alpha$ subordinate to\/~$\bW$.
 In this setting: \par
\textup{(a)} One has\/ $\Ext^{\fX,n}(\fF,\fQ)=0$ for any flat
contraherent cosheaf of contramodules\/ $\fF$ on\/ $\fX$, any locally
cotorsion\/ $\bW$\+locally contraherent cosheaf of contramodules\/
$\fQ$ on\/ $\fX$, and any integer $n\ge1$.
 Consequently, all flat contraherent cosheaves of contramodules on\/
$\fX$ are antilocally flat (with respect to any open
covering\/~$\bW$). \par
\textup{(b)} A contraherent cosheaf of contramodules\/ $\fG$ on\/ $\fX$
is flat if and only if\/ $\fG$ is a direct summand of a finitely
iterated extension of the direct images of flat contraherent cosheaves
of contramodules from\/~$\fU_\alpha$.
\end{cor}

\begin{proof}
 This is the formal scheme generalization
of~\cite[Corollary~6.2.9]{Pcosh}.
 Part~(b): the ``if'' implication holds by
Lemma~\ref{from-semi-separated-flat-direct-image}.
 The ``only if'' implication follows immediately from
Lemma~\ref{coflasque-flat-precover}(b) and
Corollary~\ref{flat-proj-lct-preenvelope-cor}(b).
 Since $\fP$ is an injective object in $\fX\Ctrh^\fl$, the short
exact sequence from Lemma~\ref{coflasque-flat-precover}(b) splits
and $\fG$ is a direct summand of~$\fF$.

 Part~(a): according to
Corollary~\ref{flat-proj-lct-preenvelope-cor}(b), any flat contraherent
cosheaf of contramodules on $\fX$ has a finite coresolution by
projective locally cotorsion contraherent cosheaves of contramodules.
 As the $\Ext$ groups in the exact categories $\fX\Lcth_\bW$ and
$\fX\Lcth^\lct_\bW$ agree by
Corollary~\ref{minus-ctr-bctr-derived-categories-agree}(d),
the assertion follows.

 Alternatively, one can avoid using
Corollary~\ref{minus-ctr-bctr-derived-categories-agree} in the proof
of part~(a) for $n=1$ by arguing as follows.
 The groups $\Ext^1$ computed in various exact categories in question
agree simply because these full subcategories are closed under
extensions in one another; so let us denote these groups by
$\Ext^{\fX,1}({-},{-})$.
 In view of part~(b), it suffices to show that
$\Ext^{\fX,1}(\fj_!\fH,\fQ)=0$ for any affine open formal subscheme
$\fU\subset\fX$ subordinate to $\bW$ with the open immersion morphism
$\fj\:\fU\rarrow\fX$ and any flat contraherent cosheaf of
contramodules $\fH$ on~$\fU$.
 The argument is based on $\Ext^1$\+adjunction considerations in
the spirit of~\cite[Lemma~1.7(d) or the opposite version of
Lemma~1.7(b)]{Pal}.

 Let $0\rarrow\fQ\rarrow\fM\rarrow\fj_!\fH\rarrow0$ be an (admissible)
short exact sequence in $\fX\Lcth_\bW$.
 Applying the exact functor~$\fj^!$, we obtain a short exact sequence
$0\rarrow\fj^!\fQ\rarrow\fj^!\fM\rarrow\fH\rarrow0$ in $\fU\Ctrh$,
which splits by Lemma~\ref{cotorsion-contramodules-lemma}(2)
(as $\fj^!\fQ\in\fU\Ctrh^\lct$).
 Let $\fH\rarrow\fj^!\fM$ be a section for the split epimorphism
$\fj^!\fM\rarrow\fH$ in $\fU\Ctrh$.
 By adjunction~\eqref{cosheaves-direct-image-restriction-adjunction}
or~\eqref{cosheaves-direct-inverse-adjunction}, we have
the corresponding morphism $\fj_!\fH\rarrow\fM$ in $\fX\Lcth_\bW$,
which splits the epimorphism $\fM\rarrow\fj_!\fH$.

 Part~(a) for $n\ge2$ can be subsequently deduced from the $n=1$ case
using Corollary~\ref{projective-contrah-on-finite-Krull-dim}(a,c) and
the fact that the full subcategory $\fX\Ctrh^\fl$ is closed under
kernels of admissible epimorphisms in $\fX\Lcth_\bW$.
\end{proof}

 The argument in the following lemma, proposition, and corollary is
our version of~\cite[Theorem~1.1]{AB}.
 See~\cite[Remark~10.10]{Pns} for a discussion.

\begin{lem} \label{cokernel-lct-preenvelope-lemma}
 Let\/ $\fX$ be a Noetherian formal scheme of finite Krull dimension
with an open covering\/~$\bW$.
 Consider the class of all\/ $\bW$\+locally contraherent cosheaves of
contramodules\/ $\fM$ on\/ $\fX$ for which there exists an admissible
short exact sequence\/ $0\rarrow\fM\rarrow\fP\rarrow\fF\rarrow0$ in\/
$\fX\Lcth_\bW$ with\/ $\fP\in\fX\Lcth^\lct_\bW$ and\/
$\fF\in\fX\Ctrh^\fl$.
 Then the class of all such objects\/ $\fM$ is closed under cokernels
of admissible monomorphisms in\/~$\fX\Lcth_\bW$.
\end{lem}

\begin{proof}
 Let $0\rarrow\fM\rarrow\fP\rarrow\fF\rarrow0$ and $0\rarrow\fN\rarrow
\fQ\rarrow\fG\rarrow0$ be two (admissible) short exact sequences in
$\fX\Lcth_\bW$ with $\fP$, $\fQ\in\fX\Lcth_\bW^\lct$ and $\fF$,
$\fG\in\fX\Ctrh^\fl$.
 Let $f\:\fM\rarrow\fN$ be a morphism in $\fX\Lcth_\bW^\lct$.
 By Corollary~\ref{flat-contraherent-cosheaves-cor}(a), we have
$\Ext^{\fX,1}(\fF,\fQ)=0$.
 Therefore, the composition $\fM\rarrow\fN\rarrow\fQ$ factorizes
through the admissible monomorphism $\fM\rarrow\fP$, and we have
the commutative diagram of a morphism of short exact sequences
\begin{equation} \label{morphism-of-short-exact-sequences-diag-I}
\begin{gathered}
 \xymatrix{
  0 \ar[r] & \fM \ar[r] \ar[d]_f & \fP \ar[r] \ar[d] & \fF
  \ar[r] \ar[d]^g & 0 \\
  0 \ar[r] & \fN \ar[r] & \fQ \ar[r] & \fG \ar[r] & 0
 }
\end{gathered}
\end{equation}

 By Corollary~\ref{flat-proj-lct-preenvelope-cor}(a), there exists
a short exact sequence $0\rarrow\fF\overset t\rarrow\fT\rarrow\fH
\rarrow0$ in $\fX\Ctrh^\fl$ with $\fT\in\fX\Ctrh^\lct_\prj$.
 Adding the trivial short exact sequence $0\rarrow0\rarrow\fT\rarrow\fT
\rarrow0$ as an additional direct summand to the lower line
of~\eqref{morphism-of-short-exact-sequences-diag-I}, we obtain another
commutative diagram of a morphism of short exact sequences
\begin{equation} \label{morphism-of-short-exact-sequences-diag-II}
\begin{gathered}
 \xymatrix{
  0 \ar[r] & \fM \ar[r] \ar[d]_f & \fP \ar[r] \ar[d]^r & \fF
  \ar[r] \ar[d]^{(g,t)} & 0 \\
  0 \ar[r] & \fN \ar[r] & \fQ\oplus\fT \ar[r] & \fG\oplus\fT \ar[r] & 0
 }
\end{gathered}
\end{equation}
 Now $t$~is an admissible monomorphism in $\fX\Ctrh^\fl$, hence
$(g,t)$ is an admissible monomorphism in $\fX\Ctrh^\fl$ as well.
 Assuming that $f$~is an admissible monomorphism in $\fX\Lcth_\bW$,
it follows that $r$~is also an admissible monomorphism in
$\fX\Lcth_\bW$.
 Actually, we have $\fP\in\fX\Lcth^\lct_\bW$ and
$\fQ\oplus\fT\in\fX\Lcth^\lct_\bW$.
 Any admissible monomorphism in $\fX\Lcth$ between objects from
$\fX\Lcth^\lct_\bW$ is an admissible monomorphism in
$\fX\Lcth^\lct_\bW$, so $r$~is an admissible monomorphism
in $\fX\Lcth^\lct_\bW$.

 Passing to the cokernels in the termwise admissible monomorphism
of short exact
sequences~\eqref{morphism-of-short-exact-sequences-diag-II},
we obtain the desired representation of $\fN/\fM=\coker(f)$ as
the kernel of an admissible epimorphism in $\fX\Lcth_\bW$ acting
from an object of $\fX\Lcth^\lct_\bW$ to an object of $\fX\Ctrh^\fl$.
\end{proof}

 The result of the following proposition was mentioned as an open
question in the more restrictive context of (nonformal) schemes
in~\cite[Section~6.2]{Pcosh}.
 Now we are able to prove it in the more general context of formal
schemes.

\begin{prop} \label{finite-krulldim-locally-cotorsion-preenvelope}
 Let\/ $\fX$ be a Noetherian formal scheme of finite Krull dimension
$D$ with an open covering\/~$\bW$. \par
\textup{(a)} Any\/ $\bW$\+locally contraherent cosheaf of
contramodules\/ $\fM$ on\/ $\fX$ can be included in an (admissible)
short exact sequence\/ $0\rarrow\fM\rarrow\fP\rarrow\fF\rarrow0$
in the exact category\/ $\fX\Lcth_\bW$, where\/ $\fP$ is a locally
cotorsion\/ $\bW$\+locally contraherent cosheaf of contramodules and\/
$\fF$ is a flat contraherent cosheaf of contramodules on\/~$\fX$. \par
\textup{(b)} The full subcategory of locally cotorsion\/ $\bW$\+locally
contraherent cosheaves of contramodules\/ $\fX\Lcth^\lct_\bW$ is
coresolving in the exact category of $\bW$\+locally contraherent
cosheaves of contramodules\/ $\fX\Lcth_\bW$.
 The coresolution dimension of any object of\/ $\fX\Lcth_\bW$ with
respect to the coresolving subcategory\/ $\fX\Lcth^\lct_\bW$ does not
exceed~$D$.
\end{prop}

\begin{proof}
 Part~(a): by Lemma~\ref{coflasque-locally-cotorsion-preenvelope}(a),
the assertion of the proposition holds for any coflasque contraherent
cosheaf of contramodules $\fM=\fE$.
 Furthermore, the full subcategory $\fX\Ctrh_\cfq$ is resolving in
$\fX\Lcth_\bW$, as it was pointed out in the proof of
Corollary~\ref{minus-ctr-bctr-derived-categories-agree}(b\+-c).
 Moreover, the resolution dimension of any object of $\fX\Lcth_\bW$
with respect to the resolving subcategory $\fX\Ctrh_\cfq\subset
\fX\Lcth_\bW$ does not exceed the Krull dimension $D$ of
the formal scheme $\fX$, by
Lemma~\ref{finite-Krull-dim-co-flasque-co-resol-dim}(b).
 Thus the cosheaf $\fM$ has a finite resolution by cosheaves from
$\fX\Ctrh_\cfq$ in the exact category $\fX\Lcth_\bW$.
 Now the existence of the desired short exact sequence
$0\rarrow\fM\rarrow\fP\rarrow\fF\rarrow0$ follows by iterated
application of Lemma~\ref{cokernel-lct-preenvelope-lemma}.

 The proof of part~(b) is similar to that of
Lemma~\ref{coflasque-locally-cotorsion-preenvelope}(b); the point
is that the full subcategory $\fX\Lcth^\lct_\bW$ is cogenerating
in $\fX\Lcth_\bW$ by part~(a).
\end{proof}

\begin{cor} \label{finite-krulldim-flat-precover}
 Let\/ $\fX$ be a Noetherian formal scheme of finite Krull dimension
with an open covering\/~$\bW$.
 Then any\/ $\bW$\+locally contraherent cosheaf of contramodules\/ $\fM$
on\/ $\fX$ can be included in an (admissible) short exact sequence\/
$0\rarrow\fP\rarrow\fF\rarrow\fM\rarrow0$ in the exact category\/
$\fX\Lcth_\bW$, where\/ $\fP$ is a locally cotorsion\/ $\bW$\+locally
contraherent cosheaf of contramodules and\/ $\fF$ is a flat contraherent
cosheaf of contramodules on\/~$\fX$.
\end{cor}

\begin{proof}
 This is yet another application of a general version of
the Salce lemma, as in~\cite[Lemma~B.1.1(a)]{Pcosh}.
 The argument is similar to the proof of
Lemma~\ref{coflasque-flat-precover}(a) and based on
Proposition~\ref{finite-krulldim-locally-cotorsion-preenvelope}(a)
together with
Corollary~\ref{projective-contrah-on-finite-Krull-dim}(a,c).
\end{proof}

\begin{cor} \label{finite-krulldim-flat-cosheaves-characterized}
 Let\/ $\fX$ be a Noetherian formal scheme of finite Krull dimension
with an open covering\/~$\bW$.
 In this setting: \par
\textup{(a)} A\/ $\bW$\+locally contraherent cosheaf of contramodules\/
$\fF$ on\/ $\fX$ is flat if and only if\/ $\Ext^{\fX,1}(\fF,\fP)=0$ for
all locally cotorsion\/ $\bW$\+locally contraherent cosheaves of
contramodules\/ $\fP$ on\/ $\fX$, and if and only if\/
$\Ext^{\fX,n}(\fF,\fP)=0$ for all locally cotorsion\/ $\bW$\+locally
contraherent cosheaves of contramodules\/ $\fP$ on\/ $\fX$ and all
$n\ge1$. \par
\textup{(b)} A\/ $\bW$\+locally contraherent cosheaf of contramodules
on\/ $\fX$ is flat if and only if it is antilocally flat.
 Consequently, the full subcategory of antilocally flat\/ $\bW$\+locally
contraherent cosheaves of contramodules in the exact category\/
$\fX\Lcth$ does not depend on the open covering\/~$\bW$.
\end{cor}

\begin{proof}
 Part~(a): the ``only if'' implications hold by
Corollary~\ref{flat-contraherent-cosheaves-cor}(a).
 The ``if'' implications follow from
Corollary~\ref{finite-krulldim-flat-precover} by virtue of
the argument proving Lemma~\ref{direct-summand-lemma} below.

 Part~(b): notice first that the full subcategory of flat
$\bW$\+locally contraherent cosheaves of contramodules in $\fX\Lcth$
does not depend on $\bW$ by
Corollary~\ref{W-flat-implies-flat-and-coflasque-cor}(b), and all
flat contraherent cosheaves on contramodules on $\fX$ are antilocally
flat by Corollary~\ref{flat-contraherent-cosheaves-cor}(a).
 The assertion that all antilocally flat $\bW$\+locally contraherent
cosheaves of contramodules are flat follows from part~(a) by virtue of
Lemma~\ref{Hom-from-to-short-exact-sequences-lemma}(b) below,
which is applicable since the full subcategory $\fX\Lcth_\bW^\lct$
is cogenerating in $\fX\Lcth_\bW$ by
Proposition~\ref{finite-krulldim-locally-cotorsion-preenvelope}.
\end{proof}

\begin{cor} \label{finite-krulldim-lct-cosheaves-characterized}
 Let\/ $\fX$ be a Noetherian formal scheme of finite Krull dimension
with an open covering\/~$\bW$.
 Then, for any\/ $\bW$\+locally contraherent cosheaf of contramodules\/
$\fP$ on\/ $\fX$, the following conditions are equivalent:
\begin{enumerate}
\item the functor\/ $\Hom^\fX({-},\fP)$ takes short exact sequences of
flat contraherent cosheaves of contramodules on\/ $\fX$ to short exact
sequences of abelian groups;
\item one has\/ $\Ext^{\fX,1}(\fF,\fP)=0$ for all flat contraherent
cosheaves of contramodules\/ $\fF$ on\/~$\fX$;
\item one has\/ $\Ext^{\fX,n}(\fF,\fP)=0$ for all flat contraherent
cosheaves of contramodules\/ $\fF$ on\/ $\fX$ and all integers $n\ge1$;
\item $\fP$~is a locally cotorsion\/ $\bW$\+locally contraherent
cosheaf of contramodules on\/~$\fX$.
\end{enumerate}
\end{cor}

\begin{proof}
 The proof is similar to that of
Corollary~\ref{locally-cotorsion-characterized-cor}, and based on
Corollary~\ref{flat-contraherent-cosheaves-cor}(a) and
Proposition~\ref{finite-krulldim-locally-cotorsion-preenvelope}(a)
together with a weak version of
Corollary~\ref{projective-contrah-on-finite-Krull-dim}(a,c)
or~\ref{finite-krulldim-flat-precover}.
\end{proof}

 In the terminology of Section~\ref{appx-cotorsion-pairs-subsecn}
below, the results of
Proposition~\ref{finite-krulldim-locally-cotorsion-preenvelope}(a)
and Corollaries~\ref{finite-krulldim-flat-precover},
\ref{finite-krulldim-flat-cosheaves-characterized}(a),
and~\ref{finite-krulldim-lct-cosheaves-characterized}
can be summarized by the following formulation.
 For any Noetherian formal scheme $\fX$ of finite Krull dimension with
an open covering $\bW$, the pair of classes (flat contraherent
cosheaves of contramodules on~$\fX$, locally cotorsion $\bW$\+locally
contraherent cosheaves of contramodules on~$\fX$) is a \emph{hereditary
complete cotorsion pair} in the exact category of $\bW$\+locally
contraherent cosheaves of contramodules $\fX\Lcth_\bW$.
 Even for (nonformal) schemes, this formulation is an improvement upon
the result of~\cite[Corollary~6.2.11]{Pcosh}, which only provided
such a cotorsion pair in the exact category of \emph{coflasque}
contraherent cosheaves.

\begin{rem} \label{flat-and-antilocally-flat-remark}
 For the reader's benefit and as an improvement
upon~\cite[Remark~6.2.10]{Pcosh}, let us summarize our results
concerning flat and antilocally flat contraherent cosheaves
of contramodules.
 The definitions of both the concepts were given in
Section~\ref{antilocally-flat-subsecn}.

 On a semi-separated Noetherian formal scheme, the class of antilocally
flat $\bW$\+locally contraherent cosheaves of contramodules
does not depend on an open covering~$\bW$
(Corollary~\ref{antilocally-flat-independent-on-covering}),
and all antilocally flat contraherent cosheaves of contramodules
are flat (Corollary~\ref{antilocally-flat-independent-on-covering}).

 On a locally Noetherian formal scheme $\fX$ of finite Krull dimension,
any $\bW$\+flat $\bW$\+locally contraherent cosheaf of contramodules
is flat and contraherent
(Corollary~\ref{W-flat-implies-flat-and-coflasque-cor}(b)).
 So the class of $\bW$\+flat $\bW$\+locally contraherent cosheaves
of contramodules on $\fX$ does not depend on an open covering~$\bW$.

 On a Noetherian formal scheme of finite Krull dimension, the class
of antilocally flat $\bW$\+locally contraherent cosheaves of
contramodules does not depend on an open covering $\bW$, and
coincides with the class of flat contraherent cosheaves of contramodules
(Corollary~\ref{finite-krulldim-flat-cosheaves-characterized}(b)).

 On a semi-separated Noetherian formal scheme, a description of
antilocally flat contraherent cosheaves of contramodules is provided
by Corollary~\ref{antilocally-flat-corollary}(c), and several
equivalent characterizations can be found in
Corollary~\ref{antilocally-flat-characterized-cor}(a).

 On a Noetherian formal scheme of finite Krull dimension, a very
similar description of (antilocally) flat contraherent cosheaves of
contramodules is provided by
Corollary~\ref{flat-contraherent-cosheaves-cor}(b), and several
equivalent characterizations can be found in
Corollary~\ref{finite-krulldim-flat-cosheaves-characterized}(a).

 Similarly to the discussion in Section~\ref{antilocally-flat-subsecn}, 
the result of Corollary~\ref{flat-contraherent-cosheaves-cor}(b) can be
rephrased by saying that the class of (antilocally) flat contraherent
cosheaves of contramodules on Noetherian formal schemes of finite
Krull dimension is \emph{antilocal} in a sense
resembling~\cite[Section~4]{Pal}.

 The assertion that on a Noetherian formal scheme of finite Krull
dimension the classes of flat and antilocally flat contraherent
cosheaves of contramodules coincide is our main improvement
upon~\cite[Remark~6.2.10]{Pcosh} in the context of the present remark
(on top of the generalization from schemes to formal schemes).
\end{rem}

\begin{cor} \label{flat-contrah-direct-images}
 Let\/ $\fX$ be a locally Noetherian formal scheme, $\fY$ be
a locally Noetherian formal scheme of finite Krull dimension, and\/
$\ff\:\fY\rarrow\fX$ be a quasi-compact morphism of formal schemes.
 Then \par
\textup{(a)} the functor of direct image of cosheaves of
contramodule\/ $\fO$\+modules\/~$\ff_!$ takes flat contraherent
cosheaves of contramodules on\/ $\fY$ to contraherent cosheaves
of contramodules on\/ $\fX$, and induces an exact functor\/
$\ff_!\:\fY\Ctrh^\fl\rarrow\fX\Ctrh$ between these exact categories;
\par
\textup{(b)} if the morphism of formal schemes\/~$\ff$ is tight and
flat, then the functor of direct image of cosheaves of contamodule\/
$\fO$\+modules\/~$\ff_!$ takes flat contraherent cosheaves of
contramodules on\/ $\fY$ to flat contraherent cosheaves of
contramodules on\/ $\fX$, and induces an exact functor\/
$\ff_!\:\fY\Ctrh^\fl\rarrow\fX\Ctrh^\fl$ between these exact
categories; \par
\textup{(c)} if the formal scheme\/ $\fY$ is Noetherian (i.~e.,
quasi-compact) and the morphism of formal schemes\/~$\ff$ is tight
and very flat, then the functor of direct image of cosheaves of
contramodule\/ $\fO$\+modules\/~$\ff_!$ takes projective contraherent
cosheaves of contramodles on\/ $\fY$ to projective contraherent
cosheaves of contramodules on\/~$\fX$.
\end{cor}

\begin{proof}
 This is the formal scheme generalization
of~\cite[Corollary~6.2.12]{Pcosh}, and the arguments are similar.

 Part~(a) is a particular case of
Corollary~\ref{co-flasque-direct-image-cor}(b) in view of
Corollary~\ref{W-flat-implies-flat-and-coflasque-cor}(c).

 In part~(b), one can assume the formal scheme $\fX$ is affine;
then the formal scheme $\fY$ is Noetherian of finite Krull dimension.
 Using part~(a), it suffices to show that the contramodule
$\fO_\fX(\fX)$\+module $(\ff_!\fF)[\fX]=\ff[\fY]$ is flat for any
flat contraherent cosheaf $\fF$ on~$\fY$.
 For this purpose, consider a finite coresolution of the cosheaf $\fF$
by objects from $\fY\Ctrh^\lct_\prj$ in $\fY\Ctrh^\lct$ (which exists
by Corollary~\ref{flat-proj-lct-preenvelope-cor}(b)), and apply
the functor of cosections $\Delta(\fY,{-})$ to such a coresolution.
 By Corollary~\ref{W-flat-implies-flat-and-coflasque-cor}(c)
and Lemma~\ref{coflasque-adjusted-to-cosections}(c), the resulting
sequence remains exact.
 In view of Corollary~\ref{proj-lct-direct-images}(b), we obtain
a finite coresolution of the contramodule $\fO_\fX(\fX)$\+module
$\fF[\fY]$ by flat (and cotorsion) contramodule $\fO_\fX(\fX)$\+modules.
 Thus the contramodule $\fO_\fX(\fX)$\+module $\fF[\fY]$ is also flat.

 Part~(c) follows from part~(a) or~(b) together with
Corollary~\ref{projective-contrah-on-finite-Krull-dim}(a,c),
exactness of the inverse image
functor~\eqref{inverse-image-lcth-contramod-nocovering},
and the adjunction~\eqref{cosheaves-direct-inverse-adjunction}.
\end{proof}

\subsection{Homology of contraherent cosheaves of contramodules~II}
\label{homology-of-cosheaves-II-subsecn}
 This section follows up on the discussion in
Section~\ref{homology-of-cosheaves-I-subsecn}.

 Let $\fX$ be a locally Noetherian formal scheme.
 Then the left derived functor $\boL_*\Delta(\fX,{-})\:\fX\Lcth^\lct
\rarrow\fO_\fX(\fO)\Modl$ of the functor of global cosections
$\Delta(\fX,{-})$ is defined as usual using projective resolutions
in the exact category $\fX\Lcth^\lct$ (which exist by
Theorem~\ref{loc-Noetherian-projective-lct-cosh-contramods}(a)).
 It follows from the same theorem that the derived functors of
$\Delta(\fX,{-})$ computed in the exact category $\fX\Lcth^\lct_\bW$
for a fixed open covering $\bW$ and in the whole exact category
$\fX\Lcth^\lct$ agree.

 Similarly to~\cite[Section~6.3]{Pcosh}, the groups/modules
$\boL_*\Delta(\fX,\fM)$ are called the \emph{homology groups/modules}
of a locally cotorsion locally contraherent cosheaf of contramodules
$\fM$ on a locally Noetherian formal scheme~$\fX$.

 Let us show that $\boL_n\Delta(\fX,\fE)=0$ for any coflasque
locally cotorsion contraherent cosheaf of contramodules $\fE$ on $\fX$
and all integers $n\ge1$.
 By Corollary~\ref{projective-lct-of-contramods-local-coflasque}(b),
any projective locally cotorsion contraherent cosheaf of contramodules
on $\fX$ is coflasque.
 In view of Lemma~\ref{coflasque-adjusted-to-cosections}(b), any
resolution of an object of $\fX\Ctrh^\lct_\cfq$ by objects from
$\fX\Ctrh^\lct_\cfq$ in the exact category $\fX\Lcth^\lct$ is
actually exact in the exact category $\fX\Lcth^\lct_\cfq$.
 By part~(c) of the same lemma, the functor $\Delta(\fX,{-})$ preserves
exactness of such exact sequences.
 Hence one can compute the derived functor $\boL_*\Delta(\fX,{-})$
using resolutions by coflasque locally cotorsion contraherent
cosheaves of contramodules.

 Now let $\fX$ be a Noetherian formal scheme of finite Krull dimension.
 Then the left derived functor $\boL_*\Delta(\fX,{-})\:\fX\Lcth
\rarrow\fO_\fX(\fO)\Modl$ of the functor of global sections
$\Delta(\fX,{-})$ is defined as usual using projective resolutions
in the exact category $\fX\Lcth$ (which exist by
Corollary~\ref{projective-contrah-on-finite-Krull-dim}(a)).
 It follows from the same corollary that the derived functors of
$\Delta(\fX,{-})$ computed in the exact category $\fX\Lcth_\bW$ for
a fixed open covering $\bW$ and in the whole exact category $\fX\Lcth$
agree.
 Furthermore, in view of
Corollary~\ref{projective-contrah-on-finite-Krull-dim}(c)
and Lemma~\ref{coflasque-adjusted-to-cosections}(b\+-c) one can
compute the derived functor $\boL_*\Delta(\fX,{-})$ using coflasque
resolutions.
 The groups/modules $\boL_*\Delta(\fX,\fM)$ are called
the \emph{homology groups/modules} of a locally contraherent cosheaf
of contramodules $\fM$ on a Noetherian formal scheme $\fX$ of finite
Krull dimension.

 It is clear from the discussion above (in particular, concerning
the coflasque resolutions) that the constructions of the derived
functors $\boL_*\Delta(\fX,{-})$ on the exact categories
$\fX\Lcth^\lct$ and $\fX\Lcth$ above agree when they are both
applicable.
 They also obviously agree with the constructions given
in Section~\ref{homology-of-cosheaves-I-subsecn}.

 Using injective coresolutions, one can similarly define
the right derived functor of global sections of quasi-coherent torsion
sheaves $\boR^*\Gamma(\fX,{-})\:\fX\Tors\rarrow\fO_\fX(\fX)\Modl$ for
a locally Noetherian formal scheme~$\fX$.
 The abelian category of quasi-coherent torsion sheaves $\fX\Tors$ has
enough injective objects by
Theorem~\ref{loc-Noetherian-injective-qcoh-torsion}(a),
and all of them are flasque sheaves by
Corollary~\ref{injective-torsion-local-flasque}(b); hence the derived
functor $\boR^*\Gamma(\fX,{-})$ agrees with the conventional cohomology
of sheaves of abelian groups or $\fO_\fX$\+modules on~$\fX$.
 So we have $\boR^n\Gamma(\fX,\cM)\simeq H^n(\fX,\cM)$ for all
$\cM\in\fX\Tors$ and $n\ge0$.

\subsection{Special inverse images}
\label{special-inverse-images-subsecn}
 This section is a formal scheme generalization
of~\cite[Section~7.3]{Pcosh}.
 We denote by $\cO_X$ the structure sheaf of a (nonformal) scheme~$X$.
 So $\fO_\fX=\cO_X$ if $\fX=X$ is a scheme.
 
 Let $\ff\:\fY\rarrow\fX$ be a tight morphism of locally Noetherian
formal schemes.
 Then, given a closed subscheme of definition $\overline X\subset\fX$,
the related closed subscheme of definition $\overline Y=
\overline X\times_\fX\fY\subset\fY$ is constructed as follows.

 Let $\fV\subset\fY$ and $\fU\subset\fX$ be a pair of affine open
formal subschemes such that $\ff(\fV)\subset\fU$.
 So $\fU=\Spf\fR$ and $\fV=\Spf\fS$ for some complete, separated
adic Noetherian rings $\fR$ and~$\fS$.
 Let $\ff|^\fU_\fV\:\fV\rarrow\fU$ be the induced morphism of affine
Noetherian formal schemes and $f_{\fU,\fV}\:\fR\rarrow\fS$ be
the corresponding continuous morphism of adic Noetherian rings,
as in Section~\ref{morphisms-of-formal-schemes-subsecn}.
 By assumption, the ring homomorphism $f_{\fU,\fV}$ is tight.

 Put $\overline U=\overline X\cap\fU$; so $\overline U$ is a closed
subscheme of definition in~$\fU$.
 Let $\fI\subset\fR$ be the ideal of definition corresponding to
the closed subscheme of definition $\overline U\subset\fU$.
 Since the ring homomorphism $f_{\fU,\fV}\:\fR\rarrow\fS$ is
continuous and tight, the ideal $\fJ=\fS f_{\fU,\fV}(\fI)$ is
an ideal of definition in~$\fS$.
 Let $\overline V\subset\fV$ be the closed subscheme of definition
corresponding to the ideal of definition $\fJ\subset\fS$.
 Then the closed subscheme of definition $\overline Y\subset\fY$
is defined by the rule that $\overline Y\cap\fV=\overline V$ for
all affine open formal subschemes $\fV\subset\fY$ and $\fU\subset\fX$
as above.

\medskip

 Let $\fii\:\fZ\rarrow\fX$ be a tight affine morphism of locally
Noetherian formal schemes.
 Let $\overline X\subset\fX$ be a closed subscheme of definition;
put $\overline Z=\overline X\times_\fX\fZ\subset\fZ$ as above.
 Let $\fU\subset\fX$ be an affine open formal subscheme; then
$\fV=\fii^{-1}(\fU)$ is an affine open formal subscheme in~$\fZ$.
 Put $\overline U=\overline X\cap\fU$ and $\overline V=
\overline Z\cap\fV$.

 By Corollary~\ref{finite-maps-of-adic-topological-rings-cor},
the ring $\fO_\fZ(\fV)$ is a finitely generated module over
the ring $\fO_\fX(\fU)$ if and only if the ring
$\cO_{\overline Z}(\overline V)$ is a finitely generated module over
the ring $\cO_{\overline X}(\overline U)$.
 If this is the case for all affine open formal subschemes
$\fU\subset\fX$, the morphism of formal schemes $\fii\:\fZ\rarrow\fX$
is said to be \emph{finite}.

 So $\fii\:\fZ\rarrow\fX$ is a finite morphism of formal schemes
if and only if $\bar i\:\overline Z=\overline X\times_\fX\fZ\rarrow
\overline X$ is a finite morphism of schemes (for any chosen closed
subscheme of definition $\overline X\subset\fX$, or equivalently,
for all such closed subschemes of definition).
 Being a finite morphism of schemes is a local
property~\cite[Lemma Tag~01WI]{SP}; so it suffices to check
the finiteness condition above for affine open subschemes $\overline U$
belonging to any chosen affine open covering of the scheme
$\overline X$, or equivalently, for affine open formal subschemes $\fU$
belonging to any chosen affine open covering of the formal scheme~$\fX$.

 In particular, by
Corollary~\ref{surjective-maps-of-adic-topological-rings-cor},
the ring homomorphism $\fO_\fX(\fU)\rarrow\fO_\fZ(\fV)$ is surjective
if and only if the ring homomorphism $\cO_{\overline X}(\overline U)
\rarrow\cO_{\overline Z}(\overline V)$ is surjective.
 If this is the case for all affine open formal subschemes
$\fU\subset\fX$, the morphism of formal schemes $\fii\:\fZ\rarrow\fX$
is said to be a \emph{closed immersion}.

 So $\fii\:\fZ\rarrow\fX$ is a closed immersion of formal schemes
if and only if $\bar i\:\overline Z=\overline X\times_\fX\fZ\rarrow
\overline X$ is a closed immersion of schemes (for any chosen closed
subscheme of definition $\overline X\subset\fX$, or equivalently,
for all such closed subschemes of definition).
 Being a closed immersion  of schemes is a local
property~\cite[Lemma Tag~01QO]{SP}; so it suffices to check
the surjectivity condition above for affine open subschemes
$\overline U$ belonging to any chosen affine open covering of the scheme
$\overline X$, or equivalently, for affine open formal subschemes $\fU$
belonging to any chosen affine open covering of the formal scheme~$\fX$.

 Let us emphasize that any finite morphism (in particular, any closed
immersion) of locally Noetherian formal schemes is a tight affine
morphism by the definition.

\begin{lem} \label{finite-morphism-preimage-tensor-isomorphism-lemma}
 Let\/ $\fk\:\fZ\rarrow\fX$ be a finite morphism of locally Noetherian
formal schemes.
 Let\/ $\fU'\subset\fU$ be two affine open formal subschemes in\/ $\fX$
and\/ $\fk^{-1}(\fU')\subset\fk^{-1}(\fU)$ be their preimages
in\/~$\fZ$.
 Then the induced ring map
\begin{equation} \label{finite-morphism-preimage-tensor-isom-formula}
 \fO_\fX(\fU')\ot_{\fO_\fX(\fU)}\fO_\fZ(\fk^{-1}(\fU))
 \lrarrow\fO_\fZ(\fk^{-1}(\fU'))
\end{equation}
is an isomorphism of adic topological rings.
\end{lem}

\begin{proof}
 All the topological ring homomorphisms involved are continuous and
tight by assumption and by
Lemma~\ref{tight-flat-adic-ring-base-change-context}(a)
or~\ref{open-immersion-adic-ring-base-change-context}(a).
 So the notation for the adic completion functor $\Lambda$
is unambigous.
 According to Section~\ref{morphisms-of-formal-schemes-subsecn},
the topological ring
map~\eqref{finite-morphism-preimage-tensor-isom-formula} becomes
an isomorphism after the functor $\Lambda$ is applied.
 It remains to point out that the left-hand side
of~\eqref{finite-morphism-preimage-tensor-isom-formula} is a finitely
generated module over the complete, separated adic Noetherian ring
$\fO_\fX(\fU')$, and refer to
Lemma~\ref{finitely-generated-modules-are-complete-separated}(b).
\end{proof}

 Let $\fk\:\fZ\rarrow\fX$ be a finite morphism of locally Noetherian
formal schemes.
 Then the functor of special inverse image of quasi-coherent torsion
sheaves
\begin{equation} \label{qcoh-tors-special-inverse-image}
 \fk^!\:\fX\Tors\lrarrow\fZ\Tors
\end{equation}
is constructed as follows.

 Let $\cM$ be a quasi-coherent torsion sheaf on~$\fX$.
 Then the quasi-coherent torsion sheaf $\fk^!\cM$ on $\fZ$ is defined
by the rule
\begin{equation} \label{spec-inv-image-of-qcoh-tors-defining-rule}
 (\fk^!\cM)(\fV)=\fO_\fZ(\fV)\ot_{\fO_\fZ(\fk^{-1}(\fU))}
 \Hom_{\fO_\fX(\fU)}\bigl(\fO_\fZ(\fk^{-1}(\fU)),\cM(\fU)\bigr)
\end{equation}
for any pair of affine open formal subschemes $\fV\subset\fZ$ and
$\fU\subset\fX$ such that $\fk(\fV)\subset\fU$.

 So, in particular, denoting by $k_{\fU,\fk^{-1}(\fU)}\:
\fO_\fX(\fU)\rarrow\fO_\fZ(\fk^{-1}(\fU))$ the finite map of adic
Noetherian rings corresponding to the morphism of affine Noetherian
formal schemes $\fk|^\fU_{\fk^{-1}(\fU)}\:\fk^{-1}(\fU)\rarrow\fU$,
we have
$$
 \Hom_{\fO_\fX(\fU)}\bigl(\fO_\fZ(\fk^{-1}(\fU)),\cM(\fU)\bigr)=
 k_{\fU,\fk^{-1}(\fU)}^\diamond(\cM(\fU))
$$
by Lemma~\ref{finpres-topol-ring-map-co-contra-extensn-of-scalars}(a).
 This proves that $\Hom_{\fO_\fX(\fU)}(\fO_\fZ(\fk^{-1}(\fU)),\cM(\fU))$
is a torsion module over $\fO_\fZ(\fk^{-1}(\fU))$.
 Consequently, the right-hand side
of~\eqref{spec-inv-image-of-qcoh-tors-defining-rule}
is a torsion $\fO_\fZ(\fV)$\+module by
Lemma~\ref{tight-co-extension-of-scalars}(a).

 The following computation from~\cite[Section~7.3]{Pcosh} shows that
the torsion $\fO_\fZ(\fV)$\+mod\-ule $(\fk^!\cM)(\fV)$ is well-defined
by the formula~\eqref{spec-inv-image-of-qcoh-tors-defining-rule}.
 Let $\fU'\subset\fU$ be two affine open formal subschemes in~$\fX$.
 Then we have {\hbadness=1625
\begin{multline*}
 \Hom_{\fO_\fX(\fU')}\bigl(\fO_\fZ(\fk^{-1}(\fU')),\cM(\fU')\bigr) \\
 \simeq\Hom_{\fO_\fX(\fU')}
 \bigl(\fO_\fX(\fU')\ot_{\fO_\fX(\fU)}\fO_\fZ(\fk^{-1}(\fU)),\>
 \fO_\fX(\fU')\ot_{\fO_\fX(\fU)}\cM(\fU)\bigr) \\
 \simeq\Hom_{\fO_\fX(\fU)}\bigl(\fO_\fZ(\fk^{-1}(\fU)),\>
 \fO_\fX(\fU')\ot_{\fO_\fX(\fU)}\cM(\fU)\bigr) \\
 \simeq\fO_\fX(\fU')\ot_{\fO_\fX(\fU)}
 \Hom_{\fO_\fX(\fU)}\bigl(\fO_\fZ(\fk^{-1}(\fU)),\cM(\fU)\bigr) \\
 \simeq\bigl(\fO_\fX(\fU')\ot_{\fO_\fX(\fU)}
 \fO_\fZ(\fk^{-1}(\fU))\bigr)\ot_{\fO_\fZ(\fk^{-1}(\fU))}
 \Hom_{\fO_\fX(\fU)}\bigl(\fO_\fZ(\fk^{-1}(\fU)),\cM(\fU)\bigr) \\
 \simeq\fO_\fZ(\fk^{-1}(\fU'))\ot_{\fO_\fZ(\fk^{-1}(\fU))}
 \Hom_{\fO_\fX(\fU)}\bigl(\fO_\fZ(\fk^{-1}(\fU)),\cM(\fU)\bigr).
\end{multline*}
 Here} the first and the final (fifth) isomorphisms hold by
Lemma~\ref{finite-morphism-preimage-tensor-isomorphism-lemma},
while the third isomorphism (or the combination of the third and
the fourth isomorphisms) is provided by
Lemma~\ref{Hom-into-tensor-and-injective-tensored-with-flat}(a)
or Lemma~\ref{co-extension-of-scalars-tensor-product-lemma}(a) or~(b).

 The computation in the previous paragraph proves that
the formula~\eqref{spec-inv-image-of-qcoh-tors-defining-rule}
produces the same torsion $\fO_\fZ(\fV)$\+module $(\fk^!\cM)(\fV)$
when applied to two pairs of affine open formal subschemes
$(\fV,\fU)$ and $(\fV,\fU')$ with $\fk(\fV)\subset\fU'\subset\fU$.
 For a treatment of the remaining non-semi-separatedness issue,
we refer to Section~\ref{inverse-images-of-qcoh-tors-subsecn} above
or~\cite[Section~3.3]{Pcosh}.

 The special inverse image functor $\fk^!\:\fX\Tors\rarrow\fZ\Tors$
is right adjoint to the (exact) direct image functor
$\fk_*\:\fZ\Tors\rarrow\fX\Tors$, for the same reasons as explained
in~\cite[Section~7.3]{Pcosh}.
 The key point is that it suffices to define a morphism of
quasi-coherent torsion sheaves on $\fZ$ on the modules of sections
over the affine open formal subschemes $\fk^{-1}(\fU)\subset\fZ$,
where $\fU$ ranges over the affine open formal subschemes of~$\fX$,
with the obvious compatibilities between such morphisms on the modules
of sections; essentially because morphisms of quasi-coherent torsion
sheaves on $\fk^{-1}(\fU)$ correspond bijectively to morphisms of
the related modules of sections over $\fk^{-1}(\fU)$.

 It follows that the functor~$\fk^!$ takes injective quasi-coherent
torsion sheaves on $\fX$ to injective quasi-coherent torsion sheaves
on~$\fZ$ (the same conclusion can be obtained using
Corollary~\ref{local-section-criterion-of-injectivity-cor}).
 The special inverse image functor~$\fk^!$ also preserves infinite
direct sums, as one can see from the construction.

\medskip

 Let $\fii\:\fZ\rarrow\fX$ be a closed immersion of locally Noetherian
formal schemes.
 Let $\bW$ be an open covering of $\fX$ and $\bT$ be an open covering
of $\fZ$ such that the morphism~$\fii$ is $(\bW,\bT)$\+coaffine.
 Then the functor of special inverse image of $\bW$\+flat
$\bW$\+locally contraherent cosheaves of contramodules
\begin{equation} \label{lcth-W-T-flat-special-inverse-image}
 \fii^*\:\fX\Lcth_\bW^\fl\lrarrow\fZ\Lcth_\bT^\fl
\end{equation}
is constructed as follows.

 Let $\fF$ be a $\bW$\+flat $\bW$\+locally contraherent cosheaf of
contramodules on~$\fX$.
 Then the $\bT$\+flat $\bT$\+locally contraherent cosheaf of
contramodules $\fii^*\fF$ on $\fZ$ is defined by the rule
\begin{equation} \label{spec-inv-image-of-W-flat-lcth-defining-rule}
 (\fii^*\fF)[\fV]=\Hom_{\fO_\fZ(\fii^{-1}(\fU))}\bigl(\fO_\fZ(\fV),\>
 \fO_\fZ(\fii^{-1}(\fU))\ot_{\fO_\fX(\fU)}\fF[\fU]\bigr)
\end{equation}
for any pair of affine open formal subschemes $\fV\subset\fZ$ 
subordinate to $\bT$ and $\fU\subset\fX$ subordinate to $\bW$ such that
$\fk(\fV)\subset\fU$.

 So, in particular, denoting by $i_{\fU,\fii^{-1}(\fU)}\:
\fO_\fX(\fU)\rarrow\fO_\fZ(\fii^{-1}(\fU))$ the surjective map of adic
Noetherian rings corresponding to the morphism of affine Noetherian
formal schemes $\fii|^\fU_{\fii^{-1}(\fU)}\:\fii^{-1}(\fU)\rarrow\fU$,
we have
$$
 \fO_\fZ(\fii^{-1}(\fU))\ot_{\fO_\fX(\fU)}\fF[\fU]=
 i_{\fU,\fii^{-1}(\fU)}^\#(\fF[\fU])
$$
by Lemma~\ref{finpres-topol-ring-map-co-contra-extensn-of-scalars}(b).
 This proves that $\fO_\fZ(\fii^{-1}(\fU))\ot_{\fO_\fX(\fU)}\fF[\fU]$
is a contramodule module over $\fO_\fZ(\fii^{-1}(\fU))$.
 Moreover, $\fO_\fZ(\fii^{-1}(\fU))\ot_{\fO_\fX(\fU)}\fF[\fU]$ is
a flat contraadjusted contramodule module over $\fO_\fZ(\fii^{-1}(\fU))$
by Corollary~\ref{flat-contramodules-contraextension-of-scalars} and
Lemma~\ref{contraadjusted-extension-of-scalars-lemma}(b), since
$\fF[\fU]$ is a flat contraadjusted contramodule $\fO_\fX(\fU)$\+module.
 Consequently, the right-hand side
of~\eqref{spec-inv-image-of-W-flat-lcth-defining-rule} is a flat
contraadjusted contramodule $\fO_\fZ(\fV)$\+module by
Corollary~\ref{colocalization-of-flat-contramodules-cor}.

 The following computation from~\cite[Section~7.3]{Pcosh} shows that
the contramodule $\fO_\fZ(\fV)$\+module $(\fii^*\fF)[\fV]$ is
well-defined by
the formula~\eqref{spec-inv-image-of-W-flat-lcth-defining-rule}.
 Let $\fU'\subset\fU$ be two affine open formal subschemes in~$\fX$.
 Then we have
\begin{multline*}
 \fO_\fZ(\fii^{-1}(\fU'))\ot_{\fO_\fX(\fU')}\fF[\fU'] \\
 \simeq\bigl(\fO_\fZ(\fii^{-1}(\fU))
 \ot_{\fO_\fX(\fU)}\fO_\fX(\fU')\bigr)\ot_{\fO_\fX(\fU')}
 \Hom_{\fO_\fX(\fU)}(\fO_\fX(\fU'),\fF[\fU]) \\
 \simeq\fO_\fZ(\fii^{-1}(\fU))\ot_{\fO_\fX(\fU)}
 \Hom_{\fO_\fX(\fU)}(\fO_\fX(\fU'),\fF[\fU]) \\
 \simeq\Hom_{\fO_\fX(\fU)}\bigl(\fO_\fX(\fU'),\>
 \fO_\fZ(\fii^{-1}(\fU))\ot_{\fO_\fX(\fU)}\fF[\fU]\bigr) \\
 \simeq\Hom_{\fO_\fZ(\fii^{-1}(\fU))}
 \bigl(\fO_\fZ(\fii^{-1}(\fU))\ot_{\fO_\fX(\fU)}\fO_\fX(\fU'),\>
 \fO_\fZ(\fii^{-1}(\fU))\ot_{\fO_\fX(\fU)}\fF[\fU]\bigr) \\
 \simeq\Hom_{\fO_\fZ(\fii^{-1}(\fU))}\bigl(\fO_\fZ(\fii^{-1}(\fU')),\>
 \fO_\fZ(\fii^{-1}(\fU))\ot_{\fO_\fX(\fU)}\fF[\fU]\bigr).
\end{multline*}
 Here the first and the final (fifth) isomorphisms hold by
Lemma~\ref{finite-morphism-preimage-tensor-isomorphism-lemma},
while the third isomorphism (or the combination of the third
and the fourth isomorphisms) is provided by
Lemma~\ref{contraadjusted-extension-of-scalars-lemma}(c).

 The computation in the previous paragraph proves that
the formula~\eqref{spec-inv-image-of-W-flat-lcth-defining-rule}
produces the same contramodule $\fO_\fZ(\fV)$\+module $(\fii^*\fF)[\fV]$
when applied to two pairs of affine open formal subschemes
$(\fV,\fU)$ and $(\fV,\fU')$ with $\fk(\fV)\subset\fU'\subset\fU$.
 For a treatment of the remaining non-semi-separatedness issue,
we refer to Sections~\ref{inverse-images-of-qcoh-tors-subsecn}\+-%
\ref{inverse-images-of-ctrh-contramod-subsecn} above
or~\cite[Section~3.3]{Pcosh}.

 For any $(\bW,\bT)$\+coaffine closed immersion of locally Noetherian
formal schemes $\fii\:\fZ\rarrow\fX$, the special inverse image functor
$\fii^*\:\fX\Lcth_\bW^\fl\rarrow\fZ\Lcth_\bT^\fl$ preserves infinite
direct products.
 Assuming that the morphism~$\fii$ is $(\bW,\bT)$\+affine and
$(\bW,\bT)$\+coaffine, there is a natural adjunction isomorphism
of abelian groups
\begin{equation} \label{lcth-cntrmod-specia-inverse-direct-adjunction}
 \Hom^\fZ(\fii^*\fF,\fQ)\simeq\Hom^\fX(\fF,\fii_!\fQ)
\end{equation}
for any $\bW$\+flat $\bW$\+locally contraherent cosheaf of
contramodules $\fF$ on $\fX$ and any $\bT$\+locally contraherent
cosheaf of contramodules $\fQ$ on~$\fZ$.
 We refer to~\cite[Section~7.3]{Pcosh} for the explanation
of the adjunction.
 The key point is that it suffices to define a morphism of
$\bT$\+locally contraherent cosheaves of contramodules on $\fZ$ on
the modules of cosections over the affine open formal subschemes
$\fii^{-1}(\fU)\subset\fZ$, where $\fU$ ranges over the affine open
formal subschemes of $\fX$ subordinate to~$\bW$; essentially because
morphisms of contraherent cosheaves of contramodules on $\fii^{-1}(\fU)$
correspond bijectively to morphisms of the related modules of cosections
over $\fii^{-1}(\fU)$.
 See also the last paragraph at the end of this section below.

 Notice that, as a special property of all closed immersion
morphisms~$\fii$, for any open covering $\bT$ of the formal scheme
$\fZ$ there exists an open covering $\bW$ of the formal scheme $\fX$
such that  the morphism~$\fii$ is $(\bW,\bT)$\+affine (since 
the topology on $\fZ$ is the restriction of the topology on~$\fX$).
 Moreover, the formal scheme $\fZ$ has an open covering $\bT'$
subordinate to $\bT$ such that the morphism~$\fii$ is
$(\bW,\bT')$\+affine \emph{and} $(\bW,\bT')$\+coaffine (it suffices
to let $\bT'$ consist of all the preimages $\fii^{-1}(\fW)$ in $\fZ$ of
the affine open formal subschemes $\fW\subset\fX$ subordinate to~$\bW$).

 Furthermore, whenever the locally Noetherian formal scheme $\fX$ has
finite Krull dimension, we have $\fX\Lcth_\bW^\fl=\fX\Ctrh^\fl$ and
$\fZ\Lcth_\bT^\fl=\fZ\Ctrh^\fl$ by
Corollary~\ref{W-flat-implies-flat-and-coflasque-cor}(b).
 In this case, it follows that the adjunction
isomorphism~\eqref{lcth-cntrmod-specia-inverse-direct-adjunction}
holds for any flat contraherent cosheaf of contramodules $\fF$
on $\fX$ and any locally contraherent cosheaf of contramodules $\fQ$
on~$\fZ$.

\medskip

 Let $\fk\:\fZ\rarrow\fX$ be a finite morphism of locally Noetherian
formal schemes.
 Then the functor of special inverse image of projective locally
cotorsion contraherent cosheaves of contramodules
\begin{equation} \label{lct-prj-special-inverse-image}
 \fk^*\:\fX\Ctrh^\lct_\prj\lrarrow\fZ\Ctrh^\lct_\prj
\end{equation}
is constructed as follows.

 Let $\fF$ be a projective locally cotorsion contraherent cosheaf of
contramodules on~$\fX$.
 Then the projective locally cotorsion contraherent cosheaf of
contramodules $\fk^*\fF$ on $\fZ$ is defined by the rule
\begin{equation} \label{spec-inv-image-of-lct-prj-defining-rule}
 (\fk^*\fF)[\fV]=\Hom_{\fO_\fZ(\fk^{-1}(\fU))}\bigl(\fO_\fZ(\fV),\>
 \fO_\fZ(\fk^{-1}(\fU))\ot_{\fO_\fX(\fU)}\fF[\fU]\bigr)
\end{equation}
for any pair of affine open formal subschemes $\fV\subset\fZ$ and
$\fU\subset\fX$ such that $\fk(\fV)\subset\fU$.

 So, in particular, in the above notation $k_{\fU,\fk^{-1}(\fU)}$
for the finite morphism of adic Noetherian rings $\fO_\fX(\fU)\rarrow
\fO_\fZ(\fk^{-1}(\fU))$, we have
$$
 \fO_\fZ(\fk^{-1}(\fU))\ot_{\fO_\fX(\fU)}\fF[\fU]=
 k_{\fU,\fk^{-1}(\fU)}^\#(\fF[\fU])
$$
by Lemma~\ref{finpres-topol-ring-map-co-contra-extensn-of-scalars}(b).
 This proves that $\fO_\fZ(\fk^{-1}(\fU))\ot_{\fO_\fX(\fU)}\fF[\fU]$
is a contramodule module over $\fO_\fZ(\fk^{-1}(\fU))$.
 Moreover, $\fO_\fZ(\fk^{-1}(\fU))\ot_{\fO_\fX(\fU)}\fF[\fU]$ is
a flat cotorsion contramodule module over $\fO_\fZ(\fk^{-1}(\fU))$
by Lemma~\ref{cotorsion-extension-of-scalars-lemma}(a), since
$\fF[\fU]$ is a flat cotorsion contramodule $\fO_\fX(\fU)$\+module
by Corollary~\ref{projective-lct-of-contramods-local-coflasque}(a)
or~\ref{W-flat-lct-coincide-with-projective-lct}.
 Consequently, the right-hand side
of~\eqref{spec-inv-image-of-lct-prj-defining-rule} is a flat cotorsion
contramodule $\fO_\fZ(\fV)$\+module by
Proposition~\ref{colocalization-of-cotorsion-contramodule-prop} and
Corollary~\ref{colocalization-of-flat-contramodules-cor}.

 The contramodule $\fO_\fZ(\fV)$\+module $\fk^*(\fF)[\fV]$ is
well-defined by
the formula~\eqref{spec-inv-image-of-lct-prj-defining-rule} due to
the natural isomorphism of $\fO_\fZ(\fk^{-1}(\fU'))$\+modules
$$
 \fO_\fZ(\fk^{-1}(\fU'))\ot_{\fO_\fX(\fU')}\fF[\fU']
 \simeq\Hom_{\fO_\fZ(\fk^{-1}(\fU))}\bigl(\fO_\fZ(\fk^{-1}(\fU')),\>
 \fO_\fZ(\fk^{-1}(\fU))\ot_{\fO_\fX(\fU)}\fF[\fU]\bigr),
$$
which holds for any two affine open formal subschemes $\fU'\subset\fU$
in~$\fX$.
 This isomorphism is provable by the same computation as in the case
of a closed immersion morphism~$\fii$ above.
 One needs to use Lemma~\ref{cotorsion-extension-of-scalars-lemma}(c).
 Once again, for a treatment of the remaining non-semi-separatedness
issue, we refer to
Sections~\ref{inverse-images-of-qcoh-tors-subsecn}\+-%
\ref{inverse-images-of-ctrh-contramod-subsecn} above
or~\cite[Section~3.3]{Pcosh}.

 All the discussion above only defines the contramodule
$\fO_\fZ(\fV)$\+modules $\fk^*(\fF)[\fV]$, in such a way that
the contraherence axiom is satisfied, for affine open formal subschemes
$\fV\subset\fZ$ for which there exists an affine open formal subscheme
$\fU\subset\fX$ such that $\fk(\fV)\subset\fU$.
 However, under the latter assumption we have seen that the contramodule
$\fO_\fZ(\fV)$\+module $\fk^*(\fF)[\fV]$ is flat and cotorsion.
 It remains to point out that any locally contraherent cosheaf of
contramodules that is locally flat and locally cotorsion on a locally
Noetherian formal scheme $\fZ$ is globally contraherent and projective
locally cotorsion by
Corollary~\ref{projective-lct-of-contramods-local-coflasque}(d)
or~\ref{W-flat-lct-coincide-with-projective-lct}.

 For any finite morphism of locally Noetherian formal schemes
$\fk\:\fZ\rarrow\fX$, the special inverse image functor
$\fk^*\:\fX\Ctrh^\lct_\prj\rarrow\fZ\Ctrh^\lct_\prj$ preserves infinite
direct products.
 Assuming that the morphism~$\fk$ is $(\bW,\bT)$\+affine, there is
a natural adjunction isomorphism of abelian groups
\begin{equation} \label{lct-prj-special-inverse-direct-adjunction}
 \Hom^\fZ(\fk^*\fF,\fQ)\simeq\Hom^\fX(\fF,\fk_!\fQ)
\end{equation}
for any projective locally cotorsion contraherent cosheaf of
contramodules $\fF$ on $\fX$ and any $\bT$\+locally contraherent
cosheaf of contramodules $\fQ$ on~$\fZ$.

 Once again, we refer to~\cite[Section~7.3]{Pcosh} for the explanation
of the adjunction.
 As above, the key point is that it suffices to define a morphism of
$\bT$\+locally contraherent cosheaves of contramodules on $\fZ$ on
the modules of cosections over the affine open formal subschemes
$\fk^{-1}(\fU)\subset\fZ$, where $\fU$ ranges over the affine open
formal subschemes of $\fX$ subordinate to~$\bW$, with the obvious
compatibilities between such morphisms on the modules of cosections.
 Essentially, this follows from the result about topology bases in
Section~\ref{co-sheaves-of-modules-subsecn} together with the fact that
morphisms of contraherent cosheaves of contramodules on $\fk^{-1}(\fU)$
correspond bijectively to morphisms of the related modules of cosections
over $\fk^{-1}(\fU)$,
by Lemma~\ref{contrah-cosheaves-of-contramods-lemma}(b).

\subsection{Dilute quasi-coherent torsion sheaves}
\label{dilute-subsecn}
 In this section and the next one we return to the setting of
semi-separated Noetherian formal scheme~$\fX$.
 The aim is to extend to formal schemes the theory suggested
in~\cite[Section~4.1]{M-n} and taken up in~\cite[Section~4.2]{Pcosh}.
 We follow the exposition in~\cite{Pcosh}.

 Given a semi-separated Noetherian formal scheme $\fX$, we are
interested in open formal subschemes $\fW\subset\fX$ for which
the open immersion morphism $\fj\:\fW\rarrow\fX$ is affine.
 Notice that if $\fX=\bigcup_\alpha\fU_\alpha$ is an affine open
covering of $\fX$, then $\fW=\bigcup_\alpha(\fW\cap\fU_\alpha)$
is an affine open covering of~$\fW$.
 In particular, if $\fX$ is affine, then so is~$\fW$.
 
 Let $\ff\:\fY\rarrow\fX$ be a morphism of semi-separated Noetherian
formal schemes.
 The fact that, for any open formal subscheme $\fW\subset\fX$ with
an affine open immersion morphism $\fj\:\fW\rarrow\fX$, the open
formal subscheme $\fT=\ff^{-1}(\fW)\subset\fY$ has an affine open
immersion morphism $\fj'\:\fT\rarrow\fY$ plays a key role in
the arguments in this and the next section.
 Furthermore, for any affine open formal subscheme $\fU\subset\fX$,
the open immersion morphism $\fU\rarrow\fX$ is affine.

 We will say that a quasi-coherent torsion sheaf $\cM$ on $\fX$ is
\emph{dilute} if, for every open formal subscheme $\fW\subset\fX$
with affine open immersion morphism $\fj\:\fW\rarrow\fX$, one has
$\boR^n\Gamma(\fW,\fj^*\cM)=0$ for all $n\ge1$.
 Clearly, if this is the case for all $\fW$, then the quasi-coherent
torsion sheaf $\fj^*\cM=\cM|_\fW$ is also dilute (as the compositions
of affine morphisms are affine).
 Furthermore, any quasi-coherent torsion sheaf on an affine formal
scheme $\fU$ is dilute (as the functor $\Gamma(\fU,{-})$ is exact
on $\fU\Tors$).

 We will denote the full subcategory of dilute quasi-coherent torsion
sheaves by $\fX\Tors^\dil\subset\fX\Tors$.
 Obviously, the full subcategory $\fX\Tors^\dil$ is closed under
extensions and cokernels of monomorphisms in $\fX\Tors$; so it inherits
an exact category structure from the abelian exact structure of
$\fX\Tors$.

 The tightness, flatness, and affineness assumptions in the following
lemma can be dropped, as we will see below in this section.

\begin{lem} \label{tight-flat-affine-direct-images-are-dilute}
 Let\/ $\ff\:\fY\rarrow\fX$ be a tight, flat, affine morphism of
semi-separated Noetherian formal schemes.
 Then the direct image functor\/ $\ff_*\:\fY\Tors\rarrow\fX\Tors$
takes dilute quasi-coherent torsion sheaves on\/ $\fY$ to dilute
quasi-coherent torsion sheaves on\/~$\fX$.
\end{lem}

\begin{proof}
 This is the Noetherian formal scheme version
of~\cite[Lemma~4.2.1]{Pcosh}.
 Similarly to the base change argument in~\cite{Pcosh}, it suffices to
show that one $\boR^n\Gamma(\fY,\cN)\simeq\boR^n\Gamma(\fX,\ff_*\cN)$
for any quasi-coherent torsion sheaf $\cN$ on $\fY$ and all $n\ge0$.
 This isomorphism follows from exactness of the functor $\ff_*$
\,\eqref{direct-image-qcoh-tors} and the fact that the functor~$\ff_*$
takes injective objects to injective objects by
Lemma~\ref{injective-torsion-sheaves-direct-image-lem}.
\end{proof}

\begin{cor} \label{fin-iter-exten-of-dirim-from-affopens-dilute}
 Let\/ $\fX$ be a semi-separated Noetherian formal scheme with
a finite affine open covering\/ $\fX=\bigcup_\alpha\fU_\alpha$.
 Then any finitely iterated extension of the direct images of
quasi-coherent torsion sheaves from\/ $\fU_\alpha$ is a dilute
quasi-coherent torsion sheaf on\/~$\fX$.
\end{cor}

\begin{proof}
 This is the Noetherian formal scheme version
of~\cite[Corollary~4.2.2]{Pcosh}.
 The assertion follows from
Lemma~\ref{tight-flat-affine-direct-images-are-dilute} and
the remarks preceding the lemma.
 The assumption that the formal scheme $\fX$ is semi-separated is
essential here.
\end{proof}

\begin{lem} \label{injectives-and-flasques-are-dilute}
 Let\/ $\fX$ be a semi-separated Noetherian formal scheme.
 Then \par
\textup{(a)} any injective quasi-coherent torsion sheaf on\/ $\fX$
is dilute; \par
\textup{(b)} any flasque quasi-coherent torsion sheaf on\/ $\fX$
is dilute.
\end{lem}

\begin{proof}
 Part~(a): this is our version of~\cite[Lemma~4.2.3(a)]{Pcosh}.
 The assertion is deduced by comparing
Lemma~\ref{tight-flat-affine-direct-images-are-dilute} or
Corollary~\ref{fin-iter-exten-of-dirim-from-affopens-dilute}
with Lemma~\ref{injective-torsion-sheaves-on-qcomp}.
 Alternatively, one can use
Corollary~\ref{injective-torsion-local-flasque}(a),
or obtain part~(a) as a special case of part~(b) using
Corollary~\ref{injective-torsion-local-flasque}(b).

 Part~(b): the point is that the restriction of a flasque sheaf of
abelian groups to an open subset is a flasque sheaf of abelian groups
again.
 So the assertion follows from the discussion in the last paragraph
of Section~\ref{homology-of-cosheaves-II-subsecn}.
\end{proof}

\begin{lem} \label{direct-sums-of-dilute}
 For any semi-separated Noetherian formal scheme\/ $\fX$,
the full subcategory of dilute quasi-coherent torsion sheaves is
closed under infinite direct sums in\/ $\fX\Tors$.
\end{lem}

\begin{proof}
 This is the (easy) Noetherian formal scheme version
of~\cite[Corollary~4.2.4]{Pcosh}.
 The point is that the derived functor $\boR^*\Gamma(\fX,{-})$ preserves
infinite direct sums.
 Indeed, the global sections functor\/ $\Gamma(\fX,{-})$ preserves
infinite direct sums (since $\fX$ is quasi-compact; see
Section~\ref{qcoh-torsion-sheaves-subsecn}), and infinite direct sums
of injective quasi-coherent torsion sheaves on $\fX$ are injective
(by Corollary~\ref{qcoh-tors-direct-sums-of-injectives-cor}).
 Furthermore, the restriction functor $\fj^*\:\cM\longmapsto\cM|_\fW$
also preserves infinite direct sums.
\end{proof}

 It follows from Lemma~\ref{injectives-and-flasques-are-dilute}(a)
and the discussion in the beginning of this section that the full
subcategory of dilute quasi-coherent torsion sheaves $\fX\Tors^\dil$
is coresolving in $\fX\Tors$ for any semi-separated Noetherian formal
scheme~$\fX$.
 One also has $\boR^n\Gamma(\fX,\cM)=0$ for any dilute quasi-coherent
torsion sheaf $\cM$ on $\fX$ an all $n\ge1$, by the definition.
 Therefore, the derived functor $\boR^*\Gamma(\fX,{-})$ can be
computed using dilute coresolutions.

\begin{lem} \label{affine-direct-images-are-dilute}
 Let\/ $\ff\:\fY\rarrow\fX$ be an affine morphism of semi-separated
Noetherian formal schemes.
 Then the direct image functor\/ $\ff_*\:\fY\Tors\rarrow\fX\Tors$
takes dilute quasi-coherent torsion sheaves on\/ $\fY$ to dilute
quasi-coherent torsion sheaves on\/~$\fX$.
\end{lem}

\begin{proof}
 This is the Noetherian formal scheme version
of~\cite[Lemma~4.2.5]{Pcosh}.
 Similarly to the proof of
Lemma~\ref{tight-flat-affine-direct-images-are-dilute}
(see~\cite[first paragraph of the proof of Lemma~4.2.1]{Pcosh}),
it suffices to prove an isomorphism $\boR^n\Gamma(\fY,\cN)\simeq
\boR^n\Gamma(\fX,\ff_*\cN)$ for all quasi-coherent torsion sheaves
$\cN$ on $\fY$ and all $n\ge0$.

 Similarly to~\cite[proof of Lemma~4.2.5]{Pcosh}, one observes
that the \v Cech coresolution~\eqref{Cech-sequence-of-torsion-sheaves}
is a coresolution of a quasi-coherent torsion sheaf $\cM$ on $\fX$
by dilute quasi-coherent torsion sheaves (in view of
Lemma~\ref{tight-flat-affine-direct-images-are-dilute}
or Corollary~\ref{fin-iter-exten-of-dirim-from-affopens-dilute}).
 So the functor $\boR^*\Gamma(\fX,{-})$ can be computed using \v Cech
coresolutions with respect to any affine open covering.
 The argument finishes exactly as in~\cite[proof of Lemma~4.2.5]{Pcosh}.
\end{proof}

\begin{lem} \label{direct-images-of-injectives-are-dilute}
 Let\/ $\ff\:\fY\rarrow\fX$ be a morphism of semi-separated Noetherian
formal schemes.
 Then the direct image functor\/ $\ff_*\:\fY\Tors\rarrow\fX\Tors$
takes injective quasi-coherent torsion sheaves on\/ $\fY$ to
dilute quasi-coherent torsion sheaves on\/~$\fX$.
\end{lem}

\begin{proof}
 This is the Noetherian formal scheme version
of~\cite[Lemma~4.2.6]{Pcosh}.
 The assertion follows from
Lemmas~\ref{injective-torsion-sheaves-on-qcomp}
and~\ref{affine-direct-images-are-dilute}.
 The point is that, for any affine open formal subscheme
$\fV\subset\fY$, the composition $\fV\rarrow\fY\rarrow\fX$ is
an affine morphism of semi-separated Noetherian formal schemes
$\fV\rarrow\fX$ (since any morphism from an affine formal scheme
to a semi-separated one is affine).
 Alternatively, one can observe that injective quasi-coherent torsion
sheaves on $\fY$ are flasque by
Corollary~\ref{injective-torsion-local-flasque}(b), the direct
images of flasque sheaves of abelian groups are flasque, and
flasque quasi-coherent torsion sheaves on $\fX$ are dilute by
Lemma~\ref{injectives-and-flasques-are-dilute}(b).
\end{proof}

\begin{cor} \label{dilute-direct-images}
 For any morphism of semi-separated Noetherian formal schemes\/
$\ff\:\fY\rarrow\fX$, the functor of direct image of quasi-coherent
torsion sheaves\/ $\ff_*\:\fY\Tors\allowbreak\rarrow\fX\Tors$ takes
dilute quasi-coherent torsion sheaves on\/ $\fY$ to dilute
quasi-coherent torsion sheaves on\/ $\fX$, and induces an exact
functor\/ $\ff_*\:\fY\Tors^\dil\rarrow\fX\Tors^\dil$.
\end{cor}

\begin{proof}
 This is the Noetherian formal scheme version
of~\cite[Corollary~4.2.7]{Pcosh}, and the proof is the same.
 One needs to use Lemmas~\ref{injectives-and-flasques-are-dilute}(a)
and~\ref{direct-images-of-injectives-are-dilute}, and the fact that
the cokernels of monomorphisms of dilute quasi-coherent torsion
sheaves are dilute. 
\end{proof}

\subsection{Codilute contraherent cosheaves of contramodules}
\label{codilute-subsecn}
 We did not need codilute contraherent cosheaves for the purposes
of the book manuscript~\cite{Pcosh}, because antilocal contraherent
cosheaves could be used instead.
 In the context of formal schemes, we do not seem to have
antilocal contraherent cosheaves of contramodules, so we introduce
the codilute contraherent cosheaves of contramodules instead.

 We continue with the setting of Section~\ref{dilute-subsecn};
so $\fX$ is a semi-separated Noetherian formal scheme.
 We will say that a locally contraherent cosheaf of contramodules
$\fP$ on $\fX$ is \emph{codilute} if, for every open formal subscheme
$\fW\subset\fX$ with affine open immersion morphism
$\fj\:\fW\rarrow\fX$, one has $\boL_n\Delta(\fW,\fj^!\fP)=0$
for all $n\ge1$.
 Clearly, if this is the case for all $\fW$, then the locally
contraherent cosheaf $\fj^!\fP=\fP|_\fW$ is also codilute.

\begin{lem} \label{dilute-lcth-are-globally-contraherent}
\textup{(a)} A locally contraherent cosheaf of contramodules on
an affine Noetherian formal scheme is codilute if and only if
it is contraherent. \par
\textup{(b)} Any codilute locally contraherent cosheaf of contramodules
on a semi-separated Noetherian formal scheme is (globally) contraherent.
\end{lem}

\begin{proof}
 Part~(a) follows directly from
Corollary~\ref{homological-criterion-of-contraherence-cor}.
 Notice that if $\fU$ is an affine formal scheme and $\fj\:\fW\rarrow
\fX$ is an affine morphism of formal schemes, then the formal scheme
$\fW$ is affine.

 Part~(b): let $\fP$ be a codilute locally contraherent cosheaf of
contramodules on a semi-separated Noetherian formal scheme $\fX$, and
let $\fU\subset\fX$ be an affine open formal subscheme.
 Then the open immersion morphism $\fj\:\fU\rarrow\fX$ is affine,
so the locally contraherent cosheaf of contramodules $\fj^!\fP$ on $\fU$
is codilute.
 By part~(a), it follows that $\fj^!\fP$ is a contraherent cosheaf
of contramodules on $\fU$, as desired.
\end{proof}

\begin{lem} \label{coflasque-are-dilute}
\textup{(a)} Any coflasque locally cotorsion contraherent cosheaf of
contramodules on a semi-separated Noetherian formal scheme is codilute.
\par
\textup{(b)} Any coflasque contraherent cosheaf of contramodules on
a semi-separated Noetherian formal scheme of finite Krull dimension
is codilute.
\end{lem}

\begin{proof}
 The restriction of a coflasque cosheaf of abelian groups onto any open
subset is coflasque.
 So the assertions follow from the facts that $\boL_n\Delta(\fX,\fE)=0$
for all $n\ge1$ whenever either $\fX$ is a locally Noetherian formal
scheme and $\fE$ is a coflasque locally cotorsion contraherent cosheaf
of contramodules (for part~(a)), or $\fX$ is a Noetherian formal scheme
of finite Krull dimension and $\fE$ is a coflasque contraherent cosheaf
of contramodules (for part~(b)).
 Those were explained in Section~\ref{homology-of-cosheaves-II-subsecn}.
\end{proof}

 We will denote the full subcategory of codilute contraherent cosheaves
of contramodules by $\fX\Ctrh_\cdl\subset\fX\Ctrh\subset\fX\Lcth$.
 We also put $\fX\Ctrh^\lct_\cdl=\fX\Ctrh_\cdl\cap\fX\Lcth^\lct$;
so $\fX\Ctrh^\lct_\cdl$ is the full subcategory of codilute locally
cotorsion contraherent cosheaves of contramodules.
 Obviously, the full subcategory $\fX\Ctrh_\cdl$ is closed under
extensions and kernels of epimorphisms in $\fX\Ctrh$ and in $\fX\Lcth$;
so it inherits an exact category structure.
 The same applies to the full subcategory $\fX\Ctrh^\lct_\cdl\subset
\fX\Ctrh^\lct\subset\fX\Lcth^\lct$.

 The tightness, very flatness, and affineness assumptions
in the following lemma can be dropped, as we will see below
in this section.

\begin{lem} \label{tight-very-flat-affine-direct-images-are-codilute}
 Let\/ $\ff\:\fY\rarrow\fX$ be a tight, very flat, affine morphism of
semi-separated Noetherian formal schemes.
 Then the direct image functor\/ $\ff_!\:\fY\Ctrh\rarrow\fX\Ctrh$
takes codilute contraherent cosheaves of contramodules on\/ $\fY$ to 
codilute contraherent cosheaves of contramodules on\/~$\fX$.
\end{lem}

\begin{proof}
 This is dual-analogous to
Lemma~\ref{tight-flat-affine-direct-images-are-dilute}.
 In view of the same base change argument, the assertion follows from
the isomorphism~\eqref{veryflat-affine-dirim-cosheaf-homology-iso}
from Section~\ref{homology-of-cosheaves-I-subsecn}.
\end{proof}

\begin{cor} \label{fin-iter-exten-of-dirim-from-affopens-codilute}
 Let\/ $\fX$ be a semi-separated Noetherian formal scheme with
a finite affine open covering\/ $\fX=\bigcup_\alpha\fU_\alpha$.
 Then any finitely iterated extension of the direct images of
contraherent cosheaves of contramodules from\/ $\fU_\alpha$ is
a codilute contraherent cosheaf of contramodules on\/~$\fX$.
\end{cor}

\begin{proof}
 This is dual-analogous to
Corollary~\ref{fin-iter-exten-of-dirim-from-affopens-dilute}.
 The assertion follows from the ``if'' implication
of Lemma~\ref{dilute-lcth-are-globally-contraherent}(a) together
Lemma~\ref{tight-very-flat-affine-direct-images-are-codilute} and
the remarks preceding the latter lemma.
\end{proof}

\begin{lem} \label{antilocally-flat-are-codilute}
 Any antilocally flat contraherent cosheaf of contramodules on
a semi-separated Noetherian formal scheme is codilute.
\end{lem}

\begin{proof}
 This is dual-analogous to
Lemma~\ref{injectives-and-flasques-are-dilute}(a).
 The assertion is deduced by comparing
Corollary~\ref{fin-iter-exten-of-dirim-from-affopens-codilute}
with Corollary~\ref{antilocally-flat-corollary}(c)
(see also Corollary~\ref{antilocally-flat-independent-on-covering}).
\end{proof}

 As particular cases of Lemma~\ref{antilocally-flat-are-codilute},
all projective contraherent cosheaves of contramodules and all
projective locally cotorsion contraherent cosheaves of contramodules
on a semi-separated Noetherian formal scheme are codilute.

\begin{lem} \label{products-of-codilute}
 For any semi-separated Noetherian formal scheme\/ $\fX$, the full
subcategory of codilute contraherent cosheaves of contramodules is
closed under infinite direct products in\/ $\fX\Ctrh$.
\end{lem}

\begin{proof}
 This is dual-analogous to
Lemma~\ref{direct-sums-of-dilute}.
 The point is that the derived functor $\boL_*\Delta(\fX,{-})$
preserves infinite products.
 Indeed, the global cosections functor\/ $\Delta(\fX,{-})$ preserves
infinite products (since $\fX$ is quasi-compact; see
Section~\ref{contraherent-cosheaves-of-contramods-subsecn}), one can
compute the derived functor $\boL_*\Delta(\fX,{-})$ using antilocally
flat resolutions (by Section~\ref{homology-of-cosheaves-I-subsecn}),
and infinite products of antilocally flat contraherent cosheaves of
contramodules on $\fX$ are antilocally flat
(by Corollary~\ref{antilocally-flat-products}).
 Furthermore, the restriction functor $\fj^!\:\fP\longmapsto\fP|_\fW$
also preserves infinite products.
\end{proof}

 It follows from Lemma~\ref{antilocally-flat-are-codilute} and
the discussion above in this section that the full
subcategory of codilute contraherent cosheaves of contramodules
$\fX\Ctrh_\cdl$ is resolving in $\fX\Lcth$ for any semi-separated
Noetherian formal scheme~$\fX$.
 One also has $\boL_n\Delta(\fX,\fP)=0$ for any codilute contraherent
cosheaf of contramodules $\fP$ on $\fX$ an all $n\ge1$, by
the definition.
 Therefore, the derived functor of global cosections
$\boL_*\Delta(\fX,{-})$ on the exact category of locally contraherent
cosheaves of contramodules $\fX\Lcth$ can be computed using codilute
resolutions.

\begin{lem} \label{affine-direct-images-are-codilute}
 Let\/ $\ff\:\fY\rarrow\fX$ be an affine morphism of semi-separated
Noetherian formal schemes.
 Then the direct image functor\/ $\ff_!\:\fY\Ctrh\rarrow\fX\Ctrh$
takes codilute contraherent cosheaves of contramodules on\/ $\fY$ to
codilute contraherent cosheaves of contramodules on\/~$\fX$.
\end{lem}

\begin{proof}
 This is dual-analogous to
Lemma~\ref{affine-direct-images-are-dilute}.
 The question reduces to showing that $\boL_n\Delta(\fY,\fQ)\simeq
\boL_n\Delta(\fX,\ff_!\fQ)$ for all contraherent cosheaves of
contramodules $\fQ$ on $\fY$ and all $n\ge0$.

 Let $\fX=\bigcup_\alpha\fU_\alpha$ be a finite affine open
covering of~$\fX$.
 For any contraherent cosheaf of contramodules $\fP$ on $\fX$,
the \v Cech resolution~\eqref{Cech-sequence-of-cosheaves-of-contramods}
is a resolution of $\fP$ by codilute contraherent cosheaves of
contramodules
(by Lemma~\ref{tight-very-flat-affine-direct-images-are-codilute}
or Corollary~\ref{fin-iter-exten-of-dirim-from-affopens-codilute}).
 So one can compute the derived functor $\boL_*\Delta(\fX,{-})$
on the exact category $\fX\Ctrh$ using the \v Cech resolutions with
respect to finite affine open coverings.
 (Given an open covering $\bW$ of $\fX$ and a finite affine open
covering $\fX=\bigcup_\alpha\fU_\alpha$ subordinate to $\bW$,
the same applies to the derived functor $\boL_*\Delta(\fX,{-})$
on the exact category $\fX\Lcth_\bW$.)

 Put $\fV_\alpha=\ff^{-1}(\fU_\alpha)\subset\fY$; then
$\fY=\bigcup_\alpha\fV_\alpha$ is a finite affine open covering
of $\fY$, and one can compute the derived functor
$\boL_*\Delta(\fY,{-})$ using the \v Cech resolutions for this
covering.
 It remains to observe that the complex of abelian groups
$C_\bu(\{\fU_\alpha\},\ff_!\fQ)$
\,\eqref{Cech-complex-of-abelian-groups-for-copresheaf} obtained by
applying the functor $\Delta(\fX,{-})$ to the \v Cech resolution
of $\ff_!\fQ$ with respect to the open covering
$\fX=\bigcup_\alpha\fU_\alpha$ is naturally isomorphic to the complex
of abelian groups $C_\bu(\{\fV_\alpha\},\fQ)$ obtained by applying
the functor $\Delta(\fY,{-})$ to the \v Cech resolution of $\fQ$ with
respect to the open covering $\fY=\bigcup_\alpha\fV_\alpha$.
\end{proof}

\begin{lem} \label{direct-images-of-projectives-are-codilute}
 Let\/ $\ff\:\fY\rarrow\fX$ be a morphism of semi-separated Noetherian
formal schemes.
 Then \par
\textup{(a)} the functor of direct image of cosheaves of\/
$\fO$\+modules\/ $\ff_!\:(\fY,\fO_\fY)\Cosh\rarrow(\fX,\fO_\fX)\Cosh$
takes projective contraherent cosheaves of contramodules on\/ $\fY$ to
codilute contraherent cosheaves of contramodules on\/~$\fX$; \par
\textup{(b)} the functor of direct image of cosheaves of\/
$\fO$\+modules\/ $\ff_!\:(\fY,\fO_\fY)\Cosh\rarrow(\fX,\fO_\fX)\Cosh$
takes projective locally cotorsion contraherent cosheaves of
contramodules on\/ $\fY$ to codilute locally cotorsion contraherent
cosheaves of contramodules on\/~$\fX$.
\end{lem}

\begin{proof}
 This is dual-analogous to
Lemma~\ref{direct-images-of-injectives-are-dilute}.
 The assertion follows from the description of projective (or
projective locally cotorsion) contraherent cosheaves of contramodules
provided by Lemma~\ref{projective-contrah-on-qcomp-ssep} (or
Lemma~\ref{projective-contrah-lct-on-qcomp-ssep}) together with
the result of Lemma~\ref{affine-direct-images-are-codilute}.
 Once again the point is that any morphism from an affine formal
scheme to a semi-separated formal scheme is affine.
 Notice also that, for an affine morphism of formal schemes,
the functor of direct image of contraherent cosheaves of
contramodules~\eqref{ctrh-contramod-affine-direct-image} takes
locally cotorsion contraherent cosheaves of contramodules to
locally cotorsion contraherent cosheaves of
contramodules~\eqref{ctrh-lct-contramod-affine-direct-image}.
\end{proof}

\begin{lem} \label{direct-images-of-codilute-are-contraherent}
 Let\/ $\fY$ be a semi-separated Noetherian formal scheme, $\fU$ be
an affine Noetherian formal scheme, and\/ $\ff\:\fY\rarrow\fU$ be
a morphism of formal schemes.
 Then \par
\textup{(a)} the functor of direct image of cosheaves of\/
$\fO$\+modules\/ $\ff_!\:(\fY,\fO_\fY)\Cosh\rarrow(\fU,\fO_\fU)\Cosh$
takes codilute contraherent cosheaves of contramodules on\/ $\fY$ to
contraherent cosheaves of contramodules on\/ $\fU$, and induces
an exact functor\/ $\ff_!\:\fY\Ctrh_\cdl\rarrow\fU\Ctrh$;
{\hbadness=1200\par}
\textup{(b)} the functor of direct image of cosheaves of\/
$\fO$\+modules\/ $\ff_!\:(\fY,\fO_\fY)\Cosh\rarrow(\fU,\fO_\fU)\Cosh$
takes codilute locally cotorsion contraherent cosheaves of contramodules
on\/ $\fY$ to locally cotorsion contraherent cosheaves of contramodules
on\/ $\fX$, and induces an exact functor\/ $\ff_!\:\fY\Ctrh^\lct_\cdl
\rarrow\fU\Ctrh^\lct$.
\end{lem}

\begin{proof}
 Let $\fY=\bigcup_\alpha\fV_\alpha$ be a finite affine open covering,
and let $\fQ$ be a codilute contraherent cosheaf of contramodules
on~$\fY$.
 Then the \v Cech exact
sequence~\eqref{Cech-sequence-of-cosheaves-of-contramods} with respect
to the open covering $\{\fV_\alpha\}$ is a finite resolution of
a codilute contraherent cosheaf of contramodules $\fQ$ by codilute
contraherent cosheaves of contramodules on~$\fY$ (as per the proof of
Lemma~\ref{affine-direct-images-are-codilute}).
 So the \v Cech complex of $\fO_\fU(\fU)$\+modules
$C_\bu(\{\fV_\alpha\},\fQ)\rarrow\fQ[\fY]\rarrow0$ obtained by applying
the functor $\Delta(\fY,{-})$ to this \v Cech resolution of $\fQ$
is exact.

 All the terms of the resulting \v Cech complex of
$\fO_\fU(\fU)$\+modules above, except perhaps the rightmost one,
are contraadjusted contramodule $\fO_\fU(\fU)$\+modules, by
Lemma~\ref{restriction-of-scalars-contraadjusted} (or even
cotorsion contramodule $\fO_\fU(\fU)$\+modules, if $\fQ$ is a codilute
locally cotorsion contraherent cosheaf of contramodules in the context
of part~(b), by Lemma~\ref{restriction-of-scalars-cotorsion}).
 Hence it follows that the $\fO_\fU(\fU)$\+module $\fQ[\fY]$ is
a contraadjusted contramodule (resp., a cotorsion contramodule).
 The contraherence axiom for the cosheaf $\ff_!\fQ$ on $\fU$ is
established by the argument from the next-to-last paragraph of
the proof of Corollary~\ref{co-flasque-direct-image-cor}(b\+-c) above.

 Finally, exactness of the functor $\ff_!\:\fY\Ctrh_\cdl\rarrow\fU\Ctrh$
follows from the fact that the degree~$0$ homology of a complex of
abelian groups having no homology in all the other degrees except~$0$
is an exact functor.
\end{proof}

\begin{cor} \label{codilute-direct-images}
 Let\/ $\ff\:\fY\rarrow\fX$ be a morphism of semi-separated Noetherian
formal schemes.
 Then \par
\textup{(a)} the functor of direct image of cosheaves of\/
$\fO$\+modules\/ $\ff_!\:(\fY,\fO_\fY)\Cosh\rarrow(\fX,\fO_\fX)\Cosh$
takes codilute contraherent cosheaves of contramodules on\/ $\fY$ to
codilute contraherent cosheaves of contramodules on\/ $\fX$, and induces
an exact functor\/ $\ff_!\:\fY\Ctrh_\cdl\rarrow\fX\Ctrh_\cdl$; \par
\textup{(b)} the functor of direct image of cosheaves of\/
$\fO$\+modules\/ $\ff_!\:(\fY,\fO_\fY)\Cosh\rarrow(\fX,\fO_\fX)\Cosh$
takes codilute locally cotorsion contraherent cosheaves of contramodules
on\/ $\fY$ to codilute locally cotorsion contraherent cosheaves of
contramodules on\/ $\fX$, and induces an exact functor\/
$\ff_!\:\fY\Ctrh^\lct_\cdl\rarrow\fX\Ctrh^\lct_\cdl$.
\end{cor}

\begin{proof}
 This is dual-analogous to Corollary~\ref{dilute-direct-images}.
 This is also a more general version
of~\cite[Corollary~4.6.3(a\+-b)]{Pcosh}.
 The assertion and proof bear some similarity to those of
Corollaries~\ref{ssep-flat-direct-image-of-antilocally-flat-cor}
and~\ref{co-flasque-direct-image-cor}(b\+-c) (we have already seen
some of that in the proof of the preceding
Lemma~\ref{direct-images-of-codilute-are-contraherent}).

 Let us prove part~(a).
 It follows from Lemma~\ref{direct-images-of-codilute-are-contraherent}
that the direct images of codilute contraherent cosheaves of
contramodules are contraherent, and the functor
$\ff_!\:\fY\Ctrh_\cdl\rarrow\fX\Ctrh_\cdl$ is exact.
 It remains to prove that the functor~$\ff_!$ takes codilute contraherent
cosheaves of contramodules on $\fY$ to codilute contraherent cosheaves
of contramodules on~$\fX$.
 By the very same base change argument from~\cite[first paragraph of
the proof of Lemma~4.2.1]{Pcosh}, it suffices to show that
$\boL_n\Delta(\fX,\ff_!\fQ)=0$ for any codilute contraherent cosheaf
of contramodules $\fQ$ on~$\fY$.

 Let $\fF_\bu\rarrow\fQ$ be a projective resolution of $\fQ$ in
the exact category $\fY\Ctrh$.
 Then the complex $\fF_\bu\rarrow\fQ\rarrow0$ is exact in the exact
category $\fY\Ctrh_\cdl$ (since the kernels of admissible epimorphisms
of codilute contraherent cosheaves of contramodules are codilute, and
projective contraherent cosheaves of contramodules on $\fY$ are
codilute by Lemma~\ref{antilocally-flat-are-codilute}).
 According to the previous paragraph, it follows that the complex
$\ff_!\ff_\bu\rarrow\ff_!\fQ\rarrow0$ is exact in $\fX\Ctrh$.
 By Lemma~\ref{direct-images-of-projectives-are-codilute},
\,$\ff_!\fF_\bu$ is a complex of codilute contraherent cosheaves of
contramodules on~$\fX$; so $\ff_!\fF_\bu$ is a codilute resolution
of $\ff_!\fQ$ in $\fX\Ctrh$.
 Thus the complex $\Delta(\fX,\ff_!\fF_\bu)$ computes
$\boL_*\Delta(\fX,\ff_!\fQ)$.
 Finally, the natural isomorphism of complexes of abelian groups
$\Delta(\fX,\ff_!\fF_\bu)\simeq\Delta(\fY,\fF_\bu)$ induces
an isomorphism of the cohomology groups $\boL_n\Delta(\fX,\ff_!\fQ)
\simeq\boL_n\Delta(\fY,\fQ)$, and we have $\boL_n\Delta(\fY,\fQ)=0$
for $n\ge1$.
\end{proof}

\subsection{Nonaffine projection formulas for the contratensor product}
\label{projection-contratensor-subsecn}
 In this section we extend the results of
Proposition~\ref{contratensor-projection-formula-prop}(b\+-c) to
nonaffine morphisms of nonaffine formal schemes, under certain
assumptions.

\begin{lem} \label{projection-map-for-contratensor-trian-diag-lemma}
 Let\/ $\fg\:\fW\rarrow\fY$ and\/ $\ff\:\fY\rarrow\fX$ be a composable
pair of flat tight quasi-compact morphisms of locally Noetherian
formal schemes.
 Let\/ $\cM$ be a quasi-coherent torsion sheaf on\/ $\fX$ and\/
$\fR$ be a cosheaf of\/ $\fO_\fW$\+modules on\/~$\fW$.
 Then two natural morphisms in\/ $\fX\Tors$ provided by
Proposition~\ref{contratensor-projection-formula-prop}(a) and
the direct image under\/~$\ff$ of a natural morphism in\/ $\fY\Tors$
provided by the same proposition form a commutative triangular
diagram
$$
 \xymatrix{
  \cM\ocn_\fX\ff_!\fg_!\fR \ar[rr] \ar[rd]
  && \ff_*\fg_*(\fg^*\ff^*\cM\ocn_\fW\fR) \\
  & \ff_*(\ff^*\cM\ocn_\fY\fg_!\fR) \ar[ru]
 }
$$
of morphisms of quasi-coherent torsion sheaves on\/~$\fX$.
\end{lem}

\begin{proof}
 Instead of spelling out tiresome details of verification of
the commutativity of the diagram, we suggest a simpler alternative
construction of the morphism of quasi-coherent torsion sheaves from
Proposition~\ref{contratensor-projection-formula-prop}(a),
which is only applicable in the case of a \emph{flat} tight
quasi-compact morphism of formal schemes $\ff\:\fY\rarrow\fX$.
 Hence the flatness assumption in the formulation of the present lemma
(which is otherwise likely to be unnecessary).

 For any quasi-coherent torsion sheaf $\cM$ on $\fX$ and any
cosheaf of $\fO_\fY$\+modules $\fQ$ on $\fY$, the construction
of Proposition~\ref{contratensor-projection-formula-prop}(a)
provides a natural morphism $\cM\ocn_\fX\ff_!\fQ\rarrow
\ff_*(\ff^*\cM\ocn_\fY\fQ)$
\,\eqref{contratensor-projection-formula-map-eqn}
of quasi-coherent torsion sheaves on~$\fX$.
 In view of the adjunction~\eqref{qcoh-tors-inverse-direct-adjunction},
the datum of a morphism~\eqref{contratensor-projection-formula-map-eqn}
is equivalent to the datum of a morphism
\begin{equation} \label{contratensor-inverse-image-map-eqn}
 \ff^*(\cM\ocn_\fX\ff_!\fQ)\lrarrow\ff^*\cM\ocn_\fY\fQ
\end{equation}
of quasi-coherent torsion sheaves on~$\fY$.
 Let us explain how to construct the natural
morphism~\eqref{contratensor-inverse-image-map-eqn}.

 According to the construction in
Section~\ref{contratensor-product-subsecn}, we have
\begin{multline*}
 \ff^*(\cM\ocn_\fX\ff_!\fQ)=
 \ff^*\varinjlim\nolimits_{\fU\subset\fX}
 \bigl((\fj_*\fj^*\cM)\ot_{\fO_\fX(\fU)}(\ff_!\fQ)[\fU]\bigr)
 \\ \simeq
 \varinjlim\nolimits_{\fU\subset\fX}\ff^*
 \bigl((\fj_*\fj^*\cM)\ot_{\fO_\fX(\fU)}(\ff_!\fQ)[\fU]\bigr)
 \\ \simeq
 \varinjlim\nolimits_{\fU\subset\fX}
 \bigl((\ff^*\fj_*\fj^*\cM)\ot_{\fO_\fX(\fU)}
 (\ff_!\fQ)[\fU]\bigr),
\end{multline*}
where $\fU$ ranges over all the affine open formal subschemes in $\fX$
and $\fj\:\fU\rarrow\fX$ denotes the open immersion morphism.
 Here we are using the facts that the functor~$\ff^*$, being a left
adjoint functor, commutes with all colimits, as well as with all
tensor products of quasi-coherent torsion sheaves with modules over
rings (in the sense of Section~\ref{contraherent-Hom-subsecn}).

 Consider the open formal subscheme $\ff^{-1}(\fU)\subset\fY$.
 Denote by $\fj'\:\ff^{-1}(\fU)\rarrow\fY$ the identity open immersion
morphism and by $\ff'\:\ff^{-1}(\fU)\rarrow\fU$ the natural flat tight
quasi-compact morphism.
 Then we have
$$
 \ff^*\fj_*\fj^*\cM\simeq\fj'_*\ff'{}^*\fj^*\cM
 \simeq \fj'_*\fj'{}^*\ff^*\cM
$$
by Lemma~\ref{tors-openimmers-dir-tight-inv-basechange-lemma}(b).
 Furthermore, by the definition and by
Lemma~\ref{cosheaf-cosections-as-inductive-image-lemma},
$$
 (\ff_!\fQ)[\fU]=\fQ[\ff^{-1}(\fU)]\simeq
 \varinjlim\nolimits_{\fV\subset\ff^{-1}(\fU)}\fQ[\fV],
$$
where $\fV$ ranges over all the affine open formal subschemes in
$\ff^{-1}(\fU)$.
 Hence
$$
 \ff^*(\cM\ocn_\fX\ff_!\fQ)\simeq
 \varinjlim\nolimits_{\fU\subset\fX}
 \varinjlim\nolimits_{\fV\subset\ff^{-1}(\fU)}
 \bigl((\fj'_*\fj'{}^*\ff^*\cM)\ot_{\fO_\fX(\fU)}\fQ[\fV]\bigr).
$$

 Set $\cN=\ff^*\cM$.
 Consider the poset $\bD$ of all pairs of affine open formal subschemes
$\fU\subset\fX$ and $\fV\subset\fY$ such that $\ff(\fV)\subset\fU$
(ordered by inclusion of both the affine open formal subschemes in
the pair).
 Then the diagram of affine open formal subschemes $\fV\subset\fY$
indexed by $\bD$ is exhaustive in $\fY$, therefore by construction
$$
 \cN\ocn_\fY\fQ\simeq\varinjlim\nolimits_{\fU\subset\fX}
 \varinjlim\nolimits_{\fV\subset\ff^{-1}(\fU)}
 \bigl((\fk_*\fk^*\cN)\ocn_{\fO_\fY(\fV)}\fQ[\fV]\bigr),
$$
where $\fk\:\fV\rarrow\fY$ is the open immersion morphism.
 Finally, the obvious diagram of adjunction morphisms of quasi-coherent
torsion sheaves $\fj'_*\fj'{}^*\cN\rarrow\fk_*\fk^*\cN$ on $\fY$,
together with the diagram of ring homomorphisms $\fO_\fX(\fU)\rarrow
\fO_\fY(\fV)$, both indexed by $\bD$, induce the desired morphism of
the direct limits over~$\bD$
\begin{multline*}
 \ff^*(\cM\ocn_\fX\ff_!\fQ)\simeq
 \varinjlim\nolimits_{\fU\subset\fX}
 \varinjlim\nolimits_{\fV\subset\ff^{-1}(\fU)}
 \bigl((\fj'_*\fj'{}^*\cN)\ot_{\fO_\fX(\fU)}\fQ[\fV]\bigr) \\
 \lrarrow\varinjlim\nolimits_{\fU\subset\fX}
 \varinjlim\nolimits_{\fV\subset\ff^{-1}(\fU)}
 \bigl((\fk_*\fk^*\cN)\ot_{\fO_\fY(\fV)}\fQ[\fV]\bigr)
 \simeq\ff^*\cM\ocn_\fY\fQ
\end{multline*}
in $\fY\Tors$.
\end{proof}

\begin{lem} \label{contratensor-exact-on-flats-and-alf}
\textup{(a)} Let\/ $\fX$ be a Noetherian formal scheme of finite Krull
dimension.
 Then, for any fixed quasi-coherent torsion sheaf\/ $\cM\in\fX\Tors$,
the contratensor product functor\/ $\cM\ocn_\fX{-}\,\:\fX\Ctrh^\fl
\rarrow\fX\Tors$ is exact as a functor on the exact category of flat
contraherent cosheaves of contramodules\/ $\fX\Ctrh^\fl$. \par
\textup{(b)} Let\/ $\fX$ be a Noetherian formal scheme of finite Krull
dimension.
 Then the functor of contratensor product of a flasque quasi-coherent
torsion sheaf with a flat contraherent cosheaf of contramodules\/
$\ocn_\fX\:\fX\Tors^\fq\times\fX\Ctrh^\fl\rarrow\fX\Tors$ is exact
as a functor of both the arguments from the exact category of
flasque quasi-coherent torsion sheaves\/ $\fX\Tors^\fq$ and from
the exact category of flat contraherent cosheaves of contramodules\/
$\fX\Ctrh^\fl$. \par
\textup{(c)} Let\/ $\fX$ be a semi-separated Noetherian formal scheme.
 Then the functor of contratensor product with an antilocally flat
contraherent cosheaf of contramodules\/
$\ocn_\fX\:\fX\Tors\times\fX\Ctrh_\alf\rarrow\fX\Tors$ is exact as
a functor of both the arguments from the abelian category of
quasi-coherent torsion sheaves\/ $\fX\Tors$ and from the exact category
of antilocally flat contraherent cosheaves of contramodules\/
$\fX\Ctrh_\alf$. 
\end{lem}

\begin{proof}
 Part~(a): let $\fF\in\fX\Ctrh^\fl$ be a flat contraherent cosheaf of
contramodules on~$\fX$.
 There are enough injective objects in $\fX\Tors$, hence in order to
prove that the functor $\fF\longmapsto\cM\ocn_\fX\fF$ is exact,
it suffices to check that the functor
$$
 \fF\longmapsto\Hom_\fX(\cM\ocn_\fX\fF,\,\cJ)
$$
is exact on the exact category $\fX\Ctrh^\fl$ for any fixed
quasi-coherent torsion sheaf $\cM$ and any fixed injective
quasi-coherent torsion sheaf $\cJ$ on~$\fX$.
 Now, according to the adjunction
isomorphism~\eqref{fHom-contratensor-adjunction-eqn}, we have
\begin{equation} \label{fHom-contratensor-adjunction-rewritten}
 \Hom_\fX(\cM\ocn_\fX\fF,\,\cJ)\simeq
 \Hom^\fX(\fF,\fHom_\fX(\cM,\cJ)).
\end{equation}
 The cosheaf $\fHom_\fX(\cM,\cJ)$ is a locally cotorsion contraherent
cosheaf of contramodules on $\fX$ by construction (see
Section~\ref{contraherent-Hom-subsecn}).
 It remains to point out that, for any locally cotorsion contraherent
cosheaf of contramodules $\fP$ on $\fX$, the functor $\fF\longmapsto
\fHom_\fX(\fF,\fP)$ is exact on $\fX\Ctrh^\fl$ by
Corollary~\ref{flat-contraherent-cosheaves-cor}(a).

 Part~(c): similarly to the proof of part~(a), there are enough
injective objects in $\fX\Tors$, hence it suffices to check that
the functor
$$
 (\cM,\fF)\longmapsto\Hom_\fX(\cM\ocn_\fX\fF,\,\cJ)
$$
is exact as a functor of two arguments $\cM\in\fX\Tors$ and
$\fF\in\fX\Ctrh_\alf$ for any fixed injective
quasi-coherent torsion sheaf $\cJ$ on~$\fX$.
 Once again, we have the adjunction
isomorphism~\eqref{fHom-contratensor-adjunction-rewritten},
and $\fHom_\fX(\cM,\cJ)$ is a locally cotorsion contraherent
cosheaf of contramodules on $\fX$ by construction.
 By Corollary~\ref{antilocally-flat-characterized-cor}(a), the functor
$(\fF,\fP)\longmapsto\Hom^\fX(\fF,\fP)$ is exact as a functor of two
arguments $\fF\in\fX\Ctrh_\alf$ and $\fP\in\fX\Ctrh^\lct$.
 Similarly to the proof of part~(a), this is sufficient to show that
the contratensor product functor $\fF\longmapsto\cM\ocn_\fX\fF$
is exact on the exact category $\fX\Ctrh_\alf$ for any fixed object
$\cM\in\fX\Tors$.

 To prove that the contratensor product functor
$\cM\longmapsto\cM\ocn_\fX\fF$ is exact on the abelian category
$\fX\Tors$ for any fixed object $\fF\in\fX\Ctrh_\alf$, it remains
to point out that, over a semi-separated Noetherian formal scheme
$\fX$, the functor $\cM\longmapsto\fHom_\fX(\cM,\cJ)$ is exact
as a functor of a quasi-coherent torsion sheaf $\cM\in\fX\Tors$
for any fixed injective quasi-coherent torsion sheaf
$\cJ\in\fX\Tors^\inj$.
 This observation can be made immediately from the construction
of the contrareherent cosheaf of contramodules $\fHom_\fX(\cM,\cJ)$
on~$\fX$.
 The point is that, for any affine open subscheme $\fU\subset\fX$,
the open immersion morphism $\fj\:\fU\rarrow\fX$ is affine, and
the direct image functor $\fj_*\:\fU\Tors\rarrow\fX\Tors$ is
exact for an affine morphism~$\fj$.

 The proof of part~(b) is similar to that of parts~(a) and~(c).
 The point is that, for any open subscheme $\fU\subset\fX$ with
the open immersion morphism $\fj\:\fU\rarrow\fX$, the functor
$\fj_*\fj^*\:\fX\Tors^\fq\rarrow\fX\Tors^\fq$ is exact, essentially
by Corollary~\ref{co-flasque-direct-image-cor}(a).
 Hence the functor $\cE\longmapsto\fHom_\fX(\cE,\cJ)$ is exact
as a functor of a flasque quasi-coherent torsion sheaf
$\cE\in\fX\Tors^\fq$ for any fixed injective quasi-coherent torsion 
sheaf $\cJ\in\fX\Tors^\inj$.
 Then one needs to use the fact that the functor $(\fF,\fP)
\longmapsto\Hom^\fX(\fF,\fP)$ is exact as a functor of two arguments
$\fF\in\fX\Ctrh^\fl$ and $\fP\in\fX\Ctrh^\lct$ by
Corollary~\ref{flat-contraherent-cosheaves-cor}(a).
\end{proof}

\begin{lem} \label{fHom-and-contratensor-with-alf-are-co-dilute}
 Let\/ $\fX$ be a semi-separated Noetherian formal scheme
and\/ $\cM$ be a quasi-coherent torsion sheaf on\/~$\fX$.
 Then \par
\textup{(a)} for any injective quasi-coherent torsion sheaf\/ $\cJ$
on\/ $\fX$, the locally cotorsion contraherent cosheaf of
contramodules\/ $\fHom_\fX(\cM,\cJ)$ on\/ $\fX$ is codilute; \par
\textup{(b)} for any antilocally flat contraherent cosheaf of
contramodules\/ $\fF$ on\/ $\fX$, the quasi-coherent
torsion sheaf\/ $\cM\ocn_\fX\fF$ on\/ $\fX$ is dilute.
\end{lem}

\begin{proof}
 Part~(a): let $\fX=\bigcup_\alpha\fU_\alpha$ be a finite affine open
covering of~$\fX$.
 By Lemma~\ref{injective-torsion-sheaves-on-qcomp}, \,$\cJ$ is
a direct summand of a finite direct sum of the direct images of
injective quasi-coherent torsion sheaves from~$\fU_\alpha$.
 Alternatively, it follows from
Theorem~\ref{loc-Noetherian-injective-qcoh-torsion}(b) that $\cJ$ is
a finite direct sum of the direct images of injective quasi-coherent
torsion sheaves from~$\fU_\alpha$.
 Hence it suffices to consider the case of $\cJ=\fj_*\cK$, where
$\fj\:\fU\rarrow\fX$ is the immersion of an affine open formal
subscheme and $\cK\in\fU\Tors^\inj$.
 Then we have $\fHom_\fX(\cM,\fj_*\cK)\simeq
\fj_!\fHom_\fU(\fj^*\cM,\cK)$ by
Proposition~\ref{fHom-projection-formula-prop}, and it remains to
recall that the direct image of any contraherent cosheaf of
contramodules from $\fU$ is codilute on $\fX$ by
Lemma~\ref{tight-very-flat-affine-direct-images-are-codilute} or
Corollary~\ref{fin-iter-exten-of-dirim-from-affopens-codilute}.

 Part~(b): by Corollary~\ref{antilocally-flat-corollary}(c),
\,$\fF$ is a direct summand of a finitely iterated extension of
the direct images of antilocally flat contraherent cosheaves of
contramodules from~$\fU_\alpha$.
 By Lemma~\ref{contratensor-exact-on-flats-and-alf}(c), the functor
$\cM\ocn_\fX{-}\,\:\fX\Ctrh_\alf\rarrow\fX\Tors$ is exact; so it
preserves finitely iterated extensions.
 Since the class of dilute quasi-coherent torsion sheaves is closed
under extensions in $\fX\Tors$, it suffices to consider the case of
$\fF=\fj_!\fH$, where $\fj\:\fU\rarrow\fX$ is the immersion of
an affine open formal subscheme and $\fH\in\fU\Ctrh_\alf$.
 Then we have $\cM\ocn_\fX\fj_!\fH\simeq\fj_*(\fj^*\cM\ocn_\fU\fH)$ by
Proposition~\ref{contratensor-projection-formula-prop}(b) or~(c),
and it remains to recall that the direct image of any quasi-coherent
torsion sheaf from $\fU$ is dilute on $\fX$ by
Lemma~\ref{tight-flat-affine-direct-images-are-dilute} or
Corollary~\ref{fin-iter-exten-of-dirim-from-affopens-dilute}.
\end{proof}

\begin{prop} \label{alf-contratensor-projection-formula-prop}
 Let\/ $\ff\:\fY\rarrow\fX$ be a flat tight morphism of semi-separated
Noetherian formal schemes.
 Let\/ $\cM$ be a quasi-coherent torsion sheaf on\/ $\fX$ and\/ $\fG$
be an antilocally flat contraherent cosheaf of contramodules on\/~$\fY$.
 Then the natural
morphism~\eqref{contratensor-projection-formula-map-eqn} from
Proposition~\textup{\ref{contratensor-projection-formula-prop}(a)} is
an isomorphism
\begin{equation} \label{alf-contratensor-projection-formula-eqn}
 \cM\ocn_\fX\ff_!\fG\simeq\ff_*(\ff^*\cM\ocn_\fY\fG)
\end{equation}
of quasi-coherent torsion sheaves on\/~$\fX$.
\end{prop}

\begin{proof}
 This is a formal scheme version of~\cite[Lemma~8.4.2]{Pcosh}.
 The key observation is that both the left-hand side and the right-hand
side of~\eqref{alf-contratensor-projection-formula-eqn} are exact
functors of the argument $\fF\in\fX\Ctrh_\alf$.
 Concerning the left-hand side, we have an exact functor
$\ff_!\:\fY\Ctrh_\alf\rarrow\fX\Ctrh_\alf$ by
Corollary~\ref{ssep-flat-direct-image-of-antilocally-flat-cor},
and the functor $\cM\ocn_\fX{-}\,\:\fX\Ctrh_\alf\rarrow\fX\Tors$
is exact by Lemma~\ref{contratensor-exact-on-flats-and-alf}(c).

 Concerning the right-hand side, once again the functor
$\cN\ocn_\fY{-}\,\:\fY\Ctrh_\alf\rarrow\fY\Tors$ is exact
for any quasi-coherent torsion sheaf $\cN$ on $\fY$
by Lemma~\ref{contratensor-exact-on-flats-and-alf}(c).
 Furthermore, for any antilocally flat contraherent cosheaf
of contramodules $\fG$ on $\fY$, the quasi-coherent torsion sheaf
$\cN\ocn_\fY\fG$ on $\fY$ is dilute by
Lemma~\ref{fHom-and-contratensor-with-alf-are-co-dilute}(b).
 Finally, we have an exact functor $\ff_*\:\fY\Tors^\dil\rarrow
\fX\Tors^\dil$, by Corollary~\ref{dilute-direct-images}.

 Now that we know that both the sides of the desired
isomorphism~\eqref{alf-contratensor-projection-formula-eqn} are
exact functors of the argument $\fG$, we can use
Corollary~\ref{antilocally-flat-corollary}(c) to the effect that
it suffices to prove the desired isomorphism in the case of
$\fG=\fk_!\fH$, where $\fH$ is an antilocally flat contraherent
cosheaf of contramodules on~$\fV$.
 In that case, the assertion follows from the commutativity of
the triangular diagram in
Lemma~\ref{projection-map-for-contratensor-trian-diag-lemma}.
 The point is that both the morphisms of formal schemes~$\fk$ and
$\ff\circ\fk$ are tight, flat, and affine (since the formal scheme
$\fV$ is affine, and the formal schemes $\fY$ and $\fX$ are
semi-separated).
 So both the morphisms of quasi-coherent torsion sheaves
$\ff^*\cM\ocn_\fY\fk_!\fH\rarrow\fk_*(\fk^*\ff^*\cM\ocn_\fV\fH)$ and
$\cM\ocn_\fX\ff_!\fk_!\fH\rarrow\ff_*\fk_*(\ff^*\fk^*\cM\ocn_\fV\fH)$
are isomorphisms by
Proposition~\ref{contratensor-projection-formula-prop}(b).
 Hence the morphism $\cM\ocn_\fX\ff_!\fk_!\fH\rarrow
\ff_*(\ff^*\cM\ocn_\fY\fk_!\fH)$ is an isomorphism, too.
\end{proof}

 The next lemma is a partial generalization of
Lemma~\ref{fHom-contratensor-co-flasque}.

\begin{lem} \label{fHom-from-flasque-contratensor-with-alf-co-flasque}
 Let\/ $\fX$ be a Noetherian formal scheme and\/ $\cE$ be
a flasque quasi-coherent torsion sheaf on\/~$\fX$.
 In this context: \par
\textup{(a)} for any injective quasi-coherent torsion sheaf\/ $\cJ$
on\/ $\fX$, the locally cotorsion contraherent cosheaf of
contramodules\/ $\fHom_\fX(\cE,\cJ)$ on\/ $\fX$ is coflasque; \par
\textup{(b)} for any projective locally cotorsion contraherent
cosheaf of contramodules\/ $\fF$ on $\fX$, the quasi-coherent torsion
sheaf\/ $\cE\ocn_\fX\fF$ on\/ $\fX$ is flasque; \par
\textup{(c)} if the Noetherian formal scheme\/ $\fX$ has finite Krull
dimension, then for any flat contraherent cosheaf of contramodules\/
$\fF$ on\/ $\fX$, the quasi-coherent torsion sheaf\/ $\cE\ocn_\fX\fF$
on\/ $\fX$ is flasque.
\end{lem}

\begin{proof}
 The proof of part~(a) follows the argument in the proof of
Lemma~\ref{fHom-and-contratensor-with-alf-are-co-dilute}(a).
 One only needs to say that the locally cotorsion contraherent cosheaf
of contramodules $\fHom_\fU(\fj^*\cE,\cK)$ on $\fU$ is coflasque by
Lemma~\ref{fHom-contratensor-co-flasque}(a), and the direct image
functor~$\fj_!$ takes coflasque (locally cotorsion) contraherent
cosheaves of contramodules on $\fU$ to coflasque (locally cotorsion)
contraherent cosheaves of contramodules on $\fX$ by
Corollary~\ref{co-flasque-direct-image-cor}(b) or~(c).

 The proof of part~(b) is dual-analogous.
 It is clear from
Theorem~\ref{loc-Noetherian-projective-lct-cosh-contramods} that
$\fF$ is a finite direct sum of the direct images of projective
locally cotorsion contraherent cosheaves of contramodules
from $\fU_\alpha$, for any finite affine open covering
$\fX=\bigcup_\alpha\fU_\alpha$ of the formal scheme~$\fX$.
 Hence it suffices to consider the case of $\fF=\fj_!\fG$, where
$\fj\:\fU\rarrow\fX$ is the immersion of an affine open formal
subscheme and $\fG\in\fU\Ctrh^\lct_\prj$.
 Then we have $\cE\ocn_\fX\fj_!\fH\simeq\fj_*(\fj^*\cE\ocn_\fU\fH)$
by Proposition~\ref{contratensor-projection-formula-prop}(c).
 The quasi-coherent torsion sheaf $\fj^*\cE\ocn_\fU\fH$ on $\fU$
is flasque by Lemma~\ref{fHom-contratensor-co-flasque}(b), and
the direct image functor~$\fj_*$ takes flasque quasi-coherent torsion
sheaves on $\fU$ to flasque quasi-coherent torsion sheaves on $\fX$
by Corollary~\ref{co-flasque-direct-image-cor}(a).

 Part~(c): by Corollary~\ref{flat-contraherent-cosheaves-cor}(b),
\,$\fF$ is a direct summand of a finitely iterated extension of
the direct images of flat contraherent cosheaves of contramodules
from~$\fU_\alpha$.
 By Lemma~\ref{contratensor-exact-on-flats-and-alf}(a), the functor
the functor $\cE\ocn_\fX{-}\,\:\fX\Ctrh^\fl\rarrow\fX\Tors$ is exact;
so it preserves finitely iterated extensions.
 Since the class of flasque quasi-coherent torsion sheaves is closed
under extensions in $\fX\Tors$, it suffices to consider the case of
$\fF=\fj_!\fG$, where $\fG$ is a flat contraherent cosheaf of
contramodules on~$\fU$.
 Then, once again, we have $\cE\ocn_\fX\fj_!\fH\simeq
\fj_*(\fj^*\cE\ocn_\fU\fH)$ by
Proposition~\ref{contratensor-projection-formula-prop}(c).
 The argument finishes similarly to the proof of part~(b).
\end{proof}

\begin{prop} \label{flat-contratensor-projection-formula-prop}
 Let\/ $\fj\:\fY\rarrow\fX$ be an open immersion of Noetherian formal
schemes of finite Krull dimension.
 Let\/ $\cE$ be a flasque quasi-coherent torsion sheaf on\/ $\fX$ and
$\fG$ be a flat contraherent cosheaf of contramodules on\/~$\fY$.
 Then the natural
morphism~\eqref{contratensor-projection-formula-map-eqn} from
Proposition~\textup{\ref{contratensor-projection-formula-prop}(a)} is
an isomorphism
\begin{equation} \label{flat-contratensor-projection-formula-eqn}
 \cE\ocn_\fX\fj_!\fG\simeq\fj_*(\fj^*\cE\ocn_\fY\fG)
\end{equation}
of quasi-coherent torsion sheaves on\/~$\fX$.
\end{prop}

\begin{proof}
 This is the formal scheme generalization of~\cite[Lemma~8.4.3]{Pcosh}.
 Similarly to the proof of the previous proposition, we show that
both the left-hand side and the right-hand side
of~\eqref{flat-contratensor-projection-formula-eqn} are exact functors
of a flat contraherent cosheaf of contramodules $\fG\in\fY\Ctrh^\fl$
for any fixed flasque quasi-coherent torsion sheaf $\cE$ on~$\fX$.
 Concerning the left-hand side, we have an exact functor
$\fj_!\:\fY\Ctrh^\fl\rarrow\fX\Ctrh^\fl$ by
Corollary~\ref{flat-contrah-direct-images}(b), and the functor
$\cE\ocn_\fX{-}\,\:\fX\Ctrh^\fl\rarrow\fX\Tors$ is exact by
Lemma~\ref{contratensor-exact-on-flats-and-alf}(a).

 Concerning the right-hand side, once again the functor
$\fj^*\cE\ocn_\fY{-}\,\:\fY\Ctrh^\fl\rarrow\fY\Tors$ is exact by
Lemma~\ref{contratensor-exact-on-flats-and-alf}(a).
 Furthermore, for any flat contraherent cosheaf of contramodules
$\fG$ on $\fY$, the quasi-coherent torsion sheaf $\fj^*\cE\ocn_\fY\fG$
on $\fY$ is flasque by
Lemma~\ref{fHom-from-flasque-contratensor-with-alf-co-flasque}(c).
 Finally, we have an exact functor $\fj_*\:\fY\Tors^\fq\rarrow
\fX\Tors^\fq$, by Corollary~\ref{co-flasque-direct-image-cor}(a).

 Now that we know that both the sides of the desired
isomorphism~\eqref{flat-contratensor-projection-formula-eqn} are
exact functors of the argument $\fG$, we can use
Corollary~\ref{flat-contraherent-cosheaves-cor}(b) to the effect that
it suffices to prove the desired isomorphism in the case of
$\fG=\fk_!\fH$, where $\fH$ is a flat contraherent cosheaf of
contramodules on~$\fV$.
 In that case, the assertion follows from the commutativity of
the triangular diagram in
Lemma~\ref{projection-map-for-contratensor-trian-diag-lemma}.
 The point is that both the morphisms of formal schemes~$\fk$ and
$\ff\circ\fk$ are open immersion of an affine Noetherian formal
scheme~$\fV$.
 So both the morphisms of quasi-coherent torsion sheaves
$\fj^*\cE\ocn_\fY\fk_!\fH\rarrow\fk_*(\fk^*\fj^*\cE\ocn_\fV\fH)$ and
$\cE\ocn_\fX\fj_!\fk_!\fH\rarrow\fj_*\fk_*(\fj^*\fk^*\cE\ocn_\fV\fH)$
are isomorphisms by
Proposition~\ref{contratensor-projection-formula-prop}(c).
 Hence the morphism $\cE\ocn_\fX\fj_!\fk_!\fH\rarrow
\fj_*(\fj^*\cE\ocn_\fY\fk_!\fH)$ is an isomorphism, too.
\end{proof}

\subsection{Background equivalences of derived categories}
\label{background-equivalences-subsecn}
 The definitions and discussion of the exotic derived categories
under consideration in this section can be found in
Sections~\ref{appx-exotic-positselski-subsecn}\+-%
\ref{appx-homotopy-of-inj-proj-subsecn}.
 The notation with equality signs between derived category symbols
in the superindices of the letter~$\sD$ (as in the next theorem) is
explained in Section~\ref{appx-homotopy-of-inj-proj-subsecn}.

 We start with quasi-coherent torsion sheaves before passing to
locally contraherent cosheaves of contramodules.

\begin{thm} \label{qcoh-tors-coderived}
 Let\/ $\fX$ be a locally Noetherian formal scheme.
 Then the classes of Positselski-coacyclic and Becker-coacyclic
complexes in the abelian category of quasi-coherent torsion sheaves\/
$\fX\Tors$ coincide,
$$
 \Ac^\co(\fX\Tors)=\Ac^\bco(\fX\Tors).
$$
 The inclusion of additive/abelian categories\/ $\fX\Tors^\inj\rarrow
\fX\Tors$ induces an equivalence of triangulated categories
$$
 \Hot(\fX\Tors^\inj)\simeq\sD^{\co=\bco}(\fX\Tors).
$$
\end{thm}

\begin{proof}
 This is the formal scheme version of~\cite[Theorem~6.4.10(a)]{Pcosh}.
 The abelian category $\fX\Tors$ has enough injective objects by
Theorem~\ref{loc-Noetherian-injective-qcoh-torsion}(a), and
the class of injective objects is closed under infinite direct
sums in $\fX\Tors$ by
Corollary~\ref{qcoh-tors-direct-sums-of-injectives-cor}.
 So Theorem~\ref{countable-co-products-co-contra-derived-theorem}(a)
below is applicable.
\end{proof}

\begin{lem} \label{qcoh-tors-flasque-coresol-dim}
 Let\/ $\fX$ be a locally Noetherian formal scheme. 
 In this setting: \par
\textup{(a)} the full subcategory of flasque quasi-coherent torsion
sheaves\/ $\fX\Tors^\fq$ is coresolving and closed under infinite
direct sums in the abelian category\/ $\fX\Tors$; \par
\textup{(b)} if the locally Noetherian formal scheme\/ $\fX$ has finite
Krull dimension $D$, then the coresolution dimension of any object of\/
$\fX\Tors$ with respect to the coresolving subcategory\/
$\fX\Tors^\fq\subset\fX\Tors$ does not exceed~$D$.
\end{lem}

\begin{proof}
 This is the formal scheme version of~\cite[Lemma~6.3.1(a)]{Pcosh}.
 Part~(a): the full subcategory $\fX\Tors^\fq$ is closed under
extensions and cokernels of monomorphisms in $\fX\Tors$, since flasque
sheaves of abelian groups have such closedness properties (as one can
see from~\cite[Theorem and Corollary~II.3.1.2]{God}).
 For the assertion that $\fX\Tors^\fq$ is closed under infinite
direct sums in $\fX\Tors$, see Section~\ref{coflasque-subsecn}.
 The full subcategory $\fX\Tors^\fq$ is cogenerating in $\fX\Tors$
by Theorem~\ref{loc-Noetherian-injective-qcoh-torsion}(a) and
Corollary~\ref{injective-torsion-local-flasque}(b).

 Part~(b): the property of a sheaf of abelian groups to be flasque
is local by~\cite[Section~II.3.1]{God}, so in view of
Lemma~\ref{co-resol-dim-independent-of-co-resolution-lem}(b)
it suffices to consider the case of a Noetherian (i.~e.,
quasi-compact; or even Noetherian affine) formal scheme~$\fX$.
 Then Lemma~\ref{finite-Krull-dim-co-flasque-co-resol-dim}(a) is
applicable.
\end{proof}

\begin{lem} \label{qcoh-tors-dilute-coresol-dim}
 Let\/ $\fX$ be a semi-separated Noetherian formal scheme.
 In this setting: \par
\textup{(a)} the full subcategory of dilute quasi-coherent torsion
sheaves\/ $\fX\Tors^\dil$ is coresolving and closed under infinite
direct sums in the abelian category\/ $\fX\Tors$; \par
\textup{(b)} if the formal scheme\/ $\fX$ is covered by $N$ affine open
formal subschemes, then the coresolution dimension of any object of\/
$\fX\Tors$ with respect to the coresolving subcategory\/ $\fX\Tors^\dil$
does not exceed $N-1$.
\end{lem}

\begin{proof}
 This is the formal scheme version of~\cite[Lemma~4.7.1(a)]{Pcosh}.
 Part~(a): the second assertion is
Lemma~\ref{direct-sums-of-dilute}.
 The first assertion was mentioned in the discussion after that
lemma in Section~\ref{dilute-subsecn}.
 Part~(b): the \v Cech
coresolution~\eqref{Cech-sequence-of-torsion-sheaves}
is a coresolution of any quasi-coherent torsion sheaf on $\fX$
by dilute quasi-coherent torsion sheaves, as per the proof of
Lemma~\ref{affine-direct-images-are-dilute}.
\end{proof}

\begin{cor} \label{qcoh-tors-flasque-dil-derived-equivalences}
\textup{(a)} Let\/ $\fX$ be a locally Noetherian formal scheme
and\/ $\st=+$, $\co$, or\/~$\bco$ be a derived category symbol.
 Then the inclusions of additive/exact/abelian categories\/
$\fX\Tors^\inj\rarrow\fX\Tors^\fq\rarrow\fX\Tors$ induce triangulated
equivalences\/ $\Hot^\st(\fX\Tors^\inj)\simeq\sD^\st(\fX\Tors^\fq)
\simeq\sD^\st(\fX\Tors)$. \par
\textup{(b)} Let\/ $\fX$ be a locally Noetherian formal scheme of
finite Krull dimension and\/ $\st=\bb$, $+$, $-$, $\varnothing$,
$\abs+$, $\abs-$, $\abs$, $\co$, or\/~$\bco$ be a derived category
symbol.
 Then the inclusion of exact/abelian categories\/
$\fX\Tors^\fq\rarrow\fX\Tors$ induces a triangulated equivalence\/
$\sD^\st(\fX\Tors^\fq)\simeq\sD^\st(\fX\Tors)$. \par
\textup{(c)} Let\/ $\fX$ be a semi-separated Noetherian formal scheme
and\/ $\st=\bb$, $+$, $-$, $\varnothing$, $\abs+$, $\abs-$, $\abs$,
$\co$, or\/~$\bco$ be a derived category symbol.
 Then the inclusion of exact/abelian categories\/
$\fX\Tors^\dil\rarrow\fX\Tors$ induces a triangulated equivalence\/
$\sD^\st(\fX\Tors^\dil)\simeq\sD^\st(\fX\Tors)$.
\end{cor}

\begin{proof}
 Part~(a) follows from Theorem~\ref{qcoh-tors-coderived} and
Lemma~\ref{qcoh-tors-flasque-coresol-dim}(a) in view of
Propositions~\ref{infinite-co-resolutions-minus-plus-derived}(b),
\ref{infinite-co-resolutions-contra-co-derived}(b),
and~\ref{infinite-co-resolutions-becker-derived}(b).
 Part~(b): this is the formal scheme version
of~\cite[Corollary~6.3.2(a)]{Pcosh}.
 The assertion follows from part~(a) and
Lemma~\ref{qcoh-tors-flasque-coresol-dim}(b) in view of
Proposition~\ref{finite-co-resolutions}(b); see also
Proposition~\ref{finite-co-resolutions-becker-derived-I}(b).
 Part~(c): this is the formal scheme version
of~\cite[Corollary~4.7.3(a)]{Pcosh}.
 The assertion follows from
Lemma~\ref{qcoh-tors-dilute-coresol-dim}(a\+-b) in view of
Propositions~\ref{finite-co-resolutions}(b)
and~\ref{finite-co-resolutions-becker-derived-I}(b).
\end{proof}

\begin{thm} \label{lcth-lct-W-contramod-contraderived}
 Let\/ $\fX$ be a locally Noetherian formal scheme with an open
covering\/~$\bW$.
 Then the classes of Positselski-contraacyclic and Becker-contraacyclic
complexes in the exact category of locally cotorsion\/ $\bW$\+locally
contraherent cosheaves of contramodules\/ $\fX\Lcth^\lct_\bW$ coincide,
$$
 \Ac^\ctr(\fX\Lcth^\lct_\bW)=\Ac^\bctr(\fX\Lcth^\lct_\bW).
$$
 The inclusions of additive/exact categories\/ $\fX\Ctrh^\lct_\prj
\rarrow\fX\Ctrh^\lct\rarrow\fX\Lcth^\lct_\bW$ induce equivalences of
triangulated categories
$$
 \Hot(\fX\Ctrh^\lct_\prj)\simeq\sD^{\ctr=\bctr}(\fX\Ctrh^\lct)\simeq
 \sD^{\ctr=\bctr}(\fX\Lcth^\lct_\bW).
$$
\end{thm}

\begin{proof}
 This is the formal scheme version of~\cite[Theorem~6.4.10(b) and
Corollary~6.4.11]{Pcosh}.
 The exact category $\fX\Lcth^\lct_\bW$ has enough projective objects
by Theorem~\ref{loc-Noetherian-projective-lct-cosh-contramods}(a), and
the class of projective objects is closed under infinite direct
products in $\fX\Lcth^\lct_\bW$ by
Corollary~\ref{proj-lct-closed-under-products}.
 So Theorem~\ref{countable-co-products-co-contra-derived-theorem}(b)
is applicable.
 This proves the triangulated equivalences $\Hot(\fX\Ctrh^\lct_\prj)
\simeq\sD^{\ctr=\bctr}(\fX\Lcth^\lct_\bW)$ and, in particular
(for $\bW=\{\fX\}$), \,$\Hot(\fX\Ctrh^\lct_\prj)\simeq
\sD^{\ctr=\bctr}(\fX\Ctrh^\lct)$.
 Then the triangulated equivalence $\sD^{\ctr=\bctr}(\fX\Ctrh^\lct)
\simeq\sD^{\ctr=\bctr}(\fX\Lcth^\lct_\bW)$ follows.
\end{proof}

\begin{lem} \label{lcth-lct-contramod-ssep-contraherent-resol-dim}
 Let\/ $\fX$ be a locally Noetherian formal scheme with an open
covering\/~$\bW$.
 In this setting: \par
\textup{(a)} the full subcategory of locally cotorsion contraherent
cosheaves of contramodules\/ $\fX\Ctrh^\lct$ is resolving (and closed
under infinite products) in the exact category of locally cotorsion\/
$\bW$\+locally contraherent cosheaves of contramodules\/
$\fX\Lcth^\lct_\bW$; \par
\textup{(b)} if the formal scheme\/ $\fX$ is Noetherian (i.~e.,
quasi-compact) and semi-separated, then the full subcategory of
codilute locally cotorsion contraherent cosheaves of contramodules\/
$\fX\Ctrh^\lct_\cdl$ is resolving and closed under infinite products
in the exact category\/ $\fX\Lcth^\lct_\bW$; \par
\textup{(c)} if the formal scheme\/ $\fX$ is Noetherian and
semi-separated, and covered by $N$ affine open formal subschemes
subordinate to\/ $\bW$, then the resolution dimension of any object
of\/ $\fX\Lcth^\lct_\bW$ with respect to the resolving subcategory\/
$\fX\Ctrh^\lct_\cdl$ (hence also with respect to the resolving
subcategory\/ $\fX\Ctrh^\lct$) does not exceed $N-1$.
\end{lem}

\begin{proof}
 This is a formal scheme version of~\cite[Lemma~4.7.2(a)]{Pcosh}.

 Part~(a): these assertions were already mentioned in the proof of
Corollary~\ref{minus-ctr-bctr-derived-categories-agree}.
 To spell out in some more detail: The full subcategory $\fX\Ctrh^\lct$
is closed under extensions, kernels of admissible epimorphisms, and
infinite products in $\fX\Lcth^\lct_\bW$ according to the discussion in
Sections~\ref{locally-contraherent-cosheaves-of-contramods-subsecn}\+-%
\ref{locally-contraherent-lct-cosheaves-subsecn}.
 The full subcategory $\fX\Ctrh^\lct$ is generating in
$\fX\Lcth^\lct_\bW$ by
Theorem~\ref{loc-Noetherian-projective-lct-cosh-contramods}(a).
 (In the case when $\fX$ is quasi-compact and semi-separated, it
suffices to refer to
Corollary~\ref{projective-contrah-lct-on-qcomp-ssep-corollary}.)

 Part~(b): the full subcategory $\fX\Ctrh^\lct_\cdl$ is closed under
extensions and kernels of admissible epimorphisms in $\fX\Lcth^\lct_\bW$
according to Section~\ref{codilute-subsecn}.
 The full subcategory $\fX\Ctrh^\lct_\cdl$ is closed under infinite
products in $\fX\Lcth^\lct_\bW$ by Lemma~\ref{products-of-codilute}.
 The full subcategory $\fX\Ctrh^\lct_\cdl$ is generating in
$\fX\Lcth^\lct_\bW$ by Lemma~\ref{antilocally-flat-are-codilute}.

 Part~(c): for semi-separated Noetherian formal schemes, the claimed
resolution dimension bound holds due to the existence of the \v Cech
resolution~\eqref{Cech-sequence-of-cosheaves-of-contramods}, which
is a resolution of a locally cotorsion $\bW$\+locally contraherent
cosheaf of contramodules by codilute locally cotorsion contraherent
cosehaves of contramodules (as per the proof of
Lemma~\ref{affine-direct-images-are-codilute}).
\end{proof}

\begin{lem} \label{lcth-lct-contramod-coflasque-resol-dim}
 Let\/ $\fX$ be a locally Noetherian formal scheme with an open
covering\/~$\bW$.
 In this setting: \par
\textup{(a)} the full subcategory of coflasque locally cotorsion
contraherent cosheaves of contramodules\/ $\fX\Ctrh^\lct_\cfq$ is
resolving and closed under infinite products in the exact category
of locally cotorsion\/ $\bW$\+locally contraherent cosheaves of
contramodules\/ $\fX\Lcth^\lct_\bW$; \par
\textup{(b)} if the locally Noetherian formal scheme\/ $\fX$ has finite
Krull dimension $D$, then the resolution dimension of any object of\/
$\fX\Lcth^\lct_\bW$ with respect to the resolving subcategory\/
$\fX\Ctrh^\lct_\cfq$ (hence also with respect to the resolving
subcategory $\fX\Ctrh^\lct$) does not exceed~$D$.
\end{lem}

\begin{proof}
 This is the formal scheme version of~\cite[Lemma~6.3.1(b)]{Pcosh}.

 Part~(a): these assertions were also mentioned in the proof of
Corollary~\ref{minus-ctr-bctr-derived-categories-agree}.
 To spell out in some more detail: The full subcategory
$\fX\Ctrh^\lct_\cfq$ is closed under extensions and cokernels of
monomorphisms in $\fX\Lcth^\lct_\bW$ by
Lemmas~\ref{coflasque-lcth-are-globally-contraherent}
and~\ref{coflasque-adjusted-to-cosections}(a\+-b).
 For the assertion that $\fX\Ctrh^\lct_\cfq$ is closed under infinite
products in $\fX\Lcth^\lct_\bW$, see Sections~\ref{coflasque-subsecn}
and~\ref{locally-contraherent-lct-cosheaves-subsecn}.
 The full subcategory $\fX\Ctrh^\lct_\cfq$ is generating in
$\fX\Lcth^\lct_\bW$ by
Theorem~\ref{loc-Noetherian-projective-lct-cosh-contramods}(a) and
Corollary~\ref{projective-lct-of-contramods-local-coflasque}(b).

 Part~(b): the property of a cosheaf of abelian groups to be coflasque
is local by~\cite[Lemma~3.4.1(a)]{Pcosh}, so in view of
Lemma~\ref{co-resol-dim-independent-of-co-resolution-lem}(a)
it suffices to consider the case of a Noetherian (i.~e.,
quasi-compact; or even Noetherian affine) formal scheme~$\fX$.
 Then Lemma~\ref{finite-Krull-dim-co-flasque-co-resol-dim}(b) is
applicable.
\end{proof}

\begin{cor} \label{lcth-lct-contramod-coflasque-derived-equivalences}
\textup{(a)} Let\/ $\fX$ be a locally Noetherian formal scheme with
an open covering\/ $\bW$ and\/ $\st=-$, $\ctr$, or\/~$\bctr$ be
a derived category symbol.
 Then the inclusions of additive/exact categories\/ $\fX\Ctrh^\lct_\prj
\rarrow\fX\Ctrh^\lct_\cfq\rarrow\fX\Ctrh^\lct\rarrow\fX\Lcth^\lct_\bW$
induce triangulated equivalences\/ $\Hot^\st(\fX\Ctrh^\lct_\prj)\simeq
\sD^\st(\fX\Ctrh^\lct_\cfq)\simeq\sD^\st(\fX\Ctrh^\lct)\simeq
\sD^\st(\fX\Lcth^\lct_\bW)$. \par
\textup{(b)} Let\/ $\fX$ be a semi-separated Noetherian formal scheme
with an open covering\/ $\bW$ and\/ $\st=\bb$, $+$, $-$, $\varnothing$, 
$\abs+$, $\abs-$, $\abs$, $\ctr$, or\/~$\bctr$ be a derived category
symbol.
 Then the inclusions of exact categories\/ $\fX\Ctrh^\lct_\cdl\rarrow
\fX\Ctrh^\lct\rarrow\fX\Lcth^\lct_\bW$ induce triangulated
equivalences\/ $\sD^\st(\fX\Ctrh^\lct_\cdl)\simeq\sD^\st(\fX\Ctrh^\lct)
\simeq\sD^\st(\fX\Lcth^\lct_\bW)$. \par
\textup{(c)} Let\/ $\fX$ be a locally Noetherian formal scheme of
finite Krull dimension with an open covering\/ $\bW$ and\/ $\st=\bb$,
$+$, $-$, $\varnothing$, $\abs+$, $\abs-$, $\abs$, $\ctr$, or\/~$\bctr$
be a derived category symbol.
 Then the inclusions of exact categories\/ $\fX\Ctrh^\lct_\cfq\rarrow
\fX\Ctrh^\lct\rarrow\fX\Lcth^\lct_\bW$ induce triangulated
equivalences\/ $\sD^\st(\fX\Ctrh^\lct_\cfq)\simeq
\sD^\st(\fX\Ctrh^\lct)\simeq\sD^\st(\fX\Lcth^\lct_\bW)$.
\hbadness=1900
\end{cor}

\begin{proof}
 Part~(a), which is included here for the sake of completeness of
the formulation, is a part of
Corollary~\ref{minus-ctr-bctr-derived-categories-agree}(b).
 The assertions follow from
Lemmas~\ref{lcth-lct-contramod-ssep-contraherent-resol-dim}(a)
and~\ref{lcth-lct-contramod-coflasque-resol-dim}(a) in view of
Propositions~\ref{infinite-co-resolutions-minus-plus-derived}(a),
\ref{infinite-co-resolutions-minus-plus-derived}(a),
and~\ref{infinite-co-resolutions-becker-derived}(a).

 Part~(b), which a formal scheme version
of~\cite[Corollary~4.7.5(a)]{Pcosh}, follows from
Lemma~\ref{lcth-lct-contramod-ssep-contraherent-resol-dim}(a\+-c)
by virtue of Propositions~\ref{finite-co-resolutions}(a)
and~\ref{finite-co-resolutions-becker-derived-I}(a).
 Part~(c), which is a formal scheme version
of~\cite[Corollaries~6.3.2(b) and~6.3.3(a)]{Pcosh}, follows from
Lemmas~\ref{lcth-lct-contramod-ssep-contraherent-resol-dim}(a)
and~\ref{lcth-lct-contramod-coflasque-resol-dim}(a\+-b)
in view of the same two propositions below.
\end{proof}

\begin{lem} \label{finite-krulldim-ctrh-flat-all-derived-cats-agree}
 Let\/ $\fX$ be a Noetherian formal scheme of finite Krull dimension.
 Then the classes of absolutely acyclic, Positselski-coacyclic,
Positselski-contraacyclic, Becker-coacyclic, Becker-contraacyclic,
and acyclic complexes in the exact category of flat contraherent
cosheaves of contramodules\/ $\fX\Ctrh^\fl$ coincide,
\begin{multline*}
 \Ac^\abs(\fX\Ctrh^\fl)=\Ac^\co(\fX\Ctrh^\fl)=\Ac^\ctr(\fX\Ctrh^\fl)\\
 =\Ac^\bco(\fX\Ctrh^\fl)=\Ac^\bctr(\fX\Ctrh^\fl)=\Ac(\fX\Ctrh^\fl).
\end{multline*}
\end{lem}

\begin{proof}
 This is a restatement of a part of
Corollary~\ref{minus-ctr-bctr-derived-categories-agree}(a).
\end{proof}

\begin{thm} \label{lcth-W-contramod-contraderived}
 Let\/ $\fX$ be a Noetherian formal scheme of finite Krull dimension
with an open covering\/~$\bW$.
 Then the classes of Positselski-contraacyclic and Becker-contraacyclic
complexes in the exact category of\/ $\bW$\+locally contraherent
cosheaves of contramodules\/ $\fX\Lcth_\bW$ coincide,
$$
 \Ac^\ctr(\fX\Lcth_\bW)=\Ac^\bctr(\fX\Lcth_\bW).
$$
 The inclusions of additive/exact categories\/ $\fX\Ctrh_\prj\rarrow
\fX\Ctrh^\fl\rarrow\fX\Ctrh\rarrow\fX\Lcth_\bW$ induce equivalences of
triangulated categories
$$
 \Hot(\fX\Ctrh_\prj)\simeq\sD^{\abs=\varnothing}(\fX\Ctrh^\fl)
 \simeq\sD^{\ctr=\bctr}(\fX\Ctrh)\simeq
 \sD^{\ctr=\bctr}(\fX\Lcth_\bW).
$$
\end{thm}

\begin{proof}
 This is the formal scheme version of~\cite[Theorem~6.4.10(d)]{Pcosh}.
 The exact category $\fX\Lcth_\bW$ has enough projective objects,
and all of them belong to $\fX\Ctrh^\fl$, by
Corollary~\ref{projective-contrah-on-finite-Krull-dim}(a,c).
 The full subcategory $\fX\Ctrh^\fl$ is obviously closed under infinite
products in $\fX\Lcth_\bW$, and all the objects of $\fX\Ctrh^\fl$ have
finite projective dimensions in $\fX\Lcth_\bW$ by
Corollary~\ref{flat-proj-lct-preenvelope-cor}(b).
 So Theorem~\ref{countable-co-products-co-contra-derived-theorem}(b)
is applicable.
 This proves the triangulated equivalences $\Hot(\fX\Ctrh_\prj)
\simeq\sD^{\ctr=\bctr}(\fX\Lcth_\bW)$ and, in particular (for
$\bW=\{\fX\}$), \,$\Hot(\fX\Ctrh_\prj)\simeq\sD^{\ctr=\bctr}(\fX\Ctrh)$.
 Then the triangulated equivalence $\sD^{\ctr=\bctr}(\fX\Ctrh)\simeq
\sD^{\ctr=\bctr}(\fX\Lcth_\bW)$ follows.
 The triangulated equivalence $\Hot(\fX\Ctrh_\prj)\simeq
\sD^{\abs=\ctr=\bctr=\varnothing}(\fX\Ctrh^\fl)$ holds by
Corollary~\ref{minus-ctr-bctr-derived-categories-agree}(a).

 One can also observe that the full subcategory $\fX\Ctrh^\fl$ is
resolving in $\fX\Lcth_\bW$, essentially by
Corollary~\ref{projective-contrah-on-finite-Krull-dim}(a,c).
 Hence the triangulated equivalence
$\sD^\ctr(\fX\Ctrh^\fl)\simeq\sD^\ctr(\fX\Lcth_\bW)$ holds by
Proposition~\ref{infinite-co-resolutions-contra-co-derived}(a),
and the triangulated equivalence
$\sD^\bctr(\fX\Ctrh^\fl)\simeq\sD^\bctr(\fX\Lcth_\bW)$ holds by
Proposition~\ref{infinite-co-resolutions-becker-derived}(a).
\end{proof}

\begin{lem} \label{lcth-contramod-ssep-contraherent-resol-dim}
 Let\/ $\fX$ be a semi-separated Noetherian formal scheme with an open
covering\/~$\bW$.
 In this setting: \par
\textup{(a)} the full subcategory of contraherent cosheaves of
contramodules\/ $\fX\Ctrh$ is resolving (and closed under infinite
products) in the exact category of\/ $\bW$\+locally contraherent
cosheaves of contramodules\/ $\fX\Lcth_\bW$; \par
\textup{(b)} the full subcategory of codilute contraherent cosheaves of
contramodules\/ $\fX\Ctrh_\cdl$ is resolving and closed under infinite
products in the exact category of\/ $\bW$\+locally contraherent
cosheaves of contramodules\/ $\fX\Lcth_\bW$; {\hbadness=1225\par}
\textup{(c)} if the formal scheme\/ $\fX$ is covered by $N$ affine open
formal subschemes subordinate to\/ $\bW$, then the resolution dimension
of any object of\/ $\fX\Lcth_\bW$ with respect to the resolving
subcategory\/ $\fX\Ctrh_\cdl$ (hence also with respect to the resolving
subcategory\/ $\fX\Ctrh$) does not exceed $N-1$.
\end{lem}

\begin{proof}
 This is a formal scheme version of~\cite[Lemma~4.7.1(c)]{Pcosh}.
 Part~(a): The full subcategory $\fX\Ctrh$ is closed under extensions,
kernels of admissible epimorphisms, and infinite products in
$\fX\Lcth_\bW$ according to the discussion in
Section~\ref{locally-contraherent-cosheaves-of-contramods-subsecn}.
 The full subcategory $\fX\Ctrh$ is generating in
$\fX\Lcth_\bW$ by
Corollary~\ref{projective-contrah-on-qcomp-ssep-corollary}.
 The proofs of parts~(b\+-c) are similar to those of
Lemma~\ref{lcth-lct-contramod-ssep-contraherent-resol-dim}(b\+-c).
\end{proof}

\begin{cor} \label{lcth-contramod-ssep-derived-equivalences}
 Let\/ $\fX$ be a semi-separated Noetherian formal scheme
with an open covering\/~$\bW$.
 In this setting: \par
\textup{(a)} Let\/ $\st=-$, $\ctr$, or\/~$\bctr$ be a derived category
symbol.
 Then the inclusions of additive/exact categories\/ $\fX\Ctrh_\alf
\rarrow\fX\Ctrh_\cdl\rarrow\fX\Ctrh\rarrow\fX\Lcth_\bW$ induce
triangulated equivalences\/ $\sD^\st(\fX\Ctrh_\alf)\simeq
\sD^\st(\fX\Ctrh_\cdl)\simeq\sD^\st(\fX\Ctrh)\simeq
\sD^\st(\fX\Lcth_\bW)$. \par
\textup{(b)} Let\/ $\st=\bb$, $+$, $-$, $\varnothing$, $\abs+$,
$\abs-$, $\abs$, $\ctr$, or\/~$\bctr$ be a derived category symbol.
 Then the inclusions of exact categories\/ $\fX\Ctrh_\cdl\rarrow
\fX\Ctrh\rarrow\fX\Lcth_\bW$ induce triangulated equivalences\/
$\sD^\st(\fX\Ctrh_\cdl)\simeq\sD^\st(\fX\Ctrh)\simeq
\sD^\st(\fX\Lcth_\bW)$.
\end{cor}

\begin{proof}
 Part~(a): the full subcategory $\fX\Ctrh_\alf$ is resolving
in $\fX\Lcth_\bW$
by Corollaries~\ref{antilocally-flat-characterized-cor}(b)
and~\ref{antilocally-flat-corollary}(b).
 The full subcategory $\fX\Ctrh_\alf$ is also closed under infinite
products in $\fX\Lcth_\bW$ by Corollary~\ref{antilocally-flat-products}.
 The full subcategories $\fX\Ctrh_\cdl$ and $\fX\Ctrh$ in
$\fX\Lcth_\bW$ have similar properties by
Lemma~\ref{lcth-contramod-ssep-contraherent-resol-dim}(a\+-b).
 In view of
Propositions~\ref{infinite-co-resolutions-minus-plus-derived}(a),
\ref{infinite-co-resolutions-minus-plus-derived}(a),
and~\ref{infinite-co-resolutions-becker-derived}(a),
the desired assertions follow.

 Part~(b), which is a formal scheme version
of~\cite[Corollary~4.7.4(a)]{Pcosh}, follows from
Lemma~\ref{lcth-contramod-ssep-contraherent-resol-dim}(a\+-c)
by virtue of Propositions~\ref{finite-co-resolutions}(a)
and~\ref{finite-co-resolutions-becker-derived-I}(a).
\end{proof}

\begin{lem} \label{lcth-contramod-coflasque-resol-dim}
 Let\/ $\fX$ be a Noetherian formal scheme of finite Krull dimension
$D$ with an open covering\/~$\bW$.
 In this setting: \par
\textup{(a)} the full subcategory of contraherent cosheaves of
contramodules\/ $\fX\Ctrh$ is resolving and closed under infinite
products in the exact category of\/ $\bW$\+locally contraherent
cosheaves of contramodules\/ $\fX\Lcth_\bW$; \par
\textup{(b)} the full subcategory of coflasque contraherent cosheaves of
contramodules\/ $\fX\Ctrh_\cfq$ is resolving and closed under infinite
products in the exact category of\/ $\bW$\+locally contraherent
cosheaves of contramodules\/ $\fX\Lcth_\bW$; \par
\textup{(c)} the resolution dimension of any object of\/ $\fX\Lcth_\bW$
with respect to the resolving subcategory\/ $\fX\Ctrh_\cfq$ (hence also
with respect to the resolving subcategory $\fX\Ctrh$) does not
exceed~$D$.
\end{lem}

\begin{proof}
 This is the formal scheme version of~\cite[Lemma~6.3.1(c)]{Pcosh}.

 Parts~(a\+-b): these assertions were also mentioned in the proof of
Corollary~\ref{minus-ctr-bctr-derived-categories-agree}.
 To spell out in some more detail: The full subcategory $\fX\Ctrh$ is
closed under extensions, kernels of admissible epimorphisms, and
infinite products in $\fX\Lcth_\bW$ according to the discussion in
Section~\ref{locally-contraherent-cosheaves-of-contramods-subsecn}.
 The full subcategory $\fX\Ctrh_\cfq$ is closed under extensions and
cokernels of monomorphisms in $\fX\Lcth_\bW$ by
Lemmas~\ref{coflasque-lcth-are-globally-contraherent}
and~\ref{coflasque-adjusted-to-cosections}(a\+-b).
 For the assertion that $\fX\Ctrh_\cfq$ is closed under infinite
products in $\fX\Lcth_\bW$, see Section~\ref{coflasque-subsecn}.
 The full subcategory $\fX\Ctrh_\cfq$ is generating in $\fX\Lcth_\bW$
by Corollary~\ref{projective-contrah-on-finite-Krull-dim}(a,c); it
follows that the full subcategory $\fX\Ctrh\subset\fX\Lcth_\bW$ is
generating, too.

 The proof of part~(c) is similar to that of
Lemma~\ref{lcth-lct-contramod-coflasque-resol-dim}(b).
\end{proof}

\begin{cor} \label{lcth-contramod-finite-krulldim-derived-equivalences}
 Let\/ $\fX$ be a Noetherian formal scheme of finite Krull dimension
with an open covering\/~$\bW$.
 In this setting: \par
\textup{(a)} Let\/ $\st=-$, $\ctr$, or\/~$\bctr$ be a derived category
symbol.
 Then the inclusions of additive/exact categories\/ $\fX\Ctrh_\prj
\rarrow\fX\Ctrh^\fl\rarrow\fX\Ctrh_\cfq\rarrow\fX\Ctrh\rarrow
\fX\Lcth_\bW$ induce triangulated equivalences\/
$\Hot^\st(\fX\Ctrh_\prj)\simeq\sD^\st(\fX\Ctrh^\fl)\simeq
\sD^\st(\fX\Ctrh_\cfq)\simeq\sD^\st(\fX\Ctrh)\simeq
\sD^\st(\fX\Lcth_\bW)$. {\hbadness=1425\par}
\textup{(b)} Let\/ $\st=\bb$, $+$, $-$, $\varnothing$, $\abs+$,
$\abs-$, $\abs$, $\ctr$, or\/~$\bctr$ be a derived category symbol.
  Then the inclusions of exact categories\/ $\fX\Ctrh_\cfq\rarrow
\fX\Ctrh\rarrow\fX\Lcth_\bW$ induce triangulated equivalences\/
$\sD^\st(\fX\Ctrh_\cfq)\simeq\sD^\st(\fX\Ctrh)\simeq
\sD^\st(\fX\Lcth_\bW)$.
\end{cor}

\begin{proof}
 Part~(a), which is included here for the sake of completeness of
the formulation, is a part of
Corollary~\ref{minus-ctr-bctr-derived-categories-agree}(c).
 The full subcategories $\fX\Ctrh^\fl$, $\fX\Ctrh_\cfq$, and
$\fX\Ctrh$ are resolving in $\fX\Lcth_\bW$ by
Lemma~\ref{lcth-contramod-coflasque-resol-dim}(a\+-b)
and the second paragraph of the proof of
Theorem~\ref{lcth-W-contramod-contraderived}.
 In view of
Propositions~\ref{infinite-co-resolutions-minus-plus-derived}(a),
\ref{infinite-co-resolutions-minus-plus-derived}(a),
and~\ref{infinite-co-resolutions-becker-derived}(a), the desired
assertions follow.

 Part~(b), which is a formal scheme version
of~\cite[Corollaries~6.3.2(c) and~6.3.3(b)]{Pcosh}, follows from
Lemma~\ref{lcth-contramod-coflasque-resol-dim}(a\+-c) by virtue of
Propositions~\ref{finite-co-resolutions}(a)
and~\ref{finite-co-resolutions-becker-derived-I}(a).
\end{proof}

\begin{prop} \label{lcth-contramod-lcta-lct-derived-equivalences}
 Let\/ $\fX$ be a Noetherian formal scheme of finite Krull dimension
with an open covering\/ $\bW$, and let\/ $\st=\bb$, $+$, $-$,
$\varnothing$, $\abs+$, $\abs-$, $\abs$, $\ctr$, or\/~$\bctr$ be
a derived category symbol.
 In this context: \par
\textup{(a)} the inclusion of exact categories\/ $\fX\Ctrh^\lct_\cfq
\rarrow\fX\Ctrh_\cfq$ induces a triangulated equivalence\/
$\sD^\st(\fX\Ctrh^\lct_\cfq)\simeq\sD^\st(\fX\Ctrh_\cfq)$; \par
\textup{(b)} the inclusion of exact categories\/ $\fX\Lcth^\lct_\bW
\rarrow\fX\Lcth_\bW$ induces a triangulated equivalence\/
$\sD^\st(\fX\Lcth^\lct_\bW)\simeq\sD^\st(\fX\Lcth_\bW)$.
\end{prop}

\begin{proof}
 This is the formal scheme version
of~\cite[Corollary~6.4.4(a\+-b)]{Pcosh}.
 The improvement over the exposition in~\cite[Section~6.2]{Pcosh}
achieved in Section~\ref{flat-cosheaves-of-contramods-subsecn} above
allows us to prove part~(b) directly (rather than deduce it from
part~(a), as it was done in~\cite{Pcosh}).

 Part~(b): the full subcategory $\fX\Lcth^\lct_\bW$ is closed under
extensions, cokernels of admissible monomorphisms, and infinite
products in $\fX\Lcth_\bW$ according to
Section~\ref{locally-contraherent-lct-cosheaves-subsecn}.
 The full subcategory $\fX\Lcth^\lct_\bW$ is also cogenerating in
$\fX\Lcth_\bW$ by
Proposition~\ref{finite-krulldim-locally-cotorsion-preenvelope}(a);
and all objects of $\fX\Lcth_\bW$ have finite, uniformly bounded
coresolution dimensions with respect to $\fX\Lcth^\lct_\bW$ by
Proposition~\ref{finite-krulldim-locally-cotorsion-preenvelope}(b).
 Applying Proposition~\ref{finite-co-resolutions}(b), we obtain
the desired assertion for all the mentioned derived category symbols
except $\st=\bctr$.

 In the case of $\st=\bctr$, it suffices to say that
$\Ac^\bctr(\fX\Lcth^\lct_\bW)=\Ac^\ctr(\fX\Lcth^\lct_\bW)$ by
Theorem~\ref{lcth-lct-W-contramod-contraderived} and
$\Ac^\bctr(\fX\Lcth_\bW)=\Ac^\ctr(\fX\Lcth_\bW)$ by
Theorem~\ref{lcth-W-contramod-contraderived}.
 This reduces the assertion to the case of $\st=\ctr$, which is
covered by the previous paragraph.
 Alternatively, a direct argument based on
Proposition~\ref{finite-co-resolutions-becker-derived-II}(b)
is also possible.
 Then one needs to recall that the pair of classes
$(\fX\Ctrh^\fl,\allowbreak\>\fX\Lcth^\lct_\bW)$ is a hereditary
complete cotorsion pair in $\fX\Lcth_\bW$ according to the discussion
in Section~\ref{flat-cosheaves-of-contramods-subsecn}, and
there are enough projective objects in $\fX\Lcth_\bW$ by
Corollary~\ref{projective-contrah-on-finite-Krull-dim}(a).

 Part~(a) can be deduced by combining part~(b) with
Corollaries~\ref{lcth-lct-contramod-coflasque-derived-equivalences}(c)
and~\ref{lcth-contramod-finite-krulldim-derived-equivalences}(b).
 Alternatively, a direct argument similar to the proof of part~(a)
is spelled out below.  {\hbadness=1300\par}

 The full subcategory $\fX\Ctrh^\lct_\cfq$ is closed under extensions,
cokernels of admissible monomorphisms, and infinite products in
$\fX\Ctrh_\cfq$ by
Sections~\ref{locally-contraherent-lct-cosheaves-subsecn}
and~\ref{coflasque-subsecn}.
 The full subcategory $\fX\Ctrh^\lct_\cfq$ is also cogenerating in
$\fX\Ctrh_\cfq$ by
Lemma~\ref{coflasque-locally-cotorsion-preenvelope}(a);
and all objects of $\fX\Ctrh_\cfq$ have finite, uniformly bounded
coresolution dimensions with respect to $\fX\Ctrh^\lct_\cfq$ by
Lemma~\ref{coflasque-locally-cotorsion-preenvelope}(b).
 So Proposition~\ref{finite-co-resolutions}(b) covers the desired
assertion for all the mentioned derived category symbols
except $\st=\bctr$.

 In the case of $\st=\bctr$, there are the same two options as in
the proof of part~(b) above.
 One can prove the equalities
$\Ac^\bctr(\fX\Ctrh^\lct_\cfq)=\Ac^\ctr(\fX\Ctrh^\lct_\cfq)$ and
$\Ac^\bctr(\fX\Ctrh_\cfq)=\Ac^\ctr(\fX\Ctrh_\cfq)$ by applying
Theorem~\ref{countable-co-products-co-contra-derived-theorem}(b),
similarly to the proofs of 
Theorems~\ref{lcth-lct-W-contramod-contraderived}
and~\ref{lcth-W-contramod-contraderived}.
 Alternatively, a direct argument based on
Proposition~\ref{finite-co-resolutions-becker-derived-II}(b)
is possible.

 Then one needs to use the fact that the pair of classes
$(\fX\Ctrh^\fl,\allowbreak\>\fX\Ctrh^\lct_\cfq)$ is a hereditary
complete cotorsion pair in $\fX\Ctrh_\cfq$.
 This can be proved based on the results of
Section~\ref{flat-cosheaves-of-contramods-subsecn}, such as
Lemmas~\ref{coflasque-locally-cotorsion-preenvelope}(a)
and~\ref{coflasque-flat-precover}(a) together with
Corollary~\ref{flat-contraherent-cosheaves-cor}(a\+-b), similarly
to the proof of~\cite[Corollary~6.2.11]{Pcosh}.
 Alternatively one can restrict the hereditary complete cotorsion
pair in $\fX\Lcth_\bW$ mentioned above in the proof of part~(b)
to the exact subcategory $\fX\Ctrh_\cfq\subset\fX\Lcth_\bW$ using
Lemma~\ref{coflasque-adjusted-to-cosections}(a\+-b),
Corollary~\ref{W-flat-implies-flat-and-coflasque-cor}(c),
and~\cite[Lemmas~B.1.3 and~B.1.4(a)]{Pcosh}.
 One also needs to say that there are enough projective objects
in $\fX\Ctrh_\cfq$ by
Lemma~\ref{coflasque-adjusted-to-cosections}(b) and
Corollary~\ref{projective-contrah-on-finite-Krull-dim}(a,c).
\end{proof}

\begin{cor}
 Let\/ $\fX$ be a Noetherian formal scheme of finite Krull dimension
with an open covering\/ $\bW$, and let\/ $\st=\bb$, $+$, $-$,
$\varnothing$, $\abs+$, $\abs-$, $\abs$, $\ctr$, or\/~$\bctr$ be
a derived category symbol.
 Then there is a commutative diagram of triangulated equivalences
between the conventional or exotic derived categories induced by
the inclusions of exact categories
\begin{equation} \label{derived-equivs-induced-by-inclusions-diag-I}
\begin{gathered}
 \xymatrix{
  \sD^\st(\fX\Ctrh^\lct_\cfq) \ar@<-2pt>[r] \ar@<-2pt>[d]
  & \sD^\st(\fX\Ctrh^\lct) \ar@<-2pt>[r] \ar@<-2pt>@{-}[l] \ar@<-2pt>[d]
  & \sD^\st(\fX\Lcth^\lct_\bW) \ar@<-2pt>@{-}[l] \ar@<-2pt>[d] \\
  \sD^\st(\fX\Ctrh_\cfq) \ar@<-2pt>[r] \ar@<-2pt>@{-}[u]
  & \sD^\st(\fX\Ctrh) \ar@<-2pt>[r] \ar@<-2pt>@{-}[l] \ar@<-2pt>@{-}[u]
  & \sD^\st(\fX\Lcth_\bW) \ar@<-2pt>@{-}[l] \ar@<-2pt>@{-}[u]
 }
\end{gathered}
\end{equation}
\end{cor}

\begin{proof}
 The upper horizontal equivalences are the result of
Corollary~\ref{lcth-lct-contramod-coflasque-derived-equivalences}(c).
 The lower horizontal equivalences are the result of
Corollary~\ref{lcth-contramod-finite-krulldim-derived-equivalences}(b).
 The vertical equivalences are the result of
Proposition~\ref{lcth-contramod-lcta-lct-derived-equivalences}(a\+-b).
 The commutativity of
the diagram~\eqref{derived-equivs-induced-by-inclusions-diag-I}
follows from the obvious commutativity of the underlying diagram
of inclusions of exact categories.
 The point is that all the mentioned exact (sub)categories are
either resolving (in the case of the horizonal arrows) or coresolving
(in the case of the vertical arrows) in each other, with
the (co)resolution dimension always bounded by the Krull dimension
of~$\fX$, with enough projective objects and exact functors of
infinite direct products, and closed under infinite products
in each other.
 Besides, the infinite products of projective objects are either
projective (in the upper line), or have finite projective dimensions
(in the lower line).
\end{proof}

\subsection{The unbounded derived equivalence}
\label{unbounded-derived-equiv-subsecn}
 The aim of this section is to prove the following theorem
complementing the results of the previous
Section~\ref{background-equivalences-subsecn}.
 It is a formal scheme version of~\cite[Theorem~4.7.9]{Pcosh}.

\begin{thm} \label{lcth-lcta-lct-unbounded-derived-equivalence}
 Let\/ $\fX$ be a semi-separated Noetherian formal scheme with
an open covering\/~$\bW$.
 Then, for any derived category symbol\/ $\st=+$ or\/~$\varnothing$,
the inclusion of exact categories\/ $\fX\Lcth^\lct_\bW\rarrow
\fX\Lcth_\bW$ induces a triangulated equivalence\/
$\sD^\st(\fX\Lcth^\lct_\bW)\simeq\sD^\st(\fX\Lcth_\bW)$.
\end{thm}

 Let us introduce some notation and conventions.
 Given an additive category $\sE$, we denote by $\Com(\sE)$
the additive category of (unbounded) complexes in~$\sE$.
 If $\sE$ is an exact category, then we endow the additive category
$\Com(\sE)$ with the \emph{termwise exact structure}, where
the admissible short exact sequence in $\Com(\sE)$ are the short
sequences of complexes in $\sE$ that are admissible exact as short
sequences in $\sE$ in every degree.
 Whenever $\sE$ is an abelian category with the abelian exact
structure, the category $\Com(\sE)$ is also abelian, and the termwise
exact structure on $\Com(\sE)$ coincides with the abelian exact
structure.

 We start with the following lemma from adic commutative algebra.

\begin{lem} \label{contramodule-complex-cotorsion-preenvelope}
 Let $R$ be an adic topological ring and\/ $\fM^\bu$ be a complex of
quotseparated contramodule $R$\+modules.
 Then the exists a short exact sequence of complexes of quotseparated
contramodule $R$\+modules\/ $0\rarrow\fM^\bu\rarrow\fP^\bu\rarrow\fF^\bu
\rarrow0$ such that\/ $\fP^\bu$ is a complex of cotorsion quotseparated
contramodule $R$\+modules and $\fF^\bu$ is an acyclic complex in
the exact category $R\Contra_\fl$ of flat quotseparated contramodule
$R$\+modules.
\end{lem}

\begin{proof}
 The assertion actually holds for left contramodules over any
complete, separated topological associative ring $\fR$ with
a countable base of neighborhoods of zero consisting of open right
ideals (which is a generalization of the setting of the lemma
by~\cite[Proposition~1.5]{Pdc}).
 The argument is based on the results of the papers~\cite{PR,Pfltp}.
 We use the language of cotorsion pairs from
Section~\ref{appx-cotorsion-pairs-subsecn} and the language of
Becker-coacyclic complexes from Section~\ref{appx-becker-subsecn}.

 Denote, as usual, by $\fR\Contra$ the abelian category of left
contramodules over the topological ring $\fR$, and consider
the abelian category $\Com(\fR\Contra)$ of complexes in $\fR\Contra$.
 We refer to~\cite[Sections~5\+-7]{PR} for the definitions of
the full subcategories of \emph{flat left\/ $\fR$\+contramodules}
$\fR\Contra_\fl$ and \emph{cotorsion left\/ $\fR$\+contramodules}
$\fR\Contra^\cot\subset\fR\Contra$.
 By~\cite[Corollary~7.8 and its proof]{PR} together
with~\cite[Corollaries~6.8 and~6.15]{PR}, the pair of classes
$(\fR\Contra_\fl,\allowbreak\>\fR\Contra^\cot)$ is a hereditary
complete cotorsion pair generated by a set of objects in the abelian
category $\fR\Contra$.
 According to~\cite[Section~5]{PR}, the abelian category
$\fR\Contra$ is locally $\aleph_1$\+presentable (in the sense
of~\cite[Definition~1.17 and Theorem~1.20]{AR}).

 Consider the full subcategory of all complexes of cotorsion left
$\fR$\+contramodules $\Com(\fR\Contra^\cot)\subset\Com(\fR\Contra)$,
and consider the class $\Ac^\bco(\fR\Contra_\fl)\allowbreak\subset
\Com(\fR\Contra_\fl)$ of all Becker-coacyclic complexes in the exact
category of flat $\fR$\+contramodules $\fR\Contra_\fl$, which we will
also view as a full subcategory in the category of all complexes
of $\fR$\+contramodules $\Com(\fR\Contra)$.
 Similarly, let $\Ac(\fR\Contra_\fl)\subset\Com(\fR\Contra_\fl)
\subset\Com(\fR\Contra)$ be the class of all acyclic complexes
in the exact category $\fR\Contra_\fl$.
 By~\cite[Proposition~2.1]{Pfltp} applied to the locally presentable
abelian category $\sK=\fR\Contra$ with the cotorsion pair
$(\sA,\sB)=(\fR\Contra_\fl,\allowbreak\>\fR\Contra^\cot)$ from
the previous paragraph, the pair of classes
$(\Ac^\bco(\fR\Contra_\fl),\allowbreak\>\Com(\fR\Contra^\cot))$
is a hereditary complete cotorsion pair in the abelian category
of complexes $\Com(\fR\Contra)$.
 By~\cite[Corollary~4.2]{Pfltp}, all Becker-coacyclic complexes
in the exact category $\fR\Contra_\fl$ are acyclic in this exact
category, i.~e., $\Ac^\bco(\fR\Contra_\fl)\subset
\Ac(\fR\Contra_\fl)$.

 Now any special preenvelope short exact
sequence~\eqref{special-preenvelope-sequence} for the object
$\fM^\bu\in\Com(\fR\Contra)$ with respect to the mentioned
cotorsion pair in $\Com(\fR\Contra)$ is a short exact sequence
claimed in the lemma.
\end{proof}

\begin{prop} \label{lcth-complex-lct-complex-preenvelope}
 Let\/ $\fX$ be a semi-separated Noetherian formal scheme with
an open covering\/ $\bW$ and a finite affine open covering\/
$\fX=\bigcup_{\alpha=1}^N\fU_\alpha$ subordinate to\/~$\bW$.
 Then any complex of\/ $\bW$\+locally contraherent cosheaves of
contramodules\/ $\fM^\bu$ on\/ $\fX$ can be included in
an (admissible) short exact sequence\/ $0\rarrow\fM^\bu\rarrow
\fP^\bu\rarrow\fF^\bu\rarrow0$ in the exact category\/
$\Com(\fX\Lcth_\bW)$, where\/ $\fP^\bu\in\Com(\fX\Lcth^\lct_\bW)$
is a complex of locally cotorsion\/ $\bW$\+locally contraherent
cosheaves of contramodules on\/ $\fX$, while\/ $\fF^\bu$ is a finitely
iterated extension of the direct images of complexes of
contraherent cosheaves of contramodules on\/ $\fU_\alpha$
corresponding to acyclic complexes in the exact categories of
flat contraadjusted contramodule\/ $\fO_\fX(\fU_\alpha)$\+modules.
 In particular, $\fF^\bu$ is an acyclic complex in the exact
category $\fX\Ctrh_\alf$ of antilocally flat contraherent cosheaves
of contramodules on\/~$\fX$.
\end{prop}

\begin{proof}
 This is a formal scheme version of~\cite[Lemma~4.7.8]{Pcosh}.
 The proof of the first assertion of the proposition
proceeds along the lines of the proof of
Lemma~\ref{locally-cotorsion-preenvelope-lemma} above.
 Following the fruitful idea of~\cite[proof of Lemma~6.3]{ES},
we plug in whole complexes of $\bW$\+locally contraherent cosheaves of
contramodules, rather than single $\bW$\+locally contraherent cosheaves
of contramodules, into the construction from the proof of
Lemma~\ref{locally-cotorsion-preenvelope-lemma}.
 Otherwise the construction proceeds exactly in the same way,
with the only difference that, in the context of an affine
Noetherian formal scheme $\fU=\fU_\beta$, instead of the short exact
sequence of contramodule $\fO(\fU)$\+modules from
Proposition~\ref{flat-cotorsion-pair-in-quotseparated}(b) we now
use the short exact sequence of complexes of contramodule
$\fO(\fU)$\+modules from
Lemma~\ref{contramodule-complex-cotorsion-preenvelope}.

 To deduce the second assertion of the proposition, it suffices
to observe that, for an affine open formal subscheme $\fU\subset\fX$
with the open immersion morphism $\fj\:\fU\rarrow\fX$, we have
an exact functor of direct image $\fj_!\:\fU\Ctrh_\alf\rarrow
\fX\Ctrh_\alf$ (see Corollary~\ref{antilocally-flat-corollary}(c)).
 The full subcategory $\fU\Ctrh_\alf\subset\fU\Ctrh$ consists
precisely of all the contraherent cosheaves of contramodules
on $\fU$ corresponding to flat contraadjusted contramodule
$\fO_\fX(\fU)$\+modules under the equivalence of categories from
Lemma~\ref{contrah-cosheaves-of-contramods-lemma}(b).
 Finally, any finitely iterated extension of acyclic complexes in
an exact category is also an acyclic complex.
\end{proof}

\begin{proof}[Proof of
Theorem~\ref{lcth-lcta-lct-unbounded-derived-equivalence}]
 The easy case of bounded below complexes, $\st=+$, is covered
by Proposition~\ref{infinite-co-resolutions-minus-plus-derived}(b).
 Indeed, the full subcategory $\fX\Lcth^\lct_\bW$ is closed under
extensions and cokernels of admissible monomorphisms in $\fX\Lcth_\bW$
according to Section~\ref{locally-contraherent-lct-cosheaves-subsecn}.
 The full subcategory $\fX\Lcth^\lct_\bW$ is also cogenerating in
$\fX\Lcth_\bW$ by Lemma~\ref{locally-cotorsion-preenvelope-lemma}
or Corollary~\ref{antilocally-flat-corollary}(a); so the full
subcategory $\fX\Lcth^\lct_\bW$ is coresolving in $\fX\Lcth_\bW$,
and Proposition~\ref{infinite-co-resolutions-minus-plus-derived}(b)
is applicable.

 In the interesting case of two-sided unbounded complexes,
$\st=\varnothing$, the argument is based on
Proposition~\ref{lcth-complex-lct-complex-preenvelope} together
with the cotorsion periodicity theorem for contramodules.
 In view of a basic lemma about full triangulated subcategories and
triangulated Verdier quotient categories~\cite[Corollary~10.2.7]{KS},
\cite[Lemma~A.3.3(b)]{Pcosh}, \cite[Lemma~8.1(b)]{Pfltp}, in order
to prove the assertion of the theorem for $\st=\varnothing$ it suffices
to establish two properties:
\begin{enumerate}
\item for any complex $\fM^\bu$ in $\fX\Lcth_\bW$, there exists
a complex $\fP^\bu$ in $\fX\Lcth^\lct_\bW$ together with
a morphism of complexes $\fM^\bu\rarrow\fP^\bu$ whose cone is
an acyclic complex in $\fX\Lcth_\bW$;
\item any complex in $\fX\Lcth^\lct_\bW$ that is acyclic in
$\fX\Lcth_\bW$ is also acyclic in $\fX\Lcth^\lct_\bW$.
\end{enumerate}

 Now, concerning item~(1), the desired morphism of complexes
$\fM^\bu\rarrow\fP^\bu$ is provided by
Proposition~\ref{lcth-complex-lct-complex-preenvelope}.
 It suffices to observe that any acyclic complex in the exact
category $\fX\Ctrh_\alf$ is also acyclic in the exact category
$\fX\Lcth_\bW$.
 Item~(2) follows from
Lemma~\ref{complex-of-cosheaves-exactness-criterion}(a\+-b),
which is based on
Proposition~\ref{contraadjusted-and-cotorsion-periodicity-prop}(b).
\end{proof}

\begin{cor}
 Let\/ $\fX$ be a semi-separated Noetherian formal scheme with an open
covering\/ $\bW$, and let\/ $\st=+$ or\/~$\varnothing$ be a derived
category symbol.
 Then there is a commutative diagram of triangulated equivalences
between the conventional derived categories induced by the inclusions
of exact categories
\begin{equation} \label{derived-equivs-induced-by-inclusions-diag-II}
\begin{gathered}
 \xymatrix{
  \sD^\st(\fX\Ctrh^\lct_\cdl) \ar@<-2pt>[r] \ar@<-2pt>[d]
  & \sD^\st(\fX\Ctrh^\lct) \ar@<-2pt>[r] \ar@<-2pt>@{-}[l] \ar@<-2pt>[d]
  & \sD^\st(\fX\Lcth^\lct_\bW) \ar@<-2pt>@{-}[l] \ar@<-2pt>[d] \\
  \sD^\st(\fX\Ctrh_\cdl) \ar@<-2pt>[r] \ar@<-2pt>@{-}[u]
  & \sD^\st(\fX\Ctrh) \ar@<-2pt>[r] \ar@<-2pt>@{-}[l] \ar@<-2pt>@{-}[u]
  & \sD^\st(\fX\Lcth_\bW) \ar@<-2pt>@{-}[l] \ar@<-2pt>@{-}[u]
 }
\end{gathered}
\end{equation}
\end{cor}

\begin{proof}
 The upper horizontal equivalences are the result of
Corollary~\ref{lcth-lct-contramod-coflasque-derived-equivalences}(b).
 The lower horizontal equivalences are the result of
Corollary~\ref{lcth-contramod-ssep-derived-equivalences}(b).
 The middle and rightmost vertical equivalences are the result of
Theorem~\ref{lcth-lcta-lct-unbounded-derived-equivalence}.
 The commutativity of
the  diagram~\eqref{derived-equivs-induced-by-inclusions-diag-II}
follows from the obvious commutativity of the underlying diagram
of inclusions of exact categories.
 The leftmost vertical triangulated functor is a triangulated
equivalence in view of the commutativity of the diagram.
\end{proof}

\subsection{Derived direct images of torsion sheaves}
\label{derived-direct-images-of-qcoh-tors-subsecn}
 In the rest of Section~\ref{injective-projective-secn}, we denote
by~$\st$ one of the derived category symbols $\st=\bb$, $+$, $-$,
$\varnothing$, $\abs+$, $\abs-$, $\abs$, $\co$, $\ctr$, $\bco$,
or~$\bctr$.
 Given an abelian or exact category $\sE$, the complexes in $\sE$
belonging to the respective thick subcategory of trivial objects
$\Ac^\st(\sE)\subset\Hot^\st(\sE)$ will be simply called
``$\st$\+acyclic'' for brevity.

 Let $\ff\:\fY\rarrow\fX$ be a quasi-compact morphism of locally
Noetherian formal schemes.
 Applying the direct image functor $\ff_*$
\,\eqref{direct-image-qcoh-tors} to complexes in $\fY\Tors$ termwise,
one obtains the induced triangulated functor between the homotopy
categories
\begin{equation} \label{qcoh-tors-homotopy-cats-direct-image-formula}
 \ff_*\:\Hot(\fY\Tors)\lrarrow\Hot(\fX\Tors).
\end{equation}
 The aim of this
Section~\ref{derived-direct-images-of-qcoh-tors-subsecn} is
to construct right derived functors
\begin{equation} \label{qcoh-tors-derived-direct-image-formula}
 \boR\ff_*\:\sD^\st(\fY\Tors)\lrarrow\sD^\st(\fX\Tors)
\end{equation}
for $\st\ne\ctr$, $\bctr$, and discuss assumptions under which various
constructions are applicable.
 There are several approaches.
 
 Let $\fY$ be a locally Noetherian formal scheme of finite
Krull dimension.
 Then for any derived category symbol $\st\ne\ctr$, $\bctr$, we have
a triangulated equivalence $\sD^\st(\fY\Tors^\fq)\simeq
\sD^\st(\fY\Tors)$ provided by
Corollary~\ref{qcoh-tors-flasque-dil-derived-equivalences}(b).
 The functor~$\ff_*$ induces an exact functor $\ff_*\:\fY\Tors^\fq
\rarrow\fX\Tors^\fq$ by Corollary~\ref{co-flasque-direct-image-cor}(a).
 Hence, for any derived category symbol $\st\ne\co$, $\ctr$, $\bco$,
$\bctr$, the functor~$\ff_*$ takes $\st$\+acyclic complexes in
$\fY\Tors^\fq$ to $\st$\+acyclic complexes in $\fX\Tors^\fq$.
 So the derived functor $\boR\ff_*$
\,\eqref{qcoh-tors-derived-direct-image-formula} can be constructed
by applying the functor~$\ff_*$ termwise to complexes of flasque
quasi-coherent torsion sheaves on~$\fY$.

 For the derived category symbol $\st=\co$, one needs to observe
additionally that the full subcategory $\fY\Tors^\fq$ is closed
under infinite direct sums in $\fY\Tors$ according to
Section~\ref{coflasque-subsecn}, and the functor $\ff_*\:\fY\Tors
\rarrow\fX\Tors$ preserves infinite direct sums according to
Section~\ref{direct-images-of-qcoh-tors-subsecn}.
 Hence the functor~$\ff_*$ also takes Positselski-coacyclic complexes
in $\fY\Tors^\fq$ to Positselski-coacyclic complexes in $\fX\Tors^\fq$,
and the derived functor functor $\boR\ff_*$
\,\eqref{qcoh-tors-derived-direct-image-formula} for $\st=\co$ also can
be constructed by applying the functor~$\ff_*$ termwise to complexes of
flasque quasi-coherent torsion sheaves on~$\fY$.
 Finally, to handle the case $\st=\bco$, it suffices to recall
the result of Theorem~\ref{qcoh-tors-coderived}, which tells us
that for the abelian categories $\fY\Tors$ and $\fX\Tors$ the Becker
coderived category coincides with the Positselski one.

 Let $\fX$ and $\fY$ be semi-separated Noetherian formal schemes.
 Then for any derived category symbol $\st\ne\ctr$, $\bctr$, we have
a triangulated equivalence $\sD^\st(\fY\Tors^\dil)\simeq
\sD^\st(\fY\Tors)$ provided by
Corollary~\ref{qcoh-tors-flasque-dil-derived-equivalences}(c).
 The functor~$\ff_*$ induces an exact functor $\ff_*\:\fY\Tors^\dil
\rarrow\fX\Tors^\dil$ by Corollary~\ref{dilute-direct-images}.
 Hence, for any derived category symbol $\st\ne\co$, $\ctr$, $\bco$,
$\bctr$, the functor~$\ff_*$ takes $\st$\+acyclic complexes in
$\fY\Tors^\dil$ to $\st$\+acyclic complexes in $\fX\Tors^\dil$.
 So the derived functor $\boR\ff_*$
\,\eqref{qcoh-tors-derived-direct-image-formula} can be constructed
by applying the functor~$\ff_*$ termwise to complexes of dilute
quasi-coherent torsion sheaves on~$\fY$.

 For the derived category symbol $\st=\co$, one needs to observe
additionally that the full subcategory $\fY\Tors^\dil$ is closed
under infinite direct sums in $\fY\Tors$
by Lemma~\ref{direct-sums-of-dilute}, and the functor $\ff_*\:\fY\Tors
\rarrow\fX\Tors$ preserves infinite direct sums as mentioned above.
 Hence the functor~$\ff_*$ also takes Positselski-coacyclic complexes
in $\fY\Tors^\dil$ to Positselski-coacyclic complexes in
$\fX\Tors^\dil$, and the derived functor functor $\boR\ff_*$
\,\eqref{qcoh-tors-derived-direct-image-formula} for $\st=\co$ also can
be constructed by applying the functor~$\ff_*$ termwise to complexes of
dilute quasi-coherent torsion sheaves on~$\fY$.
 Finally, the case $\st=\bco$ reduces to that of $\st=\co$ by
Theorem~\ref{qcoh-tors-coderived}.

 In the general case of arbitrary locally Noetherian formal schemes
$\fX$ and $\fY$, one can still construct the derived functor
$\boR\ff_*$ \,\eqref{qcoh-tors-derived-direct-image-formula}
for \emph{some} derived category symbols~$\st$.
 For $\st=+$, $\co$, or~$\bco$, we have a triangulated equivalence
$\sD^\st(\fY\Tors^\fq)\simeq\sD^\st(\fY\Tors)$ provided by
Corollary~\ref{qcoh-tors-flasque-dil-derived-equivalences}(a).
 The same arguments as above are applicable, and the functor~$\boR\ff_*$
can be constructed by applying the functor~$\ff_*$ termwise to complexes
of flasque quasi-coherent torsion sheaves on~$\fY$.

 A tricky case is that of $\st=\varnothing$ for arbitrary locally
Noetherian formal schemes $\fX$ and~$\fY$.
 In that case, for all we know, an unbounded complex in $\fY\Tors^\fq$
that is acyclic in $\fY\Tors$ \emph{need not} be acyclic in
$\fY\Tors^\fq$.
 Nevertheless, applying the functor~$\ff_*$ termwise to any complex
in $\fY\Tors^\fq$ that is acyclic in $\fY\Tors$ produces a complex
in $\fX\Tors^\fq$ that is acyclic in $\fX\Tors$, by
Corollary~\ref{co-flasque-direct-image-of-acyclic-complexes-cor}(a).
 For this reason, the derived functor $\boR\ff_*$
\,\eqref{qcoh-tors-derived-direct-image-formula} still can be
constructed by applying the functor~$\ff_*$ termwise to complexes of
flasque quasi-coherent torsion sheaves on~$\fY$.

 Here we are also using the fact that, for every complex of
quasi-coherent torsion sheaves $\cN^\bu$ on a locally Noetherian
formal scheme $\fY$, there exists a quasi-isomorphism
$\cN^\bu\rarrow\cE^\bu$ from $\cN^\bu$ into a complex of flasque
quasi-coherent torsion sheaves $\cE^\bu$ on~$\fY$.
 This follows, e.~g., from Theorem~\ref{qcoh-tors-coderived}
(or from the existence of enough homotopy injective complexes of
injective objects in any Grothendieck abelian category).

\subsection{Derived direct images of cosheaves of contramodules}
\label{derived-direct-images-of-ctrh-contramod-subsecn}
 Let $\ff\:\fY\rarrow\fX$ be a quasi-compact morphism of locally
Noetherian formal schemes.

 When the morphism~$\ff$ is affine, we have well-defined direct
image functors $\ff_!\:\fY\Ctrh\rarrow\fX\Ctrh$
\,\eqref{ctrh-contramod-affine-direct-image} and
$\ff_!\:\fY\Ctrh^\lct\rarrow\fX\Ctrh^\lct$
\,\eqref{ctrh-lct-contramod-affine-direct-image}.
 We can apply both the functors~$\ff_!$ to complexes in the respective
exact categories $\fY\Ctrh$ and $\fY\Ctrh^\lct$ termwise, and obtain
the induced triangulated functors between the homotopy categories
\begin{align}
 \ff_!\:\Hot(\fY\Ctrh)&\lrarrow\Hot(\fX\Ctrh),
 \label{ctrh-homotopy-cats-affine-direct-image} \\
 \ff_!\:\Hot(\fY\Ctrh^\lct)&\lrarrow\Hot(\fX\Ctrh^\lct).
 \label{ctrh-lct-homotopy-cats-affine-direct-image}
\end{align}

 Furthermore, both the functors $\ff_!$
\,(\ref{ctrh-contramod-affine-direct-image}\+-%
\ref{ctrh-lct-contramod-affine-direct-image}) are exact and preserve
infinite products.
 So, for any derived category symbol $\st\ne\co$, $\bco$, $\bctr$,
both the functors $\ff_!$
\,(\ref{ctrh-homotopy-cats-affine-direct-image}\+-%
\ref{ctrh-lct-homotopy-cats-affine-direct-image}) take $\st$\+acyclic
complexes to $\st$\+acyclic complexes, and therefore induce
triangulated functors between the derived categories
\begin{align}
 \ff_!\:\sD^\st(\fY\Ctrh)&\lrarrow\sD^\st(\fX\Ctrh),
 \label{ctrh-derived-cats-affine-underived-direct-image} \\
 \ff_!\:\sD^\st(\fY\Ctrh^\lct)&\lrarrow\sD^\st(\fX\Ctrh^\lct).
 \label{ctrh-lct-derived-cats-affine-underived-direct-image}
\end{align}

 In the case of $\st=\bctr$, for the locally cotorsion contraherent
cosheaf the Becker contraderived category agrees with the Positselski
one, $\Ac^\bctr(\fY\Ctrh^\lct)=\Ac^\ctr(\fY\Ctrh^\lct)$, by
Theorem~\ref{lcth-lct-W-contramod-contraderived}.
 So we have the triangulated
functor~\eqref{ctrh-lct-derived-cats-affine-underived-direct-image}
for $\st=\bctr$ as well.
 The same applies to the triangulated
functor~\eqref{ctrh-derived-cats-affine-underived-direct-image}
for $\st=\bctr$, by virtue of
Theorem~\ref{lcth-W-contramod-contraderived}, assuming that
the formal scheme $\fY$ is Noetherian (i.~e., quasi-compact)
of finite Krull dimension.

 More generally, let $\bW$ be an open covering of $\fX$ and $\bT$
be an open covering of $\fY$ for which the morphism~$\ff$ is
$(\bW,\bT)$\+affine.
 Then the direct image
functors~(\ref{loc-ctrh-contramod-W-T-affine-direct-image}\+-%
\ref{loc-ctrh-lct-contramod-W-T-affine-direct-image}) can be applied
to complexes in $\fY\Lcth_\bT$ and $\fY\Lcth^\lct_\bT$ termwise,
producing the induced triangulated functors between
the homotopy categories
\begin{align}
 \ff_!\:\Hot(\fY\Lcth_\bT)&\lrarrow\Hot(\fX\Lcth_\bW),
 \label{lcth-homotopy-cats-affine-direct-image} \\
 \ff_!\:\Hot(\fY\Lcth^\lct_\bT)&\lrarrow\Hot(\fX\Lcth^\lct_\bW).
 \label{lcth-lct-homotopy-cats-affine-direct-image}
\end{align}
 Under the same assumptions on~$\st$ and $\fY$ as above, and based
on the same arguments, we have the induced triangulated functors between
the derived categories
\begin{align}
 \ff_!\:\sD^\st(\fY\Lcth_\bT)&\lrarrow\sD^\st(\fX\Lcth_\bW),
 \label{lcth-derived-cats-affine-underived-direct-image} \\
 \ff_!\:\sD^\st(\fY\Lcth^\lct_\bT)&\lrarrow\sD^\st(\fX\Lcth^\lct_\bW).
 \label{lcth-lct-derived-cats-affine-underived-direct-image}
\end{align}

 All the discussion above in this
Section~\ref{derived-direct-images-of-ctrh-contramod-subsecn}
is applicable to \emph{affine} morphisms~$\ff$ \emph{only}.
 When the morphism~$\ff$ is not affine, we \emph{do not} even have
well-defined functors of direct image~$\ff_!$
\,(\ref{ctrh-contramod-affine-direct-image}\+-%
\ref{ctrh-lct-contramod-affine-direct-image})
or~(\ref{ctrh-homotopy-cats-affine-direct-image}\+-%
\ref{ctrh-lct-homotopy-cats-affine-direct-image}).

 The aim of the rest of this section is to construct left derived
functors
\begin{align}
 \boL\ff_!\:\sD^\st(\fY\Ctrh)&\lrarrow\sD^\st(\fX\Ctrh),
 \label{ctrh-contramod-derived-direct-image-formula} \\
 \boL\ff_!\:\sD^\st(\fY\Ctrh^\lct)&\lrarrow\sD^\st(\fX\Ctrh^\lct)
 \label{ctrh-lct-contramod-derived-direct-image-formula}
\end{align}
for $\st\ne\co$, $\bco$ in the case of a \emph{not necessarily affine}
quasi-compact morphism of locally Noetherian formal schemes
$\ff\:\fY\rarrow\fX$.
 This means that we need to derive direct image functors that are only
partially defined.
 There are several approaches, requiring different sets of assumptions.

 We start with the locally cotorsion contraherent cosheaves of
contramodules before passing to the locally contraadjusted ones.

 Let $\fY$ be a locally Noetherian formal scheme of finite Krull
dimension.
 Then for any derived category symbol $\st\ne\co$, $\bco$, we have
a triangulated equivalence $\sD^\st(\fY\Ctrh^\lct_\cfq)\simeq
\sD^\st(\fY\Ctrh^\lct)$ provided by
Corollary~\ref{lcth-lct-contramod-coflasque-derived-equivalences}(c).
 The functor of direct image of cosheaves of $\fO$\+modules~$\ff_!$
induces an exact functor $\ff_!\:\fY\Ctrh^\lct_\cfq\rarrow
\fX\Ctrh^\lct_\cfq$ by Corollary~\ref{co-flasque-direct-image-cor}(c).
 Hence, for any derived category symbol $\st\ne\co$, $\ctr$, $\bco$,
$\bctr$, the functor~$\ff_!$ takes $\st$\+acyclic complexes in
$\fY\Ctrh^\lct_\cfq$ to $\st$\+acyclic complexes
in $\fX\Ctrh^\lct_\cfq$.
 So the derived functor $\boL\ff_!$
\,\eqref{ctrh-lct-contramod-derived-direct-image-formula} can be
constructed by applying the functor~$\ff_!$ termwise to complexes of
coflasque locally cotorsion contraherent cohseaves of contramodules
on~$\fY$.

 For the derived category symbol $\st=\ctr$, one needs to observe
additionally that the full subcategory $\fY\Ctrh^\lct_\cfq$ is closed
under infinite products in $(\fY,\fO_\fY)\Cosh$ according to
Sections~\ref{locally-contraherent-cosheaves-of-contramods-subsecn}\+-%
\ref{locally-contraherent-lct-cosheaves-subsecn}
and~\ref{coflasque-subsecn}, and the functor $\ff_!\:(\fY,\fO_\fY)\Cosh
\rarrow(\fX,\fO_\fX)\Cosh$ preserves infinite products according to
Section~\ref{direct-images-of-ctrh-contramod-subsecn}.
 Hence the functor~$\ff_!$ also takes Positselski-contraacyclic
complexes in $\fY\Ctrh^\lct_\cfq$ to Positselski-contraacyclic complexes
in $\fX\Ctrh^\lct_\cfq$, and the derived functor functor $\boL\ff_!$
\,\eqref{ctrh-lct-contramod-derived-direct-image-formula} for $\st=\ctr$
also can be constructed by applying the functor~$\ff_!$ termwise to
complexes of coflasque locally cotorsion contraherent cosheaves of
contramodules on~$\fY$.
 Finally, to handle the case $\st=\bctr$, it suffices to recall
the result of Theorem~\ref{lcth-lct-W-contramod-contraderived},
which tells us that for the exact categories $\fY\Ctrh^\lct$ and
$\fX\Ctrh^\lct$ the Becker contraderived category coincides with
the Positselski one.

 Let $\fX$ and $\fY$ be semi-separated Noetherian formal schemes.
 Then for any derived category symbol $\st\ne\co$, $\bco$, we have
a triangulated equivalence $\sD^\st(\fY\Ctrh^\lct_\cdl)\allowbreak\simeq
\sD^\st(\fY\Ctrh^\lct)$ provided by
Corollary~\ref{lcth-lct-contramod-coflasque-derived-equivalences}(b).
 The functor of direct image of cosheaves of $\fO$\+modules~$\ff_!$
induces an exact functor $\ff_!\:\fY\Ctrh^\lct_\cdl\rarrow
\fX\Ctrh^\lct_\cdl$ by Corollary~\ref{codilute-direct-images}(b).
 Hence, for any derived category symbol $\st\ne\co$, $\ctr$, $\bco$,
$\bctr$, the functor~$\ff_!$ takes $\st$\+acyclic complexes in
$\fY\Ctrh^\lct_\cdl$ to $\st$\+acyclic complexes
in $\fX\Ctrh^\lct_\cdl$.
 So the derived functor $\boL\ff_!$
\,\eqref{ctrh-lct-contramod-derived-direct-image-formula} can be
constructed by applying the functor~$\ff_!$ termwise to complexes of
codilute locally cotorsion contraherent cohseaves of contramodules
on~$\fY$.

 For the derived category symbol $\st=\ctr$, one needs to observe
additionally that the full subcategory $\fY\Ctrh^\lct_\cdl$ is closed
under infinite products in $(\fY,\fO_\fY)\Cosh$ by
Sections~\ref{locally-contraherent-cosheaves-of-contramods-subsecn}\+-%
\ref{locally-contraherent-lct-cosheaves-subsecn}
and Lemma~\ref{products-of-codilute}, and the functor
$\ff_!\:(\fY,\fO_\fY)\Cosh\rarrow(\fX,\fO_\fX)\Cosh$ preserves infinite
products as mentioned above.
 Hence the functor~$\ff_!$ also takes Positselski-contraacyclic
complexes in $\fY\Ctrh^\lct_\cdl$ to Positselski-contraacyclic complexes
in $\fX\Ctrh^\lct_\cdl$, and the derived functor functor $\boL\ff_!$
\,\eqref{ctrh-lct-contramod-derived-direct-image-formula} for $\st=\ctr$ 
also can be constructed by applying the functor~$\ff_!$ termwise to
complexes of codilute locally cotorsion contraherent cohseaves of
contramodules on~$\fY$.
 Finally, the case $\st=\bctr$ reduces to that of $\st=\ctr$ by
Theorem~\ref{lcth-lct-W-contramod-contraderived}.

 In the general case of arbitrary locally Noetherian formal schemes
$\fX$ and $\fY$, one can still construct the derived functor
$\boL\ff_!$ \,\eqref{ctrh-lct-contramod-derived-direct-image-formula}
for \emph{some} derived category symbols~$\st$.
 For $\st=-$, $\ctr$, or~$\bctr$, we have a triangulated equivalence
$\sD^\st(\fY\Ctrh^\lct_\cfq)\simeq\sD^\st(\fY\Ctrh^\lct)$ provided by
Corollary~\ref{lcth-lct-contramod-coflasque-derived-equivalences}(a).
 The same arguments as above are applicable, and the functor~$\boL\ff_!$
can be constructed by applying the functor~$\ff_!$ termwise to complexes
of coflasque locally cotorsion contraherent cohseaves of contramodules
on~$\fY$.

 A tricky case is that of $\st=\varnothing$ for arbitrary locally
Noetherian formal schemes $\fX$ and~$\fY$.
 In that case, for all we know, an unbounded complex in
$\fY\Ctrh^\lct_\cfq$ that is acyclic in $\fY\Ctrh^\lct$ \emph{need not}
be acyclic in $\fY\Ctrh^\lct_\cfq$.
 Nevertheless, applying the functor~$\ff_!$ termwise to any complex
in $\fY\Ctrh^\lct_\cfq$ that is acyclic in $\fY\Ctrh^\lct$ produces
a complex in $\fX\Ctrh^\lct_\cfq$ that is acyclic in $\fX\Ctrh^\lct$,
by Corollary~\ref{co-flasque-direct-image-of-acyclic-complexes-cor}(c).
 For this reason, the derived functor $\boL\ff_!$
\,\eqref{ctrh-lct-contramod-derived-direct-image-formula} still can be
constructed by applying the functor~$\ff_!$ termwise to complexes of
coflasque locally cotorsion contraherent cohseaves of contramodules
on~$\fY$.

 Here we are also using the fact that, for every complex of locally
cotorsion contraherent cosheaves of contramodules $\fQ^\bu$ on
a locally Noetherian formal scheme $\fY$, there exists
a quasi-isomorphism $\fE^\bu\rarrow\fQ^\bu$ into $\fQ^\bu$ from
a complex of coflasque locally cotorsion contraherent cosheaves
$\fE^\bu$ on~$\fY$.
 This follows, e.~g., from
Theorem~\ref{lcth-lct-W-contramod-contraderived}.

 Now let us pass to locally contraadjusted contraherent cosheaves of
contramodules.

 Let $\fY$ be a Noetherian formal scheme of finite Krull dimension.
 Then for any derived category symbol $\st\ne\co$, $\bco$, we have
a triangulated equivalence $\sD^\st(\fY\Ctrh_\cfq)\allowbreak\simeq
\sD^\st(\fY\Ctrh)$ provided by
Corollary~\ref{lcth-contramod-finite-krulldim-derived-equivalences}(b).
 The functor of direct image of cosheaves of $\fO$\+modules~$\ff_!$
induces an exact functor $\ff_!\:\fY\Ctrh_\cfq\rarrow
\fX\Ctrh_\cfq$ by Corollary~\ref{co-flasque-direct-image-cor}(b).
 Hence, for any derived category symbol $\st\ne\co$, $\ctr$, $\bco$,
$\bctr$, the functor~$\ff_!$ takes $\st$\+acyclic complexes in
$\fY\Ctrh_\cfq$ to $\st$\+acyclic complexes in $\fX\Ctrh_\cfq$.
 So the derived functor $\boL\ff_!$
\,\eqref{ctrh-contramod-derived-direct-image-formula} can be
constructed by applying the functor~$\ff_!$ termwise to complexes of
coflasque contraherent cohseaves of contramodules on~$\fY$.

 For the derived category symbol $\st=\ctr$, one needs to observe
additionally that the full subcategory $\fY\Ctrh_\cfq$ is closed under
infinite products in $(\fY,\fO_\fY)\Cosh$ according to
Sections~\ref{locally-contraherent-cosheaves-of-contramods-subsecn}
and~\ref{coflasque-subsecn}, and the functor $\ff_!\:(\fY,\fO_\fY)\Cosh
\rarrow(\fX,\fO_\fX)\Cosh$ preserves infinite products as mentioned 
above.
 Hence the functor~$\ff_!$ also takes Positselski-contraacyclic
complexes in $\fY\Ctrh_\cfq$ to Positselski-contraacyclic complexes
in $\fX\Ctrh_\cfq$, and the derived functor functor $\boL\ff_!$
\,\eqref{ctrh-contramod-derived-direct-image-formula} for $\st=\ctr$
also can be constructed by applying the functor~$\ff_!$ termwise to
complexes of coflasque contraherent cosheaves of contramodules on~$\fY$.
 Finally, in the case $\st=\bctr$, it suffices to recall the result of
Theorem~\ref{lcth-W-contramod-contraderived}, which tells us that for
the exact category $\fY\Ctrh$ the Becker contraderived category
coincides with the Positselski one.

 Let $\fX$ and $\fY$ be semi-separated Noetherian formal schemes.
 Then for any derived category symbol $\st\ne\co$, $\bco$, we have
a triangulated equivalence $\sD^\st(\fY\Ctrh_\cdl)\allowbreak\simeq
\sD^\st(\fY\Ctrh)$ provided by
Corollary~\ref{lcth-contramod-ssep-derived-equivalences}(b).
 The functor of direct image of cosheaves of $\fO$\+modules~$\ff_!$
induces an exact functor $\ff_!\:\fY\Ctrh_\cdl\rarrow
\fX\Ctrh_\cdl$ by Corollary~\ref{codilute-direct-images}(a).
 Hence, for any derived category symbol $\st\ne\co$, $\ctr$, $\bco$,
$\bctr$, the functor~$\ff_!$ takes $\st$\+acyclic complexes in
$\fY\Ctrh_\cdl$ to $\st$\+acyclic complexes in $\fX\Ctrh_\cdl$.
 So the derived functor $\boL\ff_!$
\,\eqref{ctrh-contramod-derived-direct-image-formula} can be
constructed by applying the functor~$\ff_!$ termwise to complexes of
codilute contraherent cohseaves of contramodules on~$\fY$.

 For the derived category symbol $\st=\ctr$, one needs to observe
additionally that the full subcategory $\fY\Ctrh_\cdl$ is closed
under infinite products in $(\fY,\fO_\fY)\Cosh$ by
Section~\ref{locally-contraherent-cosheaves-of-contramods-subsecn}
and Lemma~\ref{products-of-codilute}, and the functor
$\ff_!\:(\fY,\fO_\fY)\Cosh\rarrow(\fX,\fO_\fX)\Cosh$ preserves infinite
products.
 Hence the functor~$\ff_!$ also takes Positselski-contraacyclic
complexes in $\fY\Ctrh_\cdl$ to Positselski-contraacyclic complexes
in $\fX\Ctrh_\cdl$, and the derived functor functor $\boL\ff_!$
\,\eqref{ctrh-contramod-derived-direct-image-formula} for $\st=\ctr$ 
also can be constructed by applying the functor~$\ff_!$ termwise to
complexes of codilute contraherent cohseaves of contramodules on~$\fY$.

 When both $\fX$ and $\fY$ are Noetherian formal schemes of finite
Krull dimension, for any derived category symbol $\st\ne\co$, $\bco$
we have a commutative square diagram of compatibility
\begin{equation} \label{ctrh-lcta-lct-der-dir-im-compat-diag}
\begin{gathered}
 \xymatrix{
  \sD^\st(\fY\Ctrh^\lct) \ar@<-2pt>[r] \ar[d]^{\boL\ff_!}
  & \sD^\st(\fY\Ctrh) \ar@<-2pt>@{-}[l] \ar[d]^{\boL\ff_!} \\
  \sD^\st(\fX\Ctrh^\lct) \ar@<-2pt>[r]
  & \sD^\st(\fX\Ctrh) \ar@<-2pt>@{-}[l]
 }
\end{gathered}
\end{equation}
of the left derived direct image functors $\boL\ff_!$
\,(\ref{ctrh-contramod-derived-direct-image-formula}\+-
\ref{ctrh-lct-contramod-derived-direct-image-formula})
with the triangulated equivalences provided by
Proposition~\ref{lcth-contramod-lcta-lct-derived-equivalences}(b).

 The same commutative diagram of triangulated functors and
triangulated equivalences~\eqref{ctrh-lcta-lct-der-dir-im-compat-diag}
is also valid in the case of semi-separated Noetherian formal schemes
$\fX$ and $\fY$ for any one of the two derived category symbols
$\st=+$ or~$\varnothing$, with the triangulated equivalences provided
by Theorem~\ref{lcth-lcta-lct-unbounded-derived-equivalence}.

\subsection{Triangulated functors of inverse image}
\label{triangulated-inverse-images-subsecn}
 In the context of schemes, many derived functors of inverse image
can be constructed using resolutions by flat quasi-coherent sheaves
and coresolutions by locally injective locally contraherent cosheaves
(see~\cite[Sections~4.10\+-4.11 and~5.8]{Pcosh} for a detailed
discussion).
 In the realm of formal schemes, we \emph{do not} seem to have any flat
quasi-coherent torsion sheaves or locally injective locally contraherent
cosheaves of contramodules.
 For this reason, we restrict ourselves to underived inverse images
in this brief section.
 
 Let $\ff\:\fY\rarrow\fX$ be a flat tight morphism of locally
Noetherian formal schemes.
 Then the inverse image functor $\ff^*\:\fX\Tors\rarrow\fY\Tors$
\,\eqref{inverse-image-qcoh-tors} is exact and preserves infinite
direct sums.
 Therefore, it induces triangulated functors
\begin{equation} \label{qcoh-derived-cats-underived-inverse-image}
 \ff^*\:\sD^\st(\fX\Tors)\lrarrow\sD^\st(\fY\Tors)
\end{equation}
between the derived categories with all the derived category symbols
$\st\ne\ctr$, $\bco$, $\bctr$.
 The triangulated functor $\ff^*$
\,\eqref{qcoh-derived-cats-underived-inverse-image} is also
well-defined in the case of $\st=\bco$, as one can show by reducing
that case to the case of $\st=\co$ using
Theorem~\ref{qcoh-tors-coderived}.
 
 The triangulated functor of inverse image $\ff^*$
\,\eqref{qcoh-derived-cats-underived-inverse-image} is left adjoint
to the right derived functor of direct image $\boR\ff_*$
\,\eqref{qcoh-tors-derived-direct-image-formula} whenever both
are defined.
 This follows from the constructions and
the adjunction~\eqref{qcoh-tors-inverse-direct-adjunction}.

 Let $\bW$ be an open covering of $\fX$ and $\bT$ be an open covering
of $\fY$ such that the morphism~$\ff$ is $(\bW,\bT)$\+coaffine (in
the sense of the definition in
Section~\ref{inverse-images-of-qcoh-tors-subsecn}).
 Then the inverse image functor $\ff^!\:\fX\Lcth^\lct_\bW\rarrow
\fY\Lcth^\lct_\bT$ \,\eqref{inverse-image-lcth-lct-W-T-contramod} is
exact and preserves infinite products.
 Therefore, it induces triangulated functors
\begin{equation} \label{lcth-lct-derived-cats-underived-inverse-image}
 \ff^!\:\sD^\st(\fX\Lcth^\lct_\bW)\lrarrow\sD^\st(\fY\Lcth^\lct_\bT)
\end{equation}
between the derived categories with all the derived category symbols
$\st\ne\co$, $\bco$, $\bctr$.
 The triangulated functor $\ff^!$
\,\eqref{lcth-lct-derived-cats-underived-inverse-image} is also
well-defined in the case of $\st=\bctr$, as one can show by reducing
that case to the case of $\st=\ctr$ using
Theorem~\ref{lcth-lct-W-contramod-contraderived}.

 Whenever the left derived direct image functor $\boL\ff_!$
\,\eqref{ctrh-lct-contramod-derived-direct-image-formula} for
locally cotorsion contraherent cosheaves of contramodules can be
defined by one of the constructions from
Section~\ref{derived-direct-images-of-ctrh-contramod-subsecn},
the same construction also provides a similar left derived
direct image functor
\begin{equation}
 \label{lcth-lct-contramod-derived-direct-image-formula}
 \boL\ff_!\:\sD^\st(\fY\Lcth^\lct_\bT)\lrarrow
 \sD^\st(\fX\Lcth^\lct_\bW)
\end{equation}
acting between the derived categories of locally cotorsion
$\bT$\+locally contraherent cosheaves of contramodules on $\fY$
and locally cotorsion $\bW$\+locally contraherent cosheaves of
contramodules on~$\fX$ (for any open coverings $\bW$ of $\fX$ and
$\bT$ of~$\fY$).
 In fact, in all those cases except possibly the ``tricky' one,
the inclusions of exact categories induce triangulated equivalences
$\sD^\st(\fY\Ctrh^\lct)\simeq\sD^\st(\fY\Lcth^\lct_\bT)$
and $\sD^\st(\fX\Ctrh^\lct)\simeq\sD^\st(\fX\Lcth^\lct_\bW)$ by
the respective part of
Corollary~\ref{lcth-lct-contramod-coflasque-derived-equivalences}.

 The triangulated functor of inverse image $\ff^!$
\,\eqref{lcth-lct-derived-cats-underived-inverse-image} is right
adjoint to the left derived functor of direct image $\boL\ff_!$
\,\eqref{lcth-lct-contramod-derived-direct-image-formula} whenever
both are defined.
 This follows from the constructions and
the adjunction~\eqref{cosheaves-direct-inverse-adjunction}.

 Let $\ff\:\fY\rarrow\fX$ be a very flat tight $(\bW,\bT)$\+coaffine
morphism of locally Noetherian formal schemes.
 Then the inverse image functor $\ff^!\:\fX\Lcth_\bW\rarrow
\fY\Lcth_\bT$ \,\eqref{inverse-image-lcth-W-T-contramod} is exact
and preserves infinite products.
 Therefore, it induces triangulated functors
\begin{equation} \label{lcth-derived-cats-underived-inverse-image}
 \ff^!\:\sD^\st(\fX\Lcth_\bW)\lrarrow\sD^\st(\fY\Lcth_\bT)
\end{equation}
between the derived categories with all the derived category symbols
$\st\ne\co$, $\bco$, $\bctr$.
 In the case when the formal scheme $\fX$ is Noetherian of finite
Krull dimension, the triangulated functor $\ff^!$
\,\eqref{lcth-derived-cats-underived-inverse-image} is also
well-defined in the case of $\st=\bctr$, as one can show by reducing
that case to the case of $\st=\ctr$ using
Theorem~\ref{lcth-W-contramod-contraderived}.

 Whenever the left derived direct image functor $\boL\ff_!$
\,\eqref{ctrh-contramod-derived-direct-image-formula} for
contraherent cosheaves of contramodules can be defined by one of
the constructions from
Section~\ref{derived-direct-images-of-ctrh-contramod-subsecn},
the same construction also provides a similar left derived
direct image functor
\begin{equation}
 \label{lcth-contramod-derived-direct-image-formula}
 \boL\ff_!\:\sD^\st(\fY\Lcth_\bT)\lrarrow\sD^\st(\fX\Lcth_\bW)
\end{equation}
acting between the derived categories of $\bT$\+locally contraherent
cosheaves of contramodules on $\fY$ and $\bW$\+locally contraherent
cosheaves of contramodules on~$\fX$ (for any open coverings $\bW$ of
$\fX$ and $\bT$ of~$\fY$).
 In fact, in all those cases, the inclusions of exact categories induce
triangulated equivalences
$\sD^\st(\fY\Ctrh)\simeq\sD^\st(\fY\Lcth_\bT)$
and $\sD^\st(\fX\Ctrh)\simeq\sD^\st(\fX\Lcth_\bW)$ by
Corollaries~\ref{lcth-contramod-ssep-derived-equivalences}
and~\ref{lcth-contramod-finite-krulldim-derived-equivalences}.

 In this context, too, the triangulated functor of inverse image $\ff^!$
\,\eqref{lcth-derived-cats-underived-inverse-image} is right
adjoint to the left derived functor of direct image $\boL\ff_!$
\,\eqref{lcth-contramod-derived-direct-image-formula} whenever
both are defined.
 This follows from the constructions and
the adjunction~\eqref{cosheaves-direct-inverse-adjunction}.

\Section{Dedualizing Complexes}  \label{dedualizing-secn}

 The aim of this section is to globalize the MGM duality/equivalence
theorem, as stated in~\cite[Corollary~3.5]{Pmgm}, or more generally,
\cite[Example~5.8 and Theorem~5.10]{Pmgm}, to Noetherian formal
schemes of finite Krull dimension, not necessarily semi-separated
(or alternatively, semi-separated but of possibly infinite
Krull dimension).
 While acknowledging the prior art of the papers~\cite{AJL0,AJL2},
we follow the approach of~\cite[Section~5]{Pmgm}.
 The goal is to state and prove the MGM duality in the form of
an equivalence between the conventional or absolute derived categories
of the abelian category of quasi-coherent torsion sheaves $\fX\Tors$
and the exact category of contraherent cosheaves of contramodules
$\fX\Ctrh$ on a formal scheme~$\fX$.

\subsection{Dedualizing complexes for adic Noetherian rings}
\label{dedualizing-complexes-over-rings-subsecn}
 This section is mostly a reminder of background material
in commutative algebra related to the MGM duality.
 We follow the paper~\cite{PSY} by Porta, Shaul, and Yekutieli, as well
as the papers~\cite{Pmgm,Pdc,Ppt} by the present author.
 Our notation is closer to the notation in the latter three papers.

 Let $R$ be a commutative ring and $s\in R$ be an element.
 We denote by $K^\bu_\infty(R;s)$ the \emph{infinite dual Koszul
complex} $R\rarrow R[s^{-1}]$, with the term $R$ placed in
the cohomological degree~$0$ and the term $R[s^{-1}]$ placed in
the cohomological degree~$1$.
 Given a finite sequence of elements $s_1$,~\dots, $s_m\in R$,
the infinite dual Koszul complex $K^\bu_\infty(R;\bs)$ is
constructed as the tensor product
$$
 K^\bu_\infty(R;\bs)=K^\bu_\infty(R;s_1)\ot_R\dotsb\ot_R
 K^\bu_\infty(R;s_m).
$$
 Here $\bs$~is a shorthand notation for the sequence of elements 
$s_1$,~\dots,~$s_m$.
 Obviously, $K^\bu_\infty(R;\bs)$ is a finite complex of flat
(in fact, countably presented and very flat) $R$\+modules
concentrated in the cohomological degrees from~$0$ to~$m$.

\begin{lem} \label{infinite-dual-Koszul-has-torsion-cohomology}
 Let $R$ be a commutative ring and $s_1$,~\dots, $s_m\in R$ be
a finite sequence of elements generating an ideal $I\subset R$.
 Then the cohomology $R$\+modules of the infinite dual Koszul
complex $K^\bu_\infty(R;\bs)$ are $I$\+torsion $R$\+modules.
\end{lem}

\begin{proof}
 This is~\cite[Lemma~1.1(a)]{Pmgm}.
\end{proof}

 Now let $t_1$,~\dots, $t_n\in R$ be another finite sequence of
elements, denoted for brevity by~$\bt$.
 Let us denote the concatenation of the two sequences $s_1$,~\dots,
$s_m$, $t_1$,~\dots, $t_m$ by $\bs\sqcup\bt$.
 Then there are natural morphisms of complexes
$K^\bu_\infty(R;\>\bs\sqcup\bt)\rarrow K^\bu_\infty(R;\bs)$
and $K^\bu_\infty(R;\>\bs\sqcup\bt)\rarrow K^\bu_\infty(R;\bt)$
induced by the projection morphisms of complexes
$K^\bu_\infty(R,t_i)\rarrow R$ and $K^\bu_\infty(R,s_j)\rarrow R$
(for $1\le i\le n$ and $1\le j\le m$).

\begin{lem} \label{infinite-dual-Koszul-independent-of-generators}
 Let $R$ be a commutative ring.
 Let $s_1$,~\dots, $s_m\in R$ be a sequence of elements generating
an ideal $I\subset R$, and let $t_1$,~\dots, $t_n\in R$ be
a sequence of elements generating an ideal $J\subset R$.
 Assume that the ideals $I$ and $J$ have the same radical in the ring
$R$, that is $\sqrt{I}=\sqrt{J}$.
 Then the morphisms of complexes
$$
 K^\bu_\infty(R;\bs) \llarrow K^\bu_\infty(R;\>\bs\sqcup\bt)
 \lrarrow K^\bu_\infty(R;\bt)
$$
are quasi-isomorphisms.
\end{lem}

\begin{proof}
 This is~\cite[proof of Corollary~6.2]{PSY},
\cite[Proposition~2.20]{Yek3}, or~\cite[Lemma~2.1]{Pdc}.
\end{proof}

 Given a homomorphism of commutative rings $f\:R\rarrow S$ and
a finite sequence of elements $r_1$,~\dots, $r_m\in R$, denoted for
brevity by $\br=(r_1,\dotsc,r_m)$, we will use the notation
$f(\br)$ for the sequence of elements $f(r_1)$,~\dots,
$f(r_m)\in S$.

\begin{lem} \label{infinite-dual-Koszul-base-change}
 Let $f\:R\rarrow S$ be a homomorphism of commutative rings, and let
$r_1$,~\dots, $r_m\in R$ be a sequence of elements, denoted by\/~$\br$.
 Then the natural morphism of complexes of $R$\+modules
$K^\bu_\infty(R;\br)\rarrow K^\bu_\infty(S;f(\br))$ induces
an isomorphism of complexes of $S$\+modules
$$
 S\ot_R K^\bu_\infty(R;\br)\simeq K^\bu_\infty(S;f(\br)).
$$
\end{lem}

\begin{proof}
 Obvious.
\end{proof}

 A detailed discussion of \emph{weakly proregular} finitely generated
ideals in commutative rings can be found
in~\cite[Sections~4 and~6]{PSY}, \cite[Sections~1 and~2]{Pmgm},
\cite[Section~5]{Yek2}, \cite[Section~3]{Yek3}, \cite[Section~3]{Pdc},
or~\cite[Section~1]{Ppt}.
 Let us mention the following version of the definition based
on~\cite[Theorem~4.24]{PSY}: a finitely generated ideal $I\subset R$
with a sequence of generators $\bs=(s_1,\dotsc,s_m)$ is called weakly
proregular if, for every injective $R$\+module $J$, one has
$H^i(K^\bu_\infty(R;\bs)\ot_RJ)=0$ for all $i>0$.

 If the radicals of two finitely generated ideals $I\subset R$ and
$J\subset R$ coincide, that is $\sqrt{I}=\sqrt{J}$, then the ideal
$I$ is weakly proregular if and only if the ideal $J$ is weakly
proregular~\cite[Corollaries~6.2 and~6.3]{PSY}.
 So we will say that an adic topological ring $R$ is weakly proregular
if any one (or equivalently, all) of its ideals of definition is
weakly proregular.
 We recall once again that \emph{all ideals in Noetherian commutative
rings are weakly proregular}~\cite[Theorem~4.34]{PSY},
\cite[Section~1]{Pmgm}, \cite[Theorem~3.3]{Yek3} (this was already
mentioned in Section~\ref{prelim-contramodules-subsecn}).

 As usual, we are mainly interested in adic Noetherian rings in
this paper; but for the sake of clarity and natural generality, we
formulate the results of this section for weakly proregular adic
topological rings.
 Notice that all contramodule modules over a weakly proregular adic
topological ring are quotseparated by~\cite[Corollary~3.7]{Pdc};
see also~\cite[Remark~3.8]{Pdc}.

\begin{lem} \label{cohomology-of-infdualkoszul-completion-invariant}
 Let $R$ be a commutative ring and $s_1$,~\dots, $s_m\in R$ be a finite
sequence of elements generating a weakly proregular ideal $I\subset R$.
 Denote, as usual, by\/ $\Lambda_I(R)$ the $I$\+adic completion of
the ring $R$, and by\/ $\lambda_{I,R}\:R\rarrow\Lambda_I(R)$
the completion map.
 Then the natural morphism of complexes of $R$\+modules
$$
 K^\bu_\infty(R;\bs)\lrarrow
 K^\bu_\infty(\Lambda_I(R);\lambda_{I,R}(\bs))
$$
is a quasi-isomorphism.
\end{lem}

\begin{proof}
 This is a part of~\cite[Lemma~5.3(b)]{Pmgm} (the weak proregularity
is the running assumption in much of~\cite[Section~5]{Pmgm}).
\end{proof}

 Let $R$ be an adic topological ring and $\st=\bb$, $+$, $-$,
or~$\varnothing$ be a conventional derived category symbol.
 Then we denote by $\sD^\st_\tors(R\Modl)$ the full triangulated
subcategory of $\sD^\st(R\Modl)$ consisting of all the complexes
whose $R$\+modules of cohomology are torsion $R$\+modules.

 The following theorem is a stronger version of
Corollary~\ref{Noetherian-torsion-Ext-isomorphism-cor} above.

\begin{thm} \label{torsion-cohomology-modules-triang-equivalence}
 Let $R$ be a weakly proregular adic topological ring and\/ $\st=\bb$,
$+$, $-$, or\/~$\varnothing$ be a conventional derived category symbol.
 Then the inclusion of abelian categories $R\Modl_\tors\rarrow
R\Modl$ induces a triangulated equivalence
$$
 \sD^\st(R\Modl_\tors)\simeq\sD^\st_\tors(R\Modl).
$$
\end{thm}

\begin{proof}
 This is~\cite[Corollary~1.4]{Pmgm}.
 The weak proregularity is assumption necessary;
see~\cite[Theorem~4.1]{Pdc}.
\end{proof}

 For any finite sequence of elements $\bs=(s_1,\dotsc,s_m)$ in
a commutative ring $R$, the infinite dual Koszul complex
$K^\bu_\infty(R;\bs)$ is quasi-isomorphic to a finite complex of
projective $R$\+modules.
 This follows already from the fact that, for any element $s\in R$,
the projective dimension of the $R$\+module $R[s^{-1}]$ does not
exceed~$1$.
 Various explicit constructions of an $(m+\nobreak1)$\+term complex
of countably generated free $R$\+modules, sitting in the cohomological
degrees from~$0$ to~$m$, endowed with a natural quasi-isomorphism onto
$K^\bu_\infty(R;\bs)$ can be found in~\cite[Definition~5.1]{PSY}
and~\cite[Section~2]{Pmgm}, under notation such as
$\operatorname{Tel}^\bu(R;\bs)$ or $T^\bu(R;\bs)$.
 We will pick any one of such complexes and denote it by $T^\bu(R;\bs)$.
 So $T^\bu(R;\bs)$ is a finite complex of projective $R$\+modules
endowed with a quasi-isomorphism of complexes of $R$\+modules
$T^\bu(R;\bs)\rarrow K^\bu_\infty(R;\bs)$.

 The notion of a \emph{dedualizing complex} for a weakly proregular
finitely generated ideal $I$ in a commutative ring $R$ was introduced
in~\cite[Section~5]{Pmgm}.
 A dedualizing complex was defined as a finite complex of $I$\+torsion
$R$\+modules satisfying certain conditions~(i\+-iii).
 Then a specific construction of one particular example of a dedualizing
complex for any such ring $R$ and ideal $I$ was spelled out
in~\cite[Example~5.8]{Pmgm}.
 In our exposition in the present paper, we ignore the abstract
definition of \emph{a} dedualizing complex from~\cite[Section~5]{Pmgm},
but work with \emph{the} dedualizing complex produced by
the concrete construction of~\cite[Example~5.8]{Pmgm} instead.

 Let $R$ be a weakly proregular adic topological ring.
 Choose a sequence of elements $s_1$,~\dots, $s_m\in R$ generating
some ideal of definition $I\subset R$, and consider the infinite
dual Koszul complex $K^\bu_\infty(R;\bs)$.
 By Lemma~\ref{infinite-dual-Koszul-has-torsion-cohomology},
\,$K^\bu(R;\bs)$ is a (finite) complex of $R$\+modules with torsion
cohomology modules.

 By Lemma~\ref{infinite-dual-Koszul-independent-of-generators},
the object $K^\bu(R;\bs)\in\sD^\bb(R\Modl)$ does not depend on
the choice of a particular ideal of definition $I$ and its particular
sequence of generators~$\bs$ up to a natural isomorphism in the derived
category $\sD^\bb(R\Modl)$.
 By Theorem~\ref{torsion-cohomology-modules-triang-equivalence},
the object $K^\bu(R;\bs)\in\sD^\bb_\tors(R\Modl)$ can be represented
by a finite complex of torsion $R$\+modules, which is uniquely
defined up to a natural isomorphism in $\sD(R\Modl_\tors)$.
 
 We denote the resulting object by $B^\bu=B^\bu_R\in
\sD^\bb(R\Modl_\tors)$ and call it \emph{the dedualizing complex of
torsion modules} for a weakly proregular adic topological ring~$R$.
 The isomorphism $K^\bu(R;\bs)\simeq B^\bu_R$ in $\sD^\bb(R;\bs)$
can be represented by a fraction (``roof'') of quasi-isomorphisms of
complexes of $R$\+modules $K^\bu_\infty(R;\bs)\larrow T^\bu(R;\bs)
\rarrow B^\bu_R$.

 Recall the natural equivalence (in fact, isomorphism) of categories
$R\Modl_\tors\simeq\Lambda(R)\Modl_\tors$ mentioned in the beginning
of Section~\ref{prelim-inj-proj-flat-torsion-contra-subsecn}.
 With this isomorphism of categories in mind,
Lemma~\ref{cohomology-of-infdualkoszul-completion-invariant} claims
that the passage from $R$ to $\Lambda(R)$ does not change
the dedualizing complex.

\begin{lem} \label{T-B-Hom-tensor-quasi-isomorphism}
 Let $R$ be a weakly proregular adic topological ring and\/
$\bs=(s_1,\dotsc,s_m)$ be a finite sequence of elements generating
an ideal of definition in~$R$.
 Then \par
\textup{(a)} for any complex of injective objects $J^\bu$ in
the abelian category of torsion $R$\+modules $R\Modl_\tors$,
the quasi-isomorphism of complexes of $R$\+modules $T^\bu(R;\bs)
\rarrow B^\bu_R$ induces a quasi-isomorphism of complexes of
$R$\+modules\/ $\Hom_R(B^\bu_R,J^\bu)\rarrow
\Hom_R(T^\bu(R;\bs),J^\bu)$; \par
\textup{(b)} for any complex of flat contramodule $R$\+modules
$G^\bu$, the quasi-isomorphism of complexes of $R$\+modules
$T^\bu(R;\bs)\rarrow B^\bu_R$ induces a quasi-isomorphism of
complexes of $R$\+modules $T^\bu(R;\bs)\ot_RG^\bu\rarrow
B^\bu_R\ot_RG^\bu$.
\end{lem}

\begin{proof}
 Since the totalization of any finite acyclic complex of unbounded
complex is acyclic (and even absolutely acyclic), it suffices to
consider the cases of a one-term complex $J^\bu=J$ corresponding to
an injective object $J\in R\Modl_\tors$ in part~(a), and a one-term
complex $G^\bu=G$ corresponding to a flat contramodule $R$\+module
$G$ in part~(b).
 Then, for a Noetherian ring $R$, the assertions are obvious, as
all injective objects of $R\Modl_\tors$ are injective in $R\Modl$
by Lemma~\ref{injective-torsion-modules-characterizations}, and all
flat contramodule $R$\+modules are flat $R$\+modules by
Lemma~\ref{flat-contramodules-characterizations}.
 In the general (weakly proregular) context, part~(a) becomes
a particular case of~\cite[Lemma~5.5(b)]{Pmgm}
or~\cite[Lemma~3.6(a)]{Ppt}, and part~(b) becomes a particular case
of~\cite[Corollary~3.8]{Ppt}.
\end{proof}

\begin{cor} \label{Hom-tensor-B-preserves-acyclicity-boundedness-cor}
 Let $R$ be a weakly proregular adic topological ring.
 In this setting: \par
\textup{(a)} if a complex $J^\bu$ of injective objects in
$R\Modl_\tors$ is acyclic as a complex in $R\Modl_\tors$, then
the complex of contramodule $R$\+modules\/ $\Hom_R(B_R^\bu,J^\bu)$
is acyclic  as well; \par
\textup{(b)} if a complex $J^\bu$ of injective objects in
$R\Modl_\tors$ has bounded above (respectively, below) cohomology,
then the complex of contramodule $R$\+modules\/ $\Hom_R(B_R^\bu,J^\bu)$
also has bounded above (resp., below) cohomology; \par
\textup{(c)} if a complex $G^\bu$ of flat contramodule $R$\+modules
is acyclic as a complex in $R\Modl_\ctra$, then the complex of
torsion $R$\+modules $B_R^\bu\ot_RG^\bu$ is acyclic as well; \par
\textup{(d)} if a complex $G^\bu$ of flat contramodule $R$\+modules
has bounded below (respectively, above) cohomology, then the complex
of torsion $R$\+modules $B_R^\bu\ot_RG^\bu$ also has bounded below
(resp., above) cohomology.
\end{cor}

\begin{proof}
 Since $T^\bu(R;\bs)$ is a finite complex of projective $R$\+modules,
the assertions~(a\+-b) follow from
Lemma~\ref{T-B-Hom-tensor-quasi-isomorphism}(a), and
the assertions~(c\+-d) follow from
Lemma~\ref{T-B-Hom-tensor-quasi-isomorphism}(b).
\end{proof}

\begin{prop} \label{Hom-tensor-B-quasi-isomorphisms-prop}
 Let $R$ be a weakly proregular adic topological ring.
 Let $J^\bu$ be a complex of injective objects in $R\Modl_\tors$,
and let $G^\bu$ be a complex of flat contramodule $R$\+modules.
 In this context: \par
\textup{(a)} if $G^\bu\rarrow\Hom_R(B_R^\bu,J^\bu)$ is
a quasi-isomorphism of complexes of (contramodule) $R$\+modules,
then the composition $B_R^\bu\ot_RG^\bu\rarrow
B_R^\bu\ot_R\Hom_R(B_R^\bu,J^\bu)\rarrow J^\bu$ is a quasi-isomorphism
of complexes of (torsion) $R$\+modules; \par
\textup{(b)} if $B_R^\bu\ot_RG^\bu\rarrow J^\bu$ is
a quasi-isomorphism of complexes of (torsion) $R$\+modules, then
the composition\/ $G^\bu\rarrow\Hom_R(B_R^\bu,\>B_R^\bu\ot_RG^\bu)
\rarrow\Hom_R(B_R^\bu,J^\bu)$ is a quasi-isomorphism of complexes
of (contramodule) $R$\+modules.
\end{prop}

\begin{proof}
 Part~(a): let $\bs=(s_1,\dotsc,s_m)$ be a finite sequence of
generators of an ideal of definition in~$R$.
 Put $T^\bu=T^\bu(R;\bs)$ for brevity.
 Then we have a commutative diagram of morphisms of complexes
of $R$\+modules
$$
 \xymatrix{
  T^\bu\ot_RG^\bu \ar[r]^-{\mathrm{quis}}\ar[d]_-{\mathrm{quis}}
  & T^\bu\ot_R\Hom_R(B_R^\bu,J^\bu) \ar[r]^-{\mathrm{quis}} \ar[d]
  & T^\bu\ot_R\Hom_R(T^\bu,J^\bu) \ar[d]^-{\mathrm{quis}} \\
  B_R^\bu\ot_RG^\bu \ar[r]
  & B_R^\bu\ot_R\Hom_R(B_R^\bu,J^\bu) \ar[r]
  & J^\bu
 }
$$
 The leftmost vertical morphism is a quasi-isomorphism by
Lemma~\ref{T-B-Hom-tensor-quasi-isomorphism}(b).
 The leftmost upper horizontal morphism is a quasi-isomorphism
by the assumption of part~(a), and since $T^\bu$ is a bounded above
complex of flat $R$\+modules.
 The rightmost upper horizontal morphism is a quasi-isomorphism
by Lemma~\ref{T-B-Hom-tensor-quasi-isomorphism}(a), and since
$T^\bu$ is a bounded above complex of flat $R$\+modules.
 The rightmost vertical morphism is a quasi-isomorphism
by~\cite[Lemma~5.2(a)]{Pmgm}.
 It follows that the lower horizontal composition is
a quasi-isomorphism, too.

 Part~(b): in the notation of the proof of part~(a), we have
a commutative diagram of morphisms of complexes of $R$\+modules
$$
 \xymatrix{
  G^\bu \ar[r]\ar[d]_-{\mathrm{quis}}
  & \Hom_R(B_R^\bu,\,B_R^\bu\ot_RG^\bu) \ar[r] \ar[d]
  & \Hom_R(B_R^\bu,J^\bu) \ar[d]^-{\mathrm{quis}} \\
  \Hom_R(T^\bu,\,T^\bu\ot_RG^\bu) \ar[r]^-{\mathrm{quis}}
  & \Hom_R(T^\bu,\,B_R^\bu\ot_RG^\bu) \ar[r]^-{\mathrm{quis}}
  &\Hom_R(T^\bu,J^\bu)
 }
$$
 The rightmost vertical morphism is a quasi-isomorphism by
Lemma~\ref{T-B-Hom-tensor-quasi-isomorphism}(a).
 The leftmost vertical morphism is a quasi-isomorphism
by~\cite[Lemma~5.2(b)]{Pmgm}.
 The leftmost lower horizontal morphism is a quasi-isomorphism
by Lemma~\ref{T-B-Hom-tensor-quasi-isomorphism}(b), and since
$T^\bu$ is a bounded above complex of projective $R$\+modules.
 The rightmost lower horizontal morphism is a quasi-isomorphism
by the assumption of part~(b), and since $T^\bu$ is a bounded above
complex of projective $R$\+modules.
 It follows that the upper horizontal composition is
a quasi-isomorphism, too.
\end{proof}

\begin{lem} \label{dedualizing-complex-derived-endomorphisms}
 Let $R$ be a weakly proregular adic topological ring, $\Lambda(R)$ be
its completion, and $B^\bu\in\sD^\bb_\tors(R\Modl)$ be the dedualizing
complex for~$R$.
 Then \par
\textup{(a)} one has\/ $\Hom_{\sD^\bb(R\Modl_\tors)}(B^\bu,B^\bu[n])=0$
for all\/ $0\ne n\in\boZ$; \par
\textup{(b)} the natural action of\/ $\Lambda(R)$ on $B^\bu$ (extending
by continuity the action of $R$ on~$B^\bu$) induces a ring isomorphism
$\Lambda(R)\simeq\Hom_{\sD^\bb(R\Modl_\tors)}(B^\bu,B^\bu)$.
\end{lem}

\begin{proof}
 This is~\cite[property~(ii) in Section~5]{Pmgm}, which holds for our
construction of the dedualizing complex by~\cite[Example~5.8]{Pmgm}.
\end{proof}

 Let $f\:R\rarrow S$ be a flat tight continuous map of adic
topological rings (in the sense of the definition in
Section~\ref{adic-rings-subsecn}).
 Then, by Lemma~\ref{flat-map-of-adic-rings-direct-image-lemma}(a),
there is an exact functor of extension of scalars
$f^*\:R\Modl_\tors\rarrow S\Modl_\tors$ defined by the rule
$f^*(M)=S\ot_RM$ for all $M\in R\Modl_\tors$.
 Consequently, we have the induced triangulated functors
acting between the derived categories,
$$
 f^*\:\sD^\st(R\Modl_\tors)\lrarrow\sD^\st(S\Modl_\tors)
$$
for $\st=\bb$, $+$, $-$, or~$\varnothing$.

\begin{rem} \label{completion-preserves-not-reflects-wpr}
 To put the proof of the next lemma in a proper context, let us
mention the known fact that, for any weakly proregular adic
topological ring $R$, the adic topological ring $\Lambda(R)$ is
weakly proregular, too~\cite[Theorem~0.8 or~7.5]{Yek4}.
 For another proof, see~\cite[Example~5.2(2)]{Pper}, based
on~\cite[Lemma~5.3(b) or~5.4(b)]{Pmgm}.

 The converse assertion is not true.
 For a counterexample, pick a commutative ring $R$ with an element
$s\in R$ such that the $s$\+torsion in $R$ is nonzero and
$s$\+divisible, e.~g., take the trivial extension $R=k[s]\oplus
k[s,s^{-1}]/k[s]$ of the polynomial ring $k[s]$ in one variable~$s$
over a field~$k$ by the Pr\"ufer module $k[s,s^{-1}]/k[s]$ over
$k[s]$; and endow $R$ with the $s$\+adic topology;
cf.~\cite[Remark~3.8]{Pdc}).
\end{rem}

\begin{lem} \label{dedualizing-complex-extension-of-scalars}
 Let $f\:R\rarrow S$ be a flat tight continuous homomorphism
of weakly proregular adic topological rings, and let 
$B_R^\bu\in\sD^\bb(R\Modl_\tors)$ and $B_S^\bu\in\sD^\bb(S\Modl_\tors)$
be the dedualizing complexes of torsion modules over $R$ and~$S$.
 Then there is a natural isomorphism \emergencystretch=1em
$$
 \rho^R_S\:f^*(B_R^\bu)\simeq B_S^\bu
$$
in the derived category\/ $\sD^\bb(S\Modl_\tors)$.
\end{lem}

\begin{proof}
 First of all we notice that, in view of
Lemma~\ref{cohomology-of-infdualkoszul-completion-invariant},
we can assume the adic topological rings $R$ and $S$ to be
separated and complete.
 Then $S$ is a flat contramodule $R$\+module.
 Let $\br=(r_1,\dotsc,r_m)$ be a finite sequence of generators of
an ideal of definition in~$R$.
 Using Lemma~\ref{infinite-dual-Koszul-base-change}, it remains to
prove that the fraction (``roof'') of quasi-isomorphisms of complexes
of $R$\+modules $K^\bu_\infty(R;\br)\larrow T^\bu(R;\br)\rarrow
B_R^\bu$ induces a fraction of quasi-isomorphisms of complexes of
$S$\+modules  $S\ot_RK^\bu_\infty(R;\br)\larrow S\ot_RP^\bu
\rarrow S\ot_RB_R^\bu$.
 Now the map $S\ot_RT^\bu(R;\br)\rarrow S\ot_RB_R^\bu$ is
a quasi-isomorphism by Lemma~\ref{T-B-Hom-tensor-quasi-isomorphism}(b).
 The map $S\ot_RT^\bu(R;\br)\rarrow S\ot_RK^\bu_\infty(R;\br)$ is
a quasi-isomorphism for the simple reason that $T^\bu(R;\br)\rarrow
K^\bu_\infty(R;\br)$ is a quasi-isomorphism of bounded above
complexes of flat $R$\+modules.
\end{proof}

\subsection{Gluing complexes of quasi-coherent torsion sheaves}
\label{gluing-in-D-Tors-subsecn}
 In this section we spell out a gluing construction for objects of
the (conventional unbounded) derived categories of quasi-coherent
torsion sheaves with respect to an open covering of a Noetherian
formal scheme.
 Our discussion is inspired by a very general approach
of~\cite[Section Tag~0D65]{SP} (which in turn goes back to
an even more general~\cite[Th\'eor\`eme~3.2.4]{BBD}).
 Throughout this section, we use the constructions of the functors of
derived direct image $\boR\ff_*$
\,\eqref{qcoh-tors-derived-direct-image-formula} and
underived inverse image $\ff^*$
\,\eqref{qcoh-derived-cats-underived-inverse-image}
(for $\st=\varnothing$) from
Sections~\ref{derived-direct-images-of-qcoh-tors-subsecn}
and~\ref{triangulated-inverse-images-subsecn}.

 We work with a fixed locally Noetherian formal scheme~$\fX$.
 For an open formal subscheme $\fW\subset\fX$, we denote the open
immersion morphism by $\fj_\fW\:\fW\rarrow\fX$.
 Given two open formal subchemes $\fU$ and $\fV\subset\fX$, we use
the notation $\fk_\fU\:(\fU\cap\fV)\rarrow\fU$ and
$\fk_\fV\:(\fU\cap\fV)\rarrow\fV$ for the two open immersion morphisms
from the intersection.

\begin{lem} \label{qcoh-tors-mayer-vietoris-triangle-lemma}
 Let\/ $\fX$ be a locally Noetherian formal scheme covered by
two open formal subschemes, $\fX=\fU\cup\fV$.
 Following the notation above, we have the identity open immersion
morphisms\/ $\fj_\fU\:\fU\rarrow\fX$, \ $\fj_\fV\:\fV\rarrow\fX$,
and\/ $\fj_{\fU\cap\fV}\:\fU\cap\fV\rarrow\fX$.
 Let\/ $\cM^\bu\in\sD(\fX\Tors)$ be an object of the derived
category of quasi-coherent torsion sheaves on\/~$\fX$.
 Then there is a natural distinguished triangle of Mayer--Vietoris
\begin{equation} \label{qcoh-tors-mayer-vietoris-triangle-formula}
 \cM^\bu\lrarrow\boR\fj_\fU{}_*\fj_\fU^*\cM^\bu\oplus
 \boR\fj_\fV{}_*\fj_\fV^*\cM^\bu\lrarrow
 \boR\fj_{\fU\cap\fV}{}_*\fj_{\fU\cap\fV}^*\cM^\bu\lrarrow\cM^\bu[1]
\end{equation}
in the triangulated category\/ $\sD(\fX\Tors)$.
\end{lem}

\begin{proof}
 For any complex of quasi-coherent torsion sheaves $\cE^\bu$ on $\fX$,
there is a natural short exact sequence
\begin{equation} \label{qcoh-tors-mayer-vietoris-short-seq-formula}
 0\lrarrow\cE^\bu\lrarrow\fj_\fU{}_*\fj_\fU^*\cE^\bu\oplus
 \fj_\fV{}_*\fj_\fV^*\cM^\bu\lrarrow
 \fj_{\fU\cap\fV}{}_*\fj_{\fU\cap\fV}^*\cE^\bu\lrarrow0
\end{equation}
in the abelian category of complexes in $\fX\Tors$.
 Here the morphisms in
the sequence~\eqref{qcoh-tors-mayer-vietoris-short-seq-formula} come
from the adjunction~\eqref{qcoh-tors-inverse-direct-adjunction} of
the direct and inverse image functors, together with the isomorphisms
of compatibility of the constructions of the direct and inverse images
with compositions of morphisms of formal schemes.
 The sequence~\eqref{qcoh-tors-mayer-vietoris-short-seq-formula} is
exact, since its restrictions to $\fU$ and $\fV$ are naturally split
short exact sequences of complexes in $\fU\Tors$ and $\fV\Tors$.

 It remains to point out that, following the discussion in
Section~\ref{derived-direct-images-of-qcoh-tors-subsecn}, one can
compute the derived direct image functors~$\boR\ff_*$ for unbounded
complexes of quasi-coherent torsion sheaves by applying the direct
image functors~$\ff_*$ to arbitrary complexes of flasque quasi-coherent
torsion sheaves.
 Representing the given derived category object $\cM^\bu$ by a complex
$\cE^\bu$ of flasque quasi-coherent torsion sheaves on $\fX$, one
obtains a short exact sequence of
complexes~\eqref{qcoh-tors-mayer-vietoris-short-seq-formula},
which induces the desired distinguished
triangle~\eqref{qcoh-tors-mayer-vietoris-triangle-formula}
in $\sD(\fX\Tors)$.
 It is important for this argument that the restriction of
quasi-coherent torsion sheaves to open subschemes preserves flasqueness.
\end{proof}

\begin{cor} \label{negative-Ext-vanishing-preserved-by-gluing-cor}
 Let\/ $\fX$ be a locally Noetherian formal scheme covered by
two open formal subschemes, $\fX=\fU\cup\fV$.
 Let\/ $\cM^\bu$ and\/ $\cN^\bu\in\sD(\fX\Tors)$ be two objects of
the derived category of quasi-coherent torsion sheaves on\/~$\fX$.
 Assume that
\begin{itemize}
\item $\Hom_{\sD(\fU\Tors)}(\fj_\fU^*\cN^\bu,\fj_\fU^*\cM^\bu[i])=0$
 for all $i<0$,
\item $\Hom_{\sD(\fV\Tors)}(\fj_\fV^*\cN^\bu,\fj_\fV^*\cM^\bu[i])=0$
 for all $i<0$, and
\item $\Hom_{\sD((\fU\cap\fV)\Tors)}(\fj_{\fU\cap\fV}^*\cN^\bu,\,
 \fj_{\fU\cap\fV}^*\cM^\bu[i])=0$ for all $i<0$.
\end{itemize}
 Then\/ $\Hom_{\sD(\fX\Tors)}(\cN^\bu,\cM^\bu[i])=0$ for all $i<0$.
\end{cor}

\begin{proof}
 Applying the functor $\Hom_{\sD(\fX\Tors)}(\cN^\bu,{-}[*])$ to
the distinguished
triangle~\eqref{qcoh-tors-mayer-vietoris-triangle-formula}, we obtain
a long exact sequence of abelian groups {\hbadness=1375
\begin{multline} \label{mayer-vietoris-Hom-long-exact-sequence-I}
 \dotsb\lrarrow\Hom_{\sD(\fX\Tors)}(\cN^\bu,
 \boR\fj_{\fU\cap\fV}{}_*\fj_{\fU\cap\fV}^*\cM^\bu[i-1])\lrarrow
 \Hom_{\sD(\fX\Tors)}(\cN^\bu,\cM^\bu[i]) \\ \lrarrow
 \Hom_{\sD(\fX\Tors)}(\cN^\bu,\boR\fj_\fU{}_*\fj_\fU^*\cM^\bu[i])\oplus
 \Hom_{\sD(\fX\Tors)}(\cN^\bu,\boR\fj_\fV{}_*\fj_\fV^*\cM^\bu[i])
 \lrarrow\dotsb
\end{multline}
 It} remains to recall the natural isomorphism of abelian groups
\begin{equation} \label{qcoh-tors-derived-open-immersion-adjunction}
 \Hom_{\sD(\fX\Tors)}(\cN^\bu,\boR\fj_*\cL^\bu)\simeq
 \Hom_{\sD(\fW\Tors)}(\fj^*\cN^\bu,\cL^\bu),
\end{equation}
which holds for any open formal subscheme $\fW\subset\fX$ with
the open immersion morphism $\fj\:\fW\rarrow\fX$ and any object
$\cL^\bu\in\sD(\fW\Tors)$ according to the adjunction of the functors
$\boR\fj_*$ \,\eqref{qcoh-tors-derived-direct-image-formula} and
$\fj^*$ \,\eqref{qcoh-derived-cats-underived-inverse-image}
as per Section~\ref{triangulated-inverse-images-subsecn}.
 Rewriting
the sequence~\eqref{mayer-vietoris-Hom-long-exact-sequence-I} using
the isomorphisms~\eqref{qcoh-tors-derived-open-immersion-adjunction},
we obtain a long exact sequence of abelian groups
\begin{multline} \label{mayer-vietoris-Hom-long-exact-sequence-II}
 \dotsb\lrarrow\Hom_{\sD((\fU\cap\fV)\Tors)}(\fj_{\fU\cap\fV}^*\cN^\bu,
 \fj_{\fU\cap\fV}^*\cM^\bu[i-1])\lrarrow
 \Hom_{\sD(\fX\Tors)}(\cN^\bu,\cM^\bu[i]) \\ \lrarrow
 \Hom_{\sD(\fU\Tors)}(\fj_\fU^*\cN^\bu,\fj_\fU^*\cM^\bu[i])\oplus
 \Hom_{\sD(\fV\Tors)}(\fj_\fV^*\cN^\bu,\fj_\fV^*\cM^\bu[i]) \\
 \lrarrow\Hom_{\sD((\fU\cap\fV)\Tors)}(\fj_{\fU\cap\fV}^*\cN^\bu,
 \fj_{\fU\cap\fV}^*\cM^\bu[i])\lrarrow\dotsb,
\end{multline}
and the desired assertion immediately follows.
\end{proof}

\begin{cor} \label{uniqueness-of-morphism-of-glued-objects-cor}
 Let\/ $\fX$ be a locally Noetherian formal scheme covered by
two open formal subschemes, $\fX=\fU\cup\fV$.
 Let $f\:\cN^\bu\rarrow\cM^\bu$ be a morphism in the derived
category of quasi-coherent torsion sheaves on\/~$\fX$.
 Assume that\/ $\fj_\fU^*(f)=0=\fj_\fV^*(f)$, and further that
\begin{itemize}
\item $\Hom_{\sD((\fU\cap\fV)\Tors)}(\fj_{\fU\cap\fV}^*\cN^\bu,\,
 \fj_{\fU\cap\fV}^*\cM^\bu[-1])=0$.
\end{itemize}
 Then $f=0$.
\end{cor}

\begin{proof}
 Consider the long exact sequence of abelian
groups~\eqref{mayer-vietoris-Hom-long-exact-sequence-II} for $i=0$.
 By assumption, the image of $f\in\Hom_{\sD(\fX\Tors)}(\cM^\bu,\cN^\bu)$
under the differential in this long exact sequence vanishes.
 On the other hand, the differential mapping into
$\Hom_{\sD(\fX\Tors)}(\cM^\bu,\cN^\bu)$ comes from a zero group,
by the other assumption.
\end{proof}

\begin{cor} \label{existence-of-morphism-of-glued-objects-cor}
 Let\/ $\fX$ be a locally Noetherian formal scheme covered by
two open formal subschemes, $\fX=\fU\cup\fV$.
 So we have the identity open immersion morphisms\/
$\fk_\fU\:\fU\cap\nobreak\fV\rarrow\fU$ and\/
$\fk_\fV\:\fU\cap\nobreak\fV\rarrow\fV$.
 Let\/ $\cM^\bu$ and\/ $\cN^\bu\in\sD(\fX\Tors)$ be two objects of
the derived category of quasi-coherent torsion sheaves on\/~$\fX$.
 Let $g\:\fj_\fU^*\cN^\bu\rarrow\fj_\fU^*\cM^\bu$ be a morphism in\/
$\sD(\fU\Tors)$ and $h\:\fj_\fV^*\cN^\bu\rarrow\fj_\fV^*\cM^\bu$ be
a morphism in\/ $\sD(\fV\Tors)$.
 Assume that\/ $\fk_\fU^*(g)=\fk_\fV^*(h)$ in\/
$\sD((\fU\cap\nobreak\fV)\Tors)$.
 Then there exists a morphism $f\:\cN^\bu\rarrow\cM^\bu$ in\/
$\sD(\fX\Tors)$ such that\/ $\fj_\fU^*(f)=g$ and\/ $\fj_\fV^*(f)=h$.
\end{cor}

\begin{proof}
 This is our version of the second half of~\cite[Lemma Tag~08DG]{SP}.
 Once again, we consider the long exact sequence of abelian
groups~\eqref{mayer-vietoris-Hom-long-exact-sequence-II} for $i=0$.
 By assumption, the image of the element $(g,h)\in
\Hom_{\sD(\fU\Tors)}(\fj_\fU^*\cN^\bu,\fj_\fU^*\cM^\bu)\oplus
\Hom_{\sD(\fV\Tors)}(\fj_\fV^*\cN^\bu,\fj_\fV^*\cM^\bu)$ under
the differential in this long exact sequence vanishes.
 Hence there exists an element
$f\in\Hom_{\sD(\fX\Tors)}(\cN^\bu,\cM^\bu)$ such that $(g,h)$ is
the image of~$f$ under the differnential in the long exact sequence.
\end{proof}

\begin{lem} \label{gluing-from-covering-by-two-opens-lemma}
 Let\/ $\fX$ be a locally Noetherian formal scheme covered by
two open formal subschemes, $\fX=\fU\cup\fV$.
 Let\/ $\cA^\bu\in\sD(\fU\Tors)$ and\/ $\cB^\bu\in\sD(\fV\Tors)$ be two
derived category objects.
 Suppose given an isomorphism\/ $\sigma\:\fk_\fU^*\cA^\bu\simeq
\fk_\fV^*\cB^\bu$ in\/ $\sD((\fU\cap\nobreak\fV)\Tors)$.
 Then there exists an object\/ $\cM^\bu\in\sD(\fX\Tors)$ together with
isomorphisms\/ $\rho^\fX_\fU\:\fj_\fU^*\cM^\bu\simeq\cA^\bu$ in
$\sD(\fU\Tors)$ and\/ $\rho^\fX_\fV\:\fj_\fV^*\cM^\bu\simeq\cB^\bu$
in $\sD(\fV\Tors)$ such that\/
$\fk_\fV^*(\rho^\fX_\fV)=\sigma\circ\fk_\fU^*(\rho^\fX_\fU)$
in\/ $\sD((\fU\cap\nobreak\fV)\Tors)$.
\end{lem}

\begin{proof}
 This is our version of the first half of~\cite[Lemma Tag~08DG]{SP}.
 Identifying $\fk_\fU^*\cA^\bu$ with $\fk_\fV^*\cB^\bu$ by means
of~$\sigma$, put $\fk_\fU^*\cA^\bu=\cC^\bu=\fk_\fV^*\cB^\bu\in
\sD((\fU\cap\fV)\Tors)$.
 The adjunction
isomorphism~\eqref{qcoh-tors-derived-open-immersion-adjunction}
provides natural morphisms $\cA^\bu\rarrow
\boR\fk_\fU{}_*\fk_\fU^*\cA^\bu=\boR\fk_\fU{}_*\cC^\bu$ in
$\sD(\fU\Tors)$ and $\cB^\bu\rarrow\boR\fk_\fV{}_*\fk_\fV^*\cB^\bu=
\boR\fk_\fV{}_*\cC^\bu$ in $\sD(\fV\Tors)$.
 Applying the functor $\boR\fj_\fU{}_*$ to the former morphism and
the functor $\boR\fj_\fV{}_*$ to the latter one, and combining
them together, we obtain a morphism
\begin{equation} \label{gluing-morphism-for-taking-cocone}
 \boR\fj_\fU{}_*\cA^\bu\oplus\boR\fj_\fV{}_*\cB^\bu\lrarrow
 \boR\fj_{\fU\cap\fV}{}_*\cC^\bu
\end{equation}
in $\sD(\fX\Tors)$.
 Let us extend the morphism~\eqref{gluing-morphism-for-taking-cocone}
to a distinguished triangle in $\sD(\fX\Tors)$,
\begin{equation} \label{gluing-cocone-triangle}
 \cM^\bu\lrarrow\boR\fj_\fU{}_*\cA^\bu\oplus
 \boR\fj_\fV{}_*\cB^\bu\lrarrow\boR\fj_{\fU\cap\fV}{}_*\cC^\bu
 \lrarrow\cM^\bu[1]
\end{equation}
 Applying the triangulated functors $\fj_\fU^*$ and~$\fj_\fV^*$ to
the morphism~\eqref{gluing-morphism-for-taking-cocone} produces
split epimorphisms with natural splittings.
 Hence, applying the same triangulated functors to
the triangle~\eqref{gluing-cocone-triangle}, one obtains the desired
isomorphisms $\rho^\fX_\fU\:\fj_\fU^*\cM^\bu\simeq\cA^\bu$ and
$\rho^\fX_\fV\:\fj_\fV^*\cM^\bu\simeq\cB^\bu$.
 Finally, applying the triangulated functor $\fj_{\fU\cap\fV}^*$
to the triangle~\eqref{gluing-cocone-triangle}, the vanishing of
the composition of the leftmost and the middle arrows in the resulting
split triangle translates into the desired compatibility equation
on the isomorphisms $\fk_\fV^*(\rho^\fX_\fV)=
\sigma\circ\fk_\fU^*(\rho^\fX_\fU)$.
\end{proof}

 In the next four propositions, we consider the following setting.
 Let $\fX$ be a locally Noetherian formal scheme and $\bC$ be
an open covering of~$\fX$.
 So $\bC=\{\fU\}$ is a set of some open formal subschemes
$\fU\subset\fX$ covering~$\fX$.
 We will assume the following two conditions to be satisfied:
\begin{enumerate}
\renewcommand{\theenumi}{\roman{enumi}}
\item there exists a finite open covering
$\fX=\bigcup_{\alpha=1}^n\fU_\alpha$ of the formal scheme $\fX$
by open formal subschemes $\fU_\alpha\in\bC$;
\item one has $\fU\cap\fV\in\bC$ whenever $\fU\in\bC$ and $\fV\in\bC$.
\end{enumerate}

 Furthermore, in the present section, given an open formal subscheme
$\fY\subset\fX$, we denote by $\bC|^\fY$ the set of all open formal
subschemes $\{\fV\mid \fV\in\bC,\,\fV\subset\fY\}$.
 (Notice that the meaning of this notation is \emph{different}
from that of the notation $\bW|_\fY$ from
Section~\ref{inverse-images-of-ctrh-contramod-subsecn} and
the subsequent exposition.)
 Generally speaking, $\bC|^\fY$ need not even be an open covering
of~$\fY$.
 However, assume that there exists a finite set of open subsets
$\fV_\beta\in\bC$, \, $1\le\beta\le m$, such that
$\fY=\bigcup_{\beta=1}^m\fV_\beta$.
 Then $\bC|^\fY$ is an open covering of $\fY$ satisfying
conditions~(i) and~(ii).

\begin{prop} \label{gluing-negative-Ext-vanishing-prop}
 Let\/ $\fX$ be a locally Noetherian formal scheme, and let\/
$\bC$ be an open covering of\/ $\fX$ satisfying conditions~(i)
and~(ii) above.
 Let\/ $\cM^\bu$ and\/ $\cN^\bu\in\sD(\fX\Tors)$ be two objects of
the derived category of quasi-coherent torsion sheaves on\/~$\fX$.
 Assume that
\begin{itemize}
\item $\Hom_{\sD(\fU\Tors)}(\fj_\fU^*\cN^\bu,\fj_\fU^*\cM^\bu[i])=0$
 for all\/ $\fU\in\bC$ and all $i<0$.
\end{itemize}
 Then\/ $\Hom_{\sD(\fX\Tors)}(\cN^\bu,\cM^\bu[i])=0$ for all $i<0$.
\end{prop}

\begin{proof}
 This is a simple argument by induction on the cardinality~$n$ of
the finite open covering $\fX=\bigcup_{\alpha=1}^n\fU_\alpha$ from
condition~(i).
 For $n=1$, there is nothing to prove.
 For $n\ge2$, put $\fY=\bigcup_{\beta=1}^{n-1}\fU_\beta\subset\fX$
and $\fU=\fU_n$.
 Applying the induction assumption to the formal scheme $\fY$ with
its open covering $\bC|^\fY$, we conclude that
$$
 \Hom_{\sD(\fY\Tors)}(\fj_\fY^*\cN^\bu,\fj_\fY^*\cM^\bu[i])=0
 \qquad\text{for all $i<0$}.
$$
 Similarly, applying the induction assumption to the formal scheme
$\fY\cap\fU$ with its open covering $\bC|^{\fY\cap\fU}$ and its finite
subcovering $\fY\cap\fU=\bigcup_{\beta=1}^{n-1}(\fU_\beta\cap\fU)$,
we see that
$$
 \Hom_{\sD((\fY\cap\fU)\Tors)}
 (\fj_{\fY\cap\fU}^*\cN^\bu,\fj_{\fY\cap\fU}^*\cM^\bu[i])=0
 \qquad\text{for all $i<0$}.
$$
 Now it remains to apply
Corollary~\ref{negative-Ext-vanishing-preserved-by-gluing-cor}
to the formal scheme $\fX$ with its open covering $\fX=\fY\cup\fU$.
\end{proof}

\begin{prop} \label{gluing-morphism-uniqueness-prop}
 Let\/ $\fX$ be a locally Noetherian formal scheme, and let\/
$\bC$ be an open covering of\/ $\fX$ satisfying conditions~(i)
and~(ii) above.
 Let\/ $f\:\cN^\bu\rarrow\cM^\bu$ be a morphism in the derived
category\/ $\sD(\fX\Tors)$ of quasi-coherent torsion sheaves on\/~$\fX$.
 Assume that\/ $\fj_\fU^*(f)=0$ for all\/ $\fU\in\bC$, and further that
\begin{itemize}
\item $\Hom_{\sD(\fU\Tors)}(\fj_\fU^*\cN^\bu,\fj_\fU^*\cM^\bu[i])=0$
 for all\/ $\fU\in\bC$ and all $i<0$.
\end{itemize}
 Then $f=0$.
\end{prop}

\begin{proof}
 The proof uses induction on the cardinality~$n$ of the finite open
covering from condition~(i) similar to the one in the proof
of Proposition~\ref{gluing-negative-Ext-vanishing-prop}.
 As in that proof, for $n=1$ the assertion is obvious; and
for $n\ge2$, we put $\fY=\bigcup_{\beta=1}^{n-1}\fU_\beta\subset\fX$
and $\fU=\fU_n$.
 Applying the induction assumption to the formal scheme $\fY$ with
its open covering $\bC|^\fY$ and the morphism $\fj_\fY^*(f)$ in
$\sD(\fY\Tors)$, we see that $\fj_\fY^*(f)=0$.
 By assumption, we also have $\fj_\fU^*(f)=0$.
 Furthermore, $\Hom_{\sD((\fY\cap\fU)\Tors)}(\fj_{\fY\cap\fU}^*\cN^\bu,
\>\fj_{\fY\cap\fU}^*\cM^\bu[-1])=0$ by
Proposition~\ref{gluing-negative-Ext-vanishing-prop} applied to
the formal scheme $\fY\cap\fU$ applied to the formal scheme
$\fY\cap\fU$ with its open covering $\bC|^{\fY\cap\fU}$ and its finite
subcovering $\fY\cap\fU=\bigcup_{\beta=1}^{n-1}
(\fU_\beta\cap\nobreak\fU)$.
 Appling Corollary~\ref{uniqueness-of-morphism-of-glued-objects-cor}
to the formal scheme $\fX$ with its open covering $\fX=\fY\cup\fU$,
we conclude that $f=0$.
\end{proof}

 Given two open formal subschemes $\fU$ and $\fV\subset\fX$ such that
$\fV\subset\fU$, we denote by $\fh_{\fV,\fU}\:\fV\rarrow\fU$
the identity open immersion morphism.

\begin{prop} \label{gluing-morphism-existence-prop}
 Let\/ $\fX$ be a locally Noetherian formal scheme, and let\/
$\bC$ be an open covering of\/ $\fX$ satisfying conditions~(i)
and~(ii) above.
 Let\/ $\cM^\bu$ and\/ $\cN^\bu\in\sD(\fX\Tors)$ be two objects of
the derived category of quasi-coherent torsion sheaves on\/~$\fX$.
 Suppose that we are given
\begin{itemize}
\item a morphism $f_\fU\:\fj_\fU^*\cN^\bu\rarrow\fj_\fU^*\cM^\bu$
in\/ $\sD(\fU\Tors)$ for every\/ $\fU\in\bC$ such that
\item $\fh_{\fV,\fU}^*(f_\fU)=f_\fV$ for all\/ $\fU$, $\fV\in\bC$
with\/ $\fV\subset\fU$.
\end{itemize}
 Assume further that
\begin{itemize}
\item $\Hom_{\sD(\fU\Tors)}(\fj_\fU^*\cN^\bu,\fj_\fU^*\cM^\bu[i])=0$
 for all\/ $\fU\in\bC$ and all $i<0$.
\end{itemize}
 Then there exists a unique morphism $f\:\cN^\bu\rarrow\cM^\bu$ in\/
$\sD(\fX\Tors)$ such that
\begin{itemize}
\item $\fj_\fU^*(f)=f_\fU$ for all\/ $\fU\in\bC$.
\end{itemize}
\end{prop}

\begin{proof}
 The uniqueness assertion is provided by
Proposition~\ref{gluing-morphism-uniqueness-prop}.
 To prove the existence, we once again argue by induction on
the cardinality~$n$ of the finite open covering from condition~(i),
similarly to the proofs of
Propositions~\ref{gluing-negative-Ext-vanishing-prop}
and~\ref{gluing-morphism-uniqueness-prop}.
 As in the previous two proofs, for $n=1$ the assertion is obvious; and
for $n\ge2$, we put $\fY=\bigcup_{\beta=1}^{n-1}\fU_\beta\subset\fX$
and $\fU=\fU_n$.

 Applying the induction assumption to the formal scheme $\fY$ with
its open covering $\bC|^\fY$ and the objects $\fj_\fY^*\cM^\bu$
and $\fj_\fY^*\cN^\bu\in\sD(\fY\Tors)$, we obtain a unique morphism
$g\:\fj_\fY^*\cN^\bu\rarrow\fj_\fY^*\cM^\bu$ in $\sD(\fY\Tors)$ such
that $\fh_{\fV,\fY}^*(g)=f_\fV$ for all $\fV\in\bC$, \,$\fV\subset\fY$.
 Now the morphism $l'=\fh_{\fY\cap\fU,\fY}^*(g)$ in
$\sD((\fY\cap\fU)\Tors)$ has the property that
$\fh_{\fV,\fY\cap\fU}^*(l')=f_\fV$ for all $\fV\in\bC$,
\,$\fV\subset\fY\cap\fU$.
 The morphism $l''=\fh_{\fY\cap\fU,\fU}^*(f_\fU)$ in
in $\sD((\fY\cap\fU)\Tors)$ has the same property,
$\fh_{\fV,\fY\cap\fU}^*(l'')=f_\fV$ for all $\fV\in\bC$,
\,$\fV\subset\fY\cap\fU$.
 By Proposition~\ref{gluing-morphism-uniqueness-prop} applied to
the formal scheme $\fY\cap\fU$ with its open covering
$\bC|^{\fY\cap\fU}$ and its finite subcovering
$\fY\cap\fU=\bigcup_{\beta=1}^{n-1}(\fU_\beta\cap\fU)$, such a morphism
$l\:\fj^*_{\fY\cap\fU}(\cN^\bu)\rarrow\fj^*_{\fY\cap\fU}(\cN^\bu)$ in
$\sD((\fY\cap\fU)\Tors)$ is unique.
 Hence we have $l'=l''$, i.~e., $\fh_{\fY\cap\fU,\fY}^*(g)=
\fh_{\fY\cap\fU,\fU}^*(f_\fU)$.

 By Corollary~\ref{existence-of-morphism-of-glued-objects-cor} applied
to the formal scheme $\fX$ with its open covering $\fX=\fY\cup\fU$,
there exists a morphism $f\:\cN^\bu\rarrow\cM^\bu$ in $\sD(\fX\Tors)$
such that $\fj_\fY^*(f)=g$ and $\fj_\fU^*(f)=f_\fU$.
 It remains to show that $\fj_\fW^*(f)=f_\fW$ for all $\fW\in\bC$.
 By construction, we know that this equation holds whenever
$\fW\subset\fY$ or $\fW\subset\fU$, but this is not yet enough
to establish the properties of~$f$ claimed in the proposition.
 In the general case, the desired equation is provable by applying
Proposition~\ref{gluing-morphism-uniqueness-prop} to the formal
scheme $\fW$ with its covering $\bC'$ consisting of all the open
formal subschemes $\fV\in\bC$ such that either $\fV\subset\fW\cap\fY$
or $\fV\subset\fW\cap\fU$.
 It is important that the open covering $\bC'$ of the formal scheme
$\fW$ has a finite subcovering
$\fW=\bigcup_{\alpha=1}^n(\fW\cap\fU_\alpha)$.

 Alternatively, one can apply
Corollary~\ref{uniqueness-of-morphism-of-glued-objects-cor} to
the formal scheme $\fW$ with its open covering $\fW=(\fW\cap\fY)
\cup(\fW\cap\fU)$.
 Notice that $\Hom_{\sD((\fW\cap\fY\cap\fU)\Tors)}
(\fj_{\fW\cap\fY\cap\fU}^*\cN^\bu,\>\allowbreak
\fj_{\fW\cap\fY\cap\fU}^*\cM^\bu[-1])=0$
by Proposition~\ref{gluing-negative-Ext-vanishing-prop} applied to
the formal scheme $\fW\cap\fY\cap\fU$ with its open covering
$\bC|^{\fW\cap\fY\cap\fU}$ and its finite subcovering
$\fW\cap\fY\cap\fU=\bigcup_{\beta=1}^{n-1}(\fW\cap\fU_\beta\cap\fU)$.
 So Corollary~\ref{uniqueness-of-morphism-of-glued-objects-cor}
is indeed applicable.
\end{proof}

\begin{prop} \label{main-gluing-prop}
 Let\/ $\fX$ be a locally Noetherian formal scheme, and let\/
$\bC$ be an open covering of\/ $\fX$ satisfying conditions~(i)
and~(ii) above.
 Suppose given
\begin{itemize}
\item a derived category object\/ $\cM_\fU^\bu\in\sD(\fU\Tors)$
for every open formal subscheme\/ $\fU\in\bC$, and
\item an isomorphism\/ $\rho^\fU_\fV\:\fh_{\fV,\fU}^*\cM_\fU^\bu\simeq
\cM_\fV^\bu$ in\/ $\sD(\fV\Tors)$ for every pair of open formal
subschemes\/ $\fU$, $\fV\in\bC$ such that\/ $\fV\subset\fU$.
\end{itemize}
 Assume that 
\begin{itemize}
\item the equation
$\rho^\fU_\fW=\rho^\fV_\fW\circ\fh_{\fW,\fV}^*(\rho^\fU_\fV)$
holds for all triples of open formal subschemes\/ $\fU$, $\fV$,
$\fW\in\bC$ such that\/ $\fW\subset\fV\subset\fU$.
\item $\Hom_{\sD(\fU\Tors)}(\cM_\fU^\bu,\cM_\fU^\bu[i])=0$ for all\/
$\fU\in\bC$ and all $i<0$.
\end{itemize}
 Then there exists a unique, up to a unique isomorphism,
derived category object\/ $\cM^\bu\in\sD(\fX\Tors)$ together with
\begin{itemize}
\item isomorphisms\/ $\rho_\fU^\fX:\fj_\fU^*\cM^\bu\simeq\cM_\fU^\bu$
in\/ $\sD(\fU\Tors)$, given for all\/ $\fU\in\bC$, for which
\item the equation
$\rho^\fX_\fV=\rho^\fU_\fV\circ\fh_{\fV,\fU}^*(\rho^\fX_\fU)$
holds for all pairs of open formal subschemes\/ $\fU$, $\fV\in\bC$
such that\/ $\fV\subset\fU$.
\end{itemize}
 In addition, one has
\begin{itemize}
\item $\Hom_{\sD(\fX\Tors)}(\cM^\bu,\cM^\bu[i])=0$ for all $i<0$.
\end{itemize}
\end{prop}

\begin{proof}
 The uniqueness assertion is provided by
Proposition~\ref{gluing-morphism-existence-prop}.
 To prove the existence, we again argue by induction on
the cardinality~$n$ of the finite open covering from condition~(i).
 As in the previous three proofs, for $n=1$ the assertion is obvious;
and for $n\ge2$, we put $\fY=\bigcup_{\beta=1}^{n-1}\fU_\beta
\subset\fX$ and $\fU=\fU_n$.

 Applying the induction assumption to the formal scheme $\fY$ with its
open covering $\bC|^\fY$ and its finite subcovering
$\fY=\bigcup_{\beta=1}^{n-1}\fU_\beta$, we obtain a derived category
object $\cM_\fY^\bu\in\sD(\fY\Tors)$ together with isomorphisms
$\rho^\fY_\fV\:\fh_{\fV,\fY}^*\cM_\fY^\bu\simeq\cM_\fV^\bu$ in
$\sD(\fV\Tors)$ for all open formal subschemes $\fV\in\bC$,
\,$\fV\subset\fY$.
 Moreover, as a part of the induction assumption, the equations
$\rho^\fY_\fW=\rho^\fV_\fW\circ\fh_{\fW,\fV}^*(\rho^\fY_\fV)$ hold for
all pairs of open formal subschemes\/ $\fV$, $\fW\in\bC$ such that\/
$\fW\subset\fV\subset\fY$.

 Now the object ${}'\cM^\bu_{\fY\cap\fU}=
\fh_{\fY\cap\fU,\fY}^*\cM_\fY^\bu\in\sD((\fY\cap\fU)\Tors)$ is endowed
with naturally arising isomorphisms ${}'\!\rho^{\fY\cap\fU}_\fV\:
\fh_{\fV,\fY\cap\fU}^*{}'\cM^\bu_{\fY\cap\fU}\simeq\cM_\fV^\bu$
for all $\fV\in\bC$, \,$\fV\subset\fY\cap\fU$, satisfying
the compatibility equations as above.
 The same applies to the object ${}''\cM^\bu_{\fY\cap\fU}=
\fh_{\fY\cap\fU,\fU}^*\cM_\fU^\bu\in\sD((\fY\cap\fU)\Tors)$,
which is endowed with naturally arising isomorphisms
${}''\!\rho^{\fY\cap\fU}_\fV\:\fh_{\fV,\fY\cap\fU}^*
{}''\cM^\bu_{\fY\cap\fU}\simeq\cM_\fV^\bu$ for all $\fV\in\bC$,
\,$\fV\subset\fY\cap\fU$, satisfying the same compatibility equations.
 Hence Proposition~\ref{gluing-morphism-existence-prop} provides
a unique isomorphism $\sigma\:\fh_{\fY\cap\fU,\fY}^*\cM_\fY^\bu
\simeq \fh_{\fY\cap\fU,\fU}^*\cM_\fU^\bu$ such that
${}''\!\rho^{\fY\cap\fU}_\fV\circ\sigma={}'\!\rho^{\fY\cap\fU}_\fV$
for all for all $\fV\in\bC$, \,$\fV\subset\fY\cap\fU$.

 Applying Lemma~\ref{gluing-from-covering-by-two-opens-lemma}, we
obtain the desired derived category object $\cM^\bu\in\sD(\fX\Tors)$
together with isomorphisms $\fj_\fY^*\cM^\bu\simeq\cM_\fY^\bu$
in $\sD(\fY\Tors)$ and $\fj_\fU^*\cM^\bu\simeq\cM_\fU^\bu$ in
$\sD(\fU\Tors)$ satisfying the compatibility equation with~$\sigma$.
 This provides the isomorphisms $\rho_\fW^\fX:\fj_\fW^*\cM^\bu\simeq
\cM_\fW^\bu$ for all open formal subschemes $\fW\in\bC$ such that
either $\fW\subset\fY$ or $\fW\subset\fU$.
 In the general case of $\fW\in\bC$, the desired isomorphim
$\rho_\fW^\fX$ is constructed by applying
Proposition~\ref{gluing-morphism-existence-prop} to the formal
scheme $\fW$ with its covering $\bC'$ consisting of all the open
formal subschemes $\fV\in\bC$ such that either $\fV\subset\fW\cap\fY$
or $\fV\subset\fW\cap\fU$.
 This is a one step higher analogue of the argument in the proof of
Proposition~\ref{gluing-morphism-existence-prop}.
 Alternatively, one can apply
Corollary~\ref{existence-of-morphism-of-glued-objects-cor} to
the formal scheme $\fW$ with its open covering $\fW=(\fW\cap\fY)
\cup(\fW\cap\fU)$.
\end{proof}

\begin{thm} \label{main-gluing-theorem}
 Let\/ $\fX$ be a Noetherian formal scheme with an open
covering\/~$\bW$.
 Suppose given
\begin{itemize}
\item a derived category object\/ $\cM_\fU^\bu\in\sD(\fU\Tors)$
for every affine open formal subscheme\/ $\fU$ subordinate
to\/ $\bW$ in $\fX$, and
\item an isomorphism\/ $\rho^\fU_\fV\:\fh_{\fV,\fU}^*\cM_\fU^\bu\simeq
\cM_\fV^\bu$ in\/ $\sD(\fV\Tors)$ for every pair of affine open formal
subschemes\/ $\fV\subset\fU$ subordinate to\/ $\bW$ in\/~$\fX$.
\end{itemize}
 Assume that 
\begin{itemize}
\item the equation
$\rho^\fU_\fW=\rho^\fV_\fW\circ\fh_{\fW,\fV}^*(\rho^\fU_\fV)$
holds for all triples of affine open formal subschemes\/ $\fW\subset\fV
\subset\fU$ subordinate to\/ $\bW$ in\/~$\fX$;
\item $\Hom_{\sD(\fU\Tors)}(\cM_\fU^\bu,\cM_\fU^\bu[i])=0$ for all\/
affine open formal subschemes\/ $\fU\subset\fX$ subordinate to\/ $\bW$
and all $i<0$.
\end{itemize}
 Then there exists a unique, up to a unique isomorphism,
derived category object\/ $\cM^\bu\in\sD(\fX\Tors)$ together with
\begin{itemize}
\item isomorphisms\/ $\rho_\fU^\fX:\fj_\fU^*\cM^\bu\simeq\cM_\fU^\bu$
in\/ $\sD(\fU\Tors)$, given for all affine open formal subschemes\/
$\fU\subset\fX$ subordinate to\/ $\bW$, for which
\item the equation
$\rho^\fX_\fV=\rho^\fU_\fV\circ\fh_{\fV,\fU}^*(\rho^\fX_\fU)$
holds for all pairs of affine open formal subschemes\/ $\fV\subset\fU$
subordinate to\/ $\bW$ in\/~$\fX$.
\end{itemize}
 In addition, one has
\begin{itemize}
\item $\Hom_{\sD(\fX\Tors)}(\cM^\bu,\cM^\bu[i])=0$ for all $i<0$.
\end{itemize}
\end{thm}

\begin{proof}
 The proof is based on Proposition~\ref{main-gluing-prop} and
proceeds in two steps.
 Firstly, let us assume that the formal scheme $\fX$ is semi-separated.
 Then we simply apply Proposition~\ref{main-gluing-prop} to
the open covering $\bC$ of $\fX$ consisting of all the affine open
formal subschemes $\fU\subset\fX$ subordinate to~$\bW$.
 We have $\fU\cap\fV\in\bC$ for all $\fU$, $\fV\in\bC$, since $\fX$
is semi-separated; and the open covering $\bC$ of $\fX$ has a finite
subcovering, since $\fX$ is quasi-compact by the theorem's assumption.
 This proves both the existence and the uniqueness of the desired
object $\cM^\bu\in\sD(\fX\Tors)$.

 In the general case of a non-semi-separated Noetherian formal
scheme $\fX$, we consider the open covering $\bC$ of $\fX$ consisting
of all the semi-separated open formal subschemes $\fY\subset\fX$.
 The intersection of any two semi-separated open formal subschemes
in $\fX$ is obviously semi-separated.
 The previous paragraph, applied to $\fY$, provides a derived category
object $\cM_\fY^\bu\in\sD(\fY\Tors)$.
 Given a pair of semi-separated open formal subschemes $\fT\subset\fY$
in $\fX$, the uniqueness assertion in the previous paragraph provides
a uniquely defined isomorphism $\rho^\fY_\fT\:\fh_{\fT,\fY}^*
\cM_\fY^\bu\simeq\cM_\fT^\bu$.
 The equations $\rho^\fY_\fW=\rho^\fT_\fW\circ\fh_{\fW,\fT}^*
(\rho^\fY_\fT)$ for all triples of semi-separated open formal
subschemes\/ $\fW\subset\fT\subset\fY$ in $\fX$ can be checked using
Proposition~\ref{gluing-morphism-uniqueness-prop}.
 So we can apply Proposition~\ref{main-gluing-prop} again, proving
the theorem in the general case.
\end{proof}

\subsection{Dedualizing complexes on Noetherian formal schemes}
\label{dedualizing-complexes-on-formal-schemes-subsecn}
 The following corollary is obtained by combining the final lemmas from
Section~\ref{dedualizing-complexes-over-rings-subsecn} with the gluing
theorem from Section~\ref{gluing-in-D-Tors-subsecn}.

 For every affine Noetherian formal scheme $\fU=\Spf\fR$, we denote
by $\cB_\fU^\bu\in\sD^\bb(\fU\Tors)$ the derived category object
corresponding to the dedualizing complex $B_\fR^\bu$ for the adic
Noetherian ring $\fR$ constructed in
Section~\ref{dedualizing-complexes-over-rings-subsecn}.
 Given an open immersion of affine Noetherian formal schemes
$\fh\:\fV\rarrow\fU$, where $\fU=\Spf\fR$ and $\fV=\Spf\fS$,
consider the related formal open immersion morphism of adic
Noetherian rings $h\:\fR\rarrow\fS$ (in the terminology of
Section~\ref{prelim-formal-open-immersions-coverings-subsecn}).
 Denote by $\rho^\fU_\fV\:\fh^*(\cB_\fU^\bu)\simeq\cB_\fV^\bu$
the isomorphism in $\sD^\bb(\fV\Tors)$ corresponding to
the isomorphism $\rho^\fR_\fS\:h^*(B_\fR^\bu)\simeq B_\fS^\bu$
in $\sD^\bb(\fS\Modl_\tors)$ constructed in
Lemma~\ref{dedualizing-complex-extension-of-scalars}.

 Keeping the notation of Section~\ref{gluing-in-D-Tors-subsecn},
we consider a Noetherian formal scheme $\fX$, and for any open
formal subscheme $\fU\subset\fX$, denote by $\fj_\fU\:\fU\rarrow\fX$
the identity open immersion morphism.
 Given two open formal subschemes $\fV\subset\fU$ in $\fX$, we
denote their identity open immersion morphism by $\fh_{\fV,\fU}\:
\fV\rarrow\fU$.

\begin{cor} \label{nonaffine-dedualizing-complex-cor}
 For any Noetherian formal scheme\/ $\fX$, there exists a unique,
up to a unique isomorphism, derived category object\/
$\cB^\bu=\cB_\fX^\bu\in\sD^\bb(\fX\Tors)$ together with
\begin{itemize}
\item isomorphisms\/ $\rho_\fU^\fX:\fj_\fU^*\cB_\fX^\bu\simeq
\cM_\fU^\bu$ in\/ $\sD^\bb(\fU\Tors)$, given for all affine open
formal subschemes\/ $\fU\subset\fX$,
\end{itemize}
for which
\begin{itemize}
\item the equation
$\rho^\fX_\fV=\rho^\fU_\fV\circ\fh_{\fV,\fU}^*(\rho^\fX_\fU)$
holds for all pairs of affine open formal subschemes\/ $\fV\subset\fU$
in\/~$\fX$.
\end{itemize}
 Moreover, for any chosen open covering\/ $\bW$ of the formal scheme\/
$\fX$, a datum of the isomorphisms\/~$\rho^\fX_\fU$ for all affine
open subschemes\/ $\fU$ subordinate to\/ $\bW$ in\/ $\fX$, satifying
the equations above for all pairs of affine open subschemes\/
$\fV\subset\fU$ subordinate to\/ $\bW$, is sufficient to determine
a derived category object\/ $\cB_\fX^\bu$ uniquely up to a unique
isomorphism.
 In addition, one has
\begin{itemize}
\item $\Hom_{\sD^\bb(\fX\Tors)}(\cB_\fX^\bu,\cB_\fX^\bu[i])=0$
for all $i<0$.
\end{itemize}
\end{cor}

\begin{proof}
 The assumptions of Theorem~\ref{main-gluing-theorem} hold for
the trivial open covering $\bW_\fX=\{\fX\}$, the derived category
objects $\cB_\fU^\bu$, and the derived category
isomorphisms~$\rho^\fU_\fV$ (where $\fV\subset\fU$ are affine open
formal subschemes in~$\fX$) by
Lemma~\ref{dedualizing-complex-derived-endomorphisms}(a) and
(the proof of) Lemma~\ref{dedualizing-complex-extension-of-scalars}.
 Applying Theorem~\ref{main-gluing-theorem}, we obtain all
the assertions of the corollary, except for the ``moreover''
clause, which is based on the uniqueness assertion of
Theorem~\ref{main-gluing-theorem} for the given open covering~$\bW$.
\end{proof}

 For any Noetherian formal scheme $\fX$, we call the derived
category object $\cB^\bu=\cB_\fX^\bu\in\sD^\bb(\fX\Tors)$ provided
by Corollary~\ref{nonaffine-dedualizing-complex-cor}
\emph{the dedualizing complex of quasi-coherent torsion sheaves
on\/~$\fX$}.
 The dedualizing complex $\cB_\fX^\bu$ plays the crucial role in
the construction of the MGM equivalence in the next sections.

 The next corollary refers to the triangulated functor of inverse
image $\ff^*\:\sD^\bb(\fX\Tors)\allowbreak\rarrow\sD^\bb(\fY\Tors)$
\,\eqref{qcoh-derived-cats-underived-inverse-image} defined in
Section~\ref{triangulated-inverse-images-subsecn} for any
flat tight morphism of locally Noetherian formal schemes
$\ff\:\fY\rarrow\fX$.
 Similarly to the notation above, for any flat tight morphism of
affine Noetherian formal schemes $\fg\:\fV\rarrow\fU$, where
$\fU=\Spf\fR$ and $\fV=\Spf\fS$, consider the related flat tight
continuous homomorphism of adic Noetherian rings $g\:\fR\rarrow\fS$,
and denote by $\rho^\fU_\fV\:\fg^*(\cB_\fU^\bu)\simeq\cB_\fV^\bu$
the isomorphism in $\sD^\bb(\fV\Tors)$ corresponding to
the isomorphism $\rho^\fR_\fS\:g^*(B_\fR^\bu)\simeq B_\fS^\bu$
in $\sD^\bb(\fS\Modl_\tors)$ constructed in
Lemma~\ref{dedualizing-complex-extension-of-scalars}.

\begin{cor} \label{flat-inverse-image-of-dedualizing-is-dedualizing}
 Let\/ $\ff\:\fY\rarrow\fX$ be a flat tight morphism of Noetherian
formal schemes.
 Then there exists a unique isomorphism $\rho^\fX_\fY\:
\ff^*\cB_\fX^\bu\simeq\cB_\fY^\bu$ in the derived category\/
$\sD^\bb(\fY\Tors)$ such that
\begin{itemize}
\item for all affine open formal subschemes\/ $\fU\subset\fX$ and\/
$\fV\subset\fY$ such that\/ $\ff(\fV)\subset\fU$, denoting by\/ $\fg\:
\fV\rarrow\fU$ the induced flat tight morphism of affine Noetherian
formal schemes, and denoting by\/ $\fj\:\fU\rarrow\fX$ and\/
$\fk\:\fV\rarrow\fY$ the open immersion morphisms, we have
a commutative diagram
$$
 \xymatrix{
  \fk^*\cB_\fY^\bu \ar@{=}[rrr]^{\rho^\fY_\fV}
  \ar@{=}[dd]_{\fk^*(\rho^\fX_\fY)}
  &&& \cB_\fV^\bu \ar@{=}[dd]^{\rho^\fU_\fV} \ar@{=}[llldd] \\ \\
  \fk^*\ff^*\cB_\fX^\bu \ar@{=}[r] & \fg^*\fj^*\cB_\fX^\bu
  \ar@{=}[rr]_{\fg^*(\rho^\fX_\fU)} && \fg^*\cB_\fU^\bu
 }
$$
of isomorphisms in the derived category\/ $\sD^\bb(\fV\Tors)$.
\end{itemize}
\end{cor}

\begin{proof}
 Pick a pair of open coverings $\bW$ of $\fX$ and $\bT$ of $\fY$
such that the morphism~$\ff$ is $(\bW,\bT)$\+coaffine (in the sense of
the definition in Section~\ref{inverse-images-of-qcoh-tors-subsecn}).
 For any pair of affine open subschemes $\fU\subset\fX$ subordinate
to $\bW$ and $\fV\subset\fY$ subordinate to $\bT$, consider the 
the diagonal isomorphisms shown on the diagram, constructed from
the condition of commutativity of the lower right triangle.
 Use Proposition~\ref{gluing-morphism-uniqueness-prop} in order to
show that the diagonal isomorphism does not depend on the choice
of $\fU$ for any fixed~$\fV$.
 Then check that the complex $\ff^*\cB_\fX^\bu\in\sD^\bb(\fY\Tors)$,
together with the collection of the diagonal isomorphisms for all
affine open subschemes $\fV\subset\fY$ subordinate to $\bT$,
satisfies the equations in the definition of a dedualizing complex
on~$\fY$.
 Finally, use the ``moreover'' uniqueness assertion of
Corollary~\ref{nonaffine-dedualizing-complex-cor}.
\end{proof}

\begin{rem}
 A simper alternative construction of the dedualizing complex
$\cB_\fX^\bu$ on a Noetherian formal scheme $\fX$ was suggested
to the author by Amnon Yekutieli.
 It is based on the results of his paper~\cite{Yek0} and uses
homological algebra in the abelian category $(\fX,\fO_\fX)\Sh$ of
sheaves of $\fO_\fX$\+modules on~$\fX$.

 The construction works as follows.
 Let $0\rarrow\fO_\fX\rarrow\cJ^\bu$ be an injective coresolution
of the structure sheaf $\fO_\fX$ in the Grothendieck abelian
category $(\fX,\fO_\fX)\Sh$.
 Applying to $\cJ^\bu$ the functor of maximal discrete/torsion
subsheaf $\underline{\Gamma}_{\mathrm{disc}}$ defined
in~\cite[Definition~3.7]{Yek0}, one obtains a complex of discrete
sheaves of $\fO_\fX$\+modules on $\fX$ that can be denoted by
$\boR\underline{\Gamma}_{\mathrm{disc}}(\fO_X)=
\underline{\Gamma}_{\mathrm{disc}}(\cJ^\bu)$.
 According to~\cite[Theorem~4.8]{Yek0} (based on~\cite[Proposition~4.2
and Lemma~4.7]{Yek0}), the complex
$\boR\underline{\Gamma}_{\mathrm{disc}}(\fO_X)$ is quasi-isomorphic
to a uniquely defined object of the derived category (in our notation)
$\sD^+(\fX\Tors)$.
 (See also~\cite[Proposition~5.3.1]{AJL} for an alternative formulation
and proof.)
 The latter object is actually cohomologically bounded
by~\cite[Lemma~4.5]{Yek0}, so it can be truncated to obtain
an object of the bounded derived category $\sD^\bb(\fX\Tors)$.

 The resulting object of $\sD^\bb(\fX\Tors)$ is naturally isomorphic
to our dedualizing complex $\cB_\fX^\bu$, as one can see
from~\cite[Lemma~4.5]{Yek0} and/or~\cite[Example~5.8]{Pmgm}.
\end{rem}

\subsection{MGM equivalence for formal schemes of finite
Krull dimension~I} \label{mgm-equiv-fin-krulldim-I-subsecn}
 The following theorem is the first main result of this paper.
 It is our version of the MGM equivalence on formal schemes,
interpreted as a triangulated equivalence between the derived
categories of the abelian category of quasi-coherent torsion sheaves
and the exact category of contraherent cosheaves of contramodules
in the spirit of~\cite[Corollary~3.5 or Theorem~5.10]{Pmgm}
and~\cite[Corollary~8.7]{Ppt}.
 For a generalization including the absolute derived categories,
see Main Theorem~\ref{fin-krulldim-mgm-main-theorem-II} below.

\begin{mainthm} \label{fin-krulldim-mgm-main-theorem-I}
 Let\/ $\fX$ be a (not necessarily semi-separated) Noetherian formal
scheme of finite Krull dimension, and let\/ $\cB_\fX^\bu\in
\sD^\bb(\fX\Tors)$ be the dedualizing complex of\/ $\fX$ as per
Section~\ref{dedualizing-complexes-on-formal-schemes-subsecn}.
 Let\/ $\st=\bb$, $+$, $-$, or\/~$\varnothing$ be a conventional
derived category symbol.
 Then there is a natural equivalence of triangulated categories
\begin{equation} \label{fin-krulldim-mgm-equivalence-eqn-I}
 \boR\fHom_\fX(\cB_\fX^\bu,{-})\:\sD^\st(\fX\Tors)\,\simeq\,
 \sD^\st(\fX\Ctrh):\!\cB_\fX^\bu\ocn_\fX^\boL{-}
\end{equation}
given by derived functors of the functor\/ $\fHom_\fX(\cB_\fX^\bu,{-})$
from Section~\ref{contraherent-Hom-subsecn} and the contratensor
product functor\/ $\cB_\fX^\bu\ocn_\fX{-}$ from
Section~\ref{contratensor-product-subsecn}.
\end{mainthm}

\begin{proof}
 This is a far-reaching formal scheme generalization
of~\cite[Theorem~6.8.1]{Pcosh}, and the proof is similar.
 The right-hand side of
the desired equivalence~\eqref{fin-krulldim-mgm-equivalence-eqn-I}
can be constructed in several alternative ways, as it is
the derived category of many exact categories naturally embedded
into one another, as per the commutative diagram of triangulated
equivalences~\eqref{derived-equivs-induced-by-inclusions-diag-I}.
 For any open covering $\bW$ of the formal scheme $\fX$, we will
construct a triangulated equivalence $\sD(\fX\Tors)\simeq
\sD(\fX\Lcth^\lct_\bW)$, and then show that it takes the full
subcategories $\sD^\st(\fX\Tors)\subset\sD(\fX\Tors)$ and
$\sD^\st(\fX\Lcth^\lct_\bW)\subset\sD(\fX\Lcth^\lct_\bW)$ into
one another for all symbols $\st=\bb$, $+$, or~$-$.

 The notation $\boR\fHom_\fX(\cB_\fX^\bu,{-})$ and
$\cB_\fX^\bu\ocn_\fX^\boL{-}$
in~\eqref{fin-krulldim-mgm-equivalence-eqn-I} needs to
be taken with a grain of salt, as in some sense our construction is
that of \emph{two-sided} derived functors.
 The point is that, in addition to coresolving the complexes of
quasi-coherent torsion sheaves and resolving the complexes of (locally
cotorsion locally) contraherent cosheaves of contramodules representing
the objects of the derived categories involved, we also need to
coresolve the dedualizing complex $\cB_\fX^\bu$ itself.

 Specifically, by
Corollary~\ref{qcoh-tors-flasque-dil-derived-equivalences}(b),
the derived category object $\cB_\fX^\bu\in\sD^\bb(\fX\Tors)$ can be
repreresented by a finite complex of flasque quasi-coherent torsion
sheaves $\cE^\bu$ on~$\fX$.
 We fix such a complex $\cE^\bu$ for the time being.

 Let us define the functor $\boR\fHom_\fX(\cB_\fX^\bu,{-})\:
\sD(\fX\Tors)\rarrow\sD(\fX\Ctrh^\lct)$.
 It follows, e.~g., from Theorem~\ref{qcoh-tors-coderived} (or from
the existence of enough homotopy injective complexes of injective
objects in any Grothendieck abelian category, cf.\
Section~\ref{derived-direct-images-of-qcoh-tors-subsecn}) that
for any complex $\cM^\bu$ in $\fX\Tors$ there exists a complex
$\cJ^\bu$ in $\fX\Tors^\inj$ together with a quasi-isomorphism of
complexes of quasi-coherent torsion sheaves $\cM^\bu\rarrow\cJ^\bu$.
 By the definition, we let the functor $\boR\fHom_\fX(\cB_\fX^\bu,{-})$
assign to $\cM^\bu$ the totalization of the bicomplex
$\fHom_\fX(\cE^\bu,\cJ^\bu)$ in $\fX\Ctrh^\lct$.
 See the construction in Section~\ref{contraherent-Hom-subsecn}.

 The following lemma shows that the derived category object
$\boR\fHom_\fX(\cB_\fX^\bu,\cM^\bu)\in\sD(\fX\Ctrh^\lct)$ is
well-defined for $\cM^\bu\in\sD(\fX\Tors)$.

\begin{lem} \label{derived-fHom-takes-acyclics-to-acyclics}
 For any complex\/ $\cJ^\bu$ in\/ $\fX\Tors^\inj$ that is acyclic
in the abelian category\/ $\fX\Tors$, the complex\/
$\fHom_\fX(\cE^\bu,\cJ^\bu)$ is acyclic in the exact category\/
$\fX\Ctrh^\lct$.
\end{lem}

\begin{proof}
 By Theorem~\ref{loc-Noetherian-injective-qcoh-torsion}(b),
any injective quasi-coherent torsion sheaf $\cI$ on $\fX$ can be
presented as an infinite direct sum $\cI=\bigoplus_{x\in\fX}
\iota_x{}_*\widetilde I_x$, where $I_x$ are some injective
$\widehat\fm_{x,\fX}$\+torsion $\widehat\fO_{x,\fX}$\+modules.
 Let $\fX=\bigcup_{\alpha=1}^N\fU_\alpha$ be a finite affine open
covering of~$\fX$.
 Set $\fS_\beta=\fU_\beta\setminus\bigcup_{\alpha<\beta}\fU_\alpha$
for $1\le\beta\le N$, and consider the finite direct sum decomposition
$\cI=\bigoplus_{\beta=1}^N\cI_\beta$, where
$\cI_\beta=\bigoplus_{z\in\fS_\beta}\iota_z{}_*\widetilde I_z$.
 We are interested in the finite decreasing filtration on $\cI$
associated with this direct sum decomposition, $\cI_{\ge\alpha}=
\bigoplus_{\beta=\alpha}^N\cI_\beta$ for all $1\le\alpha\le N+1$.

 By Lemma~\ref{sheaves-lower-upper-sets-fully-faithful-lemma}(a),
the decreasing filtration $(\cI_{\ge\alpha})_{\alpha=1}^{N+1}$
is preserved by all morphisms of injective quasi-coherent torsion
sheaves $\cI$ on~$\fX$.
 Hence we obtain a termwise split filtration $\cJ^\bu_{\ge\alpha}$
on any complex $\cJ^\bu$ in $\fX\Tors^\inj$.
 The associated quotient complexes $\cJ^\bu_\alpha=
\cJ^\bu_{\ge\alpha}/\cJ^\bu_{\ge\alpha+1}$ are isomorphic to
the direct images $\fj_\alpha{}_*\cK^\bu_\alpha$ of uniquely defined
complexes of injective quasi-coherent torsion sheaves $\cK^\bu_\alpha$
on~$\fU_\alpha$ (where $\fj_\alpha\:\fU_\alpha\rarrow\fX$ denote
the identity open immersion morphisms).

 The next observation is that the complexes $\cK^\bu_\alpha$ are
acyclic in $\fU_\alpha\Tors$ and the complexes $\cJ^\bu_\alpha$ are
acyclic in $\fX\Tors$ for all $1\le\alpha\le N$ whenever the complex
$\cJ^\bu$ in $\fX\Tors^\inj$ is acyclic in $\fX\Tors$.
 Indeed, for $\alpha=1$ we have $\cK^\bu_1\simeq\fj_1^*\cJ^\bu$,
which is an acyclic complex in $\fU_1\Tors$ since the functor
$\fj_1^*\:\fX\Tors\rarrow\fU_1\Tors$ is exact.
 Furthermore, the complex $\cJ^\bu_1=\fj_1{}_*\cK^\bu_1$ is acyclic in
$\fX\Tors$ by Lemma~\ref{finite-Krull-dim-co-flasque-co-resol-dim}(a)
with Corollary~\ref{co-flasque-direct-image-cor}(a), or directly by
the more general
Corollary~\ref{co-flasque-direct-image-of-acyclic-complexes-cor}(a).
 Hence the complex $\cJ^\bu_{\ge2}=\ker(\cJ^\bu\to\cJ^\bu_1)$ is
acyclic in $\fX\Tors$, too.
 Now we have $\cK^\bu_2\simeq\fj_2^*\cJ^\bu_{\ge2}$.
 Proceeding by induction, one proves the claimed acyclicity of
the complexes $\cK^\bu_\alpha$ and $\cJ^\bu_\alpha$ for all
$1\le\alpha\le N$.

 Now we have $\fHom_\fX(\cE^\bu,\fj_\alpha{}_*\cK^\bu_\alpha)\simeq
\fj_\alpha{}_!\fHom_{\fU_\alpha}(\fj_\alpha^*\cE^\bu,\cK^\bu_\alpha)$
by Proposition~\ref{fHom-projection-formula-prop}.
 The complex $\fHom_{\fU_\alpha}(\fj_\alpha^*\cE^\bu,\cK^\bu_\alpha)$
of locally cotorsion contraherent cosheaves of contramodules on
the affine Noetherian formal scheme $\fU_\alpha$ corresponds,
under the equivalence of exact categories from
Lemma~\ref{lct-contrah-cosheaves-of-contramods-on-affine-lemma},
to the complex of cotorsion contramodule $\fO_\fX(\fU_\alpha)$\+modules
$\Hom_{\fO_\fX(\fU_\alpha)}(\cE^\bu(\fU_\alpha),
\cK^\bu_\alpha(\fU_\alpha))$ (by the definition of $\fHom$).
 Since the complex of injective objects $\cK^\bu_\alpha$ is acyclic
in $\fU_\alpha\Tors$, while $\cE^\bu(\fU_\alpha)$ is a dedualizing
complex of torsion $\fO_\fX(\fU_\alpha)$\+modules, it follows by virtue
of Corollary~\ref{Hom-tensor-B-preserves-acyclicity-boundedness-cor}(a)
that the complex $\Hom_{\fO_\fX(\fU_\alpha)}(\cE^\bu(\fU_\alpha),
\cK^\bu_\alpha(\fU_\alpha))$ is acyclic in $\fO_\fX(\fU_\alpha)\Contra$.
 Using the result of Proposition~\ref{adic-raynaud-gruson-prop}(b)
or the much more general result of
Proposition~\ref{contraadjusted-and-cotorsion-periodicity-prop}(b),
one can see that the complex
$\Hom_{\fO_\fX(\fU_\alpha)}(\cE^\bu(\fU_\alpha),
\cK^\bu_\alpha(\fU_\alpha))$ is acyclic in
$\fO_\fX(\fU_\alpha)\Contra^\cot$.
 Thus the complex
$\fHom_{\fU_\alpha}(\fj_\alpha^*\cE^\bu,\cK^\bu_\alpha)$ is
acyclic in $\fU_\alpha\Ctrh^\lct$.

 Furthermore, $\fj_\alpha^*\cE^\bu$ is a complex of flasque
quasi-coherent torsion sheaves on~$\fU_\alpha$.
 By Lemma~\ref{fHom-contratensor-co-flasque}(a), it follows that
$\fHom_{\fU_\alpha}(\fj_\alpha^*\cE^\bu,\cK^\bu_\alpha)$ is
a complex of coflasque locally cotorsion contraherent cosheaves of
contramodules on~$\fU_\alpha$.
 Applying Lemma~\ref{finite-Krull-dim-co-flasque-co-resol-dim}(b)
with Corollary~\ref{co-flasque-direct-image-cor}(c), or directly
the more general
Corollary~\ref{co-flasque-direct-image-of-acyclic-complexes-cor}(c),
we conclude that $\fHom_\fX(\cE^\bu,\cJ^\bu_\alpha)\simeq
\fj_\alpha{}_!\fHom_{\fU_\alpha}(\fj_\alpha^*\cE^\bu,\cK^\bu_\alpha)$
is an acyclic complex in $\fX\Ctrh^\lct$.
 As the class of acyclic complexes in an exact category is closed
under extensions, it follows that the complex
$\fHom_\fX(\cE^\bu,\cJ^\bu)$ is also acyclic in $\fX\Ctrh^\lct$.
\end{proof}

 Let us define the functor $\cB_\fX^\bu\ocn_\fX^\boL{-}\,\:
\sD(\fX\Lcth^\lct_\bW)\rarrow\sD(\fX\Tors)$.
 It follows from Theorem~\ref{lcth-lct-W-contramod-contraderived} that
for any complex $\fP^\bu$ in $\fX\Lcth^\lct_\bW$ there exists
a complex $\fF^\bu$ in $\fX\Ctrh^\lct_\prj$ together with
a quasi-isomorphism $\fF^\bu\rarrow\fP^\bu$ of complexes in
$\fX\Lcth^\lct_\bW$ (indeed, any Positselski-contraacyclic complex in
an exact category with exact functors of infinite product is acyclic).
 By the definition, we let the functor $\cB_\fX^\bu\ocn_\fX^\boL{-}$
assign to $\fP^\bu$ the totalization of the bicomplex
$\cE^\bu\ocn_\fX\fF^\bu$ in $\fX\Tors$.
 See the construction in Section~\ref{contratensor-product-subsecn}.

 The next lemma shows that the derived category object
$\cB_\fX^\bu\ocn_\fX^\boL\fP^\bu\in\sD(\fX\Tors)$ is well-defined
for $\fP^\bu\in\sD(\fX\Lcth^\lct_\bW)$.

\begin{lem} \label{derived-contratensor-takes-acyclics-to-acyclics}
 For any complex\/ $\fF^\bu$ in\/ $\fX\Ctrh^\lct_\prj$ that is
acyclic in the exact category\/ $\fX\Lcth^\lct_\bW$, the complex\/
$\cE^\bu\ocn_\fX\fF^\bu$ is acyclic in the abelian category\/
$\fX\Tors$.
\end{lem}

\begin{proof}
 The argument is dual-analogous to the proof of
Lemma~\ref{derived-fHom-takes-acyclics-to-acyclics}.
 By Theorem~\ref{loc-Noetherian-projective-lct-cosh-contramods}(b),
any projective locally cotorsion contraherent cosheaf of
contramodules $\fH$ on $\fX$ can be presented as an infinite direct
product $\fH=\prod_{x\in\fX}\iota_x{}_!\widecheck H_x$,
where $H_x$ are some free $\widehat\fm_{x,\fX}$\+contramodule
$\widehat\fO_{x,\fX}$\+modules.
 Let $\fX=\bigcup_{\alpha=1}^N\fU_\alpha$ be a finite affine open
covering of~$\fX$.
 Set $\fS_\beta=\fU_\beta\setminus\bigcup_{\alpha<\beta}\fU_\alpha$
for $1\le\beta\le N$, and consider the finite direct sum decomposition
$\fH=\bigoplus_{\beta=1}^N\fH_\beta$, where
$\fH_\beta=\prod_{z\in\fS_\beta}\iota_z{}_!\widecheck H_z$.
 We are interested in the finite increasing filtration on $\fH$
associated with this direct sum decomposition, $\fH_{\le\alpha}=
\bigoplus_{\gamma=1}^\alpha\fH_\gamma$ for all $0\le\alpha\le N$.

 By Lemma~\ref{cosheaves-lower-upper-sets-fully-faithful-lemma}(a),
the increasing filtration $(\fH_{\le\alpha})_{\alpha=0}^N$ is
preserved by all morphisms of projective locally cotorsion
contraherent cosheaves of contramodules $\fH$ on~$\fX$.
 Hence we obtain a termwise split filtration $\fF^\bu_{\le\alpha}$
on any complex $\fF^\bu$ in $\fX\Ctrh^\lct_\prj$.
 The associated quotient complexes $\fF^\bu_\alpha=
\fF^\bu_{\le\alpha}/\fF^\bu_{\le\alpha-1}$ are isomorphic to
the direct images $\fj_\alpha{}_!\fG_\alpha^\bu$ of uniquely
defined complexes of projective locally cotorsion contraherent
cosheaves of contramodules $\fG_\alpha^\bu$ on~$\fU_\alpha$.

 The next observation is that the complexes $\fG_\alpha^\bu$ are
acyclic in $\fU_\alpha\Ctrh$ and the complexes $\fF_\alpha^\bu$ are
acyclic in $\fF\Lcth^\lct_\bW$ for all $1\le\alpha\le N$ whenever
the complex $\fF^\bu$ in $\fX\Ctrh^\lct_\prj$ is acyclic in
$\fX\Lcth^\lct_\bW$.
 Indeed, for $\alpha=1$ we have $\fG_1^\bu\simeq\fj_1^!\fF^\bu$,
which is an acyclic complex in $\fU_1\Ctrh^\lct$ since the functor
$\fj_1^!\:\fX\Lcth^\lct_\bW\rarrow\fU_1\Ctrh^\lct$ is exact.
 Furthermore, the complex $\fF_1^\bu=\fj_1{}_!\fG_1^\bu$ is acyclic
in $\fX\Ctrh^\lct$ by
Lemma~\ref{finite-Krull-dim-co-flasque-co-resol-dim}(b)
with Corollary~\ref{co-flasque-direct-image-cor}(c), or directly by
the more general
Corollary~\ref{co-flasque-direct-image-of-acyclic-complexes-cor}(c).
 Hence the complex $\fF^\bu/\fF^\bu_1$ is acyclic in
$\fX\Lcth^\lct_\bW$ as well.
 Now we have $\fG_2^\bu\simeq\fj_2^!(\fF^\bu/\fF^\bu_1)$.
 Proceeding by induction, one proves the claimed acyclicity of
the complexes $\fG_\alpha^\bu$ and $\fF_\alpha^\bu$ for all
$1\le\alpha\le N$.

 Now we have $\cE^\bu\ocn_\fX\fj_\alpha{}_!\fG_\alpha^\bu\simeq
\fj_\alpha{}_*(\fj_\alpha^*\cE^\bu\ocn_{\fU_\alpha}\fG_\alpha^\bu)$
by Proposition~\ref{contratensor-projection-formula-prop}(c).
 (Here and below, one can also reference the more general result of
Proposition~\ref{flat-contratensor-projection-formula-prop} instead
of Proposition~\ref{contratensor-projection-formula-prop}(c).)
 The complex $\fj_\alpha^*\cE^\bu\ocn_{\fU_\alpha}\fG_\alpha^\bu$ of
quasi-coherent torsion sheaves on the affine Noetherian formal
scheme $\fU_\alpha$ corresponds, under the natural equivalence of
abelian categories $\fU_\alpha\Tors\simeq\fO_\fX(\fU_\alpha)\Tors$,
to the complex of torsion $\fO_\fX(\fU_\alpha)$\+modules
$\cE^\bu(\fU_\alpha)\ot_{\fO_\fX(\fU_\alpha)}
\fG_\alpha^\bu[\fU_\alpha]$ (by the definition of the contratensor
product functor~$\ocn$).
 Since the complex of flat cotorsion contramodule
$\fO_\fX(\fU_\alpha)$\+modules $\fG_\alpha^\bu[\fU_\alpha]$ is
acyclic in $\fO_\fX(\fU_\alpha)\Contra$ (cf.\ the equivalence
of exact categories from
Lemma~\ref{lct-contrah-cosheaves-of-contramods-on-affine-lemma}),
while $\cE^\bu(\fU_\alpha)$ is a dedualizing complex of torsion
$\fO_\fX(\fU_\alpha)$\+modules, it follows by virtue of
Corollary~\ref{Hom-tensor-B-preserves-acyclicity-boundedness-cor}(c)
that the complex $\cE^\bu(\fU_\alpha)\ot_{\fO_\fX(\fU_\alpha)}
\fG_\alpha^\bu[\fU_\alpha]$ is acyclic in $\fO_\fX(\fU_\alpha)\Tors$.
 Thus the complex $\fj_\alpha^*\cE^\bu\ocn_{\fU_\alpha}\fG_\alpha^\bu$
is acyclic in $\fU_\alpha\Tors$.

 Furthermore, $\fj_\alpha^*\cE^\bu$ is a complex of flasque
quasi-coherent torsion sheaves on~$\fU_\alpha$.
 By Lemma~\ref{fHom-contratensor-co-flasque}(b), it follows that
$\fj_\alpha^*\cE^\bu\ocn_{\fU_\alpha}\fG_\alpha^\bu$ is a complex of
flasque quasi-coherent torsion sheaves on $\fU_\alpha$ as well.
 Applying Lemma~\ref{finite-Krull-dim-co-flasque-co-resol-dim}(a)
with Corollary~\ref{co-flasque-direct-image-cor}(a), or directly
the more general
Corollary~\ref{co-flasque-direct-image-of-acyclic-complexes-cor}(a),
we conclude that $\cE^\bu\ocn_\fX\fF_\alpha^\bu\simeq
\fj_\alpha{}_*(\fj_\alpha^*\cE^\bu\ocn_{\fU_\alpha}\fG_\alpha^\bu)$
is an acyclic complex in $\fX\Tors$.
 Recall that the filtration $\fF^\bu_{\le\alpha}$ on the complex
$\fF^\bu$ is termwise split, so it is preserved by the functor
$\cE^\bu\ocn_\fX{-}$.
 As the class of acyclic complexes in an abelian category is closed
under extensions, it follows that the complex $\cE^\bu\ocn_\fX\fF^\bu$
is also acyclic in $\fX\Tors$.
\end{proof}

 We have constructed the desired derived functors
\begin{equation} \label{derived-fHom-from-dedualizing}
 \boR\fHom_\fX(\cB^\bu_\fX,{-})\:\sD(\fX\Tors)\lrarrow
 \sD(\fX\Lcth^\lct_\bW)
\end{equation}
and
\begin{equation} \label{derived-contratensor-with-dedualizing}
 \cB^\bu_\fX\ocn_\fX^\boL{-}\,\:\sD(\fX\Lcth^\lct_\bW)\lrarrow
 \sD(\fX\Tors).
\end{equation}

 Let us show that the functor~\eqref{derived-fHom-from-dedualizing}
does not depend on the choice of a finite complex of flasque
quasi-coherent torsion sheaves $\cE^\bu$ representing the derived
category object~$\cB_\fX^\bu$.
 For this purpose, we point out that, for any acyclic complex $\cL^\bu$
in $\fX\Tors^\fq$ and any injective quasi-coherent torsion sheaf $\cI$
on $\fX$, the complex of locally cotorsion contraherent cosheaves of
contramodules $\fHom_\fX(\cL^\bu,\cI)$ on $\fX$ is acyclic
by construction.
 In other words, the functor $\fHom_\fX({-},\cI)\:
(\fX\Tors^\fq)^\sop\rarrow\fX\Ctrh^\lct$ is exact.
 Indeed, for any open formal subscheme $\fU\subset\fX$ with the open
immersion morphism $\fj\:\fU\rarrow\fX$, the functor $\fj_*\fj^*\:
\fX\Tors^\fq\rarrow\fX\Tors^\fq$ is exact
by Corollary~\ref{co-flasque-direct-image-cor}(a).
 Alternatively, one can refer to the more general result of
Corollary~\ref{co-flasque-direct-image-of-acyclic-complexes-cor}(a).

 To prove that the functor~\eqref{derived-contratensor-with-dedualizing}
does not depend on the choice of $\cE^\bu$, we need to check that,
in the same notation, the complex $\cL^\bu\ocn_\fX\fH$ is acyclic
in $\fX\Tors$ for any projective locally cotorsion contraherent
cosheaf of contramodules $\fH$ on~$\fX$.
 Here it suffices to refer to
Lemma~\ref{contratensor-exact-on-flats-and-alf}(b).
 Alternatively, it follows from
Theorem~\ref{loc-Noetherian-projective-lct-cosh-contramods}(b) that
the cosheaf $\fH$ has the form $\fH=\bigoplus_\alpha\fj_\alpha{}_!
\fG_\alpha$, where $\fG_\alpha$ are some projective locally cotorsion
contraherent cosheaves of contramodules on~$\fU_\alpha$ (in
the notation of the proofs of the lemmas above).
 By Proposition~\ref{contratensor-projection-formula-prop}(c),
we have $\cL^\bu\ocn_\fX\fj_\alpha{}_!\fG_\alpha\simeq
\fj_\alpha{}_*(\fj_\alpha^*\cL^\bu\ocn_{\fU_\alpha}\fG_\alpha)$.
 Now $\fj_\alpha^*\cL^\bu\ocn_{\fU_\alpha}\fG_\alpha$ is an acyclic
complex of quasi-coherent torsion sheaves on $\fU_\alpha$ by
construction (since the $\fO_\fX(\fU_\alpha)$\+module
$\fG_\alpha[\fU_\alpha]$) is flat), and it is also a complex of
flasque quasi-coherent torsion sheaves on $\fU_\alpha$ by
Lemma~\ref{fHom-contratensor-co-flasque}(b).
 As above, it remains to refer either to
Lemma~\ref{finite-Krull-dim-co-flasque-co-resol-dim}(a)
with Corollary~\ref{co-flasque-direct-image-cor}(a), or directly to
Corollary~\ref{co-flasque-direct-image-of-acyclic-complexes-cor}(a).

 It is a standard fact that
the adjunction~\eqref{fHom-contratensor-adjunction-eqn} leads to
an adjunction of the two
derived functors~\eqref{derived-fHom-from-dedualizing}
and~\eqref{derived-contratensor-with-dedualizing} constructed as above
(cf.~\cite[proof of Lemma in Section~8.3]{Psemi}).
 Quite generally, this is applicable to any derived functors in
the sense of Deligne~\cite[Sections~1.2.1\+-1.2.2]{Del}.
 Our next aim is to show that the related adjunction morphisms are
isomorphisms in the derived categories.

 Let us check that the adjunction morphism $\cB_\fX^\bu\ocn_\fX^\boL
\fHom_\fX(\cE^\bu,\cJ^\bu)\rarrow\cJ^\bu$ is an isomorphism in
$\sD(\fX\Tors)$ for any complex $\cJ^\bu$ in $\fX\Tors^\inj$.
 Due to the existence of the finite filtration $\cJ^\bu_{\ge\alpha}$
constructed in the proof of
Lemma~\ref{derived-fHom-takes-acyclics-to-acyclics}, we can assume
that $\cJ^\bu=\fj_*\cK^\bu$ for some affine open formal subscheme
$\fU\subset\fX$ with the open immersion morphism $\fj\:\fU\rarrow\fX$
and some complex $\cK^\bu$ in $\fU\Tors$.
 Then $\fHom_\fX(\cE^\bu,\cJ^\bu)\simeq\fj_!\fHom_\fU(\fj^*\cE^\bu,
\cK^\bu)$ by Proposition~\ref{fHom-projection-formula-prop}.

 Let $\fG^\bu$ be a complex in $\fU\Ctrh^\lct_\prj$ endowed with
a quasi-isomorphism $\fG^\bu\rarrow\fHom_\fU(\fj^*\cE^\bu,\cK^\bu)$
of complexes in $\fU\Ctrh^\lct$.
 Then $\fj_!\fG^\bu$ is a complex in $\fX\Ctrh^\lct_\prj$ endowed
with a quasi-isomorphism $\fj_!\fG^\bu\rarrow
\fj_!\fHom_\fU(\fj^*\cE^\bu,\cK^\bu)$ of complexes in $\fX\Ctrh^\lct$
(as both $\fG^\bu$ and $\fHom_\fU(\fj^*\cE^\bu,\cK^\bu)$ are
complexes of coflasque locally cotorsion contraherent cosheaves of
contramodules on $\fU$, by Lemma~\ref{fHom-contratensor-co-flasque}(a)).
 Hence we have $\cB_\fX^\bu\ocn_\fX^\boL\fHom_\fX(\cE^\bu,\cJ^\bu)=
\cE^\bu\ocn_\fX\fj_!\fG^\bu\simeq\fj_*(\fj^*\cE^\bu\ocn_\fU\fG^\bu)$
by Proposition~\ref{contratensor-projection-formula-prop}(c).
 Both $\cK^\bu$ and $\fj^*\cE^\bu\ocn_\fU\fG^\bu$ are complexes of
flasque quasi-coherent torsion sheaves on~$\fU$ (by
Lemma~\ref{fHom-contratensor-co-flasque}(b)), so it remains to
show that the natural morphism $\fj^*\cE^\bu\ocn_\fU\fG^\bu
\rarrow\cK^\bu$ is a quasi-isomorphism in $\fX\Tors$.
 This is the result of
Proposition~\ref{Hom-tensor-B-quasi-isomorphisms-prop}(a).

 Dual-analogously, let us check that the adjunction morphism
$\boR\fHom_\fX(\cB_\fX^\bu,\>\allowbreak\cE^\bu\ocn_\fX\nobreak\fF^\bu)
\allowbreak\rarrow\fF^\bu$ is an isomorphism in $\sD(\fX\Ctrh^\lct)$
for any complex $\fF^\bu$ in $\fX\Ctrh^\lct_\prj$.
 Due to the existence of the finite filtration $\fF^\bu_{\le\alpha}$
constructed in the proof of
Lemma~\ref{derived-contratensor-takes-acyclics-to-acyclics},
we can assume that $\fF^\bu=\fj_!\fG^\bu$ for some complex
$\fG^\bu$ in $\fU\Ctrh^\lct_\prj$, in the notation above.
 Then $\cE^\bu\ocn_\fX\fF^\bu\simeq\fj_*(\fj^*\cE^\bu\ocn_\fU\fG^\bu)$
by Proposition~\ref{contratensor-projection-formula-prop}(c).

 Let $\cK^\bu$ be a complex in $\fU\Tors^\inj$ endowed with
a quasi-isomorphism $\cK^\bu\rarrow\fj^*\cE^\bu\ocn_\fU\fG^\bu$ of
complexes in $\fU\Tors$.
 Then $\fj_*\cK^\bu$ is a complex in $\fX\Tors$ endowed with
a quasi-isomorphism $\fj_*\cK^\bu\rarrow
\fj_*(\fj^*\cE^\bu\ocn_\fU\fG^\bu)$ of complexes in $\fX\Tors$
(as both $\cK^\bu$ and $\fj^*\cE^\bu\ocn_\fU\fG^\bu$ are complexes
of flasque quasi-coherent torsion sheaves, by
Lemma~\ref{fHom-contratensor-co-flasque}(b)).
 Hence we have $\boR\fHom_\fX(\cB_\fX^\bu,\>\allowbreak
\cE^\bu\ocn_\fX\nobreak\fF^\bu)=\fHom_\fX(\cE^\bu,\fj_*\cK^\bu)
\simeq\fj_!\fHom_\fU(\fj^*\cE^\bu,\cK^\bu)$ by
Proposition~\ref{fHom-projection-formula-prop}.
 Both $\fG^\bu$ and $\fHom_\fU(\fj^*\cE^\bu,\cK^\bu)$ are complexes
of coflasque locally cotorsion contraherent cosheaves of contramodules
on~$\fU$ (by Lemma~\ref{fHom-contratensor-co-flasque}(a)), so it
remains to show that the natural morphism $\fG^\bu\rarrow
\fHom_\fU(\fj^*\cE^\bu,\cK^\bu)$ is a quasi-isomorphism in
$\fU\Ctrh^\lct$.
 By Proposition~\ref{Hom-tensor-B-quasi-isomorphisms-prop}(b),
the map of complexes of $\fO_\fX(\fU)$\+modules
$\fG^\bu[\fU]\rarrow\fHom_\fU(\fj^*\cE^\bu,\cK^\bu)[\fU]$ is
a quasi-isomorphism.
 Finally, we refer to Proposition~\ref{adic-raynaud-gruson-prop}(b)
or the much more general result of
Proposition~\ref{contraadjusted-and-cotorsion-periodicity-prop}(b)
for the assertion that any quasi-isomorphism in
$\fO_\fX(\fU)\Contra$ between two complexes in
$\fO_\fX(\fU)\Contra^\cot$ is a quasi-isomorphism in
$\fO_\fX(\fU)\Contra^\cot$. {\emergencystretch=1em\par}

 We have constructed the desired triangulated
equivalence~\eqref{fin-krulldim-mgm-equivalence-eqn-I} for the unbounded
derived categories, $\st=\varnothing$.
 It remains to show that the derived functors
$\boR\fHom_\fX(\cB_\fX^\bu,{-})$ \,\eqref{derived-fHom-from-dedualizing}
and $\cB_\fX^\bu\ocn_\fX^\boL{-}$
\,\eqref{derived-contratensor-with-dedualizing} take bounded above
(respectively, below) complexes in $\fX\Tors$ to bounded above
(resp., below) complexes in $\fX\Lcth^\lct_\bW$ and vice versa, up
to quasi-isomorphism in the respective abelian or exact category.

 If a complex $\cM^\bu$ in $\fX\Tors$ is bounded below, then it has
a bounded below injective coresolution $\cJ^\bu$, and the complex
$\fHom_\fX(\cE^\bu,\cJ^\bu)$ in $\fX\Ctrh^\lct$ is also bounded below.
 The case of a bounded above complex $\cM^\bu$ in $\fX\Tors$
requires a bit more work, and
Lemma~\ref{co-flasque-open-subscheme-adjunction-morphism-lemma}(a)
is helpful in this context.

\begin{lem} \label{derived-fHom-takes-bounded-to-bounded}
 For any bounded above (respectively, two-sided bounded) complex\/
$\cM^\bu$ in\/ $\fX\Tors$, the derived category object\/
$\boR\fHom_\fX(\cB_\fX^\bu,\cM^\bu)\in\sD(\fX\Lcth_\bW^\lct)$
can be represented by a bounded above (resp., two-sided bounded)
complex in\/ $\fX\Lcth_\bW^\lct$.
\end{lem}

\begin{proof}
 Let $\cM^\bu$ be a bounded above complex in $\fX\Tors$.
 Using Corollary~\ref{qcoh-tors-flasque-dil-derived-equivalences}(b),
one can find a bounded above complex in $\cN^\bu$ in $\fX\Tors^\fq$
together with a quasi-isomorphism $\cM^\bu\rarrow\cN^\bu$ of
complexes in $\fX\Tors$.
 Suppose that the complex $\cN^\bu$ is concentrated in
the cohomological degrees not exceeding some integer~$n$.
 Let $\cJ^\bu$ be a complex in $\fX\Tors^\inj$ such that there is
a quasi-isomorphism $\cN^\bu\rarrow\cJ^\bu$ of complexes in
$\fX\Tors$ (or equivalently, in $\fX\Tors^\fq$).

 Consider the finite filtration $\cJ^\bu_{\ge\alpha}$ on $\cJ^\bu$
constructed in the proof of
Lemma~\ref{derived-fHom-takes-acyclics-to-acyclics}.
 Following that proof, one shows by induction on~$\alpha$
that every complex $\cJ^\bu_\alpha$ is quasi-isomophic in $\fX\Tors$
to a complex of flasque quasi-coherent torsion sheaves on $\fX$
concentrated in the cohomological degrees~$\le n$, while
the complex $\cK^\bu_\alpha$ is quasi-isomorphic in $\fU_\alpha\Tors$
to a complex of flasque quasi-coherent torsion sheaves on $\fU_\alpha$
concentrated in the cohomological degrees~$\le n$.

 Specifically, the complex $\cK_1^\bu\simeq\fj_1^*\cJ^\bu$ is
quasi-isomorphic to $\fj_1^*\cN^\bu$, and the complex $\cJ_1^\bu$
is quasi-isomorphic to $\cN_1^\bu=\fj_1{}_*\fj_1^*\cN^\bu$.
 By Lemma~\ref{co-flasque-open-subscheme-adjunction-morphism-lemma}(a),
it follows that the complex $\cJ^\bu_{\ge2}=\ker(\cJ^\bu\to\cJ^\bu_1)$
is quasi-isomorphic to the complex $\cN^\bu_{\ge2}=
\ker(\cN^\bu\to\cN^\bu_1)$, which is also a complex of flasque
quasi-coherent torsion sheaves on $\fX$ concentrated in
the cohomological degrees~$\le n$.
 Then the complex $\cK_2^\bu$ is quasi-isomorphic to
$\fj_2^*\cN^\bu_{\ge2}$, the complex $\cJ_2^\bu$ is quasi-isomorphic
to $\fj_2{}_*\fj_2^*\cN^\bu_{\ge2}$, etc.

 Now the complex of cotorsion contramodule
$\fO_\fX(\fU_\alpha)$\+modules $\Hom_{\fO_\fX(\fU_\alpha)}
(\cE^\bu(\fU_\alpha),\allowbreak\cK^\bu_\alpha(\fU_\alpha))$ is
quasi-isomorphic in $\fO_\fX(\fU_\alpha)\Contra$ to a bounded above
complex of contramodule $\fO_\fX(\fU_\alpha)$\+modules by
Corollary~\ref{Hom-tensor-B-preserves-acyclicity-boundedness-cor}(b).
 Any bounded above complex of contramodule
$\fO_\fX(\fU_\alpha)$\+modules is quasi-isomorphic
in $\fO_\fX(\fU_\alpha)\Contra$ to a bounded above complex of
cotorsion contramodule $\fO_\fX(\fU_\alpha)$\+modules by
Proposition~\ref{adic-raynaud-gruson-prop}(b).
 If two complexes in $\fO_\fX(\fU_\alpha)\Contra^\cot$ are isomorphic
in $\sD(\fO_\fX(\fU_\alpha)\Contra)$, then they are also isomorphic
in $\sD(\fO_\fX(\fU_\alpha)\Contra^\cot)$ in view of
Proposition~\ref{finite-co-resolutions}(b).
 Thus the complex $\fHom_{\fU_\alpha}(\fj_\alpha^*\cE^\bu,
\cK^\bu_\alpha)$ of locally cotorsion contraherent cosheaves of
contramodules on the affine Noetherian formal scheme $\fU_\alpha$
is isomorphic in $\sD(\fU_\alpha\Ctrh^\lct)$ to a bounded above
complex in $\fU_\alpha\Ctrh^\lct$.

 By
Corollary~\ref{lcth-lct-contramod-coflasque-derived-equivalences}(c),
the latter complex is isomorphic in $\sD^-(\fU_\alpha\Ctrh^\lct)$
to a bounded above complex in $\fU_\alpha\Ctrh^\lct_\cfq$.
 By the same corollary, it follows that the complex
$\fHom_{\fU_\alpha}(\fj_\alpha^*\cE^\bu,\cK^\bu_\alpha)$ in
$\fU_\alpha\Ctrh^\lct_\cfq$ (as per the proof of
Lemma~\ref{derived-fHom-takes-acyclics-to-acyclics}) is isomorphic
in $\sD(\fU_\alpha\Ctrh^\lct_\cfq)$ to a bounded above complex in
$\fU_\alpha\Ctrh^\lct_\cfq$.
 Using Corollary~\ref{co-flasque-direct-image-cor}(c),
we arrive to the conclusion that the complex
$\fHom_\fX(\cE^\bu,\cJ^\bu_\alpha)\simeq
\fj_\alpha{}_!\fHom_{\fU_\alpha}(\fj_\alpha^*\cE^\bu,\cK^\bu_\alpha)$
in $\fX\Ctrh^\lct$ is isomorphic in $\sD(\fX\Ctrh^\lct)$ to
a bounded above complex in $\fX\Ctrh^\lct$.

 Since $\sD^-(\fX\Ctrh^\lct)$ is a full triangulated subcategory
in $\sD(\fX\Ctrh^\lct)$, it finally follows that the complex
$\fHom_\fX(\cE^\bu,\cJ^\bu)$ is isomorphic in $\sD(\fX\Ctrh^\lct)$
to a bounded above complex in $\fX\Ctrh^\lct$.
 The case of a two-sided bounded complex $\cM^\bu$ in $\fX\Tors$
can be treated similarly.
\end{proof}

 Dual-analogously, if a complex $\fP^\bu$ in $\fX\Lcth^\lct_\bW$ is
bounded above, then it has a bounded above projective resolution
$\fF^\bu$, and the complex $\cE^\bu\ocn_\fX\fF^\bu$ in $\fX\Tors$ is
also bounded above.
 The case of a bounded below complex $\fP^\bu$ in $\fX\Lcth^\lct_\bW$
is a bit more complicated, and 
Lemma~\ref{co-flasque-open-subscheme-adjunction-morphism-lemma}(c)
is helpful in this context.

\begin{lem} \label{derived-contratensor-takes-bounded-to-bounded}
 For any bounded below (respectively, two-sided bounded) complex\/
$\fP^\bu$ in\/ $\fX\Lcth^\lct_\bW$, the derived category object\/
$\cB_\fX^\bu\ocn_\fX^\boL\fP^\bu\in\sD(\fX\Tors)$ can be represented
by a bounded below (resp., two-sided bounded) complex in\/ $\fX\Tors$.
\end{lem}

\begin{proof}
 Let $\fP^\bu$ be a bounded below complex in $\fX\Lcth^\lct_\bW$.
 Using
Corollary~\ref{lcth-lct-contramod-coflasque-derived-equivalences}(c),
one can find a bounded below complex $\fQ^\bu$ in $\fX\Ctrh^\lct_\cfq$
together with a quasi-isomorphism $\fQ^\bu\rarrow\fP^\bu$ of
complexes in $\fX\Lcth^\lct_\bW$.
 Suppose that the complex $\fQ^\bu$ is concentrated in
the cohomological degrees not smaller than some integer~$n$.
 Let $\fF^\bu$ be a complex in $\fX\Ctrh^\lct_\prj$ such that there is
a quasi-isomorphism $\fF^\bu\rarrow\fQ^\bu$ of complexes in
$\fX\Lcth^\lct_\bW$ (or equivalently, in $\fX\Ctrh^\lct_\cfq$).

 Consider the finite filtration $\fF^\bu_{\le\alpha}$ on $\fF^\bu$
constructed in the proof of
Lemma~\ref{derived-contratensor-takes-acyclics-to-acyclics}.
 Following that proof, one shows by induction on~$\alpha$ that
every complex $\fF^\bu_\alpha$ is quasi-isomorphic in $\fX\Ctrh^\lct$
to a complex of coflasque locally cotorsion contraherent cosheaves
of contramodules on $\fX$ concentrated in the cohomological
degrees~$\ge n$, while the complex $\fG^\bu_\alpha$ is quasi-isomorphic
in $\fU_\alpha\Ctrh^\lct$ to a complex of coflasque locally cotorsion
contraherent cosheaves of contramodules on $\fU_\alpha$ concentrated
in the cohomological degrees~$\ge n$.

 Specifically, the complex $\fG_1^\bu\simeq\fj_1^!\fF^\bu$ is
quasi-isomorphic to $\fj_1^!\fQ^\bu$, and the complex $\fF_1^\bu$ is
quasi-isomorphic to $\fQ_1^\bu=\fj_1{}_!\fj_1^!\fQ^\bu$.
 By Lemma~\ref{co-flasque-open-subscheme-adjunction-morphism-lemma}(c),
it follows that the complex $\fF^\bu_{\ge2}=\fF^\bu/\fF_1^\bu$ is
quasi-isomorphic to the complex $\fQ^\bu_{\ge2}=\fQ^\bu/\fQ_1^\bu$,
which is also a complex of coflasque locally cotorsion contraherent
cosheaves of contramodules on $\fX$ concentrated in the cohomological
degrees~$\ge n$.
 Then the complex $\fG_2^\bu$ is quasi-isomorphic to
$\fj_2^!\fQ^\bu_{\ge2}$, the complex $\fF_2^\bu$ is quasi-isomorphic
to $\fj_2{}_!\fj_2^!\fQ^\bu_{\ge2}$, etc.

 Now the complex of torsion $\fO_\fX(\fU_\alpha)$\+modules
$\cE^\bu(\fU_\alpha)\ot_{\fO(\fU_\alpha)}\fG^\bu_\alpha[\fU_\alpha]$
is quasi-isomorphic in $\fO_\fX(\fU_\alpha)\Tors$ to a bounded below
complex of torsion $\fO_\fX(\fU_\alpha)$\+modules by
Corollary~\ref{Hom-tensor-B-preserves-acyclicity-boundedness-cor}(d).
 Thus the complex $\fj_\alpha^*\cE^\bu\ocn_{\fU_\alpha}\fG^\bu_\alpha$
of quasi-coherent torsion sheaves on the affine Noetherian formal
scheme $\fU_\alpha$ is isomorphic in $\sD(\fX\Tors)$ to a bounded
below complex in $\fX\Tors$.

 By Corollary~\ref{qcoh-tors-flasque-dil-derived-equivalences}(b),
the latter complex is isomorphic in $\sD^+(\fX\Tors)$ to a bounded
below complex in $\fX\Tors^\fq$.
 By the same corollary, it follows that the complex
$\fj_\alpha^*\cE^\bu\ocn_{\fU_\alpha}\fG^\bu_\alpha$ in
$\fU_\alpha\Tors^\fq$ (as per the proof of
Lemma~\ref{derived-contratensor-takes-acyclics-to-acyclics}) is
isomorphic in $\sD(\fU_\alpha\Tors^\fq)$ to a bounded below complex
in $\fU_\alpha\Tors^\fq$.
 Using Corollary~\ref{co-flasque-direct-image-cor}(a), we arrive to
the conclusion that the complex $\cE^\bu\ocn_\fX\fF^\bu_\alpha\simeq
\fj_\alpha{}_*(\fj_\alpha^*\cE^\bu\ocn_{\fU_\alpha}\fG^\bu_\alpha)$
in $\fX\Tors$ is isomorphic in $\sD(\fX\Tors)$ to a bounded below
complex in $\fX\Tors$.

 Since $\sD^+(\fX\Tors)$ is a full triangulated subcategory in
$\sD(\fX\Tors)$, it finally follows that the complex
$\cE^\bu\ocn_\fX\fF^\bu$ is isomorphic in $\sD(\fX\Tors)$ to
a bounded below complex in $\fX\Tors$.
 The case of a two-sided bounded complex $\fP^\bu$ in
$\fX\Lcth^\lct_\bW$ can be treated similarly.
\end{proof}

 The proof of Main Theorem~\ref{fin-krulldim-mgm-main-theorem-I}
is finished.
\end{proof}

\subsection{MGM equivalence for formal schemes of finite
Krull dimension~II} \label{mgm-equiv-fin-krulldim-II-subsecn}
 The following theorem extends the result of
Main Theorem~\ref{fin-krulldim-mgm-main-theorem-I} to absolute
derived categories, in the spirit of~\cite[Theorem~4.8.1]{Pcosh},
\cite[Theorems~4.9 and~5.10]{Pmgm}, and~\cite[Corollary~8.7]{Ppt}.
 On top of the generalization of~\cite[Theorem~6.8.1]{Pcosh} to formal 
schemes, our theorem in this section also extends the latter result
from~\cite{Pcosh} to the absolute derived category symbols.

\begin{mainthm} \label{fin-krulldim-mgm-main-theorem-II}
 Let\/ $\fX$ be a (not necessarily semi-separated) Noetherian formal
scheme of finite Krull dimension, and let\/ $\cB_\fX^\bu\in
\sD^\bb(\fX\Tors)$ be the dedualizing complex of\/ $\fX$ as per
Section~\ref{dedualizing-complexes-on-formal-schemes-subsecn}.
 Let\/ $\st=\bb$, $+$, $-$, $\varnothing$, $\abs+$, $\abs-$,
or\/~$\abs$ be a conventional or absolute derived category symbol.
 Then there is a natural equivalence of triangulated categories
\begin{equation} \label{fin-krulldim-mgm-equivalence-eqn-II}
 \boR\fHom_\fX(\cB_\fX^\bu,{-})\:\sD^\st(\fX\Tors)\,\simeq\,
 \sD^\st(\fX\Ctrh):\!\cB_\fX^\bu\ocn_\fX^\boL{-}
\end{equation}
given by derived functors of the functor\/ $\fHom_\fX(\cB_\fX^\bu,{-})$
from Section~\ref{contraherent-Hom-subsecn} and the contratensor
product functor\/ $\cB_\fX^\bu\ocn_\fX{-}$ from
Section~\ref{contratensor-product-subsecn}.
\end{mainthm}

\begin{proof}
 Just as for the equivalence~\eqref{fin-krulldim-mgm-equivalence-eqn-I}
in Theorem~\ref{fin-krulldim-mgm-main-theorem-I}, the right-hand side
of the desired equivalence~\eqref{fin-krulldim-mgm-equivalence-eqn-II}
can be constructed in several alternative ways, because it is
the derived category of many exact categories naturally embedded
into one another, as per the commutative diagram of triangulated
equivalences~\eqref{derived-equivs-induced-by-inclusions-diag-I}.
 As in the proof of Theorem~\ref{fin-krulldim-mgm-main-theorem-I},
we prefer to work with the locally cotorsion $\bW$\+locally
contraherent cosheaves of contramodules.
 For any open covering $\bW$ of the formal scheme $\fX$, we will
construct the triangulated equivalences $\sD^\st(\fX\Tors)\simeq
\sD^\st(\fX\Lcth^\lct_\bW)$ for all the conventional and absolute
derived category symbols~$\st$.

 Having Theorem~\ref{fin-krulldim-mgm-main-theorem-I} and all
the lemmas proved in the arguments
in Section~\ref{mgm-equiv-fin-krulldim-I-subsecn} at our disposal,
the proof of Theorem~\ref{fin-krulldim-mgm-main-theorem-II} becomes
a straightforward application of the general technique
of~\cite[Appendix~B]{Pmgm} or~\cite[Appendix~A]{Pps}.

 In the notation of~\cite[Appendix~A]{Pps}, we put
\begin{alignat*}{4}
     \sA =\fX\Tors &&&=\sE &&\supset \sJ&&=\fX\Tors^\inj, \\
 \sB =\fX\Lcth^\lct_\bW &&&=\sF &&\supset \sP&&=\fX\Ctrh^\lct_\prj.
\end{alignat*}
 Consider the adjoint pair of DG\+functors
\begin{alignat*}{2}
 \Psi=\fHom_\fX(\cB_\fX^\bu,{-})&\:\sC^+(\sJ)&&\lrarrow\sC^+(\sB), \\
   \Phi=\cB_\fX^\bu\ocn_\fX{-}\,&\:\sC^-(\sP)&&\lrarrow\sC^-(\sA).
\end{alignat*}

 Whenever a complex $\cJ^\bu$ in $\fX\Tors^\inj$ is a coresolution
of an object of $\fX\Tors$, the complex $\Psi(\cJ^\bu)$ in
$\sC^+(\fX\Lcth^\lct_\bW)$ is concentrated in the cohomological
degrees~$\ge-l_2$ and quasi-isomorphic to a complex in
$\fX\Lcth^\lct_\bW$ concentrated in the cohomological degrees
$-l_2\le m\le l_1$, for suitable fixed integers $l_1$ and~$l_2$,
essentially by Lemma~\ref{derived-fHom-takes-bounded-to-bounded}.
 Whenever a complex $\fP^\bu$ in $\fX\Ctrh^\lct_\prj$ is
a resolution of an object of $\fX\Lcth^\lct_\bW$, the complex
$\Phi(\fP^\bu)$ in $\sC^-(\fX\Tors)$ is concentrated in
the cohomological degrees~$\le l_2$ and quasi-isomorphic to
a complex in $\fX\Tors$ concentrated in the cohomological degrees
$-l_1\le m\le l_2$ for suitable fixed integers $l_1$ and~$l_2$,
essentially by
Lemma~\ref{derived-contratensor-takes-bounded-to-bounded}.
 The discussion of truncations of complexes quasi-isomorphic to
bounded complexes in exact categories from~\cite[Section~4.1]{PRoos}
is relevant here.

 So the assumptions of~\cite[Sections~A.1 and~A.3]{Pps} are
satisfied, and the constructions of~\cite[Sections~A.2\+-A.3]{Pps}
provide the desired derived functors
\begin{equation} \label{absolute-derived-fHom-from-dedualizing}
 \boR\Psi=\boR\fHom_\fX(\cB^\bu_\fX,{-})\:\sD^\st(\fX\Tors)\lrarrow
 \sD^\st(\fX\Lcth^\lct_\bW)
\end{equation}
and
\begin{equation} \label{absolute-derived-contratensor-with-dedualizing}
 \boL\Phi=\cB^\bu_\fX\ocn_\fX^\boL{-}\,\:\sD^\st(\fX\Lcth^\lct_\bW)
 \lrarrow\sD^\st(\fX\Tors)
\end{equation}
for all the conventional and absolute derived category symbols
$\st=\bb$, $+$, $-$, $\varnothing$, $\abs+$, $\abs-$, or~$\abs$.
 According to~\cite[Section~A.4]{Pps}, the functor $\boR\Psi$ is
right adjoint to~$\boL\Phi$.

 Finally, the result of~\cite[first assertion of Theorem~A.5]{Pps}
allows to deduce the assertion that $\boR\Psi$ and $\boL\Phi$ are
mutually inverse equivalences from the same assertion in the case
of $\st=\bb$, which is covered by
Theorem~\ref{fin-krulldim-mgm-main-theorem-I}.
 Alternatively, the conditions in~\cite[second assertion of
Theorem~A.5]{Pps} are satisfied according to the arguments in the proof
of Theorem~\ref{fin-krulldim-mgm-main-theorem-I}, implying that
the functors $\boR\Psi$ \,\eqref{absolute-derived-fHom-from-dedualizing}
and $\boL\Phi$ \,\eqref{absolute-derived-contratensor-with-dedualizing}
are mutually inverse equivalences for all the conventional and absolute
derived category symbols.
\end{proof}

\subsection{MGM equivalence for semi-separated formal schemes}
\label{mgm-equiv-semi-separated-subsecn}
 The following MGM equivalence theorem is a version of
Theorem~\ref{fin-krulldim-mgm-main-theorem-I} for semi-separated
Noetherian formal schemes of possibly infinite Krull dimension.

\begin{mainthm} \label{semi-separated-mgm-main-theorem}
 Let\/ $\fX$ be a semi-separated Noetherian formal scheme, and let\/
$\cB_\fX^\bu\in\sD^\bb(\fX\Tors)$ be the dedualizing complex of\/ $\fX$
as per Section~\ref{dedualizing-complexes-on-formal-schemes-subsecn}.
 Let\/ $\st=+$ or\/~$\varnothing$ be a derived category symbol.
 Then there is a natural equivalence of triangulated categories
\begin{equation} \label{semi-separated-mgm-equivalence-eqn}
 \boR\fHom_\fX(\cB_\fX^\bu,{-})\:\sD^\st(\fX\Tors)\,\simeq\,
 \sD^\st(\fX\Ctrh):\!\cB_\fX^\bu\ocn_\fX^\boL{-}
\end{equation}
given by derived functors of the functor\/ $\fHom_\fX(\cB_\fX^\bu,{-})$
from Section~\ref{contraherent-Hom-subsecn} and the contratensor
product functor\/ $\cB_\fX^\bu\ocn_\fX{-}$ from
Section~\ref{contratensor-product-subsecn}.
\end{mainthm}

\begin{proof}
 The argument is similar to but much simpler than the proof of
Theorem~\ref{fin-krulldim-mgm-main-theorem-I}.
 Once again, the right-hand side of
the desired equivalence~\eqref{semi-separated-mgm-equivalence-eqn}
can be constructed in several alternative ways, as it is
the derived category of many exact categories naturally embedded
into one another, as per the commutative diagram of triangulated
equivalences~\eqref{derived-equivs-induced-by-inclusions-diag-II}.
 In order to prove the theorem, for any open covering $\bW$ of
the formal scheme $\fX$ one constructs a triangulated equivalence
$\sD(\fX\Tors)\simeq\sD(\fX\Lcth^\lct_\bW)$, and then shows that
it takes the full subcategories $\sD^+(\fX\Tors)\subset\sD(\fX\Tors)$
and $\sD^+(\fX\Lcth^\lct_\bW)\subset\sD(\fX\Lcth^\lct_\bW)$ into
one another.

 The major simplification as compared to the proof of
Theorem~\ref{fin-krulldim-mgm-main-theorem-I} is that, for any
affine open formal subscheme $\fU\subset\fX$, the open immersion
morphism $\fj\:\fU\rarrow\fX$ is affine now.
 Accordingly, the direct image functor $\fj_*\:\fU\Tors\rarrow\fX\Tors$
is exact, while the direct image functor $\fj_!\:\fU\Ctrh\rarrow
\fX\Ctrh$ is \emph{everywhere defined and} exact.
 For this reason, all the considerations of (co)flasqueness and all
the references to the results of Section~\ref{coflasque-subsecn},
which were needed in the proof of
Theorem~\ref{fin-krulldim-mgm-main-theorem-I} in order to deal with
the issue of nonexact/partially defined direct images from
the open immersions of affine open formal subschemes, are no longer
necessary in the context of the proof of our present theorem,
and can be simply dropped.

 In particular, one no longer needs to replace the dedualizing
complex $\cB_\fX^\bu$ with a quasi-isomorphic finite complex of
flasque quasi-coherent torsion sheaves $\cE^\bu$ for the construction
of the derived functors $\boR\fHom_\fX(\cB_\fX^\bu,{-})$ and
$\cB_\fX^\bu\ocn_\fX^\boL{-}$ in the context of our present theorem.
 Any finite complex $\cE^\bu$ in $\fX\Tors$ representing
the derived category object $\cB^\bu_\fX\in\sD^\bb(\fX\Tors)$ can be
plugged into the same construction as in the proof of
Theorem~\ref{fin-krulldim-mgm-main-theorem-I} in order to produce
the desired derived functors~\eqref{derived-fHom-from-dedualizing}
and~\eqref{derived-contratensor-with-dedualizing} in our present
context.

 Also, if one so prefers, one can replace all the references to
Proposition~\ref{contratensor-projection-formula-prop}(c)
in the proof of Theorem~\ref{fin-krulldim-mgm-main-theorem-I} with
references to part~(b) of the same proposition in the context of
the proof of Theorem~\ref{semi-separated-mgm-main-theorem}.
 Of course, in both the contexts one can also reference the more
general results of
Propositions~\ref{alf-contratensor-projection-formula-prop}
and/or~\ref{flat-contratensor-projection-formula-prop} instead
of the respective parts of
Proposition~\ref{contratensor-projection-formula-prop}.

 There are two reasons why our present theorem is only stated for
two derived category symbols~$+$ or~$\varnothing$, while the other
two conventional derived category symbols~$\bb$ and~$-$ appear to
be unavailable for MGM equivalence over infinite-dimensional
semi-separated Noetherian formal schemes.
 Firstly, the commutative diagram of triangulated
equivalences~\eqref{derived-equivs-induced-by-inclusions-diag-II},
based as it is on
Theorem~\ref{lcth-lcta-lct-unbounded-derived-equivalence}, only holds
for $\st=+$ or~$\varnothing$ in the case of infinite Krull dimension.
 Secondly, the proof of
Lemma~\ref{derived-fHom-takes-bounded-to-bounded} uses
Proposition~\ref{adic-raynaud-gruson-prop}(b), so it only works
for finite Krull dimension.
 Lemma~\ref{derived-contratensor-takes-bounded-to-bounded},
on the other hand, remains perfectly valid in the context of
Theorem~\ref{semi-separated-mgm-main-theorem}.

 Without the coflasqueness assumption, the convenient assertion of
Lemma~\ref{co-flasque-open-subscheme-adjunction-morphism-lemma}(c)
is no longer available for the purposes of the proof of
Lemma~\ref{derived-contratensor-takes-bounded-to-bounded},
but this is not a problem.
 The direct image functors are exact, and one can simply use
the cones of the adjunction morphisms instead of the cokernels in
the constructions of the complexes $\fQ^\bu_{\ge2}$,
\,$\fQ^\bu_{\ge3}$,~\dots, $\fQ^\bu_{\ge N}$.
 The cone of a morphism of bounded below complexes is still bounded
below, even if the bound shifts by~$1$ at every step.
 There is only a finite number $N-1$ of such steps in the situation
at hand.
\end{proof}

\begin{rem}
 In both the contexts of the proofs of
Theorems~\ref{fin-krulldim-mgm-main-theorem-I}
and~\ref{semi-separated-mgm-main-theorem}, one could try to avoid
the passage to the locally cotorsion $\bW$\+locally contraherent
cosheaves of contramodules altogether, and instead work with the derived
category of locally contraadjusted $\bW$\+locally contraherent cosheaves
of contramodules $\sD(\bW\Lcth_\bW)$ throughout the whole argument.
 This means that, instead of a complex of projective locally
cotorsion contraherent cosheaves of contramodules $\fF^\bu$ as in
the arguments above, one would have to apply the functor
$\cE^\bu\ocn_\fX{-}$ to a complex of projective locally contraadjusted
contraherent cosheaves of contramodules $\fF^\bu$ or (which is
generally more advisable) to a complex of flat or antilocally flat
contraherent cosheaves of contramodules $\fF^\bu$ on~$\fX$.

 Then the explicit classification provided by
Theorem~\ref{loc-Noetherian-projective-lct-cosh-contramods}(b), on which
the proof of Lemma~\ref{derived-contratensor-takes-acyclics-to-acyclics}
is based, becomes inapplicable, and one has to devise another argument
instead.
 Presumably, such an argument would require a proper undestanding of
the structure of complexes of flat or antilocally flat contraherent
cosheaves of contramodules \emph{that are acyclic as complexes of\/
$\bW$\+locally contraherent cosheaves of contramodules}, i.~e., of
complexes in $\fX\Ctrh^\fl$ or $\fX\Ctrh_\alf$ that are acyclic as
complexes in $\fX\Lcth_\bW$ or in $X\Ctrh$.

 In the case of a (nonformal) Noetherian scheme $X$ of finite Krull
dimension, such an argument is partly spelled out in~\cite[proof of
Theorem~6.10.1]{Pcosh}, based on the description of the complexes
in $X\Ctrh^\fl$ acyclic in $X\Ctrh$ provided
by~\cite[Theorem~6.9.4(b)]{Pcosh}.
 It looks likely that this approach could be extended to Noetherian
formal schemes $\fX$ of finite Krull dimension without great
difficulties.

 The analogue of~\cite[Theorem~6.9.4(b)]{Pcosh} for quasi-compact
semi-separated (nonformal) schemes $X$ can be found
in~\cite[Corollary~5.3.5(c)]{Pcosh}, which provides a description
of the complexes in $X\Ctrh_\alf$ that are acyclic in the exact
category of antilocal contraherent cosheaves $X\Ctrh_{\mathsf{al}}$
(or equivalently, in the exact category of contraherent
cosheaves $X\Ctrh$).
 Extending this result to semi-separated Noetherian formal schemes
$\fX$ looks like a harder problem.
 Solving it might open the way to deducing the assertion of
Theorem~\ref{semi-separated-mgm-main-theorem} for the derived
category symbols $\st=\bb$ and~$-$.

 To summarize, the proofs of the lemmas in
Section~\ref{mgm-equiv-fin-krulldim-I-subsecn} work by reducing
questions about a Noetherian formal scheme $\fX$ to similar questions
about its affine open formal subschemes $\fU_\alpha \subset\fX$ forming
a covering.
 However, the contraherent $\fHom(\cE,{-})$ and the contratensor
product $\cE\ocn{-}$ are essentially global constructions, having
\emph{no} compatibility with the restrictions of quasi-coherent
torsion sheaves or contraherent cosheaves of contramodules to open
subschemes.
 Instead, these functors have compatibility properties with
the direct images of quasi-coherent sheaves and contraherent cosheaves
of contramodules, as per the results of
Section~\ref{projection-formulas-subsecn}.

 For this reason, any local (or colocal) descriptions of the classes
of complexes under consideration in the proof of
Theorem~\ref{fin-krulldim-mgm-main-theorem-I} appear to be useless
of the purposes of checking the lemmas, and \emph{antilocal}
descriptions of the relevant classes of complexes (in the sense of
the discussion in the paper~\cite[Section~4]{Pal}) are needed instead.
 The classification theorems for injective quasi-coherent torsion
sheaves and locally cotorsion contraherent cosheaves of contramodules,
Theorems~\ref{loc-Noetherian-injective-qcoh-torsion}(b)
and~\ref{loc-Noetherian-projective-lct-cosh-contramods}(b),
are used in the proofs in
Section~\ref{mgm-equiv-fin-krulldim-I-subsecn} as a foundation for
an elementary ad hoc approach to establishing the requisite
antilocality.

 Other and more general techniques for proving the antilocality have
been developed in~\cite{Pal}, \cite[Appendix~B and Chapter~5]{Pcosh},
etc.\ (see also Corollaries~\ref{antilocally-flat-corollary}(c)
and~\ref{flat-contraherent-cosheaves-cor}(b) above), but
antilocality of the class of all complexes in $\fX\Ctrh_\alf$ that are
acyclic in $\fX\Ctrh$ (cf.~\cite[Example~6.9]{Pal}) appears to be
a complicated case, and the most difficult antilocality results are
proved in~\cite[Chapter~5]{Pcosh} using the passage from
the antilocal contraherent cosheaves to the contraadjusted
quasi-coherent sheaves together with the small object argument in
the Grothendieck category of complexes of quasi-coherent sheaves.
 Passing from contraherent cosheaves of contramodules to quasi-coherent
torsion sheaves, without a flatness or antilocal flatness assumption
on the cosheaves, seems like a harder task.
\end{rem}

\subsection{Compatibility with flat direct images}
\label{mgm-flat-direct-images-compat-subsecn}
 In this section we prove several versions
of~\cite[Theorem~7.5.1]{Pcosh} for Noetherian formal schemes.

\begin{thm} \label{ssep-mgm-flat-direct-images-compatible-thm}
 Let\/ $\ff\:\fY\rarrow\fX$ be a flat morphism of semi-separated
Noetherian formal schemes.
 Then the equivalences of triangulated categories\/ $\sD(\fY\Tors)
\simeq\sD(\fY\Ctrh)$ and\/ $\sD(\fX\Tors)\simeq\sD(\fX\Ctrh)$ from
Main Theorem~\textup{\ref{semi-separated-mgm-main-theorem}} transform
the right derived functor of direct image of quasi-coherent torsion
sheaves\/ $\boR\ff_*$ \,\eqref{qcoh-tors-derived-direct-image-formula}
into the left derived functor of direct image of contraherent cosheaves
of contramodules\/
$\boL\ff_!$~\,\eqref{ctrh-contramod-derived-direct-image-formula}.
 In other words, there is a commutative square diagram of triangulated
functors and triangulated equivalences
\begin{equation} \label{ssep-mgm-flat-direct-images-compat-diag}
\begin{gathered}
 \xymatrix{
  \sD(\fY\Tors) \ar@{=}[r] \ar[d]_{\boR\ff_*}
  & \sD(\fY\Ctrh) \ar[d]^{\boL\ff_!} \\
  \sD(\fX\Tors) \ar@{=}[r] & \sD(\fX\Ctrh)
 }
\end{gathered}
\end{equation}
\end{thm}

\begin{proof}
 This theorem can be viewed as a formal scheme version
of~\cite[Theorem~4.10.1]{Pcosh} (for flat morphisms~$\ff$).
 Taking the commutative
diagram~\eqref{ctrh-lcta-lct-der-dir-im-compat-diag} into account,
we need to prove that the triangulated equivalences
$\sD(\fY\Tors)\simeq\sD(\fY\Ctrh^\lct)$ and
$\sD(\fX\Tors)\simeq\sD(\fX\Ctrh^\lct)$ constructed in the proof of
Theorem~\ref{semi-separated-mgm-main-theorem} transform the derived
functor~\eqref{qcoh-tors-derived-direct-image-formula} into the derived
functor~\eqref{ctrh-lct-contramod-derived-direct-image-formula}.

 Let $\cE^\bu$ be any finite complex in $\fX\Tors$ representing
the derived category object $\cB_\fX^\bu\in\sD^\bb(\fX\Tors)$.
 Then, by
Corollary~\ref{flat-inverse-image-of-dedualizing-is-dedualizing},
the finite complex $\ff^*\cE^\bu$ in $\fY\Tors$ represents
the derived category object $\cB_\fY^\bu\in\sD^\bb(\fY\Tors)$.

 Let $\cJ^\bu$ be a complex of injective objects in $\fY\Tors$.
 Then, following the proof of
Theorem~\ref{semi-separated-mgm-main-theorem}, the complex
$\fHom_\fY(\ff^*\cE^\bu,\cJ^\bu)$ represents the derived category
object $\boR\fHom_\fY(\cB_\fY^\bu,\cJ^\bu)\in\sD(\fY\Ctrh^\lct)$.
 By Lemma~\ref{fHom-and-contratensor-with-alf-are-co-dilute}(a),
\,$\fHom_\fY(\ff^*\cE^\bu,\cJ^\bu)$ is a complex of codilute
locally cotorsion contraherent cosheaves of contramodules on~$\fY$.
 Following the construction of the derived functor $\boL\ff_!$
\,\eqref{ctrh-lct-contramod-derived-direct-image-formula} in
Section~\ref{derived-direct-images-of-ctrh-contramod-subsecn},
the complex $\ff_!\fHom_\fY(\ff^*\cE^\bu,\cJ^\bu)$ represents
the derived category object $\boL\ff_!\boR\fHom_\fY(\cB_\fY^\bu,\cJ^\bu)
\in\sD(\fX\Ctrh^\lct)$.

 On the other hand, the complex $\ff_*\cJ^\bu$ represents
the derived category object $\boR\ff_*\cJ^\bu\in\sD(\fX\Tors)$.
 Furthermore, $\ff_*\cJ^\bu$ is a complex in $\fX\Tors^\inj$ by
Lemma~\ref{injective-torsion-sheaves-direct-image-lem}.
 Following the proof of
Theorem~\ref{semi-separated-mgm-main-theorem}, the complex
$\fHom_\fX(\cE^\bu,\ff_*\cJ^\bu)$ represents the derived category
object $\boR\fHom_\fX(\cB_\fX^\bu,\boR\ff_*\cJ^\bu)
\in\sD(\fX\Ctrh^\lct)$.
 Now the isomorphism $\ff_!\fHom_\fY(\ff^*\cE^\bu,\cJ^\bu)\simeq
\fHom_\fX(\cE^\bu,\ff_*\cJ^\bu)$ provided by
Proposition~\ref{fHom-projection-formula-prop} represents the desired
isomorphism of triangulated functors establishing commutativity
of the diagram~\eqref{ssep-mgm-flat-direct-images-compat-diag}.

 Alternatively, let $\fF^\bu$ be a complex of projective objects in
$\fX\Ctrh^\lct$.
 Then, following the proof of
Theorem~\ref{semi-separated-mgm-main-theorem}, the complex
$\ff^*\cE^\bu\ocn_\fY\fF^\bu$ represents the derived category object
$\cB_\fY^\bu\ocn_\fY^\boL\fF^\bu\in\sD(\fY\Tors)$.
 By Lemma~\ref{fHom-and-contratensor-with-alf-are-co-dilute}(b),
$\ff^*\cE^\bu\ocn_\fY\fF^\bu$ is a complex of dilute quasi-coherent
torsion sheaves on~$\fY$.
 Following the construction of the derived functor $\boR\ff_*$
\,\eqref{qcoh-tors-derived-direct-image-formula} in
Section~\ref{derived-direct-images-of-qcoh-tors-subsecn}, the complex
$\ff_*(\ff^*\cE^\bu\ocn_\fY\fF^\bu)$ represents the derived category
object $\boR\ff_*(\cB_\fY^\bu\ocn_\fY^\boL\fF^\bu)\in\sD(\fX\Tors)$.
 
 On the other hand, the complex $\ff_!\fF^\bu$ represents the derived
category object $\boL\ff_!\fF^\bu\in\sD(\fX\Ctrh^\lct)$.
 Furthermore, $\ff_!\fF^\bu$ is a complex in $\fX\Ctrh^\lct_\prj$
by Corollary~\ref{proj-lct-direct-images}(b).
 Following the proof of
Theorem~\ref{semi-separated-mgm-main-theorem}, the complex
$\cE^\bu\ocn_\fX\ff_!\ff^\bu$ represents the derived category object
$\cB_\fX^\bu\ocn_\fX^\boL\boL\ff_!\fF^\bu\in\sD(\fX\Tors)$.
 Now the isomorphism $\ff_*(\ff^*\cE^\bu\ocn_\fY\fF^\bu)\simeq
\cE^\bu\ocn_\fX\ff_!\fF^\bu$ provided by
Proposition~\ref{alf-contratensor-projection-formula-prop} represents
the desired isomorphism of triangulated functors establishing 
commutativity of
the diagram~\eqref{ssep-mgm-flat-direct-images-compat-diag}.
\end{proof}

 The next theorem is a version of the previous one for Noetherian
formal schemes of finite Krull dimension instead of semi-separated
Noetherian formal schemes.

\begin{thm} \label{fin-krulldim-mgm-flat-direct-images-compatible-thm}
 Let\/ $\ff\:\fY\rarrow\fX$ be a flat morphism of Noetherian formal
schemes of finite Krull dimension.
 Then the equivalences of triangulated categories\/ $\sD(\fY\Tors)
\simeq\sD(\fY\Ctrh)$ and\/ $\sD(\fX\Tors)\simeq\sD(\fX\Ctrh)$ from
Main Theorem~\textup{\ref{fin-krulldim-mgm-main-theorem-I}} transform
the right derived functor of direct image of quasi-coherent torsion
sheaves\/ $\boR\ff_*$ \,\eqref{qcoh-tors-derived-direct-image-formula}
into the left derived functor of direct image of contraherent cosheaves
of contramodules\/
$\boL\ff_!$~\,\eqref{ctrh-contramod-derived-direct-image-formula}.
 In other words, there is a commutative square diagram of triangulated
functors and triangulated equivalences
\begin{equation} \label{fin-kdim-mgm-flat-direct-images-compat-diag}
\begin{gathered}
 \xymatrix{
  \sD(\fY\Tors) \ar@{=}[r] \ar[d]_{\boR\ff_*}
  & \sD(\fY\Ctrh) \ar[d]^{\boL\ff_!} \\
  \sD(\fX\Tors) \ar@{=}[r] & \sD(\fX\Ctrh)
 }
\end{gathered}
\end{equation}
\end{thm}

\begin{proof}
 This is the formal scheme generalization
of~\cite[Theorem~7.5.1]{Pcosh}, and our proof is
an improvement/simplification upon the one in~\cite{Pcosh}.
 Similarly to the proof of the previous theorem, taking the commutative
diagram~\eqref{ctrh-lcta-lct-der-dir-im-compat-diag} into account,
we need to prove that the triangulated equivalences
$\sD(\fY\Tors)\simeq\sD(\fY\Ctrh^\lct)$ and
$\sD(\fX\Tors)\simeq\sD(\fX\Ctrh^\lct)$ constructed in the proof of
Theorem~\ref{fin-krulldim-mgm-main-theorem-I} transform the derived
functor~\eqref{qcoh-tors-derived-direct-image-formula} into the derived
functor~\eqref{ctrh-lct-contramod-derived-direct-image-formula}.

 Let $\cE_\fX^\bu$ be a finite complex of flasque quasi-coherent
torsion sheaves on $\fX$ representing the derived category object
$\cB_\fX^\bu\in\sD^\bb(\fX\Tors)$.
 Then, by
Corollary~\ref{flat-inverse-image-of-dedualizing-is-dedualizing},
the finite complex $\ff^*\cE_\fX^\bu$ in $\fY\Tors$ represents
the derived category object $\cB_\fY^\bu\in\sD^\bb(\fY\Tors)$; but
$\ff^*\cE_\fX^\bu$ \emph{need not} be a complex of \emph{flasque}
quasi-coherent torsion sheaves on~$\fY$.
 Choose a finite complex $\cE_\fY^\bu$ of flasque quasi-coherent
torsion sheaves on $\fY$ together with a quasi-isomorphism
$\ff^*\cE_\fX^\bu\rarrow\cE_\fY^\bu$ of complexes in $\fY\Tors$.

 For any complex $\cJ^\bu$ in $\fY\Tors^\inj$, there is a natural
isomorphism of complexes of locally cotorsion contraherent cosheaves
of contramodules $\ff_!\fHom_\fY(\ff^*\cE_\fX^\bu,\cJ^\bu)\simeq
\fHom_\fX(\cE_\fX^\bu,\ff_*\cJ^\bu)$ on $\fX$ provided by
Proposition~\ref{fHom-projection-formula-prop}.
 Composing this isomorphism with the morphism
$\ff_!\fHom_\fY(\cE_\fY^\bu,\cJ^\bu)\rarrow
\ff_!\fHom_\fY(\ff^*\cE_\fX^\bu,\cJ^\bu)$ induced by
the quasi-isomorphism $\ff^*\cE_\fX^\bu\rarrow\cE_\fY^\bu$,
we obtain a natural morphism $\ff_!\fHom_\fY(\cE_\fY^\bu,\cJ^\bu)
\rarrow\fHom_\fX(\cE_\fX^\bu,\ff_*\cJ^\bu)$ of complexes of coflasque
locally cotorsion contraherent cosheaves of contramodules on~$\fX$.
 Here $\fHom_\fY(\cE_\fY^\bu,\cJ^\bu)$ and
$\fHom_\fX(\cE_\fX^\bu,\ff_*\cJ^\bu)$ are complexes of coflasque
locally cotorsion contraherent cosheaves of contramodules by
Lemma~\ref{fHom-from-flasque-contratensor-with-alf-co-flasque}(a),
and it follows that $\ff_!\fHom_\fY(\cE_\fY^\bu,\cJ^\bu)$ is
a complex of coflasque locally cotorsion contraherent cosheaves of
contramodules by Corollary~\ref{co-flasque-direct-image-cor}(c).
{\hbadness=1900\par}

 We observe that the complex $\ff_!\fHom_\fY(\cE_\fY^\bu,\cJ^\bu)$
represents the derived category object
$\boL\ff_!\boR\fHom_\fY(\cB_\fY^\bu,\cJ^\bu)\in\sD(\fX\Ctrh^\lct)$,
while the complex $\fHom_\fX(\cE_\fX^\bu,\ff_*\cJ^\bu)$ represents
the derived category object $\boR\fHom_\fX(\cB_\fX^\bu,\boR\ff_*\cJ^\bu)
\in\sD(\fX\Ctrh^\lct)$ (following the construction of
the triangulated equivalence in
Theorem~\ref{fin-krulldim-mgm-main-theorem-I} and the constructions
of the derived functors of direct
image~\eqref{qcoh-tors-derived-direct-image-formula}
and~\eqref{ctrh-lct-contramod-derived-direct-image-formula} in
Sections~\ref{derived-direct-images-of-qcoh-tors-subsecn}\+-%
\ref{derived-direct-images-of-ctrh-contramod-subsecn}).
 It is helpful to keep in mind that $\ff_*\cJ^\bu$ is a complex of
injective quasi-coherent torsion sheaves on~$\fX$.
 We have constructed a natural transformation
\begin{equation} \label{fHom-and-flat-direct-images-natural-transform}
 \boL\ff_!\,\boR\fHom_\fY(\cB_\fY^\bu,{-})\lrarrow
 \boR\fHom_\fX(\cB_\fX^\bu,\boR\ff_*({-})).
\end{equation}
of functors $\sD(\fY\Tors)\rarrow\sD(\fX\Ctrh^\lct)$.
 For the argument below, one needs to check that such natural
transformations are compatible with the compositions of flat
morphisms~$\ff$ of Noetherian formal schemes of finite Krull
dimension.

 Now one can cover the formal scheme $\fX$ with affine open formal
subschemes $\fU_\alpha$ and the formal scheme $\fY$ with affine open
formal subschemes $\fV_\beta$ in such a way that for every
index~$\beta$ there exists an index~$\alpha$ with
$\ff(\fV_\beta)\subset\fU_\alpha$.
 The triangulated category $\sD(\fY\Tors)$ is generated by
the derived direct images of objects from $\sD(\fV_\beta\Tors)$,
as one can see from the proof of
Lemma~\ref{derived-fHom-takes-acyclics-to-acyclics}.
 So it suffices to show that the natural
transformation~\eqref{fHom-and-flat-direct-images-natural-transform}
is an isomorphism for the open immersions $\fV_\beta\rarrow\fY$
and $\fU_\alpha\rarrow\fX$, and for the flat morphisms of affine
Noetherian formal schemes $\fV_\beta\rarrow\fU_\alpha$.
 Finally, the case of any flat morphism of semi-separated Noetherian
formal schemes $\ff\:\fV\rarrow\fU$ is covered by
Theorem~\ref{ssep-mgm-flat-direct-images-compatible-thm}.
 In the case of an open immersion of Noetherian formal schemes
$\fj\:\fY\rarrow\fX$, one can take $\cE_\fY^\bu=\fj^*\cE_\fX^\bu$,
and then the map~\eqref{fHom-and-flat-direct-images-natural-transform}
is an isomorphism by construction.

 Alternatively, for any complex $\fF^\bu$ in $\fY\Ctrh^\lct_\prj$,
there is a natural morphism of complexes of quasi-coherent torsion
sheaves $\cE_\fX^\bu\ocn_\fX\ff_!\fF^\bu\rarrow
\ff_*(\ff^*\cE_\fX^\bu\ocn_\fY\fF^\bu)$ on $\fY$ provided by
Proposition~\ref{contratensor-projection-formula-prop}(a).
 Composing this morphism with the morphism $\ff_*(\ff^*\cE_\fX^\bu
\ocn_\fY\fF^\bu)\rarrow\ff_*(\cE_\fY^\bu\ocn_\fY\fF^\bu)$ induced
by the quasi-isomorphism $\ff^*\cE_\fX^\bu\rarrow\cE_\fY^\bu$, we
obtain a natural morphism $\cE_\fX^\bu\ocn_\fX\ff_!\fF^\bu\rarrow
\ff_*(\cE_\fY^\bu\ocn_\fY\fF^\bu)$ of complexes of flasque
quasi-coherent torsion sheaves on~$\fX$.
 Here $\ff^*\cE_\fX^\bu\ocn_\fY\fF^\bu$ and
$\cE_\fX^\bu\ocn_\fX\ff_!\fF^\bu$ are complexes of flasque
quasi-coherent sheaves by 
Lemma~\ref{fHom-from-flasque-contratensor-with-alf-co-flasque}(b) or~(c)
(it is helpful to keep in mind that $\ff_!\fF^\bu$ is a complex of
projective locally cotorsion contraherent cosheaves of contramodules
on $\fX$ by Corollary~\ref{proj-lct-direct-images}(b)), and it follows
that $\ff_*(\cE_\fY^\bu\ocn_\fY\fF^\bu)$ is a complex of flasque
quasi-coherent torsion sheaves by
Corollary~\ref{co-flasque-direct-image-cor}(a).

 We observe that the complex $\ff_*(\cE_\fY^\bu\ocn_\fY\fF^\bu)$
represents the derived category object
$\boR\ff_*(\cB_\fY^\bu\ocn_\fY^\boL\fF^\bu)\in\sD(\fX\Tors)$, while
the complex $\cE_\fX^\bu\ocn_\fX\ff_!\fF^\bu$ represents the derived
category object $\cB_\fX^\bu\ocn_\fX^\boL\boL\ff_!\fF^\bu\in
\sD(\fX\Tors)$ (following the construction of the triangulated
equivalence in Theorem~\ref{fin-krulldim-mgm-main-theorem-I} and
the constructions of the derived functors of direct
image~\eqref{qcoh-tors-derived-direct-image-formula}
and~\eqref{ctrh-lct-contramod-derived-direct-image-formula} in
Sections~\ref{derived-direct-images-of-qcoh-tors-subsecn}\+-%
\ref{derived-direct-images-of-ctrh-contramod-subsecn}).
 We have constructed a natural transformation
\begin{equation} \label{ctrtnsr-and-flat-direct-images-nat-transform}
 \cB_\fX^\bu\ocn_\fX^\boL\boL\ff_!\fF^\bu\lrarrow
 \boR\ff_*(\cB_\fY^\bu\ocn_\fY^\boL\fF^\bu).
\end{equation}
of functors $\sD(\fY\Ctrh^\lct)\rarrow\sD(\fX\Tors)$.
 Once again, for the argument below, one needs to check that such
natural transformations are compatible with the compositions of flat
morphisms~$\ff$.

 Choose open coverings $\fX=\bigcup_\alpha\fU_\alpha$ and
$\fY=\bigcup_\beta\fV_\beta$ as above.
 The triangulated category $\sD(\fY\Ctrh^\lct)$ is generated by
the derived direct images of objects from $\sD(\fV_\beta\Ctrh^\lct)$,
as one can see from the proof of
Lemma~\ref{derived-contratensor-takes-acyclics-to-acyclics}.
 So it suffices to show that the natural
transformation~\eqref{ctrtnsr-and-flat-direct-images-nat-transform}
is an isomorphism for the open immersions $\fV_\beta\rarrow\fY$
and $\fU_\alpha\rarrow\fX$, and for the flat morphisms of affine
Noetherian formal schemes $\fV_\beta\rarrow\fU_\alpha$.
 Finally, the case of any flat morphism of semi-separated Noetherian
formal schemes $\ff\:\fV\rarrow\fU$ is covered by
Theorem~\ref{ssep-mgm-flat-direct-images-compatible-thm}.
 In the case of an open immersion of Noetherian formal schemes
$\fj\:\fY\rarrow\fX$, one can take $\cE_\fY^\bu=\fj^*\cE_\fX^\bu$.
 Then the map~\eqref{ctrtnsr-and-flat-direct-images-nat-transform}
is an isomorphism by
Proposition~\ref{flat-contratensor-projection-formula-prop}
and by construction.
\end{proof}

 To end, we formulate a more general version of the previous theorem
applicable to all conventional and absolute derived category symbols.

\begin{thm} \label{absderived-mgm-flat-direct-images-compatible-thm}
 Let\/ $\ff\:\fY\rarrow\fX$ be a flat morphism of Noetherian formal
schemes of finite Krull dimension, and let\/ $\st=\bb$, $+$, $-$,
$\varnothing$, $\abs+$, $\abs-$, or\/~$\abs$ be a conventional or
absolute derived category symbol.
 Then the equivalences of triangulated categories\/ $\sD^\st(\fY\Tors)
\simeq\sD^\st(\fY\Ctrh)$ and\/ $\sD^\st(\fX\Tors)\simeq
\sD^\st(\fX\Ctrh)$ from
Main Theorem~\textup{\ref{fin-krulldim-mgm-main-theorem-II}} transform
the right derived functor of direct image of quasi-coherent torsion
sheaves\/ $\boR\ff_*$ \,\eqref{qcoh-tors-derived-direct-image-formula}
into the left derived functor of direct image of contraherent cosheaves
of contramodules\/
$\boL\ff_!$~\,\eqref{ctrh-contramod-derived-direct-image-formula}.
 In other words, there is a commutative square diagram of triangulated
functors and triangulated equivalences
\begin{equation} \label{absderived-mgm-flat-direct-images-compat-diag}
\begin{gathered}
 \xymatrix{
  \sD^\st(\fY\Tors) \ar@{=}[r] \ar[d]_{\boR\ff_*}
  & \sD^\st(\fY\Ctrh) \ar[d]^{\boL\ff_!} \\
  \sD^\st(\fX\Tors) \ar@{=}[r] & \sD^\st(\fX\Ctrh)
 }
\end{gathered}
\end{equation}
\end{thm}

\begin{proof}
 Similarly to the proofs of the previous theorems, taking
the commutative diagram~\eqref{ctrh-lcta-lct-der-dir-im-compat-diag}
into account, one can work with the locally cotorsion contraherent
cosheaves of contramodules instead of the locally contraadjusted ones.
 The case of $\st=\varnothing$ is the result of
Theorem~\ref{fin-krulldim-mgm-flat-direct-images-compatible-thm}.
 The conventional derived categories $\sD^\bb(\sE)$, \,$\sD^+(\sE)$,
and $\sD^-(\sE)$ (for any exact category~$\sE$) are full subcategories
of the conventional unbounded derived category $\sD(\sE)$, so
the assertion for the conventional derived category symbols follows
from Theorem~\ref{fin-krulldim-mgm-flat-direct-images-compatible-thm}
(and the constructions of the derived equivalences in
Theorem~\ref{fin-krulldim-mgm-main-theorem-I}).
 Finally, the claimed assertion for all the conventional or absolute
derived category symbols follows from the case of $\st=\bb$ in view
of the constructions of the derived equivalences in
Theorem~\ref{fin-krulldim-mgm-main-theorem-II} and the constructions
of the derived functors~\eqref{qcoh-tors-derived-direct-image-formula}
and~\eqref{ctrh-contramod-derived-direct-image-formula}.
\end{proof}

\Section{Dualizing Complexes}  \label{dualizing-secn}

\subsection{Dualizing complexes for adic Noetherian rings}
\label{dualizing-complexes-for-adic-Noetherian-rings-subsecn}
 Let $R$ be a Noetherian commutative ring.
 A finite complex of injective $R$\+modules $D^\bu$ is said to be
a \emph{dualizing complex for~$R$} \,\cite[Section~V.2]{HarRD},
\cite[Section Tag~0A7A]{SP} if the following two conditions
are satisfied:
\begin{enumerate}
\renewcommand{\theenumi}{\roman{enumi}}
\item the cohomology $R$\+modules $H^i(D^\bu)$ are finitely
generated for all $i\in\boZ$;
\item the natural homothety map $R\rarrow\Hom_R(D^\bu,D^\bu)$
is an isomorphism.
\end{enumerate}
 Notice that, unlike many other expositions including~\cite{HarRD}
and~\cite{SP}, we prefer to include the two-sided boundedness and
the termwise injectivity condition into the definition of
a dualizing complex.

 Any Noetherian commutative ring admitting a dualizing complex has
finite Krull dimension~\cite[Corollary~V.7.2]{HarRD},
\cite[Lemma Tag~0A80]{SP}; in fact, the Krull dimension is bounded
by the length of a dualizing complex.
 The following locality lemma helps one feel more comfortable.

\begin{lem} \label{dualizing-complex-of-modules-local-property-lemma}
 Let $R\rarrow S_\alpha$, \,$1\le\alpha\le N$ be a open covering
(in the terminology of
Section~\ref{prelim-formal-open-immersions-coverings-subsecn}).
 Then a complex of $R$\+modules $D^\bu$ is a dualizing complex
(respectively, quasi-isomorphic to a dualizing complex) if and only
if the complex of $S_\alpha$\+modules $S\ot_RD^\bu$ is a dualizing
complex (resp., quasi-isomorphic to a dualizing complex) for
every index\/~$\alpha$.
\end{lem}

\begin{proof}
 This is~\cite[Corollary~V.2.3]{HarRD}
or~\cite[Lemmas Tags~0A7G and~0A7H]{SP}.
 The underlying fact that being an injective module over a Noetherian
commutative ring is a local property~\cite[Section~II.7]{HarRD}
is also used here.
\end{proof}

\begin{lem} \label{finite-map-Hom-complex-of-modules-dualizing}
 Let $R\rarrow S$ be a finite homomorphism of Noetherian rings
(i.~e., the ring $S$ is a finitely generated $R$\+module).
 Let $D^\bu$ be a dualizing complex for~$R$.
 Then the complex of $S$\+modules\/ $\Hom_R(S,D^\bu)$ is
a dualizing complex for~$S$.
\end{lem}

\begin{proof}
 This is~\cite[Proposition~V.2.4]{HarRD} or~\cite[Lemma Tag~0AX0]{SP}.
\end{proof}

 Now we pass to adic Noetherian rings.

\begin{lem} \label{dualizing-complex-of-torsion-modules-lemma}
 Let $R$ be an adic Noetherian ring and\/ $\cD^\bu$ be a finite complex
of injective torsion $R$\+modules.
 Then the following three conditions are equivalent:
\begin{enumerate}
\item there exists an ideal of definition $I\subset R$ such that
the complex of $R/I$\+modules\/ $\Hom_R(R/I,\cD^\bu)$ is
a dualizing complex for the ring~$R/I$; \
\item there exists an ideal of definition $I\subset R$ such that,
for every $n\ge1$, the complex of $R/I^n$\+modules\/
$\Hom_R(R/I^n,\cD^\bu)$ is a dualizing complex for the ring~$R/I^n$;
\item for every ideal of definition $J\subset R$, the complex of
$R/J$\+modules\/ $\Hom_R(R/J,\allowbreak\cD^\bu)$ is a dualizing
complex for the ring~$R/J$.
\end{enumerate}
\end{lem}

\begin{proof}
 (3)~$\Longrightarrow$~(1) Obvious.

 (1)~$\Longrightarrow$~(2) This is~\cite[Lemma~D.1.3]{Pcosh}.
 
 (2)~$\Longrightarrow$~(3) Follows immediately from
Lemma~\ref{finite-map-Hom-complex-of-modules-dualizing}.
\end{proof}

 Let $R$ be an adic Noetherian ring.
 A finite complex of injective torsion $R$\+modules $\cD^\bu$ is
said to be a \emph{dualizing complex of torsion
$R$\+modules}~\cite[Section~D.1]{Pcosh} if it satisfies any one
of the equivalent conditions of
Lemma~\ref{dualizing-complex-of-torsion-modules-lemma}.

 It follows from
Lemma~\ref{finite-map-Hom-complex-of-modules-dualizing} that if
$D^\bu$ is a dualizing complex for the Noetherian commutative ring
(viewed as an abstract ring with the topology disregarded),
then its maximal torsion subcomplex $\cD^\bu=\Gamma(D^\bu)\subset
D^\bu$ (in the notation of
Section~\ref{prelim-inj-proj-flat-torsion-contra-subsecn}) is
a dualizing complex of torsion $R$\+modules.

\begin{cor} \label{dualizing-complex-of-torsion-modules-local-cor}
 Let $R\rarrow S_\alpha$, \,$1\le\alpha\le N$ be a formal open
covering (in the terminology of
Section~\ref{prelim-formal-open-immersions-coverings-subsecn})
of an adic Noetherian ring $R$ by adic Noetherian rings~$S_\alpha$.
 Then a complex of torsion $R$\+modules $\cD^\bu$ is a dualizing
complex (respectively, quasi-isomorphic to a dualizing complex as
a complex in $R\Tors$) if and only if the complex of
$S_\alpha$\+modules $S\ot_RD^\bu$ is a dualizing complex
(resp., quasi-isomorphic to a dualizing complex as a complex in
$S_\alpha\Tors$) for every index\/~$\alpha$.
\end{cor}

\begin{proof}
 The assertion follows from the definitions and
Lemma~\ref{dualizing-complex-of-modules-local-property-lemma}.
 It is helpful to keep in mind that the functors of extension of
scalars $S_\alpha\ot_R{-}\,\:R\Tors\rarrow S_\alpha\Tors$ are exact
(by Lemma~\ref{flat-map-of-adic-rings-direct-image-lemma}(a)),
jointly faithful, and preserve, as well as jointly reflect, injectivity
(by Theorem~\ref{locality-of-injective-torsion-modules-thm}).
\end{proof}

\begin{lem} \label{adic-Noetherian-hom-tensor-flatness-injectivity}
 Let $R$ be an adically Noetherian adic topological ring (in
the sense of Section~\ref{adic-rings-subsecn}).
 In this setting: \par
\textup{(a)} for any injective torsion $R$\+modules\/ $\cI$ and\/
$\cJ\in R\Tors^\inj$, the quotseparated contramodule $R$\+module\/
$\Hom_R(\cI,\cJ)$ is flat and cotorsion; \par
\textup{(b)} for any injective torsion $R$\+module\/ $\cI$ and any
flat contramodule $R$\+module\/ $\fF$, the torsion $R$\+module\/
$\cI\ot_R\fF$ is injective.
\end{lem}

\begin{proof}
 This is our version of~\cite[Lemma~1.7.1]{Pcosh}.
 In part~(a), the assertion that $\Hom_R(\cI,\cJ)$ is a cotorsion
quotseparated contramodule $R$\+module is a special case of
Lemma~\ref{flat-cotorsion-injective-tensor-Hom}(c).
 All the other claims are covered by~\cite[Lemma~10.2]{Ppt}, which
also provides more general assertions than stated in the lemma.

 Alternatively, in the particular case when the ring $R$ is Noetherian,
it suffices to recall that a torsion $R$\+module is injective in
$R\Tors$ if and only if it is injective in $R\Modl$
(Lemma~\ref{injective-torsion-modules-characterizations}),
a contramodule $R$\+module is flat if and only if it is flat
as an $R$\+module (Lemma~\ref{flat-contramodules-characterizations}),
and a contramodule $R$\+module is cotorsion if and only if it is
cotorsion as an $R$\+module
(Lemma~\ref{cotorsion-contramodules-over-Noetherian}).
 So the desired assertions follow from~\cite[Lemmas~1.3.2(b)
and~1.7.1]{Pcosh}.
\end{proof}

\begin{lem} \label{adic-Noetherian-dualizing-adjunctions-quasi-isom}
 Let $R$ be an adic Noetherian ring and\/ $\cD^\bu$ be a dualizing
complex of torsion $R$\+modules.
 In this setting: \par
\textup{(a)} for any injective torsion $R$\+module\/ $\cJ$,
the natural morphism of finite complexes of injective torsion
$R$\+modules\/ $\cD^\bu\ot_R\Hom_R(\cD^\bu,\cJ)\rarrow\cJ$ is
a quasi-isomorphism of complexes in $R\Tors$ (in fact, a homotopy
equivalence of complexes in $R\Tors^\inj$); \par
\textup{(b)} for any flat contramodule $R$\+module\/ $\fF$,
the natural morphism of finite complexes of flat contramodule
$R$\+modules\/ $\fF\rarrow\Hom_R(\cD^\bu,\>\cD^\bu\ot_R\fF)$ is
a quasi-isomorphism of complexes in the abelian category $R\Contra$
(and also in the exact category $R\Contra_\fl$).
\end{lem}

\begin{proof}
 This is a part of~\cite[proof of Theorem~D.1.4]{Pcosh}.
 In part~(a), it is helpful to keep in mind that any quasi-isomorphism
of bounded below complexes of injective objects in an abelian
(or exact) category is a homotopy equivalence.
 In part~(b), it is similarly helpful to keep in mind that any
quasi-isomorphism in $R\Contra$ of bounded above complexes in
$R\Contra_\fl$ is a quasi-isomorphism in $R\Contra_\fl$ (because
$R\Contra_\fl$ is a resolving subcategory in $R\Contra$).
 The assertions of the lemma follow
from~\cite[Lemmas~C.4.1 and~D.1.2]{Pcosh}.
 In the proof of part~(b), one also needs to use the result of
Lemma~\ref{flat-quotseparated-contramodules-separated} above
(cf.~\cite[Lemma~D.1.5(b)]{Pcosh}).
 For generalizations, see~\cite[Theorems~D.5.10 and~E.2.9]{Pcosh}
and~\cite[Lemma~4.2 and Corollaries~11.4\+-11.5]{Ppt}.
\end{proof}

\subsection{Dualizing complexes on Noetherian formal schemes}
\label{dualizing-complexes-on-formal-schemes-subsecn}
 This very short section consists, essentially, of one definition.

 Let $\fX$ be a locally Noetherian formal scheme.
 As a preliminary remark we recall that, for every injective
quasi-coherent torsion sheaf $\cJ$ on $\fX$ and every affine open formal
subscheme $\fU\subset\fX$, the torsion $\fO_\fX(\fU)$\+module $\cJ(\fU)$
is injective (by Corollary~\ref{injective-torsion-local-flasque}(a)).
 The converse assertion also holds (by part~(d) of the same corollary).

 Let $\cD^\bu$ be a finite complex of injective quasi-coherent torsion
sheaves on~$\fX$.
 We will say that $\cD^\bu$ is a \emph{dualizing complex} for $\fX$ if,
for every affine open formal subscheme $\fU\subset\fX$, the complex
of torsion $\fO_\fX(\fU)$\+modules $\cD^\bu(\fU)$ is a dualizing
complex of torsion modules over the adic Noetherian
ring~$\fO_\fX(\fU)$ (in the sense of the definition in
Section~\ref{dualizing-complexes-for-adic-Noetherian-rings-subsecn}).
 In view of
Corollary~\ref{dualizing-complex-of-torsion-modules-local-cor}, it
suffices to check this condition for affine open formal subschemes
$\fU$ belonging to any chosen affine open covering of the locally
Noetherian formal scheme~$\fX$.

 So it follows from the definitions and the discussion in
Section~\ref{dualizing-complexes-for-adic-Noetherian-rings-subsecn}
that any locally Noetherian formal scheme admitting a dualizing
complex has finite Krull dimension.

 An equivalent definition and an illuminating discussion of dualizing
complexes on Noetherian formal schemes (including, in particular,
a proof of uniqueness of a dualizing complex) can be found
in~\cite[Section~5]{Yek0}.

 In the more general context of ind-Noetherian
ind-schemes~\cite[Sections~4.1\+-4.2]{Psemten}, dualizing complexes can
be unbounded (typically, unbounded below), and demonstrate complicated
behavior.
 In particular, a dualizing complex on an ind-Noetherian ind-scheme can
quite easily be acyclic~\cite[Remark~4.6 and Section~11(7)]{Psemten}.
 Dualizing complexes on Noetherian formal schemes are easiler to
comprehend.

\subsection{Co-contra correspondence}
\label{co-contra-subsecn}
 The following derived co-contra correspondence theorem is
the second main result of this paper.

 We refer to the introductions to the papers~\cite{Pmgm,Pps} for
a general discussion of the philosophy of co-contra correspondence,
including the roles of dedualizing complexes (as in the MGM equivalence
theorems, see Main Theorems~\ref{fin-krulldim-mgm-main-theorem-I},
\ref{fin-krulldim-mgm-main-theorem-II},
and~\ref{semi-separated-mgm-main-theorem} above) as compared to
the roles of dualizing complexes (as in the following main theorem).
 Let us emphasize that the MGM equivalence is an equivalence between
the \emph{conventional} or \emph{absolute} derived categories,
while the main theorem in the present section is an equivalence
between the \emph{coderived category} on one side and
the \emph{contraderived category} on the other side.

\begin{mainthm} \label{nonssep-derived-co-contra-main-thm}
 Let\/ $\fX$ be a (not necessarily semi-separated) Noetherian formal
scheme with a dualizing complex\/~$\cD^\bu$ (as defined in
Section~\ref{dualizing-complexes-on-formal-schemes-subsecn}).
 Then there is a natural equivalence of triangulated categories
\begin{equation} \label{nonssep-derived-co-contra-eqn}
 \boR\fHom_\fX(\cD^\bu,{-})\:\sD^\co(\fX\Tors)\,\simeq\,
 \sD^\ctr(\fX\Ctrh):\!\cD^\bu\ocn^\boL_\fX{-}
\end{equation}
given by derived functors of the functor\/ $\fHom_\fX(\cD^\bu,{-})$
from Section~\ref{contraherent-Hom-subsecn} and the contratensor
product functor\/ $\cD^\bu\ocn_\fX{-}$ from
Section~\ref{contratensor-product-subsecn}.
\end{mainthm}

\begin{proof}
 This is a formal scheme generalization of~\cite[Theorem~6.8.2]{Pcosh},
and the proof is similar.
 The right-hand side of
the equivalence~\eqref{nonssep-derived-co-contra-eqn} can be
constructed in several alternative ways, as it is the contraderived
category of many exact categories naturally embedded into one
another, as per the commutative diagram of triangulated
equivalences~\eqref{derived-equivs-induced-by-inclusions-diag-I}.
 Moreover, for all the exact categories appearing in this
commutative diagram, the Becker contraderived category coincides
with the Positselski one, according to
Theorems~\ref{lcth-lct-W-contramod-contraderived}
and~\ref{lcth-W-contramod-contraderived}.
 In particular, for any open covering $\bW$ of the formal scheme $\fX$,
we have $\sD^\bctr(\fX\Lcth^\lct_\bW)=\sD^\ctr(\fX\Lcth^\lct_\bW)\simeq
\sD^\ctr(\fX\Lcth_\bW)=\sD^\bctr(\fX\Lcth_\bW)$, as well as
$\sD^{\ctr=\bctr}(\fX\Ctrh^\lct)\simeq
\sD^{\ctr=\bctr}(\fX\Lcth^\lct_\bW)$ and
$\sD^{\ctr=\bctr}(\fX\Ctrh)\simeq\sD^{\ctr=\bctr}(\fX\Lcth_\bW)$.
 All of these triangulated categories are also naturally equivalent
to the homotopy categories $\Hot(\fX\Ctrh^\lct_\prj)$
and $\Hot(\fX\Ctrh_\prj)$, as well as to the conventional or absolute
derived category $\sD^{\abs=\varnothing}(\fX\Ctrh^\fl)$, as per
the same two theorems and
Lemma~\ref{finite-krulldim-ctrh-flat-all-derived-cats-agree}.

 The left-hand side of
the equivalence~\eqref{nonssep-derived-co-contra-eqn} can be
also constructed in some alternative ways,
as in Theorem~\ref{qcoh-tors-coderived} and
Corollary~\ref{qcoh-tors-flasque-dil-derived-equivalences}(a).
 In particular, the Becker coderived category coincides with
the Positselski one, $\sD^\bco(\fX\Tors)\simeq\sD^\co(\fX\Tors)$,
according to the mentioned theorem.
 The same Theorem~\ref{qcoh-tors-coderived} also provides
a triangulated equivalence $\Hot(\fX\Tors^\inj)\simeq
\sD^{\co=\bco}(\fX\Tors)$.

 Let us show that the functors $\fHom_\fX(\cD^\bu,{-})$ and
$\cD^\bu\ocn_\fX{-}$ are mutually inverse equivalences between
the homotopy categories $\Hot^\st(\fX\Tors^\inj)$ and
$\Hot^\st(\fX\Ctrh^\lct_\prj)$ for any symbol $\st=\bb$, $+$, $-$,
or~$\varnothing$.
 First of all, it needs to be checked that the two functors
land in the categories of complexes with the terms belonging to
the desired additive subcategories.

\begin{lem} \label{fHom-contratensor-with-dualizing-land-where-needed}
 Let\/ $\fX$ be a Noetherian formal scheme and\/ $\cI$ be an injective
quasi-coherent torsion sheaf on\/~$\fX$.
 Then \par
\textup{(a)} for any injective quasi-coherent torsion sheaf\/ $\cJ$
on\/ $\fX$, the contraherent cosheaf of contramodules\/
$\fHom_\fX(\cI,\cJ)$ on\/ $\fX$ is projective locally cotorsion; \par
\textup{(b)} for any projective locally cotorsion contraherent cosheaf
of contramodules\/ $\fF$ on\/ $\fX$, the quasi-coherent torsion sheaf\/
$\cI\ocn_\fX\fF$ on\/ $\fX$ is injective.
\end{lem}

\begin{proof}
 Let\/ $\fX=\bigcup_{\alpha=1}^N\fU_\alpha$ be a finite affine open
covering of~$\fX$.
 Denote by $\fj_\alpha\:\fU_\alpha\rarrow\fX$ the identity open
immersion morphisms.
 It is clear from Theorem~\ref{loc-Noetherian-injective-qcoh-torsion}(b)
that a quasi-coherent torsion sheaf $\cJ$ on $\fX$ is injective if
and only if it has the form
$\cJ=\bigoplus_{\alpha=1}^N\fj_\alpha{}_*\cK_\alpha$, where $\cK_\alpha$
are some injective quasi-coherent torsion sheaves on~$\fU_\alpha$.
 Dual-analogously, it follows from
Theorem~\ref{loc-Noetherian-projective-lct-cosh-contramods}(b) that
a cosheaf of contramodule $\fO_X$\+modules $\fF$ on $\fX$ is
projective locally cotorsion if and only it has the form
$\fF=\bigoplus_{\alpha=1}^N\fj_\alpha{}_!\fG_\alpha$, where $\fG_\alpha$
are some projective locally cotorsion contraherent cosheaves of
contramodules on~$\fU_\alpha$.

 Part~(a): for any quasi-coherent torsion sheaf $\cI$ and any injective
quasi-coherent torsion sheaf $\cJ$ on $\fX$, the cosheaf
$\fHom_\fX(\cI,\cJ)$ on $\fX$ is a locally cotorsion contraherent
cosheaf of contramodules by construction (see
Section~\ref{contraherent-Hom-subsecn}).
 We need to prove that the object $\fHom_\fX(\cI,\cJ)\in\fX\Ctrh^\lct$
is projective when both $\cI$ and $\cJ$ are injective in $\fX\Tors$.
 According to the previous paragraph, it suffices to consider
the case of $\cJ=\fj_*\cK$, where $\cK$ is an injective quasi-coherent
torsion sheaf on an affine open formal subscheme $\fU\subset\fX$
with the open immersion morphism $\fj\:\fU\rarrow\fX$.
 Then we have $\fHom_\fX(\cI,\fj_*\cK)\simeq
\fj_!\fHom_\fU(\fj^*\cI,\cK)$ by
Proposition~\ref{fHom-projection-formula-prop}.
 By construction, the locally cotorsion contraherent cosheaf of
contramodules $\fHom_\fU(\fj^*\cI,\cK)$ on $\fU$ corresponds to
the cotorsion contramodule $\fO_\fX(\fU)$\+module
$\Hom_{\fO_\fX(\fU)}(\cI(\fU),\cK(\fU))$ under the equivalence of
categories from
Lemma~\ref{lct-contrah-cosheaves-of-contramods-on-affine-lemma}.

 Now $\cK(\fU)$ is an injective torsion $\fO_\fX(\fU)$\+module
by assumption, and $\cI(\fU)$ is an injective torsion
$\fO_\fX(\fU)$\+module by
Corollary~\ref{local-section-criterion-of-injectivity-cor}.
 By Lemma~\ref{adic-Noetherian-hom-tensor-flatness-injectivity}(a),
it follows that $\Hom_{\fO_\fX(\fU)}(\cI(\fU),\cK(\fU))$ is
a flat cotorsion contramodule $\fO_\fX(\fU)$\+module.
 According to the proof of
Lemma~\ref{projective-contrah-lct-on-qcomp-ssep}, this means that
$\fHom_\fU(\fj^*\cI,\cK)$ is a projective locally cotorsion
contraherent cosheaf of contramodules on $\fU$, and it remains
to refer to Corollary~\ref{proj-lct-direct-images}(b) or to
the remarks in the first paragraph of this proof.

 Part~(b): according to the remarks in the first paragraph of
the proof, it suffices to consider the case of $\fF=\fj_!\fG$, where
$\fG$ is a projective locally cotorsion contraherent cosheaf on
an affine open formal subscheme $\fU\subset\fX$ with the open immersion
morphism $\fj\:\fU\rarrow\fX$.
 Then we have $\cI\ocn_\fX\fj_!\fG\simeq\fj_*(\fj^*\cI\ocn_\fU\fG)$
by Proposition~\ref{contratensor-projection-formula-prop}(c).
 By construction, the quasi-coherent torsoin sheaf
$\fj^*\cI\ocn_\fU\fG$ on $\fU$ corresponds to the torsion
$\fO_\fU(\fX)$\+module $\cI(\fU)\ot_{\fO_\fU(\fX)}\fG[\fU]$ under
the equivalence of categories $\fU\Tors\simeq\fO_\fX(\fU)\Tors$.

 Now $\fG[\fU]$ is a flat cotorsion contramodule
$\fO_\fX(\fU)$\+module according to the proof of
Lemma~\ref{projective-contrah-lct-on-qcomp-ssep}, and
$\cI(\fU)$ is an injective torsion $\fO_\fX(\fU)$\+module as
mentioned above.
 By Lemma~\ref{adic-Noetherian-hom-tensor-flatness-injectivity}(b),
it follows that $\cI(\fU)\ot_{\fO_\fU(\fX)}\fG[\fU]$ is an injective
torsion $\fO_\fX(\fU)$\+module.
 Hence $\fj^*\cI\ocn_\fU\fG$ is an injective quasi-coherent torsion
sheaf on~$\fU$, and it remains to refer to
Lemma~\ref{injective-torsion-sheaves-direct-image-lem}
or to the remarks in the first paragraph of this proof.
\end{proof}

 We have constructed the desired pair of functors
\begin{equation} \label{homotopy-fHom-from-dualizing}
 \fHom_\fX(\cD^\bu,{-})\:\Hot^\st(\fX\Tors^\inj)\lrarrow
 \Hot^\st(\fX\Ctrh^\lct_\prj)
\end{equation}
and
\begin{equation} \label{homotopy-contratensor-with-dualizing}
 \cD^\bu\ocn_\fX{-}\,\:\Hot^\st(\fX\Ctrh^\lct_\prj)\lrarrow
 \Hot^\st(\fX\Tors^\inj).
\end{equation}
 It is clear from the adjunction
isomorphism~\eqref{fHom-contratensor-adjunction-eqn} that
the functor~\eqref{homotopy-contratensor-with-dualizing} is
left adjoint to the functor~\eqref{homotopy-fHom-from-dualizing}.
 It remains to check that the adjunction morphisms are isomorphisms
in the homotopy categories.

\begin{lem} \label{formal-scheme-dualizing-adjunctions-quasi-isom}
 Let\/ $\fX$ be a Noetherian formal scheme with a dualizing
complex\/~$\cD^\bu$.
 Then \par
\textup{(a)} for any injective quasi-coherent torsion sheaf\/ $\cJ$
on\/ $\fX$, the adjunction morphism\/ $\cD^\bu\ocn_\fX
\fHom_\fX(\cD^\bu,\cJ)\rarrow\cJ$ is a quasi-isomorphism of
finite complexes in\/ $\fX\Tors$ (in fact, a homotopy equivalence
of finite complexes in\/ $\fX\Tors^\inj$); \par
\textup{(b)} for any projective locally cotorsion contraherent cosheaf
of contramodules\/ $\fF$ on\/ $\fX$, the adjunction morphism\/
$\fF\rarrow\fHom_\fX(\cD^\bu,\>\cD^\bu\ocn_\fX\fF)$ is
a quasi-isomorphism of finite complexes in\/ $\fX\Ctrh^\lct$ (in fact,
a homotopy equivalence of finite complexes in\/ $\fX\Ctrh^\lct_\prj$).
\end{lem}

\begin{proof} 
 Part~(a): in view of the first paragraph of the proof of
Lemma~\ref{fHom-contratensor-with-dualizing-land-where-needed},
it suffices to consider the case of $\cJ=\fj_*\cK$, where $\cK$ is
an injective quasi-coherent torsion sheaf on an affine open formal
subscheme $\fU\subset\fX$ with the open immersion morphism
$\fj\:\fU\rarrow\fX$.
 By Propositions~\ref{fHom-projection-formula-prop}
and~\ref{contratensor-projection-formula-prop}(c), we have
$$
 \cD^\bu\ocn_\fX\fHom_\fX(\cD^\bu,\fj_*\cK)\simeq
 \cD^\bu\ocn_\fX\fj_!\fHom_\fU(\fj^*\cD^\bu,\cK)\simeq
 \fj_*(\fj^*\cD^\bu\ocn_\fU\fHom_\fU(\fj^*\cD^\bu,\cK)).
$$
 Now the adjunction morphism $\cK\rarrow\fj^*\cD^\bu\ocn_\fU
\fHom_\fU(\fj^*\cD^\bu,\cK)$ of finite complexes of quasi-coherent
torsion sheaves on an affine Noetherian formal scheme $\fU$ is
a homotopy equivalence by
Lemma~\ref{adic-Noetherian-dualizing-adjunctions-quasi-isom}(a).
 The direct image functor~$\fj_*$ takes preserves homotopy
equivalences (as any additive functor does).

 Part~(b): in view of the first paragraph of the proof of
Lemma~\ref{fHom-contratensor-with-dualizing-land-where-needed},
it suffices to consider the case of $\fF=\fj_!\fG$, where $\fG$ is
a projective locally cotorsion contraherent cosheaf of
contramodules on~$\fU$.
 By Propositions~\ref{contratensor-projection-formula-prop}(c)
and~\ref{fHom-projection-formula-prop}, we have
$$
 \fHom_\fX(\cD^\bu,\,\cD^\bu\ocn_\fX\fj_!\fG)\simeq
 \fHom_\fX(\cD^\bu,\,\fj_*(\fj^*\cD^\bu\ocn_\fU\fG))\simeq
 \fj_!\fHom_\fU(\fj^*\cD^\bu,\>\fj^*\cD^\bu\ocn_\fU\fG).
$$
 Now the global cosections functor ${-}[\fU]$ from
Lemma~\ref{lct-contrah-cosheaves-of-contramods-on-affine-lemma}
takes the adjunction morphism $\fG\rarrow\fHom_\fU(\fj^*\cD^\bu,
\>\fj^*\cD^\bu\ocn_\fU\fG)$ of finite complexes of locally cotorsion
contraherent cosheaves of contramodules on an affine Noetherian formal
scheme $\fU$ to a quasi-isomorphism of finite complexes in the exact
category $\fO_\fX(\fU)\Contra_\fl$, by
Lemma~\ref{adic-Noetherian-dualizing-adjunctions-quasi-isom}(b).
 It remains to notice that the flat cotorsion contramodule
$\fO_\fX(\fU)$\+modules are the injective objects of
the exact category $\fO_\fX(\fU)\Contra_\fl$, by
Lemma~\ref{cotorsion-pair-kernel-as-injectives-projectives}(a).
 Any quasi-isomorphism of bounded below complexes of injective
objects in an exact category is a homotopy equivalence, so
the adjunction morphism $\fG\rarrow\fHom_\fU(\fj^*\cD^\bu,
\>\fj^*\cD^\bu\ocn_\fU\fG)$ is actually a homotopy equivalence
of finite complexes of locally cotorsion contraherent cosheaves of
contramodules on~$\fU$.
 Finally, the direct image functor~$\fj_!$, being an additive functor,
preserves homotopy equivalences.
\end{proof}

 The assertions that, for any complexes $\cJ^\bu$ in $\fX\Tors^\inj$
and $\fF^\bu$ in $\fX\Ctrh^\lct_\prj$, the adjunction morphisms
$\cD^\bu\ocn_\fX\fHom_\fX(\cD^\bu,\cJ^\bu)\rarrow\cJ^\bu$ and
$\fF^\bu\rarrow\fHom_\fX(\cD^\bu,\>\cD^\bu\ocn_\fX\fF^\bu)$ are
homotopy equivalences follow purely formally from the respective parts
of Lemma~\ref{formal-scheme-dualizing-adjunctions-quasi-isom}.
 The point is that the totalization of any finite acyclic complex of
complexes in an exact category is an (absolutely) acyclic complex.
 In particular, this applies to the split exact categories
$\fX\Tors^\inj$ and $\fX\Ctrh^\lct_\prj$.

 We have proved the theorem already, but let us offer an additional
discussion of a perhaps more informative alternative approach.
 Alternatively, one can prove that the functors
$\fHom_\fX(\cD^\bu,{-})$ and $\cD^\bu\ocn_\fX{-}$ induce mutually
inverse equivalences between the homotopy category
$\Hot^\st(\fX\Tors^\inj)$ and the (absolute) derived category
$\sD^\st(\fX\Ctrh^\fl)$, for any symbol $\st=\bb$, $+$, $-$,
or~$\varnothing$ (or~$\abs+$, $\abs-$, or~$\abs$).
 The argument is based on the following lemmas.

\begin{lem} \label{over-fin-krulldim-injective-contratensor-with-flat}
 Let\/ $\fX$ be a Noetherian formal scheme of finite Krull dimension.
 Then, for any injective quasi-coherent torsion sheaf\/ $\cI$ and
any flat contraherent cosheaf of contramodules $\fF$ on\/ $\fX$,
the quasi-coherent torsion sheaf\/ $\cI\ocn_\fX\fF$ on\/ $\fX$ is
injective.
\end{lem}

\begin{proof}
 By Corollary~\ref{flat-contraherent-cosheaves-cor}(b),
\,$\fF$ is a direct summand of a finitely iterated extension of
the direct images of flat contraherent cosheaves of contramodules
from affine open formal subschemes $\fU\subset\fX$.
 In view of Lemma~\ref{contratensor-exact-on-flats-and-alf}(a),
it suffices to show that, for any flat contraherent cosheaf of
contramodules $\fG$ on $\fU$, the quasi-coherent torsion sheaf
$\cI\ocn_\fX\fj_!\fG$ on $\fX$ is injective (where $\fj\:\fU\rarrow\fX$
is the open immersion morphism).

 Indeed, we have $\cI\ocn_\fX\fj_!\fG\simeq
\fj_*(\fj^*\cI\ocn_\fU\fG)$ by
Proposition~\ref{contratensor-projection-formula-prop}(c).
 The quasi-coherent torsion sheaf $\fj^*\cI\ocn_\fU\fG$ on $\fU$ is
injective by Corollary~\ref{local-section-criterion-of-injectivity-cor}
and Lemma~\ref{adic-Noetherian-hom-tensor-flatness-injectivity}(b),
and it remains to point out that the direct image functor $\fj_*\:
\fU\Tors\rarrow\fX\Tors$ takes injective objects to injective objects
by Lemma~\ref{injective-torsion-sheaves-direct-image-lem}.
\end{proof}

 It is clear from
Lemma~\ref{over-fin-krulldim-injective-contratensor-with-flat}
that the functor $\cD^\bu\ocn_\fX{-}$ takes all complexes in
$\fX\Ctrh^\fl$ to complexes in $\fX\Tors$.
 It also follows easily from
Lemma~\ref{contratensor-exact-on-flats-and-alf}(a) together with
Lemma~\ref{over-fin-krulldim-injective-contratensor-with-flat}
that the same functor takes acyclic complexes in $\fX\Ctrh^\fl$ to
contractible complexes in $\fX\Tors$.
 Hence we obtain the induced triangulated functor
\begin{equation} \label{induced-contratensor-with-dualizing-on-flats}
 \cD^\bu\ocn_\fX{-}\,\:\sD^\st(\fX\Ctrh^\fl)\lrarrow
 \Hot^\st(\fX\Tors^\inj).
\end{equation}
 It follows from the adjunction
isomorphism~\eqref{fHom-contratensor-adjunction-eqn} that
the functor~\eqref{homotopy-fHom-from-dualizing}, viewed as
a functor to the derived category $\sD^\st(\fX\Ctrh^\fl)$, is right
adjoint to
the functor~\eqref{induced-contratensor-with-dualizing-on-flats}.
 It remains to check that the related adjunction morphisms are
isomorphisms; and for the adjunction morphism in the homotopy category
$\Hot^\st(\fX\Tors)$ we know that already.
 The next lemma takes care of the adjunction morphism in
the derived category $\sD^\st(\fX\Ctrh^\fl)$.

\begin{lem} \label{formal-scheme-dualizing-adjunction-quasi-isom-II}
 Let\/ $\fX$ be a Noetherian formal scheme with a dualizing
complex\/~$\cD^\bu$.
 Then, for any flat contraherent cosheaf of contramodules\/ $\fF$
on\/ $\fX$, the adjunction morphism\/ $\fF\rarrow\fHom_\fX(\cD^\bu,
\>\cD^\bu\ocn_\fX\fF)$ is a quasi-isomorphism of finite complexes
in\/ $\fX\Ctrh^\fl$.
\end{lem}

\begin{proof}
 In any exact category, a morphism of complexes admitting a finite
filtration whose associated quotient morphisms are quasi-isomorphisms
is a quasi-isomorphism itself.
 Hence, in view of Corollary~\ref{flat-contraherent-cosheaves-cor}(b)
and Lemma~\ref{contratensor-exact-on-flats-and-alf}(a),
it suffices to consider the case of $\fF=\fj_!\fG$, where
$\fj\:\fU\rarrow\fX$ is the open immersion of an affine open
formal subscheme $\fU\subset\fX$ and $\fG$ is a flat contraherent
cosheaf of contramodules on~$\fU$.

 In this case, we argue similarly to the proof of
Lemma~\ref{formal-scheme-dualizing-adjunctions-quasi-isom}(b).
 Any complex of flat contraadjusted contramodule
$\fO_\fX(\fU)$\+modules that is acyclic as a complex of flat
contramodule $\fO_\fX(\fU)$\+modules is also acyclic as a complex of
flat contraadjusted contramodule $\fO_\fX(\fU)$\+modules, so it follows
from Lemma~\ref{adic-Noetherian-dualizing-adjunctions-quasi-isom}(b)
that the adjunction morphism $\fG\rarrow\fHom_\fU(\fj^*\cD^\bu,
\>\fj^*\cD^\bu\ocn_\fU\fG)$ is a quasi-isomorphism in $\fU\Ctrh^\fl$.
 It remains to point out the (existence and) exactness of
the direct image functor $\fj_!\:\fU\Ctrh^\fl\rarrow\fX\Ctrh^\fl$,
see Lemma~\ref{from-semi-separated-flat-direct-image} or
Corollary~\ref{flat-contrah-direct-images}(b).
\end{proof}

 Similarly to the argument above, the assertion that, for any
complex $\fF^\bu$ in $\fX\Ctrh^\fl$, the adjunction morphism
$\fF^\bu\rarrow\fHom_\fX(\cD^\bu,\>\cD^\bu\ocn_\fX\fF^\bu)$ is
an isomorphism in $\sD^\st(\fX\Ctrh^\fl)$ follows purely formally
from Lemma~\ref{formal-scheme-dualizing-adjunction-quasi-isom-II}.
 Once again, the totalization of any finite acyclic complex of
complexes in the exact category $\fX\Ctrh^\fl$ (as in any exact
category) is absolutely acyclic.
\end{proof}

\subsection{Weakly smooth morphisms of Noetherian formal schemes}
 The exposition in this and the next section follows the relevant
parts of~\cite[Chapter~10]{Psemten}, specialized from ind-Noetherian
ind-schemes to Noetherian formal schemes.
 However, unlike in~\cite[Chapter~10]{Psemten}, we do not presume our
formal schemes to be semi-separated.

 A homomorphism of commutative rings $f\:R\rarrow S$ is said to be
\emph{of finite type} if $S$ is a finitely generated commutative
$R$\+algebra.

\begin{lem} \label{finite-type-nilpotent-reduction}
 Let $f\:R\rarrow S$ be a homomorphism of commutative rings and
$I\subset R$ be a nilpotent ideal.
 Then the ring homomorphism~$f$ is of finite type if and only if
the induced ring homomorphism $\bar f\:R/I\rarrow S/Sf(I)$ is
of finite type.
\end{lem}

\begin{proof}
 The ``only if'' assertion is obvious.
 To prove the ``if'', pick a finite set of generators
$\bar x_1$,~\dots, $\bar x_m\in S/Sf(I)$ of the $R/I$\+algebra
$S/Sf(I)$, and lift them to some elements $x_1$,~\dots,  $x_m\in S$.
 Let $S'$ be the $R$\+subalgebra in $S$ spanned by
$x_1$,~\dots,~$x_m$.
 Then we have $S=S'+Sf(I)$.
 Denote by $M$ the $R$\+module $S/S'$; then $M=IM$, and it follows
that $M=0$.
 Thus $S'=S$ and $S$ is a finitely generated $R$\+algebra.
\end{proof}

 The following lemma uses the terminology of open coverings of
commutative rings; see
Section~\ref{prelim-formal-open-immersions-coverings-subsecn}.

\begin{lem} \label{finite-type-local-property}
\textup{(a)} Let $f\:R\rarrow S$ be a homomorphism of commutative rings
and $g_\alpha\:S\rarrow T_\alpha$ be an open covering.
 Then the ring homomorphism~$f$ is of finite type if and only if
the ring homomorphisms $g_\alpha\circ f$ are of finite type for
all the indices\/~$\alpha$. \par
\textup{(b)} Let $f\:R\rarrow S$ be a homomorphism of commutative rings
and $g_\alpha\:R\rarrow T_\alpha$ be an open covering.
 Then the ring homomorphism~$f$ is of finite type if and only if
the induced ring homomorphisms $f'_\alpha\:T_\alpha\rarrow 
T_\alpha\ot_RS$ are of finite type for all the indices\/~$\alpha$.
\end{lem}

\begin{proof}
 See~\cite[D\'efinition~I.6.2.1 and Proposition~I.6.2.5]{EGAI}
or~\cite[Lemma Tag~01T2]{SP}.
\end{proof}

\begin{cor} \label{adic-topological-rings-map-of-finite-type}
 Let $f\:R\rarrow S$ be a tight continuous homomorphism of adic
topological rings.
 Then the following conditions are equivalent:
\begin{enumerate}
\item there exists an ideal of definition $I\subset R$ such that
the ring homomorphism $\bar f\:R/I\rarrow S/Sf(I)$ is of finite type;
\item for every ideal of definition $I\subset R$, the ring
homomorphism $\bar f\:R/I\rarrow S/Sf(I)$ is of finite type;
\item for every pair of open ideals $I\subset R$ and $J\subset S$
such that $f(I)\subset J$, the ring homomorphism $\bar f\:R/I\rarrow
S/J$ is of finite type.
\end{enumerate}
\end{cor}

\begin{proof}
 (1)~$\Longleftrightarrow$~(2) Follows from
Lemma~\ref{finite-type-nilpotent-reduction}.

 (2)~$\Longleftrightarrow$~(3) Obvious.
\end{proof}

 We will say that a tight continuous homomorphism of adic topological
rings $f\:R\rarrow S$ is \emph{of finite type} (as a map of adic
topological rings) if it satisfies the equivalent conditions of
Corollary~\ref{adic-topological-rings-map-of-finite-type}.

\begin{cor} \label{topological-finite-type-local-property}
\textup{(a)} Let $f\:R\rarrow S$ be a tight continuous homomorphism
of adic topological rings and $g_\alpha\:S\rarrow T_\alpha$ be
a formal open covering.
 Then the topological ring map~$f$ is of finite type if and only if
the topological ring maps $g_\alpha\circ f$ are of finite type for
all the indices\/~$\alpha$. \par
\textup{(b)} Let $f\:R\rarrow S$ be a tight continuous homomorphism
of adic topological rings and $g_\alpha\:R\rarrow T_\alpha$ be
a formal open covering.
 Put $W_\alpha=T_\alpha\ot_RS$, and endow the rings $W_\alpha$ with
the tensor product topologies.
 Then the topological ring map~$f$ is of finite type if and only if
the induced topological ring maps $f'_\alpha\:T_\alpha\rarrow 
W_\alpha$ are of finite type for all the indices\/~$\alpha$.
\end{cor}

\begin{proof}
 Follows from the definitions and
Lemma~\ref{finite-type-local-property}.
\end{proof}

 Let $\ff\:\fY\rarrow\fX$ be a tight morphism of locally Noetherian
formal schemes.
 We will say that the morphism~$\ff$ is \emph{locally of finite
type} (cf.~\cite[Section~I.6.2]{EGAI}, \cite[Section~01T0]{SP}) if,
for every pair of affine open formal subschemes $\fU\subset\fX$ and
$\fV\subset\fY$ such that $\ff(\fV)\subset\fU$, the tight continuous
map of adic Noetherian rings $f_{\fU,\fV}\:\fO_\fX(\fU)\rarrow
\fO_\fY(\fV)$ is of finite type in the sense of the definition above.
 It follows from Corollary~\ref{topological-finite-type-local-property}
that a morphism~$\ff$ is of finite type if and only if there exist
two affine open coverings $\fX=\bigcup_\alpha\fU_\alpha$ and
$\fY=\bigcup_\beta\fV_\beta$, indexed by the same set of
indices~$\alpha$, such that $\ff(\fV_\alpha)\subset\fU_\alpha$ for every
index~$\alpha$ and the continuous ring map~$f_{\fU_\alpha,\fV_\alpha}$
is of finite type for every~$\alpha$.

 We will say that a morphism of locally Noetherian formal schemes
$\ff\:\fY\rarrow\fX$ is \emph{of finite type}
(cf.~\cite[Section~I.6.3]{EGAI}, \cite[Section~01T0]{SP}) if it is
locally of finite type and quasi-compact.

\medskip

 Let $f\:R\rarrow S$ be a homomorphism of commutative rings.
 By the \emph{fiber rings} of~$f$ we mean commutative rings of
the form $k\ot_RS$, where $k$~is a field and the $R$\+module
structure on~$k$ comes from a ring map $\kappa\:R\rarrow k$ such that
\emph{no proper subfield of~$k$ contains the image of\/~$\kappa$.}
 In other words, this means that $k$~is the residue field of a prime
ideal in~$R$.
 We are interested in ring homomorphisms~$f$ such that \emph{all
the fibers of~$f$ are regular Noetherian rings of bounded Krull
dimension}.

\begin{lem} \label{regular-fibers-local-property}
 Let $f\:R\rarrow S$ be a homomorphism of commutative rings and
$d$~be an integer. \par 
\textup{(a)} Let $g_\alpha\:S\rarrow T_\alpha$ be an open covering.
 Then the fiber rings of~$f$ are regular Noetherian rings of
Krull dimension\/~$\le\nobreak d$ if and only if the fiber rings
of $g_\alpha\circ f$ are regular Noetherian rings of Krull
dimension\/~$\le\nobreak d$ for all the indices\/~$\alpha$. \par
\textup{(b)} Let $g_\alpha\:R\rarrow T_\alpha$ be an open covering.
 Then the fiber rings of~$f$ are regular Noetherian rings of
Krull dimension\/~$\le\nobreak d$ if and only if the fiber rings of
the induced ring homomorphisms $f'_\alpha\:T_\alpha\rarrow 
T_\alpha\ot_RS$ are regular Noetherian rings of Krull
dimension\/~$\le\nobreak d$ for all the indices\/~$\alpha$.
\end{lem}

\begin{proof}
 Obvious.
\end{proof}

 Let $f\:R\rarrow S$ be a tight continuous homomorphism of adic
topological rings.
 In this context, by the \emph{fiber rings} of~$f$ we will mean
commutative rings of the form $k\ot_RS$, where $k$~is a residue
field of an \emph{open} prime ideal in~$R$.
 In other words, $k$~is a field with the discrete topology and
$\kappa\:R\rarrow k$ is a continuous ring homomorphism such that
no proper subfield of~$k$ contains the image of~$\kappa$.

 Let $I\subset R$ be an ideal of definition.
 Clearly, the fiber rings of the topological ring map
$f\:R\rarrow S$ are precisely the fiber rings of the induced
homomorphism of discrete commutative rings $\bar f\:R/I\rarrow S/Sf(I)$.

\begin{cor} \label{topological-regular-fibers-local-property}
 Let $f\:R\rarrow S$ be a tight continuous homomorphism of adic
toplogical rings and $d\ge0$ be an integer. \par 
\textup{(a)} Let $g_\alpha\:S\rarrow T_\alpha$ be a formal open
covering.
 Then the fiber rings of the topological ring map~$f$ are regular
Noetherian rings of Krull dimension\/~$\le\nobreak d$ if and only if
the fiber rings of the topological ring maps $g_\alpha\circ f$ are
regular Noetherian rings of Krull dimension\/~$\le\nobreak d$ for all
the indices\/~$\alpha$. \par
\textup{(b)} Let $g_\alpha\:R\rarrow T_\alpha$ be a formal open
covering.
 Put $W_\alpha=T_\alpha\ot_RS$, and endow the rings $W_\alpha$ with
the tensor product topologies.
 Then the fiber rings of the topological ring map~$f$ are regular
Noetherian rings of Krull dimension\/~$\le\nobreak d$ if and only if
the fiber rings of the induced topological ring maps
$f'_\alpha\:T_\alpha\rarrow W_\alpha$ are regular Noetherian rings of
Krull dimension\/~$\le\nobreak d$ for all the indices\/~$\alpha$.
\end{cor}

\begin{proof}
 Follows from Lemma~\ref{regular-fibers-local-property}.
\end{proof}

 Let $\ff\:\fY\rarrow\fX$ be a tight morphism of locally Noetherian
formal schemes, and let $d\ge0$ be an integer.
 We will say that the morphism~$\ff$ has \emph{regular fibers of Krull
dimension\/~$\le\nobreak d$} if, for every pair of affine open formal
subschemes $\fU\subset\fX$ and $\fV\subset\fY$ such that $\ff(\fV)
\subset\fU$, the fiber rings of the tight continuous map of adic
Noetherian rings $f_{\fU,\fV}\:\fO_\fX(\fU)\rarrow\fO_\fY(\fV)$ are
regular Noetherian rings of Krull dimension~$\le\nobreak d$.
 It follows from
Corollary~\ref{topological-regular-fibers-local-property}
that a morphism~$\ff$ has regular fibers of Krull
dimension~$\le\nobreak d$ if and only if there exist two affine open
coverings $\fX=\bigcup_\alpha\fU_\alpha$ and $\fY=
\bigcup_\beta\fV_\beta$, indexed by the same set of indices~$\alpha$,
such that $\ff(\fV_\alpha)\subset\fU_\alpha$ for every index~$\alpha$
and the fiber rings of the tight continuous ring
map~$f_{\fU_\alpha,\fV_\alpha}$ are regular Noetherian rings of
Krull dimension~$\le\nobreak d$ for every~$\alpha$.

 We will say that a continuous ring map $f\:R\rarrow S$ is
\emph{weakly smooth} if $f$~is tight and flat with regular fibers of
bounded Krull dimension (i.~e., of Krull dimension~$\le\nobreak d$ for
some finite integer~$d$).
 Similarly, a morphism of locally Noetherian formal schemes
$\ff\:\fY\rarrow\fX$ is said to be \emph{weakly smooth} if $\ff$~is
tight and flat with regular fibers of bounded Krull dimension.

\subsection{Inverse image of dualizing complex under weakly smooth
morphism}
 In this section we spell out the formulations and proofs of
the formal scheme versions (essentially, particular cases)
of~\cite[Lemmas~10.20 and~10.21]{Psemten}, based on the results
of~\cite[Propositions~10.2(b) and~10.14, and Lemma~10.15]{Psemten}.

\begin{lem} \label{extension-of-scalars-injective-dimension}
 Let $R$ and $S$ be Noetherian commutative rings, and let
$f\:R\rarrow S$ be ring homomorphism of finite type such that $S$
is a flat $R$\+module and all the fiber rings of~$f$ are regular of
Krull dimension\/~$\le\nobreak d$.
 Let $J$ be an injective $R$\+module.
 Then the injective dimension of the $S$\+module $S\ot_RJ$ does not
exceed~$d$.
\end{lem}

\begin{proof}
 Obviously, the underlying $R$\+module of the $S$\+module $S\ot_RJ$ is
injective (as the tensor product of a flat and an injective module
over a Noetherian commutative ring).
 Now we claim that, for any $S$\+module $K$, if the underlying
$R$\+module of $K$ is injective, then the injective dimension of $K$
as an $S$\+module does not exceed~$d$.
 This is a particular case of~\cite[Proposition~10.2(b)]{Psemten}.
\end{proof}

\begin{lem} \label{adic-extension-of-scalars-injective-dimension}
 Let $R$ and $S$ be adic Noetherian rings, and let $f\:R\rarrow S$ be
a flat tight continuous ring map of finite type with regular fiber
rings of Krull dimension\/~$\le\nobreak d$.
 Let\/ $\cJ$ be an injective torsion $R$\+module.
 Then the injective dimension of the torsion $S$\+module
$f^*(\cJ)=S\ot_R\cJ$ does not exceed~$d$.
\end{lem}

\begin{proof}
 Let $0\rarrow S\ot_R\cJ\rarrow\cK^0\rarrow\cK^1\rarrow\dotsb
\rarrow\cK^{d-1}\rarrow\cK\rarrow0$ be an exact sequence in
$S\Modl_\tors$ with injective torsion $S$\+modules $\cK^i$,
\,$0\le i\le d-1$.
 We need to show that the torsion $S$\+module $\cK$ is injective,
and for this purpose, it suffices to check that the $S/Sf(I)$\+module
$\Hom_S(S/Sf(I),\cK)\simeq\Hom_R(R/I,\cK)$ is injective for every ideal
of definition $I\subset R$.
 Clearly, we have $\Hom_R(R/I,\>S\ot_R\nobreak\cJ)\allowbreak\simeq
S\ot_R\Hom_R(R/I,\cJ)$.
 As a special case of~\cite[Lemma~10.15 or~10.18]{Psemten}, the complex
of $S/Sf(I)$\+modules $0\rarrow S\ot_R\Hom_R(R/I,\cJ)\rarrow
\Hom_R(R/I,\cK^0)\rarrow\Hom_R(R/I,\cK^1)\rarrow\dotsb\rarrow
\Hom_R(R/I,\cK^{d-1})\rarrow\Hom_R(R/I,\cK)\rarrow0$ is acyclic.
 The $S/Sf(I)$\+modules $\Hom_R(R/I,\cK^i)$ are injective for all
$0\le i\le d-1$.
 By Lemma~\ref{extension-of-scalars-injective-dimension} applied to
the ring homomorphism $\bar f\:R/I\rarrow S/Sf(I)$ and
the $R/I$\+module $\Hom_R(R/I,\cJ)$, it follows that the $S$\+module
$\Hom_R(R/I,\cK)$ is injective.
\end{proof}

\begin{prop} \label{adic-smooth-inverse-image-of-dualizing-prop}
 Let $R$ and $S$ be adic Noetherian rings, and let $f\:R\rarrow S$ be
a weakly smooth continuous ring map of finite type.
 Let\/ $\cD_R^\bu$ be a dualizing complex of torsion $R$\+modules.
 Then \par
\textup{(a)} there exists a finite complex of injective torsion
$S$\+modules\/ $\cD_S^\bu$ together with a quasi-isomorphism
$S\ot_R\cD_R^\bu\rarrow\cD_S^\bu$ of finite complexes in
$S\Modl_\tors$; \par
\textup{(b)} any such complex\/ $\cD_S^\bu$ is a dualizing complex of
torsion $S$\+modules.
\end{prop}

\begin{proof}
 Part~(a) follows from
Lemma~\ref{adic-extension-of-scalars-injective-dimension}.
 Part~(b) is a special case of~\cite[Lemma~10.21]{Psemten}.
 Notice that, in the context of~\cite{Psemten}, the dualizing
complexes on ind-Noetherian ind-schemes are usually unbounded below,
and quite often acyclic.
 For this reason, they have to be considered as objects of
the \emph{coderived category} of quasi-coherent torsion sheaves,
i.~e., viewed up to morphisms of quasi-coherent torsion sheaves with
\emph{coacyclic} cones.
 Hence the role of the morphisms with coacyclic cones
in~\cite[Lemma~10.21]{Psemten} and the other cited results
from~\cite{Psemten}.
 However, in our present context, all the dualizing complexes are
finite; and a bounded below complex is acyclic if and only if it
is coacyclic~\cite[Lemma~2.1]{Psemi}, \cite[Lemma~A.1.3]{Pcosh}.
 So we simply speak of quasi-isomorphisms instead of morphisms with
coacyclic cones.
\end{proof}

\begin{cor} \label{formsch-smooth-inverse-image-of-dualizing-cor}
 Let\/ $\fX$ and\/ $\fY$ be locally Noetherian formal schemes, and
let\/ $\ff\:\fY\rarrow\fX$ be a weakly smooth morphism locally of
finite type.
 Let\/ $\cD_\fX^\bu$ be a dualizing complex on\/~$\fX$.
 Then \par
\textup{(a)} there exists a finite complex of injective quasi-coherent
torsion sheaves\/ $\cD_\fY^\bu$ on\/ $\fY$ together with
a quasi-isomorphism\/ $\ff^*\cD_\fX^\bu\rarrow\cD_\fY^\bu$ of finite
complexes in\/ $\fY\Tors$; \par
\textup{(b)} any such complex\/ $\cD_\fY^\bu$ is a dualizing complex 
on\/~$\fY$.
\end{cor}

\begin{proof}
 Part~(a) follows from
Lemma~\ref{adic-extension-of-scalars-injective-dimension}.
 Part~(b) follows from
Proposition~\ref{adic-smooth-inverse-image-of-dualizing-prop}(b).
\end{proof}

\subsection{Compatibility with smooth direct images}
\label{nonssep-co-contra-smooth-direct-images-compat-subsecn}
 The following theorem is the formal scheme generalization
of~\cite[Theorem~7.6.3]{Pcosh}.

\begin{thm} \label{nonssep-co-contra-smooth-direct-images-compat-thm}
 Let\/ $\fX$ and\/ $\fY$ be Noetherian formal schemes, and let\/
$\ff\:\fY\rarrow\fX$ be a weakly smooth morphism of finite type.
 Let\/ $\cD_\fX^\bu$ be a dualizing complex for\/ $\fX$ and\/
$\cD_\fY^\bu$ be the corresponding dualizing complex for\/ $\fY$ as per
Corollary~\textup{\ref{formsch-smooth-inverse-image-of-dualizing-cor}}.
 Then the related equivalences of triangulated categories\/
$\sD^\co(\fY\Tors)\simeq\sD^\ctr(\fY\Ctrh)$ and\/ $\sD^\co(\fX\Tors)
\simeq\sD^\ctr(\fX\Ctrh)$ from
Main Theorem~\textup{\ref{nonssep-derived-co-contra-main-thm}} transform
the right derived functor of direct image of quasi-coherent torsion
sheaves\/ $\boR\ff_*$ \,\eqref{qcoh-tors-derived-direct-image-formula}
into the left derived functor of direct image of contraherent cosheaves
of contramodules\/
$\boL\ff_!$~\,\eqref{ctrh-contramod-derived-direct-image-formula}.
 In other words, there is a commutative square diagram of triangulated
functors and triangulated equivalences
\begin{equation} \label{co-contra-smooth-direct-images-compat-diag}
\begin{gathered}
 \xymatrix{
  \sD^\co(\fY\Tors) \ar@{=}[r] \ar[d]_{\boR\ff_*}
  & \sD^\ctr(\fY\Ctrh) \ar[d]^{\boL\ff_!} \\
  \sD^\co(\fX\Tors) \ar@{=}[r] & \sD^\ctr(\fX\Ctrh)
 }
\end{gathered}
\end{equation}
\end{thm}

\begin{proof}
 Similarly to the proofs of the theorems in
Section~\ref{mgm-flat-direct-images-compat-subsecn}, taking
the commutative diagram~\eqref{ctrh-lcta-lct-der-dir-im-compat-diag}
into account, we need to prove that the triangulated equivalences
$\sD^\co(\fY\Tors)\simeq\sD^\ctr(\fY\Ctrh^\lct)$ and
$\sD^\co(\fX\Tors)\simeq\sD^\ctr(\fX\Ctrh^\lct)$ constructed in
the proof of Theorem~\ref{nonssep-derived-co-contra-main-thm} transform
the derived functor~\eqref{qcoh-tors-derived-direct-image-formula}
into the derived
functor~\eqref{ctrh-lct-contramod-derived-direct-image-formula}.
 The latter equivalences were constructed in terms of the injective
and projective resolutions as the equivalences of (bounded or
unbounded) homotopy categories $\Hot^\st(\fY\Tors^\inj)\simeq
\Hom^\st(\fY\Ctrh^\lct_\prj)$ provided by the functors
$\fHom_\fY(\cD_\fY^\bu,{-})$ and $\cD_\fY^\bu\ocn_\fY{-}$, and
similarly for~$\fX$.
 In particular, it was shown in
Lemma~\ref{fHom-contratensor-with-dualizing-land-where-needed} that
such functors of contraherent $\fHom$ and contraherensor product
take complexes in $\fY\Tors^\inj$ to complexes in $\fY\Ctrh^\lct_\prj$
and back.

 Furthermore, following the constructions in
Sections~\ref{derived-direct-images-of-qcoh-tors-subsecn}\+-
\ref{derived-direct-images-of-ctrh-contramod-subsecn}, we notice that
any complexes in $\fY\Tors^\inj$ and $\fY\Ctrh^\lct_\prj$ are
adjusted to the derived functors $\boR\ff_*$ and $\boL\ff_!$
acting between the co/contraderived categories.
 Since the morphism~$\ff$ is tight and flat, the direct image functors
$\ff_*$ and~$\ff_!$ take $\fY\Tors^\inj$ into $\fX\Tors^\inj$ and
$\fY\Ctrh^\lct_\prj$ into $\fX\Ctrh^\lct_\prj$ (by
Corollary~\ref{proj-lct-direct-images}(b)).

 For any complex $\cJ^\bu$ in $\fY\Tors^\inj$, there is a natural
isomorphism of complexes of locally cotorsion contraherent cosheaves
of contramodules $\ff_!\fHom_\fY(\ff^*\cD_\fX^\bu,\cJ^\bu)\simeq
\fHom_\fX(\cD_\fX^\bu,\ff_*\cJ^\bu)$ provided by
Proposition~\ref{fHom-projection-formula-prop}.
 Composing this isomorphism with the morphism
$\ff_!\fHom_\fY(\cD_\fY^\bu,\cJ^\bu)\rarrow
\ff_!\fHom_\fY(\ff^*\cD_\fX^\bu,\cJ^\bu)$ induced by
the quasi-isomorphism $\ff^*\cD_\fX^\bu\rarrow\cD_\fY^\bu$, we
obtain a natural morphism
\begin{equation} \label{fHom-and-smooth-direct-images-natural-transform}
 \ff_!\fHom_\fY(\cD_\fY^\bu,\cJ^\bu)\lrarrow
 \fHom_\fX(\cD_\fX^\bu,\ff_*\cJ^\bu)
\end{equation}
of complexes in $\fX\Ctrh^\lct_\prj$.
 Dual-analogously, by
Proposition~\ref{contratensor-projection-formula-prop}(a),
for any complex $\fF^\bu$ in $\fY\Ctrh^\lct_\prj$ there is a natural
morphism of complexes of quasi-coherent torsion sheaves
$\cD_\fX^\bu\ocn_\fX\ff_!\fF^\bu\rarrow
\ff_*(\ff^*\cD_\fX^\bu\ocn_\fY\fF^\bu)$ on~$\fX$.
 Composing it with the morphism $\ff_*(\ff^*\cD_\fX^\bu\ocn_\fY\fF^\bu)
\rarrow\ff_*(\cD_\fY^\bu\ocn_\fY\fF^\bu)$ induced by
the quasi-isomorphism $\ff^*\cD_\fX^\bu\rarrow\cD_\fY^\bu$, we obtain
a natural morphism
\begin{equation} \label{ctrtnsr-and-smooth-direct-images-nat-transform}
 \cD_\fX^\bu\ocn_\fX\ff_!\fF^\bu\lrarrow
 \ff_*(\cD_\fY^\bu\ocn_\fY\fF^\bu)
\end{equation}
of complexes in $\fX\Tors^\inj$.
 It is important for our argument below that the natural
transformations~(\ref{fHom-and-smooth-direct-images-natural-transform}%
\+-\ref{ctrtnsr-and-smooth-direct-images-nat-transform})
are compatible with the compositions of weakly smooth morphisms~$\ff$ 
of finite type between Noetherian formal schemes of finite Krull
dimension.
 The two natural transformations are also compatible with
the adjunction~\eqref{fHom-contratensor-adjunction-eqn}, in
the sense explained in~\cite[Remark~7.6.4]{Pcosh}.

 As usual, it is helpful to observe that the totalization of any
finite acyclic complex in an exact category is an (absolutely)
acyclic complex.
 In particular, this applies to the split exact categories
$\fX\Tors^\inj$ and $\fX\Ctrh^\lct_\prj$.
 Hence, in order to prove that
the morphism~\eqref{fHom-and-smooth-direct-images-natural-transform} is
a homotopy equivalence for any complex $\cJ^\bu$ in $\fY\Tors$,
it suffices to check
that~\eqref{fHom-and-smooth-direct-images-natural-transform}
is a homotopy equivalence for any one-term complex
$\cJ^\bu=\cJ\in\fY\Tors^\inj$.
 Furthermore, a finite complex in $\fX\Ctrh^\lct_\prj$ is contractible
if and only if it is acyclic in $\fX\Ctrh^\lct$ (or in $\fX\Ctrh$).
 So it is suffices to prove
that~\eqref{fHom-and-smooth-direct-images-natural-transform} is
a quasi-isomorphism in $\fX\Ctrh$ for $\cJ^\bu=\cJ$.

 Dual-analogously, in order to prove that
the morphism~\eqref{ctrtnsr-and-smooth-direct-images-nat-transform}
is a homotopy equivalence for any complex $\fF^\bu$ in
$\fX\Ctrh^\lct_\prj$, it suffices to check
that~\eqref{ctrtnsr-and-smooth-direct-images-nat-transform}
is a homotopy equivalence for any one-term complex
$\fF^\bu=\fF\in\fY\Ctrh^\lct_\prj$.
 Furthermore, a finite complex in $\fX\Tors^\inj$ is contractible if
and only if it is acyclic in $\fX\Tors$.
 So it suffices to prove
that~\eqref{ctrtnsr-and-smooth-direct-images-nat-transform} is
a quasi-isomorphism in $\fX\Tors$ for $\fF^\bu=\fF$.

 We argue similarly to the proof of
Theorem~\ref{fin-krulldim-mgm-flat-direct-images-compatible-thm}.
 Cover the formal scheme $\fX$ with affine open formal
subschemes $\fU_\alpha$ and the formal scheme $\fY$ with affine open
formal subschemes $\fV_\beta$ such that for every index~$\beta$ there
exists an index~$\alpha$ with $\ff(\fV_\beta)\subset\fU_\alpha$.

 Let us show that
the morphism~\eqref{fHom-and-smooth-direct-images-natural-transform}
is a homotopy equivalence.
 One can use the fact that any injective quasi-coherent torsion
sheaf on $\fY$ is a finite direct sum of the direct images of injective
quasi-coherent torsion sheaves from~$\fV_\beta$.
 Alternatively, one can say that any complex $\cJ^\bu$ in
$\fY\Tors^\inj$ has a finite split filtration whose successive quotients
are direct images of complexes of injective quasi-coherent torsion
sheaves from~$\fV_\beta$, according to the proof of
Lemma~\ref{derived-fHom-takes-acyclics-to-acyclics}.
 In any case, the conclusion is that it suffices to show that
the morphism~\eqref{fHom-and-smooth-direct-images-natural-transform}
is a homotopy equivalence for complexes $\cJ^\bu$ of the form
$\cJ^\bu=\fk_\beta{}_*\cK^\bu$, where $\cK^\bu$ is a complex in
$\fV_\beta\Tors^\inj$ and $\fk_\beta\:\fV_\beta\rarrow\fX$ is
the open immersion morphism.

 Thus it suffices to show that the natural
transformation~\eqref{fHom-and-smooth-direct-images-natural-transform}
is a homotopy equivalence for the open immersions $\fV_\beta\rarrow\fY$
and $\fU_\alpha\rarrow\fX$, and for the weakly smooth morphisms of
finite type of affine Noetherian formal schemes
$\fV_\beta\rarrow\fU_\alpha$.
 In the case of an open immersion~$\ff$, the morphism
$\ff^*\cD_\fY^\bu\rarrow\cD_\fX^\bu$ is a homotopy equivalence;
hence, by construction,
the morphism~\eqref{fHom-and-smooth-direct-images-natural-transform}
is a homotopy equivalence, too.

 Let us treat the case of a weakly smooth morphism $\ff\:\fY\rarrow\fX$
of finite type between semi-separated Noetherian formal schemes $\fY$
and~$\fX$.
 In this case, for any injective quasi-coherent torsion sheaf
$\cJ$ on $\fY$, the functor $\fHom_\fY({-},\cJ)\:\fY\Tors\rarrow
\fY\Ctrh^\lct$ is exact, as explained in the proof of
Lemma~\ref{contratensor-exact-on-flats-and-alf}(c).
 So the morphism $\fHom_\fY(\cD_\fY^\bu,\cJ)\rarrow
\fHom_\fY(\ff^*\cD_\fX^\bu,\cJ)$ is a quasi-isomorphism in
$\fY\Ctrh^\lct$.
 Furthermore, by
Lemma~\ref{fHom-and-contratensor-with-alf-are-co-dilute}(a), the latter
morphism is a quasi-isomorphism of finite complexes of codilute locally
cotorsion contraherent cosheaves of contramodules on~$\fY$.
 By Corollary~\ref{codilute-direct-images}(b), it follows
that~\eqref{fHom-and-smooth-direct-images-natural-transform} is
a quasi-isomorphism in $\fX\Ctrh^\lct$ for $\cJ^\bu=\cJ$.
 According to the arguments above, this is sufficient.

 Alternatively, let us show that
the morphism~\eqref{ctrtnsr-and-smooth-direct-images-nat-transform}
is a homotopy equivalence.
 One can use the fact that any object of $\fY\Ctrh^\lct_\prj$ is
a finite direct sum of the direct images of objects from
$\fV_\beta\Ctrh^\lct_\prj$.
 Alternatively, one can say that any complex $\fF^\bu$ in
$\fY\Ctrh^\lct_\prj$ has a finite split filtration whose successive
quotients are direct images of complexes of projective locally
cotorsion contraherent cosheaves of contramodules
from~$\fV_\beta$, according to the proof of
Lemma~\ref{derived-contratensor-takes-acyclics-to-acyclics}.
 In any case, the conclusion is that it suffices to show that
the morphism~\eqref{ctrtnsr-and-smooth-direct-images-nat-transform}
is a homotopy equivalence for complexes $\fF^\bu$ of the form
$\fF^\bu=\fk_\beta{}_!\fG^\bu$, where $\fG^\bu$ is a complex in
$\fV_\beta\Ctrh^\lct_\prj$.

 Thus it is sufficient to show that the natural
transformation~\eqref{ctrtnsr-and-smooth-direct-images-nat-transform}
is a homotopy equivalence for the open immersions $\fV_\beta\rarrow\fY$
and $\fU_\alpha\rarrow\fX$, and for the weakly smooth morphisms of
finite type of affine Noetherian formal schemes
$\fV_\beta\rarrow\fU_\alpha$.
 In the case of an open immersion~$\ff$, the morphism
$\cD_\fX^\bu\ocn_\fX\ff_!\fF^\bu\rarrow
\ff_*(\ff^*\cD_\fX^\bu\ocn_\fY\fF^\bu)$ is an isomorphism
by Proposition~\ref{contratensor-projection-formula-prop}(c)
or~\ref{flat-contratensor-projection-formula-prop}, and the morphism
$\ff^*\cD_\fY^\bu\rarrow\cD_\fX^\bu$ is a homotopy equivalence.
 Hence, by construction,
the morphism~\eqref{ctrtnsr-and-smooth-direct-images-nat-transform}
is a homotopy equivalence, too. 

 Let us treat the case of a weakly smooth morphism $\ff\:\fY\rarrow\fX$
of finite type between semi-separated Noetherian formal schemes $\fY$
and~$\fX$.
 In this case, the morphism $\cD_\fX^\bu\ocn_\fX\ff_!\fF^\bu\rarrow
\ff_*(\ff^*\cD_\fX^\bu\ocn_\fY\fF^\bu)$ is an isomorphism by
Proposition~\ref{alf-contratensor-projection-formula-prop}.
 Furthermore, for any projective locally cotorsion (or just antilocally
flat) contraherent cosheaf of contramodules $\fF$ on $\fY$, 
the functor ${-}\ocn_\fY\fF\:\fY\Tors\rarrow\fY\Tors$ is exact by
Lemma~\ref{contratensor-exact-on-flats-and-alf}(c).
 So the morphism $\ff^*\cD_\fX^\bu\ocn_\fY\fF\rarrow
\cD_\fY^\bu\ocn_\fY\fF$ is a quasi-isomorphism in $\fY\Tors$.
 Finally, by
Lemma~\ref{fHom-and-contratensor-with-alf-are-co-dilute}(b), the latter
morphism is a quasi-isomorphism of finite compltexes of dilute
quasi-coherent torsion sheaves on~$\fY$.
 By Corollary~\ref{dilute-direct-images}, it follows
that~\eqref{ctrtnsr-and-smooth-direct-images-nat-transform} is
a quasi-isomorphism in $\fX\Tors$ for $\fF^\bu=\fF$.
 Once again, according to the arguments above, this is sufficient.
\end{proof}

\Section{Reminder of Contraherent Cosheaves on Schemes}
\label{reminder-secn}

 This section is a partial recollection of those definitions and
constructions from the book manuscript~\cite[Chapters~2\+-4]{Pcosh}
that are not fully covered by
Sections~\ref{contraherent-cosheaves-of-contramods-secn}\+-%
\ref{injective-projective-secn} of this paper, yet are needed
in the sequel.

\subsection{Locally contraherent cosheaves on schemes}
\label{lcth-on-schemes-subsecn}
 In the rest of this paper, we generally presume our schemes to
be locally Noetherian.
 This allows us to consider every scheme as a (locally Noetherian)
formal scheme, making all the preceding exposition in
Sections~\ref{contraherent-cosheaves-of-contramods-secn}\+-%
\ref{dualizing-secn} applicable to schemes as particular cases of
formal schemes.

 Still, sometimes we mention the Noetherianity assumptions on
schemes, for emphasis.
 On the other hand, we suppress all the mentions of quasi-separatedness
assumptions (as all locally Noetherian schemes and all morphisms of
locally Noetherian schemes are always quasi-separated).

 We denote by $\cO_X$ the structure sheaf of a scheme~$X$.
 The category of quasi-coherent sheaves on $X$ is denoted by
$X\Qcoh$; so we have $\fX\Tors=X\Qcoh$ when $\fX=X$ is a scheme.
 In particular, $\fX\Tors^\inj=X\Qcoh^\inj\subset X\Qcoh$ is the full
subcategory of injective quasi-coherent sheaves.

 Let $\bW$ be an open covering of a scheme~$X$.
 Specializing to schemes the definitions of locally contraherent
cosheaves of contramodules and locally cotorsion locally contraherent
cosheaves of contramodules from
Sections~\ref{locally-contraherent-cosheaves-of-contramods-subsecn}\+-%
\ref{locally-contraherent-lct-cosheaves-subsecn}, we obtain
the exact categories of (\emph{locally contraadjusted})
\emph{$\bW$\+locally contraherent cosheaves} $\fX\Lcth_\bW=X\Lcth_\bW$
and \emph{locally cotorsion\/ $\bW$\+locally contraherent cosheaves}
$\fX\Lcth_\bW^\lct=X\Lcth_\bW^\lct$ on a scheme~$X$.

 In line with our general notation system, we put
$X\Ctrh=X\Lcth_{\{X\}}$ and $X\Ctrh^\lct=X\Lcth_{\{X\}}^\lct$;
these are the exact categories of \emph{contraherent cosheaves} and
\emph{locally cotorsion contraherent cosheaves} on~$X$.
 Similarly, $X\Lcth=\bigcup_\bW X\Lcth_\bW$ and $X\Lcth^\lct=
\bigcup_\bW X\Lcth_\bW^\lct$ are the exact categories of
\emph{locally contraherent cosheaves} and \emph{locally cotorsion
locally contraherent cosheaves} on~$X$.
 All of these categories are full subcategories in the category of
cosheaves of $\cO_X$\+modules on~$X$ (which was defined in
Section~\ref{co-sheaves-of-modules-subsecn}),
$$
 X\Lcth^\lct\subset X\Lcth\subset(X,\cO_X)\Cosh.
$$
 We refer to~\cite[Sections~2.2 and~3.1\+-3.2]{Pcosh} for a more
detailed and specific discussion of contraherent and locally
contraherent cosheaves on schemes.

 As usual, we denote by $\Hom_X({-},{-})$ the groups of morphisms
in the abelian category $X\Qcoh$ and by $\Hom^X({-},{-})$ the groups
of morphisms in the exact category $X\Lcth$.
 The groups of morphisms in the abelian category $(X,\cO_X)\Sh$ of
sheaves of $\cO_X$\+modules on $X$ are denoted by
$\Hom_{\cO_X}({-},{-})$, while the groups of morphisms in
the additive category $(X,\cO_X)\Cosh$ are denoted by
$\Hom^{\cO_X}({-},{-})$.

\subsection{Locally injective locally contraherent cosheaves}
\label{loc-injective-on-schemes-subsecn}
 The important notions of a \emph{flat} and a \emph{very flat}
quasi-coherent sheaf do not seem to make sense in the context of
quasi-coherent torsion sheaves on formal schemes, so they need to be
defined separately for schemes.

 A quasi-coherent sheaf $\cF$ on a scheme $X$ is said to be
\emph{flat} if the $\cO_X(U)$\+module $\cF(U)$ is flat for every
affine open subscheme $U\subset X$.
 It suffices to check this condition for affine open subschemes $U$
belonging to any chosen affine open covering of~$X$
\,\cite[Proposition~IV.2.1.2]{EGAIV2}, \cite[Lemma Tag~01U4]{SP},
\cite[Lemma~1.8.3]{Pcosh}.

 A quasi-coherent sheaf $\cF$ on $X$ is said to be \emph{very flat}
if the $\cO_X(U)$\+module $\cF(U)$ is very flat for every affine
open subscheme $U\subset X$ (see the discussion of very flat modules
in Section~\ref{prelim-veryflat-and-contraadjusted-subsecn} above).
 Once again, it suffices to check this condition for affine open
subschemes $U$ belonging to any chosen affine open covering of~$X$
\,\cite[Lemma~1.2.6(a) and Section~1.10]{Pcosh},
\cite[Example~2.5]{Pal}.

 We denote the full subcategory of flat quasi-coherent sheaves by
$X\Qcoh_\fl\subset X\Qcoh$ and the full subcategory of very flat
quasi-coherent sheaves by $X\Qcoh_\vfl\subset X\Qcoh$.
 Both the full subcategories $X\Qcoh_\fl$ and $X\Qcoh_\vfl$ are
closed under extensions, kernels of epimorphisms, infinite direct
sums, and tensor products in $X\Qcoh$.
 In particular, both the full subcategories $X\Qcoh_\vfl\subset
X\Qcoh_\fl$ inherit exact category structures from the abelian
exact structure of $X\Qcoh$.
 We refer to~\cite[Sections~1.10 and~4.1]{Pcosh} for further details
(see also~\cite[Section~3.2]{PRoos}).

 Let $\bW$ be an open covering of a scheme~$X$.
 A $\bW$\+locally contraherent cosheaf $\fJ$ on $X$ is said to be
\emph{locally injective} if the $\cO_X(U)$\+module $\fJ[U]$ is
injective for every affine open subscheme $U\subset X$ subordinate
to~$\bW$.
 It suffices to check this condition for affine open subschemes $U$
belonging to any chosen affine open covering of $X$ subordinate
to~$\bW$ \,\cite[Lemma~1.3.6(b) and Section~3.1]{Pcosh},
\cite[Example~3.7]{Pal}.

 The full subcategory of locally injective $\bW$\+locally contraherent
cosheaves on $X$ is denoted by $X\Lcth^\lin_\bW\subset X\Lcth_\bW$.
 The category of locally injective contraherent cosheaves on $X$ is
denoted by $X\Ctrh^\lin=X\Lcth^\lin_{\{X\}}$.
 It follows from the previous paragraph that $X\Ctrh^\lin=
X\Ctrh\cap X\Lcth^\lin_\bW$.
 Moreover, for any open covering $\bW'$ of $X$ such that $\bW'$ is
subordinate to $\bW$, one has $X\Lcth^\lin_\bW=X\Lcth_\bW\cap
X\Lcth^\lin_{\bW'}$.
 The category of locally injective locally contraherent cosheaves of
contramodules on $X$ is denoted by $X\Lcth^\lin=
\bigcup_\bW X\Lcth^\lin_\bW$.

 By the definition, one has $X\Lcth^\lin_\bW\subset X\Lcth^\lct_\bW
\subset X\Lcth_\bW$.
 Hence $X\Ctrh^\lin\subset X\Ctrh^\lct\subset X\Ctrh$ and
$X\Lcth^\lin\subset X\Lcth^\lct\subset X\Lcth$.

 Clearly, the full subcategory $X\Lcth^\lin_\bW$ is closed under
extensions, cokernels of admissible monomorphisms, and infinite
products in $X\Lcth^\lct_\bW$ and in $X\Lcth_\bW$.
 So the full subcategories $X\Ctrh^\lin$, \,$X\Lcth^\lin_\bW$, and
$X\Lcth^\lin$ inherit exact category structures from the ambient
exact categories $X\Ctrh$, \,$X\Lcth_\bW$, and $X\Lcth$.
 A morphism in $X\Lcth^\lin_\bW$ is an admissible monomorphism in
$X\Lcth^\lin_\bW$ if and only if it is an admissible monomorphism in
$X\Lcth_\bW$.

 For an affine scheme $U$, the exact category of locally injective
contraherent cosheaves $U\Ctrh^\lin$ is split exact.
 The functor $\fJ\longmapsto\fJ[U]$ extablishes an equivalence between
the additive category $U\Ctrh^\lin$ and the additive category of
injective $\cO_U(U)$\+modules.
 We refer to~\cite[Sections~2.2 and~3.1\+-3.2]{Pcosh} for further
details (see also~\cite[Section~4.4]{Pns}).

\subsection{Direct images}  \label{direct-images-for-schemes-subsecn}
 Specializing to schemes the definitions of flat and very flat morphisms
of formal schemes in Section~\ref{morphisms-of-formal-schemes-subsecn},
one obtains the definitions of \emph{flat} and \emph{very flat}
morphisms of schemes $f\:Y\rarrow X$.
 The notion of a flat morphism of schemes is, of course, very well
known~\cite[Section~IV.2.1.1]{EGAIV2}, \cite[Definition Tag~01U3]{SP}.
 We refer to~\cite[Sections~1.8 and~1.10]{Pcosh} for a more detailed
and specific discussion of very flat morphisms of schemes.

 Let $f\:Y\rarrow X$ be a quasi-compact morphism of schemes.
 Specializing to schemes the construction of the direct image
functor~\eqref{direct-image-qcoh-tors}, one can obtain the (very
well known) functor of direct image of quasi-coherent sheaves
\begin{equation} \label{direct-image-qcoh}
 f_*\:Y\Qcoh\lrarrow X\Qcoh.
\end{equation}
 For an affine morphism~$f$, the functor~$f_*$ is exact.
 Whenever the morphism~$f$ is affine and flat, the functor~$f_*$ takes
flat quasi-coherent sheaves on $Y$ to flat quasi-coherent sheaves
on~$X$.
 Whenever the morphism~$f$ is affine and very flat, the functor~$f_*$
takes very flat quasi-coherent sheaves on $Y$ to very flat
quasi-coherent sheaves on~$X$ \,\cite[Lemma~1.2.3(b) and
Section~1.10]{Pcosh}.

 Let $f\:Y\rarrow X$ be an affine morphism of schemes.
 Specializing to schemes the constructions of the direct image
functors~(\ref{ctrh-contramod-affine-direct-image}\+-%
\ref{ctrh-lct-contramod-affine-direct-image}), one can obtain
the exact functors of direct image of locally contraadjusted and
locally cotorsion contraherent cosheaves
\begin{align}
 f_!\:Y\Ctrh &\lrarrow X\Ctrh, \label{ctrh-affine-direct-image} \\
 f_!\:Y\Ctrh^\lct &\lrarrow X\Ctrh^\lct.
 \label{ctrh-lct-affine-direct-image}
\end{align}
 Whenever the affine morphism~$f$ is flat,
the functors~(\ref{ctrh-affine-direct-image}\+-%
\ref{ctrh-lct-affine-direct-image}) take locally injective
contraherent cosheaves on $Y$ to locally injective contraherent
cosheaves on $X$, so there is an exact functor of direct image
of locally injective contraherent cosheaves
\begin{equation} \label{ctrh-lin-affine-direct-image}
 f_!\:Y\Ctrh^\lin\lrarrow X\Ctrh^\lin.
\end{equation}
 We refer to~\cite[Section~2.3]{Pcosh} for a more detailed and
specific discussion of the direct images of contraherent cosheaves
on schemes.

 More generally, let $\bW$ be an open covering of $X$ and $\bT$ be
an open covering of $Y$ such that the morphism $f\:Y\rarrow X$ is
$(\bW,\bT)$\+affine (in the sense of the definition
in Section~\ref{direct-images-of-ctrh-contramod-subsecn} above
or in~\cite[Section~3.3]{Pcosh}).
 Specializing to schemes the constructions of the direct image
functors~(\ref{loc-ctrh-contramod-W-T-affine-direct-image}\+-%
\ref{loc-ctrh-lct-contramod-W-T-affine-direct-image}), one can obtain
the exact functors of direct image of locally contraadjusted and
locally cotorsion locally contraherent cosheaves
\begin{align}
 f_!\:Y\Lcth_\bT &\lrarrow X\Lcth_\bW,
 \label{loc-ctrh-W-T-affine-direct-image} \\
 f_!\:Y\Lcth^\lct_\bT &\lrarrow X\Lcth^\lct_\bW.
 \label{loc-ctrh-lct-W-T-affine-direct-image}
\end{align}
 Whenever the $(\bW,\bT)$\+affine morphism~$f$ is flat,
the functors~(\ref{loc-ctrh-W-T-affine-direct-image}\+-%
\ref{loc-ctrh-lct-W-T-affine-direct-image}) take locally injective
$\bT$\+locally contraherent cosheaves on $Y$ to locally injective
$\bW$\+locally contraherent cosheaves on $X$, so there is an exact
functor of direct image of locally injective locally contraherent
cosheaves
\begin{equation} \label{loc-ctrh-lin-W-T-affine-direct-image}
 f_!\:Y\Lcth^\lin_\bT\lrarrow X\Lcth^\lin_\bW.
\end{equation}
 We refer to~\cite[Section~3.3]{Pcosh} for a more detailed and
specific discussion of the direct images of locally contraherent
cosheaves on schemes.

\subsection{Inverse images}  \label{inverse-images-for-schemes-subsecn}
 Let $f\:Y\rarrow X$ be a morphism of schemes.
 Specializing to schemes the construction of the inverse image
functor~\eqref{inverse-image-qcoh-tors}, one can obtain the (very
well known) functor of inverse image of quasi-coherent sheaves
\begin{equation} \label{inverse-image-qcoh}
 f^*\:X\Qcoh\rarrow Y\Qcoh.
\end{equation}
 For a flat morphism~$f$, the functor~$f^*$ is exact.
 For any morphism~$f$, the functor~$f^*$ takes flat quasi-coherent
sheaves on $X$ to flat quasi-coherent sheaves on~$Y$.
 For any morphism~$f$, the functor~$f^*$ also takes very flat
quasi-coherent sheaves on $X$ to very flat quasi-coherent sheaves
on~$Y$ \,\cite[Lemma~1.2.2(b) and Section~1.10]{Pcosh}.

 Let $\bW$ be an open covering of $X$ and $\bT$ be an open covering
of $Y$ such that the morphism $f\:Y\rarrow X$ is $(\bW,\bT)$\+coaffine
(in the sense of the definition in
Section~\ref{inverse-images-of-qcoh-tors-subsecn} above
or in~\cite[Sections~2.3 and~3.3]{Pcosh}).
 Let us first assume that the morphism~$f$ is very flat.
 Specializing to schemes the construction of the inverse image
functor~\eqref{inverse-image-lcth-W-T-contramod}, one can obtain
the exact functor of inverse image of locally contraadjusted locally
contraherent cosheaves
\begin{equation} \label{inverse-image-lcth-W-T}
 f^!\:X\Lcth_\bW\lrarrow Y\Lcth_\bT.
\end{equation}
 Passing to the inductive limits of categories with respect to
the refinements of the open coverings $\bW$ and $\bT$, one obtains
the inverse image functor
(cf.~\eqref{inverse-image-lcth-contramod-nocovering})
\begin{equation} \label{inverse-image-lcth-nocovering}
 f^!\:X\Lcth\lrarrow Y\Lcth.
\end{equation}

 More generally, let $f\:Y\rarrow X$ be a flat $(\bW,\bT)$\+coaffine
morphism of schemes.
 Specializing to schemes the construction of the inverse image
functor~\eqref{inverse-image-lcth-lct-W-T-contramod}, one can obtain
the exact functor of inverse image of locally cotorsion locally
contraherent cosheaves
\begin{equation} \label{inverse-image-lcth-lct-W-T}
 f^!\:X\Lcth^\lct_\bW\lrarrow Y\Lcth^\lct_\bT.
\end{equation}
 Passing to the inductive limits of categories with respect to
the refinements of the open coverings $\bW$ and $\bT$, one obtains
the inverse image functor
(cf.~\eqref{inverse-image-lcth-lct-contramod-nocovering})
\begin{equation} \label{inverse-image-lcth-lct-nocovering}
 f^!\:X\Lcth^\lct\lrarrow Y\Lcth^\lct.
\end{equation}

 Even more generally, let $f\:Y\rarrow X$ be an arbitrary
$(\bW,\bT)$\+coaffine morphism of schemes.
 Then the construction very similar to the one in
Section~\ref{inverse-images-of-ctrh-contramod-subsecn} provides
an exact functor of inverse image of locally injective locally
contraherent cosheaves
\begin{equation} \label{inverse-image-lcth-lin-W-T}
 f^!\:X\Lcth^\lin_\bW\lrarrow Y\Lcth^\lin_\bT.
\end{equation}
 Given a locally injective $\bW$\+locally contraherent cosheaf $\fJ$
on $X$, the locally injective $\bT$\+locally contraherent cosheaf
$f^!\fJ$ on $Y$ is defined by the rule
$$
 (f^!\fJ)[V]=\Hom_{\cO_X(U)}(\cO_Y(V),\fJ[U])
$$
for any pair of affine open subschemes $V\subset Y$ subordinate to
$\bT$ and $U\subset X$ subordinate to $\bW$ such that $f(V)\subset U$.

 One needs to check that, for a fixed affine open subscheme $V$,
the $\cO_Y(V)$\+module $\Hom_{\cO_X(U)}(\cO_Y(V),\fJ[U])$ does not
depend on the choice of an affine open subscheme $U$ up to a natural
isomorphism.
 The case of a semi-separated scheme $X$ is easier than the general
case.
 See the discussion in
Sections~\ref{inverse-images-of-qcoh-tors-subsecn}\+-%
\ref{inverse-images-of-ctrh-contramod-subsecn} above.

 Passing to the inductive limits of categories with respect to
the refinements of the open coverings $\bW$ and $\bT$, one obtains
the inverse image functor
\begin{equation} \label{inverse-image-lcth-lin-nocovering}
 f^!\:X\Lcth^\lin\lrarrow Y\Lcth^\lin.
\end{equation}
 We refer to~\cite[Sections~2.3 and~3.3]{Pcosh} for a more detailed
and specific discussion of the inverse images of contraherent and
locally contraherent cosheaves on schemes.

 For any morphism of schemes $f:Y\rarrow X$, the adjunction
isomorphism~\cite[formula~(3.4) in Section~3.3]{Pcosh}
\begin{equation} \label{cosh-on-schemes-lin-direct-inverse-adjunction}
 \Hom^{\cO_Y}(\fQ,f^!\fJ)\simeq\Hom^{\cO_X}(f_!\fQ,\fJ)
\end{equation}
of groups of morphisms in the categories $(Y,\cO_Y)\Cosh$ and
$(X,\cO_X)\Cosh$ holds for every locally injective locally
contraherent cosheaf $\fJ$ on $X$ and every cosheaf of
$\cO_Y$\+modules $\fQ$ on~$Y$.
 Similarly, for any flat morphism of schemes~$f$, the adjunction
isomorphism~\cite[formula~(3.5) in Section~3.3]{Pcosh}
\begin{equation} \label{cosh-on-schemes-lct-lcta-direct-inverse-adjn}
 \Hom^{\cO_Y}(\fQ,f^!\fP)\simeq\Hom^{\cO_X}(f_!\fQ,\fP)
\end{equation}
of groups of morphisms in the categories $(Y,\cO_Y)\Cosh$ and
$(X,\cO_X)\Cosh$ holds for every locally cotorsion locally
contraherent cosheaf $\fP$ on $X$ and every cosheaf of
$\cO_Y$\+modules $\fQ$ on~$Y$.
 For a very flat morphism of schemes~$f$, the adjnction
isomorphism~\eqref{cosh-on-schemes-lct-lcta-direct-inverse-adjn} holds
for every locally contraherent cosheaf $\fP$ on $X$ and every
cosheaf of $\cO_Y$\+modules $\fQ$ on~$Y$
(cf.~\eqref{cosheaves-direct-inverse-adjunction}).

\subsection{Contraadjusted and cotorsion quasi-coherent sheaves}
\label{cta-and-cot-qcoh-sheaves-subsecn}
 This section is an extraction from~\cite[Section~4.1]{Pcosh}.

 Let $X$ be a scheme.
 We denote by $\Ext_X^*({-},{-})$ the $\Ext$ groups in the abelian
category $X\Qcoh$ of quasi-coherent sheaves on~$X$.
 A quasi-coherent sheaf $\cP$ on $X$ is said to be
\emph{contraadjusted}~\cite[Section~2.5]{Pcosh},
\cite[Section~3.5]{PRoos} if $\Ext_X^1(\cF,\cP)=0$ for all very flat
quasi-coherent sheaves $\cF$ on~$X$ (in the sense of the definition
in Section~\ref{loc-injective-on-schemes-subsecn}).

 A quasi-coherent sheaf $\cP$ on an affine scheme $U$ is contraadjusted
if and only if the corresponding $\cO_U(U)$\+module $\cP(U)$ is
contraadjusted.
 The full subcategory $X\Qcoh^\cta$ of contraadjusted quasi-coherent
sheaves is closed under extensions and direct summands in $X\Qcoh$,
so $X\Qcoh^\cta$ inherits an exact category structure from the abelian
exact structure of $X\Qcoh$.
 The class of contraadjusted quasi-coherent sheaves is preserved by
the direct images with respect to quasi-compact morphisms of
schemes~\cite[Section~2.5]{Pcosh}.

\begin{prop} \label{cta-qcoh-sheaves-prop}
 Let $X$ be quasi-compact semi-separated scheme with a finite affine 
open covering $X=\bigcup_{\alpha=1}^N U_\alpha$.
 Then a quasi-coherent sheaf on $X$ is contraadjusted if and only if
it is a direct summand of a finitely iterated extension of the direct
images of contraadjusted quasi-coherent sheaves from\/~$U_\alpha$.
\end{prop}

\begin{proof}
 This is essentially~\cite[Corollary~4.1.4(c)]{Pcosh}.
\end{proof}

 The result of Proposition~\ref{cta-qcoh-sheaves-prop} can be rephrased
by saying that the class of contraadjusted quasi-coherent sheaves on
quasi-compact semi-separated schemes is~\emph{antilocal} in the sense
similar to~\cite[Section~4]{Pal}.

\begin{prop} \label{vfl-cta-cotorsion-pair-in-qcoh}
 Let $X$ be a quasi-compact semi-separated scheme.
 Then \par
\textup{(a)} for any quasi-coherent sheaf\/ $\cM$ on $X$ there exists
a short exact sequence of quasi-coherent sheaves\/ $0\rarrow\cP
\rarrow\cF\rarrow\cM\rarrow0$ on $X$ with a very flat quasi-coherent
sheaf\/ $\cF$ and a contraadjusted quasi-coherent sheaf\/~$\cP$; \par
\textup{(b)} for any quasi-coherent sheaf\/ $\cM$ on $X$ there exists
a short exact sequence of quasi-coherent sheaves\/ $0\rarrow\cM
\rarrow\cP\rarrow\cF\rarrow0$ on $X$ with a contraadjusted
quasi-coherent sheaf\/ $\cP$ and a very flat quasi-coherent
sheaf\/~$\cF$.
\end{prop}

\begin{proof}
 This is~\cite[Corollary~4.1.4(a\+-b)]{Pcosh}.
\end{proof}

\begin{cor} \label{cta-qcoh-sheaves-characterized}
 Let $X$ be a quasi-compact semi-separated scheme and\/ $\cP$ be
a quasi-coherent sheaf on~$X$.
 Then the following conditions are equivalent:
\begin{enumerate}
\item the quasi-coherent sheaf\/ $\cP$ is contraadjusted;
\item the functor\/ $\Hom_X({-},\cP)$ takes short exact sequences
of very flat quasi-coherent sheaves to short exact sequences of
abelian groups;
\item one has\/ $\Ext_X^n(\cF,\cP)=0$ for all very flat quasi-coherent
sheaves\/ $\cF$ on $X$ and all $n\ge1$.
\end{enumerate}
 Consequently, the class of contraadjusted quasi-coherent sheaves is
closed under cokernels of monomorphisms in $X\Qcoh$.
\end{cor}

\begin{proof}
 This is~\cite[Corollary~4.1.2]{Pcosh}.
 The assertions follow from the existence of enough very flat
quasi-coherent sheaves (a weak form of
Proposition~\ref{vfl-cta-cotorsion-pair-in-qcoh}(a)) by virtue of
Lemma~\ref{Hom-from-to-short-exact-sequences-lemma}(a) below.
\end{proof}

 In the terminology of Section~\ref{appx-cotorsion-pairs-subsecn} below,
the results of Proposition~\ref{vfl-cta-cotorsion-pair-in-qcoh}
and Corollary~\ref{cta-qcoh-sheaves-characterized}
can be summarized by saying that, on any quasi-compact semi-separated
scheme $X$, the pair of classes (very flat quasi-coherent sheaves,
contraadjusted quasi-coherent sheaves) is a \emph{hereditary complete
cotorsion pair} in the abelian category $X\Qcoh$.
 The direct summand lemma (Lemma~\ref{direct-summand-lemma} below)
is also relevant here.

 A quasi-coherent sheaf $\cP$ on a scheme $X$ is said to be
\emph{cotorsion}~\cite{EE}, \cite[Section~2.5]{Pcosh} if
$\Ext_X^1(\cF,\cP)=0$ for all flat quasi-coherent sheaves $\cF$ on~$X$.
 Clearly, any cotorsion quasi-coherent sheaf in contraadjusted.

 A quasi-coherent sheaf $\cP$ on an affine scheme $U$ is cotorsion if
and only if the corresponding $\cO_U(U)$\+module $\cP(U)$ is cotorsion.
 The full subcategory $X\Qcoh^\cot$ of cotorsion quasi-coherent sheaves
is closed under extensions and direct summands in $X\Qcoh$, so
$X\Qcoh^\cot$ inherits an exact category structure from the abelian
exact structure of $X\Qcoh$.
 The class of cotorsion quasi-coherent sheaves is preserved by
the direct images with respect to quasi-compact morphisms of
schemes~\cite[Section~2.5]{Pcosh}.

\begin{prop} \label{cot-qcoh-sheaves-prop}
 Let $X$ be quasi-compact semi-separated scheme with a finite affine 
open covering $X=\bigcup_{\alpha=1}^N U_\alpha$.
 Then a quasi-coherent sheaf on $X$ is cotorsion if and only if
it is a direct summand of a finitely iterated extension of the direct
images of cotorsion quasi-coherent sheaves from\/~$U_\alpha$.
\end{prop}

\begin{proof}
 This is essentially~\cite[Corollary~4.1.11(c)]{Pcosh}.
\end{proof}

 Similarly to the case of contraadjusted quasi-coherent sheaves,
the result of Proposition~\ref{cot-qcoh-sheaves-prop} can be rephrased
by saying that the class of cotorsion quasi-coherent sheaves on
quasi-compact semi-separated schemes is~\emph{antilocal}.

\begin{prop} \label{flat-cot-cotorsion-pair-in-qcoh}
 Let $X$ be a quasi-compact semi-separated scheme.
 Then \par
\textup{(a)} for any quasi-coherent sheaf\/ $\cM$ on $X$ there exists
a short exact sequence of quasi-coherent sheaves\/ $0\rarrow\cP
\rarrow\cF\rarrow\cM\rarrow0$ on $X$ with a flat quasi-coherent sheaf\/
$\cF$ and a cotorsion quasi-coherent sheaf\/~$\cP$; \par
\textup{(b)} for any quasi-coherent sheaf\/ $\cM$ on $X$ there exists
a short exact sequence of quasi-coherent sheaves\/ $0\rarrow\cM
\rarrow\cP\rarrow\cF\rarrow0$ on $X$ with a cotorsion quasi-coherent
sheaf\/ $\cP$ and a flat quasi-coherent sheaf\/~$\cF$.
\end{prop}

\begin{proof}
 This is~\cite[Corollary~4.1.11(a\+-b)]{Pcosh}.
\end{proof}

\begin{cor} \label{cot-qcoh-sheaves-characterized}
 Let $X$ be a quasi-compact semi-separated scheme and\/ $\cP$ be
a quasi-coherent sheaf on~$X$.
 Then the following conditions are equivalent:
\begin{enumerate}
\item the quasi-coherent sheaf\/ $\cP$ is cotorsion;
\item the functor\/ $\Hom_X({-},\cP)$ takes short exact sequences
of flat quasi-coherent sheaves to short exact sequences of
abelian groups;
\item one has\/ $\Ext_X^n(\cF,\cP)=0$ for all flat quasi-coherent
sheaves on $X$ and all $n\ge1$.
\end{enumerate}
 Consequently, the class of cotorsion quasi-coherent sheaves is
closed under cokernels of monomorphisms in $X\Qcoh$.
\end{cor}

\begin{proof}
 This is~\cite[Corollary~4.1.9]{Pcosh}.
 The assertions follow from the existence of enough flat
quasi-coherent sheaves by virtue of
Lemma~\ref{Hom-from-to-short-exact-sequences-lemma}(a) below.
\end{proof}

 In the terminology of Section~\ref{appx-cotorsion-pairs-subsecn} below,
the assertions of Proposition~\ref{flat-cot-cotorsion-pair-in-qcoh}
and Corollary~\ref{cot-qcoh-sheaves-characterized}
can be summarized by saying that, on any quasi-compact semi-separated
scheme $X$, the pair of classes (flat quasi-coherent sheaves,
cotorsion quasi-coherent sheaves) is a \emph{hereditary complete
cotorsion pair} in the abelian category $X\Qcoh$.
 This result goes back to~\cite[Corollary~4.2]{EE} together
with~\cite[Section~2.4]{M-n} or~\cite[Section~3.2]{M-th}
(see also~\cite[Lemma~A.1]{EP}).

 Let $X$ be a quasi-compact semi-separated scheme.
 Here is the notation for the intersections of full subcategories
in $X\Qcoh$ introduced in~\cite[Section~4.1]{Pcosh}:
\begin{alignat*}{2}
 &X\Qcoh^\cta_\vfl &&=X\Qcoh_\vfl\cap X\Qcoh^\cta; \\
 &X\Qcoh^\cta_\fl &&=X\Qcoh_\fl\cap X\Qcoh^\cta; \\
 &X\Qcoh^\cot_\fl &&=X\Qcoh_\fl\cap X\Qcoh^\cot.
\end{alignat*}
 Notice that $X\Qcoh^\cta_\vfl$ is the full subcategory of projective
objects in the exact category $X\Qcoh^\cta$ and simultaneously
the full subcategory of injective objects in the exact category
$X\Qcoh_\vfl$, by
Lemma~\ref{cotorsion-pair-kernel-as-injectives-projectives} below.
 Similarly, $X\Qcoh^\cot_\fl$ is the full subcategory of projective
objects in the exact category $X\Qcoh^\cot$ and simultaneously
the full subcategory of injective objects in the exact category
$X\Qcoh_\fl$, by the same lemma.
 So the full subcategories $X\Qcoh^\cta_\vfl$ and $X\Qcoh^\cot_\fl$
inherit trivial/split exact structures from the abelian exact
structure of $X\Qcoh$.

 The inherited exact structure on the full subcategory
$X\Qcoh^\cta_\fl\subset X\Qcoh$ is nontrivial, on the other hand.
 The exact category $X\Qcoh^\cta_\fl$ has both enough injective
and enough projective objects.
 The full subcategory of injective objects in $X\Qcoh^\cta_\fl$ is
$X\Qcoh^\cot_\fl$, while the full subcategory of projective objects
in $X\Qcoh^\cta_\fl$ is $X\Qcoh^\cta_\vfl$
(by Lemma~\ref{nested-cotorsion-pairs-projectives-injectives}).

\begin{prop} \label{flat-cta-qcoh-sheaves-prop}
 Let $X$ be quasi-compact semi-separated scheme with a finite affine 
open covering $X=\bigcup_{\alpha=1}^N U_\alpha$.
 In this setting: \par
\textup{(a)} a quasi-coherent sheaf on $X$ is very flat and
contraadjusted if and only if it is a direct summand of a finite
direct sum of the direct images of very flat contraadjusted
quasi-coherent sheaves from~$U_\alpha$; \par
\textup{(b)} a quasi-coherent sheaf on $X$ is flat and contraadjusted
if and only if it is a direct summand of a finitely iterated extension
of the direct images of flat contraadjusted quasi-coherent sheaves
from~$U_\alpha$; \par
\textup{(c)} a quasi-coherent sheaf on $X$ is flat and cotorsion if
and only if it is a direct summand of a finite direct sum of
the direct images of flat cotorsion quasi-coherent sheaves
from~$U_\alpha$.
\end{prop}

\begin{proof}
 Part~(a) is~\cite[Lemma~4.1.5]{Pcosh}.
 Part~(b) is~\cite[Lemma~4.1.6]{Pcosh}.
 Part~(c) is~\cite[Lemma~4.1.12]{Pcosh}.
\end{proof}

\begin{cor} \label{qcomp-ssep-direct-image-of-cta-cot-cor}
 Let $f\:Y\rarrow X$ be a morphism of quasi-compact semi-separated
schemes.
 Then \par
\textup{(a)} the functor of direct image of quasi-coherent sheaves
$f_*\:Y\Qcoh\rarrow X\Qcoh$ takes contraadjusted quasi-coherent
sheaves on $Y$ to contraadjusted quasi-coherent sheaves on $X$, and
induces an exact functor $f_*\:Y\Qcoh^\cta\rarrow Y\Qcoh^\cta$; \par
\textup{(b)} the functor of direct image of quasi-coherent sheaves
$f_*\:Y\Qcoh\rarrow X\Qcoh$ takes cotorsion quasi-coherent
sheaves on $Y$ to cotorsion quasi-coherent sheaves on $X$, and
induces an exact functor $f_*\:Y\Qcoh^\cot\rarrow Y\Qcoh^\cot$.
\end{cor}

\begin{proof}
 This is~\cite[Corollary~4.1.14]{Pcosh}.
\end{proof}

\begin{cor} \label{qcomp-ssep-direct-image-of-flat-cta-cor}
 Let $f\:Y\rarrow X$ be a morphism of quasi-compact semi-separated
schemes.
 In this context: \par
\textup{(a)} if the morphism~$f$ is very flat (in the sense of
Section~\ref{direct-images-for-schemes-subsecn}), then the functor of
direct image of quasi-coherent sheaves~$f_*$ takes very flat
contraadjusted quasi-coherent sheaves on $Y$ to very flat contraadjusted
quasi-coherent sheaves on~$X$; \par
\textup{(b)} if the morphism~$f$ is flat, then the functor of
direct image of quasi-coherent sheaves~$f_*$ takes flat contraadjusted
quasi-coherent sheaves on $Y$ to flat contraadjusted quasi-coherent
sheaves on~$X$; \par
\textup{(c)} if the morphism~$f$ is flat, then the functor of
direct image of quasi-coherent sheaves~$f_*$ takes flat cotorsion
quasi-coherent sheaves on $Y$ to flat cotorsion quasi-coherent
sheaves on~$X$.
\end{cor}

\begin{proof}
 Using Proposition~\ref{flat-cta-qcoh-sheaves-prop}
and Corollary~\ref{qcomp-ssep-direct-image-of-cta-cot-cor},
one reduces the desired assertions to the case of an affine
morphism~$f$, which is covered by the discussion in the beginning of
Section~\ref{direct-images-for-schemes-subsecn}.
 The point is that if $Y=\bigcup_{\beta=1}^M V_\beta$ is a finite
affine open covering of $Y$, then the compositions $V_\beta\rarrow Y
\rarrow X$ are affine morphisms with the same flatness or very flatness
properties as the morphism~$f$ (since any morphism from an affine
scheme to a semi-separated one is affine).
\end{proof}

\subsection{Antilocal contraherent cosheaves} \label{antilocal-subsecn}
 Let $X$ be a quasi-compact semi-sep\-a\-rated scheme with an open
covering~$\bW$.
 Specializing the discussion of projective contraherent cosheaves of
contramodules and projective locally cotorsion contraherent cosheaves
of contramodules from Section~\ref{semi-separated-inj-proj-subsecn} to
the case of a (nonformal) scheme $\fX=X$, we obtain descriptions of
the full subcategory $X\Ctrh_\prj$ of projective objects in
$X\Lcth_\bW$ and the full subcategory $X\Ctrh^\lct_\prj$ of
projective objects in $X\Lcth^\lct_\bW$.

 Furthermore, specializing the definition and description of the full
subcategory of antilocally flat contraherent cosheaves of contramodules
from Section~\ref{antilocally-flat-subsecn}, we obtain the definition
of description of the full subcategory $X\Ctrh_\alf\subset X\Lcth_\bW$
of \emph{antilocally flat} contraherent cosheaves on~$X$.

 The rest of this section is an extraction
from~\cite[Sections~4.3 and~4.6]{Pcosh}.

 Let $X$ be a scheme with an open covering~$\bW$.
 A $\bW$\+locally contraherent cosheaf $\fP$ on $X$ is said to be
\emph{antilocal} if, for any short exact sequence $0\rarrow\fI\rarrow
\fJ\rarrow\fK\rarrow0$ of locally injective $\bW$\+locally contraherent
cosheaves on $X$ (in the sense of the definition in
Section~\ref{loc-injective-on-schemes-subsecn}), the sequence of
abelian groups $0\rarrow\Hom^X(\fP,\fI)\rarrow \Hom^X(\fP,\fJ)\rarrow
\Hom^X(\fP,\fK)\rarrow0$ is exact.
 Clearly, any antilocally flat $\bW$\+locally contraherent cosheaf on
$X$ is antilocal.

 All contraherent cosheaves on an affine scheme $U$ with the open
covering $\bW_U=\{U\}$ are antilocal (since all short exact sequences
in $U\Ctrh^\lin\simeq\cO_U(U)\Modl^\inj$ are split).
 It follows from the adjunction
isomorphism~\eqref{cosh-on-schemes-lin-direct-inverse-adjunction}
and the formula~\eqref{inverse-image-lcth-lin-W-T}
that, for any $(\bW,\bT)$\+affine $(\bW,\bT)$\+coaffine morphism of
schemes $f\:Y\rarrow X$, the functor of direct image of $\bT$\+locally
contraherent cosheaves $f_!$
\,\eqref{loc-ctrh-W-T-affine-direct-image} takes antilocal
$\bT$\+locally contraherent cosheaves on $Y$ to antilocal
$\bW$\+locally contraherent cosheaves on~$X$.

\begin{prop} \label{antilocal-prop}
 Let $X$ be a quasi-compact semi-separated scheme. \par
\textup{(a)} Let\/ $\bW$ be an open covering of $X$ and
$X=\bigcup_{\alpha=1}^N U_\alpha$ be a finite affine open covering
subordinate to\/~$\bW$.
 Then a\/ $\bW$\+locally contraherent cosheaf on $X$ is antilocal if
and only if it is a direct summand of a finitely iterated extension
of the direct images of contraherent cosheaves from~$U_\alpha$. \par
\textup{(b)} All antilocal\/ $\bW$\+locally contraherent cosheaves on
$X$ are (globally) contraherent.
 The full subcategory of antilocal\/ $\bW$\+locally contraherent
cosheaves in the exact category $X\Lcth_\bW$ does not depend on an open
covering\/~$\bW$.
\end{prop}

\begin{proof}
 Part~(a) is~\cite[Corollary~4.3.5(c)]{Pcosh}.
 Part~(b) is~\cite[Corollary~4.3.7]{Pcosh}.
\end{proof}

 Let $X$ be a quasi-compact semi-separated scheme.
 Then the full subcategory of antilocal contraherent cosheaves on $X$
is denoted by $X\Ctrh_\al\subset X\Ctrh\subset X\Lcth_\bW$.
 The notation $X\Ctrh^\lct_\al$ stands for the intersection
$X\Ctrh^\lct_\al=X\Ctrh_\al\cap X\Lcth^\lct_\bW$.

\begin{prop} \label{antilocal-lct-prop}
 Let $X$ be a quasi-compact semi-separated scheme with a finite affine
open covering $X=\bigcup_{\alpha=1}^N U_\alpha$.
 Then a locally cotorsion ($\bW$\+locally) contraherent cosheaf on $X$
is antilocal if and only if its a direct summand of a finitely iterated
extension of the direct images of locally cotorsion contraherent
cosheaves from~$U_\alpha$.
\end{prop}

\begin{proof}
 This is~\cite[Corollary~4.3.6(c)]{Pcosh}.
\end{proof}

 The results of Propositions~\ref{antilocal-prop}(a)
and~\ref{antilocal-lct-prop} can be rephrased by saying that the classes
of antilocal contraherent cosheaves and antilocal locally cotorsion
contraherent cosheaves on quasi-compact semi-separated schemes are
\emph{antilocal} in the sense similar to~\cite[Section~4]{Pal}
(just as the terminology suggests).

\begin{prop} \label{al-lin-cotorsion-pair-in-lcth-W}
 Let $X$ be a quasi-compact semi-separated scheme with an open
covering\/~$\bW$.
 Then \par
\textup{(a)} for any\/ $\bW$\+locally contraherent cosheaf\/ $\fM$
on $X$ there exists an (admissible) short exact sequence\/
$0\rarrow\fM\rarrow\fJ\rarrow\fP\rarrow0$ in $X\Lcth_\bW$ with
a locally injective\/ $\bW$\+locally contraherent cosheaf\/ $\fJ$
and an antilocal contraherent cosheaf\/~$\fP$; \par
\textup{(b)} for any\/ $\bW$\+locally contraherent cosheaf\/ $\fM$
on $X$ there exists a short exact sequence\/ $0\rarrow\fJ\rarrow\fP
\rarrow\fM\rarrow0$ in $X\Lcth_\bW$ with an antilocal contraherent
cosheaf\/ $\fP$ and a locally injective\/ $\bW$\+locally contraherent
cosheaf\/~$\fJ$.
\end{prop}

\begin{proof}
 This is~\cite[Corollary~4.3.5(a\+-b)]{Pcosh}.
\end{proof}

 Let $X$ be a quasi-compact semi-separated scheme with an open
covering~$\bW$.
 It is explained in~\cite[Section~4.3]{Pcosh} that the $\Ext$ groups
computed in all our exact categories of (locally) contraherent
cosheaves on $X$, including $X\Ctrh$, $X\Ctrh^\lct$, $X\Ctrh^\lin$,
as well as $X\Lcth_\bW$, $X\Lcth^\lct_\bW$, $X\Lcth^\lin_\bW$ and
$X\Lcth$, $X\Lcth^\lct$, $X\Lcth^\lin$ agree with each other.
 The reason is that all those full exact subcategories are either
resolving or coresolving in each other, so the results of
Proposition~\ref{infinite-co-resolutions-minus-plus-derived} below
are applicable.
 We denote these $\Ext$ groups by $\Ext^{X,*}({-},{-})$.

\begin{cor} \label{antilocal-characterized}
 Let $X$ be a quasi-compact semi-separated scheme with an open
covering\/~$\bW$, and let\/ $\fP$ be a\/ $\bW$\+locally contraherent
cosheaf on\/~$X$.
 Then the following conditions are equivalent:
\begin{enumerate}
\item the cosheaf\/ $\fP$ is (contraherent and) antilocal;
\item $\Ext^{X,1}(\fP,\fJ)=0$ for all locally injective\/
$\bW$\+locally contraherent cosheaves\/ $\fJ$ on~$X$; \par
\item $\Ext^{X,n}(\fP,\fJ)=0$ for all locally injective\/
$\bW$\+locally contraherent cosheaves\/ $\fJ$ on $X$ and all $n\ge1$.
\end{enumerate}
 Consequently, the class of antilocal contraherent cosheaves is
closed under kernels of admissible epimorphisms in $X\Lcth_\bW$.
\end{cor}

\begin{proof}
 This is~\cite[Corollary~4.3.3]{Pcosh}.
 The assertions follow from the existence of enough locally injective
$\bW$\+locally contraherent cosheaves in $X\Lcth_\bW$ (a weak form
of Proposition~\ref{al-lin-cotorsion-pair-in-lcth-W}(a)) by virtue of
Lemma~\ref{Hom-from-to-short-exact-sequences-lemma}(b) below.
\end{proof}

\begin{prop} \label{al-lin-cotorsion-pair-in-lcth-lct-W}
 Let $X$ be a quasi-compact semi-separated scheme with an open
covering\/~$\bW$.
 Then \par
\textup{(a)} for any locally cotorsion\/ $\bW$\+locally contraherent
cosheaf\/ $\fM$ on $X$ there exists an (admissible) short exact
sequence\/ $0\rarrow\fM\rarrow\fJ\rarrow\fP\rarrow0$ in
$X\Lcth^\lct_\bW$ with a locally injective\/ $\bW$\+locally
contraherent cosheaf\/ $\fJ$ and an antilocal locally cotorsion
contraherent cosheaf\/~$\fP$; \par
\textup{(b)} for any locally cotorsion\/ $\bW$\+locally contraherent
cosheaf\/ $\fM$ on $X$ there exists a short exact sequence\/ $0\rarrow
\fJ\rarrow\fP\rarrow\fM\rarrow0$ in $X\Lcth^\lct_\bW$ with an antilocal
locally cotorsion contraherent cosheaf\/ $\fP$ and a locally
injective\/ $\bW$\+locally contraherent cosheaf\/~$\fJ$.
\end{prop}

\begin{proof}
 This is~\cite[Corollary~4.3.6(a\+-b)]{Pcosh}.
\end{proof}

 In the terminology of Section~\ref{appx-cotorsion-pairs-subsecn} below,
the results of Proposition~\ref{al-lin-cotorsion-pair-in-lcth-W}
and Corollary~\ref{antilocal-characterized} can be summarized by saying
that, on any quasi-compact semi-separated scheme $X$ with an open
covering $\bW$, the pair of classes (antilocal contraherent cosheaves,
locally injective $\bW$\+locally contraherent cosheaves) is
a \emph{hereditary complete cotorsion pair} in the exact category 
$X\Lcth_\bW$.
 The direct summand lemma (Lemma~\ref{direct-summand-lemma} below)
is also relevant here.

 Similarly, the results of
Proposition~\ref{al-lin-cotorsion-pair-in-lcth-lct-W}
and Corollary~\ref{antilocal-characterized} can be summarized by saying
that the pair of classes (antilocal locally cotorsion contraherent
cosheaves, locally injective $\bW$\+locally contraherent cosheaves) is
hereditary complete cotorsion pair in the exact category
$X\Lcth^\lct_\bW$.

 Let $X$ be a quasi-compact semi-separated scheme with an open
covering~$\bW$.
 Then the notation $X\Ctrh^\lin_\al$ stands for the intersection
$X\Ctrh^\lin_\al=X\Ctrh_\al\cap X\Lcth^\lin_\bW$.

 Notice that $X\Ctrh^\lin_\al$ is the full subcategory of projective
objects in the exact category $X\Lcth^\lin_\bW$ and simultaneously
the full subcategory of injective objects in the exact category
$X\Ctrh_\al$, by
Lemma~\ref{cotorsion-pair-kernel-as-injectives-projectives} below.
 So the full subcategory $X\Ctrh^\lin_\al$ inherits a trivial/split
exact structure from the exact category $X\Lcth_\bW$.

\begin{prop} \label{antilocal-lin-prop}
 Let $X$ be a quasi-compact semi-separated scheme with a finite affine
open covering $X=\bigcup_{\alpha=1}^N U_\alpha$.
 Then a locally injective ($\bW$\+locally) contraherent cosheaf on $X$
is antilocal if and only if its a direct summand of a finite direct sum
of the direct images of locally injective contraherent cosheaves
from~$U_\alpha$.
\end{prop}

\begin{proof}
 This is~\cite[Corollary~4.3.8]{Pcosh}.
\end{proof}

\begin{cor} \label{antilocal-direct-products-cor}
 Let $X$ be a quasi-compact semi-separated scheme with an open
covering\/~$\bW$.
 Then the full subcategories of antilocal contraherent cosheaves
$X\Ctrh_\al$, \,$X\Ctrh^\lct_\al$, and $X\Ctrh^\lin_\al$ are closed
under infinite direct products in the categories $X\Lcth_\bW$,
\,$X\Lcth^\lct_\bW$, and $X\Lcth^\lin_\bW$.
\end{cor}

\begin{proof}
 This is~\cite[Corollary~4.3.10]{Pcosh}.
\end{proof}

\begin{cor} \label{inverse-image-of-antilocal-cor}
 Let $f\:Y\rarrow X$ be an affine morphism of quasi-compact
semi-separated schemes.
 In this setting: \par
\textup{(a)} the functor of inverse image of locally injective
locally contraherent cosheaves $f^!\:X\Lcth^\lin\rarrow Y\Lcth^\lin$
\,\eqref{inverse-image-lcth-lin-nocovering} takes the full subcategory
$X\Ctrh^\lin_\al$ into the full subcategory $Y\Ctrh^\lin_\al$; \par
\textup{(b)} if the morphism~$f$ is flat, then the functor of inverse
image of locally cotorsion locally contraherent cosheaves
$f^!\:X\Lcth^\lct\rarrow Y\Lcth^\lct$
\,\eqref{inverse-image-lcth-lct-nocovering} takes the full subcategory
$X\Ctrh^\lct_\al$ into the full subcategory $Y\Ctrh^\lct_\al$; \par
\textup{(c)} if the morphism~$f$ is very flat, then the functor of
inverse image of locally contraherent cosheaves $f^!\:X\Lcth\rarrow
Y\Lcth$ \,\eqref{inverse-image-lcth-nocovering} takes the full 
subcategory $X\Ctrh_\al$ into the full subcategory $Y\Ctrh_\al$.
\end{cor}

\begin{proof}
 This is~\cite[Corollary~4.3.11]{Pcosh}.
\end{proof}

\begin{cor} \label{direct-image-of-antilocal-cor}
 Let $f\:Y\rarrow X$ be a morphism of quasi-compact semi-separated
schemes.
 In this setting: \par
\textup{(a)} the functor of direct image of cosheaves of\/
$\cO$\+modules\/ $f_!\:(Y,\cO_Y)\Cosh\rarrow(X,\cO_X)\Cosh$ takes
antilocal contraherent cosheaves on $Y$ to antilocal contraherent
cosheaves on $X$, and induces an exact functor\/ $f_!\:Y\Ctrh_\al
\rarrow X\Ctrh_\al$; \par
\textup{(b)} the functor of direct image of cosheaves of\/
$\cO$\+modules\/ $f_!\:(Y,\cO_Y)\Cosh\rarrow(X,\cO_X)\Cosh$ takes
antilocal locally cotorsion contraherent cosheaves on $Y$ to antilocal
locally cotorsion contraherent cosheaves on $X$, and induces
an exact functor\/ $f_!\:Y\Ctrh^\lct_\al\rarrow X\Ctrh^\lct_\al$; \par
\textup{(c)} if the morphism~$f$ is flat, then the functor of direct
image of cosheaves of\/ $\cO$\+modules\/ $f_!\:(Y,\cO_Y)\Cosh\rarrow
(X,\cO_X)\Cosh$ takes antilocal locally injective contraherent
cosheaves on $Y$ to antilocal locally injective contraherent cosheaves
on~$X$.
\end{cor}

\begin{proof}
 This is~\cite[Corollary~4.6.3(a,b,d)]{Pcosh}.
\end{proof}

\subsection{Contrahereable cosheaves and the contraherator}
\label{contraherator-subsecn}
 This section is an extraction from~\cite[Section~3.5]{Pcosh}.
 It is a cosheaf version of the construction of the \emph{coherator}
functor from~\cite[Appendix~B]{TT} and~\cite[Proposition Tag~077P]{SP}.

 Let $X$ be a scheme with an open covering~$\bW$.
 A cosheaf of $\cO_X$\+modules $\fP$ on $X$ is said to be
\emph{$\bW$\+locally derived contrahereable} if the following
two conditions hold:
\begin{enumerate}
\setcounter{enumi}{3}
\renewcommand{\theenumi}{\roman{enumi}$^\circ$}
\item for any affine open subscheme $U\subset X$ subordinate to $\bW$
and any finite affine open covering $U=\bigcup_{\alpha=1}^N U_\alpha$
of the affine open subscheme $U$, the homological \v Cech sequence
(cf.~\eqref{Cech-complex-of-abelian-groups-for-copresheaf})
\begin{multline} \label{Cech-sequence-exactness-condition}
 0\lrarrow\fP[U_1\cap\dotsb\cap U_N]\lrarrow\dotsb \\ \lrarrow
 \bigoplus\nolimits_{1\le\alpha<\beta\le N}\fP[U_\alpha\cap U_\beta]
 \lrarrow\bigoplus\nolimits_{\alpha=1}^N\fP[U_\alpha]\lrarrow
 \fP[U]\lrarrow0
\end{multline}
is exact, and
\renewcommand{\theenumi}{\roman{enumi}}
\item for any affine open subscheme $U\subset X$ subordinate to $\bW$,
the $\cO_X(U)$\+module $\fP[U]$ is contraadjusted.
\end{enumerate}
 We call the condition~(iv$^\circ$) the \emph{exactness condition}
and condition~(v) the \emph{contraadjustedness condition}.

 Notice that the contraadjustedness condition~(v) from the definition
above is equivalent to the contraadjustedness condition~(v) from
Section~\ref{contraherent-cosheaves-of-contramods-subsecn} restricted
from formal schemes to schemes.
 On other hand, the exactness condition~(iv$^\circ$) is weaker than
the scheme version of the contraherence condition~(iv) from
Section~\ref{contraherent-cosheaves-of-contramods-subsecn}, provided
that the contraadjustedness condition is assumed.
 This follows from~\cite[Lemma~1.2.6(b)]{Pcosh} or
Lemma~\ref{contraadj-homol-Cech-sequence-lemma} above.

 In view of Remarks~\ref{on-co-sheaf-axioms-remark}(1\+-2), for any
copresheaf of $\cO_X$\+modules $\fP$ on the topology base $\bB$ of
affine open subschemes in $X$ subordinate to $\bW$, the exactness
condition~(iv$^\circ$) implies the cosheaf
axiom~\eqref{topology-base-cosheaf-axiom}.
 The point is that the exactness of the rightmost fragment
of~\eqref{Cech-sequence-exactness-condition} is precisely the cosheaf
axiom for the open covering $U=\bigcup_{\alpha=1}^N U_\alpha$.
 So (iv$^\circ$)~is a stronger version of the cosheaf axiom.

 A cosheaf of $\cO_X$\+modules is called \emph{derived contrahereable}
if it is $\bW_X$\+locally derived contrahereable for the trivial open
covering $\bW_X=\{X\}$.
 A cosheaf of $\cO_X$\+modules is called \emph{locally derived
contrahereable} if it is $\bW$\+locally derived contrahereable for
some open covering $\bW$ of~$X$.

 Any $\bW$\+locally contraherent cosheaf on $X$ is $\bW$\+locally
derived contrahereable.
 Conversely, by~\cite[Lemma~3.2.2]{Pcosh} or
Proposition~\ref{homological-criterion-of-contraherence-prop} above,
if a $\bW$\+locally contraherent cosheaf is derived contrahereable,
then it is contraherent.

 A $\bW$\+locally derived contrahereable cosheaf $\fP$ on $X$ is
called \emph{locally cotorsion} if the $\cO_X(U)$\+module $\fP[U]$ is
cotorsion for every affine open subscheme $U\subset X$ subordinate
to~$\bW$.
 Given an open covering $\bW'$ subordinate to $\bW$, it follows from
condition~(iv$^\circ$) that a $\bW$\+locally derived contrahereable
cosheaf on $X$ is locally cotorsion if and only if it is locally
cotorsion \emph{as a\/ $\bW'$\+locally derived contrahereable cosheaf}.

 Similarly, a $\bW$\+locally derived contrahereable cosheaf $\fJ$ on $X$
is called \emph{locally injective} if the $\cO_X(U)$\+module $\fP[U]$
is injective for every affine open subscheme $U\subset X$ subordinate
to~$\bW$.
 It follows from condition~(iv$^\circ$) that a $\bW$\+locally derived
contrahereable cosheaf on $X$ is locally injective if and only if it is
locally injective \emph{as a\/ $\bW'$\+locally derived contrahereable
cosheaf}. 

 A short sequence of $\bW$\+locally derived contrahereable cosheaves
$0\rarrow\fL\rarrow\fM\rarrow\fN\rarrow0$ on $X$ is said to be
(\emph{admissible}) \emph{exact} if the short sequence of
$\cO_X(U)$\+modules $0\rarrow\fL[U]\rarrow\fM[U]\rarrow\fN[U]\rarrow0$
is exact for every affine open subscheme $U\subset X$ subordinate
to~$\bW$.
 A short sequence of $\bW$\+locally derived contrahereable cosheaves
on $X$ is exact if and only if it is exact as a short sequence of
$\bW'$\+locally derived contrahereable cosheaves.
 The category of $\bW$\+locally derived contrahereable cosheaves on
$X$, endowed with the class of short exact sequences defined above,
becomes an exact category.

\medskip

 Let $\fQ$ be a derived contrahereable cosheaf on an affine scheme~$U$.
 In this simple case, the \emph{contraherator} $\fC\fQ$ of the cosheaf
$\fQ$ is defined as the contraherent cosheaf on $U$ corresponding
to the contraadjusted $\cO_U(U)$\+module $\fQ[U]$ under the equivalence
of categories from~\cite[Corollary~2.2.3]{Pcosh} or
Lemma~\ref{contrah-cosheaves-of-contramods-lemma}(b).
 We use the notation similar to that of
Section~\ref{projective-lct-cosheaves-subsecn}: the contraherent
cosheaf on $U$ corresponding to a contraadjusted $\cO_U(U)$\+module
$Q$ is denoted by $\widecheck Q\in U\Ctrh$.
 So, by the definition, we have $\fC\fQ=\widecheck{\fQ[U]}$.
 There is a natural morphism of derived contrahereable cosheaves
$\fQ\rarrow\fC\fQ$ on~$U$ \,\cite[Lemma~2.2.4 and Section~3.5]{Pcosh}.

 Let $X$ be a quasi-compact semi-separated scheme with an open
covering $\bW$, and let $X=\bigcup_{\alpha=1}^N U_\alpha$ be a finite
affine open covering subordinate to~$\bW$.
 Let $\fP$ be a $\bW$\+locally derived contrahereable cosheaf on~$X$.
 The \emph{contraherator complex} $\fC_\bu(\{U_\alpha\},\fP)$ of $\fP$
is following finite \v Cech complex of (globally) contraherent
cosheaves on~$X$:
\begin{multline} \label{contraherator-complex}
 0\lrarrow j_{1,\dotsc,N}{}_!\fC(\fP|_{U_1\cap\dotsb\cap U_N})
 \lrarrow\dotsb \\
 \lrarrow\bigoplus\nolimits_{1\le\alpha<\beta\le N}
 j_{\alpha,\beta}{}_!\fC(\fP|_{U_\alpha\cap U_\beta})\lrarrow
 \bigoplus\nolimits_{\alpha=1}^N j_\alpha{}_!\fC(\fP|_{U_\alpha}),
\end{multline}
where $j_{\alpha_1,\dotsc,\alpha_k}\:U_{\alpha_1}\cap\dotsb\cap
U_{\alpha_k}\rarrow X$ denotes the identity open immersion morphism
and the notation $\fC(\fQ)$ was introduced above.
 The term $\bigoplus_{\alpha=1}^N j_\alpha{}_!\fC(\fP|_{U_\alpha})$
is placed in the homological degree~$0$.
 The construction of the differentials in
the complex~\eqref{contraherator-complex} is explained
in~\cite[Section~3.5]{Pcosh}.

 For a locally cotorsion $\bW$\+locally derived contrahereable
cosheaf $\fP$ on $X$, the complex $\fC_\bu(\{U_\alpha\},\fP)$
\,\eqref{contraherator-complex} is a complex of locally cotorsion
contraherent cosheaves on~$X$.
 For a locally injective $\bW$\+locally derived contrahereable
cosheaf $\fJ$ on $X$, the complex $\fC_\bu(\{U_\alpha\},\fJ)$
is a complex of locally injective contraherent cosheaves on~$X$.

\begin{lem} \label{derived-contraherator-well-defined}
 Let $X$ be a quasi-compact semi-separated scheme with an open
covering\/~$\bW$.
 In this setting: \par
\textup{(a)} Let\/ $\fP$ be a\/ $\bW$\+locally derived contrahereable
cosheaf on~$X$.
 Then the object of the bounded derived category\/ $\sD^\bb(X\Ctrh)$
represented by the complex\/ $\fC_\bu(\{U_\alpha\},\fP)$ in
the exact category of contraherent cosheaves on $X$ does not depend
ot the choice of a finite affine open covering
$X=\bigcup_\alpha U_\alpha$ subordinate to\/~$\bW$. \par
\textup{(b)} Let\/ $\fP$ be a locally cotorsion\/ $\bW$\+locally derived
contrahereable cosheaf on~$X$.
 Then the object of the bounded derived category\/
$\sD^\bb(X\Ctrh^\lct)$ represented by the complex\/
$\fC_\bu(\{U_\alpha\},\fP)$ in the exact category of locally cotorsion
contraherent cosheaves on $X$ does not depend ot the choice of a finite
affine open covering $X=\bigcup_\alpha U_\alpha$ subordinate to\/~$\bW$.
\par
\textup{(b)} Let\/ $\fJ$ be a locally injective\/ $\bW$\+locally derived
contrahereable cosheaf on~$X$.
 Then the object of the bounded derived category\/
$\sD^\bb(X\Ctrh^\lin)$ represented by the complex\/
$\fC_\bu(\{U_\alpha\},\fJ)$ in the exact category of locally injective
contraherent cosheaves on $X$ does not depend ot the choice of a finite
affine open covering $X=\bigcup_\alpha U_\alpha$ subordinate to\/~$\bW$.
\end{lem}

\begin{proof}
 Part~(a) is~\cite[Lemma~3.5.1]{Pcosh}, and parts~(b\+-c) are
provable similarly.
 Alternatively, parts~(b\+-c) actually follow from part~(a), since
the functors $\sD^\bb(X\Ctrh^\lin)\rarrow\sD^\bb(X\Ctrh^\lct)\rarrow
\sD^\bb(X\Ctrh)$ are fully faithful by
Proposition~\ref{infinite-co-resolutions-minus-plus-derived}(b).
 Notice that the full exact subcategories $X\Lcth^\lin_\bW\subset
X\Lcth^\lct_\bW\subset X\Lcth_\bW$ (and in particular, $X\Ctrh^\lin
\subset X\Ctrh^\lct\subset X\Ctrh$) are coresolving in each other
according to the discussion in
Section~\ref{loc-injective-on-schemes-subsecn} together with
the result of~\cite[Lemma~4.3.2(a) or Corollary~4.3.5(a)]{Pcosh}
or~\cite[Proposition~6.1]{Pns}.
\end{proof}

 A $\bW$\+locally derived contrahereable cosheaf $\fP$ on $X$ is called
\emph{$\bW$\+locally contrahereable} if the derived category object
represented by the complex $\fC_\bu(\{U_\alpha\},\fP)$ is isomorphic
in $\sD^\bb(X\Lcth_\bW)$ to some $\bW$\+locally contraherent cosheaf
on~$X$ (viewed as a complex concentrated in degree~$0$).
 If this is the case, then the mentioned $\bW$\+locally contraherent
cosheaf on $X$ is called the (\emph{$\bW$\+local}) \emph{contraherator}
of $\fP$ and denoted by~$\fC\fP$.
 A derived contrahereable cosheaf $\fP$ on $X$ is called
\emph{contrahereable} if it is locally contrahereable with respect
to the trivial open covering $\bW_X=\{X\}$.

 Any derived contrahereable cosheaf $\fQ$ on an affine scheme $U$ is
contrahereable, and the cosheaf $\fC\fQ$ constructed above specifically
in the affine scheme case is the contraherator of $\fQ$ in the sense of
the definition in the previous paragraph~\cite[Section~3.5]{Pcosh}.
 So our notation and terminology is consistent.

 Any $\bW$\+locally contraherent cosheaf $\fP$ on a quasi-compact
semi-separated scheme $X$ is $\bW$\+locally contrahereable, and
one has $\fC\fP=\fP$.
 The contraherator complex~\eqref{contraherator-complex} coincides
with (the scheme special case of) the contraherent \v Cech
resolution~\eqref{Cech-sequence-of-cosheaves-of-contramods} of
a $\bW$\+locally contraherent cosheaf $\fP$ in this case.

 The full subcategory of $\bW$\+locally contrahereable cosheaves is
closed under extensions, kernels of admissible epimorphisms, and
infinite products in the exact category of $\bW$\+locally derived
contrahereable cosheaves; so it inherits an exact category structure.
 The contraherator $\fC$ is an exact functor from the exact category
of $\bW$\+locally contrahereable cosheaves to that of $\bW$\+locally
contraherent cosheaves on~$X$.

 All of the above applies to locally cotorsion and locally injective
$\bW$\+locally derived contrahereable cosheaves just as well.
 The full subcategories of locally cotorsion and locally injective
$\bW$\+locally derived contrahereable cosheaves are closed under
extensions, cokernels of admissible monomorphisms, and infinite
products in the exact category of $\bW$\+locally derived contrahereable
cosheaves.

 A locally cotorsion $\bW$\+locally derived contrahereable cosheaf is
called $\bW$\+locally contrahereable if it is $\bW$\+locally
contahereable as a locally contraadjusted $\bW$\+locally derived
contrahereable cosheaf.
 Similarly, a locally injective $\bW$\+locally derived contrahereable
cosheaf is called $\bW$\+locally contrahereable if it is $\bW$\+locally
contahereable as a locally contraadjusted $\bW$\+locally derived
contrahereable cosheaf.
 The contraherator functor $\fC$ takes locally cotorsion $\bW$\+locally
contrahereable cosheaves to locally cotorsion $\bW$\+locally
contraherent cosheaves, and locally injective $\bW$\+locally
contrahereable cosheaves to locally injective $\bW$\+locally
contraherent cosheaves.

\medskip

 For any $(\bW,\bT)$\+affine morphism of schemes $f\:Y\rarrow X$,
the functor of direct image of cosheaves of $\cO$\+modules
$f_!\:(Y,\cO_Y)\Cosh\rarrow(X,\cO_X)\Cosh$ takes $\bT$\+locally
derived contrahereable cosheaves on $Y$ to $\bW$\+locally derived
contrahereable cosheaves on~$X$.
 For a $(\bW,\bT)$\+affine morphism~$f$ of quasi-compact
semi-separated schemes, the direct image functor~$f_!$ also commutes
with the construction of the contraherator complex (if one produces
the finite affine open covering of $Y$ subordinate to $\bT$ as
the preimage of any given finite affine open covering of $X$ subordinate
to~$\bW$) \cite[Section~3.5]{Pcosh}.
 It follows that the direct image functor~$f_!$ takes $\bT$\+locally
contrahereable cosheaves on $Y$ to $\bW$\+locally contrahereable
cosheaves on $X$, and commutes with the contraherator functor~$\fC$.
 So, for any $\bT$\+locally contrahereable cosheaf $\fQ$ on~$Y$, one
has a natural isomorphism
\begin{equation} \label{contraherator-commutes-with-W-T-affine-dirim}
 f_!(\fC\fQ)\simeq\fC(f_!\fQ)
\end{equation}
of $\bW$\+locally contraherent cosheaves on~$X$.

 For any $\bW$\+locally contrahereable cosheaf $\fP$ and any
$\bW$\+locally contraherent cosheaf $\fQ$ on a quasi-compact
semi-separated scheme $X$, there is a natural adjunction
isomorphism~\cite[formula~(3.9) in Section~3.5]{Pcosh}
\begin{equation} \label{contraherator-adjunction}
 \Hom^X(\fC\fP,\fQ)\simeq\Hom^{\cO_X}(\fP,\fQ)
\end{equation}
of the abelian groups of morphisms in the additive categories
$X\Lcth_\bW$ and $(X,\cO_X)\Cosh$.
 In other words, the contraherator functor $\fC$ is left adjoint to
the identity inclusion of the category of $\bW$\+locally
contraherent cosheaves into the category of $\bW$\+locally
contrahereable ones.

\subsection{The $\Cohom$ functor}  \label{Cohom-subsecn}
 This section is an extraction from~\cite[Sections~2.4 and~3.6]{Pcosh}.
 The functor $\Cohom_X$ from a quasi-coherent sheaf to a locally
contraherent cosheaf on a scheme $X$ is the dual analogue of the functor
of tensor product~$\ot_{\cO_X}$ of quasi-coherent sheaves on~$X$.

 The construction of the functor $\Cohom_X$ is relatively
straightforward and uncomplicated, but there is a caveat.
 As many functors into exact (rather than abelian) categories,
the functor $\Cohom_X$ is only \emph{partially defined}.
 See the discussion in~\cite[Section~5.6]{Pphil} and~\cite[Section~0.10
of the Introduction]{Pcosh}.

 So $\Cohom_X$ is a partially defined functor of two arguments
$$
 \Cohom_X\:X\Qcoh^\sop\times X\Lcth_\bW\dasharrow X\Lcth_\bW.
$$
 An itemized list of four alternative options of sufficient conditions
on a quasi-coherent sheaf $\cM$ and a $\bW$\+locally contraherent
cosheaf $\fP$ under which a $\bW$\+locally contraherent cosheaf
$\Cohom_X(\cM,\fP)$ is well-defined can be found
in~\cite[Section~5.5]{Pphil}.

 Let us spell out some details.
 More generally, a $\bW$\+locally derived contrahereable cosheaf
$\Cohom_X(\cM,\fP)$ can be defined for a quasi-coherent sheaf $\cM$
and a $\bW$\+locally derived contrahereable cosheaf $\fP$ on $X$
under certain assumptions.
 In all the cases, the defining formula for the cosheaf
$\Cohom_X(\cM,\fP)$ is
\begin{equation} \label{cohom-defining-formula}
 \Cohom_X(\cM,\fP)[U]=\Hom_{\cO_X(U)}(\cM(U),\fP[U])
\end{equation}
for all affine open subschemes $U\subset X$ subordinate to~$\bW$.
 One only needs to impose conditions on $\cM$ and $\fP$ guaranteeing
that the contraherence axiom~(iv) from
Section~\ref{contraherent-cosheaves-of-contramods-subsecn}
or the exactness condition~(iv$^\circ$) from
Section~\ref{contraherator-subsecn} is satisfied for
$\Cohom_X(\cM,\fP)$, together with the contraadjustedness axiom~(v).

 In fact, the contraherence axiom~(iv) does not require any
adjustedness assumptions on $\cM$ and~$\fP$.
 This axiom is satisfied for $\Cohom_X(\cM,\cP)$ and a given pair of
affine open subschemes $V\subset U$ in $X$ whenever it is satisfied
for the cosheaf $\fP$ and the same pair of affine open subschemes
$U\subset V$.
 (The quasi-coherence of the sheaf $\cM$ is presumed here.)
 However, some adjustedness assumptions on $\cM$ and $\fP$ are needed
in order to satisfy conditions~(iv$^\circ$) and/or~(v).

 Let $\cF$ be a very flat quasi-coherent sheaf (in the sense of
the definition in Section~\ref{loc-injective-on-schemes-subsecn}) and
$\fP$ be a $\bW$\+locally derived contrahereable cosheaf on~$X$
(in the sense of the definition in Section~\ref{contraherator-subsecn}).
 Then the $\bW$\+locally derived contrahereable cosheaf
$\Cohom_X(\cF,\fP)$ on $X$, given by the rule $U\longmapsto
\Hom_{\cO_X(U)}(\cF(U),\fP[U])$ \,\eqref{cohom-defining-formula},
is well-defined.
 If $\fP$ is a $\bW$\+locally contraherent cosheaf, then
$\Cohom_X(\cF,\fP)$ is a $\bW$\+locally contraherent cosheaf on~$X$.

 Let $\cF$ be a flat quasi-coherent sheaf and $\fP$ be a locally
cotorsion $\bW$\+locally derived contrahereable cosheaf on~$X$.
 Then the locally cotorsion $\bW$\+locally derived contrahereable
cosheaf $\Cohom_X(\cF,\fP)$ on $X$, given by the rule $U\longmapsto
\Hom_{\cO_X(U)}(\cF(U),\fP[U])$ \,\eqref{cohom-defining-formula},
is well-defined.
 If $\fP$ is a locally cotorsion $\bW$\+locally contraherent cosheaf,
then $\Cohom_X(\cF,\fP)$ is a locally cotorsion $\bW$\+locally
contraherent cosheaf on~$X$.

 Let $\cM$ be a quasi-coherent sheaf and $\fJ$ be a locally injective
$\bW$\+locally derived contrahereable cosheaf on~$X$.
 Then the locally cotorsion $\bW$\+locally derived contrahereable
cosheaf $\Cohom_X(\cM,\fJ)$ on $X$, given by the rule $U\longmapsto
\Hom_{\cO_X(U)}(\cM(U),\fJ[U])$ \,\eqref{cohom-defining-formula},
is well-defined.
 If $\fJ$ is a locally injective $\bW$\+locally contraherent cosheaf,
then $\Cohom_X(\cM,\fJ)$ is a locally cotorsion $\bW$\+locally
contraherent cosheaf on~$X$.

 In particular, let $\cF$ be a flat quasi-coherent sheaf and $\fJ$ be
a locally injective $\bW$\+locally derived contrahereable cosheaf
on~$X$.
 Then the locally injective $\bW$\+locally derived contrahereable
cosheaf $\Cohom_X(\cF,\fJ)$ on $X$, given by the rule $U\longmapsto
\Hom_{\cO_X(U)}(\cF(U),\fJ[U])$ \,\eqref{cohom-defining-formula},
is well-defined.
 If $\fJ$ is a locally injective $\bW$\+locally contraherent cosheaf,
then $\Cohom_X(\cF,\fJ)$ is a locally injective $\bW$\+locally
contraherent cosheaf on~$X$.

\subsection{The contraherent tensor product functor}
 This section is an extraction from~\cite[Section~3.7]{Pcosh}.
 The construction of the \emph{contraherent tensor product} of
a quasi-coherent sheaf and a locally contraherent cosheaf on
a scheme $X$ is dual-analogous to the contruction of
the quasi-coherent internal $\cHom$ sheaf from one quasi-coherent
sheaf on $X$ into another.
 Let us start with a brief discussion of the latter construction.

 Let $\cM$ and $\cN$ be two sheaves of $\cO_X$\+modules on
a scheme~$X$.
 Then the sheaf of $\cO_X$\+modules $\cHom_{\cO_X}(\cM,\cN)$ is
defined by the rule $\cHom_{\cO_X}(\cM,\cN)(Y)=
\Hom^{\cO_Y}(\cM|_Y,\cN|_Y)$ for any open subscheme $Y\subset X$.
 Here, as usual, $\Hom^{\cO_Y}$ denotes the group of morphisms in
the abelian category $(Y,\cO_Y)\Sh$.
 
 In the case of two quasi-coherent sheaves $\cM$ and $\cN$ on $X$,
the same sheaf of $\cO_X$\+modules $\cHom_{\cO_X}(\cM,\cN)$ on $X$
can be defined by the rule $\cHom_{\cO_X}(\cM,\cN)(U)=
\Hom_{\cO_X(U)}(\cM(U),\cN(U))$ for any affine open subscheme
$U\subset X$.
 Indeed, the abelian category $U\Qcoh$ is equivalent to $\cO_X(U)\Modl$,
and $U\Qcoh$ is a full subcategory in $(U,\cO_U)\Sh$.

 For a coherent sheaf $\cM$ and a quasi-coherent sheaf $\cN$ on
a locally Noetherian scheme $X$, the sheaf of $\cO_X$\+modules
$\cHom_{\cO_X}(\cM,\cN)$ is quasi-coherent.
 For two quasi-coherent sheaves $\cM$ and $\cN$ on $X$, the sheaf
of $\cO_X$\+modules $\cHom_{\cO_X}(\cM,\cN)$ is usually \emph{not}
quasi-coherent.

 Given two quasi-coherent sheaves $\cM$ and $\cN$ on $X$,
the \emph{quasi-coherent internal Hom sheaf} $\cHom_{X\qc}(\cM,\cN)
\in X\Qcoh$ is defined as the quasi-coherent sheaf representing
the functor $\cL\longmapsto\Hom_X(\cL\ot_{\cO_X}\cM,\>\cN)$
on the category $X\Qcoh$.
 So, by the definition, one has a natural isomorphism of abelian
groups $\Hom_X(\cL,\cHom_{X\qc}(\cM,\cN))\simeq
\Hom_X(\cL\ot_{\cO_X}\cM,\>\cN)$.

 The quasi-coherent sheaf $\cHom_{X\qc}(\cM,\cN)$ on $X$ can be
constructed by applying the coherator functor~$\cQ$
\,\cite[Appendix~B]{TT}, \cite[Proposition Tag~077P]{SP} to
the sheaf of $\cO_X$\+modules $\cHom_{\cO_X}(\cM,\cN)$.
 So one has $\cHom_{X\qc}(\cM,\cN)=\cQ\cHom_{\cO_X}(\cM,\cN)$.

\medskip

 The notion of a \emph{coadjusted} module over a commutative ring $R$,
defined in~\cite[Section~1.7]{Pcosh}, makes most sense in the case
when $R$ is Noetherian (cf.~\cite[Remarks~1.7.12 and~8.12.2]{Pcosh}).
 An $R$\+module $K$ is said to be coadjusted if the tensor product
functor $K\ot_R{-}$ takes contraadjusted $R$\+modules to
contraadjusted $R$\+modules.

 All finitely generated $R$\+modules are coadjusted (since the class
of contraadjusted $R$\+modules is closed under finite direct sums
and quotients).
 Over a Noetherian commutative ring $R$, all injective modules are
coadjusted (since infinite direct sums of injective modules are
injective and all quotient modules of injective modules are
contraadjusted) \cite[Section~1.7]{Pcosh}.

 The full subcategory of coadjusted $R$\+modules $R\Modl^\coa$ is closed
under extensions and quotients in $R\Modl$.
 Hence the full subcategory $R\Modl^\coa$ inherits an exact category
structure from the abelian exact structure of $R\Modl$.

 A quasi-coherent sheaf $\cK$ on a scheme $X$ is said to be
\emph{coadjusted} if the $\cO_X(U)$\+module $\cK(U)$ is coadjusted
for every affine open subscheme $U\subset X$.
 In view of~\cite[Lemma~1.7.11]{Pcosh}, it suffices to check this
condition for affine open subschemes $U$ belonging to any chosen
affine open covering of the scheme~$X$.

 So any coherent sheaf on a locally Noetherian scheme is coadjusted.
 Any injective quasi-coherent sheaf on a locally Noetherian scheme
is coadjusted as well (by~\cite[Theorem~II.7.18]{HarRD},
\cite[Theorem~A.3]{EP}, or
Corollary~\ref{local-section-criterion-of-injectivity-cor} above).

\medskip

 Let $\cM$ be a quasi-coherent sheaf and $\fP$ be a cosheaf of
$\cO_X$\+modules on a scheme~$X$.
 The copresheaf of $\cO_X$\+modules $\cM\ot_X\fP$ on the topology
base $\bB$ of all affine open subschemes $U\subset X$ is defined by
the rule
$$
 (\cM\ot_X\fP)[U]=\cM(U)\ot_{\cO_X(U)}\fP[U].
$$
 One can easily construct the corestriction maps
$(\cM\ot_X\fP)[V]\rarrow(\cM\ot_X\fP)[U]$ for all pairs of affine open
subschemes $V\subset U$ in $X$ (the key observation is that
$\cM(V)\ot_{\cO_X(V)}\fP[V]\simeq\cM(U)\ot_{\cO_X(U)}\fP[V]$), and
it is also easy to see that $\cM\ot_X\fP$ is a cosheaf of
$\cO_X$\+modules on~$\bB$ \,\cite[Section~3.7]{Pcosh}.
 So, in view of the discussion in
Section~\ref{co-sheaves-of-modules-subsecn}, we obtain a cosheaf
of $\cO_X$\+modules $\cM\ot_X\fP$ on whole of~$X$.

 The natural adjunction isomorphism of abelian
groups~\cite[formula~(3.14) in Section~3.7]{Pcosh}
\begin{equation} \label{Cohom-cosheaf-tensor-adjunction}
 \Hom^{\cO_X}(\fP,\Cohom_X(\cM,\fJ))\simeq
 \Hom^{\cO_X}(\cM\ot_X\fP,\>\cJ)
\end{equation}
holds for any quasi-coherent sheaf $\cM$, any cosheaf of
$\cO_X$\+modules $\fP$, and any locally derived contrahereable
cosheaf $\cJ$ on $X$ such that $\cM$ and $\cJ$ satisfy one of
the sufficient conditions from Section~\ref{Cohom-subsecn}
guaranteeing that the locally derived contrahereable cosheaf
$\Cohom_X(\cM,\fJ)$ on $X$ is well-defined.

 Let $\bW$ be an open covering of~$X$.
 Specializing the definition of a $\bW$\+flat cosheaf of contramodule
$\fO_\fX$\+modules given at the end of
Section~\ref{antilocally-flat-subsecn} to the case of a (nonformal)
scheme $\fX=X$, we obtain the definition of a $\bW$\+flat cosheaf
of $\cO_X$\+modules on~$X$.
 See~\cite[Sections~3.4, 3.7, 4.4, and~6.2]{Pcosh} for a more detailed
descussion of flat (contraherent) cosheaves specific to the case of
schemes.

 In the following two cases the cosheaf $\cM\ot_X\fP$ satisfies
the exactness condition~(iv$^\circ$) from
Section~\ref{contraherator-subsecn}:
\begin{itemize}
\item if $\fF$ is a $\bW$\+flat cosheaf of $\cO_X$\+modules satisfying
condition~(iv$^\circ$) for all affine open subschemes $U\subset X$
subordinate to $\bW$, and $\cM$ is an arbitrary quasi-coherent sheaf
on $X$, then the cosheaf $\cM\ot_X\fF$ satisfies condition~(iv$^\circ$) 
for all such affine open subschemes $U\subset X$;
\item if $\fP$ is a cosheaf of $\cO_X$\+modules satisfying
condition~(iv$^\circ$) for all affine open subschemes $U\subset X$
subordinate to $\bW$, and $\cF$ is a flat quasi-coherent sheaf on $X$,
then the cosheaf $\cF\ot_X\fP$ satisfies condition~(iv$^\circ$) 
for all such affine open subschemes $U\subset X$.
\end{itemize}

 In the case of a coadjusted quasi-coherent sheaf $\cK$ and
a $\bW$\+flat $\bW$\+locally derived contrahereable cosheaf $\fF$
on $X$, the cosheaf of $\cO_X$\+modules $\cK\ot_X\fF$ satisfies both
the conditions~(iv$^\circ$) and~(v) from
Section~\ref{contraherator-subsecn}.
 So $\cK\ot_X\fF$ is a $\bW$\+locally derived contrahereable cosheaf
on $X$ is this case.

 Over a locally Noetherian scheme $X$, the tensor product $\cM\ot_X\fP$
is a locally injective $\bW$\+locally derived contrahereable cosheaf
on $X$ in the following two cases~\cite[Section~3.7]{Pcosh}:
\begin{itemize}
\item for any $\bW$\+flat $\bW$\+locally derived contrahereable cosheaf
$\fF$ and any injective quasi-coherent sheaf $\cJ$ on $X$,
the tensor product $\cJ\ot_X\fF$ is a locally injective
$\bW$\+locally derived contrahereable cosheaf on~$X$;
\item for any locally injective $\bW$\+locally derived contrahereable
cosheaf $\fJ$ and any flat quasi-coherent sheaf $\cF$ on $X$,
the tensor product $\cF\ot_X\fJ$ is a locally injective
$\bW$\+locally derived contrahereable cosheaf on~$X$.
\end{itemize}

 Let $\cK$ be a coadjusted quasi-coherent sheaf and $\fF$ be
a $\bW$\+flat $\bW$\+locally derived contrahereable cosheaf on~$X$.
 Assuming that the scheme $X$ is quasi-compact and quasi-separated,
the \emph{derived contraherent tensor product} $\cK\ot_{X\ct}^\boL\fF
\in\sD^\bb(X\Ctrh)$ is defined as the derived category object 
represented by the contraherator complex of contraherent cosheaves
$\fC_\bu(\{U_\alpha\},\>\cK\ot_X\fF)$ constructed in
Section~\ref{contraherator-subsecn}.
 By Lemma~\ref{derived-contraherator-well-defined}(a), the object
$\cK\ot_{X\ct}^\boL\fF$ of the derived category $\sD(X\Ctrh)$ does
not depend on the choice of a finite affine open covering
$X=\bigcup_{\alpha=1}^N U_\alpha$ subordinate to~$\bW$.

 Whenever the object $\cK\ot_{X\ct}^\boL\fF\in\sD^\bb(X\Ctrh)$ turns out
to be isomorphic to an object from $X\Lcth_\bW$ in the derived
category of $\bW$\+locally contraherent cosheaves $\sD^\bb(X\Lcth_\bW)$,
one says that the \emph{contraherent tensor product} of $\cK$ and $\fF$
is defined.
 The resulting object is denoted by $\cK\ot_{X\ct}\fF\in X\Lcth_\bW$.
 In other words, we put $\cK\ot_{X\ct}\fF=\fC(\cK\ot_X\fF)$ whenever
the $\bW$\+locally derived contrahereable cosheaf $\cK\ot_X\fF$ on $X$
turns out to be actually $\bW$\+locally contrahereable (in the notation
and terminology of Section~\ref{contraherator-subsecn}).

 Over a locally Noetherian scheme $X$, the tensor product
$\cM\ot_X\fF$ is actually a $\bW$\+locally contraherent cosheaf
for any \emph{coherent} sheaf $\cM$ and any $\bW$\+flat $\bW$\+locally
contraherent cosheaf $\fF$ on~$X$ \,\cite[Corollary~1.7.4(c)
and Section~3.7]{Pcosh}.
 Whenever the cosheaf $\fF$ is locally cotorsion, so is the cosheaf
$\cM\ot_X\fF$ \,\cite[Corollary~1.7.5(a) and Section~3.7]{Pcosh}.
 So, assuming additionally that the scheme $X$ is quasi-compact and
semi-separated, one has $\cM\ot_{X\ct}\fF=\cM\ot_X\fF$ in these cases.

\subsection{The contraherent $\fHom$ functor~II}
\label{contraherent-Hom-II-subsecn}
 This section in an extraction from~\cite[Section~2.5]{Pcosh}.
 The functor of contraherent $\fHom_X$ between two quasi-coherent
sheaves on a scheme $X$ has \emph{no dual analogue within the realm
of sheaves}.
 Rather, the functor $\fHom_X$ is dual-analogous to the functor
of contratensor product~$\ocn_X$, which will be discussed in
the next Section~\ref{contratensor-product-II-subsecn}.

 Similarly to the functor $\Cohom_X$ from
Section~\ref{Cohom-subsecn} above, the functor $\fHom_X$ is only
\emph{partially defined}.
 See the discussion of partially defined functors into exact categories
in~\cite[Section~5.6]{Pphil} and~\cite[Section~0.10 of
the Introduction]{Pcosh}.

 So $\fHom_X$ is a partially defined functor of two arguments
$$
 \fHom_X\:X\Qcoh^\sop\times X\Qcoh\dasharrow X\Ctrh.
$$
 An itemized list of four alternative options of sufficient conditions
on quasi-coherent sheaves $\cM$ and $\cP$ on a scheme $X$ under which
a contraherent cosheaf $\fHom_X(\cM,\cP)$ is well-defined can be
found in~\cite[Section~8.7]{Pphil}.
 (However, the required semi-separatedness or quasi-separatedness
assumptions on the scheme $X$ are not mentioned in that list,
as quasi-compactness and semi-separatedness of the scheme are
generally presumed throughout~\cite[Section~8]{Pphil}.)

 Out of the four cases listed in~\cite[Section~8.7]{Pphil}, only one
set of assumptions for the definition of the functor $\fHom_X$ seems
make sense for quasi-coherent torsion sheaves and contraherent cosheaves
of contramodules on Noetherian formal schemes~$\fX$.
 This case was discussed in Section~\ref{contraherent-Hom-subsecn}.
 Specializing the construction presented in
Section~\ref{contraherent-Hom-subsecn} to the case of a (nonformal)
scheme $\fX=X$, we obtain the definition of the locally cotorsion
contraherent cosheaf $\fHom_X(\cM,\cJ)$ for any quasi-coherent sheaf
$\cM$ and any injective quasi-coherent sheaf $\cJ$ on~$X$.

 Let us spell out some details.
 In all cases, the defining formula for the cosheaf $\fHom_X(\cM,\cP)$
is similar to the one in Section~\ref{contraherent-Hom-subsecn},
\begin{equation} \label{fHom-defining-formula}
 \fHom_X(\cM,\cP)[U]=\Hom_X(j_*j^*\cM,\cP)
\end{equation}
for any affine open subscheme $U\subset X$ with the identity open
immersion morphism $j\:U\rarrow X$.
 One only needs to impose conditions on $\cM$ and $\cP$ guaranteeing
that the contraherence and contraadjustedness axioms~(iv\+-v)
from Section~\ref{contraherent-cosheaves-of-contramods-subsecn}
are satisfied for $\fHom_X(\cM,\cP)$.

 In fact, the contraherence axiom~(iv) does not require any
adjustedness assumptions on $\cM$ and~$\cP$.
 For any quasi-coherent sheaves $\cM$ and $\cP$ on $X$ and any
pair of affine open subschemes $V\subset U$ in $X$, the contraherence
axiom is satisfied for $\fHom_X(\cM,\cP)$.
 However, some adjustedness assumptions on $\cM$ and $\fP$ are needed
in order to satisfy the contraadjustedness axiom~(v).

 Let $\cF$ be a very flat quasi-coherent sheaf (in the sense of
the definition in Section~\ref{loc-injective-on-schemes-subsecn}),
and let $\cP$ be a contraadjusted quasi-coherent sheaf on~$X$
(in the sense of the definition in
Section~\ref{cta-and-cot-qcoh-sheaves-subsecn}).
 Assume that the scheme $X$ is semi-separated.
 Then the contraherent cosheaf $\fHom_X(\cF,\cP)$ on $X$, given
by the rule $U\longmapsto\Hom_X(j_*j^*\cF,\cP)$
\,\eqref{fHom-defining-formula}, is well-defined.

 Let $\cF$ be a flat quasi-coherent sheaf and $\cP$ be a cotorsion
quasi-coherent sheaf on~$X$ (in the sense of the definition in
Section~\ref{cta-and-cot-qcoh-sheaves-subsecn}).
 Assume that the scheme $X$ is semi-separated.
 Then the locally cotorsion contraherent cosheaf $\fHom_X(\cF,\cP)$
on $X$, given by the rule $U\longmapsto\Hom_X(j_*j^*\cF,\cP)$
\,\eqref{fHom-defining-formula}, is well-defined.

 In particular, let $\cF$ be a flat quasi-coherent sheaf and $\cJ$
be an injective quasi-coherent sheaf on~$X$.
 Assume that the scheme $X$ is semi-separated.
 Then the locally injective contraherent cosheaf $\fHom_X(\cF,\cJ)$
on $X$, given by the rule $U\longmapsto\Hom_X(j_*j^*\cF,\cJ)$
\,\eqref{fHom-defining-formula}, is well-defined.

 The following proposition complements the result of
Proposition~\ref{cosections-of-fHom-prop}.

\begin{prop}
 Let $X$ be a semi-separated scheme and $Y$ be a quasi-compact
open subscheme in $X$ with an affine open immersion morphism
$f\:Y\rarrow X$.
 Assume that either \par
\textup{(a)} $\cF$ is a very flat quasi-coherent sheaf and\/ $\cP$
is a contraadjusted quasi-coherent sheaf on~$X$; or \par
\textup{(b)} $\cF$ is a flat quasi-coherent sheaf and\/ $\cP$ is
a cotorsion quasi-coherent sheaf on~$X$. \par
\noindent Then there is a natural isomorphism of\/ $\cO_X(Y)$\+modules
$$
 \fHom_X(\cF,\cP)[Y]\simeq\Hom_X(f_*f^*\cF,\cP).
$$
\end{prop}

\begin{proof}
 This is~\cite[Lemma~2.5.2(a\+-b)]{Pcosh}.
\end{proof}

\subsection{The contratensor product functor~II}
\label{contratensor-product-II-subsecn}
 As mentioned in Section~\ref{contraherent-Hom-II-subsecn},
the functor of contraherensor product~$\ocn_X$ of a quasi-coherent sheaf
and a cosheaf of $\cO_X$\+modules on a scheme $X$ is dual-analogous
to the contraherent $\fHom$ functor.

 Both $\fHom_X$ and~$\ocn_X$ are essentially sheaf-cosheaf phenomena,
i.~e., constructions connecting sheaves with cosheaves.
 They are also essentially global constructions, in the sense that
these functors \emph{do not} seem to have any compatibility properties
with the functors of restriction to open subschemes.
 Instead, the functors $\fHom_X$ and~$\ocn_X$ have good
\emph{antilocality} properties, in the sense of compatibility with
the \emph{direct image} functors (see
Sections~\ref{projection-formulas-subsecn}
and~\ref{projection-contratensor-subsecn} above
and Section~\ref{projection-formulas-III-subsecn} below).

 Let $\cM$ be a quasi-coherent sheaf and $\fP$ be a cosheaf of
$\cO_X$\+modules on a scheme~$X$.
 Specializing the construction of the quasi-coherent torsion sheaf
of contratensor product from Section~\ref{contratensor-product-subsecn}
to the case of a (nonformal) scheme $\fX=X$, we obtain the construction
of the quasi-coherent sheaf of \emph{contrahensor product}
$$
 \cM\ocn_X\fP=\varinjlim\nolimits_{U\subset X}
 \bigl((j_*j^*\cM)\ot_{\cO_X(U)}\fP[U]\bigr)
$$
on~$X$.
 Here the inductive limit in the category of quasi-coherent sheaves
$X\Qcoh$ is taken over the (nondirected) poset of all affine open
subschemes $U\subset X$, and $j\:U\rarrow X$ denotes the identity open
immersion morphism.

 So the contratensor product is a functor
$$
 \ocn_X\:X\Qcoh\times(X,\cO_X)\Cosh\lrarrow X\Qcoh.
$$
 We refer to~\cite[Section~2.6]{Pcosh} for a more detailed and specific
discussion of the contratensor product functor over a scheme.

 Specializing the discussion of diagrams of affine open formal
subschemes from Section~\ref{contratensor-product-subsecn}, we obtain
the definition of a diagram of affine open subchemes, and particularly
the notion of an \emph{exhaustive} diagram of affine open subschemes
in~$X$.
 The contratensor product $\cM\ocn_X\fP$ can be computed as
the inductive limit over any exhaustive diagram of affine open
subschemes in~$X$.

\subsection{Projection formulas~III}
\label{projection-formulas-III-subsecn}
 This section is an extraction from~\cite[Sections~3.8 and~8.4]{Pcosh}.
 It complements the discussion in
Sections~\ref{projection-formulas-subsecn}
and~\ref{projection-contratensor-subsecn} in the context of schemes
(as opposed to formal schemes).

 Let $X$ be a scheme with an open covering $\bW$ and $Y$ be a scheme
with an open covering~$\bT$.
 Let $f\:Y\rarrow X$ be a flat $(\bW,\bT)$\+affine morphism of schemes.
 Then, for any quasi-coherent sheaf $\cM$ on $X$ and any locally
injective $\bT$\+locally contraherent cosheaf $\fK$ on $Y$, there
is a natural isomorphism of locally cotorsion $\bW$\+locally
contraherent cosheaves on~$X$
\begin{equation} \label{Cohom-into-flat-direct-image-of-lin}
 \Cohom_X(\cM,f_!\fK)\simeq f_!\Cohom_Y(f^*\cM,\fK).
\end{equation}
 This is~\cite[formula~(3.22) in Section~3.8]{Pcosh}.

\Section{Locally Injective Ind-Contraherent Ind-Cosheaves}
\label{ind-cosheaves-secn}

\subsection{Formal schemes as ind-schemes}
\label{formal-schemes-as-ind-schemes-subsecn}
 Introductory discussions of \emph{ind-schemes} can be found
in the notes~\cite[Section~1]{Rich}, in the unfinished
manuscript~\cite[Section~7.11]{BD2}, and in
the book~\cite[Chapter~1]{Psemten}.

 Let $\AffSch$ denote the category of affine schemes, $\CSch$ denote
the category of concentrated (i.~e., quasi-compact quasi-separated)
schemes, $\NoethAffSch$ denote the category of Noetherian affine
schemes, and $\NoethSch$ denote the category of Noetherian schemes.
 We denote by $\Sch$ the ambient category of (arbitrary) schemes and
by $\Sets$ the category of sets.

 Ind-schemes are represented by filtered inductive systems of
schemes $(X_\gamma\to X_\delta)_{\gamma<\delta\in\Gamma}$ indexed by
directed posets~$\Gamma$.
 To any such inductive system, the contravariant functor
$\CSch^\sop\rarrow\Sets$ defined by the rule $Z\longmapsto
\varinjlim_{\gamma\in\Gamma}\Mor_\Sch(Z,X_\gamma)$ is assigned.
 Here $\Mor_\sC({-},{-})$ is the notation for the set of morphisms in
a (nonadditive) category $\sC$, and the inductive limit over
$\gamma\in\Gamma$ is taken in the category of sets.

 The category of ind-schemes is defined as the full subcategory in
the category of contravariant functors $\AffSch^\sop\rarrow\Sets$
consisting of the functors representable by filtered inductive systems
of schemes.
 For the reasons explained in~\cite[Section~1.2]{Psemten}, the same
category can be defined as the full subcategory in the category of
contravariant functors $\CSch^\sop\rarrow\Sets$ consisting of
the functors representable by filtered inductive systems of schemes.
 The ind-scheme represented by a filtered inductive system of schemes
$(X_\gamma)_{\gamma\in\Gamma}$ is denoted by
$\fX=\ilim_{\gamma\in\Gamma}X_\gamma$.

 An ind-scheme is said to be \emph{ind-Noetherian} if it can be
represented by a filtered inductive system
$(X_\gamma)_{\gamma\in\Gamma}$ of Noetherian schemes~$X_\gamma$.
 The category of ind-Noetherian ind-schemes can be equivalently
defined as the full subcategory in the category of contravariant
functors $\NoethAffSch^\sop\rarrow\Sets$ consisting of the functors
representable by filtered inductive systems of Noetherian schemes.
 It can be also equivalently defined as the full subcategory in
the category of contravariant functors $\NoethSch^\sop\rarrow\Sets$
consisting of the functors representable by filtered inductive systems
of Noetherian schemes.

 An ind-scheme $\fX$ is said to be \emph{strict} if it can be
represented by a filtered inductive system $(X_\gamma\to
X_\delta)_{\gamma<\delta\in\Gamma}$ such that the transition
morphisms $i_{\gamma\delta}\:X_\gamma\rarrow X_\delta$ are closed
immersions of schemes.
 An \emph{$\aleph_0$\+ind-scheme} is an ind-scheme $\fX$ that can
be represented by a filtered inductive system
$(X_\gamma)_{\gamma\in\Gamma}$ indexed by a \emph{countable}
directed poset $\Gamma$, or equivalently, by an inductive system
$(X_\gamma)_{\gamma=0}^\infty$ indexed by the directed poset of
nonnegative integers $\Gamma=\boZ_{\ge0}$.

 In this paper, we are only interested in strict ind-Noetherian
$\aleph_0$\+ind-schemes.
 Any strict ind-Noetherian $\aleph_0$\+ind-scheme $\fX$ can be
represented as $\fX=\ilim_{\gamma\ge0}X_\gamma$, where
$(X_\gamma)_{\gamma=0}^\infty$ is an inductive system of closed
immersions of Noetherian schemes indexed by the poset of nonnegative
integers.
 We refer to~\cite[Chapter~1 and Section~4.1]{Psemten} for a further
discussion.

 Any Noetherian formal scheme $\fX$ can be viewed as a strict
ind-Noetherian $\aleph_0$\+ind-scheme in the following way.
 The contravariant functor $\NoethSch^\sop\rarrow\Sets$ corresponding
to $\fX$ is simply the functor of morphisms $Z\longmapsto
\Mor_\NoethFormSch(Z,\fX)$ in the category of Noetherian formal
schemes $\NoethFormSch$.
 Notice that any Noetherian scheme $Z$ can be considered as
a Noetherian formal scheme.

 Let $\overline X\subset\fX$ be a closed subscheme of definition
in $\fX$ corresponding to a coherent sheaf of ideals of definition
$\fI\subset\fO_\fX$ (as per
Section~\ref{basics-of-formal-schemes-subsecn}).
 For every integer $n\ge1$, denote by $\overline X^{(n)}\subset\fX$
the closed subscheme of definition corresponding to the coherent sheaf
of ideals of definition $\fI^n\subset\fO_\fX$.
 Then the ind-scheme corresponding to $\fX$ can be represented by
the countable filtered inductive system of closed immersions of
Noetherian schemes $\overline X=\overline X^{(1)}\rarrow
\overline X^{(2)}\rarrow\overline X^{(3)}\rarrow\dotsb$.

 The category of Noetherian formal schemes is a full subcategory of
the category of strict ind-Noetherian $\aleph_0$\+ind-schemes.
 By abuse of notation, we will simply write
$\fX=\ilim_{n\ge1}\overline X^{(n)}$ for a Noetherian formal
scheme~$\fX$.

\subsection{Flat and very flat pro-quasi-coherent pro-sheaves}
\label{flat-vfl-pro-qcoh-pro-scheaves-subsecn}
 The notion of a pro-quasi-coherent pro-sheaf on an ind-scheme was
suggested in~\cite[Section~7.11.3]{BD2} and taken up
in~\cite[Chapter~3]{Psemten}.
 The full subcategory of \emph{flat} pro-quasi-coherent pro-sheaves is
better behaved, in that it has a natural exact category structure.
 In this section, we specialize the discussion
from~\cite[Sections~3.1, 3.4, and~3.5]{Psemten} to the case of
a Noetherian formal scheme~$\fX$.

 Let $\fX$ be a Noetherian formal scheme.
 A \emph{flat pro-quasi-coherent pro-sheaf} $\rF$ on $\fX$ is
the following set of data:
\begin{enumerate}
\renewcommand{\theenumi}{\roman{enumi}}
\item to every closed subscheme of definition $Y\subset\fX$,
a flat quasi-coherent sheaf $\rF^{(Y)}$ on $Y$ is assigned;
\item to every pair of closed subschemes of definition
$Y$, $Z\subset \fX$, \,$Z\subset Y$, with the related closed
immersion morphism $i_{ZY}\:Z\rarrow Y$, an isomorphism
$i_{ZY}^*\rF^{(Y)}\simeq\rF^{(Z)}$ of flat quasi-coherent sheaves
on $Z$ is assigned;
\end{enumerate}
for which the following compatibility equation holds:
\begin{enumerate}
\renewcommand{\theenumi}{\roman{enumi}}
\setcounter{enumi}{2}
\item for every triple of closed subschemes of definition
$Y$, $Z$, $T\subset\fX$, \,$T\subset Z\subset Y$, the triangular
diagram of isomorphisms $i_{TY}^*\rF^{(Y)}\simeq
i_{TZ}^*\rF^{(Z)}\simeq\rF^{(T)}$ is commutative in $T\Qcoh_\fl$.
\end{enumerate}

 Let $(X_\gamma)_{\gamma\in\Gamma}$ be an inductive system of closed
subschemes of definition in $\fX$, indexed by a directed poset $\Gamma$,
such that $\fX=\ilim_{\gamma\in\Gamma}X_\gamma$.
 Then, in order to produce a flat pro-quasi-coherent pro-sheaf $\rF$ on
$\fX$, it suffices to specify flat quasi-coherent sheaves
$\rF^{(X_\gamma)}$ on $X_\gamma$ for all $\gamma\in\Gamma$, together
with isomorphisms $i_{X_\gamma X_\delta}^*\rF^{(X_\delta)}\simeq
\rF^{(X_\gamma)}$ in $X_\gamma\Qcoh_\fl$ for every $\gamma<\delta\in
\Gamma$, such that the compatibility equation~(iii) holds for all
$T=X_\beta$, \,$Z=X_\gamma$, \,$Y=X_\delta$, \ $\beta<\gamma<\delta
\in\Gamma$.
 The flat quasi-coherent sheaves $\rF^{(Y)}$ for all other closed
subschemes of definition $Y\subset\fX$ and the related
isomorphisms~(ii) can then be uniquely recovered so that
condition~(iii) is satisfied for all closed subschemes of
definition in~$\fX$.

 We denote the category of flat pro-quasi-coherent pro-sheaves on
a Noetherian formal scheme $\fX$ by $\fX\Proqc_\fl$.
 Following~\cite[Section~3.4]{Psemten}, we say that a short sequence of
flat pro-quasi-coherent pro-sheaves $0\rarrow\rF\rarrow\rG\rarrow\rH
\rarrow0$ on $\fX$ is (\emph{admissible}) \emph{exact} if, for every
closed subscheme of definition $Z\subset\fX$, the short sequence of flat
quasi-coherent sheaves $0\rarrow\rF^{(Z)}\rarrow\rG^{(Z)}\rarrow
\rH^{(Z)}\rarrow0$ is exact in the exact category $Z\Qcoh_\fl$.
 Given a filtered inductive system $(X_\gamma)_{\gamma\in\Gamma}$ of
closed subschemes of definition in $\fX$ such that
$\fX=\ilim_{\gamma\in\Gamma}X_\gamma$, it suffices to check
the exactness condition above for the closed subschemes $Z=X_\gamma$
with $\gamma\in\Gamma$.

 By~\cite[Proposition~3.5]{Psemten}, the rule stated in the previous
paragraph defines an exact category structure on the additive
category $\fX\Proqc_\fl$.

 A flat pro-quasi-coherent pro-sheaf $\rF$ on $\fX$ is said to be
\emph{very flat} if, for every closed subscheme of definition
$Z\subset\fX$, the flat quasi-coherent sheaf $\rF^{(Z)}$ on $Z$ is
very flat (in the sense of the definition in
Section~\ref{loc-injective-on-schemes-subsecn}).
 Given a filtered inductive system $(X_\gamma)_{\gamma\in\Gamma}$ of
closed subschemes of definition in $\fX$ such that
$\fX=\ilim_{\gamma\in\Gamma}X_\gamma$, it suffices to check
the very flatness condition for the closed subschemes $Z=X_\gamma$
with $\gamma\in\Gamma$.

 We denote the full subcategory of very flat pro-quasi-coherent
pro-sheaves by $\fX\Proqc_\vfl\subset\fX\Proqc_\fl$.
 The full subcategory $\fX\Proqc_\vfl$ is closed under extensions
and kernels of admissible epimorphisms in $\fX\Proqc_\fl$, so
$\fX\Proqc_\vfl$ inherits an exact category structure from
$\fX\Proqc_\fl$.

\begin{ex} \label{flat-pro-sheaves-on-affine=flat-contramodules}
 Let $\fU=\Spf\fR$ be an affine Noetherian formal scheme.
 So $\fR=\fO_\fU(\fU)$ is a complete, separated adic Noetherian ring.
 For any open ideal of definition $\fI\subset\fR$, denote by
$V_\fI=\Spec\fR/\fI\subset\fU$ the related closed subscheme of
definition in~$\fU$.
 So the correspondence $\fI\longmapsto V_\fI$ is an ordered bijection
between the directed posets of open ideals of definition in $\fR$
and closed subschemes of definition in~$\fU$.
 Then the rule $\fF\longmapsto\rF$ with $\rF^{(V_\fI)}=
\widetilde{\fF/\fI\fF}$ defines an equivalence of exact categories
$$
 \fR\Contra_\fl\simeq\fU\Proqc_\fl
$$
between the exact category $\fR\Contra_\fl$ of flat contramodule
$\fR$\+modules $\fF$ and the exact category $\fU\Proqc_\fl$ of flat
pro-quasi-coherent pro-sheaves $\rF$ on~$\fU$.
 Here $\fR\Contra_\fl\allowbreak\subset\fR\Contra$ is the full
subcategory of flat contramodules in the abelian category $\fR\Contra$,
endowed with the exact category structure inherited from the abelian
exact structure of $\fR\Contra$.
 The notation $\widetilde{\fF/\fI\fF}$ stands for the flat
quasi-coherent sheaf on the affine scheme $V_\fI=\Spec\fR/\fI$
corresponding to the flat $\fR/\fI$\+module $\fF/\fI\fF$.
 This equivalence of exact categories is a special case of
the discussion in~\cite[Examples~3.8(1\+-2)]{Psemten}.
 Another reference is~\cite[Lemma~8.6(b)]{Pflcc}.
 The results of
Lemmas~\ref{flat-quotseparated-contramodules-separated}
and~\ref{flat-quotseparated-contramodules-well-behaved}(b) above
are relevant here.

 The restriction of the equivalence of exact categories
$\fR\Contra_\fl\simeq\fU\Proqc_\fl$ to the full subcategories of
very flat contramodule $\fR$\+modules (in the sense of the definition
in Section~\ref{prelim-veryflat-and-contraadjusted-subsecn}) and
very flat pro-quasi-coherent pro-sheaves is an equivalence of exact
categories
$$
 \fR\Contra_\vfl\simeq\fU\Proqc_\vfl.
$$
\end{ex}

 For any Noetherian formal scheme $\fX$, the exact category
of flat pro-quasi-coherent pro-sheaves $\fX\Proqc_\fl$ has exact
functors of infinite direct sum~\cite[Section~3.5]{Psemten}.
 For all closed subschemes of definition $Z\subset\fX$, the functors
$\rF\longmapsto\rF^{(Z)}\:\fX\Proqc_\fl\rarrow Z\Qcoh_\fl$ preserve
infinite direct sums.
 The full exact subcategory of very flat pro-quasi-coherent pro-sheaves
$\fX\Proqc_\vfl$ is closed under infinite direct sums in the exact
category $\fX\Proqc_\fl$.

 We denote the groups of morphisms in the category $\fX\Proqc_\fl$ by
$\Hom_{(\fX)}({-},{-})$.

\subsection{Locally injective locally ind-contraherent ind-cosheaves}
\label{lin-ind-lcth-in-cosheaves-secn}
 The notion of a locally injective (locally) ind-contraherent
ind-cosheaf on a formal scheme or ind-scheme is dual-analogous to
that of a flat or very flat pro-quasi-coherent pro-sheaf (which
was discussed in the previous
Section~\ref{flat-vfl-pro-qcoh-pro-scheaves-subsecn}).

 Let $\fX$ be a Noetherian formal scheme with an open covering~$\bW$.
 For any closed subscheme of definition $Z\subset\fX$, the underlying
topological space of $Z$ coincides with the underlying topological
space of~$\fX$.
 With this fact in mind, we will denote simply by $\bW$ the same open
covering viewed as an open covering of~$Z$.

 Notice that, for any bijective closed immersion of Noetherian schemes
$i\:Z\rarrow Y$ and any open covering $\bW$ of $Y$, viewed also
as an open covering of $Z$, the morphism~$i$ is $(\bW,\bW)$\+coaffine
with respect to the open covering $\bW$ of $Y$ and the open covering
$\bW$ of~$Z$.
 This follows from the fact that a Noetherian scheme is affine if and
only if its maximal reduced closed subscheme is affine
(see~\cite[Exercise~III.3.1]{HarAG}, \cite[Lemma Tag~01YQ]{SP},
and/or Section~\ref{basics-of-formal-schemes-subsecn} above).

 A \emph{locally injective\/ $\bW$\+locally ind-contraherent
ind-cosheaf} $\gJ$ on $\fX$ is the following set of data:
\begin{enumerate}
\renewcommand{\theenumi}{\roman{enumi}}
\item to every closed subscheme of definition $Y\subset\fX$,
a locally injective $\bW$\+locally contraherent cosheaf $\gJ_{(Y)}$
on~$Y$ (in the sense of the definition in
Section~\ref{loc-injective-on-schemes-subsecn}) is assigned;
\item to every pair of closed subschemes of definition
$Y$, $Z\subset \fX$, \,$Z\subset Y$, with the related closed
immersion morphism $i_{ZY}\:Z\rarrow Y$, an isomorphism
$\gJ_{(Z)}\simeq i_{ZY}^!\gJ_{(Y)}$ of locally injective $\bW$\+locally
contraherent cosheaves on~$Z$ (with the $f^!$~notation referring to
formula~\eqref{inverse-image-lcth-lin-W-T}) is assigned;
\end{enumerate}
for which the following compatibility equation holds:
\begin{enumerate}
\renewcommand{\theenumi}{\roman{enumi}}
\setcounter{enumi}{2}
\item for every triple of closed subschemes of definition
$Y$, $Z$, $T\subset\fX$, \,$T\subset Z\subset Y$, the triangular
diagram of isomorphisms $\gJ_{(T)}\simeq i_{TZ}^!\gJ_{(Z)}\simeq
i_{TY}^!\gJ_{(Y)}$ is commutative in $T\Lcth^\lin_\bW$.
\end{enumerate}

 Let $(X_\gamma)_{\gamma\in\Gamma}$ be an inductive system of closed
subschemes of definition in $\fX$, indexed by a directed poset $\Gamma$,
such that $\fX=\ilim_{\gamma\in\Gamma}X_\gamma$.
 Then, in order to produce a locally injective $\bW$\+locally
ind-contraherent ind-cosheaf $\gJ$ on $\fX$, it suffices to specify
locally injective $\bW$\+locally contraherent cosheaves
$\gJ_{(X_\gamma)}$ on $X_\gamma$ for all $\gamma\in\Gamma$,
together with isomorphisms $\gJ_{(X_\gamma)}\simeq
i_{X_\gamma X_\delta}^!\gJ_{(X_\delta)}$ in $X_\gamma\Lcth^\lin_\bW$
for every $\gamma<\delta\in\Gamma$, such that the compatibility
equation~(iii) holds for all $T=X_\beta$, \,$Z=X_\gamma$,
\,$Y=X_\delta$, \ $\beta<\gamma<\delta \in\Gamma$.
 The locally injective $\bW$\+locally contraherent cosheaves $\gJ_{(Y)}$
for all other closed subschemes of definition $Y\subset\fX$ and
the related isomorphisms~(ii) can then be uniquely recovered so that
condition~(iii) is satisfied for all closed subschemes of
definition in~$\fX$.

 We denote the category of locally injective $\bW$\+locally
ind-contraherent ind-cosheaves on a Noetherian formal scheme $\fX$
by $\fX\Indlc^\lin_\bW$.
 A locally injective $\bW_{\fX}$\+locally ind-contraherent ind-cosheaf 
on $\fX$ for the trivial open covering $\bW_{\fX}=\{\fX\}$ is said to
be a \emph{locally injective ind-contraherent ind-cosheaf} on~$\fX$.
 The category of locally injective ind-contraherent ind-cosheaves
on $\fX$ is denoted by $\fX\Indct^\lin=\fX\Indlc^\lin_{\{\fX\}}$.
 The direct limit (essentialy, the direct union) of the categories
$\fX\Indlc^\lin_\bW$ over refinements of the open coverings $\bW$ is
denoted by $\fX\Indlc^\lin=\bigcup_\bW\Indlc^\lin_\bW$.
 The objects of the category $\fX\Indlc^\lin$ are called
\emph{locally injective locally ind-contraherent ind-cosheaves}
on~$\fX$.

 A short sequence of locally injective $\bW$\+locally ind-contraherent
ind-cosheaves $0\rarrow\gJ\rarrow\gK\rarrow\gL\rarrow0$ on $\fX$ is said
to be (\emph{admissible}) \emph{exact} if, for every closed subscheme
of definition $Z\subset\fX$, the short sequence of locally injective
$\bW$\+locally contraherent cosheaves $0\rarrow\gJ_{(Z)}\rarrow\gK_{(Z)}
\rarrow\gL_{(Z)}\rarrow0$ is exact in the exact category 
$Z\Lcth^\lin_\bW$.
 Given a filtered inductive system $(X_\gamma)_{\gamma\in\Gamma}$ of
closed subschemes of definition in $\fX$ such that
$\fX=\ilim_{\gamma\in\Gamma}X_\gamma$, it suffices to check
the exactness condition above for the closed subschemes $Z=X_\gamma$
with $\gamma\in\Gamma$.

\begin{lem} \label{indlc-lin-W-is-an-exact-category-lem}
 For any Noetherian scheme\/ $\fX$ with an open covering $\bW$,
the category\/ $\fX\Indlc^\lin_\bW$ of locally injective\/
$\bW$\+locally ind-contraherent ind-cosheaves on\/ $\fX$, endowed
with the class of admissible short exact sequences defined above,
is an exact category.
\end{lem}

\begin{proof}
 Similar to the proof of~\cite[Proposition~3.5]{Psemten}.
\end{proof}

 In particular, in the case of the trivial open covering
$\bW_\fX=\{\fX\}$, we obtain an exact category structure on
the category of locally injective ind-contraherent ind-cosheaves
$\fX\Indct^\lin$ on a Noetherian formal scheme~$\fX$.

 Let $\bW'$ be an open covering of $\fX$ subordinate to~$\bW$.
 Then a short sequence in $\fX\Indlc^\lin_\bW$ is exact in
$\fX\Indlc^\lin_{\bW'}$ if and only if it is exact in
$\fX\Indlc^\lin_\bW$.
 This follows from the fact that a short exact sequence in
$Z\Lcth^\lin_\bW$ is exact in $Z\Lcth^\lin_{\bW'}$ if and only if
it is exact in $Z\Lcth^\lin_\bW$ (see
Sections~\ref{locally-contraherent-cosheaves-of-contramods-subsecn}
and~\ref{loc-injective-on-schemes-subsecn} above,
or~\cite[Lemma~1.4.3(a) and Section~2.2]{Pcosh}).

 Passing to the direct limit of exact categories over the open
coverings $\bW$, we obtain an exact category structure on
the exact category of locally injective locally ind-contraherent
ind-cosheaves $\fX\Indlc^\lin$.
 So a short sequence is exact in $\fX\Indlc^\lin$ if and only if
it is exact in $\fX\Indlc^\lin_\bW$ for some open covering $\bW$
of~$\fX$.

\begin{ex} \label{lin-ind-ctrh-on-affine=injective-torsion-modules}
 Let $\fU=\Spf\fR$ be an affine Noetherian formal scheme.
 We keep the notation of
Example~\ref{flat-pro-sheaves-on-affine=flat-contramodules}.

 In this context, the rule $\cJ\longmapsto\gJ$ with
$\gJ_{(V_\fI)}=\widecheck{({}_\fI\cJ)}$, where ${}_\fI\cJ=
\Hom_\fR(\fR/\fI,\cJ)$, defines an equivalence of categories
$$
 \fR\Tors^\inj\simeq\fU\Indct^\lin
$$
between the additive category $\fR\Tors^\inj$ of injective torsion
$\fR$\+modules $\cJ$ and the additive category $\fU\Indct^\lin$ of
locally injective ind-contraherent ind-cosheaves $\gJ$ on~$\fU$.
 Here the notation $\widecheck{({}_\fI\cJ)}$ stands for the locally
injective contraherent cosheaf on the affine scheme
$V_\fI=\Spec\fR/\fI$ corresponding to the injective $\fR/\fI$\+module
$\Hom_\fR(\fR/\fI,\cJ)$ under the equivalence of additive categories
$\fR/\fI\Modl^\inj\simeq V_\fI\Ctrh^\lin$ from the end of
Section~\ref{loc-injective-on-schemes-subsecn}.

 It follows that, for an affine Noetherian formal scheme $\fU$,
the exact structure on the category $\fU\Indct^\lin$ is trivial,
i.~e., every (admissible) exact sequence of locally injective
ind-contraherent ind-cosheaves on $\fU$ is split exact.
 This assertion only applies to locally injective ind-contraherent
ind-cosheaves and \emph{not} to locally injective \emph{locally}
ind-contraherent ind-cosheaves, of course.
\end{ex}

 For any Noetherian formal scheme $\fX$ with an open covering $\bW$,
the exact category of locally injective $\bW$\+locally
ind-contraherent ind-cosheaves $\fX\Indlc^\lin_\bW$ has exact functors
of infinite direct product.
 For all closed subschemes of definition $Z\subset\fX$, the functors
$\gJ\longmapsto\gJ_{(Z)}\:\fX\Indlc^\lin_\bW\rarrow Z\Lcth^\lin_\bW$
preserve infinite products.

 We denote the groups of morphisms in the category $\fX\Indlc^\lin$
by $\Hom^{(\fX)}({-},{-})$.

\subsection{Inverse and direct images}
\label{pro-sheaves-ind-cosheaves-inverse-direct-images-subsecn}
 In this section we specialize the discussion of the inverse and
direct images of pro-quasi-coherent pro-sheaves on ind-schemes
from~\cite[Section~3.3]{Psemten} to the case of Noetherian formal
schemes, and extend it to the dual-analogous context of locally
injective locally ind-contraherent ind-cosheaves.

 Let $\ff\:\fY\rarrow\fX$ be a morphism of Noetherian formal schemes.
 Then, for every closed subscheme of definition $\overline Y\subset\fY$,
there exists a closed subscheme of definition $\overline X\subset\fX$
such that the composition $\overline Y\rarrow\fY\overset\ff\rarrow\fX$
factorizes as $\overline Y\overset{\bar f}\rarrow\overline X\rarrow\fX$.

 The functor of inverse image of flat pro-quasi-coherent pro-sheaves
\begin{equation} \label{flat-pro-sheaves-inverse-image}
 \ff^*\:\fX\Proqc_\fl\lrarrow\fY\Proqc_\fl
\end{equation}
is defined by the rule $(\ff^*\rF)^{(\overline Y)}=
\bar f^*(\rF^{(\overline X)})$ for all $\rF\in\fX\Proqc_\fl$,
where $\bar f\:\overline Y\rarrow\overline X$ is the induced morphism of
closed subschemes of definition in $\fY$ and $\fX$ as above.
 The inverse image functor $\ff^*$
\,\eqref{flat-pro-sheaves-inverse-image} is exact (as a functor
of exact categories).

 Furthermore, the functor $\ff^*$
\,\eqref{flat-pro-sheaves-inverse-image} takes very flat
pro-quasi-coherent pro-sheaves on $\fX$ to very flat
pro-quasi-coherent pro-sheaves on $\fY$, so there is the induced
exact functor of inverse image of very flat pro-quasi-coherent
pro-sheaves
\begin{equation} \label{very-flat-pro-sheaves-inverse-image}
 \ff^*\:\fX\Proqc_\vfl\lrarrow\fY\Proqc_\vfl.
\end{equation}

 In particular, for a closed subscheme of definition $Z\subset\fX$,
denote by $\fii_{Z,\fX}\:Z\rarrow\fX$ the identity closed immersion,
viewed as a morphism of Noetherian formal schemes.
 Then, for any flat pro-quasi-coherent pro-sheaf $\rF$ on $\fX$, one
has $\fii_{Z,\fX}^*(\rF)=\rF^{(Z)}$.

 Let $\bW$ be an open covering of $\fX$ and $\bT$ be an open covering
of $\fY$ such that the morphism $\ff\:\fY\rarrow\fX$ is
$(\bW,\bT)$\+coaffine.
 Then the functor of inverse image of locally injective $\bW$\+locally
ind-contraherent ind-cosheaves
\begin{equation} \label{lin-ind-lcth-cosheaves-W-T-inverse-image}
 \ff^!\:\fX\Indlc^\lin_\bW\lrarrow\fY\Indlc^\lin_\bT
\end{equation}
is defined by the rule $(\ff^!\gJ)_{(\overline Y)}=
\bar f^!(\gJ_{(\overline X)})$ for all $\gJ\in\fX\Indlc^\lin_\bW$,
in the same notation for the morphism $\bar f\:\overline Y
\rarrow\overline X$.
 The notation $\bar f^!\:\overline X\Lcth^\lin_\bW\rarrow
\overline Y\Lcth^\lin_\bT$ stands for the inverse image
functor~\eqref{inverse-image-lcth-lin-W-T}.
 The inverse image functor $\ff^!$
\,\eqref{lin-ind-lcth-cosheaves-W-T-inverse-image} is exact (as
a functor of exact categories).

 Passing to the inductive limits of categories with respect to
the refinements of the open coverings $\bW$ and $\bT$, one obtains
the exact inverse image functor
\begin{equation} \label{lin-ind-lcth-cosheaves-inverse-image-nocovers}
 \ff^!\:\fX\Indlc^\lin\lrarrow\fY\Indlc^\lin.
\end{equation}

 In particular, for any closed subscheme of definition $Z\subset\fX$
and any locally injective locally ind-contraherent ind-cosheaf $\gJ$ on
$\fX$, one has $\fii_{Z,\fX}^!(\gJ)=\gJ_{(Z)}$ in the notation
$\fii_{Z,\fX}\:Z\rarrow\fX$ for the identity closed immersion
(as above).

\medskip

 Now let $\ff\:\fY\rarrow\fX$ be a tight morphism of Noetherian
formal schemes.
 Then, given a closed subscheme of definition $\overline X\subset\fX$,
we consider the related closed subscheme of definition $\overline Y=
\overline X\times_\fX\fY\subset\fY$ constructed in the beginning of
Section~\ref{special-inverse-images-subsecn}.

 Obviously, the composition $\overline Y\rarrow\fY\overset\ff\rarrow\fX$
factorizes as $\overline Y\overset{\bar f}\rarrow\overline X\rarrow\fX$.
 So our new notation $\overline X$ and $\overline Y$ is a particular
case of the previous one.
 If the morphism $\ff\:\fY\rarrow\fX$ is affine, flat, or very flat,
then the morphism $\bar f\:\overline Y\rarrow\overline X$ is also
affine, flat, or very flat, respectively.

 Let $(X_\gamma)_{\gamma\in\Gamma}$ be an inductive system of closed
subschemes of definition in $\fX$, indexed by a directed poset $\Gamma$,
such that $\fX=\ilim_{\gamma\in\Gamma}X_\gamma$.
 Then the construction above produces an inductive system
$(Y_\gamma)_{\gamma\in\Gamma}$ of closed subschemes of definition
$Y_\gamma=X_\gamma\times_\fX\fY\subset\fY$, indexed by the same
directed poset $\Gamma$, such that $\fY=
\ilim_{\gamma\in\Gamma}Y_\gamma$.
 Furthermore, one has $Y_\gamma=X_\gamma\times_{X_\delta}Y_\delta$
for all pairs of indices $\gamma<\delta\in\Gamma$.

\begin{lem} \label{bijective-closed-immersion-affine-base-change-lemma}
 Let $f_2\:Y_2\rarrow X_2$ be an affine morphism of Noetherian schemes
and $i_{12}\:X_1\rarrow X_2$ be a bijective closed immersion of schemes.
 Put $Y_1=X_1\times_{X_2}Y_2$; so we have a pullback diagram of schemes
$$
 \xymatrix{
  Y_1 \ar@{>->}[r]^{k_{12}} \ar[d]_{f_1} & Y_2 \ar[d]^{f_2} \\
  X_1 \ar@{>->}[r]_{i_{12}} & X_2
 }
$$
with an affine morphism $f_1\:Y_2\rarrow X_2$ and a bijective closed
immersion $k_{12}\:Y_1\rarrow Y_2$.
 Let\/ $\bW$ be an open covering of $X_2$, viewed also as an open
covering of $X_1$, and let\/ $\bT$ be an open covering of $Y_2$,
viewed also as an open covering of $Y_1$, such that the morphism~$f_2$
(or equivalently, the morphism~$f_1$) is $(\bW,\bT)$\+affine.
 Then \par
\textup{(a)} there is a natural base change isomorphism
$i_{12}^*f_2{}_*\simeq f_1{}_*k_{12}^*$ of functors
$Y_2\Qcoh\rarrow X_1\Qcoh$; \par
\textup{(b)} assuming additionally that the morphism~$f_2$ (hence
also~$f_1$) is flat, there is a natural base change isomorphism
$i_{12}^!f_2{}_!\simeq f_1{}_!k_{12}^!$ of functors
$Y_2\Lcth^\lin_\bT\rarrow X_1\Lcth^\lin_\bW$.
\end{lem}

\begin{proof}
 The assumption that the morphism~$i_{12}$ is a bijective closed
immersion is only relevant insofar as it simplifies the discussion
of the open coverings $\bW$ and $\bT$ describing the extension of
the locality of locally contraherent cosheaves.
 Otherwise $i_{12}$~could be an arbitrary morphism of schemes.
 Under the bijective closed immersion assumption of the lemma,
the morphism~$f_2$ is $(\bW,\bT)$\+affine if and only if
the morphism~$f_1$ is $(\bW,\bT)$\+affine, since a Noetherian scheme
is affine if and only if its maximal reduced closed subscheme is
affine~\cite[Exercise~III.3.1]{HarAG}, \cite[Lemma Tag~01YQ]{SP}.
 Part~(a) is a special case of~\cite[Lemma Tag~02KG]{SP}
or~\cite[Lemma~3.3(a)]{Psemten}.
 Part~(b) is a special case of the discussion of base change
in~\cite[Section~3.3]{Pcosh}.
\end{proof}

 Let $\ff\:\fY\rarrow\fX$ be a flat tight affine morphism of
Noetherian formal schemes.
 The functor of direct image of flat pro-quasi-coherent pro-sheaves
\begin{equation} \label{flat-pro-sheaves-direct-image}
 \ff_*\:\fY\Proqc_\fl\lrarrow\fX\Proqc_\fl
\end{equation}
is defined by the rule $(\ff_*\rG)^{(\overline X)}=
\bar f_*(\rG^{(\overline Y)})$ for all $\rG\in\fY\Proqc_\fl$, where
the notation is $\overline Y=\overline X\times_\fX\fY$ as above.
 Lemma~\ref{bijective-closed-immersion-affine-base-change-lemma}(a)
implies that the flat pro-quasi-coherent pro-sheaf $\ff_*\rG$ on $\fX$
is well-defined by this formula.
 The direct image functor $\ff_*$
\,\eqref{flat-pro-sheaves-direct-image} is exact and right
adjoint to the inverse image functor $\ff^*$
\,\eqref{flat-pro-sheaves-inverse-image}.

 For a very flat tight affine morphism of Noetherian formal schemes
 $\ff\:\fY\rarrow\fX$, the functor $\ff_*$
\,\eqref{flat-pro-sheaves-direct-image} takes very flat
pro-quasi-coherent pro-sheaves on $\fY$ to very flat
pro-quasi-coherent pro-sheaves on $\fX$, so there is the induced
exact functor of direct image of very flat pro-quasi-coherent
pro-sheaves
\begin{equation} \label{very-flat-pro-sheaves-direct-image}
 \ff_*\:\fY\Proqc_\vfl\lrarrow\fX\Proqc_\vfl.
\end{equation}

 Let $\fX$ be a Noetherian formal scheme with an open covering $\bW$
and $\fY$ be a Noetherian formal scheme with an open covering~$\bT$.
 Let $\ff\:\fY\rarrow\fX$ be a flat tight $(\bW,\bT)$\+affine morphism 
of formal schemes.
 Then the functor of direct image of locally injective $\bT$\+locally
ind-contraherent ind-cosheaves
\begin{equation} \label{lin-ind-lcth-cosheaves-W-T-direct-image}
 \ff_!\:\fY\Indlc^\lin_\bT\lrarrow\fX\Indlc^\lin_\bW
\end{equation}
is defined by the rule $(\ff_!\gK)_{(\overline X)}=
\bar f_!(\gK_{(\overline Y)})$ for all $\gK\in\fY\Indlc^\lin_\bT$,
where the notation is $\overline Y=\overline X\times_\fX\fY$ as above.
 The notation $\bar f_!\:\overline Y\Lcth^\lin_\bT\rarrow
\overline X\Lcth^\lin_\bW$ stands for the direct image
functor~\eqref{loc-ctrh-lin-W-T-affine-direct-image}.
 Lemma~\ref{bijective-closed-immersion-affine-base-change-lemma}(b)
implies that the locally injective $\bW$\+locally ind-contraherent
ind-cosheaf $\ff_!\gK$ on $\fX$ is well-defined by the rule above.
 The direct image functor
$\ff_!$ \,\eqref{lin-ind-lcth-cosheaves-W-T-direct-image} is exact.

 In particular, for any flat tight affine morphism of Noetherian
formal schemes $\ff\:\fY\rarrow\fX$, we have the exact functor of
direct image of locally injective ind-contraherent ind-cosheaves
\begin{equation} \label{lin-ind-ctrh-cosheaves-direct-image}
 \ff_!\:\fY\Indct^\lin\lrarrow\fX\Indct^\lin.
\end{equation}

 For any flat tight $(\bW,\bT)$\+affine $(\bW,\bT)$\+coaffine
morphism of formal schemes $\ff\:\fY\rarrow\fX$, the direct image
functor $\ff_!$ \,\eqref{lin-ind-lcth-cosheaves-W-T-direct-image}
is left adjoint to the inverse image functor
$\ff^!$~\,\eqref{lin-ind-lcth-cosheaves-W-T-inverse-image}.
 This follows from the adjunction isomorphism
of~\cite[formula~(3.2) in Section~3.3]{Pcosh} or
formula~\eqref{cosh-on-schemes-lin-direct-inverse-adjunction} above.

\subsection{Contraadjusted and cotorsion flat pro-quasi-coherent
pro-sheaves}
 We recall a couple of ring-theoretic lemmas
from~\cite[Section~1.7]{Pcosh}.

\begin{lem} \label{tensor-with-fg-cot-cta-module}
\textup{(a)} Let $R$ be a commutative ring, $M$ be a finitely
generated $R$\+module, and $C$ be a contraadjusted $R$\+module.
 Then the $R$\+module $M\ot_RC$ is contraadjusted. \par
\textup{(b)} Let $R$ be a coherent commutative ring, $M$ be
a finitely presented $R$\+module, and $C$ be a flat cotorsion
$R$\+module.
 Then the $R$\+module $M\ot_RC$ is cotorsion.
\end{lem}

\begin{proof}
 Part~(a) is~\cite[Corollary~1.7.4(a)]{Pcosh}
(cf.\ Lemma~\ref{tensor-with-finitely-presented-lemma}(c) above).
 Part~(b) is~\cite[Corollary~1.7.5(a)]{Pcosh}
(cf.\ Lemma~\ref{Hom-from-flat-to-flat-cotorsion} above).
\end{proof}

\begin{lem} \label{cta-cot-modules-stable-under-quotreductions}
\textup{(a)} Let $R$ be a commutative ring, $I\subset R$ be an ideal,
and $C$ be a contraadjusted $R$\+module.
 Then the $R/I$\+module $C/IC$ is contraadjusted. \par
\textup{(b)} Let $R\rarrow S$ be a finite morphism of Noetherian
commutative rings (i.~e., the ring $R$ is Noetherian and
the $R$\+module $S$ is finitely generated).
 Let $C$ be a flat cotorsion $R$\+module.
 Then the $S$\+module $S\ot_RC$ is (flat and) cotorsion.
\end{lem}

\begin{proof}
 Part~(a) is~\cite[Lemma~1.7.7(b)]{Pcosh} (cf.\
Lemma~\ref{contraadjusted-contramodules-lemma} above).
 Part~(b) is~\cite[Lemma~1.7.9(a)]{Pcosh}.
\end{proof}

\begin{cor} \label{tensor-product-with-coh-sheaf-cta-cot}
 Let $X$ be a semi-separated Noetherian scheme and\/ $\cM$ be
a coherent sheaf on~$X$.
 Then \par
\textup{(a)} for any flat contraadjusted quasi-coherent sheaf\/ $\cP$ 
on $X$, the quasi-coherent sheaf\/ $\cM\ot_{\cO_X}\cP$ on $X$ is
contraadjusted; \par
\textup{(b)} for any flat cotorsion quasi-coherent sheaf\/ $\cP$ 
on $X$, the quasi-coherent sheaf\/ $\cM\ot_{\cO_X}\cP$ on $X$ is
cotorsion.
\end{cor}

\begin{proof}
 Part~(a): by Proposition~\ref{flat-cta-qcoh-sheaves-prop}(b),
\,$\cP$ is a direct summand of a finitely iterated extension of
the direct images of flat contraadjusted quasi-coherent sheaves
from affine open subschemes of~$X$.
 So, let $U\subset X$ be an affine open subscheme with the open
immersion morphism $j\:U\rarrow X$.
 Then the direct image functor $j_*\:U\Qcoh\rarrow X\Qcoh$ preserves
both flatness (since the morphism~$j$ is flat and affine) and
contraadjustedness of quasi-coherent sheaves (according to
the beginning of Section~\ref{cta-and-cot-qcoh-sheaves-subsecn}
or Corollary~\ref{qcomp-ssep-direct-image-of-cta-cot-cor}(a)).

 Furthermore, the functor $\cM\ot_{\cO_X}{-}$ preserves short exact
sequences of flat quasi-coherent sheaves on $X$ (for any quasi-coherent
sheaf $\cM$ on~$X$), and the class of contraadjusted quasi-coherent
sheaves on $X$ is closed under extensions in $X\Qcoh$ (also by
the beginning of Section~\ref{cta-and-cot-qcoh-sheaves-subsecn}).
 Finally, for any flat contraadjusted quasi-coherent sheaf $\cQ$
on $U$ we have $\cM\ot_{\cO_X}j_*\cQ\simeq j_*(j^*\cM\ot_{\cO_U}\cQ)$
by the projection formula, and the quasi-coherent sheaf
$j^*\cM\ot_{\cO_U}\cQ$ on $U$ is contraadjusted by
Lemma~\ref{tensor-with-fg-cot-cta-module}(a).
 The desired assertion follows from these observations.

 The proof of part~(b) is similar but simpler (since one does not
need to deal with finitely iterated extensions, but only with finite
direct sums).
 The argument is based on
Proposition~\ref{flat-cta-qcoh-sheaves-prop}(c),
Lemma~\ref{tensor-with-fg-cot-cta-module}(b), and the properties
of cotorsion quasi-coherent sheaves mentioned in
Section~\ref{cta-and-cot-qcoh-sheaves-subsecn} (one can also refer to
Corollary~\ref{qcomp-ssep-direct-image-of-cta-cot-cor}(b)).
\end{proof}

\begin{cor} \label{cta-cot-sheaves-finite-pullbacks}
\textup{(a)} Let $i\:Z\rarrow X$ be a closed immersion of quasi-compact
semi-separated schemes.
 Then the inverse image functor $i^*\:X\Qcoh\rarrow Z\Qcoh$ takes
flat contraadjusted quasi-coherent sheaves on $X$ to flat
contraadjusted quasi-coherent sheaves on~$Z$. \par
\textup{(b)} Let $k\:Z\rarrow X$ be a finite morphism of semi-separated
Noetherian schemes.
 Then the inverse image functor $k^*\:X\Qcoh\rarrow Z\Qcoh$ takes
flat cotorsion quasi-coherent sheaves on $X$ to flat cotorsion
quasi-coherent sheaves on~$Z$.
\end{cor}

\begin{proof}
 Part~(a): we use Proposition~\ref{flat-cta-qcoh-sheaves-prop}(b)
again.
 Let $U\subset X$ be an affine open subscheme with the open
immersion morphism $j\:U\rarrow X$.
 Put $V=i^{-1}(V)\subset Z$, and denote by $j'\:V\rarrow Z$
and $i'\:V\rarrow U$ the natural morphisms.
 Then by~\cite[Lemma Tag~02KG]{SP} or~\cite[Lemma~3.3(a)]{Psemten} (see
also Lemma~\ref{bijective-closed-immersion-affine-base-change-lemma}(a)
above), we have a base change isomorphism $i^*j_*\simeq j'_*i'{}^*$
of functors $U\Qcoh\rarrow Z\Qcoh$.
 It is helpful to keep in mind that the morphism~$j$ is affine.

 The direct image functor $j'_*\:V\Qcoh\rarrow Z\Qcoh$ preserves
both flatness and contraadjustedness of quasi-coherent sheaves, as per
the proof of the previous
Proposition~\ref{tensor-product-with-coh-sheaf-cta-cot}(a).
 The inverse image functor $i'{}^*\:U\Qcoh\rarrow V\Qcoh$ for
the closed immersion of affine schemes~$i'$ also preserves both
the flatness (obviously) and contraadjustedness of quasi-coherent
sheaves, by Lemma~\ref{cta-cot-modules-stable-under-quotreductions}(a).
 Finally, the functor $i^*\:X\Qcoh\rarrow Z\Qcoh$ preserves short
exact sequences of flat quasi-coherent sheaves (for any morphism
of schemes~$i$), and both the classes of flat and contraadjusted
quasi-coherent sheaves on $Z$ are closed under extensions in $Z\Qcoh$.

 The proof of part~(b) is similar but simpler.
 The argument is based on
Proposition~\ref{flat-cta-qcoh-sheaves-prop}(c)
and Lemma~\ref{cta-cot-modules-stable-under-quotreductions}(b).
\end{proof}

 The following lemma is easy.

\begin{lem} \label{closed-subscheme-adjunction-morphism-in-qcoh-lemma}
 Let $i\:Z\rarrow X$ be a closed immersion of schemes and\/
$\cI\subset\cO_X$ be the corresponding quasi-coherent sheaf
of ideals on~$X$.
 Let\/ $\cF$ be a quasi-coherent sheaf on~$X$.
 In this context: \par
\textup{(a)} The adjunction morphism\/ $\cF\rarrow i_*i^*\cF$ is
an epimorphism in $X\Qcoh$ with the kernel\/ $\cI\cF\subset\cF$.
 So there is a short exact sequence of quasi-coherent sheaves
$0\rarrow\cI\cF\rarrow\cF\rarrow i_*i^*\cF\rarrow0$ on~$X$. \par
\textup{(b)} If the quasi-coherent sheaf\/ $\cF$ on $X$ is flat,
then the natural epimorphism of quasi-coherent sheaves\/
$\cI\ot_{\cO_X}\cF\rarrow\cI\cF$ is an isomorphism.
 So the short exact sequence from~\textup{(a)} takes the form\/
$0\rarrow\cI\ot_{\cO_X}\cF\rarrow\cF\rarrow i_*i^*\cF\rarrow0$.
\end{lem}

\begin{proof}
 Part~(a): the assertion is local in $X$, so it reduces to the case
of affine schemes, for which it means the following.
 Let $R\rarrow S$ be a surjective homomorphism of commutative rings
with the kernel ideal $I\subset R$, and let $F$ be an $R$\+module.
 Then the adjunction morphism of $R$\+modules $F\rarrow S\ot_RF$ is
surjective with the kernel submodule $IF\subset F$.
 This is obvious.

 Part~(b): applying the exact functor $\cF\ot_{\cO_X}{-}$ to
the short exact sequence $0\rarrow\cI\rarrow\cO_X\rarrow i_*\cO_Z
\rarrow0$ in $X\Qcoh$ produces the desired short exact sequence.
\end{proof}

 Let $\fX$ be a Noetherian formal scheme.
 We say that a flat pro-quasi-coherent pro-sheaf $\rP$ on $\fX$ is
\emph{contraadjusted} if the flat quasi-coherent sheaf $\rP^{(Z)}$
on $Z$ is contraadjusted for every closed subscheme of definition
$Z\subset\fX$.

 Given a filtered inductive system $(X_\gamma)_{\gamma\in\Gamma}$
of closed subschemes of definition in $\fX$ such that
$\fX=\ilim_{\gamma\in\Gamma}X_\gamma$, it suffices to check the latter
contraadjustedness condition for the closed subschemes $Z=X_\gamma
\subset\fX$, \,$\gamma\in\Gamma$.
 This follows from Corollary~\ref{cta-cot-sheaves-finite-pullbacks}(a).
 
 Similarly, we say that a flat pro-quasi-coherent pro-sheaf $\rP$ on
$\fX$ is \emph{cotorsion} if the flat quasi-coherent sheaf $\rP^{(Z)}$
on $Z$ is cotorsion for every closed subscheme of definition
$Z\subset\fX$.

 Given a filtered inductive system $(X_\gamma)_{\gamma\in\Gamma}$
of closed subschemes of definition in $\fX$ such that
$\fX=\ilim_{\gamma\in\Gamma}X_\gamma$, it suffices to check the latter
cotorsion condition for the closed subschemes $Z=X_\gamma\subset\fX$,
\,$\gamma\in\Gamma$.
 This follows from Corollary~\ref{cta-cot-sheaves-finite-pullbacks}(b).

 We denote the full subcategories of contraadjusted and cotorsion
flat pro-quasi-coherent pro-sheaves by $\fX\Proqc_\fl^\cta$
and $\fX\Proqc_\fl^\cot\subset\fX\Proqc_\fl$.
 Similarly, the full subcategory of very flat contraadjusted
pro-quasi-coherent pro-sheaves is denoted by $\fX\Proqc_\vfl^\cta
=\fX\Proqc_\vfl\cap\fX\Proqc^\cta\subset\fX\Proqc_\fl$.

\begin{exs} \label{cta-cot-flat-pro-sheaves-on-affine}
 Let $\fU=\Spf\fR$ be an affine Noetherian formal scheme.
 In this context:

 (1)~The equivalence of exact categories $\fR\Contra_\fl\simeq
\fU\Proqc_\fl$ from
Example~\ref{flat-pro-sheaves-on-affine=flat-contramodules}
identifies the full subcategory of flat cotorsion contramodule
$\fR$\+modules in $\fR\Contra_\fl$ with the full subcategory
of flat cotorsion pro-quasi-coherent pro-sheaves on~$\fU$,
$$
 \fR\Contra_\fl^\cot\simeq\fU\Proqc_\fl^\cot.
$$
 This is a restatement of
Lemma~\ref{flat-cotorsion-contramods-in-terms-of-reductions}\,%
(1)\,$\Leftrightarrow$\,(4).

 (2)~The same equivalence of exact categories $\fR\Contra_\fl\simeq
\fU\Proqc_\fl$ identifies the full subcategory of flat contraadjusted
contramodule $\fR$\+modules in $\fR\Contra_\fl$ with the full
subcategory of flat contraadjusted pro-quasi-coherent pro-sheaves
on~$\fU$,
$$
 \fR\Contra_\fl^\cta\simeq\fU\Proqc_\fl^\cta.
$$
 This follows from
Lemma~\ref{contraadjusted-contramodules-lemma}\,%
(4)\,$\Leftrightarrow$\,(8)\,$\Leftrightarrow$\,(9) together with
Lemma~\ref{flat-quotseparated-contramodules-separated}.

 (3)~As a particular case of~(2), the equivalence of exact categories
$\fR\Contra_\vfl\simeq\fU\Proqc_\vfl$ from
Example~\ref{flat-pro-sheaves-on-affine=flat-contramodules}
identifies the full subcategory of very flat contraadjusted
contramodule $\fR$\+modules in $\fR\Contra_\vfl$ with the full
subcategory of very flat contraadjusted pro-quasi-coherent pro-sheaves
on~$\fU$,
$$
 \fR\Contra_\vfl^\cta\simeq\fU\Proqc_\vfl^\cta.
$$
\end{exs} 

 Let $\fX$ be a Noetherian formal scheme.
 Clearly, all the three full subcategories $\fX\Proqc_\fl^\cot
\subset\fX\Proqc_\fl^\cta$ and $\fX\Proqc_\vfl^\cta\subset
\fX\Proqc_\fl^\cta$ are closed under extensions in $\fX\Proqc_\fl$,
so they inherit exact category structures.
 It will follow from the next lemma that all the short exact sequences
in $\fX\Proqc_\fl^\cot$ and $\fX\Proqc_\vfl^\cta$ are actually split;
but the exact category structure on $\fX\Proqc_\fl^\cta$ is nontrivial.

 We denote the $\Ext$ groups in the exact category
$\fX\Proqc_\fl$ by $\Ext_{(\fX)}^*({-},{-})$.

\begin{lem} \label{pro-qc-pro-sheaves-Ext1-vanishing-lemma}
 Let\/ $\fX$ be a semi-separated Noetherian formal scheme. \par
\textup{(a)} For any very flat pro-quasi-coherent pro-sheaf\/ $\rF$
and any flat contraadjusted pro-quasi-coherent pro-sheaf\/ $\rP$ on\/
$\fX$, one has\/ $\Ext_{(\fX)}^1(\rF,\rP)=0$; \par
\textup{(b)} For any flat pro-quasi-coherent pro-sheaf\/ $\rF$ and
any flat cotorsion pro-quasi-coherent pro-sheaf\/ $\rP$ on\/ $\fX$,
one has\/ $\Ext_{(\fX)}^1(\rF,\rP)=0$.
 So every flat cotorsion pro-quasi-coherent pro-sheaf on\/ $\fX$ is
an injective object of the exact category\/ $\fX\Proqc_\fl$.
\end{lem}

\begin{proof}
 Part~(a): let $(X_n)_{n=1}^\infty$ be an inductive system of closed
subschemes of definition in $\fX$, indexed by the ordered set of
positive integers, such that $\fX=\ilim_{n\ge1}X_n$.
 For example, one can take $X_n=\overline X^{(n)}$ in the notation of
Section~\ref{formal-schemes-as-ind-schemes-subsecn}.
 Denote by $i_{n,n+1}\:X_n\rarrow X_{n+1}$ the transition morphisms.
 For every flat pro-quasi-coherent pro-sheaf $\rG$ on $\fX$, put
$\rG^{(n)}=\rG^{(X_n)}$ for brevity.

 Let $0\rarrow\rP\rarrow\rG\rarrow\rF\rarrow0$ be a short exact
sequence in $\fX\Proqc_\fl$.
 For every $n\ge1$, we have a split short exact sequence
$0\rarrow\rP^{(n)}\rarrow\rG^{(n)}\rarrow\rF^{(n)}\rarrow0$
in $X_n\Qcoh_\fl$.
 We need to contruct a sequence of splittings $s_n\:\rF^{(n)}\rarrow
\rG^{(n)}$ that would be compatible as $n$~varies, in the sense
that $i_{n,n+1}^*(s_{n+1})=s_n$ for all $n\ge1$.
 Proceeding by induction, assume that a compatible system of
splittings $s_1$,~\dots, $s_n$ has been constructed already.
 Let us construct a splitting $s_{n+1}$ compatible with~$s_n$.

 Pick an arbitrary splitting $s_{n+1}'\:\rF^{(n+1)}\rarrow\rG^{(n+1)}$
of the split short exact sequence $0\rarrow\rP^{(n+1)}\rarrow\rG^{(n+1)}
\rarrow\rF^{(n+1)}\rarrow0$ in $X_{n+1}\Qcoh_\fl$.
 Then the difference $s_n-i_{n,n+1}^*(s'_{n+1})\:\rF^{(n)}\rarrow
\rG^{(n)}$ factorizes through the admissible monomorphism
$\rP^{(n)}\rarrow\rG^{(n)}$; so we obtain a morphism of
quasi-coherent sheaves $t_n\:\rF^{(n)}\rarrow\rP^{(n)}$.
 We claim that any such morphism~$t_n$ can be lifted to a morphism
$t_{n+1}\:\rF^{(n+1)}\rarrow\rP^{(n+1)}$, i.~e., one has
$t_n=i_{n,n+1}^*(t_{n+1})$.
 If this is the case, then it suffices to put
$s_{n+1}=s'_{n+1}+t_{n+1}$ (where the composition
$\rF^{(n+1)}\xrightarrow{t_{n+1}}\rP^{(n+1)}\rarrow\rG^{(n+1)}$
is still denoted by $t_{n+1}$ by abuse of notation).

 By adjunction, the given morphism
$t_n\:i_{n,n+1}^*\rF^{(n+1)}\rarrow i_{n,n+1}^*\rP^{(n+1)}$
in $X_n\Qcoh$ corresponds to a morphism
$\tilde t_n\:\rF^{(n+1)}\rarrow i_{n,n+1}{}_*i_{n,n+1}^*\rP^{(n+1)}$
in $X_{n+1}\Qcoh$.
 We need to show that the morphism~$\tilde t_n$ factorizes through
the natural adjunction (epi)morphism $\rP^{(n+1)}\rarrow
i_{n,n+1}{}_*i_{n,n+1}^*\rP^{(n+1)}$ in $X_{n+1}\Qcoh$.

 By Lemma~\ref{closed-subscheme-adjunction-morphism-in-qcoh-lemma}(b),
we have a short exact sequence $0\rarrow\cI\ot_{\cO_{X_{n+1}}}
\rP^{(n+1)}\rarrow\rP^{(n+1)}\rarrow i_{n,n+1}{}_*i_{n,n+1}^*\rP^{(n+1)}
\rarrow0$ in $X_{n+1}\Qcoh$ (since the quasi-coherent sheaf
$\rP^{(n+1)}$ on $X_{n+1}$ is flat).
 Here $\cI\subset\cO_{X_{n+1}}$ denotes the quasi-coherent sheaf of
ideals on $X_{n+1}$ corresponding to the closed subscheme
$i_{n,n+1}(X_n)\subset X_{n+1}$.
 Finally, the quasi-coherent sheaf $\cI\ot_{\cO_{X_{n+1}}}\rP^{(n+1)}$
on $X_{n+1}$ is contraadjusted by
Corollary~\ref{tensor-product-with-coh-sheaf-cta-cot}(a), while
the quasi-coherent sheaf $\rF^{(n+1)}$ on $X_{n+1}$ is very flat
by assumption.
 Thus the morphism of quasi-coherent sheaves
$\tilde t_n\:\rF^{(n+1)}\rarrow i_{n,n+1}{}_*i_{n,n+1}^*\rP^{(n+1)}$
can be lifted to a morphism of quasi-coherent sheaves
$t_{n+1}\:\allowbreak\rF^{(n+1)}\allowbreak\rarrow\rP^{(n+1)}$
on $X_{n+1}$, and we are done.

 The proof of part~(b) is similar and based on
Corollary~\ref{tensor-product-with-coh-sheaf-cta-cot}(b).
\end{proof}

 Let $\ff\:\fY\rarrow\fX$ be a flat tight affine morphism of
Noetherian formal schemes.
 Then it follows from the construction together with the discussion of
the direct images of contraadjusted quasi-coherent sheaves in
the beginning of Section~\ref{cta-and-cot-qcoh-sheaves-subsecn} or in
Corollary~\ref{qcomp-ssep-direct-image-of-cta-cot-cor}(a)
or~\ref{qcomp-ssep-direct-image-of-flat-cta-cor}(b) that
the direct image functor $\ff_*\:\fY\Proqc_\fl\rarrow\fX\Proqc_\fl$
\,\eqref{flat-pro-sheaves-direct-image} takes flat contraadjusted
pro-quasi-coherent pro-sheaves on $\fY$ to flat contraadjusted
pro-quasi-coherent pro-sheaves on~$\fX$.
 Similarly, it follows from the discussion of direct
images of cotorsion quasi-coherent sheaves in
Section~\ref{cta-and-cot-qcoh-sheaves-subsecn} or from
Corollary~\ref{qcomp-ssep-direct-image-of-cta-cot-cor}(b)
or~\ref{qcomp-ssep-direct-image-of-flat-cta-cor}(c)
that the direct image functor $\ff_*$
\,\eqref{flat-pro-sheaves-direct-image} takes flat cotorsion
pro-quasi-coherent pro-sheaves on $\fY$ to flat cotorsion
pro-quasi-coherent pro-sheaves on~$\fX$.

\begin{cor} \label{enough-injectives-in-flat-vfl-pro-sheaves}
 Let\/ $\fX$ be a semi-separated Noetherian formal scheme with a finite
affine open covering\/ $\fX=\bigcup_{\alpha=1}^N\fU_\alpha$.
 In this setting: \par
\textup{(a)} There are enough injective objects in the exact category\/
$\fX\Proqc_\fl$.
 A flat pro-quasi-coherent pro-sheaf is an injective object in\/
$\fX\Proqc_\fl$ if and only if it is a flat cotorsion pro-quasi-coherent
pro-sheaf.
 A flat pro-quasi-coherent pro-sheaf on\/ $\fX$ is cotorsion if and
only if it is a direct summand of a finite direct sum of the direct
images of flat cotorsion pro-quasi-coherent pro-sheaves
from\/~$\fU_\alpha$. \par
\textup{(b)} There are enough injective objects in the exact category\/
$\fX\Proqc_\vfl$.
 A very flat pro-quasi-coherent pro-sheaf is an injective object in\/
$\fX\Proqc_\vfl$ if and only if it is a very flat contraadjusted
pro-quasi-coherent pro-sheaf.
 A flat pro-quasi-coherent pro-sheaf on\/ $\fX$ is very flat and
contraadjusted if and only if it is a direct summand of a finite direct
sum of the direct images of very flat contraadjusted pro-quasi-coherent
pro-sheaves from\/~$\fU_\alpha$.
\end{cor}

\begin{proof}
 Part~(a): first of all, flat cotorsion pro-quasi-coherent pro-sheaves
are injective objects in $\fX\Proqc_\fl$ by
Lemma~\ref{pro-qc-pro-sheaves-Ext1-vanishing-lemma}(b), and
the direct images of flat cotorsion pro-quasi-coherent pro-sheaves
from $\fU_\alpha$ are flat cotorsion pro-quasi-coherent pro-sheaves
on $\fX$ according to the paragraph preceding the present corollary.
 Thus it remains to prove that every object of $\fX\Proqc_\fl$ is
an admissible subobject of an object of the form
$\bigoplus_{\alpha=1}^N\fj_\alpha{}_*\rQ_\alpha$, where
$\fj_\alpha\:\fU_\alpha\rarrow\fX$ are the open immersion morphisms
and $\rQ_\alpha$ are flat cotorsion pro-quasi-coherent pro-sheaves
on~$\fU_\alpha$.

 Let us first consider the case of an affine Noetherian formal
scheme~$\fU$.
 We need to prove that every object $\fU\Proqc_\fl$ is an admissible
subobject of a flat cotorsion pro-quasi-coherent pro-sheaf on~$\fU$.
 Let $\fR=\fO_\fU(\fU)$ be the complete, separated adic Noetherian ring
for which $\fU=\Spf\fR$.
 By Examples~\ref{flat-pro-sheaves-on-affine=flat-contramodules}
and~\ref{cta-cot-flat-pro-sheaves-on-affine}(1), we have an equivalence
of exact categories $\fU\Proqc_\fl\simeq\fR\Contra_\fl$ identifying
the full subcategory $\fU\Proqc_\fl^\cot\subset\fU\Proqc_\fl$ with
the full subcategory $\fR\Contra_\fl^\cot\subset\fR\Contra_\fl$.
 Now it remains to refer to
Proposition~\ref{flat-cotorsion-pair-in-quotseparated}(b) together
with the fact that the full subcategory $\fR\Contra_\fl$ is closed
under extensions in $\fR\Contra$ (by
Lemma~\ref{flat-quotseparated-contramodules-well-behaved}(a)
or Corollary~\ref{flat-contramodules-cor}).
 Lemma~\ref{cotorsion-pair-kernel-as-injectives-projectives}(a)
is relevant here.

 In the general case of a semi-separated Noetherian formal scheme $\fX$,
let $\rF$ be a flat pro-quasi-coherent pro-sheaf on~$\fX$.
 By the previous paragraph, for every index~$\alpha$ there exists
a flat cotorsion pro-quasi-coherent pro-sheaf $\rQ_\alpha$ on
$\fU_\alpha$ together with an admissible monomorphism
$\fj_\alpha^*\rF\rarrow\rQ_\alpha$ in $\fU\Proqc_\fl$.
 By the adjunction of the inverse image functors~$\fj_\alpha^*$ with
the direct image functors~$\fj_\alpha{}_*$ (see
Section~\ref{pro-sheaves-ind-cosheaves-inverse-direct-images-subsecn}),
we have the related morphism $g\:\rF\rarrow\bigoplus_{\alpha=1}^N
\fj_\alpha{}_*\rQ_\alpha$.
 The morphism~$g$ is an admissible monomorphism of flat
pro-quasi-coherent pro-sheaves in restriction to every $\fU_\alpha$,
and it follows easily that $g$~is an admissible monomorphism
in $\fX\Proqc_\fl$.

 The proof of part~(b) is similar.
 All very flat contraadjusted pro-quasi-coherent pro-sheaves are
injective objects of $\fX\Proqc_\vfl$ by
Lemma~\ref{pro-qc-pro-sheaves-Ext1-vanishing-lemma}(a).

 In the case of an affine Noetherian formal scheme $\fU$,
by Examples~\ref{flat-pro-sheaves-on-affine=flat-contramodules}
and~\ref{cta-cot-flat-pro-sheaves-on-affine}(3), we have an equivalence
of exact categories $\fU\Proqc_\vfl\simeq\fR\Contra_\vfl$ identifying
the full subcategory $\fU\Proqc_\vfl^\cta\subset\fU\Proqc_\vfl$ with
the full subcategory $\fR\Contra_\vfl^\cta\subset\fR\Contra_\vfl$.
 Now it remains to refer to
Proposition~\ref{very-flat-cotorsion-pair-in-quotseparated}(b) together
with the fact that the full subcategory $\fR\Contra_\vfl$ is closed
under extensions in $\fR\Contra$ (as mentioned in
Section~\ref{prelim-veryflat-and-contraadjusted-subsecn}).
 Once again,
Lemma~\ref{cotorsion-pair-kernel-as-injectives-projectives}(a)
is relevant.
\end{proof}

\begin{lem} \label{flat-pro-sheaves-very-flat-precover}
 Let\/ $\fX$ be a semi-separated Noetherian formal scheme with
a finite affine open covering\/ $\fX=\bigcup_{\alpha=1}^N\fU_\alpha$.
 Then any flat pro-quasi-coherent pro-sheaf\/ $\rG$ on\/ $\fX$ can be
included in an (admissible) short exact sequence\/ $0\rarrow\rP
\rarrow\rF\rarrow\rG\rarrow0$ in the exact category\/ $\fX\Proqc_\fl$,
where\/ $\rF$ is a very flat pro-quasi-coherent pro-sheaf on\/ $\fX$
and\/ $\rP$ is a finitely iterated extension of the direct images of
flat contraadjusted pro-quasi-coherent pro-sheaves from\/~$\fU_\alpha$.
\end{lem}

\begin{proof}
 The argument is similar to~\cite[proof of Lemma~4.1.1]{Pcosh} and
dual-analogous to the proof of
Lemma~\ref{locally-cotorsion-preenvelope-lemma} above.

 Arguing by induction on $1\le\beta\le N$, we consider the open
formal subscheme $\fV=\bigcup_{\alpha<\beta}\fU_\alpha$ with
the identity open immersion morphism $\fh\:\fV\rarrow\fX$.
 Assume that we have constructed a short exact sequence
$0\rarrow\rQ\rarrow\rK\rarrow\rG\rarrow0$ of flat pro-quasi-coherent
pro-sheaves on $\fX$ such that the restriction $\fh^*\rK$ of
the flat pro-quasi-coherent pro-sheaf $\rK$ to the open formal
subscheme $\fV\subset\fX$ is very flat, while the flat
pro-quasi-coherent pro-sheaf $\rQ$ on $\fX$ is a finitely iterated
extension of the direct images of flat contraadjusted
pro-quasi-coherent pro-sheaves from the affine open formal subschemes
$\fU_\alpha\subset\fX$, \,$\alpha<\beta$.
 When $\beta=1$, it suffices to take $\rK=\rG$ and $\rQ=0$ for
the induction base.
 Set $\fU=\fU_\beta$, and denote by $\fj\:\fU\rarrow\fX$ the identity
open immersion morphism.

 The assertion of the lemma for the affine formal scheme $\fU$ has
to be established first.
 In this case, by
Examples~\ref{flat-pro-sheaves-on-affine=flat-contramodules}
and~\ref{cta-cot-flat-pro-sheaves-on-affine}(2), we have an equivalence
of exact categories $\fU\Proqc_\fl\simeq\fR\Contra_\fl$ identifying
the full subcategory $\fU\Proqc_\vfl\subset\fU\Proqc_\fl$ with
the full subcategory $\fR\Contra_\vfl\subset\fR\Contra_\fl$ and
the full subcategory $\fU\Proqc_\fl^\cta\subset\fU\Proqc_\fl$ with
the full subcategory $\fR\Contra_\fl^\cta\subset\fR\Contra_\fl$.
 So we can refer to
Proposition~\ref{very-flat-cotorsion-pair-in-quotseparated}(a) together
with the fact that the full subcategory $\fR\Contra_\fl$ is closed
under kernels of epimorphisms in $\fR\Contra$ (by
Lemma~\ref{flat-quotseparated-contramodules-well-behaved}(a)).
 The lemma about restrictions of complete cotorsion pairs to
exact subcategories~\cite[Lemma~1.5(a)]{Pal},
\cite[Lemma~B.1.4(a)]{Pcosh} is relevant here.
 This collection of references implies the existence of a short exact
sequence $0\rarrow\rR\rarrow\rH\rarrow\fj^*\rK\rarrow0$ in
$\fU\Proqc_\fl$ with a very flat pro-quasi-coherent pro-sheaf $\rH$
and a flat contraadjusted pro-quasi-coherent pro-sheaf $\rR$ on~$\fU$.

 Consider the direct image $0\rarrow\fj_*\rR\rarrow\fj_*\rH\rarrow
\fj_*\fj^*\rK\rarrow0$ of the short exact sequence $0\rarrow\rR
\rarrow\rH\rarrow\fj^*\rK\rarrow0$ with respect to the affine open
immersion~$\fj$, and take its pullback with respect to the adjunction
morphism $\rK\rarrow\fj_*\fj^*\rK$ in the exact category
$\fX\Proqc_\fl$.
 Similarly to~\cite[proof of Lemma~4.1.1]{Pcosh}, one shows that,
in the resulting short exact sequence $0\rarrow\fj_*\rR\rarrow\rL
\rarrow\rK\rarrow0$, the flat pro-quasi-coherent pro-sheaf $\rL$
on $\fX$ is very flat in restriction to $\fU\cup\fV$.
 The locality of the very flatness property of quasi-coherent
sheaves (see~\cite[Lemma~1.2.6(a) and Section~1.10]{Pcosh} and
Section~\ref{loc-injective-on-schemes-subsecn} above) needs to be
used here, together with the fact that the full subcategory of
very flat pro-quasi-coherent pro-sheaves is closed under extensions
and kernels of admissible epimorphisms in the exact category of
flat pro-quasi-coherent pro-sheaves.
 One also has to use the preservation of the very flatness by
the direct images of flat pro-quasi-coherent pro-sheaves with respect
to affine open immersions of formal
schemes~\eqref{very-flat-pro-sheaves-direct-image}.

 Finally, the composition of the admissible epimorphisms of flat
pro-quasi-coherent pro-sheaves $\rL\rarrow\rK\rarrow\rG$ on $\fX$
is an admissible epimorphism whose kernel is an extension of
the flat pro-quasi-coherent pro-sheaves $\rQ$ and $\fj_*\rR$ on $\fX$,
hence also a fintitely iterated extension of the direct images of
flat contraadjusted pro-quasi-coherent pro-sheaves from the affine
open formal subschemes $\fU_\alpha\subset\fX$, \,$\alpha\le\beta$.
 This finishes the argument for the induction step.
\end{proof}

\begin{cor} \label{flat-pro-sheaves-flat-ctra-preenvelope}
 Let\/ $\fX$ be a semi-separated Noetherian formal scheme with
a finite affine open covering\/ $\fX=\bigcup_{\alpha=1}^N\fU_\alpha$.
 Then any flat pro-quasi-coherent pro-sheaf\/ $\rG$ on\/ $\fX$ can be
included in a short exact sequence\/ $0\rarrow\rG\rarrow\rP\rarrow\rF
\rarrow0$ in the exact category\/ $\fX\Proqc_\fl$, where\/ $\rF$ is
a very flat pro-quasi-coherent pro-sheaf on\/ $\fX$ and\/ $\rP$ is
a finitely iterated extension of the direct images of flat
contraadjusted pro-quasi-coherent pro-sheaves from\/~$\fU_\alpha$.
\end{cor}

\begin{proof}
 The assertion follows from
Lemma~\ref{flat-pro-sheaves-very-flat-precover} and
Corollary~\ref{enough-injectives-in-flat-vfl-pro-sheaves}(a)
by virtue of (the construction from proof of) the Salce
lemma~\cite{Sal}, \cite[Lemma~B.1.1(b)]{Pcosh}.
 A brief discussion can be found in
Section~\ref{appx-cotorsion-pairs-subsecn} below.
\end{proof}

\begin{cor} \label{cta-flat-pro-sheaves-cor}
 Let\/ $\fX$ be a semi-separated Noetherian formal scheme with
a finite affine open covering\/ $\fX=\bigcup_{\alpha=1}^N\fU_\alpha$.
 In this setting: \par
\textup{(a)} For any flat pro-quasi-coherent pro-sheaf\/ $\rG$ on\/
$\fX$ there exists a short exact sequence\/ $0\rarrow\rP\rarrow\rF
\rarrow\rG\rarrow0$ in\/ $\fX\Proqc_\fl$ with a very flat
pro-quasi-coherent pro-sheaf\/ $\rF$ and a flat contraadjusted
pro-quasi-coherent pro-sheaf\/~$\rP$. \par
\textup{(b)} For any flat pro-quasi-coherent pro-sheaf\/ $\rG$ on\/
$\fX$ there exists a short exact sequence\/ $0\rarrow\rG\rarrow\rP
\rarrow\rF\rarrow0$ in\/ $\fX\Proqc_\fl$ with a flat contraadjusted
pro-quasi-coherent pro-sheaf\/ $\rP$ and a very flat
pro-quasi-coherent pro-sheaf\/~$\rF$. \par
\textup{(c)} A flat pro-quasi-coherent pro-sheaf on\/ $\fX$ is
contraadjusted if and only if it is a direct summand of a finitely
iterated extension of the direct images of flat contraadjusted
pro-quasi-coherent pro-sheaves from\/~$\fU_\alpha$.
\end{cor}

\begin{proof}
 This is a formal scheme version of~\cite[Corollary~4.1.4]{Pcosh}.
 The ``if'' implication in part~(c) follows from the preservation of
contraadjustedness of flat pro-quasi-coherent pro-sheaves by
the direct images with respect to flat tight affine morphisms of
Noetherian formal schemes (as per the remarks above in this section).
 After the ``if'' assertion of part~(c) is shown, part~(a) follows
from Lemma~\ref{flat-pro-sheaves-very-flat-precover} and part~(b)
follows from Corollary~\ref{flat-pro-sheaves-flat-ctra-preenvelope}.
 The ``only if'' implication in part~(c) follows from
Corollary~\ref{flat-pro-sheaves-flat-ctra-preenvelope} and
Lemma~\ref{pro-qc-pro-sheaves-Ext1-vanishing-lemma}(a) by virtue of
the direct summand lemma (essentially, Lemma~\ref{direct-summand-lemma}
below).
\end{proof}

 The result of Corollary~\ref{cta-flat-pro-sheaves-cor}(c) can be
rephrased by saying that the class of flat contraadjusted
pro-quasi-coherent pro-sheaves on semi-separated Noetherian formal
schemes is \emph{antilocal} in the sense similar
to~\cite[Section~4]{Pal}.

 By Lemma~\ref{flat-pro-sheaves-very-flat-precover}, the full
subcategory $\fX\Proqc_\vfl$ is resolving in $\fX\Proqc_\fl$.
 By Corollaries~\ref{enough-injectives-in-flat-vfl-pro-sheaves}(a)
and~\ref{flat-pro-sheaves-flat-ctra-preenvelope}, the full
subcategories $\fX\Proqc_\fl^\cot$ and $\fX\Proqc_\fl^\cta$ are
coresolving is $\fX\Proqc_\fl$.
 By Corollaries~\ref{enough-injectives-in-flat-vfl-pro-sheaves}(b),
the full subcategory $\fX\Proqc_\vfl^\cta$ is coresolving in
$\fX\Proqc_\vfl$.
 It follows by virtue of
Proposition~\ref{infinite-co-resolutions-minus-plus-derived} that
the $\Ext$ groups in all the mentioned exact (sub)categories agree
with each other.

\begin{cor} \label{cta-flat-pro-sheaves-characterized}
 Let\/ $\fX$ be a semi-separated Noetherian formal scheme and\/ $\rP$
be a flat pro-quasi-coherent pro-sheaf on\/~$\fX$.
 Then the following conditions are equivalent:
\begin{enumerate}
\item the flat pro-quasi-coherent pro-sheaf\/ $\rP$ is contraadjusted;
\item the functor\/ $\Hom_{(\fX)}({-},\rP)$ takes short exact
sequences of very flat pro-quasi-coherent pro-sheaves to short exact
sequences of abelian groups;
\item one has\/ $\Ext^1_{(\fX)}(\rF,\rP)=0$ for all very flat
pro-quasi-coherent pro-sheaves\/ $\rF$ on\/~$\fX$;
\item one has\/ $\Ext^n_{(\fX)}(\rF,\rP)=0$ for all very flat
pro-quasi-coherent pro-sheaves\/ $\rF$ on\/ $\fX$ and all $n\ge1$.
\end{enumerate}
\end{cor}

\begin{proof}
 (1)~$\Longrightarrow$~(3) This is
Lemma~\ref{pro-qc-pro-sheaves-Ext1-vanishing-lemma}(a).

 (3)~$\Longrightarrow$~(1) Follows from
Corollary~\ref{cta-flat-pro-sheaves-cor}(b) by virtue of the direct
summand lemma (Lemma~\ref{direct-summand-lemma} below).

 (2)~$\Longleftrightarrow$~(3)~$\Longleftrightarrow$~(4) Follows from
the facts that there are enough very flat pro-quasi-coherent
pro-sheaves in $\fX\Proqc_\fl$ (a weak form of
Corollary~\ref{cta-flat-pro-sheaves-cor}(a)) and the full subcategory
$\fX\Proqc_\vfl$ is closed under kernels of admissible epimorphisms
in $\fX\Proqc_\fl$ by virtue of
Lemma~\ref{Hom-from-to-short-exact-sequences-lemma}(a) below.
\end{proof}

 In the terminology of Section~\ref{appx-cotorsion-pairs-subsecn} below,
the results of Corollaries~\ref{cta-flat-pro-sheaves-cor}(a\+-b)
and~\ref{cta-flat-pro-sheaves-characterized} can be summarized by
saying that, for any semi-separated Noetherian formal scheme $\fX$,
the pair of classes (very flat pro-quasi-coherent pro-sheaves on $\fX$,
flat contraadjusted pro-quasi-coherent pro-sheaves on $\fX$) is
a \emph{hereditary complete cotorsion pair} in the exact category
of flat pro-quasi-coherent pro-sheaves $\fX\Proqc_\fl$.
 The direct summand lemma (Lemma~\ref{direct-summand-lemma})
is also relevant here.

 By Lemmas~\ref{cotorsion-pair-kernel-as-injectives-projectives}(b)
and~\ref{cotorsion-pair-the-same-projectives-injectives}(b) (using
also Corollary~\ref{enough-injectives-in-flat-vfl-pro-sheaves}(a)),
it follows that the exact category of flat contraadjusted
pro-quasi-coherent pro-sheaves $\fX\Proqc_\fl^\cta$ on a semi-separated
Noetherian formal scheme $\fX$ has both enough projective and enough
injective objects.
 The full subcategory of very flat contraadjusted pro-quasi-coherent
pro-sheaves $\fX\Proqc_\vfl^\cta$ is the full subcategory of projective
objects in $\fX\Proqc_\fl^\cta$, while the full subcategory of flat
cotorsion pro-quasi-coherent pro-sheaves $\fX\Proqc_\fl^\cot$ is
the full subcategory of injective objects in $\fX\Proqc_\fl^\cta$.

\begin{lem} \label{flat-contraadjusted-base-change-lemma}
 Let $f_2\:Y_2\rarrow X_2$ be a flat morphism of semi-separated
Noetherian schemes and $i_{12}\:X_1\rarrow X_2$ be a closed immersion
of schemes.
 Put $Y_1=X_1\times_{X_2}Y_2$; so we have a pullback diagram of schemes
$$
 \xymatrix{
  Y_1 \ar@{>->}[r]^{k_{12}} \ar[d]_{f_1} & Y_2 \ar[d]^{f_2} \\
  X_1 \ar@{>->}[r]_{i_{12}} & X_2
 }
$$
with a flat morphism $f_1\:Y_2\rarrow X_2$ and a closed immersion
$k_{12}\:Y_1\rarrow Y_2$.
 Then, for any quasi-coherent sheaf\/ $\cP$ on $Y_2$, there is
an natural base change comparison morphism
\begin{equation} \label{flat-closedimm-qcoh-basechange-morphism}
 i_{12}^*f_2{}_*\cP\lrarrow f_1{}_*k_{12}^*\cP
\end{equation}
of quasi-coherent sheaves on~$X_1$.
 For any flat contraadjusted quasi-coherent sheaf\/ $\cP$ on $Y_2$,
the morphism~\eqref{flat-closedimm-qcoh-basechange-morphism} is
an isomorphism.
 So there is a natural base change isomorphism
$i_{12}^*f_2{}_*\simeq f_1{}_*k_{12}^*$ of functors
$Y_2\Qcoh_\fl^\cta\rarrow X_1\Qcoh_\fl^\cta$.
\end{lem}

\begin{proof}
 The construction of the base change comparison
morphism~\eqref{flat-closedimm-qcoh-basechange-morphism} works
quite generally for any quasi-compact morphism of schemes~$f_2$ and
any morphism of schemes~$i_{12}$.
 The desired morphism~\eqref{flat-closedimm-qcoh-basechange-morphism}
corresponds by adjunction to a morphism
$f_1^*i_{12}^*f_2{}_*\cP\simeq k_{12}^*f_2^*f_2{}_*\cP\lrarrow
k_{12}^*\cP$ in $Y_1\Qcoh$, which is obtainable by applying
the inverse image functor~$k_{12}^*$ to the adjunction morphism
$f_2^*f_2{}_*\cP\rarrow\cP$.
 The same morphism~\eqref{flat-closedimm-qcoh-basechange-morphism}
corresponds by adjunction to a morphism
$f_2{}_*\cP\rarrow f_2{}_*k_{12}{}_*\cP\simeq
i_{12}{}_*f_1{}_*k_{12}^*\cP$ in $X_2\Qcoh$, which is obtainable
by applying the direct image functor~$f_2{}_*$ to the adjunction
morphism $\cP\rarrow k_{12}{}_*k_{12}^*\cP$.

 To prove the second assertion of the lemma, notice that we have
exact functors $f_2{}_*\:Y_2\Qcoh_\fl^\cta\rarrow X_2\Qcoh_\fl^\cta$
by Corollaries~\ref{qcomp-ssep-direct-image-of-cta-cot-cor}(a)
and~\ref{qcomp-ssep-direct-image-of-flat-cta-cor}(b), and
$i_{12}^*\:X_2\Qcoh_\fl^\cta\rarrow X_1\Qcoh^\fl_\cta$ by
Corollary~\ref{cta-cot-sheaves-finite-pullbacks}(a).
 Similarly, we also have exact functors
$f_1{}_*\:Y_1\Qcoh_\fl^\cta\rarrow X_1\Qcoh_\fl^\cta$ and
$k_{12}^*\:Y_2\Qcoh_\fl^\cta\rarrow Y_1\Qcoh^\fl_\cta$.

 By Proposition~\ref{flat-cta-qcoh-sheaves-prop}(b), \,$\cP$ is
a direct summand of a finitely iterated extension of the direct
images of flat contraadjusted quasi-coherent sheaves from
affine open subschemes $V\subset Y_2$.
 Since a morphism of finitely filtered objects (e.~g., in an exact
category) is an isomorphism whenever it is an isomorphism on
the successive quotients, the question reduces to the case of
$\cP=j_*\cQ$, where $\cQ$ is a flat contraadjusted quasi-coherent
sheaf on $V$ and $j\:V\rarrow Y_2$ is the open immersion morphism.

 Now both the morphisms $j\:V\rarrow Y_2$ and $f_2j\:V\rarrow X_2$
are flat and affine.
 Hence it suffices to prove the second assertion of the lemma
in the case of a flat affine morphism~$f_2$, which is covered
by~\cite[Lemma Tag~02KG]{SP} or~\cite[Lemma~3.3(a)]{Psemten} (see also
Lemma~\ref{bijective-closed-immersion-affine-base-change-lemma}(a)
above).
\end{proof}

 Let $\ff\:\fY\rarrow\fX$ be a flat tight morphism of semi-separated
Noetherian formal schemes.
 Then the exact functor of direct image of flat contraadjusted
pro-quasi-coherent pro-sheaves
\begin{equation} \label{flat-cta-pro-sheaves-direct-image}
 \ff_*\:\fY\Proqc_\fl^\cta\lrarrow\fX\Proqc_\fl^\cta
\end{equation}
is defined by the rule $(\ff_*\rQ)^{(\overline X)}=
\bar f_*(\rQ^{(\overline Y)})$ for all $\rQ\in\fY\Proqc_\fl^\cta$,
where the notation is $\overline Y=\overline X\times_\fX\fY$ as
in the beginning of Section~\ref{special-inverse-images-subsecn}
and in the constructions of the direct image functors in
Section~\ref{pro-sheaves-ind-cosheaves-inverse-direct-images-subsecn}.
 Lemma~\ref{flat-contraadjusted-base-change-lemma} implies that
the flat contraadjusted pro-quasi-coherent pro-sheaf $\ff_*\rG$
on $\fX$ is well-defined by this formula.

 By Corollary~\ref{qcomp-ssep-direct-image-of-flat-cta-cor}(c),
the direct image functor $\ff_*$
\,\eqref{flat-cta-pro-sheaves-direct-image} takes flat cotorsion
pro-quasi-coherent pro-sheaves on $\fY$ to flat cotorsion
pro-quasi-coherent pro-sheaves on $\fX$, so we also have the direct
image functor
\begin{equation} \label{flat-cot-pro-sheaves-direct-image}
 \ff_*\:\fY\Proqc_\fl^\cot\lrarrow\fX\Proqc_\fl^\cot.
\end{equation}
 For a very flat tight morphism of Noetherian formal schemes
$\ff\:\fY\rarrow\fX$, the direct image functor $\ff_*$
\,\eqref{flat-cta-pro-sheaves-direct-image} also takes very flat
contraadjusted pro-quasi-coherent pro-sheaves on $\fY$ to very flat
contraadjusted pro-quasi-coherent pro-sheaves on $\fX$ by
Corollary~\ref{qcomp-ssep-direct-image-of-flat-cta-cor}(a),
so we have the direct image functor
\begin{equation} \label{vfl-cta-pro-sheaves-direct-image}
 \ff_*\:\fY\Proqc_\vfl^\cta\lrarrow\fX\Proqc_\vfl^\cta.
\end{equation}

 Let $\ff\:\fY\rarrow\fX$ be a flat tight morphism of semi-separated
Noetherian formal schemes.
 Then, for any flat pro-quasi-coherent pro-sheaf $\rF$ on $\fX$ and
any flat contraadjusted pro-quasi-coherent pro-sheaf $\rQ$ on $\fY$,
there is a natural adjunction isomorphism of abelian groups
\begin{equation} \label{cta-pro-sheaves-inverse-direct-adjunction}
 \Hom_{(\fY)}(\ff^*\rF,\rQ)\simeq\Hom_{(\fX)}(\rF,\ff_*\rQ).
\end{equation}
 This follows from the adjunctions of the functors of inverse image
and direct image of quasi-coherent sheaves with respect to
the morphisms of schemes $\bar f\:\overline Y\rarrow\overline X$
(cf.\ formula~\eqref{qcoh-tors-inverse-direct-adjunction}).

\subsection{The flat-very flat derived equivalence}
 The aim of this section is to prove the following theorem,
which is a formal scheme version of the theorem of Estrada and
Sl\'avik~\cite[Corollary~6.2]{ES} (see
also~\cite[Theorem~5.5.16]{Pcosh}).

\begin{thm} \label{estrada-slavik-theorem-for-formal-schemes}
 Let\/ $\fX$ be a semi-separated Noetherian formal scheme.
 Then, for any derived category symbol\/ $\st=-$ or\/~$\varnothing$,
the inclusion of exact categories\/ $\fX\Proqc_\vfl\rarrow
\fX\Proqc_\fl$ induces an equivalence of the conventional derived
categories
$$
 \sD^\st(\fX\Proqc_\vfl)\simeq\sD^\st(\fX\Proqc_\fl).
$$
\end{thm}

 The argument is similar in spirit to the proof of
Theorem~\ref{lcth-lcta-lct-unbounded-derived-equivalence} above.
 We start with the following lemma from adic commutative algebra.

\begin{lem} \label{contramodule-complex-very-flat-precover}
 Let $R$ be an adic topological ring and\/ $\fG^\bu$ be a complex
of flat quotseparated contramodule $R$\+modules.
 Then there exists a short exact sequence of complexes of
flat quotseparated contramodule $R$\+modules\/ $0\rarrow\fB^\bu
\rarrow\fP^\bu\rarrow\fG^\bu\rarrow0$ such that\/ $\fP^\bu$ is
a complex of very flat quotseparated contramodule $R$\+modules
and\/ $\fB^\bu$ is an acyclic complex in the exact category of flat
quotseparated contramodule $R$\+modules.
\end{lem}

\begin{proof}
 The formulation of the lemma can be strengthened: in fact, one can
choose $\fP^\bu$ to be a complex of projective quotseparated
contramodule $R$\+modules rather than merely very flat ones.
 In this form, the assertion actually holds for left contramodules
over any complete, separated topological associative ring with
a countable base of neighborhoods of zero consisting of open
two-sided ideals (which is a generalization of the setting of
the lemma by~\cite[Proposition~1.5]{Pdc}).
 The argument is based on the results of the papers~\cite{PS4,Pbc}.
 We use the language of cotorsion pairs from
Section~\ref{appx-cotorsion-pairs-subsecn} and the language of
Becker-contraacyclic complexes from Section~\ref{appx-becker-subsecn}.

 Denote, as usual, by $\fR\Contra$ the abelian category of left
contramodules over the topological ring~$\fR$.
 The abelian category $\fR\Contra$ is locally $\aleph_1$\+presentable
(in the sense of~\cite[Definition~1.17 and Theorem~1.20]{AR}) according
to~\cite[Section~5]{PR}.
 The free left $\fR$\+contramodule $\fR$ with one generator is
a projective generator of $\fR\Contra$.

 Hence the result of~\cite[Theorem~7.3]{PS4} or~\cite[Theorem~7.2]{Pbc}
is applicable to $\fR\Contra$.
 This theorem claims, in particular, that for any complex of left
$\fR$\+con\-tra\-mod\-ules $\fM^\bu$ there exists a short exact sequence
of complexes of left $\fR$\+con\-tra\-mod\-ules $0\rarrow\fB^\bu\rarrow
\fP^\bu\rarrow\fM^\bu\rarrow0$ with a complex of projective left
$\fR$\+con\-tra\-mod\-ules $\fP^\bu$ and a Becker-contraacyclic complex
of left $\fR$\+con\-tra\-mod\-ules~$\fB^\bu$.

 Take $\fM^\bu=\fG^\bu$ to be a complex of flat left
$\fR$\+contramodules (in the sense of~\cite[Sections~5\+-7]{PR}
or~\cite[Section~E.1]{Pcosh}).
 Then, since all projective left $\fR$\+contramodules are
flat~\cite[Lemma~6.9]{PR}, \cite[Lemma~E.1.6]{Pcosh} and
the class of flat left $\fR$\+contramodules is closed under kernels
of epimorphisms in $\fR\Contra$ \,\cite[Corollaries~6.8 and~6.15]{PR},
\cite[Lemma~E.1.4 and Corollary~E.1.7]{Pcosh}, the complex $\fB^\bu$
is also a complex of flat left $\fR$\+contramodules.
 Finally, by~\cite[Theorem~6.1 or~11.1]{Pbc}, any Becker-contraacyclic
complex of flat left $\fR$\+contramodules is acyclic in the exact
category of flat left $\fR$\+contramodules $\fR\Contra_\fl$.
\end{proof}

\begin{prop} \label{flat-pro-sheaf-complexes-vfl-precover}
 Let\/ $\fX$ be a semi-separated Noetherian formal scheme with a finite
affine open covering\/ $\fX=\bigcup_{\alpha=1}^N \fU_\alpha$.
 Then any complex of flat pro-quasi-coherent pro-sheaves\/ $\rG^\bu$
on\/ $\fX$ can be included in a (termwise admissible) short exact
sequence\/ $0\rarrow\rB^\bu\rarrow\rF^\bu\rarrow\rG^\bu\rarrow0$ of
complexes in the exact category\/ $\fX\Proqc_\fl$, where\/ $\rF^\bu$
is a complex of very flat pro-quasi-coherent pro-sheaves on\/ $\fX$,
while\/ $\rB^\bu$ is a finitely iterated extension of the direct
images of acyclic complexes (in the exact categories) of flat
pro-quasi-coherent pro-sheaves from\/~$\fU_\alpha$.
 In particalar, $\rB^\bu$ is an acyclic complex in the exact category\/
$\fX\Proqc_\fl$.
\end{prop}

\begin{proof}
 This is a weak formal scheme version of~\cite[Lemma~5.1.9(a)]{Pcosh}
or~\cite[Proposition~2.2 and Corollary~4.4]{Pfltp}.
 The proof of the first assertion of the proposition proceeds
along the lines of the proof of
Lemma~\ref{flat-pro-sheaves-very-flat-precover} above.
 Following the idea of~\cite[proof of Lemma~6.3]{ES}, we plug in
whole complexes of flat pro-quasi-coherent pro-sheaves, rather than
single flat pro-quasi-coherent pro-sheaves, into the the construction
from the proof of Lemma~\ref{flat-pro-sheaves-very-flat-precover}.
 Otherwise, the construction proceeds exactly in the same way, with
the only difference that, in the context of an affine Noetherian
formal scheme $\fU=\fU_\beta$, instead of the short exact sequence
of contramodule $\fO(\fU)$\+modules from
Proposition~\ref{very-flat-cotorsion-pair-in-quotseparated}(a),
we now use the short exact sequence of complexes of (flat)
contramodule $\fO(\fU)$\+modules from
Lemma~\ref{contramodule-complex-very-flat-precover}.

 To deduce the second assertion of the proposition, it suffices to
observe that, for an affine open subscheme $\fU\subset\fX$ with
the open immersion morphism $\fj\:\fU\rarrow\fX$, we have an exact
functor of direct image $\fj_*\:\fU\Proqc_\fl\rarrow\fX\Proqc_\fl$
\,\eqref{flat-pro-sheaves-direct-image}.
 Furthermore, any finitely iterated extension of acyclic complexes in
an exact category is an acyclic complex.
\end{proof}

\begin{proof}[Proof of
Theorem~\ref{estrada-slavik-theorem-for-formal-schemes}]
 The easy case of bounded above complexes, $\st=-$, is covered by
Proposition~\ref{infinite-co-resolutions-minus-plus-derived}(a).
 Indeed, the full subcategory $\fX\Proqc_\vfl$ is resolving in
$\fX\Proqc_\fl$, essentially by
Lemma~\ref{flat-pro-sheaves-very-flat-precover}
or Corollary~\ref{cta-flat-pro-sheaves-cor}(a); so
Proposition~\ref{infinite-co-resolutions-minus-plus-derived}(a)
is applicable.

 In the interesting case of two-sided unbounded complexes,
$\st=\varnothing$, the argument is based on
Proposition~\ref{flat-pro-sheaf-complexes-vfl-precover} together
with the flat and projective periodicity theorem for modules over
associative rings (due to Benson--Goodearl~\cite[Theorem~2.5]{BG}
and Neeman~\cite[Remark~2.15]{Neem}).
 More specifically, we use \v St\!'ov\'\i\v cek's flat and very flat
periodicity theorem for modules over commutative rings or quasi-coherent
sheaves on schemes~\cite[Proposition~4.4]{ES},
\cite[Lemma~5.5.15]{Pcosh}, which implies a similar flat and very flat
periodicity assertion for flat pro-quasi-coherent pro-sheaves:
 Every complex in $\fX\Proqc_\vfl$ that is acyclic in $\fX\Proqc_\fl$
is also acyclic in $\fX\Proqc_\vfl$.

 In view of a basic lemma about full triangulated subcategories and
triangulated Verdier quotient categories~\cite[Corollary~10.2.7]{KS},
\cite[Lemma~A.3.3(b)]{Pcosh}, \cite[Lemma~8.1(b)]{Pfltp}, in order
to prove the assertion of the theorem for $\st=\varnothing$ it suffices
to establish two properties:
\begin{enumerate}
\item for any complex $\rG^\bu$ in $\fX\Proqc_\fl$, there exists
a complex $\rF^\bu$ in $\fX\Proqc_\vfl$ together with a morphism of
complexes $\rF^\bu\rarrow\rG^\bu$ whose cone is an acyclic complex
in $\fX\Proqc_\fl$;
\item any complex in $\fX\Proqc_\vfl$ that is acyclic in
$\fX\Proqc_\fl$ is also acyclic in $\fX\Proqc_\fl$.
\end{enumerate}
 Now item~(1) is provided by
Proposition~\ref{flat-pro-sheaf-complexes-vfl-precover}, and
item~(2) was explained in the previous paragraph.
\end{proof}

\subsection{Antilocal locally injective ind-contraherent ind-cosheaves}
 Let $\fX$ be a semi-separated Noetherian formal scheme.
 We say that a locally injective ind-contraherent ind-cosheaf $\gP$
on $\fX$ is \emph{antilocal} if the locally injective contraherent
cosheaf $\gP_{(Z)}$ on $Z$ is antilocal for every closed subscheme
of definition $Z\subset\fX$.

 Given a filtered inductive system $(X_\gamma)_{\gamma\in\Gamma}$
of closed subschemes of definition in $\fX$ such that
$\fX=\ilim_{\gamma\in\Gamma}X_\gamma$, it suffices to check the latter
antilocality condition for the closed subschemes $Z=X_\gamma\subset\fX$,
\,$\gamma\in\Gamma$.
 This follows from Corollary~\ref{inverse-image-of-antilocal-cor}(a).

 We denote the full subcategory of antilocal locally injective
ind-contraherent ind-cosheaves by $\fX\Indct^\lin_\al\subset
\fX\Indct^\lin\subset\fX\Indlc^\lin$.

 Over an affine scheme $U$, all contraherent cosheaves are antilocal
(as mentioned in the beginning of Section~\ref{antilocal-subsecn}).
 Consequently, over an affine Noetherian formal scheme $\fU$, all
locally injective ind-contraherent ind-cosheaves are antilocal, so one
has $\fU\Indct^\lin_\al=\fU\Indct^\lin$.
 Of course, this only applies to locally injective ind-contraherent
ind-cosheaves and \emph{not} to locally injective \emph{locally}
ind-contraherent ind-cosheaves on~$\fU$.

\begin{lem} \label{cohom-into-antilocal-lin-is-antilocal}
 Let $X$ be a quasi-compact semi-separated scheme.
 Then, for any quasi-coherent sheaf\/ $\cM$ and any antilocal locally
injective contraherent cosheaf\/ $\fP$ on $X$, the locally cotorsion
contraherent cosheaf\/ $\Cohom_X(\cM,\fP)$ on $X$ is antilocal.
\end{lem}

\begin{proof}
 This lemma is a dual analogue of
Corollary~\ref{tensor-product-with-coh-sheaf-cta-cot}.
 By Proposition~\ref{antilocal-lin-prop}, the cosheaf $\fP$ is
a direct summand of a finite direct sum of the direct images of
locally injective contraherent cosheaves from affine open subschemes
of~$X$.
 So, let $U\subset X$ be an affine open subscheme with the open
immersion morphism $j\:U\rarrow X$.
 It suffices to consider the case of a cosheaf $\fP$ of the form
$\fP=j_!\fK$, where $\fK$ is a locally injective contraherent
cosheaf on~$\fU$.
 In this case, we have $\Cohom_X(\cM,j_!\fK)\simeq
j_!\Cohom_U(j^*\cM,\fK)$ by the projection
formula~\eqref{Cohom-into-flat-direct-image-of-lin}.
 Now $\Cohom_U(j^*\cM,\fK)$ is a locally cotorsion contraherent
cosheaf on $U$ by Section~\ref{Cohom-subsecn}, and the cosheaf
$j_!\Cohom_U(j^*\cM,\fK)$ on $X$ is locally cotorsion contraherent
by formula~\eqref{ctrh-lct-affine-direct-image} and antilocal by
the discussion in the beginning of Section~\ref{antilocal-subsecn}
or by Corollary~\ref{direct-image-of-antilocal-cor}(a) or~(b).
\end{proof}

\begin{lem} \label{closed-subscheme-adjunction-morphism-in-lcth-lemma}
 Let $i\:Z\rarrow X$ be a closed immersion of schemes and\/
$\cI\subset\cO_X$ be the corresponding quasi-coherent sheaf
of ideals on~$X$.
 Let\/ $\bW$ be an open covering of $X$ and\/ $\fJ$ be a locally
injective\/ $\bW$\+locally contraherent cosheaf on~$X$.
 Then the adjunction morphism $i_!i^!\fJ\rarrow\fJ$ is an admissible
monomorphisms in the exact category $X\Lcth^\lct_\bW$ of locally
cotorsion\/ $\bW$\+locally contraherent cosheaves on~$X$.
 Moreover, there is a natural (admissible) short exact sequence\/
$0\rarrow i_!i^!\fJ\rarrow\fJ\rarrow\Cohom_X(\cI,\fJ)\rarrow0$
in the exact category $X\Lcth^\lct_\bW$.
\end{lem}

\begin{proof}
 This is dual-analogous to
Lemma~\ref{closed-subscheme-adjunction-morphism-in-qcoh-lemma}.
 Let $\bT$ be an open covering of $Z$ such that the morphism
$i\:Z\rarrow X$ is $(\bW,\bT)$\+affine and $(\bW,\bT)$ coaffine
(for example, one can take $\bT$ to be the collection of all
the preimages $i^{-1}(U)\subset Z$ of affine open subschemes
$U\subset X$ subordinate to~$\bW$).
 By formula~\eqref{inverse-image-lcth-lin-W-T}, we have
a locally injective $\bT$\+locally contraherent cosheaf $i^!\fJ$ on~$Z$.
 By formula~\eqref{ctrh-lct-affine-direct-image}, \,$i_!i^!\fJ$ is
a locally cotorsion $\bW$\+locally contraherent cosheaf on~$X$.
 By formula~\eqref{cosh-on-schemes-lin-direct-inverse-adjunction},
we have a natural adjunction morphism $i_!i^!\fJ\rarrow\fJ$
in $X\Lcth^\lct_\bW$.

 According to Section~\ref{Cohom-subsecn}, we have a locally
cotorsion $\bW$\+locally contraherent cosheaf $\Cohom_X(\cI,\fJ)$
on~$X$.
 The natural inclusion morphism of quasi-coherent sheaves
$\cI\rarrow\cO_X$ induces a morphism $\fJ\rarrow\Cohom_X(\cI,\fJ)$
in $X\Lcth^\lct_\bW$.

 Finally, the assertion that the short sequence in
$0\rarrow i_!i^!\fJ\rarrow\fJ\rarrow\Cohom_X(\cI,\fJ)\rarrow0$ is exact
is local in $X$, so it suffices to check it for affine schemes,
for which it means the following.
 Let $R\rarrow S$ be a surjective homomorphism of commutative rings
with the kernel ideal $I\subset R$, and let $J$ be an injective
$R$\+module.
 Then there is a natural short exact sequence of (cotorsion)
$R$\+modules $0\rarrow\Hom_R(S,J)\rarrow J\rarrow\Hom_R(I,J)\rarrow0$.
 This is obvious.
\end{proof}

 Let $\fX$ be a Noetherian formal scheme with an open covering $\bW$
and an open covering $\bW'$ subordinate to~$\bW$.
 Then the full subcategory $\fX\Indlc^\lin_\bW$ is closed under
extensions in $\fX\Indlc^\lin_{\bW'}$ (since, for any closed subscheme
of definition $Z\subset\fX$, the full subcategory $Z\Lcth^\lin_\bW$
is closed under extensions in the exact category $Z\Lcth^\lin_{\bW'}$,
as per~\cite[Section~2.2 and Corollary~3.2.4]{Pcosh} or
Sections~\ref{locally-contraherent-cosheaves-of-contramods-subsecn}
and~\ref{loc-injective-on-schemes-subsecn} above).
 Accordingly, the groups $\Ext^1$ computed in the exact categories
$\fX\Indlc^\lin_\bW$ agree as the open covering $\bW$ varies.
 We denote the resulting $\Ext^1$ groups by $\Ext^{(\fX),1}({-},{-})$.

\begin{lem} \label{ind-lcth-ind-cosheaves-Ext1-vanishing-lemma}
 Let\/ $\fX$ be a semi-separated Noetherian formal scheme with
an open covering\/~$\bW$.
 Then, for any antilocal locally injective ind-contraherent
ind-cosheaf\/ $\gP$ and any locally injective\/ $\bW$\+locally
ind-contraherent ind-cosheaf\/ $\gJ$ on\/ $\fX$, one has\/
$\Ext^{(\fX),1}(\gP,\gJ)=0$.
 So every antilocal locally injective ind-contraherent ind-cosheaf\/
$\gP$ on\/ $\fX$ is a projective object of the exact category\/
$\fX\Indlc^\lin_\bW$.
\end{lem}

\begin{proof}
 This is dual-analogous to
Lemma~\ref{pro-qc-pro-sheaves-Ext1-vanishing-lemma}.
 Let $(X_n)_{n=1}^\infty$ be an inductive system of closed subschemes
of definition in $\fX$, indexed by the ordered set of positive integers,
such that $\fX=\ilim_{n\ge1}X_n$.
 Denote by $i_{n,n+1}\:X_n\rarrow X_{n+1}$ the transition morphisms.
 For every locally injective locally ind-contraherent ind-cosheaf $\gK$
on $\fX$, put $\gK_{(n)}=\gK_{(X_n)}$ for brevity.
 So $\gK_{(n)}\in X_n\Lcth^\lin$.

 Let $0\rarrow\gJ\rarrow\gK\rarrow\gP\rarrow0$ be a short exact
sequence in $\fX\Indlc^\lin_\bW$.
 For every $n\ge1$, we have a split short exact sequence $0\rarrow
\gJ_{(n)}\rarrow\gK_{(n)}\rarrow\gP_{(n)}\rarrow0$ in
$X_n\Lcth^\lin_\bW$ (by Proposition~\ref{antilocal-prop}(b)
and Corollary~\ref{antilocal-characterized}).
 We need to construct a sequence of splittings $s_n\:\gP_{(n)}\rarrow
\gK_{(n)}$ that would be compatible as $n$~varies, in the sense that
$i_{n,n+1}^!(s_{n+1})=s_n$ for all $n\ge1$.
 Proceeding by induction, assumt that a compatible system of
splittings $s_1$,~\dots, $s_n$ has been constructed already.
 Let us construct a splitting $s_{n+1}$ compatible with~$s_n$.

 Pick an arbitrary splitting $s_{n+1}'\:\gP_{(n+1)}\rarrow\gK_{(n+1)}$
of the split short exact sequence $0\rarrow\gJ_{(n+1)}\rarrow\gK_{(n+1)}
\rarrow\gP_{(n+1)}\rarrow0$ in $X_{n+1}\Lcth^\lin_\bW$.
 Then the difference $s_n-i_{n,n+1}^!(s'_{n+1})\:\gP_{(n)}\rarrow
\gK_{(n)}$ factorizes through the admissible monomorphism
$\gJ_{(n)}\rarrow\gK_{(n)}$; so we obtain a morphism of locally
contraherent cosheaves $t_n\:\gP_{(n)}\rarrow\gJ_{(n)}$.
 We claim that any such morphism can be lifted to a morphism
$t_{n+1}\:\gP_{(n+1)}\rarrow\gJ_{(n+1)}$, i.~e., one has
$t_n=i_{n,n+1}^!(t_{n+1})$.
 If this is the case, then it suffices to put $s_{n+1}=s'_{n+1}+t_{n+1}$
(where the composition $\gP_{(n+1)}\xrightarrow{t_{n+1}}\gJ_{(n+1)}
\rarrow\gK_{(n+1)}$ is still denoted by $t_{n+1}$ by abuse of notation).

 By adjunction~\eqref{cosh-on-schemes-lin-direct-inverse-adjunction},
the given morphism $t_n\:i_{n,n+1}^!\gP_{(n+1)}\rarrow
i_{n,n+1}^!\gJ_{(n+1)}$ in $X_n\Lcth^\lin_\bW$ corresponds to
a morphism $\tilde t_n\:i_{n,n+1}{}_!i_{n,n+1}^!\gP_{(n+1)}\rarrow
\gJ_{(n+1)}$ in $X_n\Lcth^\lct_\bW$.
 We need to show that the morphism~$\tilde t_n$ factorizes through
the natural adjunction (mono)morphism $i_{n,n+1}{}_!i_{n,n+1}^!
\gP_{(n+1)}\rarrow\gP_{n+1}$ in $X_{n+1}\Lcth^\lct_\bW$.

 By Lemma~\ref{closed-subscheme-adjunction-morphism-in-lcth-lemma},
we have a short exact sequence $0\rarrow i_{n,n+1}{}_!i_{n,n+1}^!
\gP_{(n+1)}\rarrow\gP_{(n+1)}\rarrow\Cohom_{X_{n+1}}(\cI,\gP_{(n+1)})
\rarrow0$ in the exact category $X_{n+1}\Lcth^\lct_\bW$.
 Here $\cI\subset\cO_{X_{n+1}}$ denotes the quasi-coherent sheaf of
ideals on $X_{n+1}$ corresponding to the closed subscheme
$i_{n,n+1}(X_n)\subset X_{n+1}$.
 Finally, the locally cotorsion contraherent cosheaf
$\Cohom_{X_{n+1}}(\cI,\gP_{(n+1)})$ on $X_{n+1}$ is antilocal by
Lemma~\ref{cohom-into-antilocal-lin-is-antilocal}, while
the $\bW$\+locally contraherent cosheaf $\gJ_{(n+1)}$ on $X_{n+1}$
is locally injective by the definition of the exact category
$\fX\Indlc^\lin_\bW$.
 Thus, by Corollary~\ref{antilocal-characterized}(2), the morphism
of locally contraherent cosheaves $\tilde t_n\:i_{n,n+1}{}_!i_{n,n+1}^!
\gP_{(n+1)}\rarrow\gJ_{(n+1)}$ can be extended to a morphism of
locally contraherent cosheaves $t_{n+1}\:\gP_{(n+1)}\rarrow\gJ_{(n+1)}$
on $X_{n+1}$, as desired.
\end{proof}

 Let $\ff\:\fY\rarrow\fX$ be a flat tight affine morphism of
semi-separated Noetherian formal schemes.
 Then it follows from the construction together with the discussion of
the direct images of antilocal contraherent cosheaves in the beginning
of Section~\ref{antilocal-subsecn} or in
Corollary~\ref{direct-image-of-antilocal-cor}(a,c) that
the direct image functor $\ff_!\:\fY\Indct^\lin\rarrow\fX\Indct^\lin$
\,\eqref{lin-ind-ctrh-cosheaves-direct-image} takes antilocal locally
injective ind-contraherent ind-cosheaves on $\fY$ to antilocal locally
injective ind-contraherent ind-cosheaves on~$\fX$.

\begin{cor} \label{enough-projectives-in-lin-ind-cosheaves}
 Let\/ $\fX$ be a semi-separated Noetherian formal scheme with an open
covering\/ $\bW$ and a finite affine open covering\/
$\fX=\bigcup_{\alpha=1}^N\fU_\alpha$ subordinate to\/~$\bW$.
 Then there are enough projective objects in the exact category\/
$\fX\Indlc^\lin_\bW$.
 All the projective objects of\/ $\fX\Indlc^\lin_\bW$ belong to\/
$\fX\Indct^\lin$, and the full subcategory of projective objects in\/
$\fX\Indlc^\lin_\bW$ does not depend on an open covering\/~$\bW$.
 A locally injective\/ $\bW$\+locally ind-contraherent ind-cosheaf is
a projective object in\/ $\fX\Indlc^\lin_\bW$ if and only if it is
an antilocal locally injective ind-contraherent ind-cosheaf.
 A locally injective ind-contraherent ind-cosheaf on\/ $\fX$ is
antilocal if and only if is a direct summand of a finite direct sum
of the direct images of locally injective ind-contraherent ind-cosheaves
from\/~$\fU_\alpha$.
\end{cor}

\begin{proof}
 First of all, antilocal locally injective ind-contraherent
ind-cosheaves are projective objects in $\fX\Indlc^\lin_\bW$ by
Lemma~\ref{ind-lcth-ind-cosheaves-Ext1-vanishing-lemma}, and
the direct images of locally injective ind-contraherent ind-cosheaves
from $\fU_\alpha$ are antilocal locally injective ind-contraherent
ind-cosheaves on $\fX$ according to the paragraph preceding the present
corollary.
 Thus it remains to prove that every object of $\fX\Indlc^\lin_\bW$
is an admissible quotient object of an object of the form
$\bigoplus_{\alpha=1}^N\fj_\alpha{}_!\gQ_\alpha$, where
$j_\alpha\:\fU_\alpha\rarrow\fX$ are the open immersion morphisms and
$\gQ_\alpha$ are locally injective ind-contraherent ind-cosheaves
on~$\fU_\alpha$.

 Specifically, let $\gJ$ be a locally injective $\bW$\+locally
ind-contraherent ind-cosheaf on~$\fX$.
 Then the adjunction of the direct image functors~$\fj_\alpha{}_!$ with
the inverse image functors~$\fj_\alpha^!$ (see
Section~\ref{pro-sheaves-ind-cosheaves-inverse-direct-images-subsecn})
provides a natural morphism $g\:\bigoplus_{\alpha=1}^N\fj_\alpha{}_!
\fj_\alpha^!\gJ\rarrow\gJ$.
 The morphism~$g$ is an admissible epimorphism of locally injective
ind-contraherent ind-cosheaves in restriction to every~$\fU_\alpha$.
 In view the discussion of the locality of admissible epimorphisms of
(locally injective locally) contraherent cosheaves
in~\cite[Lemma~1.4.3(b), Remark~2.2.5, and Section~3.2]{Pcosh}
or Sections~\ref{locally-contraherent-cosheaves-of-contramods-subsecn}
and~\ref{loc-injective-on-schemes-subsecn} above), it follows
that $g$~is an admissible epimorphism in $\fX\Indlc^\lin_\bW$.
\end{proof}

 Let $\fX$ be a semi-separated Noetherian formal scheme with an open
covering $\bW$ and an open covering $\bW'$ subordinate to~$\bW$.
 By Corollary~\ref{enough-projectives-in-lin-ind-cosheaves} together
with the discussion of the locality of admissible epimorphisms of
(locally injective locally) contraherent cosheaves mentioned in
its proof, the full subcategory $\fX\Indlc^\lin_\bW$ is resolving
in $\fX\Indlc^\lin_{\bW'}$.
 By virtue of
Proposition~\ref{infinite-co-resolutions-minus-plus-derived}(a),
it follows that the $\Ext$ groups in all the exact categories
$\fX\Indlc^\lin_\bW$ agree with each other as the open covering
$\bW$ varies.
 We denote these $\Ext$ groups by $\Ext^{(\fX),*}({-},{-})$.

\begin{lem} \label{flat-affine-antiloc-base-change-lemma}
 Let $f_2\:Y_2\rarrow X_2$ be a flat morphism of quasi-compact
semi-separated schemes and $i_{12}\:X_1\rarrow X_2$ be an affine
morphism of schemes.
 Put $Y_1=X_1\times_{X_2}Y_2$; so we have a pullback diagram of schemes
$$
 \xymatrix{
  Y_1 \ar[r]^{k_{12}} \ar[d]_{f_1} & Y_2 \ar[d]^{f_2} \\
  X_1 \ar[r]_{i_{12}} & X_2
 }
$$
with a flat morphism $f_1\:Y_2\rarrow X_2$ and a bijective closed
immersion $k_{12}\:Y_1\rarrow Y_2$.
 Then, for any antilocal locally injective contraherent cosheaf\/ $\fP$
on $Y_2$, there is a natural base change isomorphism
\begin{equation} \label{flat-affine-antiloc-basechange-isom}
 f_1{}_!k_{12}^!\fP\simeq i_{12}^!f_2{}_!\fP
\end{equation}
of antilocal locally injective contraherent cosheaves on~$X_1$.
 In other words, there is a natural base change isomorphism
$f_1{}_!k_{12}^!\simeq i_{12}^!f_2{}_!$ of functors
$Y_2\Ctrh^\lin_\al\rarrow X_1\Ctrh^\lin_\al$.
\end{lem}

\begin{proof}
 This is dual-analogous to
Lemma~\ref{flat-contraadjusted-base-change-lemma}.
 Notice that we have a direct image functor $f_2{}_!\:
Y_2\Ctrh^\lin_\al\rarrow X_2\Ctrh^\lin_\al$ by
Corollary~\ref{direct-image-of-antilocal-cor}(c) and an inverse image
functor $i_{12}^!\:X_2\Ctrh^\lin_\al\rarrow X_1\Ctrh^\lin_\al$ by
Corollary~\ref{inverse-image-of-antilocal-cor}(a).
 Similarly, we have a direct image functor $f_1{}_!\:
Y_1\Ctrh^\lin_\al\rarrow X_1\Ctrh^\lin_\al$ and an inverse image functor
$k_{12}^!\:Y_2\Ctrh^\lin_\al\rarrow Y_1\Ctrh^\lin_\al$.
 So both the left-hand side and the right-hand side
of~\eqref{flat-affine-antiloc-basechange-isom} are well-defined.

 The natural morphism from the left-hand side to the right-hand side
of~\eqref{flat-affine-antiloc-basechange-isom} corresponds by
adjunction~\eqref{cosh-on-schemes-lin-direct-inverse-adjunction}
or~\eqref{cosh-on-schemes-lct-lcta-direct-inverse-adjn}
to a morphism $k_{12}^!\fP\rarrow k_{12}^!f_2^!f_2{}_!\fP
\simeq f_1^!i_{12}^!f_2{}_!\fP$ in $Y_1\Lcth^\lin$,
which is obtainable by applying the inverse image functor
$k_{12}^!\:Y_2\Lcth^\lin\rarrow Y_1\Lcth^\lin$
\,\eqref{inverse-image-lcth-lin-nocovering} to the adjunction
morphism $\fP\rarrow f_2^!f_2{}_!\fP$ in $Y_2\Lcth^\lin$.
 The same natural morphism from the left-hand side to the right-hand
side of~\eqref{flat-affine-antiloc-basechange-isom} corresponds by
adjunction~\eqref{cosh-on-schemes-lin-direct-inverse-adjunction}
to a morphism $i_{12}{}_!f_1{}_!k_{12}^!\fP\simeq
f_2{}_!k_{12}{}_!k_{12}^!\fP\rarrow f_2{}_!\fP$ in
$X_2\Ctrh^\lct$, which is obtainable by applying the direct image
functor $f_2{}_!\:Y_2\Ctrh^\lct_\al\rarrow Y_1\Ctrh^\lct_\al$
from Corollary~\ref{direct-image-of-antilocal-cor}(b) to
the adjunction morphism $k_{12}{}_!k_{12}^!\fP\rarrow\fP$
in $Y_2\Ctrh^\lct_\al$.

 It remains to prove that the morphism in question is an isomorphism.
 By Proposition~\ref{antilocal-lin-prop}, \,$\fP$ is a direct summand
of a finite direct sum of the direct images of locally injective
contraherent cosheaves from affine open subschemes $U\subset Y_2$.
 Hence the question reduces to the case of $\fP=j_!\fQ$, where $\fQ$ is
a locally injective contraherent cosheaf on $V$ and $j\:V\rarrow Y_2$
is the open immersion morphism.

 Now both the morphisms $j\:V\rarrow Y_2$ and $f_2j\:V\rarrow X_2$
are flat and affine.
 Hence it suffices to prove the desired isomorphism in the case of
a flat affine morphism~$f_2$, which is covered by the discussion
in~\cite[Section~3.3]{Pcosh}.
\end{proof}

 Let $\ff\:\fY\rarrow\fX$ be a flat tight morphism of semi-separated
Noetherian formal schemes.
 Then the functor of direct image of antilocal locally injective
ind-contraherent ind-cosheaves
\begin{equation} \label{antilocal-lin-ind-cosheaves-direct-image}
 \ff_!\:\fY\Indct^\lin_\al\lrarrow\fX\Indct^\lin_\al
\end{equation}
is defined by the rule $(\ff_!\gQ)_{(\overline X)}=
\bar f_!(\gQ_{(\overline Y)})$ for all $\gQ\in\fY\Indct^\lin_\al$,
where the notation is $\overline Y=\overline X\times_\fX\fY$ as
in the beginning of Section~\ref{special-inverse-images-subsecn}
and in the constructions of the direct image functors in
Section~\ref{pro-sheaves-ind-cosheaves-inverse-direct-images-subsecn}.
 The notation $\bar f_!\:\overline Y\Ctrh^\lin_\al\rarrow
\overline X\Ctrh^\lin_\al$ stands for the direct image functor
from Corollary~\ref{direct-image-of-antilocal-cor}(c).
 Lemma~\ref{flat-affine-antiloc-base-change-lemma} implies that
the antilocal locally injective ind-contraherent ind-cosheaf $\ff_!\gQ$
on $\fX$ is well-defined by the rule above.

 For any locally injective locally ind-contraherent ind-cosheaf $\gJ$
on $\fX$ and any antilocal locally injective ind-contraherent
ind-cosheaf $\gQ$ on $\fY$, there is a natural adjunction isomorphism
of abelian groups
\begin{equation} \label{al-ind-cosheaves-direct-inverse-adjunction}
 \Hom^{(\fY)}(\gQ,\ff^!\gJ)\simeq\Hom^{(\fX)}(\ff_!\gQ,\gJ).
\end{equation}
 This follows from the adjunctions of the functors of direct image of
antilocal locally injective contraherent cosheaves and locally
injective locally contraherent cosheaves with respect to the morphisms
of schemes $\bar f\:\overline Y\rarrow\overline X$
(see formula~\eqref{cosh-on-schemes-lin-direct-inverse-adjunction}).

\appendix

\Section{Category-Theoretic Background}

 This appendix contains almost no proofs, but only definitions,
formulations, and references.
 Our exposition is mostly based on~\cite[Appendices~A and~B]{Pcosh}.
 Other relevant reference sources are~\cite{Kel}, \cite{Neem0},
\cite{Psemi}, \cite{Bec}, \cite{SaoSt}, \cite{Sto-ICRA}, \cite{Sto},
\cite{PS4}, \cite{Pksurv}, \cite{BHP}, \cite{PS6}, \cite{Pal},
and~\cite{Pedg}.

\subsection{Resolving and coresolving subcategories}
\label{appx-resolving-subsecn}
 For background material on exact categories (in the sense of
Quillen), we suggest the survey paper~\cite{Bueh}.
 In the present paper, we are mostly interested in
\emph{idempotent-complete} exact categories~\cite[Section~6]{Bueh},
though this assumption is rarely needed, and the assumption of
\emph{weak idempotent-completeness}~\cite[Section~7]{Bueh} is often
more relevant.

 Given an exact category $\sK$ and a full subcategory $\sE\subset\sK$,
we always presume by default that $\sE$ is endowed with
the \emph{inherited exact structure} in which the (admissible) short
exact sequences in $\sE$ are the short exact sequences in $\sK$ with
the terms belonging to~$\sE$.
 More precisely, the full subcategory $\sE$ is said to \emph{inherit
an exact structure} from $\sK$ if the class of short exact sequences
defined in the previous sentence forms an exact structure on~$\sE$.
 A characterization of full subcategories in exact categories
inheriting an exact structure can be found in~\cite[Theorem~2.6]{DS}
or~\cite[Lemma~4.21]{Pedg}.
 Any full subcategory closed under extensions inherits an exact
category structure.

 Let $\sE$ be an exact category and $\sF\subset\sE$ be a full
subcategory.
 One says that the full subcategory $\sF$ is
\emph{self-generating}~\cite[Section~1]{BHP}, \cite[Section~7]{PS6},
\cite[Section~A.2]{Pcosh} in $\sE$ if for every admissible epimorphism
$E\rarrow F$ in $\sE$ with $F\in\sF$ there exists a morphism
$G\rarrow E$ in $\sE$ with $G\in\sF$ such that the composition
$G\rarrow F$ is an admissible epimorphism in~$\sE$.
 The full subcategory $\sF\subset\sE$ is said to be \emph{generating}
if for every object $E\in\sE$ there exists an object $G\in\sF$
together with an admissible epimorphism $G\rarrow E$ in~$\sE$.
 Clearly, any generating full subcategory in $\sE$ is self-generating.

 A self-generating full subcategory $\sF\subset\sE$ is said to be
\emph{self-resolving}~\cite[Section~7.1]{Pedg}, \cite[Section~9]{PS6},
\cite[Section~A.2]{Pcosh} if $\sF$ is closed under extensions and
kernels of admissible epimorphisms in~$\sE$.
 Similarly, a generating full subcategory $\sF\subset\sE$ that is
closed under extensions and kernels of admissible epimorphisms in $\sE$
is said to be \emph{resolving}~\cite[Section~2]{Sto},
\cite[Section~A.3]{Pcosh}.
 Any resolving full subcategory is self-resolving.

 Dually, let $\sC\subset\sE$ be another full subcategory.
 One says that the full subcategory $\sC$ is \emph{self-cogenerating}
in $\sE$ if for every admissible monomorphism $D\rarrow E$ in $\sE$ with
$D\in\sC$ there exist a morphism $E\rarrow C$ in $\sE$ with $C\in\sC$
such that the composition $D\rarrow C$ is an admissible monomorphism
in~$\sE$.
 The full subcategory $\sC\subset\sE$ is said to be \emph{cogenerating}
if for every object $E\in\sE$ there exists an object $C\in\sC$
together with an admissible monomorphism $E\rarrow C$ in~$\sE$.
 Clearly, any cogenerating full subcategory in $\sE$ is
self-cogenerating.

 A self-cogenerating full subcategory $\sC\subset\sE$ is said to be
\emph{self-coresolving} if $\sF$ is closed under extensions and
cokernels of admissible monomorphisms in~$\sE$.
 Similarly, a cogenerating full subcategory $\sF\subset\sE$ that is
closed under extensions and cokernels of admissible monomorphisms
in $\sE$ is said to be \emph{coresolving}.
 Any coresolving full subcategory is self-coresolving.

\begin{lem} \label{intersection-co-resolving-lemma}
 Let\/ $\sE$ be an exact category and\/ $\sE'\subset\sE$ be a full
subcategory inheriting an exact category structure. \par
\textup{(a)} Let\/ $\sF\subset\sE$ be a resolving subcategory in\/
$\sE$ such that the intersection\/ $\sF'=\sE'\cap\sF$ is a generating
subcategory in\/~$\sE'$.
 Then\/ $\sF'$ is a resolving subcategory in\/~$\sE'$. \par
\textup{(b)} Let\/ $\sC\subset\sE$ be a coresolving subcategory in\/
$\sE$ such that the intersection\/ $\sC'=\sE'\cap\sC$ is a cogenerating
subcategory in\/~$\sE'$.
 Then\/ $\sC'$ is a coresolving subcategory in\/~$\sE'$.
\end{lem}

\begin{proof}
 Obvious.
\end{proof}

 Notice that any full subcategory closed under kernels of epimorphisms
or cokernels of monomorphisms in a weakly idempotent-complete exact
category is weakly idempotent-complete.

 Let $\sF$ be a resolving subcategory in a weakly idempotent-complete
exact category~$\sE$.
 Let $d\ge-1$ be an integer.
 An object $E\in\sE$ is said to have \emph{$\sF$\+resolution
dimension\/~$\le\nobreak d$}~\cite[Section~2]{Sto} if there exists
a finite $(d\nobreak+2)$\+term exact sequence $0\rarrow F_d\rarrow
\dotsb\rarrow F_0\rarrow E\rarrow0$ in $\sE$ with $F_i\in\sF$ for
all $0\le i\le d$.
 We denote the $\sF$\+resolution dimension of an object $E\in\sE$
by $\rsd_{\sF/\sE} E\in\boZ_{\ge-1}\cup\{\infty\}$.

 Dually, let $\sC\subset\sF$ be a coresolving subcategory.
 An object $E\in\sE$ is said to have \emph{$\sC$\+coresolution
dimension\/~$\le\nobreak d$} if there exists a finite
$(d\nobreak+2)$\+term exact sequence $0\rarrow E\rarrow C^0\rarrow
\dotsb\rarrow C^d\rarrow0$ in $\sE$ with $C^i\in\sC$ for all
$0\le i\le d$.
 We denote the $\sC$\+coresolution dimension of an object $E\in\sE$
by $\crd_{\sC/\sE} E\in\boZ_{\ge-1}\cup\{\infty\}$.

 The following lemma expresses the assertion that the (co)resolution
dimension does not depend on the choice of a (co)resolution.

\begin{lem} \label{co-resol-dim-independent-of-co-resolution-lem}
 Let\/ $\sE$ be a weakly idempotent-complete exact category. \par
\textup{(a)} Let\/ $\sF\subset\sE$ be a resolving subcategory and
$E\in\sE$ be an object of\/ $\sF$\+resolution dimension\/~$\le d$.
 Let\/ $0\rarrow F\rarrow F_{d-1}\rarrow\dotsb\rarrow F_0\rarrow E
\rarrow0$ be a $(d\nobreak+2)$\+term exact sequence in\/ $\sE$
with $F_i\in\sF$ for all\/ $0\le i\le d-1$.
 Then one has $F\in\sF$. \par
\textup{(b)} Let\/ $\sC\subset\sE$ be a coresolving subcategory and
$E\in\sE$ be an object of\/ $\sC$\+coresolution dimension\/~$\le d$.
 Let\/ $0\rarrow E\rarrow C^0\rarrow\dotsb\rarrow C^{d-1}\rarrow C
\rarrow0$ be a $(d\nobreak+2)$\+term exact sequence in\/ $\sE$
with $C^i\in\sC$ for all\/ $0\le i\le d-1$.
 Then one has $C\in\sC$.
\end{lem}

\begin{proof}
 Various versions of this assertion can be found, on increasing levels
of generality, in~\cite[Lemma~2.1]{Zhu}, \cite[Proposition~2.3(1)]{Sto},
and~\cite[Corollary~A.5.2]{Pcosh}.
\end{proof}

\begin{cor} \label{intersection-co-resolution-dimension-cor}
 In the setting of Lemma~\textup{\ref{intersection-co-resolving-lemma}},
assume that both the additive categories\/ $\sE$ and\/ $\sE'$ are
weakly idempotent-complete.
 Then \par
\textup{(a)} in the context of
Lemma~\textup{\ref{intersection-co-resolving-lemma}(a)}, one has\/
$\rsd_{\sF/\sE}E'=\rsd_{\sF'/\sE'}E'$ for every object $E'\in\sE'$; \par
\textup{(b)} in the context of
Lemma~\textup{\ref{intersection-co-resolving-lemma}(b)}, one has\/
$\crd_{\sC/\sE}E'=\crd_{\sC'/\sE'}E'$ for every object $E'\in\sE'$.
\end{cor}

\begin{proof}
 This is~\cite[Corollary~A.5.5]{Pcosh}.
 The assertions follow from
Lemma~\ref{co-resol-dim-independent-of-co-resolution-lem}.
\end{proof}

\subsection{Cotorsion pairs} \label{appx-cotorsion-pairs-subsecn}
 Let $\sE$ be an exact category.
 As usual, we denote by $\Ext_\sE^*({-},{-})$ the Yoneda Ext groups
computed in the exact category~$\sE$.
 
 Let $\sF$ and $\sC\subset\sE$ be two classes of objects in~$\sE$.
 The notation $\sF^{\perp_1}\subset\sE$ stands for the class of all
objects $X\in\sE$ such that $\Ext_\sE^1(F,X)=0$ for all $F\in\sF$.
 Similarly, $\sF^{\perp_{\ge1}}\subset\sE$ is the class of all objects
$X\in\sE$ such that $\Ext_\sE^i(F,X)=0$ for all $F\in\sF$ and all
$i\ge1$.
 Dually, one denotes by ${}^{\perp_1}\sC\subset\sE$ the class of all
objects $Y\in\sE$ such that $\Ext_\sE^1(Y,C)=0$ for all $C\in\sC$.
 Similarly, ${}^{\perp_{\ge1}}\sC\subset\sE$ is the class of all objects
$Y\in\sE$ such that $\Ext_\sE^i(Y,C)=0$ for all $C\in\sC$ and all
$i\ge1$.

 A pair of classes of objects $(\sF,\sC)$ in $\sE$ is said to be
a \emph{cotorsion pair}~\cite{Sal} if $\sC=\sF^{\perp_1}$ and
$\sF={}^{\perp_1}\sC$.

\begin{lem} \label{hereditary-cotorsion-pairs-lemma}
 Let $(\sF,\sC)$ be a cotorsion pair in an exact category\/~$\sE$.
 Assume that the class\/ $\sF$ is generating in\/ $\sE$, while
the class\/ $\sC$ is cogenerating in\/~$\sE$.
 Then the following four conditions are equivalent:
\begin{enumerate}
\item the class\/ $\sF$ is closed under kernels of admissible
epimorphisms in\/~$\sE$;
\item the class\/ $\sC$ is closed under cokernels of admissible
monomorphisms in\/~$\sE$;
\item $\Ext^2_\sE(\sF,\sC)=0$ for all $F\in\sF$ and $C\in\sC$;
\item $\Ext^i_\sE(\sF,\sC)=0$ for all $F\in\sF$, \,$C\in\sC$, and
$i\ge1$.
\end{enumerate}
\end{lem}

\begin{proof}
 See~\cite[Lemma~4.25]{SaoSt}, \cite[Lemma~6.17]{Sto-ICRA},
\cite[Lemma~1.4]{PS4}, \cite[Lemma~1.3]{BHP}, \cite[Section~1 and
Lemma~1.1]{Pal}, or~\cite[Lemma~7.1]{PS6}.
\end{proof}

 A cotorsion pair $(\sF,\sC)$ in $\sE$ is said to be \emph{hereditary}
if it satisfies the assumptions and any one of the equivalent
conditions of Lemma~\ref{hereditary-cotorsion-pairs-lemma}.

 The following lemma complements
Lemma~\ref{hereditary-cotorsion-pairs-lemma} in a way relevant
for our exposition in~\cite{Pcosh} and~\cite{Pdomc}, as well as
in the present paper.

\begin{lem} \label{Hom-from-to-short-exact-sequences-lemma}
\textup{(a)} Let\/ $\sF$ be a generating (or more generally,
self-generating) class of objects closed under kernels of admissible
epimorphisms in an exact category\/~$\sE$.
 Let $C\in\sE$ be an object.
 Then one has $C\in\sF^{\perp_1}$ if and only if
$C\in\sF^{\perp_{\ge1}}$, and if and only if the functor\/
$\Hom_\sE({-},C)$ takes all short exact sequences in\/ $\sE$ with
the terms belonging to\/ $\sF$ to short exact sequences of abelian
groups. \par
\textup{(b)} Let\/ $\sC$ be a cogenerating (or more generally,
self-cogenerating) class of objects closed under cokernels of
admissible monomorphisms in an exact category\/~$\sE$.
 Let $F\in\sE$ be an object.
 Then one has $F\in{}^{\perp_1}\sC$ if and only if
$F\in{}^{\perp_{\ge1}}\sC$, and if and only if the functor\/
$\Hom_\sE(F,{-})$ takes all short exact sequences in\/ $\sE$ with
the terms belonging to\/ $\sC$ to short exact sequences of abelian
groups.
\end{lem}

\begin{proof}
 This is clear from the proofs of
Lemma~\ref{hereditary-cotorsion-pairs-lemma} as per the references
above (see, e.~g., \cite[proof of Lemma~1.3]{BHP}).
\end{proof}

 Let $(\sF,\sC)$ be a cotorsion pair in an exact category~$\sE$.
 Given an object $E\in\sE$, by a \emph{special precover sequence}
(for $E$ with respect to $(\sF,\sC)$) one means an (admissible)
short exact sequence
\begin{equation} \label{special-precover-sequence}
 0\lrarrow C'\lrarrow F\lrarrow E\lrarrow0
\end{equation}
in $\sE$ with objects $F\in\sF$ and $C'\in\sC$.
 Dually, a \emph{special preenvelope sequence} is a short exact
sequence
\begin{equation} \label{special-preenvelope-sequence}
 0\lrarrow E\lrarrow C\lrarrow F'\lrarrow0
\end{equation}
in $\sE$ with objects $C\in\sC$ and $F'\in\sF$.
 Collectively, the short exact
sequences~(\ref{special-precover-sequence}\+-%
\ref{special-preenvelope-sequence}) are referred to as
the \emph{approximation sequences}.

 A cotorsion pair $(\sF,\sC)$ in $\sE$ is said to be \emph{complete}
if every object $E\in\sE$ admits a special precover sequence and
a special preenvelope sequence.
 In fact, if the class $\sC$ is cogenerating, then existence of
special precover sequences with respect to $(\sF,\sC)$ implies
existence of special preenvelope sequences.
 Dually, if the class $\sF$ is generating, then existence of
special preenvelope sequences implies existence of special precover
sequences.
 These assertions are known as the \emph{Salce lemmas}~\cite{Sal};
see~\cite[Lemma~B.1.1]{Pcosh} for a more general formulation.
 Another reference is~\cite[Lemma~1.1]{PS4}.

 Given a class of objects $\sA$ in an ambient additive/exact
category $\sE$, we denote by $\sA^\oplus\subset\sE$ the class of
all direct summands of the object from $\sA$ in~$\sE$.
 The following direct summand lemma is convenient.

\begin{lem} \label{direct-summand-lemma}
 Let\/ $\sE$ be an exact category and $(\sF,\sC)$ be
an\/ $\Ext^1$\+orthogonal pair of classes of objects in\/ $\sE$,
that is, $\Ext^1_\sE(F,C)=0$ for all $F\in\sF$ and $C\in\sC$.
 Assume that a special precover short exact
sequence~\eqref{special-precover-sequence} and a special preenvelope
short exact sequence~\eqref{special-preenvelope-sequence} with
$F$, $F'\in\sF$ and $C$, $C'\in\sC$ exist for every object $E\in\sE$.
 Then the pair of classes $(\sF^\oplus,\sC^\oplus)$ is a complete
cotorsion pair in\/~$\sE$.
 In other words, one has\/ $\sF^{\perp_1}=\sC^\oplus$ and\/
${}^{\perp_1}\sC=\sF^\oplus$.
\end{lem}

\begin{proof}
 This easy lemma can be found in\cite[Lemma~1.2]{PS4},
\cite[Lemma~B.1.2]{Pcosh} or~\cite[Lemma~1.4]{Pal}.
\end{proof}

 We refer to~\cite[Section~11]{Bueh} for a discussion of projective
and injective objects in exact categories.
 Given an exact category $\sE$, we denote by $\sE_\prj\subset\sE$
the full subcategory of projective objects and by $\sE^\inj\subset\sE$
the full subcategory of injective objects in~$\sE$.

 For a cotorsion pair $(\sF,\sC)$ in an exact category $\sE$,
the intersection of the two classes $\sF\cap\sC$ is called
the \emph{kernel} (or in a different terminology, the \emph{core})
of the cotorsion pair $(\sF,\sC)$.

\begin{lem} \label{cotorsion-pair-kernel-as-injectives-projectives}
 Let $(\sF,\sC)$ be a complete cotorsion pair in an exact
category\/~$\sE$.
 In this context: \par
\textup{(a)} The exact category\/ $\sF$ has enough injective
objects.
 The kernel\/ $\sF\cap\sC$ of the cotorsion pair $(\sF,\sC)$ is
the class of all injective objects in\/~$\sF$. \par
\textup{(b)} The exact category\/ $\sC$ has enough projective
objects.
 The kernel\/ $\sF\cap\sC$ of the cotorsion pair $(\sF,\sC)$ is
the class of all projective objects in\/~$\sC$.
\end{lem}

\begin{proof}
 This easy and well-known lemma can be found
in~\cite[Lemma~1.9]{Pfltp}.
\end{proof}

\begin{lem} \label{cotorsion-pair-the-same-projectives-injectives}
 Let $(\sF,\sC)$ be a hereditary cotorsion pair in an exact
category\/~$\sE$.
 In this context: \par
\textup{(a)} The classes of projective objects in the exact categories\/
$\sF$ and\/ $\sE$ coincide, $\sE_\prj=\sF_\prj$.
 There are enough projective objects in\/ $\sF$ if and only if there are
enough projective objects in\/~$\sE$. \par
\textup{(b)} The classes of injective objects in the exact categories\/
$\sC$ and\/ $\sE$ coincide, $\sE^\inj=\sC^\inj$.
 There are enough injective objects in\/ $\sC$ if and only if there are
enough injective objects in\/~$\sE$.
\end{lem}

\begin{proof}
 This easy and well-known lemma is a special case
of~\cite[Lemma~7.1]{Pfltp}.
\end{proof}

\begin{lem} \label{nested-cotorsion-pairs-projectives-injectives}
 Let $(\sF_1,\sC_1)$ and $(\sF_2,\sC_2)$ be two \emph{nested} hereditary
complete cotorsion pairs in an exact category\/ $\sE$, that is,
$\sF_1\subset\sF_2$, or equivalently, $\sC_1\supset\sC_2$.
 Then the exact category\/ $\sF_2\cap\sC_1$ has both enough projective
and enough injective objects.
 The kernel\/ $\sF_1\cap\sC_1$ of the cotorsion pair $(\sF_1,\sC_1)$ is
the full subcategory of projective objects in\/ $\sF_2\cap\sC_1$, while
the kernel\/ $\sF_2\cap\sC_2$ of the cotorsion pair $(\sF_2,\sC_2)$ is
the full subcategory of injective objects in\/ $\sF_2\cap\sC_1$.
\end{lem}

\begin{proof}
 This is a special case of~\cite[Lemma~4.5]{Pfltp}.
\end{proof}

\subsection{Exotic derived categories in the sense of Positselski}
\label{appx-exotic-positselski-subsecn}
 Given an additive category $\sE$, we denote by $\Hot(\sE)$
the homotopy category of~$\sE$.
 So the objects of $\sE$ are cochain complexes in $\sE$, while
the morphisms in $\Hot(\sE)$ are the closed morphisms of degree~$0$
up to the cochain homotopy.

 As usual, $\Hot^\bb(\sE)\subset\Hot(\sE)$ denotes the full subcategory
in $\sE$ formed by bounded complexes, while $\Hot^+(\sE)$
and $\Hot^-(\sE)\subset\Hot(\sE)$ are the full subcategories of
complexes that are bounded below and bounded above, respectively.

 There are eleven conventional and exotic derived category symbols
$\star=\bb$, $+$, $-$, $\varnothing$ (the empty symbol), $\abs+$,
$\abs-$, $\abs$, $\co$, $\ctr$, $\bco$, and~$\bctr$.
 The notation $\Hot^\st(\sE)$ stands for
\begin{itemize}
\item $\Hot^\bb(\sE)$ when $\st=\bb$;
\item $\Hot^+(\sE)$ when $\st=+$ or~$\abs+$;
\item $\Hot^-(\sE)$ when $\st=-$ or~$\abs-$;
\item $\Hot(\sE)$ when $\st=\varnothing$, $\abs$, $\co$, $\ctr$,
$\bco$, or~$\bctr$.
\end{itemize}
 So $\Hot^\st(\sE)$ are triangulated categories.

 Now let $\sE$ be an exact category.
 The definitions of the conventional derived categories
$\Hot^\st(\sE)=\Hot^\st(\sE)/\Ac^\st(\sE)$ for $\st=\bb$, $+$, $-$,
or~$\varnothing$ can be found in~\cite{Neem0}, \cite{Kel},
or~\cite[Section~10]{Bueh}.
 Following the terminology in~\cite[Section~2.1]{Psemi}
and~\cite[Section~A.1]{Pcosh}, we make a distinction between
a possibly smaller class of \emph{exact} complexes and a possibly
wider class of \emph{acyclic} complexes in~$\sE$.
 The exact complexes are obtained by splicing (admissible) short
exact sequences.
 The acyclic complexes are the complexes homotopy equivalent to
exact complexes, or equivalently, the direct summands of exact
complexes.
 The notation $\Ac^\st(\sE)$ stands for the class of $\st$\+bounded
acyclic complexes.
 A morphism of complexes is called a \emph{quasi-isomorphism} if its
cone is acyclic.
 For an idempotent-complete exact category $\sE$, the classes of
exact complexes and acyclic complexes in $\sE$
coincide~\cite[Lemma~1.2]{Neem0}, \cite[Lemma~10.7]{Bueh}.

 The \emph{exotic derived categories in the sense of Positselski}
(with the symbols $\abs+$, $\abs-$, $\abs$, $\co$, and~$\ctr$)
are defined below in this section, following~\cite[Section~A.1]{Pcosh}.
 The \emph{coderived} and \emph{contraderived categories in
the sense of Becker} (with the symbols $\bco$ and~$\bctr$) are
defined in the next Section~\ref{appx-becker-subsecn},
following~\cite[Section~B.7]{Pcosh}.
 We refer to~\cite[Remark~9.2]{PS4} or~\cite[Section~7]{Pksurv} for
a historical and philosophical discussion.

 The following definitions go back
to~\cite[Sections~2.1 and~4.1]{Psemi}; see also~\cite{KLN}.
 A survey exposition can be found in~\cite[Sections~7 and~9.3]{Pksurv}.

 Let $\sE$ be an exact category.
 We are interested in (degreewise admissible) \emph{short exact
sequences of complexes} $0\rarrow K^\bu\rarrow L^\bu\rarrow M^\bu
\rarrow0$ in~$\sE$.
 Any such short exact sequence of complexes can be viewed as
a bicomplex with three rows, so the total complex $\Tot(K^\bu\to
L^\bu\to M^\bu)$ can be attached to it.

 A complex in $\sE$ is said to be \emph{absolutely acyclic} if
it belongs to the minimal thick subcategory of $\Hot(\sE)$
containing the totalizations of short exact sequences in~$\sE$.
 The full subcategory of absolutely acyclic complexes is denoted
by $\Ac^\abs(\sE)\subset\Hot(\sE)$.
 All absolutely acyclic complexes are always acyclic, so one has
$\Ac^\abs(\sE)\subset\Ac(\sE)$.
 The triangulated Verdier quotient category
$$
 \sD^\abs(\sE)=\Hot(\sE)/\Ac^\abs(\sE)
$$
is called the \emph{absolute derived category} of~$\sE$.

 The definitions from the previous paragraph, restricted to the realms
of bounded below (respectively, bounded above) complexes in $\sE$,
produce the thick subcategory $\Ac^{\abs+}(\sE)\subset\Hot^+(\sE)$
(resp., $\Ac^{\abs-}(\sE)\subset\Hot^-(\sE)$).
 The related triangulated Verdier quotient categories
$$
 \sD^\st(\sE)=\Hot(\sE)/\Ac^\st(\sE), \qquad
 \text{$\st=\abs+$ or~$\abs-$}
$$
are called the \emph{bounded below} and the \emph{bounded above
absolute derived categories} of~$\sE$.
 In fact, the inclusions $\Hot^\pm(\sE)\rarrow\Hot(\sE)$ induce
fully faithful triangulated functors $\sD^{\abs\pm}(\sE)\rarrow
\sD^\abs(\sE)$ \,\cite[Lemma~A.1.2]{Pcosh}, and it follows that any
absolutely acyclic complex that happens to be bounded below (resp.,
above) is also absolutely acyclic in the world of complexes bounded
below (resp., above), i.~e., belongs to $\Ac^{\abs+}(\sE)$
(resp., $\Ac^{\abs-}(\sE)$).

 The following definition is often given under the assumptions that
the exact category $\sE$ has exact functors of infinite direct sums
(i.~e., all infinite direct sums exist and preserve admissible short
exact sequences), but this is not strictly necessary.
 A complex in $\sE$ is said to be \emph{coacyclic} (\emph{in the sense
of Positselski}) or \emph{Positselski-coacyclic} if it belongs to
the minimal thick subcategory of $\Hot(\sE)$ containing
the totalizations of short exact sequences in $\sE$ and closed under
infinite direct sums.
 The thick subcategory of Positselski-coacyclic complexes is denoted
by $\Ac^\co(\sE)\subset\Hot(\sE)$; and the related triangulated
Verdier quotient category
$$
 \sD^\co(\sE)=\Hot(\sE)/\Ac^\co(\sE)
$$
is called the \emph{coderived category} (\emph{in the sense of 
Positselski}) or the \emph{Positselski coderived category} of
an exact category~$\sE$.

 Obviously, one always has $\Ac^\abs(\sE)\subset\Ac^\co(\sE)$.
 When the exact category $\sE$ has exact functors of infinite direct
sums, the inclusion $\Ac^\co(\sE)\subset\Ac(\sE)$ also holds.
 So there are natural triangulated Verdier quotient functors
$$
 \xymatrix{
  \sD^\abs(\sE) \ar@{->>}[r] & \sD^\co(\sE) \ar@{->>}[r] &
  \sD(\sE)
 }
$$
in this case.
 Furthermore, under the same assumptions, the full subcategory of
bounded below complexes in $\Ac^\co(\sE)$ is equivalent to
the conventional derived category $\sD^+(\sE)$
\,\cite[Lemma~A.1.3(a)]{Pcosh}.

 The next definition is the opposite version of the previous one.
 A complex in $\sE$ is said to be \emph{contraacyclic} (\emph{in
the sense of Positselski}) or \emph{Positselski-contraacyclic} if it
belongs to the minimal thick subcategory of $\Hot(\sE)$ containing
the totalizations of short exact sequences in $\sE$ and closed under
infinite products.
 The thick subcategory of Positselski-contraacyclic complexes is denoted
by $\Ac^\ctr(\sE)\subset\Hot(\sE)$; and the related triangulated
Verdier quotient category
$$
 \sD^\ctr(\sE)=\Hot(\sE)/\Ac^\ctr(\sE)
$$
is called the \emph{contraderived category} (\emph{in the sense of
Positselski}) or the \emph{Positselski contraderived category} of
an exact category~$\sE$.
 
 Obviously, one always has $\Ac^\abs(\sE)\subset\Ac^\ctr(\sE)$.
 When the exact category $\sE$ has exact functors of infinite products
(i.~e., all infinite products exist and preserve admissible short
exact sequences in~$\sE$), the inclusion $\Ac^\ctr(\sE)\subset\Ac(\sE)$
also holds.
 So there are natural triangulated Verdier quotient functors
$$
 \xymatrix{
  \sD^\abs(\sE) \ar@{->>}[r] & \sD^\ctr(\sE) \ar@{->>}[r] &
  \sD(\sE)
 }
$$
in this case.
 Furthermore, under the same assumptions, the full subcategory of
bounded above complexes in $\Ac^\ctr(\sE)$ is equivalent to
the conventional derived category $\sD^-(\sE)$
\,\cite[Lemma~A.1.3(b)]{Pcosh}.

\begin{prop} \label{infinite-co-resolutions-minus-plus-derived}
\textup{(a)} Let\/ $\sE$ be an exact category and\/ $\sF\subset\sE$ be
a resolving subcategory.
 Then the inclusion functor\/ $\sF\rarrow\sE$ induces a triangulated
equivalence of the bounded above conventional derived categories,
$\sD^-(\sF)\simeq\sD^-(\sE)$. \par
\textup{(b)} Let\/ $\sE$ be an exact category and\/ $\sC\subset\sE$ be
a coresolving subcategory.
 Then the inclusion functor\/ $\sC\rarrow\sE$ induces a triangulated
equivalence of the bounded below conventional derived categories,
$\sD^+(\sC)\simeq\sD^+(\sE)$.
\end{prop}

\begin{proof}
 This is a well-known result.
 Parts~(a) and~(b) are dual to each other.
 Part~(b) is~\cite[Lemma~I.4.6(1)]{HarRD}, \cite[Theorem~12.1]{Kel}
or~\cite[Proposition~13.2.2(i)]{KS}.
 Part~(a) is~\cite[Proposition~A.3.1(a)]{Pcosh}.
\end{proof}

\begin{prop} \label{infinite-co-resolutions-contra-co-derived}
\textup{(a)} Let\/ $\sE$ be an exact category with exact functors of
infinite products and\/ $\sF\subset\sE$ be a resolving subcategory
closed under infinite products.
 Then the inclusion functor\/ $\sF\rarrow\sE$ induces a triangulated
equivalence of the Positselski contraderived categories,
$\sD^\ctr(\sF)\simeq\sD^\ctr(\sE)$. \par
\textup{(b)} Let\/ $\sE$ be an exact category with exact functors of
infinite direct sums and\/ $\sC\subset\sE$ be a coresolving subcategory
closed under infinite direct sums.
 Then the inclusion functor\/ $\sC\rarrow\sE$ induces a triangulated
equivalence of the Positselski coderived categories,
$\sD^\co(\sC)\simeq\sD^\co(\sE)$.
\end{prop}

\begin{proof}
 Part~(a) is~\cite[Proposition~A.3.1(b)]{Pcosh}
(see~\cite[Proposition~A.2.1]{Pcosh} for a more general result).
 Part~(b) is dual.
 A far-reaching generalization can be found in~\cite[Theorem~7.11 and
Corollary~7.12]{Pedg}.
\end{proof}

\begin{prop} \label{finite-co-resolutions}
\textup{(a)} Let\/ $\sE$ be an exact category and\/ $\sF\subset\sE$ be
a resolving subcategory such that the\/ $\sF$\+resolution dimensions
of all the objects of\/ $\sE$ are finite and bounded by a fixed integer.
 Then, for any conventional or exotic derived category symbol\/
$\st=\bb$, $+$, $-$, $\varnothing$, $\abs+$, $\abs-$, $\co$,
or\/~$\ctr$, the inclusion functor\/ $\sF\rarrow\sE$ induces
a triangulated equivalence\/ $\sD^\st(\sF)\rarrow\sD^\st(\sE)$. \par
 Here, in the case of\/ $\st=\co$, it is presumed that\/ $\sE$ is
an exact category with infinite direct sums and the resolving
subcategory\/ $\sF\subset\sE$ is closed under infinite direct sums.
 In the case of\/ $\st=\ctr$, it is presumed that\/ $\sE$ is an exact
category with infinite products and the resolving subcategory\/
$\sF\subset\sE$ is closed under infinite products. \par
\textup{(b)} Let\/ $\sE$ be an exact category and\/ $\sC\subset\sE$ be
a coresolving subcategory such that the\/ $\sC$\+coresolution dimensions
of all the objects of\/ $\sE$ are finite and bounded by a fixed integer.
 Then, for any conventional or exotic derived category symbol\/
$\st=\bb$, $+$, $-$, $\varnothing$, $\abs+$, $\abs-$, $\co$,
or\/~$\ctr$, the inclusion functor\/ $\sC\rarrow\sE$ induces
a triangulated equivalence\/ $\sD^\st(\sC)\rarrow\sD^\st(\sE)$. \par
 Here, in the case of\/ $\st=\co$, it is presumed that\/ $\sE$ is
an exact category with infinite direct sums and the coresolving
subcategory\/ $\sC\subset\sE$ is closed under infinite direct sums.
 In the case of\/ $\st=\ctr$, it is presumed that\/ $\sE$ is an exact
category with infinite products and the coresolving subcategory\/
$\sC\subset\sE$ is closed under infinite products.
\end{prop}

\begin{proof}
 In part~(a), the case $\st=-$ is covered by
Proposition~\ref{infinite-co-resolutions-minus-plus-derived}(a)
and the case $\st=\ctr$ is covered by
Proposition~\ref{infinite-co-resolutions-contra-co-derived}(a),
which hold under less restrictive assumptions.
 Dually, in part~(b), the cases $\st=+$ and $\st=\co$ hold under
less restrictive assumptions by
Propositions~\ref{infinite-co-resolutions-minus-plus-derived}(b)
and~\ref{infinite-co-resolutions-contra-co-derived}(b).

 The case $\st=\bb$ in part~(b) is~\cite[Proposition~13.2.2(ii)]{KS}.
 The case $\st=\varnothing$ in part~(b) goes back
to~\cite[Lemma~I.4.6(2)]{HarRD}.
 The case $\st=\varnothing$ in part~(a)
is~\cite[Proposition~5.14]{Sto}.

 All the assertions of part~(a) are covered
by~\cite[Proposition~A.5.8]{Pcosh}.
 Part~(b) is dual.
 For a far-reaching generalization, see~\cite[Theorem~6.6
and Corollary~6.20]{Pedg}.
\end{proof}

\subsection{Coderived and contraderived categories in the sense
of Becker} \label{appx-becker-subsecn}
 The concepts defined in this section go back to the work of
J\o rgensen~\cite{Jorg}, Krause~\cite{Kra}, Neeman~\cite{Neem},
and Becker~\cite{Bec}.
 Our exposition is based on~\cite[Section~B.7]{Pcosh}.

 For comparison between the contexts of the definitions below and
those of Section~\ref{appx-exotic-positselski-subsecn}, it is helpful
to keep in mind that the functors of infinite direct sums are exact
in any exact category with enough injective objects where infinite
direct sums exist.
 Dually, the functors of infinite products are exact in any exact
category with enough projective objects where infinite products
exist~\cite[Remark~5.2]{Pedg}.

 For any two complexes $A^\bu$ and $B^\bu$ in an additive category
$\sE$, we denote by $\Hom_\sE(A^\bu,B^\bu)$ the total complex of
the bicomplex of abelian groups $\Hom_\sE(A^i,B^j)$, \,$i$, $j\in\boZ$,
constructed by taking infinite products along the diagonals.
 So one has a natural isomorphism of abelian groups
$\Hom_{\Hot(\sE)}(A^\bu,B^\bu[n])\simeq
H^n\Hom_\sE(A^\bu,B^\bu)$ for every integer $n\in\boZ$.

 Let $\sE$ be an exact category.
 A complex $A^\bu$ is $\sE$ is said to be \emph{Becker-coacyclic} if,
for every complex of injective objects $J^\bu\in\Hot(\sE^\inj)$,
the complex of abelian groups $\Hom_\sE(A^\bu,J^\bu)$ is acyclic.
 The thick subcategory of Becker-coacyclic complexes is denoted by
$\Ac^\bco(\sE)\subset\Hot(\sE)$; and the related triangulated Verdier
quotient category
$$
 \sD^\bco(\sE)=\Hot(\sE)/\Ac^\bco(\sE)
$$
is called the \emph{coderived category in the sense of Becker} or
the \emph{Becker coderived category} of an exact category~$\sE$.

 One always has $\Ac^\abs(\sE)\subset\Ac^\co(\sE)\subset\Ac^\bco(\sE)$
\,\cite[Lemma~B.7.1(a,b)]{Pcosh}.
 So there are natural triangulated Verdier quotient functors
$$
 \xymatrix{
  \sD^\abs(\sE) \ar@{->>}[r] & \sD^\co(\sE) \ar@{->>}[r] &
  \sD^\bco(\sE).
 }
$$
 It is an \emph{open problem} whether all Becker-coacyclic complexes
are acyclic in general.
 See~\cite[Lemmas~B.7.2 and~B.7.3, and Remark~B.7.4]{Pcosh} for
a discussion.

 A complex $B^\bu$ is $\sE$ is said to be \emph{Becker-contraacyclic}
if, for every complex of projective objects $P^\bu\in\Hot(\sE_\prj)$,
the complex of abelian groups $\Hom_\sE(P^\bu,B^\bu)$ is acyclic.
 The thick subcategory of Becker-contraacyclic complexes is denoted by
$\Ac^\bctr(\sE)\subset\Hot(\sE)$; and the related triangulated Verdier
quotient category
$$
 \sD^\bctr(\sE)=\Hot(\sE)/\Ac^\bctr(\sE)
$$
is called the \emph{contraderived category in the sense of Becker} or
the \emph{Becker contraderived category} of an exact category~$\sE$.

 One always has $\Ac^\abs(\sE)\subset\Ac^\ctr(\sE)\subset\Ac^\bctr(\sE)$
\,\cite[Lemma~B.7.1(a,c)]{Pcosh}.
 So there are natural triangulated Verdier quotient functors
$$
 \xymatrix{
  \sD^\abs(\sE) \ar@{->>}[r] & \sD^\ctr(\sE) \ar@{->>}[r] &
  \sD^\bctr(\sE).
 }
$$
 Once again, it is an \emph{open problem} whether all
Becker-contraacyclic complexes are acyclic in general.
 See~\cite[Lemmas~B.7.2 and~B.7.3, and Remark~B.7.4]{Pcosh}.

\begin{prop} \label{infinite-co-resolutions-becker-derived}
\textup{(a)} Let\/ $\sE$ be an exact category with countable products
and\/ $\sF\subset\sE$ be a resolving subcategory closed under
countable products.
 Then the inclusion functor\/ $\sF\rarrow\sE$ induces a triangulated
equivalence of the Becker contraderived categories,
$\sD^\bctr(\sF)\simeq\sD^\bctr(\sE)$. \par
\textup{(b)} Let\/ $\sE$ be an exact category with countable direct sums
and\/ $\sC\subset\sE$ be a coresolving subcategory closed under
countable direct sums.
 Then the inclusion functor\/ $\sC\rarrow\sE$ induces a triangulated
equivalence of the Becker coderived categories,
$\sD^\bco(\sC)\simeq\sD^\bco(\sE)$.
\end{prop}

\begin{proof}
 Part~(a) is~\cite[Proposition~B.7.11]{Pcosh}.
 Part~(b) is dual.
\end{proof}

\begin{prop} \label{finite-co-resolutions-becker-derived-I}
\textup{(a)} Let\/ $\sE$ be a weakly idempotent-complete exact category
and\/ $\sF\subset\sE$ be a resolving subcategory closed under
direct summands.
 Assume that the\/ $\sF$\+resolution dimensions of all the objects of\/
$\sE$ are finite and bounded by a fixed integer.
 Then the inclusion functor\/ $\sF\rarrow\sE$ induces a triangulated
equivalence of the Becker contraderived categories,
$\sD^\bctr(\sF)\simeq\sD^\bctr(\sE)$. \par
\textup{(b)} Let\/ $\sE$ be a weakly idempotent-complete exact category
and\/ $\sC\subset\sE$ be a coresolving subcategory closed under
direct summands.
 Assume that the\/ $\sC$\+coresolution dimensions of all the objects
of\/ $\sE$ are finite and bounded by a fixed integer.
 Then the inclusion functor\/ $\sC\rarrow\sE$ induces a triangulated
equivalence of the Becker coderived categories,
$\sD^\bco(\sC)\simeq\sD^\bco(\sE)$.
\end{prop}

\begin{proof}
 Part~(a) is~\cite[Proposition~B.7.9]{Pcosh}.
 Part~(b) is dual.
\end{proof}

\begin{prop} \label{finite-co-resolutions-becker-derived-II}
\textup{(a)} Let\/ $\sE$ be a weakly idempotent-complete exact category
with enough injective objects and $(\sF,\sC)$ be a hereditary complete
cotorsion pair in\/~$\sE$.
 Assume that the\/ $\sF$\+resolution dimensions of all the objects
of\/ $\sE$ are finite and bounded by a fixed integer.
 Then the inclusion functor\/ $\sF\rarrow\sE$ induces a triangulated
equivalence of the Becker coderived categories,
$\sD^\bco(\sF)\simeq\sD^\bco(\sE)$. \par
\textup{(b)} Let\/ $\sE$ be a weakly idempotent-complete exact category
with enough projective objects and $(\sF,\sC)$ be a hereditary complete
cotorsion pair in\/~$\sE$.
 Assume that the\/ $\sC$\+coresolution dimensions of all the objects
of\/ $\sE$ are finite and bounded by a fixed integer.
 Then the inclusion functor\/ $\sC\rarrow\sE$ induces a triangulated
equivalence of the Becker contraderived categories,
$\sD^\bctr(\sC)\simeq\sD^\bctr(\sE)$.
\end{prop}

\begin{proof}
 Part~(b) is~\cite[Proposition~B.7.10]{Pcosh}.
 Part~(a) is dual.
\end{proof}

\subsection{Homotopy categories of injective/projective objects}
\label{appx-homotopy-of-inj-proj-subsecn}
 ``Derived categories of the second kind'' is a common name for
the coderived and contraderived categories in the sense of
Positselski and in the sense of Becker, as well as for the absolute
derived categories.
 The following theorems describe the situations when derived categories
of the second kind are well-behaved.

 The \emph{injective dimension} of an object $E$ in an exact category
$\sE$ is defined as the supremum of all integers $n\ge0$ for which
there exists an object $X\in\sE$ with $\Ext_\sE^n(X,E)\ne0$.
 If there are enough injective objects in $\sE$, then the class of
injective objects $\sE^\inj$ is coresolving in $\sE$, and
the injective dimension coincides with the coresolution dimension
with respect to~$\sE^\inj$.

 Dually, the \emph{projective dimension} of an object $E$ in an exact
category $\sE$ is defined as the supremum of all integers $n\ge0$ for
which there exists an object $Y\in\sE$ with $\Ext_\sE^n(E,Y)\ne0$.
 If there are enough projective objects in $\sE$, then the class of
projective objects $\sE_\prj$ is resolving in $\sE$, and
the projective dimension coincides with the resolution dimension
with respect to~$\sE_\prj$.

 We use the notation with equality signs between derived category
symbols in the superindices of the letter $\sD$ in order to indicate
that the respective exotic derived categories agree because
the corresponding classes of acyclic complexes coincide (as full
subcategories of the homotopy category).
 So this kind of notation (as in the theorems in this section) in itself
constitutes an assertion.

\begin{thm} \label{countable-co-products-co-contra-derived-theorem}
\textup{(a)} Let\/ $\sA$ be an exact category with countable direct sums
and enough injective objects.
 Assume that the countable direct sums of injective objects have
finite injective dimensions in\/~$\sA$.
 Then the classes of Positselski-coacyclic and Becker-coacyclic
complexes in the exact category\/ $\sA$ coincide,
$$
 \Ac^\co(\sA)=\Ac^\bco(\sA).
$$
 The inclusion of additive/exact categories\/ $\sA^\inj\rarrow\sA$
induces an equivalence of triangulated categories
$$
 \Hot(\sA^\inj)\simeq\sD^{\co=\bco}(\sA).
$$ \par
\textup{(b)} Let\/ $\sB$ be an exact category with countable direct
products and enough projective objects.
 Assume that the countable products of projective objects have
finite projective dimensions in\/~$\sB$.
 Then the classes of Positselski-contraacyclic and Becker-contraacyclic
complexes in the exact category\/ $\sB$ coincide,
$$
 \Ac^\ctr(\sB)=\Ac^\bctr(\sB).
$$
 The inclusion of additive/exact categories\/ $\sB_\prj\rarrow\sB$
induces an equivalence of triangulated categories
$$
 \Hot(\sB_\prj)\simeq\sD^{\ctr=\bctr}(\sB).
$$
\end{thm}

\begin{proof}
 See, e.~g., \cite[Proposition~7.6]{Pphil}
or~\cite[Theorem~B.7.7]{Pcosh}.
 Notice that existence of enough injective objects and countable direct
sums in $\sA$ implies exactness of such countable direct sums (and
similarly for projective objects and products);
see~\cite[Remark~5.2]{Pedg}.

 Under somewhat more restrictive assumptions, a version of part~(a)
can be obtained by applying
Propositions~\ref{infinite-co-resolutions-contra-co-derived}(b)
and~\ref{infinite-co-resolutions-becker-derived}(b) to $\sE=\sA$
and $\sC=\sA^\inj$.
 Similarly, a more restrictive version of part~(b) can be obtained by
applying Propositions~\ref{infinite-co-resolutions-contra-co-derived}(a)
and~\ref{infinite-co-resolutions-becker-derived}(a) to $\sE=\sB$
and $\sF=\sA_\prj$.
\end{proof}

 The \emph{homological dimension} of an exact category $\sE$ is
the supremum of all integers $n\ge0$ for which there exists a pair
of objects $X$, $Y\in\sE$ with $\Ext^n_\sE(X,Y)\ne0$.

\begin{lem} \label{psemi-remark21}
 Let\/ $\sE$ be an exact category of finite homological dimension.
 Then the classes of absolutely acyclic and acyclic complexes
in\/ $\sE$ coincide,
$$
 \Ac^\abs(\sE)=\Ac(\sE).
$$
 If the exact category\/ $\sE$ has exact functors of infinite direct
sums (resp., infinite products), then the class of Positselski-coacyclic
(resp., Positselski-contraacyclic) complexes in\/ $\sE$ also coincides
with the class of absolutely acyclic complexes and with the class of
acyclic complexes:
$$
 \Ac^\abs(\sE)=\Ac^\co(\sE)=\Ac(\sE) \quad\text{or}\quad
 \Ac^\abs(\sE)=\Ac^\ctr(\sE)=\Ac(\sE),
$$
respectively.
\end{lem}

\begin{proof}
 This is~\cite[Remark~2.1]{Psemi}.
\end{proof}

\begin{thm} \label{fin-homol-dim-co-contra-derived-theorem}
\textup{(a)} Let\/ $\sA$ be an exact category of finite homological
dimension with enough injective objects.
 Then the classes of absolutely acyclic, Positselski-coacyclic,
Becker-coacyclic, and acyclic complexes in the exact category\/ $\sA$
coincide,
$$
 \Ac^\abs(\sA)=\Ac^\co(\sA)=\Ac^\bco(\sA)=\Ac(\sA).
$$
 The inclusion of additive/exact categories\/ $\sA^\inj\rarrow\sA$
induces an equivalence of triangulated categories
$$
 \Hot(\sA^\inj)\simeq\sD^{\abs=\co=\bco=\varnothing}(\sA).
$$ \par
\textup{(b)} Let\/ $\sB$ be an exact category of finite homological
dimension with enough projective objects.
 Then the classes of absolutely acyclic, Positselski-contraacyclic,
Becker-contraacyclic, and acyclic complexes in the exact category\/
$\sB$ coincide,
$$
 \Ac^\abs(\sB)=\Ac^\ctr(\sB)=\Ac^\bctr(\sB)=\Ac(\sB).
$$
 The inclusion of additive/exact categories\/ $\sB_\prj\rarrow\sB$
induces an equivalence of triangulated categories
$$
 \Hot(\sB_\prj)\simeq\sD^{\abs=\ctr=\bctr=\varnothing}(\sB).
$$
\end{thm}

\begin{proof}
 See, e.~g., \cite[Proposition~7.5]{Pphil}
or~\cite[Theorem~B.7.6]{Pcosh}.

 Alternatively, to prove part~(a) one can combine the results of
Lemma~\ref{psemi-remark21}, Proposition~\ref{finite-co-resolutions}(b),
and Proposition~\ref{finite-co-resolutions-becker-derived-I}(b)
(for $\sE=\sA$ and $\sC=\sA^\inj$).
 Dually, to prove part~(b) one can combine the results of
Lemma~\ref{psemi-remark21}, Proposition~\ref{finite-co-resolutions}(a),
and Proposition~\ref{finite-co-resolutions-becker-derived-I}(a)
(for $\sE=\sB$ and $\sF=\sB_\prj$).
\end{proof}


\bigskip
\begin{thebibliography}{99}
\smallskip

\bibitem{AR}
 J.~Ad\'amek, J.~Rosick\'y.
   Locally presentable and accessible categories.
London Math.\ Society Lecture Note Series~189,
Cambridge University Press, 1994.

\bibitem{AJL0}
 L.~Alonso Tarr\'\i o, A.~Jerem\'\i as L\'opez, J.~Lipman.
   Local homology and cohomology on schemes.
\textit{Annales Sci.\ de l'\'Ecole Norm.\ Sup.~(4)} \textbf{30},
\#1, p.~1--39, 1997.
 Corrections available from
\texttt{https://www.math.purdue.edu/\textasciitilde
jlipman/papers/homologyfix.pdf}

\bibitem{AJL}
 L.~Alonso Tarr\'\i o, A.~Jerem\'\i as L\'opez, J.~Lipman.
   Duality and flat base change on formal schemes.
In: Studies in duality on Noetherian formal schemes and non-Noetherian
ordinary schemes, p.~3--90, \emph{Contemporary Math.}\ \textbf{244},
American Math.\ Society, Providence, 1999.
\texttt{arXiv:alg-geom/9708006}
 Correction, \emph{Proceedings of the American Math.\ Society}
\textbf{131}, \#2, p.~351--357, 2003.  \texttt{arXiv:math.AG/0106239}

\bibitem{AJL2}
 L.~Alonso Tarr\'\i o, A.~Jerem\'\i as L\'opez, J.~Lipman.
   Greenlees--May duality on formal schemes.
In: Studies in duality on Noetherian formal schemes and non-Noetherian
ordinary schemes, p.~93--112, \emph{Contemporary Math.}\ \textbf{244},
American Math.\ Society, Providence, 1999.

\bibitem{AJPV}
 L.~Alonso Tarr\'\i o, A.~Jerem\'\i as L\'opez,
M.~P\'erez Rodr\'\i guez, M.~J.~Vale Gonsalves.
   On the existence of a compact generator on the derived category
of a Noetherian formal scheme.
\textit{Appl.\ Categorical Struct.}\ \textbf{19}, \#6, p.~865--877,
2011.  \texttt{arXiv:0905.2063 [math.AG]}

\bibitem{AB}
 M.~Auslander, R.-O.~Buchweitz.
   The homological theory of maximal Cohen--Macaulay
approximations.
\textit{M\'emoires de la Soc.\ Math.\ de France} \textbf{38},
1989, p.~5--37.

\bibitem{BHP}
 S.~Bazzoni, M.~Hrbek, L.~Positselski.
   Fp\+projective periodicity.
\textit{Journ.\ of Pure and Appl.\ Algebra} \textbf{228}, \#3,
article ID~107497, 24~pp., 2024.  \texttt{arXiv:2212.02300 [math.CT]}

\bibitem{Bec}
 H.~Becker.
   Models for singularity categories.
\textit{Advances in Math.}\ \textbf{254}, p.~187--232, 2014.
\texttt{arXiv:1205.4473 [math.CT]}

\bibitem{Beil}
 A.~Beilinson.
   Remarks on topological algebras.
\textit{Moscow Math.\ Journ.}\ \textbf{8}, \#1, p.~1--20, 2008.
\texttt{arXiv:0711.2527 [math.QA]}

\bibitem{BBD}
  A.~Beilinson, J.~Bernstein, P.~Deligne, O.~Gabber.
    Faisceaux pervers.
\textit{Ast\'erisque} \textbf{100}, 1982.
Deuxi\`eme \'edition, Soc.\ Math.\ de France, 2018.

\bibitem{BD2}
 A.~Beilinson, V.~Drinfeld.
   Quantization of Hitchin's integrable system and Hecke
eigensheaves.  February 2000.  Available from
\texttt{https://web.ma.utexas.edu/users/benzvi/Langlands.html} or
\texttt{https://math.uchicago.edu/\textasciitilde
drinfeld/langlands.html}

\bibitem{BG}
 D.~J.~Benson, K.~R.~Goodearl.
   Periodic flat modules, and flat modules for finite groups.
\textit{Pacific Journ.\ of Math.}\ \textbf{196}, \#1, p.~45--67, 2000.

\bibitem{Bour}
 N.~Bourbaki.
   Commutative algebra.  Chapters~1--7.  Translated from the French.
Hermann, Paris; Addison-Wesley Publishing Co., Reading, MA, 1972;
Springer-Verlag, Berlin, 1989--1998.

\bibitem{Bueh}
 T.~B\"uhler.
   Exact categories.
\textit{Expositiones Math.}\ \textbf{28}, \#1, p.~1--69, 2010.
\texttt{arXiv:0811.1480 [math.HO]}

\bibitem{Del}
 P.~Deligne.
   Cohomologie \`a supports propres.
SGA4, Tome 3. \textit{Lecture Notes in Math.}\ \textbf{305},
Springer-Verlag, Berlin--Heidelberg--New York, 1973, p.~250--480.

\bibitem{DS}
 B.~Dierolf, D.~Sieg.
   Exact subcategories of the category of locally convex spaces.
\textit{Math.\ Nachrichten} \textbf{285}, \#17--18, p.~2106--2119, 2012.

\bibitem{DG}
 W.~G.~Dwyer, J.~P.~C.~Greenlees.
   Complete modules and torsion modules.
\textit{American Journ.\ of Math.}\ \textbf{124}, \#1,
p.~199--220, 2002.

\bibitem{EP}
 A.~I.~Efimov, L.~Positselski.
   Coherent analogues of matrix factorizations and relative
singularity categories.
\textit{Algebra and Number Theory} \textbf{9}, \#5, p.~1159--1292,
2015.  \texttt{arXiv:1102.0261 [math.CT]}

\bibitem{En2}
 E.~Enochs.
   Flat covers and flat cotorsion modules.
\textit{Proc.\ of the Amer.\ Math.\ Soc.}\ \textbf{92}, \#2,
p.~179--184, 1984.

\bibitem{EE}
 E.~Enochs, S.~Estrada.
   Relative homological algebra in the category of quasi-coherent
sheaves.
\textit{Advances in Math.}\ \textbf{194}, \#2, p.~284--295, 2005.

\bibitem{ES}
 S.~Estrada, A.~Sl\'avik.
   Quillen equivalent models for the derived category of flats
and the resolution property.
\textit{Journ.\ Australian Math.\ Soc.}\ \textbf{110}, \#3,
p.~302--320, 2021.  \texttt{arXiv:1708.05913 [math.AT]}.

\bibitem{GZ}
 P.~Gabriel, M.~Zisman.
   Calculus of fractions and homotopy theory.
Springer-Verlag, Berlin--Heidelberg--New York, 1967.

\bibitem{GL}
 W.~Geigle, H.~Lenzing.
   Perpendicular categories with applications to representations
and sheaves.
\textit{Journ.\ of Algebra} \textbf{144}, \#2, p.~273--343, 1991.

\bibitem{God}
 R.~Godement.
   Topologie alg\'ebrique et th\'eorie des faisceaux.
Paris, Hermann, 1958--73.

\bibitem{GM}
 J.~P.~C.~Greenlees, J.~P.~May.
   Derived functors of $I$\+adic completion and local homology.
\textit{Journ.\ of Algebra} \textbf{149}, \#2, p.~438--453, 1992.

\bibitem{EGAII}
 A.~Grothendieck, J.~Dieudonn\'e.
   \'El\'ements de g\'eom\'etrie alg\'ebrique:~II.
\'Etude globale \'el\'ementaire de quelques classes de morphismes.
\textit{Publications Mathematiques de l'IH\'ES} \textbf{8},
p.~5--222, 1961.

\bibitem{EGAIV1}
 A.~Grothendieck, J.~Dieudonn\'e.
   \'El\'ements de g\'eom\'etrie alg\'ebrique:~IV.
\'Etude locale des sch\'emas et des morphismes des sch\'emas,
Premi\`ere partie.
\textit{Publications Mathematiques de l'IH\'ES} \textbf{20},
p.~5--259, 1964.

\bibitem{EGAIV2}
 A.~Grothendieck, J.~Dieudonn\'e.
   \'El\'ements de g\'eom\'etrie alg\'ebrique~IV.
\'Etude locale des sch\'emas et des morphismes des sch\'emas,
Seconde partie.
\textit{Publications Mathematiques de l'IH\'ES} \textbf{24},
p.~5--231, 1965.

\bibitem{EGAIV4}
 A.~Grothendieck, J.~Dieudonn\'e.
   \'El\'ements de g\'eom\'etrie alg\'ebrique:~IV.
\'Etude locale des sch\'emas et des morphismes des sch\'emas,
Quatri\`eme partie.
\textit{Publications Mathematiques de l'IH\'ES} \textbf{32},
p.~5--361, 1967.

\bibitem{EGAI}
 A.~Grothendieck, J.~A.~Dieudonn\'e.
   \'El\'ements de g\'eom\'etrie alg\'ebrique.~I.
Grundlehren der mathematischen Wissenschaften, 166.
Springer-Verlag, Berlin--Heidelberg--New York, 1971.

\bibitem{HarRD}
 R.~Hartshorne.
   Residues and duality.
\textit{Lecture Notes in Math.}\ \textbf{20}, Springer, 1966.

\bibitem{HarAG}
 R.~Hartshorne.
   Algebraic geometry.
Graduate Texts in Math., 52, Springer-Verlag,
New York--Heidelberg, 1977.

\bibitem{SP}
 A.~J.~de~Jong et al.
   The Stacks Project.
Available from \texttt{https://stacks.math.columbia.edu/}

\bibitem{Jorg}
 P.~J\o rgensen.
   The homotopy category of complexes of projective modules.
\textit{Advances in Math.}\ \textbf{193}, \#1, p.~223--232, 2005.
\texttt{arXiv:math.RA/0312088}

\bibitem{KS}
 M.~Kashiwara, P.~Schapira.
   Categories and sheaves.
Grundlehren der mathematischen Wissenschaften, 332,
Springer, 2006.

\bibitem{Kel}
 B.~Keller.
   Derived categories and their uses.
In: M.~Hazewinkel, Ed., \textit{Handbook of algebra}, vol.~1,
1996, p.~671--701.

\bibitem{KLN}
 B.~Keller, W.~Lowen, P.~Nicol\'as.
   On the (non)vanishing of some ``derived'' categories of curved
dg algebras.
\textit{Journ.\ of Pure and Appl.\ Algebra} \textbf{214}, \#7,
p.~1271--1284, 2010.  \texttt{arXiv:0905.3845 [math.KT]}

\bibitem{Kra}
 H.~Krause.
   The stable derived category of a Noetherian scheme.
\textit{Compositio Math.}\ \textbf{141}, \#5, p.~1128--1162, 2005.
\texttt{arXiv:math.AG/0403526}

\bibitem{Mat}
 E.~Matlis.
   Injective modules over Noetherian rings.
\textit{Pacific Journ.\ of Math.}\ \textbf{8}, \#3, p.~511--528, 1958.

\bibitem{Mat2}
 E.~Matlis.
   The higher properties of $R$\+sequences.
\textit{Journ.\ of Algebra} \textbf{50}, \#1, p.~77--112, 1978.

\bibitem{M-n}
 D.~Murfet.
   Derived categories of quasi-coherent sheaves.
Notes, October~2006.
Available from \texttt{http://www.therisingsea.org/notes}

\bibitem{M-th}
 D.~Murfet.
   The mock homotopy category of projectives and Grothendieck duality.
Ph.~D.\ Thesis, Australian National University, September 2007.
Available from \texttt{http://www. therisingsea.org/thesis.pdf}
{\hbadness=1900\par}

\bibitem{Neem0}
 A.~Neeman.
   The derived category of an exact category.
\textit{Journ.\ of Algebra} \textbf{135}, \#2, p.~388--394, 1990.

\bibitem{Neem}
 A.~Neeman.
   The homotopy category of flat modules, and Grothendieck duality.
\textit{Inventiones Math.}\ \textbf{174}, p.~255--308, 2008.

\bibitem{PSY}
 M.~Porta, L.~Shaul, A.~Yekutieli.
   On the homology of completion and torsion.
\textit{Algebras and Represent.\ Theory} \textbf{17}, \#1, p.~31--67,
2014.  \texttt{arXiv:1010.4386 [math.AC]}.
Erratum in \textit{Algebras and Represent.\ Theory} \textbf{18},
\#5, p.~1401--1405, 2015.  \texttt{arXiv:1506.07765 [math.AC]}

\bibitem{Psemi}
 L.~Positselski.
   Homological algebra of semimodules and semicontramodules:
Semi-infinite homological algebra of associative algebraic structures.
 Appendix~C in collaboration with D.~Rumynin; Appendix~D in
collaboration with S.~Arkhipov.
 Monografie Matematyczne vol.~70, Birkh\"auser/Springer Basel, 2010. 
xxiv+349~pp. \texttt{arXiv:0708.3398 [math.CT]}

\bibitem{Pkoszul}
 L.~Positselski.
   Two kinds of derived categories, Koszul duality, and
comodule-contramodule correspondence.
\textit{Memoirs of the American Math.\ Society} \textbf{212},
\#996, 2011.  vi+133~pp.  \texttt{arXiv:0905.2621 [math.CT]}

\bibitem{Pweak}
 L.~Positselski.
   Weakly curved A${}_\infty$-algebras over a topological local ring.
\textit{M\'emoires de la Soci\'et\'e Math\'ematique de France}
\textbf{159}, 2018.  vi+206~pp.  \texttt{arXiv:1202.2697 [math.CT]}

\bibitem{Pcosh}
 L.~Positselski.
   Contraherent cosheaves on schemes.
Electronic preprint \texttt{arXiv:1209.2995v24 [math.CT]}.

\bibitem{Prev}
 L.~Positselski.
   Contramodules.
\textit{Confluentes Math.}\ \textbf{13}, \#2, p.~93--182, 2021.
\texttt{arXiv:1503.00991 [math.CT]}

\bibitem{Pmgm}
 L.~Positselski.
   Dedualizing complexes and MGM duality.
\textit{Journ.\ of Pure and Appl.\ Algebra} \textbf{220}, \#12,
p.~3866--3909, 2016.  \texttt{arXiv:1503.05523 [math.CT]}

\bibitem{Pcta}
 L.~Positselski.
   Contraadjusted modules, contramodules, and reduced cotorsion modules.
\textit{Moscow Math.\ Journ.}\ \textbf{17}, \#3, p.~385--455, 2017.
\texttt{arXiv:1605.03934 [math.CT]}

\bibitem{Pps}
 L.~Positselski.
   Pseudo-dualizing complexes and pseudo-derived categories.
\textit{Rendiconti Seminario Matematico Univ.\ Padova}
\textbf{143}, p.~153--225, 2020.  \texttt{arXiv:1703.04266 [math.CT]}

\bibitem{Pper}
 L.~Positselski.
   Abelian right perpendicular subcategories in module categories.
Electronic preprint \texttt{arXiv:1705.04960 [math.CT]}.

\bibitem{Pcoun}
 L.~Positselski.
   Flat ring epimorphisms of countable type.
\textit{Glasgow Math.\ Journ.}\ \textbf{62}, \#2, p.~383--439, 2020.
\texttt{arXiv:1808.00937 [math.RA]}

\bibitem{Pproperf}
 L.~Positselski.
   Contramodules over pro-perfect topological rings.
\textit{Forum Mathematicum} \textbf{34}, \#1, p.~1--39, 2022.
\texttt{arXiv:1807.10671 [math.CT]}

\bibitem{Pdc}
 L.~Positselski.
   Remarks on derived complete modules and complexes.
\textit{Math.\ Nachrichten} \textbf{296}, \#2, p.~811--839, 2023.
\texttt{arXiv:2002.12331 [math.AC]}

\bibitem{Psemten}
 L.~Positselski.
   Semi-infinite algebraic geometry of quasi-coherent sheaves
on ind-schemes: Quasi-coherent torsion sheaves, the semiderived
category, and the semitensor product.
Birkh\"auser/Springer Nature, Cham, Switzerland, 2023.  xix+216~pp.
\texttt{arXiv:2104.05517 [math.AG]}

\bibitem{Pedg}
 L.~Positselski.
   Exact DG\+categories and fully faithful triangulated inclusion
functors.
Electronic preprint \texttt{arXiv:2110.08237 [math.CT]}.

\bibitem{Pksurv}
 L.~Positselski.
   Differential graded Koszul duality: An introductory survey.
\textit{Bulletin of the London Math.\ Society} \textbf{55}, \#4,
p.~1551--1640, 2023.  \texttt{arXiv:2207.07063 [math.CT]}

\bibitem{Pal}
 L.~Positselski.
   Local, colocal, and antilocal properties of modules and complexes
over commutative rings.
\textit{Journ.\ of Algebra} \textbf{646}, p.~100--155, 2024.
\texttt{arXiv:2212.10163 [math.AC]}

\bibitem{Pflcc}
 L.~Positselski.
   Flat comodules and contramodules as directed colimits, and cotorsion
periodicity.
\textit{Journ.\ of Homotopy and Related Struct.}\ \textbf{19}, \#4,
p.~635--678, 2024.  \texttt{arXiv:2306.02734 [math.RA]}

\bibitem{Pphil}
 L.~Positselski.
   Philosophy of contraherent cosheaves.
Electronic preprint \texttt{arXiv:2311.14179 [math.AG]}.

\bibitem{PRoos}
 L.~Positselski.
   Roos axiom holds for quasi-coherent sheaves.
\textit{Applied Categorical Struct.}\ \textbf{34}, \#2, article no.~13,
25~pp., 2026.  \texttt{arXiv:2407.13651 [math.AG]}

\bibitem{Pbc}
 L.~Positselski.
   A contramodule generalization of Neeman's flat and projective module
theorem.
\textit{Journ.\ of Homotopy and Related Struct.}\ \textbf{20}, \#4,
p.~477--510, 2025.  \texttt{arXiv:2408.10928 [math.RA]}

\bibitem{Pdomc}
 L.~Positselski.
   $\mathcal D$\+$\Omega$ duality on the contra side.
Electronic preprint \texttt{arXiv:2504.18460v12 [math.AG]}.

\bibitem{Ptd}
 L.~Positselski.
   Torsion modules and differential operators in infinitely many
variables.
\textit{Journ.\ of Pure and Appl.\ Algebra} \textbf{230}, \#4,
article ID~108240, 44~pp., 2026.  \texttt{arXiv:2505.07739 [math.AC]}

\bibitem{PBas}
 L.~Positselski.
   A relative version of Bass' theorem about finite-dimensional
algebras.
\textit{Proceedings of the American Math.\ Society} \textbf{154}, \#6,
p.~2369--2381, 2026.  \texttt{arXiv:2507.15425 [math.RA]}

\bibitem{Pfltp}
 L.~Positselski.
   Coderived and contraderived categories for a cotorsion pair,
flat-type cotorsion pairs, and relative periodicity.
Electronic preprint \texttt{arXiv:2509.07645 [math.CT]}.

\bibitem{Ppt}
 L.~Positselski.
   Pseudo-dualizing complexes of torsion modules and
semi-infinite MGM duality.
Electronic preprint \texttt{arXiv:2511.04571 [math.AC]}.

\bibitem{Pcs}
 L.~Positselski.
   Homomorphisms of topological rings and change-of-scalar functors. 
Electronic preprint \texttt{arXiv:2603.15048 [math.RA]}.

\bibitem{Pns}
 L.~Positselski.
   The contraherent version of the theorem of Sl\'avik and
\v St\!'ov\'\i\v cek.
Electronic preprint \texttt{arXiv:2609.08491 [math.AG]}.

\bibitem{PR}
 L.~Positselski, J.~Rosick\'y.
   Covers, envelopes, and cotorsion theories in locally presentable
abelian categories and contramodule categories.
\textit{Journ.\ of Algebra} \textbf{483}, p.~83--128, 2017.
\texttt{arXiv:1512.08119 [math.CT]}

\bibitem{PSl1}
 L.~Positselski, A.~Sl\'avik.
   Flat morphisms of finite presentation are very flat.
\textit{Annali di Matematica Pura ed Applicata} \textbf{199}, \#3,
p.~875--924, 2020.  \texttt{arXiv:1708.00846 [math.AC]}

\bibitem{PS1}
 L.~Positselski, J.~\v St\!'ov\'\i\v cek.
   The tilting-cotilting correspondence.
\textit{Internat.\ Math.\ Research Notices} \textbf{2021}, \#1,
p.~189--274, 2021.  \texttt{arXiv:1710.02230 [math.CT]}

\bibitem{PS3}
 L.~Positselski, J.~\v St\!'ov\'\i\v cek.
   Topologically semisimple and topologically perfect topological rings.
\textit{Publicacions Matem\`atiques} \textbf{66}, \#2, p.~457--540,
2022.  \texttt{arXiv:1909.12203 [math.CT]}

\bibitem{PS4}
 L.~Positselski, J.~\v St\!'ov\'\i\v cek.
   Derived, coderived, and contraderived categories of locally
presentable abelian categories.
\textit{Journ.\ of Pure and Appl.\ Algebra} \textbf{226}, \#4,
article ID~106883, 2022, 39~pp.  \texttt{arXiv:2101.10797 [math.CT]}

\bibitem{PS6}
 L.~Positselski, J.~\v St\!'ov\'\i\v cek.
   Flat quasi-coherent sheaves as directed colimits, and quasi-coherent
cotorsion periodicity.
\textit{Algebras and Represent.\ Theory} \textbf{27}, \#6,
p.~2267--2293, 2024.  \texttt{arXiv:2212.09639 [math.AG]}

\bibitem{RG}
 M.~Raynaud, L.~Gruson.
   Crit\`eres de platitude et de projectivit\'e: Techniques
de ``platification'' d'un module.
\textit{Inventiones Math.}\ \textbf{13}, \#1--2, p.~1--89, 1971.

\bibitem{Rich}
 T.~Richarz.
   Basics on affine Grassmannians.
Notes available from \texttt{https://www.mathematik. tu-darmstadt.de/media/algebra/homepages/richarz/Notes\_on\_affine\_Grassmannians.pdf}
{\hfuzz=10pt\par}

\bibitem{Sal}
 L.~Salce.
   Cotorsion theories for abelian groups.
\textit{Symposia Math.}\ \textbf{XXIII},
Academic Press, London--New York, 1979, p.~11--32.

\bibitem{SaoSt}
 M.~Saor\'\i n, J.~\v St\!'ov\'\i\v cek.
   On exact categories and applications to triangulated adjoints
and model structures.
\textit{Advances in Math.}\ \textbf{228}, \#2, p.~968--1007, 2011.
\texttt{arXiv:1005.3248 [math.CT]}

\bibitem{Sim}
 A.-M.~Simon.
   Approximations of complete modules by complete big
Cohen--Macaulay modules over a Cohen--Macaulay local ring.
\textit{Algebras and Represent.\ Theory} \textbf{12}, \#2--5,
p.~385--400, 2009.

\bibitem{ST}
 A.~Sl\'avik, J.~Trlifaj.
   Very flat, locally very flat, and contraadjusted modules.
\textit{Journ.\ of Pure and Appl.\ Algebra} \textbf{220}, \#12,
p.~3910--3926, 2016.  \texttt{arXiv:1601.00783 [math.AC]}

\bibitem{Sten}
 B.~Stenstr\"{o}m.
   Rings of quotients. An Introduction to Methods of Ring Theory.
Die Grundlehren der Mathematischen Wissenschaften, Band~217.
Springer-Verlag, New York, 1975.

\bibitem{Sto-ICRA}
 J.~\v St\!'ov\'\i\v cek.
   Exact model categories, approximation theory, and cohomology of
quasi-coherent sheaves.
\textit{Advances in representation theory of algebras}, p.~297--367,
EMS Ser.\ Congr.\ Rep., Eur.\ Math.\ Soc., Z\" urich, 2013.
\texttt{arXiv:1301.5206 [math.CT]}

\bibitem{Sto}
 J.~\v St\!'ov\'\i\v cek.
   Derived equivalences induced by big cotilting modules. 
\textit{Advances in Math.}\ \textbf{263}, p.~45--87, 2014.
\texttt{arXiv:1308.1804 [math.CT]}

\bibitem{TT}
 R.~Thomason, T.~Trobaugh.
   Higher algebraic K-theory of schemes and of derived categories.
\textit{The Grothendieck Festschrift} vol.~3, p.~247--435,
Birkh\"auser, 1990.

\bibitem{Yek0}
 A.~Yekutieli.
   Smooth formal embeddings and the residue complex.
\textit{Canadian Journ.\ of Math.}\ \textbf{50}, \#4, p.~863--896, 1998.
\texttt{arXiv:alg-geom/9510007}

\bibitem{Yek}
 A.~Yekutieli.
   On flatness and completion for infinitely generated modules over
noetherian rings.
\textit{Communicat.\ in Algebra} \textbf{39}, \#11, p.~4221--4245, 2011.
\texttt{arXiv:0902.4378 [math.AC]}

\bibitem{Yek2}
 A.~Yekutieli.
   Flatness and completion revisited.
\textit{Algebras and Represent.\ Theory} \textbf{21}, \#4, p.~717--736,
2018.  \texttt{arXiv:1606.01832 [math.AC]}

\bibitem{Yek3}
 A.~Yekutieli.
   Weak proregularity, derived completion, adic flatness, and prisms.
\textit{Journ.\ of Algebra} \textbf{583}, p.~126--152, 2021.
\texttt{arXiv:2002.04901 [math.AC]}

\bibitem{Yek4}
 A.~Yekutieli.
   Derived complete complexes at weakly proregular ideals.
\textit{Journ.\ of Pure and Appl.\ Algebra} \textbf{229}, \#3,
article ID~107909, 2025, 28~pp.  \texttt{arXiv:2309.01687 [math.AC]}

\bibitem{Zhu}
 X.~Zhu.
   Resolving resolution dimensions.
\textit{Algebras and Represent.\ Theory} \textbf{16}, \#4,
p.~1165--1191, 2013.

\end{thebibliography}
\end{document}